\documentclass[12pt,reqno]{amsart}
\usepackage{amsfonts,amsrefs,latexsym,amsmath, amssymb, mathrsfs, verbatim}
\usepackage{url,color}
\usepackage{pifont}
\usepackage{upgreek}
\usepackage{fancyhdr}
\usepackage{hyperref}
\usepackage{calligra}
\usepackage{marvosym}
\usepackage[percent]{overpic}
\usepackage{pict2e}

\usepackage[hypcap=false]{caption}
\usepackage{xcolor}
\usepackage{wasysym}
\usepackage{accents}

\topmargin - .23 in

\newtheorem{theorem}{Theorem}[section]
\newtheorem{proposition}{Proposition}[section]
\newtheorem{lemma}[proposition]{Lemma}
\newtheorem{corollary}[proposition]{Corollary}

\theoremstyle{definition}
\newtheorem{definition}{Definition}[section]
\newtheorem{remark}{Remark}[section]

\theoremstyle{plain}

\DeclareMathAlphabet{\mathcalligra}{T1}{calligra}{m}{n}
\DeclareFontShape{T1}{calligra}{m}{n}{<->s*[2.2]callig15}{}

\newcommand{\Tboot}{T_{(Boot)}}

\newcommand{\Trandatasize}[1]{\mathring{\updelta}}

\newcommand{\TranminusdatasizeWithFactor}{\mathring{\updelta}_*}

\newcommand{\gt}{\underline{g}}
\newcommand{\gtinverse}{\underline{g}^{-1}}
\newcommand{\gsphere}{g \mkern-8.5mu / }
\newcommand{\ginversesphere}{\gsphere^{-1}}

\newcommand{\Euct}{e}

\newcommand{\gtancomp}{\upsilon}

\newcommand{\mytr}{{\mbox{\upshape{tr}}_{\mkern-2mu \gsphere}}}
\newcommand{\myspacetimetr}{{\mbox{\upshape{tr}}_g}}

\newcommand{\D}{\mathscr{D}}
\newcommand{\angD}{ {\nabla \mkern-14mu / \,} }
\newcommand{\angDarg}[1]{{\angD_{\mkern-3mu #1}}}

\newcommand{\angdiv}{\mbox{\upshape{div} $\mkern-17mu /$\,}}
\newcommand{\angLap}{ {\Delta \mkern-12mu / \, } }

\newcommand{\GdVar}{\upgamma}
\newcommand{\BadVar}{\underline{\upgamma}}

\newcommand{\Fullset}{\mathscr{Z}}
\newcommand{\Tanset}{\mathscr{P}}

\newcommand{\Singletan}{P}

\newcommand{\angdiff}{ {{d \mkern-9mu /} }}
\newcommand{\angdiffarg}[1]{ {d \mkern-9mu /}_{#1} }
\newcommand{\angdiffuparg}[1]{ {d \mkern-9mu /}^{#1} }

\newcommand{\Lie}{\mathcal{L}}
\newcommand{\SigmatLie}{\underline{\mathcal{L}}}
\newcommand{\angLie}{ { \mathcal{L} \mkern-10mu / } }

\newcommand{\Lgeo}{L_{(Geo)}}
\newcommand{\Lunit}{L}
\newcommand{\uLunit}{\underline{L}}
\newcommand{\uLgood}{\breve{\underline{L}}}
\newcommand{\Radunit}{X}
\newcommand{\Rad}{\breve{X}}

\newcommand{\NonRadialRad}{\Xi}
\newcommand{\XiCoordComp}{\upxi}

\newcommand{\GeoAngCoordComp}{\upupsilon}

\newcommand{\CoordAng}{\Theta}
\newcommand{\Timenormal}{N}
\newcommand{\Mult}{T}
\newcommand{\GeoAng}{Y}

\newcommand{\GeoAngFlatRadComponent}{\uprho}

\newcommand{\angJ}{ {\mathscr{J} \mkern-14mu / \, } }

\newcommand{\angpi}{ { \pi \mkern-10mu / }}

\newcommand{\angk}{ { {k \mkern-10mu /} \, } }

\newcommand{\angkdoublearg}[2]{ {{k \mkern-10mu /}_{#1 #2} \, } }
\newcommand{\angkmixedarg}[2]{ {{k \mkern-10mu /}_{#1}^{\ #2} \, } }
\newcommand{\angktriplearg}[3]{ {{k \mkern-10mu /}_{#1 #2}^{#3} \, } }
\newcommand{\angkuparg}[1]{ { {k \mkern-10mu /}^{#1} \, } }

\newcommand{\angG}{ {{G \mkern-12mu /} \, }}
\newcommand{\angGarg}[1]{ {{G \mkern-12mu /}_{\mkern 1mu #1} \, }}

\newcommand{\angGdoublearg}[2]{ {{G \mkern-12mu /}_{#1 #2} \, }}
\newcommand{\angGprime}{ {{ {G'} \mkern-16mu /} \, \, }}
\newcommand{\angGprimearg}[1]{ {{ {G'} \mkern-16mu /}_{\mkern 1mu #1} \, }}

\newcommand{\angGmixedarg}[2]{ {{G \mkern-12mu /}_{#1}^{\ #2} \, }}
\newcommand{\angGnospacemixedarg}[2]{ {{G \mkern-12mu /}_{#1}^{\mkern 2mu #2}}}

\newcommand{\angENMOM}{ {{\enmomtensor \mkern-12mu /} \, }}
\newcommand{\angENMOMarg}[1]{ {{\enmomtensor \mkern-12mu /}_{\mkern 1mu #1} \, }}

\newcommand{\angxi}{ { {\xi \mkern-9mu /}  \, }}

\newcommand{\angxiarg}[1]{ {{\xi \mkern-9mu /}_{#1}  \, }}

\newcommand{\deform}[1]{{^{(#1)} \mkern-1mu \pi}}
\newcommand{\deformarg}[3]{{^{(#1)} \mkern-1mu \pi_{#2 #3}}}
\newcommand{\deformuparg}[3]{{^{(#1)} \mkern-1mu \pi^{#2 #3}}}
\newcommand{\deformmixedarg}[3]{{^{(#1)} \mkern-1mu \pi_{#2}^{\ #3}}}

\newcommand{\angdeform}[1]{{^{(#1)} \mkern-2mu \angpi}}
\newcommand{\angdeformoneformarg}[2]{{^{(#1)} \mkern-2mu \angpi_{#2}}}
\newcommand{\angdeformoneformupsharparg}[2]{{^{(#1)} \mkern-2mu {\pi \mkern-10mu /}_{#2}^{\#}}}
\newcommand{\angdeformarg}[3]{{^{(#1)} \mkern-2mu \angpi_{#2 #3}}}

\newcommand{\Lineproject}{{\Pi \mkern-12mu / } \, }
\newcommand{\Sigmatproject}{\underline{\Pi}}

\newcommand{\vol}{\varpi}
\newcommand{\tvol}{\underline{\varpi}}
\newcommand{\conevol}{\overline{\varpi}}

\newcommand{\spherevol}{\uplambda_{{g \mkern-8.5mu /}}}
\newcommand{\argspherevol}[1]{\uplambda_{{g \mkern-8.5mu /}#1}}

\newcommand{\enzero}{\mathbb{E}}
\newcommand{\flzero}{\mathbb{F}}
\newcommand{\coercivespacetime}{\mathbb{K}}

\newcommand{\totTanmax}[1]{\mathbb{Q}_{#1}}

\newcommand{\coerciveTanspacetimemax}[1]{\mathbb{K}_{#1}}

\newcommand{\upchifullmodarg}[1]{{^{(#1)} \mkern-4mu \mathscr{X}}}
\newcommand{\upchifullmodinhom}{\mathfrak{X}}

\newcommand{\upchipartialmodarg}[1]{{^{(#1)} \mkern-4mu \widetilde{\mathscr{X}}}}
\newcommand{\upchipartialmodinhom}{\widetilde{\mathfrak{X}}}
\newcommand{\upchipartialmodinhomarg}[1]{{^{(#1)} \mkern-2mu \widetilde{\mathfrak{X}}}}

\newcommand{\waveinhom}{\mathfrak{F}}

\newcommand{\inhomleftexparg}[2]{{^{(#1)} \mkern-.5mu #2 }}

\newcommand{\upchipartialmodsourcearg}[1]{\inhomleftexparg{#1}{\mathfrak{B}}}

\newcommand{\Jcurrent}[1]{^{(#1)} \mkern-10mu \mathscr{J}}

\newcommand{\Vplus}[2]{{^{(+)} \mkern-.5mu  \mathcal{V}_{#1}^{#2}}}
\newcommand{\Vminus}[2]{{^{(-)} \mkern-.5mu  \mathcal{V}_{#1}^{#2}}}
\newcommand{\Sigmaplus}[3]{{^{(+)} \mkern-.5mu  \Sigma_{#1;#2}^{#3}}}
\newcommand{\Sigmaminus}[3]{{^{(-)} \mkern-.5mu  \Sigma_{#1;#2}^{#3}}}

\newcommand{\Contwo}{B}

\newcommand{\LateTimeLUnitMu}{\upkappa}

\newcommand{\smoothfunction}{\mathrm{f}}

\newcommand{\ThirdSmoothFunction}{\upeta}

\newcommand{\basicenergyerror}[1]{{^{(#1)} \mkern-.5mu \mathfrak{P}}}
\newcommand{\basicenergyerrorarg}[2]{{^{(#1)} \mkern-.5mu \mathfrak{P}}_{(#2)}}

\newcommand{\Cur}{\mathscr{R}}
\newcommand{\Ric}{\mbox{\upshape{Ric}}}

\newcommand{\enmomtensor}{Q}

\newcommand{\myarray}[2][]{\left(
		\begin{array}{lr}
    	 #1 \\
    	 #2 
     \end{array} \right)}

\numberwithin{equation}{subsection}

\begin{document}
\title{Stable shock formation for nearly simple outgoing plane symmetric waves 
}
\author[
JS,GH,JL,WW]{
Jared Speck$^{* \dagger}$,
Gustav Holzegel$^{** \dagger \dagger}$, 
Jonathan Luk$^{*** \dagger \dagger \dagger}$, 
Willie Wong$^{****}$}

\thanks{$^{\dagger}$JS gratefully acknowledges support from NSF grant \# DMS-1162211,
from NSF CAREER grant \# DMS-1454419,
from a Sloan Research Fellowship provided by the Alfred P. Sloan foundation,
and from a Solomon Buchsbaum grant administered by the Massachusetts Institute of Technology.
}

\thanks{$^{\dagger \dagger}$GH gratefully acknowledges support from a grant from the European Research Council.
}

\thanks{$^{\dagger \dagger \dagger}$JL gratefully acknowledges support from
NSF postdoctoral fellowship \# DMS-1204493.
}

\thanks{$^{*}$Massachusetts Institute of Technology, Cambridge, MA, USA.
\texttt{jspeck@math.mit.edu}}

\thanks{$^{**}$Imperial College, London, UK.
\texttt{g.holzegel@imperial.ac.uk}}

\thanks{$^{***}$Cambridge University, Cambridge, UK
\texttt{jluk@dpmms.cam.ac.uk}}

\thanks{$^{****}$\'{E}cole Polytechnique F\'{e}d\'{e}rale de Lausanne, Lausanne, CH; now at Michigan State University, East Lansing, Michigan, USA. \\
\indent \texttt{wongwwy@member.ams.org}}

\begin{abstract}
In an influential 1964 article,
P. Lax studied
$2 \times 2$ genuinely nonlinear strictly hyperbolic PDE systems
(in one spatial dimension). Using the method
of Riemann invariants, 
he showed that a large set of smooth initial data
lead to bounded solutions whose first spatial derivatives
blow up in finite time,
a phenomenon known as wave breaking.
In the present article, we study the Cauchy problem for
two classes of quasilinear wave equations in two spatial dimensions
that are closely related to the systems studied by Lax.
When the data have one-dimensional symmetry,
Lax's methods can be applied to the wave equations to show that
a large set of smooth initial data lead to wave breaking.
Here we study solutions with initial data that are close, as measured by
an appropriate Sobolev norm, to data belonging to a distinguished subset of 
Lax's data: the data corresponding to simple plane waves.
Our main result is that under suitable relative smallness assumptions,
the Lax-type wave breaking for simple plane waves is stable.
The key point is that we allow the data perturbations to break the symmetry.
Moreover, we give a detailed, constructive description of the 
asymptotic behavior of the solution 
all the way up to the first singularity,
which is a shock driven by the intersection
of null (characteristic) hyperplanes.
We also outline how to extend our results to 
the compressible irrotational Euler equations.
To derive our results, we use Christodoulou's 
framework for studying shock formation 
to treat a new solution regime 
in which wave dispersion is not present.

\bigskip

\bigskip

\end{abstract}

\maketitle

\centerline{\today}

\tableofcontents
\setcounter{tocdepth}{2}

\section{Introduction}
\label{S:INTRO}
In his influential article \cite{pL1964}, Lax showed that
 $2 \times 2$ genuinely nonlinear strictly
hyperbolic PDE systems\footnote{Such systems involve two unknowns in one time and one spatial dimension.} 
exhibit finite-time blowup for a large set of smooth initial data.
His approach was based on the method of Riemann invariants,
which was developed by Riemann himself in his
study \cite{bR1860} of singularity formation in compressible fluid mechanics
in one spatial dimension.
The blowup is of wave breaking type, that is, 
the solution remains bounded but its first derivatives blow up.
Lax's results are by now considered classic
and have been extended in many directions
(see the references in Subsect.~\ref{SS:PREVIOUSWORK}).
In particular, an easy modification of his approach could be used to prove finite-time blowup for solutions 
to various quasilinear wave equations in one spatial dimension: 
under suitable assumptions on the nonlinearities,
one could prove blowup by first writing the wave equation as a first-order 
system in the two characteristic derivatives of the solution
and then applying Lax's methods. In the present article, 
we study the Cauchy problem for two classes of such wave equations in two spatial dimensions,
specifically equations \eqref{E:GEOWAVE} and \eqref{E:NONGEOMETRICWAVE} below.
These equations admit plane symmetric, simple wave solutions that 
blow up in finite time (see Subsect.~\ref{SS:SPSW} for a quick proof). Lax's methods can be used to show that such solutions and their blowup are stable under small perturbations that preserve the one-dimensional plane symmetry. 
Our main result is that, under a suitable hierarchy of smallness-largeness assumptions,
these blowup-solutions are also stable under data perturbations that \emph{break the symmetry}. 
To close our proof, we must derive a sharp description of the blowup that, 
even for data with one-dimensional symmetry, 
provides more information than does Lax's approach.
In Subsect.~\ref{SS:PROOFOVERVIEW},
we explain the set of data covered by our main results in more detail.
See Subsect.~\ref{SS:SUMMARYOFMAINRESULTS} for a summary of the results 
and Theorem~\ref{T:MAINTHEOREM} for the full statement.

For some evolution equations 
in \emph{more} than one spatial dimension
that enjoy special algebraic structure,
short proofs of blowup by contradiction are known;
see Subsect.~\ref{SS:PREVIOUSWORK} for some examples.
In contrast, the typical wave equation that we study
does not have any obvious features
which suggest a short path to proving blowup.
In particular, the equations do not generally
derive from a Lagrangian,
admit coercive conserved quantities,
or have signed nonlinearities.
They do, however, enjoy a key property:
they have \emph{special null structures}
(which are distinct from the well-known \emph{null condition} of S. Klainerman). 
These null structures manifest in several ways, 
including the absence of certain terms in the equations
(as we explain in more detail in the discussion surrounding equation \eqref{E:INTROWAVEEQUATIONFRAMEDECOMPOSED})
as well as the preservation of certain good product structures
under suitable commutations and differentiations of the equations
(as we explain in Subsubsect.~\ref{SSS:DIFFERENCESFROMDISPERSIVEREGIME}).
The null structures are \emph{not visible relative to the standard coordinates}.
Thus, to expose them, we construct a dynamic 
``geometric coordinate system'' 
and a corresponding vectorfield frame\footnote{Our frame \eqref{E:INTROFRAME} is closely related to a null
frame, which is the reason that we use the phrase ``special null structures.''
To obtain what is usually called a null frame, we could replace the
vectorfield $\Rad$ in \eqref{E:INTROFRAME} with the null vectorfield
$\upmu \Lunit + 2 \Rad$. All of our results could be derived by using the
null frame in place of \eqref{E:INTROFRAME}.} 
that are adapted to the characteristics corresponding to the nonlinear flow;
see Subsect.~\ref{SS:PROOFOVERVIEW} for an overview.
We are then able to exploit the null structures 
to give a detailed, constructive description of the singularity, 
which is a shock\footnote{By a ``shock'' 
in a solution to equation \eqref{E:GEOWAVE},
we mean that the singularity is of wave-breaking type;
that is, the solution remains bounded but one of 
its first rectangular coordinate partial derivatives blows up.
By a ``shock'' in a solution to equation \eqref{E:NONGEOMETRICWAVE},
we mean that the solution and its first rectangular coordinate partial derivatives 
remain bounded but one of its second rectangular coordinate partial derivatives 
blows up. Note that in both cases, the metric $g$ remains bounded 
but one of its first rectangular coordinate partial derivatives blows up.
\label{FN:SHOCK}} 
in the regime under study. 
A key feature of the proof is that 
\emph{the solution remains regular relative to the geometric coordinates at the low derivative levels}.
The blowup occurs in the partial derivatives of the solution relative to the standard rectangular coordinates
and is tied to the degeneration of the change of variables map between
geometric and rectangular coordinates; see Subsect.~\ref{SS:PROOFOVERVIEW}
for an extended overview of these issues.

Our approach to proving shock formation is based on an extension of
the remarkable framework of Christodoulou, 
who proved \cite{dC2007} 
detailed shock formation results 
for solutions to the relativistic Euler equations
in irrotational regions of $\mathbb{R}^{1+3}$
(that is, regions with vanishing vorticity)
in a very different solution regime: the small-data dispersive regime.
In that regime, relative to a geometric coordinate system analogous 
to the one mentioned in the previous paragraph, 
the solution enjoys 
time decay\footnote{As in our work here, the blowup in the small-data dispersive regime
occurs in the rectangular coordinate partial derivatives of the solution.} 
at the low derivative levels corresponding to the dispersive nature of 
waves (see Subsubsect.~\ref{SSS:DIFFERENCESFROMDISPERSIVEREGIME} for more details).
The decay plays an important role in controlling various error terms and showing
that they do not interfere with the shock formation mechanisms.
In contrast, in the regime under study here, the solutions do not decay. This
basic feature is tied to the fact that in one spatial dimension, wave equations are essentially transport equations.\footnote{This is also true for many hyperbolic systems in one spatial dimension.}
For this reason, we must develop a new approach to controlling error terms
and to showing that the solution exists long enough for the shock to form;
see Subsubsect.~\ref{SSS:DIFFERENCESFROMDISPERSIVEREGIME} for an overview
of some of the new ideas. As we explain below in more detail, 
a key ingredient in our analysis is the propagation of a two-size-parameter
$\mathring{\upepsilon}-\mathring{\updelta}$ hierarchy all the way up to the shock.
Here and throughout, $\mathring{\updelta} > 0$ is a \emph{not necessarily small} parameter 
that corresponds to the size of derivatives in a direction
that is transversal to the characteristics and 
$\mathring{\upepsilon} \geq 0$ is a \emph{small} parameter
that corresponds to the size of derivatives in directions tangent to the characteristics.
The fact that we are able to propagate the hierarchy 
is deeply tied to the special null structures
mentioned in the previous paragraph.


We can describe the solutions 
that we study as ``nearly simple outgoing plane symmetric solutions.''
By a ``plane symmetric solution,'' we mean one that depends
only on a time coordinate $t \in \mathbb{R}$ and a single rectangular spatial coordinate $x^1 \in \mathbb{R}$.
To study nearly plane symmetric solutions, 
we consider wave equations on spacetimes with topology 
$\mathbb{R} \times \Sigma$, where $t \in \mathbb{R}$ corresponds to time, 
$(x^1,x^2) \in \Sigma := \mathbb{R} \times \mathbb{T}$ corresponds to space, and 
the torus $\mathbb{T} := [0,1)$  
(with the endpoints identified and equipped with the 
usual smooth orientation and with a corresponding
local rectangular coordinate function $x^2$) 
corresponds to the direction that is suppressed in plane symmetry.
We have made the assumption $\Sigma = \mathbb{R} \times \mathbb{T}$
mainly for technical convenience; we expect that suitable wave equations 
on other manifolds could be treated using techniques similar to the ones we use in the present article.
By a ``simple outgoing plane symmetric solution'', 
we mean a special class of plane symmetric solution with only the outgoing (moving to the right) component. Recalling the $\mathring{\upepsilon}-\mathring{\updelta}$ hierarchy that we discussed earlier, in the limit $\mathring{\upepsilon} \to 0$, the solutions that we study reduce to simple outgoing plane symmetric solutions.\footnote{Note that in the analysis of this paper, 
	the solution completely vanishes when $\mathring{\upepsilon} = 0$. 
	However, this additional restriction is not 	
	necessary (see Remark \ref{R:LARGEPSI}).}

The first class of problems that we study is the Cauchy problem for covariant 
wave equations:\footnote{Relative to arbitrary coordinates,
\eqref{E:GEOWAVE} is equivalent to 
$\partial_{\alpha}\left( \sqrt{\mbox{\upshape det} g} (g^{-1})^{\alpha \beta} \partial_{\beta} \Psi \right) = 0$.
\label{FN:COVWAVEOPINARBITRARYCOORDS}}
\begin{subequations}
\begin{align} \label{E:GEOWAVE}
	\square_{g(\Psi)} \Psi
	& = 0, \\
	(\Psi|_{\Sigma_0},\partial_t \Psi|_{\Sigma_0}) & = (\mathring{\Psi},\mathring{\Psi}_0),
\end{align}
\end{subequations}
where $\square_{g(\Psi)}$ denotes the covariant wave operator of the Lorentzian metric $g(\Psi)$
and $(\mathring{\Psi},\mathring{\Psi}_0) \in H_{\Euct}^{19}(\Sigma_0) \times H_{\Euct}^{18}(\Sigma_0)$ 
(see Remarks~\ref{R:SOBOLEVSPACESMULTIPLECOORDINATESYSTEMS} and ~\ref{R:ONTHENUMBEROFDERIVATIVES} just below)
are data with support contained in the compact subset $[0,1] \times \mathbb{T}$ of the initial Cauchy hypersurface
$\Sigma_0 := \lbrace t = 0 \rbrace \simeq \mathbb{R} \times \mathbb{T}$.
Here and throughout,\footnote{See Subsect.~\ref{SS:NOTATION} regarding our conventions for indices, and in particular
for the different roles played by Greek and Latin indices.} 
$\square_{g(\Psi)} \Psi := (g^{-1})^{\alpha \beta}(\Psi) \D_{\alpha} \D_{\beta} \Psi$,
where\footnote{Throughout we use Einstein's summation convention.}  
$\D$ is the Levi-Civita connection of $g(\Psi)$.
We assume that relative to the rectangular coordinates $\lbrace x^{\alpha} \rbrace_{\alpha = 0,1,2}$
(which we explain in more detail in Subsect.~\ref{SS:EQUATIONINRECTANGULARCOMPONENTS}),
we have $g_{\alpha \beta} (\Psi) = m_{\alpha \beta} + \mathcal{O}(\Psi)$, where 
$m_{\alpha \beta} = \mbox{\upshape diag}(-1,1,1)$ is the standard Minkowski metric
and $\mathcal{O}(\Psi)$ is an error term, 
smooth in $\Psi$ and $\lesssim |\Psi|$ in magnitude when $|\Psi|$ is small.
Above and throughout, 
$\partial_0$,
$\partial_1$,
and
$\partial_2$
denote the corresponding rectangular coordinate partial derivatives,
$x^0$ is alternate notation for the time coordinate $t$,
and similarly\footnote{Note that $\partial_t$ is \emph{not} the same as the geometric
coordinate partial derivative $\frac{\partial}{\partial t}$ appearing in equation \eqref{E:LISDDT}
and elsewhere throughout the article.} 
$\partial_t := \partial_0$.
We make further mild assumptions on the nonlinearities 
ensuring that relative to rectangular coordinates,
the nonlinear terms are effectively quadratic and 
\emph{fail to satisfy Klainerman's null condition} \cite{sk1984};
see Subsect.~\ref{SS:EQUATIONINRECTANGULARCOMPONENTS} for the details.

\begin{remark}[\textbf{Our analysis refers to more than one kind of Sobolev space}]
	\label{R:SOBOLEVSPACESMULTIPLECOORDINATESYSTEMS}
	Above and throughout, 
	$H_{\Euct}^N(\Sigma_0)$ 
	denotes the standard $N^{th}$ order Sobolev space 
	with the corresponding norm
	$\| f \|_{H_{\Euct}^N(\Sigma_0)}
			:= 
			\left\lbrace
			\sum_{|\vec{I}| \leq N}
						\int_{\Sigma_0}
							(\partial_{\vec{I}} f)^2
						\, d^2 x
			\right\rbrace^{1/2}
	$,
	where	$\partial_{\vec{I}}$ is a multi-indexed differential
	operator denoting repeated differentiation with respect
	to the rectangular spatial coordinate partial derivative vectorfields
	and $d^2 x$
	is the area form of the standard Euclidean metric $\Euct$ 
	on $\Sigma_0$, which has the form
	$\Euct := \mbox{\upshape diag}(1,1)$ relative to the rectangular coordinates.
	It is important to distinguish these $L^2-$type norms
	from the more geometric ones that we introduce in
	Subect.~\ref{SS:NORMS}; the two kinds of norms drastically differ
	near the shock.
\end{remark}

\begin{remark}[\textbf{On the number of derivatives}]
	\label{R:ONTHENUMBEROFDERIVATIVES}
	Although our analysis is not optimal regarding the number
	of derivatives, we believe that any implementation of our approach 
	requires significantly more derivatives than does
	a typical proof of existence of solutions to 
	a quasilinear wave equation based on energy methods.
	It is not clear to us whether this is a limitation of our approach
	or rather a more fundamental aspect of shock-forming solutions.
	Our derivative count is driven by 
	our energy estimate hierarchy, which is based on 
	a descent scheme in which the high-order energy
	estimates are very degenerate, with slight improvements
	in the degeneracy at each level in the descent. 
	For our proof to work, we must obtain at least several orders
	of non-degenerate energy estimates,
	which requires many derivatives.
	See Subsubsect.~\ref{SSS:INTROENERGYESTIMATES} for more details.
\end{remark}

For convenience, instead of studying the solution in the entire spacetime $\mathbb{R} \times \Sigma$,
we study only the non-trivial future portion of the solution 
that is completely determined by the portion of the data
lying to the right of the straight line 
$\lbrace x^1 = 1 - U_0 \rbrace \cap \Sigma_0$, 
where 
\begin{align} \label{E:FIXEDPARAMETER}
0 < U_0 \leq 1 
\end{align}
is a parameter,
fixed until Theorem~\ref{T:MAINTHEOREM}
(our main theorem),
and the data are non-trivial in the region
$\lbrace 1 - U_0 \leq x^1 \leq 1  \rbrace \cap \Sigma_0 := \Sigma_0^{U_0}$
of thickness $U_0$.
See Figure~\ref{F:REGION} for a picture of the setup,
where the curved null hyperplane portion
$\mathcal{P}_{U_0}^t$ and the flat null hyperplane portion $\mathcal{P}_0^t$ 
in the picture are described in detail in Subsect.~\ref{SS:PROOFOVERVIEW}.
\begin{center}
\begin{overpic}[scale=.2]{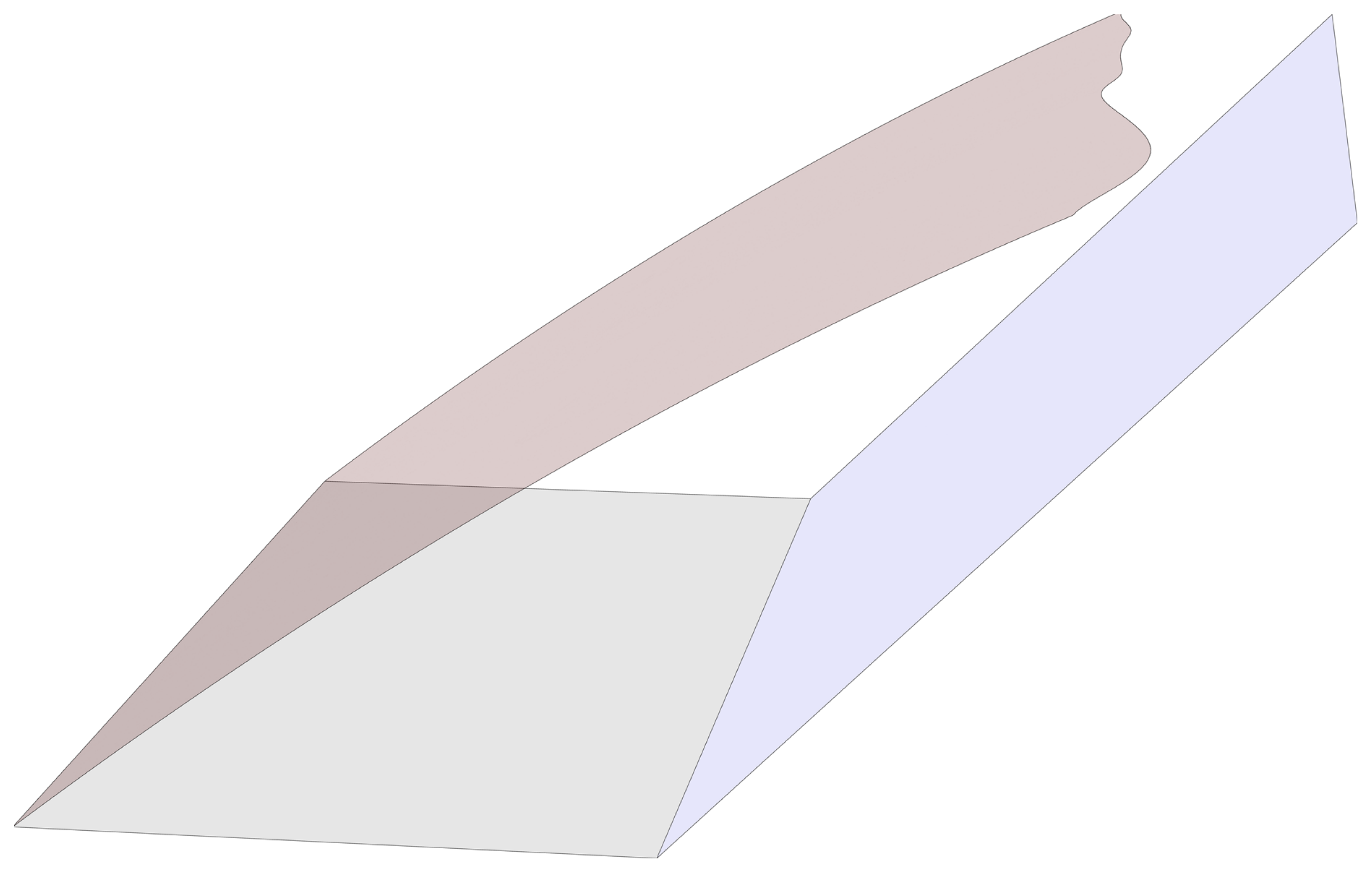} 
\put (37,33) {\large$\displaystyle \mathcal{P}_{U_0}^t$}
\put (74,33) {\large$\displaystyle \mathcal{P}_0^t$}
\put (75,12) {\large$\displaystyle \Psi \equiv 0$}
\put (20,5) {\large$\displaystyle \mbox{non-trivial data}$}
\put (60,5) {\large$\displaystyle \mbox{trivial data}$}
\put (35,18) {\large$\Sigma_0^{U_0}$}
\put (31,10) {\large$\displaystyle U_0$}
%
\put (-.8,16) {\large$\displaystyle x^2 \in \mathbb{T}$}
\put (24,-3.2) {\large$\displaystyle x^1 \in \mathbb{R}$}
\thicklines
\put (-1.1,3){\vector(.9,1){22}}
\put (.5,1.8){\vector(100,-4.5){48}}
\put (10.5,13.9){\line(.9,1){2}}
\put (52.5,11.9){\line(.43,1){1}}
\put (11.5,15){\line(100,-4.5){41.5}}
\end{overpic}
\captionof{figure}{The spacetime region under study.}
\label{F:REGION}
\end{center}

The second class of problems that we study is the Cauchy problem 
for non-covariant wave equations:
\begin{subequations}
\begin{align} \label{E:NONGEOMETRICWAVE}
	(g^{-1})^{\alpha \beta}(\partial \Phi)
	\partial_{\alpha} \partial_{\beta} \Phi
	& = 0,
		\\
	(\Phi|_{\Sigma_0},\partial_t \Phi|_{\Sigma_0}) 
	&= (\mathring{\Phi},\mathring{\Phi}_0),
		\label{E:DATAFORNONGEOMETRICWAVE}
\end{align}
\end{subequations}
where $g(\partial \Phi)$ is a Lorentzian metric with
$g_{\alpha \beta} (\partial \Phi) = m_{\alpha \beta} + \mathcal{O}(\partial \Phi)$.
We assume that the data \eqref{E:DATAFORNONGEOMETRICWAVE} are compactly supported as before, 
but we also assume one extra degree of differentiability: 
$(\mathring{\Phi},\mathring{\Phi}_0) \in H_{\Euct}^{20}(\Sigma_0) \times H_{\Euct}^{19}(\Sigma_0)$.
As we outline in Appendix~\ref{A:EXTENDINGTONONGEOMETRICWAVEEQUATIONS}, 
the second class can essentially be treated 
in the same way as the first class and thus for the remainder of the article, 
we analyze only the first class in detail.

\begin{remark}[\textbf{Special null structure}]
	\label{R:SPECIALNULLSTRUCTURE}
	As we explain in 
	Appendix~\ref{A:EXTENDINGTONONGEOMETRICWAVEEQUATIONS},
	equation \eqref{E:NONGEOMETRICWAVE} \emph{always exhibits} 
	the special null structures mentioned earlier;
	see Lemmas~\ref{L:ECOMMUTEDINTERMSOFHCON} and \ref{L:SPECIALNULLSTRUCTUREINHOMOGENEOUS}.
	However, if the null structures are ``too good,'' 
	then shocks may no longer form;
	see Footnote~\ref{FN:EXCEPTIONALLAGRANGIANS}.
\end{remark}

\subsection{Summary of the main results}
\label{SS:SUMMARYOFMAINRESULTS}
We now summarize our results. See Theorem~\ref{T:MAINTHEOREM} for the precise statement.
We also provide some extended remarks 
and preliminary comparisons to previous work;
see Subsect.~\ref{SS:PREVIOUSWORK} for
a more detailed discussion of some related work.

\begin{changemargin}{.25in}{.25in} \label{ROUGHSTATEMENT}
\noindent \textbf{Rough statement of the main results.} 
Under mild assumptions on the nonlinearities described in Subsect.~\ref{SS:EQUATIONINRECTANGULARCOMPONENTS}, 
there exists an open set\footnote{See Remark~\ref{R:EXISTENCEOFDATA} on pg.~\pageref{R:EXISTENCEOFDATA}
for a proof sketch of the existence of data to which our results apply.} 
(without symmetry assumptions) of compactly supported data 
$(\mathring{\Psi},\mathring{\Psi}_0) \in H_{\Euct}^{19}(\Sigma_0) \times H_{\Euct}^{18}(\Sigma_0)$
for equation \eqref{E:GEOWAVE}
whose corresponding solutions blow up in finite time due to the formation of a shock.
The set contains both large and small data, but each 
pair $(\mathring{\Psi},\mathring{\Psi}_0)$ belonging to the set 
is close to the data corresponding to a plane symmetric simple wave
solution;\footnote{By this,
we mean solutions that are independent of $x^2$ and that are 
constant along a family of null hyperplanes.} 
see Subsects.~\ref{SS:DATAASSUMPTIONS} and \ref{SS:SMALLNESSASSUMPTIONS}
for a precise description of our size assumptions on the data.
Finally, we provide a sharp description of the singularity and the blowup-mechanism.
Similar results hold for equation \eqref{E:NONGEOMETRICWAVE}
for an open set of data contained in $H_{\Euct}^{20}(\Sigma_0) \times H_{\Euct}^{19}(\Sigma_0)$.
\end{changemargin}

\begin{remark}[\textbf{Extending the results to higher spatial dimensions}]
\label{R:EXTENDTOHIGHERDIMENSIONS}
			Our results can be generalized to higher spatial dimensions
			(specifically, to the case of $\Sigma := \mathbb{R} \times \mathbb{T}^n$ for $n \geq 1$) 
			by making mostly straightforward modifications. 
			The only notable difference in higher dimensions is that
			one must complement the energy estimates with elliptic estimates in order to control
			some terms that completely vanish in two spatial dimensions;
			see Remark~\ref{R:ELLIPTIC}. 
\end{remark}

\begin{remark}[\textbf{Maximal development}]
	\label{R:MAXIMALDEVELOPMENT}
	We follow the solution only to the constant-time hypersurface of first blowup.
	However, with modest additional effort, our results
	could be extended to give a detailed description of a portion of the maximal development\footnote{Roughly, the maximal development 
	is the largest possible solution that is uniquely determined by the data; see, for example, \cites{jSb2016,wW2013}
	for further discussion.}
	of the data corresponding to times up to approximately twice the time of first blowup
	(see the discussion below \eqref{E:INTROKEYBLOWUPTRANSVERALFACTOR}), 
	including the shape of the boundary and the behavior of the solution along it.
	More precisely, the estimates that we prove are sufficient for invoking arguments 
	along the lines of those given in \cite{dC2007}*{Ch.~15}, 
	in which Christodoulou provided a description of the maximal development
	(without any restriction on time)
	in the context of small-data solutions to the equations of irrotational relativistic fluid mechanics
	in Minkowski spacetime.
\end{remark}

\begin{remark}[\textbf{The role of} $U_0$]
		\label{R:WHYU0}
		We have introduced the parameter $U_0$
		because one would need to vary it
		in order to extract the information concerning
		the maximal development mentioned in Remark~\ref{R:MAXIMALDEVELOPMENT}.
\end{remark}

\begin{remark}[\textbf{Extending the results to the irrotational Euler equations}]
			Our work can easily be extended to yield a class of 
			stable shock-forming 
			solutions to the irrotational Euler 
			equations (special relativistic or non-relativistic) 
			under almost any\footnote{There is precisely one exceptional equation of state
			for the irrotational relativistic Euler equations to which our results do not apply.
			The exceptional equation of state corresponds to the Lagrangian 
			$\mathscr{L} = 1 - \sqrt{1 + (m^{-1})^{\alpha \beta} \partial_{\alpha} \Phi \partial_{\beta} \Phi}$,
			where $m$ is the Minkowski metric. It is exceptional because it is the only Lagrangian for relativistic fluid mechanics 
			such that Klainerman's null condition is satisfied for perturbations near the constant states with non-zero density.
			A similar statement holds for the non-relativistic Euler equations; see \cite{dCsM2014}*{Subsect.~2.2}
			for more information.
			We note that in \cite{hL2004}, Lindblad showed that in \emph{one or more} spatial dimensions,
			the wave equation corresponding to the Lagrangian
			$\mathscr{L} = 1 - \sqrt{1 + (m^{-1})^{\alpha \beta} \partial_{\alpha} \Phi \partial_{\beta} \Phi}$
			admits global solutions whenever the data are small, smooth, and compactly supported.
			In particular, our approach to proving shock formation certainly does not apply to this equation.
			\label{FN:EXCEPTIONALLAGRANGIANS}} 
			physical equation of state.
			Extending the sharp shock formation results to solutions to the compressible Euler equations
			in regions with non-zero vorticity remains an outstanding open problem.
			The irrotational Euler equations essentially fall under the scope of equation \eqref{E:NONGEOMETRICWAVE}, but
			a few minor changes are needed; we outline them in Appendix~\ref{A:EXTENSIONTOEULER}.
			The main difference is that for the wave equations of fluid mechanics,
			we do not attempt to treat data that have a fluid-vacuum boundary, 
			along which the hyperbolicity of the equations degenerates. 
			Instead, we prove shock formation for perturbations 
			(verifying certain size assumptions)
			of the constant states with \emph{non-zero} density. In terms of a fluid potential $\Phi$,
			the constant solutions correspond to global solutions of the form $\Phi = kt$ with $k > 0$ a constant.
			In Subsect.~\ref{SS:REMARKSONSMALLNESS}, we show that
			there exist data for the irrotational relativistic Euler equations
			verifying the appropriate size assumptions needed to close the proof.
	\end{remark}

	\begin{remark}[\textbf{Additional nonlinearities that we could allow}]
		\label{R:SLIGHTEXTENSIONNOBADTERMS}
		With modest additional effort, our results could also be extended to allow for 
		$g = g(\Phi, \partial \Phi)$ in equation \eqref{E:NONGEOMETRICWAVE} where $g$ is at least linear in $\Phi$.
		That is, we could allow for
		quasilinear terms such as 
		$\Phi \cdot \partial^2 \Phi$.
		Moreover, we could also allow for the presence of semilinear
		terms verifying the \emph{strong null condition}
		(see \cite{jS2014b} for the definition)
		on RHS~\eqref{E:GEOWAVE} or \eqref{E:NONGEOMETRICWAVE}.
		In the regime close to a plane symmetric simple wave, these terms would make only a negligible contribution
		to the dynamics and in particular, they would not interfere with the shock formation processes.
		In contrast, we cannot allow for arbitrary quadratic, cubic, or even higher-order semilinear terms,
		which might highly distort the dynamics in regions where the solution's derivatives becomes large.
\end{remark}

\begin{remark}[\textbf{Possibly allowing $\Psi$ itself to be larger}]
\label{R:LARGEPSI}
For convenience, we assume in our proof that $\Psi$ (undifferentiated) is initially small
(see Subsects.~\ref{SS:DATAASSUMPTIONS} and \ref{SS:SMALLNESSASSUMPTIONS}),
and we show that the smallness is propagated all the way up to the shock.
However, we expect that with effort, one could relax this assumption
by introducing a new parameter corresponding to the $L^{\infty}$ norm of $\Psi$ itself,
which would not have to be ``very small.''
One would of course still have to assume that 
the metric $g(\Psi)$ is initially Lorentzian, which for
some nonlinearities would restrict the allowable size of the new parameter.
One would also have to make the other size assumptions on the data 
stated in Subsects.~\ref{SS:DATAASSUMPTIONS} and \ref{SS:SMALLNESSASSUMPTIONS}
and, in order to ensure that a shock forms, 
that the nonlinearities cause
the factor $G_{\Lunit \Lunit}(\Psi)$
on RHS~\eqref{E:UPMUFIRSTTRANSPORT}
to be non-vanishing.
Moreover, one would have 
to more carefully track the size of $\Psi$
throughout the evolution, especially the influence 
of the new size parameter on the evolution of other quantities. 
This would introduce new technical complications into the proof,
which we prefer to avoid.
\end{remark}

Previous work \cites{dC2007,jS2014b,sA1999a,sA2001a} in more than one spatial dimension, 
which is summarized in the survey article \cite{gHsKjSwW2016}, 
has shown shock formation in solutions to various quasilinear wave equations 
in a different regime: that of solutions
generated by small data supported in a compact subset of $\mathbb{R}^2$
or $\mathbb{R}^3$. Recently, Miao and Yu proved a related large-data
shock formation result \cite{sMpY2014} for a wave equation with cubic nonlinearities
in three spatial dimensions. In Subsect.~\ref{SS:PREVIOUSWORK}, we describe these results and others
in more detail and compare/contrast them to 
our work here. We first provide an overview of our analysis;
we provide detailed proofs starting in Sect.~\ref{S:GEOMETRICSETUP}.

At the close of this subsection,	we would like to highlight some philosophical parallels between our work
here on stable singularity formation and certain global existence results for
the Navier-Stokes equations \cites{
jYCiG2010,yYCmPpZ2014,jYCiGmP2011,yYCpZ2015}
and the Einstein-Vlasov system with a positive cosmological constant \cite{hAhR2016}.
In those works, the authors showed that a class\footnote{In \cites{jYCiG2010,yYCmPpZ2014,jYCiGmP2011,yYCpZ2015}, the symmetric solutions 
are precisely the solutions to the $2D$ Navier-Stokes equations,
which were shown by Leray \cite{jL1933} to be globally regular for data belonging to $L^2$.
In \cite{hAhR2016}, the symmetric solutions included all $\mathbb{T}^3-$Gowdy solutions
and a subset of the $\mathbb{T}^2-$ symmetric solutions (all of which are known to be future-global by \cite{JSm2011}).} 
of global smooth solutions with symmetry can be perturbed in the class of non-symmetric
solutions to produce global\footnote{More precisely, 
the solutions in \cite{hAhR2016} are only shown to be future-global.}
solutions that are approximately symmetric.\footnote{In \cite{jYCiGmP2011}, 
the perturbed solutions are allowed to be far-from-$2D$ in a certain sense, 
though the proof relies on an analyticity assumption on the data.} 
The interesting feature of these results 
is that the symmetric ``background'' solutions are allowed to be large.
Similarly, our results provide a large class of
plane symmetric shock-forming solutions 
that are orbitally stable in the class of
non-symmetric solutions.


\subsection{Overview of the analysis}
\label{SS:PROOFOVERVIEW}
We prove finite-time shock formation for solutions to \eqref{E:GEOWAVE}
for data such that initially,
$\partial_1 \Psi$ is allowed to be of any non-zero size
while,\footnote{Throughout, if $V$ is a vectorfield and $f$ is a scalar function,
then $V f := V^{\alpha} \partial_{\alpha} f$ denotes the derivative of
$f$ in the direction $V$. If $W$ is another vectorfield, then
$V W f := V^{\alpha} \partial_{\alpha} (W^{\beta} \partial_{\beta} f)$,
and similarly for higher-order differentiations.}
roughly speaking,
$\Lunit_{(Flat)} \Psi$ 
and $\partial_2 \Psi$ are relatively small.
Here and throughout,  
$\Lunit_{(Flat)} := \partial_t + \partial_1$
is a vectorfield that is null as measured
by the Minkowski metric: $m_{\alpha \beta} \Lunit_{(Flat)}^{\alpha} \Lunit_{(Flat)}^{\beta} = 0$.
We make similar size assumptions on the higher derivatives at time $0$;
see Subsects.~\ref{SS:DATAASSUMPTIONS} and \ref{SS:SMALLNESSASSUMPTIONS} for the details.

Our assumptions on the nonlinearities 
lead to Riccati-type terms $\sim (\partial_1 \Psi)^2$ in the wave equation \eqref{E:GEOWAVE},
which seem to want to drive $\partial_1 \Psi$
to blow up along the integral curves of $\Lunit_{(Flat)}$.
A caricature of this structure is:
$\Lunit_{(Flat)} \partial_1 \Psi = (\partial_1 \Psi)^2 + \mbox{\upshape Error}$. 
However, our proof does not directly rely on writing the wave equation in this form
or by proving blowup via a Riccati-type argument;
in order to make that kind of argument rigorous, 
one would have to propagate the smallness of the other directional 
derivatives of $\Psi$ 
(found in the term ``$\mbox{\upshape Error}$'')
all the way up to the singularity.
However, the rectangular coordinate partial derivatives are
inadequate for propagating the smallness near 
the singularity in more than one spatial dimension.
In fact, in the regime that we treat here, 
our arguments will suggest that generally,
$\partial_t \Psi$,
$\partial_1 \Psi$,
and
$\partial_2 \Psi$,
all blow up simultaneously
since the rectangular partial derivatives 
are generally transversal to the characteristic surfaces, whose intersection is tied to the blowup. 
These difficulties are not present in simple model problems in one spatial dimension such as Burgers' equation
$\partial_t \Psi + \Psi \partial_x \Psi = 0$;
for Burgers' equation, the blowup of $\partial_x \Psi$ is easy to derive
by commuting the equation with the coordinate derivative $\partial_x$ to obtain a Riccati ODE 
in $\partial_x \Psi$ along characteristics. 

The above discussion has alluded to a defining feature of our proof:
we avoid working with rectangular derivatives and instead
propagate the smallness of \emph{dynamic directional derivatives} of the solution,
tangent to the characteristics,
all the way up to the singularity. 
This allows us to show that the solution's
tangential derivatives do not significantly affect
the shock formation mechanisms,
which are driven by a derivative transversal to the characteristics.
Consequently, in the solution regime under study,
the shock formation mechanisms are essentially the same as in the case of exact plane symmetry.
In particular, there is partial decoupling of the solution's derivatives 
in directions tangent to the characteristics from its transversal derivatives.
We stress that \emph{this effect is not easy to see}. To uncover it, 
we develop an extension of Christodoulou's aforementioned framework \cite{dC2007} for proving shock formation;
see Subsubsect.~\ref{SSS:DIFFERENCESFROMDISPERSIVEREGIME} for a discussion of some of the new ideas that are needed.
The key ingredient in the framework of \cite{dC2007}
is an \emph{eikonal function} $u$, 
which is a solution to the \emph{eikonal equation}.
The eikonal equation is a hyperbolic PDE that depends on the spacetime metric 
$g=g(\Psi)$
and thus on the wave variable.
Specifically, in our study of equation \eqref{E:GEOWAVE},
$u$ solves the eikonal equation initial value problem
\begin{align} \label{E:INTROEIKONAL}
	(g^{-1})^{\alpha \beta}(\Psi) 
	\partial_{\alpha} u \partial_{\beta} u
	& = 0, 
	\qquad \partial_t u > 0,
		\\
	u|_{\Sigma_0}
	& = 1-x^1,
	\label{E:INTROEIKONALINITIALVALUE}
\end{align}
where $(x^1,x^2)$ are the rectangular coordinates\footnote{$x^2$ is only locally defined, 
but this is a minor detail that we typically downplay. 
We note, however, the following fact that we use throughout our analysis:
the corresponding rectangular partial derivative vectorfield $\partial_2$ can be globally defined 
so as to be non-vanishing and smooth relative to the rectangular coordinates.
} 
on $\Sigma_0 \simeq \mathbb{R} \times \mathbb{T}$.
The level sets of $u$ are null (characteristic) hyperplanes for $g(\Psi)$,
denoted by $\mathcal{P}_u$ or by $\mathcal{P}_u^t$ when they are truncated at time $t$.
 We refer to the open-at-the-top region trapped in between 
$\Sigma_0$,
$\Sigma_t$,
$\mathcal{P}_0^t$,
and
$\mathcal{P}_u^t$
as $\mathcal{M}_{t,u}$,
where $\Sigma_t$ denotes the standard flat hypersurface of constant Minkowski time.
We refer to the portion of 
$\Sigma_t$ trapped in between
$\mathcal{P}_0^t$
and
$\mathcal{P}_u^t$
as $\Sigma_t^u$.
The condition \eqref{E:INTROEIKONALINITIALVALUE} implies that 
the trace of the level sets of $u$ along $\Sigma_0$ are 
straight lines,
which we denote by $\ell_{0,u}$. For $t > 0$,
the trace of the level sets of $u$ along $\Sigma_t$
are (typically) curves\footnote{More precisely, the $\ell_{t,u}$ 
are diffeomorphic to the torus $\mathbb{T}$\label{FN:THEELLTUARETORI}.}  
$\ell_{t,u}$.
See Figure~\ref{F:SOLIDREGION} for a picture illustrating these sets and
Def.~\ref{D:HYPERSURFACESANDCONICALREGIONS} for rigorous definitions.

\begin{center}
\begin{overpic}[scale=.2]{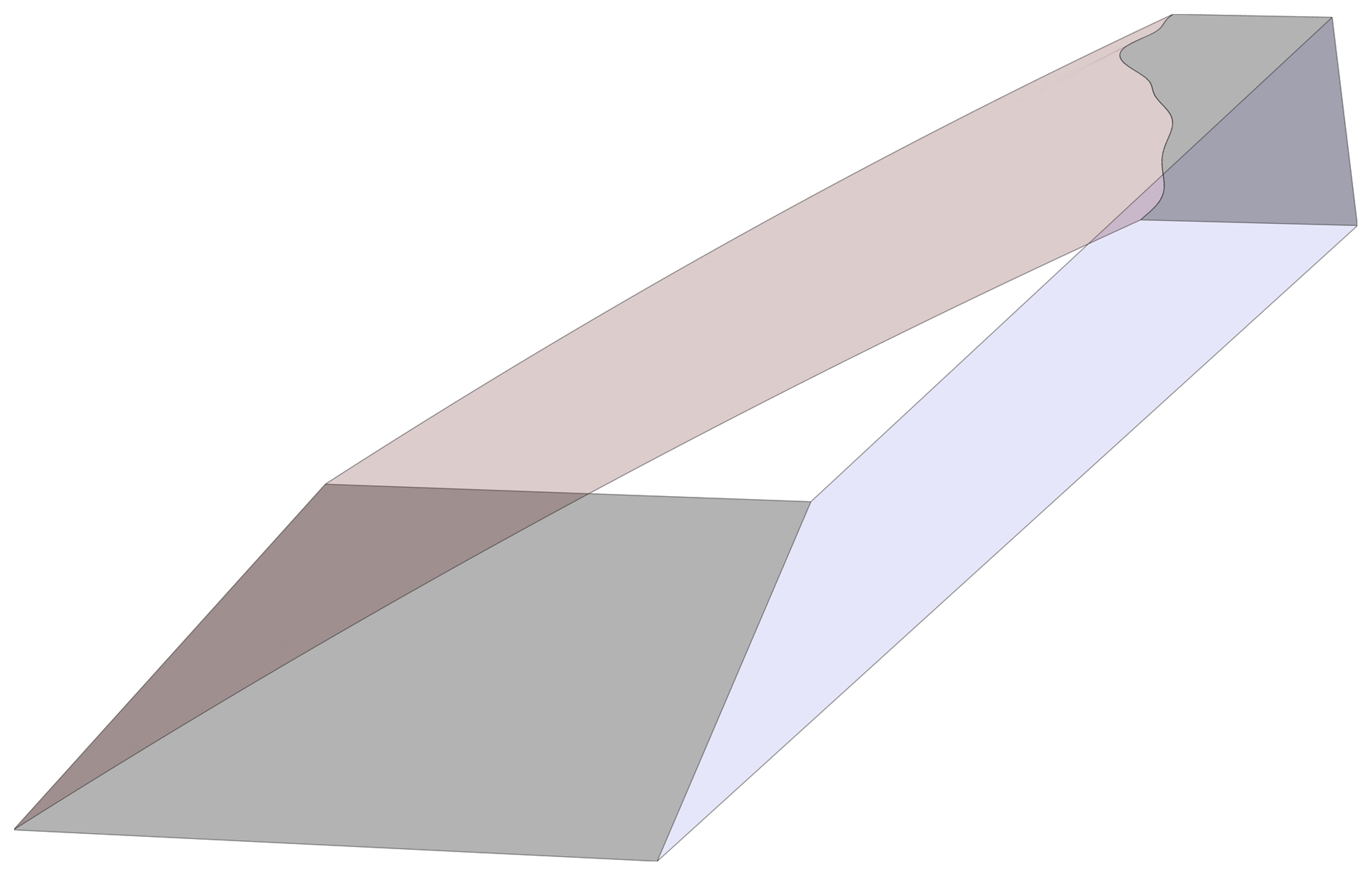} 
\put (54,31) {\large$\displaystyle \mathcal{M}_{t,u}$}
\put (42,33) {\large$\displaystyle \mathcal{P}_u^t$}
\put (74,33) {\large$\displaystyle \mathcal{P}_0^t$}
\put (75,12) {\large$\displaystyle \Psi \equiv 0$}
\put (32,17) {\large$\displaystyle \Sigma_0^u$}
\put (47,13) {\large$\displaystyle \ell_{0,0}$}
\put (12,13) {\large$\displaystyle \ell_{0,u}$}
\put (88.5,57) {\large$\displaystyle \Sigma_t^u$}
\put (93,53) {\large$\displaystyle \ell_{t,0}$}
\put (85.5,53) {\large$\displaystyle \ell_{t,u}$}
\put (-.6,16) {\large$\displaystyle x^2 \in \mathbb{T}$}
\put (22,-3) {\large$\displaystyle x^1 \in \mathbb{R}$}
\thicklines
\put (-.9,3){\vector(.9,1){22}}
\put (.7,1.8){\vector(100,-4.5){48}}
\end{overpic}
\captionof{figure}{The spacetime region and various subsets.}
\label{F:SOLIDREGION}
\end{center}

Eikonal functions $u$ can be viewed as 
coordinates dynamically adapted to the solution via a nonlinear flow.
Their use in the context of proving global results for nonlinear hyperbolic equations in more than one spatial dimension
was pioneered by Christodoulou and Klainerman in their celebrated work \cite{dCsK1993} on the stability of Minkowski spacetime.
Eikonal functions have also been used as central ingredients in proofs of 
low-regularity well-posedness for quasilinear wave equations; 
see, for example, \cites{sKiR2003,hSdT2005,sKiRjS2015,qW2014}.

From $u$, we are able to construct an assortment
of geometric quantities that can be used to derive sharp information about
the solution. The most important of these in the context of shock formation is the \emph{inverse foliation density}
\begin{align} \label{E:FIRSTUPMU}
	\upmu & := 
	\frac{-1}{(g^{-1})^{\alpha \beta}(\Psi) \partial_{\alpha} t \partial_{\beta} u} > 0,
\end{align}
where $t$ is the rectangular time coordinate. 
The quantity $1/\upmu$ measures the density of the level sets of $u$
relative to the constant-time hypersurfaces $\Sigma_t$. In our work here,
$\upmu$ is initially close to $1$ and when it vanishes,
the density becomes infinite and the level sets of $u$
(the characteristics) intersect;
see Figure~\ref{F:FRAME} below, in which we illustrate a scenario 
where $\upmu$ has become small and a shock is about to form.
In the solution regime under study, we prove that 
the rectangular components $g_{\alpha \beta}$ remain
near those of the Minkowski metric $m_{\alpha \beta} = \mbox{\upshape diag}(-1,1,1)$ all the way up to the shock. 
Thus, from \eqref{E:FIRSTUPMU}, we infer that the vanishing of $\upmu$ 
implies that some rectangular derivative of
$u$ blows up. From experience with model equations in one spatial dimension
such as Burgers' equation,
one might expect that the intersection of the characteristics
is tied to the formation of a singularity in $\Psi$.
Though it is not obvious, 
our proof in fact reveals that in the regime under study, 
$\upmu = 0$ corresponds to the blowup of the first\footnote{For equation \eqref{E:NONGEOMETRICWAVE}, the blowup occurs in the second
rectangular derivatives of $\Phi$.} rectangular derivatives of $\Psi$.
In particular, on sufficiently large time intervals,
our work affords 
\emph{a sharp description of singularity formation characterized precisely by the vanishing of $\upmu$.}

Our analysis relies on the \emph{geometric coordinates} $(t,u,\vartheta)$,
where $t = x^0$ and $u$ are as above and $\vartheta$ solves the evolution equation
$-(g^{-1})^{\alpha \beta}(\Psi) \partial_{\alpha} u \partial_{\beta} \vartheta = 0$
with $\vartheta|_{\Sigma_0} = x^2$, where $x^2$ is the local rectangular coordinate on $\mathbb{T}$.
The most important feature of the geometric coordinates is that relative to them,
the shock singularity is \emph{renormalizable}, with the possible exception of the high derivatives.\footnote{The possibility 
that the high derivatives might behave worse is a fundamental difficulty that permeates our analysis.}
More precisely, we show that
the solution and its up-to-mid-order geometric derivatives
(that is, the geometric partial derivatives
$
\frac{\partial}{\partial t}
$,
$
\frac{\partial}{\partial u}
$,
and
$
\frac{\partial}{\partial \vartheta}
$)
remain bounded in $L^{\infty}$ all the way up to the shock.
In particular, \emph{the solution's first derivatives relative to the geometric coordinates do not blow up}!
The blowup of the solution's first \emph{rectangular} partial derivatives 
is a ``low-level'' effect that could be obtained\footnote{We use a slightly different, more direct argument to prove the blowup; see Subsubsect.~\ref{SSS:INTROPOOFTHATMUVANISHES} for an overview.} 
by transforming back to the rectangular coordinates and showing
that $\upmu = 0$ causes a degeneracy in the change of variables (see Lemma~\ref{L:CHOV}).

As we alluded to in Remark~\ref{R:ONTHENUMBEROFDERIVATIVES},
the new feature that makes the proof of shock formation
more difficult than typical global results for wave equations is:
at the very high orders, 
our energies are allowed to blow up 
like $(\min_{\Sigma_t^u} \upmu)^{-p}$
as $\upmu \to 0$, where $p$ is a constant depending on the order of the energy;
see Subsubsect.~\ref{SSS:INTROENERGYESTIMATES} for an overview.
An important aspect of our proof is that 
\emph{the blowup-exponents $p$ are controlled} by certain universal\footnote{These constants are the same for all of the wave equations that we study in this article.} structural constants appearing in the equations.
The main contribution of Christodoulou in \cite{dC2007} 
was showing how to derive the degenerate 
high-order energy estimates and, crucially, 
\emph{proving that the degeneracy does not propagate 
down to the low orders.} These steps consume the majority of our effort here.

To derive estimates, 
rather than working with the geometric coordinate partial derivative frame,
we instead replace $\frac{\partial}{\partial u}$ with a 
similar vectorfield $\Rad$ that has slightly better geometric properties, 
which we describe below;
see Def.~\ref{D:RADANDXIDEFS} for the details of the construction.
That is, we rely on the following dynamic vectorfield frame,
which is depicted at two distinct points along
a fixed null hyperplane portion $\mathcal{P}_u^t$ in Figure~\ref{F:FRAME}:
\begin{align} \label{E:INTROFRAME}
	\left\lbrace
		\Lunit, \Rad, \CoordAng
	\right\rbrace.
\end{align}
The vectorfield $\Lunit =  \frac{\partial}{\partial t}$ 
is a null (that is, $g(\Lunit, \Lunit) = 0$) generator of 
$\mathcal{P}_u$ (in particular, $\Lunit$ is $\mathcal{P}_u-$tangent)
and $\CoordAng := \frac{\partial}{\partial \vartheta}$ 
is $\ell_{t,u}-$tangent with $g(\Lunit,\CoordAng) = g(\Rad,\CoordAng) = 0$.
Relative to the rectangular coordinates, we have
\begin{align}	 \label{E:LUNITGRADIENTVECTORFIELDINTRO}
	\Lunit^{\alpha}
	& = - \upmu (g^{-1})^{\alpha \beta} \partial_{\beta} u.
\end{align}
Our proof shows that all the way up to the shock,
$\Lunit$ and $\CoordAng$ remain close to their flat analogs,
which are respectively $\Lunit_{(Flat)} := \partial_t + \partial_1$
and $\partial_2$. The vectorfield $\Rad$ is transversal to $\mathcal{P}_u$,
$\Sigma_t-$tangent, $g-$orthogonal to $\ell_{t,u}$,
and, most importantly, normalized by $g(\Rad,\Rad) = \upmu^2$. 
In particular, the rectangular
components $\Rad^{\alpha}$ vanish precisely at the points where $\upmu$ vanishes (that is, at the shock points).
Our proof shows that $\Rad$ remains near $- \upmu \partial_1$ all the way up to the shock.
This is depicted in Figure~\ref{F:FRAME}, in which the vectorfield 
$\Rad$ is small in the region up top where $\upmu$ is small.

\begin{center}
\begin{overpic}[scale=.35]{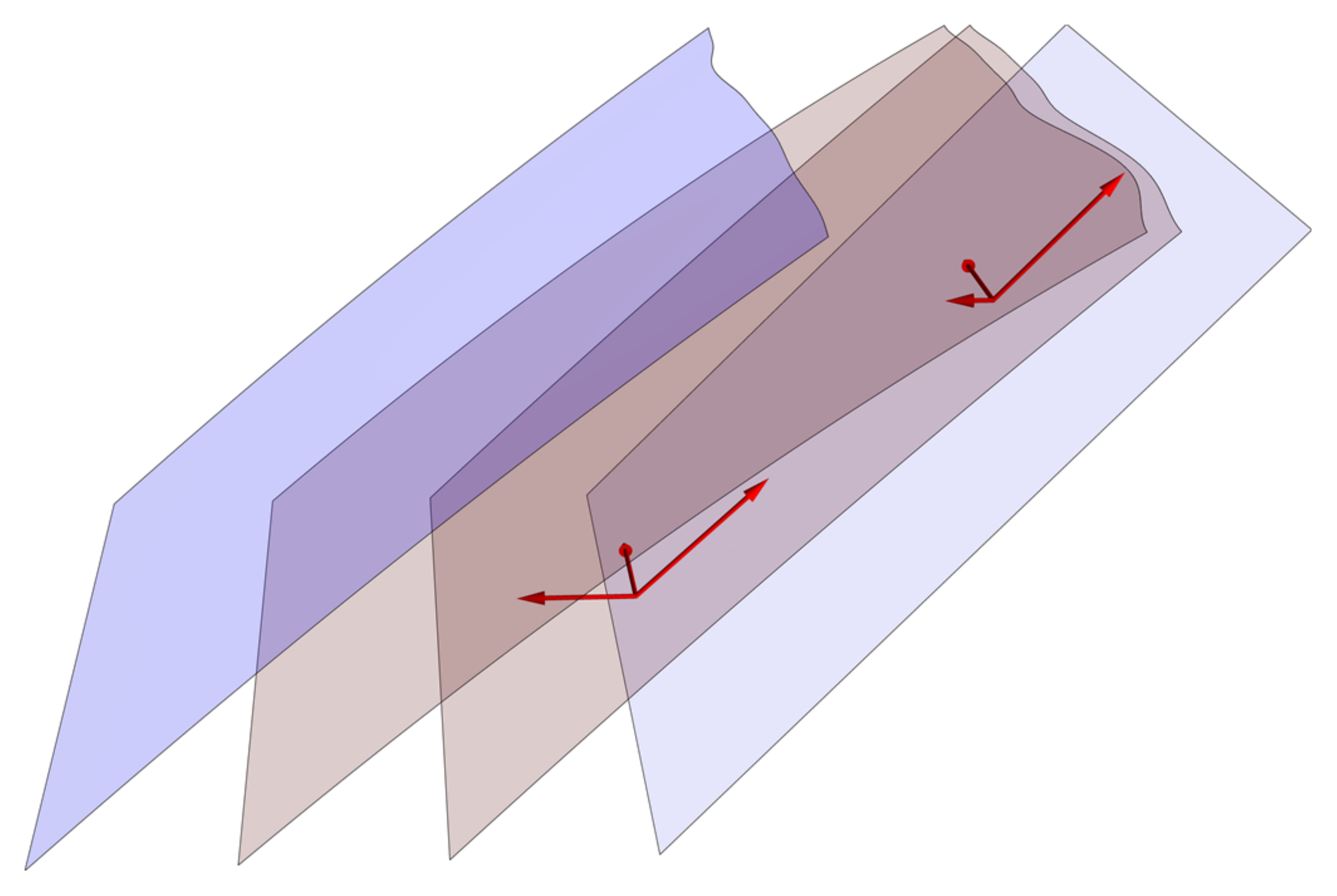} 
\put (84,55) {\large$\displaystyle \Lunit$}
\put (67,43.5) {\large$\displaystyle \Rad$}
\put (70.5,48.5) {\large$\displaystyle \CoordAng$}
\put (57,32) {\large$\displaystyle \Lunit$}
\put (35,21) {\large$\displaystyle \Rad$}
\put (45,27) {\large$\displaystyle \CoordAng$}
\put (51,13) {\large$\displaystyle \mathcal{P}_0^t$}
\put (37,13) {\large$\displaystyle \mathcal{P}_u^t$}
\put (7,13) {\large$\displaystyle \mathcal{P}_1^t$}
\put (21,10) {\large$\displaystyle \upmu \approx 1$}
\put (70,58) {\large$\displaystyle \upmu \ \mbox{\upshape small}$}
\put (75,12) {\large$\displaystyle \Psi \equiv 0$}
%
\end{overpic}
\captionof{figure}{The dynamic vectorfield frame at two distinct points in 
$\mathcal{P}_u^t$, where $0 < u < 1$.}
\label{F:FRAME}
\end{center}

Throughout the paper, we often depict $\mathcal{P}_u-$tangent derivative operators 
such as $\Lunit$ and $\CoordAng$ with the symbol $\Singletan$.
The main idea of our paper is to treat a regime in which the initial data have
pure transversal derivatives such as
$\Rad \Rad \Psi$ 
and $\Rad \Psi$ that are  
of size $\approx \mathring{\updelta} > 0$,
while all other derivatives such as 
$\Singletan \Rad \Psi$,
$\Singletan \Psi$,  
and $\Psi$ itself
are of \emph{small} size $\mathring{\upepsilon}$.
The quantity $\mathring{\updelta}$ can be either small or large,
but our required smallness of $\mathring{\upepsilon}$
depends on $\mathring{\updelta}$;
see Subsects.~\ref{SS:DATAASSUMPTIONS} and \ref{SS:SMALLNESSASSUMPTIONS} 
for the precise assumptions. 
Similar remarks apply to 
$\upmu$ and to the rectangular component functions $\Lunit^{\alpha}$ at time $0$.
To avoid lengthening the paper, 
we generally do not closely track 
the dependence of our estimates on 
$\mathring{\updelta}$. In particular, 
as we explain in Subsect.~\ref{SS:NOTATION},
we allow the ``constants'' $C$ appearing in the estimates to depend 
on $\mathring{\updelta}$.
There is one crucially important exception: 
we carefully track the dependence 
of a handful of important estimates on 
a quantity $\TranminusdatasizeWithFactor$
that is related to $\mathring{\updelta}$
and that controls the blowup-time:
\begin{align} \label{E:INTROKEYBLOWUPTRANSVERALFACTOR}
\TranminusdatasizeWithFactor
:= \frac{1}{2} \sup_{\Sigma_0^1} \left[G_{\Lunit \Lunit} \Rad \Psi \right]_-
	> 0
\end{align}
(see Def.~\ref{D:CRITICALBLOWUPTIMEFACTOR}),
where 
$
\displaystyle
G_{\Lunit \Lunit} := \frac{d}{d \Psi}g_{\alpha \beta}(\Psi) \Lunit^{\alpha} \Lunit^{\beta}
$
and $f_- = |\min \lbrace f, 0 \rbrace|$.
We explain the connection between $\TranminusdatasizeWithFactor$ and the blowup-time in
Subsubsect.~\ref{SSS:INTROPOOFTHATMUVANISHES}.
In our proof, we show that we can propagate the $\mathring{\upepsilon}-\mathring{\updelta}$ hierarchy
(in various norms)
all the way up to the time of first shock formation, 
which we show is 
$
\left\lbrace
	1 + \mathcal{O}(\mathring{\upepsilon})
\right\rbrace
\TranminusdatasizeWithFactor^{-1}
$.
We give an example of this kind of propagation in
Subsubsect.~\ref{SSS:DIFFERENCESFROMDISPERSIVEREGIME}.
In practice, when proving estimates via a bootstrap argument,
we give ourselves a margin of error by showing that we
could propagate the hierarchy for classical solutions existing up to time $2 \TranminusdatasizeWithFactor^{-1}$,
which is plenty of time for the shock to form.
Actually, our results show something stronger: \emph{no other singularities besides shocks
can form for times} $\leq 2 \TranminusdatasizeWithFactor^{-1}$.
The factor of $2$ in the previous inequality is not important 
and could be replaced with any
positive constant larger than $1$, but we would have to 
further shrink the allowable size of $\mathring{\upepsilon}$
as the size of the constant increases.

One important reason why we are able to propagate the hierarchy 
for times up to $2 \TranminusdatasizeWithFactor^{-1}$
is: relative to the frame \eqref{E:INTROFRAME}, 
the wave equation $\square_{g(\Psi)} \Psi = 0$ 
has a \textbf{miraculous structure}. Specifically, 
$\upmu \square_{g(\Psi)} \Psi = 0$
is equivalent to (see Prop.~\ref{P:GEOMETRICWAVEOPERATORFRAMEDECOMPOSED}) 
\begin{align} \label{E:INTROWAVEEQUATIONFRAMEDECOMPOSED}
	- \Lunit(\upmu \Lunit \Psi + 2 \Rad \Psi)
	+ \upmu \angLap \Psi
	& = \mathcal{N},
\end{align}
where $\angLap$ denotes the covariant Laplacian induced by $g$ along the curves $\ell_{t,u}$
and $\mathcal{N}$ denotes quadratic terms depending on $\leq 1$ derivatives of $\Psi$ and $\leq 2$ derivatives of $u$
with the following critically important null structure: 
\emph{each product in $\mathcal{N}$ contains at least one good $\mathcal{P}_u-$tangent differentiation
and thus inherits a smallness factor of $\mathring{\upepsilon}$}.
In particular, products containing quadratic or higher powers of pure transversal derivatives 
(such as $(\Rad \Psi)^2$, $(\Rad \Psi)^3$, etc.) 
are \emph{completely absent}. This good structure is related to Klainerman's null condition,
but unlike in his condition, the structure of the cubic and higher-order terms matters.
Another way to think about \eqref{E:INTROWAVEEQUATIONFRAMEDECOMPOSED} is:
by bringing $\upmu$ under the outer $\Lunit$ differentiation,
we have generated a product term of the form $- (\Lunit \upmu) \cdots$.
This leads to the cancellation of the worst term on the RHS, 
which was proportional to $\upmu^{-1} (\Rad \Psi)^2$. Put differently, 
the term $\frac{1}{2} G_{\Lunit \Lunit} \Rad \Psi$
from the RHS of equation \eqref{E:INTROMUEVOLUTION} below 
generates complete, nonlinear cancellation
of a term proportional to $\upmu^{-1} (\Rad \Psi)^2$.
This null structure survives under commutations of the wave equation with vectorfields
adapted to the eikonal function
and \emph{allows us to propagate the smallness of the 
size $\mathring{\upepsilon}$ quantities 
even though the size $\mathring{\updelta}$ quantities are allowed to be much larger}.

Our strategy of propagating the smallness of some quantities while 
simultaneously allowing derivatives transversal to the characteristics to be large
has roots in the similar approach taken 
by Christodoulou \cite{dC2009} in his celebrated proof of the formation of trapped surfaces
in solutions to the Einstein-vacuum equations
and in the related works \cites{sKiR2010b,sKiR2012,jL2013,jLiR2013,sKjLiR2014,xAjL2014,jLiR2015}.
Similar strategies have been used \cites{pYjW2016,jWpY2013,sMpLpY2014,sY2015} to prove global existence
results for semilinear wave equations verifying the null condition
in regimes that allow for large transversal derivatives.

\subsection{A short proof of blowup for plane symmetric simple waves}
\label{SS:SPSW}
We now illustrate the strategy discussed in Subsect.~\ref{SS:PROOFOVERVIEW}
by studying a model problem. Specifically, we
explain how to prove blowup for simple wave solutions 
(which we explain below)
to equation \eqref{E:GEOWAVE} in one spatial dimension.
Strictly speaking, such solutions are not covered by our main theorem
(Theorem~\ref{T:MAINTHEOREM}), but nonetheless, our model problem
provides the main idea behind the easy part of the proof of the shock formation
and the role of the smallness of the data-size parameter $\mathring{\upepsilon}$
from \eqref{E:PSIDATAASSUMPTIONS}.
That is, the solutions treated in our main theorem may be viewed
as small perturbations of solutions that are analogous to the ones
treated in this subsection.
Note that there is a difference between\footnote{In particular,
$\mbox{\upshape det} g$ depends on
the coefficients of the metric corresponding to the ``extra spatial dimensions."}
imposing plane symmetry on solutions to 
\eqref{E:GEOWAVE} in the case of two spatial dimensions
and studying equation \eqref{E:GEOWAVE} in one spatial dimension.
However, this difference is minor (as we explain at the end of this subsection)
and can be ignored here.

Specifically, we start by considering wave equations of the form
\[ \square_{g(\Psi)} \Psi = 0 \]
on $\mathbb{R}^{1+1}$.
Throughout this subsection, we denote the standard rectangular coordinates
on $\mathbb{R}^{1+1}$ by $(x^0,x^1)$. We sometimes use the alternate notation 
$(t,x) = (x^0,x^1)$.
We assume that the rectangular components of the metric verify
$g_{\alpha \beta} = g_{\alpha \beta}(\Psi) = m_{\alpha \beta} + \mathcal{O}(\Psi)$.
Here $m_{\alpha \beta} = \mbox{\upshape diag}(-1,1)$ is the standard Minkowski metric.
We assume that the data
$(\Psi|_{t=0},\partial_t \Psi|_{t=0})$
are supported in the unit interval $[0,1]$.
Also, for convenience, we make the assumption
\eqref{E:GINVERSE00ISMINUSONE}.
All of these assumptions could be significantly weakened
or eliminated, but we do not pursue those issues here.

We now let $u$ and $v$ be a pair of eikonal functions 
that increase towards the future such that
the level sets of $u$ are transversal to those of $v$.
That is, $u$ and $v$ are solutions to 
\[ (g^{-1})^{\alpha\beta}(\Psi) \partial_\alpha u \partial_\beta u 
= 
0 
= (g^{-1})^{\alpha\beta}(\Psi) \partial_\alpha v \partial_\beta v
\]
such that $\partial_t u, \partial_t v > 0$
and such that $d u$ and $d v$ are linearly independent.
For convenience, we choose 
the initial conditions $u|_{t=0} = 1 - x$,
as in \eqref{E:INTROEIKONAL}.
We also set $v|_{t=0} = x$ to be concrete.
As long as $(u,v)$ do not degenerate, we may use them as ``null coordinate'' functions 
in place of $(t,x)$.
We denote the corresponding
coordinate partial derivative vectorfields by
$
\partial_u,
\partial_v
$.

In two (spacetime) dimensions,
$g$ can be written, relative to the null
coordinates, as
$g = - \Omega^2 (du \otimes dv + dv \otimes du)$,
where $\Omega$ is a scalar-valued function.
It follows 
(see Footnote~\ref{FN:COVWAVEOPINARBITRARYCOORDS}) 
that the covariant wave equation $\square_{g(\Psi)} \Psi = 0$ is equivalent to 
\[ 
\partial_u 
\partial_v 
\Psi = 0,\]
where the nonlinearity is 
``hidden'' in the definition of $u,v$ above. 
Thus, we infer that the condition 
$
\partial_v \Psi = 0
$
is propagated by the solution if it is verified by the initial data.
We refer to such a solution $\Psi = \Psi(u)$ as a \emph{simple wave}. 
Note that the simple-wave-initial-data-assumption 
may be compared with \eqref{E:PSIDATAASSUMPTIONS} 
with $\mathring{\upepsilon} = 0$. However, we make the minor remark that the comparison
is not perfect because according to our definitions,\footnote{See Remark \ref{R:LARGEPSI} for related discussions.}
$\mathring{\upepsilon} = 0$ implies that $\Psi \equiv 0$.

For simple waves, 
$\Psi$ is constant along the level sets of $u$
and hence so are the rectangular components 
$(g^{-1})^{\alpha\beta} = (g^{-1})^{\alpha \beta}(\Psi(u))$. 
It follows that, when graphed in the
$(t,x)$ plane, the level sets of $u$ are straight lines
\emph{which are not generally parallel}.\footnote{Note that this is the same behavior seen in the characteristics associated to solutions of Burgers' equation.} 
Thus, if the characteristic velocities 
(that is, the ``slopes'' of the level sets of $u$) 
are initially not constant across different values of $u$, 
then from the compactness of the support of the data, 
we conclude that there must exist two distinct level sets of $u$ that intersect in finite time.
Clearly the rectangular derivatives 
$\partial_\beta u$ must blow up at the intersection points.
As we described in Subsect.~\ref{SS:PROOFOVERVIEW},
at such intersection points, the quantity $\upmu$ defined in \eqref{E:FIRSTUPMU}
tends to $0$. Below we explain why the vanishing of $\upmu$
is connected to the blowup of a first derivative of $\Psi$
in the direction of a vectorfield with length of order $1$.

We now compute the blowup-time by examining the quantity $1/\upmu$. 
Our goal is to explain why the blowup-time is tied
to the quantity $\TranminusdatasizeWithFactor$ 
defined in \eqref{E:INTROKEYBLOWUPTRANSVERALFACTOR}.
To this end, we define the vectorfield
$\Lunit$ as in \eqref{E:LUNITGRADIENTVECTORFIELDINTRO}
and the vectorfield $\Rad = \upmu \Radunit$ as in \eqref{E:RADDEF}.
Note that in the present context, 
$\Lunit$ 
is a scalar function multiple of 
$
\partial_v 
$.
Note also that since $\Lunit t = 1$ (see \eqref{E:LUNITOFUANDT})
and since $\Lunit$ is parallel to 
the straight line characteristics (in the $(t,x)$ plane),
it follows that $\Lunit^{\alpha} = \Lunit^{\alpha}(u)$ for $\alpha = 0,1$.
From \eqref{E:DOWNSTAIRSUPSTAIRSSRADUNITPLUSLUNITISAFUNCTIONOFPSI}
and the above discussion, we also see that
$\Radunit^{\alpha} = \Radunit^{\alpha}(u)$ for $\alpha = 0,1$.
Just below, we will derive the following evolution equation,
valid for simple waves:
\begin{align} \label{E:KEYEXPRESSIONINT}
	\Lunit \upmu 
	& = \frac{1}{2} G_{\Lunit \Lunit} \Rad \Psi,
\end{align}
where 
$
G_{\Lunit \Lunit} := \frac{d}{d \Psi}g_{\alpha \beta}(\Psi) \Lunit^{\alpha} \Lunit^{\beta}.
$
Recalling that $\Lunit = \frac{\partial}{\partial t}|_u$,
it is now clear that $\TranminusdatasizeWithFactor^{-1}$ is connected to 
the time of first vanishing of $\upmu$ (the blowup-time),
as we described in Subsect.~\ref{SS:PROOFOVERVIEW}.

We now explain why a first derivative of $\Psi$ blows up when
$\upmu$ vanishes.
To this end, we note that 
$\Lunit \upmu = \frac{1}{2} \upmu G_{\Lunit \Lunit} \Radunit \Psi$
and that by \eqref{E:RADIALVECTORFIELDSLENGTHS}, 
$g(\Radunit,\Radunit) = 1$. In particular, 
$\Radunit \Psi$ is a derivative of $\Psi$
with respect to a vectorfield of strictly positive length.
Moreover, from the above discussion, we see that 
$G_{\Lunit \Lunit}$ is constant
along the integral curves of $\Lunit$ (that is, $G_{\Lunit \Lunit} = G_{\Lunit \Lunit}(u)$).
It follows that if $\upmu$ goes to $0$ in finite time, then $|\Radunit \Psi|$ 
\emph{must blow up}.

To complete our analysis in this subsection, we will derive \eqref{E:KEYEXPRESSIONINT}.
To this end, we differentiate \eqref{E:FIRSTUPMU} to derive the following identity, 
which relies on the facts that the rectangular derivatives $\partial_\alpha t$ are constant, 
and that, by the above discussion, 
$(g^{-1})^{\alpha\beta}$ is constant along the lines of constant $u$:
\[ 
	\Lunit \upmu^{-1}
	:= 
	\Lunit^\alpha \partial_\alpha (\upmu^{-1}) = \upmu (g^{-1})^{\beta\gamma} 
	\partial_\beta t (g^{-1})^{	\alpha\delta}    
	\partial_\delta u \partial_{\alpha} \partial_{\gamma} u.
\]
Differentiating the eikonal equation 
$(g^{-1})^{\alpha\beta}(\Psi) \partial_\alpha u \partial_\beta u = 0$,
we obtain 
\[ 
	\Lunit \upmu 
	= \frac{1}{2} \upmu^3 (g^{-1})^{\beta\gamma} \partial_\beta t \partial_\alpha u \partial_\delta u \partial_\gamma (g^{-1})^{\alpha\delta},
\]
which we can simplify to 
\[ \Lunit \upmu 
= - \frac{1}{2} \upmu (g^{-1})^{\beta \gamma} \partial_\beta t \partial_\gamma \Psi G_{\Lunit \Lunit}.
\]
In the expression above, the vectorfield 
$- (g^{-1})^{\beta \gamma}\partial_\beta t$ 
is equal to $\Timenormal^\gamma$ (see \eqref{E:TIMENORMALRECTANGULAR}), 
where $\Timenormal$ is the future-directed unit normal to $\Sigma_t$.
Hence, from \eqref{E:TIMENORMAL} and the fact that $\Lunit \Psi = 0$
for simple plane waves, we obtain 
the desired key expression \eqref{E:KEYEXPRESSIONINT}.
This completes our discussion of blowup for simple plane waves.


We close this subsection by noting that similar analysis can be applied to plane symmetric
solutions to equation \eqref{E:GEOWAVE} in two spatial dimensions,
to the wave equation \eqref{E:NONGEOMETRICWAVE} 
via the discussion in Appendix \ref{A:EXTENDINGTONONGEOMETRICWAVEEQUATIONS}, 
and to the equations described in Remark \ref{R:SLIGHTEXTENSIONNOBADTERMS}.
In a coordinate system of eikonal functions $u,v$,
all of those equations take the form 
\[ 
 \partial_u \partial_v \Psi 
 = 
\mathcal{N}(\Psi,\partial \Psi) \partial_u \Psi \partial_v \Psi 
\]
for some coefficient function $\mathcal{N}$. Hence,  
for simple waves (that is, waves with 
$ 
\partial_v \Psi
\equiv 0
$), 
the above analysis carries over without any changes. 

\subsection{Overview of the main steps in the proof}
\label{SS:OVERVIEW}
We now outline the main steps in the proof of Theorem~\ref{T:MAINTHEOREM}, 
which is our main result. 
Many of the geometric ideas and insights behind these steps are contained in \cite{dC2007}. Indeed,
the main theme of the present paper is that the framework of \cite{dC2007} can be extended to 
prove shock formation in solutions to quasilinear wave equations in a regime 
different than the one treated in \cite{dC2007}: the regime
of nearly simple outgoing plane symmetric waves. For a discussion of the main new ideas in the present paper, 
see Subsubsect.~\ref{SSS:DIFFERENCESFROMDISPERSIVEREGIME}.
\begin{enumerate}
	\item We formulate the shock formation problem so that
		the fundamental dynamic quantities to be solved for are 
		$\Psi$,
		$\upmu$, and the rectangular spatial components\footnote{Note that \eqref{E:FIRSTUPMU} and 
		\eqref{E:LUNITGRADIENTVECTORFIELDINTRO}
		imply that $\Lunit^0 \equiv 1$.} 
		$\Lunit^1$, $\Lunit^2$.
		We refer to the latter three quantities as ``eikonal function quantities'' since they depend on the first rectangular 
		derivatives of $u$.
		We then derive evolution equations for $\upmu$, $\Lunit^1$, and $\Lunit^2$
		along the integral curves of the vectorfield $\Lunit$. 
		These evolution equations are essentially equivalent to the eikonal equation \eqref{E:INTROEIKONAL}.
	\item We construct a good set of vectorfields 
		$\Fullset := \lbrace \Lunit, \Rad, \GeoAng \rbrace$
		that we use to commute the wave equation and 
		also the evolution equations for the eikonal function quantities.
		From the point of view of regularity considerations,
		it is important to appreciate that the rectangular components of 
		$Z \in \Fullset$ depend on the first rectangular derivatives of $u$.
		We will explain the importance of this fact in Subsubsect.~\ref{SSS:INTROENERGYESTIMATES}
		(see especially the discussion below equation \eqref{E:VERYSCHEMATICCOMMUTED}).
		Like $\CoordAng$, the vectorfield
		$\GeoAng$ (constructed in Subsect.~\ref{SS:COMMUTATIONVECTORFIELDS})
		is tangent to the $\ell_{t,u}$,
		but it has better regularity properties than $\CoordAng$.
		We use the full commutator set $\Fullset$
		when deriving $L^{\infty}$ estimates for the derivatives of the solution.
		When deriving energy estimates, \emph{we use only the $\mathcal{P}_u-$tangent subset} 
		$\Tanset := \lbrace \Lunit, \GeoAng \rbrace$.
	\item To derive estimates, we make bootstrap assumptions on an open-at-the-top bootstrap region
		$\mathcal{M}_{\Tboot,U_0} := \cup_{s \in [0,\Tboot)} \Sigma_s^{U_0}$,
		where $0 \leq \Tboot \leq 2 \TranminusdatasizeWithFactor^{-1}$ (see \eqref{E:INTROKEYBLOWUPTRANSVERALFACTOR})
		and $\mathcal{M}_{\Tboot,U_0}$ is a spacetime subset trapped in between left-most and right-most null hyperplanes
		and the flat bottom and top hypersurfaces $\Sigma_0$ and $\Sigma_{\Tboot}$; see Figure~\ref{F:SOLIDREGION}
		on pg.~\pageref{F:SOLIDREGION}.
		We assume that $\upmu > 0$ on $\mathcal{M}_{\Tboot,U_0}$, that is, that no shocks are present. 
		We then make ``fundamental'' bootstrap assumptions about the $L^{\infty}$ norms of various low-level derivatives of $\Psi$
		with respect to vectorfields in $\Tanset$. 
		These assumptions are non-degenerate in the sense 
		that they do not lead to infinite expressions even when $\upmu = 0$.
		Using them, we derive non-degenerate $L^{\infty}$ estimates for the low-level 
		$\Fullset$ derivatives of the eikonal function quantities and other
		low-level derivatives of $\Psi$.
		Moreover, in Sect.~\ref{S:SHARPESTIMATESFORUPMU},
		we derive related but much sharper estimates for $\upmu$ and 
		some of its low-level derivatives.
		In particular, using a posteriori estimates,
		we give a precise description showing that
		\emph{$\min_{\Sigma_t^u} \upmu$ vanishes linearly in $t$
		and moreover, we connect the vanishing rate to the
		initial data quantity $\TranminusdatasizeWithFactor$ 
		defined in \eqref{E:INTROKEYBLOWUPTRANSVERALFACTOR}.}\footnote{Specifically, we show that there exists a $(t,u)-$dependent constant $\LateTimeLUnitMu$ such that
		for $0 \leq s \leq t$, we have
		$
			\min_{\Sigma_s^u} \upmu 
			\approx
			1 - \LateTimeLUnitMu s
		$;
		see \eqref{E:MUSTARBOUNDS}.
		} 
		In addition, 
		we derive related sharp estimates for certain time-integrals 
		involving degenerate factors of 
		$1/\upmu$. 
		The time integrals 
		appear in the Gronwall estimates we use to derive
		a priori energy estimates, as we describe in Step (4).
		The estimates of Sect.~\ref{S:SHARPESTIMATESFORUPMU} 
		therefore play a critical role in closing our proof.
	\item We use the $L^{\infty}$ estimates to derive 
		up-to-top order $L^2$-type (energy) estimates for
		$\Psi$ and the eikonal function quantities on $\mathcal{M}_{\Tboot,U_0}$.
		This step is difficult, in part because we must overcome
		the potential loss of a derivative 
		tied to the dependence of our commutation vectorfields on the rectangular derivatives of $u$.
		To derive the $L^2$ estimates, we commute the evolution equations with 
		only the $\mathcal{P}_u-$tangent commutators $\Singletan \in \Tanset$.
		Because of the good null structure of the wave equation 
		highlighted in \eqref{E:INTROWAVEEQUATIONFRAMEDECOMPOSED}
		and the good properties of the vectorfields in $\Tanset$,
		we do not need to commute with the transversal derivative $\Rad$ when deriving the $L^2$ estimates.
		As we have mentioned, at the high derivative levels, the energies are allowed to
		blow up in a controlled fashion near the shock,
		while at the lower derivative levels, 
		the energies remain small all the way up to the shock.
		The degeneracy of the high-order estimates 
		is tied to our approach in avoiding the derivative loss:
		we work with modified quantities that have unexpectedly 
		good regularity properties but that introduce a difficult factor of
		$1/\upmu$ into the top-order energy identities.
		This $1/\upmu$ factor is the reason that we need the sharp
		time integral estimates described Step (3);
		these sharp estimates affect the blowup-rates of our top-order energy estimates,
		which are central to the entire proof.
		We remark that the degeneracy of our high-order energy estimates reflects the ``worst-case'' behavior of $\upmu$
		along $\Sigma_t$. That is, regions where $\upmu$ is small drive the degeneracy of our high-order energy estimates
		along all of $\Sigma_t$. An added layer of complexity
		is that near the time of first shock formation,
		$\upmu$ \emph{can be large at some points while being near $0$ at others
		and thus our energy estimates along $\Sigma_t$ have to simultaneously account for both of these extremes}.
		We also highlight again the following crucially important feature of our proof: 
		\emph{we must derive non-degenerate energy estimates at the low-derivative levels}.
		From such estimates, we can recover our fundamental $L^{\infty}$ 
		bootstrap assumptions via a simple geometric Sobolev embedding result
		(see Lemma~\ref{L:SOBOLEV}).
	\item The proof that $\upmu \to 0$ and causes blowup (i.e., that the shock forms)
		before the maximum allowed bootstrap time $2 \TranminusdatasizeWithFactor^{-1}$
		is easy given the non-degenerate low-level
		$L^{\infty}$ estimates; see Subsubsect.~\ref{SSS:INTROPOOFTHATMUVANISHES}
		for an outline of the proof.
\end{enumerate}

\begin{remark}[\textbf{Straightforward bootstrap structure}]
	\label{R:LINEARORDER}
	The bootstrap structure of our proof is very simple.
	Given the simple bootstrap assumptions from Step $(3)$,
	the logic of our proof is essentially linear: the proofs of our estimates depend 
	only on previously proved estimates.
	We recover the bootstrap assumptions near the end of the proof
	of the main theorem.
\end{remark}

Steps $(1)-(3)$ involve many geometric decompositions and computations but are relatively standard.
In the remainder of Sect.~\ref{S:INTRO}, we describe
Steps $(4)$ and $(5)$ in more detail, which have some important features that
are specific to the problem of shock formation. We start with the easy Step $(5)$.

\subsubsection{Outline of the proof that the shock happens}
\label{SSS:INTROPOOFTHATMUVANISHES}
The proofs that $\upmu$ goes to $0$ and that some first rectangular derivative of $\Psi$ blows up
are easy given the non-degenerate low-level estimates.
Both of these facts are based on the following evolution equation
(derived in Lemma~\ref{L:UPMUANDLUNITIFIRSTTRANSPORT} as a consequence of the eikonal equation):
\begin{align} \label{E:INTROMUEVOLUTION}
	\Lunit \upmu
	& = \frac{1}{2} G_{\Lunit \Lunit} \Rad \Psi
		+ \mathcal{O}(\upmu \Lunit \Psi).
\end{align}
In \eqref{E:INTROMUEVOLUTION}, 
$G_{\Lunit \Lunit} := \frac{d}{d \Psi}g_{\alpha \beta}(\Psi) \Lunit^{\alpha} \Lunit^{\beta}$
and the term $\mathcal{O}(\upmu \Lunit \Psi)$ is depicted schematically. Our assumptions on the nonlinearities
ensure that in the regime under study, we have $G_{\Lunit \Lunit} \approx 1$.
Using the $\mathring{\upepsilon}-\mathring{\updelta}$ hierarchy,
we have $\Lunit(G_{\Lunit \Lunit} \Rad \Psi) = \mathcal{O}(\mathring{\upepsilon})$.
Since $\Lunit = \frac{\partial}{\partial t}$ relative to the geometric coordinates, 
we can integrate this estimate to obtain
$[G_{\Lunit \Lunit} \Rad \Psi](t,u,\vartheta) 
= 
[G_{\Lunit \Lunit} \Rad \Psi](0,u,\vartheta)
+ \mathcal{O}(\mathring{\upepsilon})
$,
where the implicit constant in $\mathcal{O}$ is allowed to depend on
the expected shock time $\TranminusdatasizeWithFactor^{-1}$
(see \eqref{E:INTROKEYBLOWUPTRANSVERALFACTOR}).
Inserting into \eqref{E:INTROMUEVOLUTION}, we obtain
\begin{align} \label{E:INTROMUEVOLUTIONESTIAMTED}
	\Lunit \upmu(t,u,\vartheta)
	& = \frac{1}{2} [G_{\Lunit \Lunit} \Rad \Psi](0,u,\vartheta)
		+ \mathcal{O}(\mathring{\upepsilon}).
\end{align}
Integrating \eqref{E:INTROMUEVOLUTIONESTIAMTED} 
and using $\upmu(0,u,\vartheta) = 1 + \mathcal{O}(\mathring{\upepsilon})$, 
we find that
\begin{align} \label{E:UPMUINTEGRATEDINTRO}
	\upmu(t,u,\vartheta)
	& = 1
	+ \frac{1}{2} [G_{\Lunit \Lunit} \Rad \Psi](0,u,\vartheta) t
	+ \mathcal{O}(\mathring{\upepsilon}).
\end{align}
From \eqref{E:INTROKEYBLOWUPTRANSVERALFACTOR} and \eqref{E:UPMUINTEGRATEDINTRO},
we see that for $0 \leq t \leq 2 \TranminusdatasizeWithFactor^{-1}$, we have
\begin{align} \label{E:INTROMUDECAYESTIMATE}
	\min_{\Sigma_t^1} \upmu
	& = 1 - \TranminusdatasizeWithFactor t 
	+ \mathcal{O}(\mathring{\upepsilon}).
\end{align}
From \eqref{E:INTROMUDECAYESTIMATE}, we see that
$\upmu$ vanishes for the first time at
$T_{Lifespan} 
=
\left\lbrace 1 + \mathcal{O}(\mathring{\upepsilon}) \right\rbrace
\TranminusdatasizeWithFactor^{-1}
$.
Moreover, the above argument can easily be extended to show that
at the points $(T_{Lifespan},u,\vartheta)$ where $\upmu$ vanishes,
the quantity
$
|\Rad \Psi|(T_{Lifespan},u,\vartheta)
$
is uniformly bounded from below, strictly away from $0$;
see inequality \eqref{E:BLOWUPPOINTINFINITE} and its proof.
Since $\sqrt{g(\Rad,\Rad)} = \upmu$, we conclude that
the derivative of $\Psi$ with respect to the $g-$unit-length
vectorfield $\Radunit := \upmu^{-1} \Rad \Psi \sim - \partial_1 \Psi$
must blow up at the points $(T_{Lifespan},u,\vartheta)$
where $\upmu$ vanishes.

\subsubsection{Energy estimates at the highest order}
\label{SSS:INTROENERGYESTIMATES}
By far, the most difficult part of the analysis 
is obtaining the high-order $L^2$ estimates of Step $(4)$.
To derive them, we use the well-known multiplier method. Specifically, 
we derive energy identities by applying the divergence theorem to the vectorfield 
$J^{\alpha} := \enmomtensor_{\ \beta}^{\alpha} \Mult^{\beta}$ 
on the region $\mathcal{M}_{t,u}$, 
where 
$\enmomtensor_{\mu \nu}[\Psi]
:= \D_{\mu} \Psi \D_{\nu} \Psi
- \frac{1}{2} g_{\mu \nu} (g^{-1})^{\alpha \beta} \D_{\alpha} \Psi \D_{\beta} \Psi
$ 
is the energy-momentum tensorfield
(see \eqref{E:ENERGYMOMENTUMTENSOR})
and
$\Mult := (1 + 2 \upmu) \Lunit + 2 \Rad$
is a timelike vectorfield\footnote{In many other works, the
      symbol $\Mult$ denotes the future-directed unit normal to $\Sigma_t$.
      In contrast, in the present article, 
      the vectorfield $\Mult$
      is not the future-directed unit normal to $\Sigma_t$.
      \label{FN:MULTNOTATION}} 
verifying 
$
g(\Mult,\Mult) 
= 
- 4 \upmu (1 + \upmu)
< 0
$;
see Prop.~\ref{P:DIVTHMWITHCANCELLATIONS} for the precise statement and
Figure~\ref{F:SOLIDREGION} for a picture illustrating the region of integration.
As we have mentioned, 
we are able to close our energy estimates by commuting the wave equation with only 
$\mathcal{P}_u-$tangent commutators $\Singletan \in \Tanset$ 
(we commute with the $\mathcal{P}_u-$transversal vectorfield 
$\Rad$ only when deriving low-level $L^{\infty}$ estimates).
Moreover, 
we do not rely on the lowest level energy identity corresponding to the non-commuted equation.
That is, we derive energy estimates for 
$\Singletan \Psi$, $\Singletan \Singletan \Psi$, etc.\
Consequently, for our data, the energies are of small size $\mathring{\upepsilon}$ at time $0$.
At the first commuted level,
the energies $\mathbb{E}[\Singletan \Psi](t,u)$ and null fluxes $\mathbb{F}[\Singletan \Psi](t,u)$ 
have the following strength (note carefully which terms contain explicit $\upmu$ weights!):
\begin{subequations}
\begin{align} \label{E:INTROENERGY}
	\enzero[\Singletan \Psi](t,u)
		& \sim
		\int_{\Sigma_t^u} 
			\upmu (\Lunit \Singletan \Psi)^2
			+ (\Rad \Singletan \Psi)^2
			+ \upmu |\angdiff \Singletan \Psi|^2
		\, d \tvol,
		\\
	\flzero[\Singletan \Psi](t,u)
		& \sim
		\int_{\mathcal{P}_u^t} 
			(\Lunit \Singletan \Psi)^2
			+ \upmu |\angdiff \Singletan \Psi|^2
		\, d \conevol.
		\label{E:INTROFLUX}
\end{align}
\end{subequations}
In \eqref{E:INTROENERGY}-\eqref{E:INTROFLUX}, 
$\angdiff \Singletan \Psi$ denotes the $\ell_{t,u}-$gradient of $\Singletan \Psi$ 
(that is, the gradient of $\Singletan \Psi$ viewed as a function of the geometric torus coordinate $\vartheta$)
and the forms $d \tvol$ 
and $d \conevol$ are constructed\footnote{$d \tvol$ is a rescaled version of the canonical form induced by $g$ on $\Sigma_t$.} 
so that they remain non-degenerate all the way up to and including the shock.
We stress that the
\emph{terms with $\upmu$ weights in \eqref{E:INTROENERGY}-\eqref{E:INTROFLUX} 
become very weak near the shock,
and they are not useful for controlling error terms that lack $\upmu$ weights}.
Since both appearances of $|\angdiff \Singletan \Psi|^2$ in \eqref{E:INTROENERGY} involve $\upmu$ weights,
we must find a different way to
control error terms proportional to $|\angdiff \Singletan \Psi|^2$ \emph{that does not rely on}
$\enzero$ or $\flzero$. To this end, 
we exploit a subtle spacetime integral $\mathbb{K}(t,u)$
with special properties first identified by Christodoulou \cite{dC2007};
we explain this in Subsubsect.~\ref{SSS:SIGNEDSPACETIME} in more detail.

With $\mathbb{E}_M$ denoting the energy corresponding to commuting the wave equation
$M$ times with elements $\Singletan \in \Tanset$, 
$\mathring{\upepsilon}$ denoting the small size of the $L^2$ quantities at time $0$,
and $\upmu_{\star}(t,u) := \min \lbrace 1, \min_{\Sigma_t^u} \upmu \rbrace$, 
we derive the following energy estimate hierarchy (see Prop.~\ref{P:MAINAPRIORIENERGY}),
valid for classical solutions when $(t,u) \in [0,2 \TranminusdatasizeWithFactor^{-1}] \times [0,U_0]$:
\begin{subequations}
\begin{align}
	\mathbb{E}_{18}(t,u)
	& \leq C \mathring{\upepsilon}^2 \upmu_{\star}^{-11.8}(t,u),
		\label{E:INTROTOPENERGY} \\
	\mathbb{E}_{17}(t,u)
	& \leq C \mathring{\upepsilon}^2 \upmu_{\star}^{-9.8}(t,u),
		\label{E:INTROONEBELOWTOPENERGY} \\
	& \cdots 
		\notag \\
	\mathbb{E}_{13}(t,u)
	& \leq C \mathring{\upepsilon}^2 \upmu_{\star}^{-1.8}(t,u),
		\\
	\mathbb{E}_{12}(t,u)
	& \leq C \mathring{\upepsilon}^2,
		\label{E:INTROFIRSTNONDEGENERATEENERGY} \\
	& \cdots 
		\notag \\
	\mathbb{E}_1(t,u)
	& \leq C \mathring{\upepsilon}^2.
	\label{E:INTROLOWESTENERGY}
\end{align}
\end{subequations}
A similar hierarchy holds for the null fluxes $\mathbb{F}$ and the spacetime integrals $\mathbb{K}$.

We now explain how to derive the top-order energy estimate 
\eqref{E:INTROTOPENERGY} and the origin of its 
degeneracy with respect to $\upmu$.
The main difficulty that one confronts in deriving \eqref{E:INTROTOPENERGY}
is that naive estimates do not work at the top order
because they lead to the loss of a derivative.
The following mantra summarizes our approach to overcoming this difficulty.
\begin{changemargin}{.25in}{.25in}
	One can gain back the derivative, \emph{but only at the expense
	of incurring a factor of $\upmu^{-1}$ in the energy identities}.
\end{changemargin}
We now flesh out these issues.
The hardest step in deriving \eqref{E:INTROTOPENERGY}
is using the $L^{\infty}$ bootstrap assumptions and the $L^{\infty}$ estimates
to obtain the following top-order energy inequality:
\begin{align} \label{E:CARICATURETOPGRONWALLREADY}
	\mathbb{E}_{18}(t,u)
	& \leq C \mathring{\upepsilon}^2
		+ 
		4
		\int_{t'=0}^t
			\left\lbrace
				\sup_{\Sigma_{t'}^u}
				\left|
					\frac{\Lunit \upmu}{\upmu}
				\right|
			\right\rbrace
				\mathbb{E}_{18}(t',u)
		\, dt'
		+ \cdots.
\end{align}
The aforementioned factor of $\upmu^{-1}$
is the one indicated on RHS~\eqref{E:CARICATURETOPGRONWALLREADY}.
The second hardest step is estimating the singular ratio
$
\displaystyle
\sup_{\Sigma_{t'}^u} 
	\left|
		\frac{\Lunit \upmu}{\upmu} 
	\right|
$
in a way that allows us to derive a Gronwall estimate from \eqref{E:CARICATURETOPGRONWALLREADY}.
To estimate the ratio, we need sharp information describing how 
$\min_{\Sigma_{t'}^u} \upmu$ goes to $0$. 
This analysis is very technical and is based on a posteriori estimates 
involving possible late-time behaviors of $\upmu$; see Sect.~\ref{S:SHARPESTIMATESFORUPMU}.
A key ingredient is that by virtue of the wave equation \eqref{E:INTROWAVEEQUATIONFRAMEDECOMPOSED}
and equation \eqref{E:INTROMUEVOLUTION},
one can show that $\Lunit \Lunit \upmu = \mathcal{O}(\mathring{\upepsilon})$, 
which implies that $\Lunit \upmu$ is approximately constant along the integral curves of $\Lunit$
on the time scale of interest.
To explain the basic idea behind the Gronwall estimates, 
let us pretend that $\upmu$ is a function of $t$ alone,
that $\upmu$ is near $0$,
and that $\Lunit \upmu < 0$.
Then recalling that $\Lunit = \frac{\partial}{\partial t}$, 
we use Gronwall's inequality 
and \eqref{E:CARICATURETOPGRONWALLREADY} to derive
$\mathbb{E}_{18}(t,u) \leq C \mathring{\upepsilon}^2 \upmu^{-4}(t) \times \cdots$.
Note that the blowup-rate $\upmu^{-4}(t)$
is determined by the numerical constant $4$ on RHS~\eqref{E:CARICATURETOPGRONWALLREADY}.
In particular, \emph{it is important that the coefficient $4$ of the dangerous integral 
is a structural constant that does not depend on the number of times that the equations are differentiated.}
We remark that the blow-up exponent on RHS~\eqref{E:INTROTOPENERGY} is $11.8$ rather than $4$
because there are other difficult error integrals
on RHS~\eqref{E:CARICATURETOPGRONWALLREADY}
(which we ignore in this introduction) 
that contribute to the top-order degeneracy. 

We now sketch how we derive inequality \eqref{E:CARICATURETOPGRONWALLREADY}
and explain the appearance of the singular factor 
$
\displaystyle
\sup_{\Sigma_{t'}^u}
				\left|
					\frac{\Lunit \upmu}{\upmu}
				\right|
$. 
To illustrate the main ideas, we commute the wave equation one time with 
a $\mathcal{P}_u-$tangent commutation vectorfield $\Singletan$ constructed in Step $(2)$
and pretend that the wave equation in $\Singletan \Psi$ represents the top-order equation.
An important fact is that the rectangular components of the vectorfields $\Singletan \in \Tanset$ depend on
$\Psi$ and $\upmu \partial u$ (see \eqref{E:LUNITGRADIENTVECTORFIELDINTRO}).
Hence, upon commuting the wave equation with $\Singletan$, 
we obtain the following schematic wave equation:
\begin{align} \label{E:VERYSCHEMATICCOMMUTED}
	\upmu \square_{g(\Psi)} \Singletan \Psi 
	& = \upmu \partial^2 (\upmu \partial u) \cdot \partial \Psi 
		+ 
		\upmu \partial (\upmu \partial u) \cdot \partial^2 \Psi
		+ \cdots
\end{align}
In \eqref{E:VERYSCHEMATICCOMMUTED}, the schematic symbol $\cdot$ 
denotes tensorial contractions
\emph{that produce products with a special structure}.
Specifically, the $\Singletan$ are designed so that the worst
imaginable error terms are completely absent on RHS~\eqref{E:VERYSCHEMATICCOMMUTED},
which is possible only because we allow $\Singletan$ to depend on $\partial u$.
In particular, a careful decomposition of RHS~\eqref{E:VERYSCHEMATICCOMMUTED}
relative to the frame \eqref{E:INTROFRAME} reveals that
the factor $\Rad \Rad \Psi$ is absent. This is important because by
signature considerations, $\Rad \Rad \Psi$ would have come with the
singular factor $1/\upmu$,
which would prevent us from deriving non-degenerate estimates at the low orders. 
Because of this structure, all terms
$\upmu \partial (\upmu \partial u) \cdot \partial^2 \Psi$ are relatively easy to control all the way up to the shock.
The main difficulty is that the factor $\upmu \partial^2 (\upmu \partial u)$ 
on RHS~\eqref{E:VERYSCHEMATICCOMMUTED}
seems to have insufficient regularity to close the estimates:
commuting the eikonal equation \eqref{E:INTROEIKONAL}, one obtains
the evolution equation $\Lunit \partial^3 u \sim \partial^3 \Psi + \cdots$,
\emph{which is inconsistent with the available regularity} (two derivatives of $\Psi$)
for solutions to \eqref{E:VERYSCHEMATICCOMMUTED}.
Clearly this difficulty propagates upon further commuting the wave equation.
In the energy estimates, this difficulty leads to error integrals that are
hard to control near the shock. 
As we will explain, the most difficult (in the sense of degeneracy created by a factor of $1/\upmu$) 
error integral\footnote{More precisely, this error integral is difficult only
when the vectorfield $\Singletan$ in \eqref{E:MODELBADINTEGRAL} is equal to the 
$\ell_{t,u}-$tangent vectorfield $\GeoAng$. The case $\Singletan = \Lunit$ is much easier to treat
because in this case, one can show that the term $\partial^2 (\upmu \partial u)$ involves at least one $\Lunit$ 
differentiation. Consequently, we can use the Raychaudhuri equation described below to algebraically replace
$\partial^2 (\upmu \partial u)$ with terms involving $\leq 2$ derivatives of $\Psi$.
\label{FN:DIFFUCULTYONLYWHENPISLUNIT}}
has the following schematic form:
\begin{align} \label{E:MODELBADINTEGRAL}
	2
	\int_{t'=0}^t
		\int_{\Sigma_{t'}^u}
			\Rad \Psi \cdot \partial^2 (\upmu \partial u) \cdot \Rad \Singletan \Psi
		\, d \tvol
	\, dt',
\end{align}
where the factor $\partial^2 (\upmu \partial u)$ in \eqref{E:MODELBADINTEGRAL} 
has a special structure that we explain just below.
It remains for us to outline why \eqref{E:MODELBADINTEGRAL} 
can be expressed as the integral on RHS~\eqref{E:CARICATURETOPGRONWALLREADY}
plus other error integrals that are similar or easier to treat.
The key fact, explained in the next paragraph, is that
$\partial^2 (\upmu \partial u) 
= \upmu^{-1} \mathsf{Modified} 
+ \upmu^{-1} G_{\Lunit \Lunit} \Rad \Singletan \Psi
+ \cdots
$,
where $G_{\Lunit \Lunit}$ is as in \eqref{E:INTROMUEVOLUTIONESTIAMTED},
$\mathsf{Modified}$ solves a good evolution equation with source terms that have an
allowable level of regularity,
and $\cdots$ denotes terms that are easy to treat.
Then observing that $\Rad \Psi \cdot \partial^2 (\upmu \partial u)$ contains the special product 
$G_{\Lunit \Lunit} \Rad \Psi$, we may  
use \eqref{E:INTROMUEVOLUTION} to substitute,
which allows us to rewrite \eqref{E:MODELBADINTEGRAL} in the form
\begin{align} \label{E:INTROREPRESENTATIVEERRORINTEGRAL}
	4
	\int_{t'=0}^t
		\int_{\Sigma_{t'}^u}
			\frac{\Lunit \upmu}{\upmu} (\Rad \Singletan \Psi)^2
		\, d \tvol
	\, dt'
	+
	2
	\int_{t'=0}^t
		\int_{\Sigma_{t'}^u}
			(\Rad \Psi)
			\frac{\mathsf{Modified}}{\upmu} 
			(\Rad \Singletan \Psi)
		\, d \tvol
	\, dt'
	+ \cdots.
\end{align}
From \eqref{E:INTROENERGY} and the first integral in \eqref{E:INTROREPRESENTATIVEERRORINTEGRAL},
we obtain the difficult integral on RHS~\eqref{E:CARICATURETOPGRONWALLREADY}.
The integral involving $\mathsf{Modified}$ in \eqref{E:INTROREPRESENTATIVEERRORINTEGRAL}
is difficult to treat,\footnote{We ignore it here; see the proofs of 
Prop.~\ref{P:KEYPOINTWISEESTIMATE}
and \ref{P:TANGENTIALENERGYINTEGRALINEQUALITIES}
and
Lemma~\ref{L:DIFFICULTTERML2BOUND} for the details.}
but the resulting estimates are similar to the ones that we have
sketched for the first integral.

We now elaborate on the special structure of the factor
$\partial^2 (\upmu \partial u)$ appearing in \eqref{E:MODELBADINTEGRAL}.
Some rather involved computations 
(see Lemmas~\ref{L:DEFORMATIONTENSORFRAMECOMPONENTS} and \ref{L:BOXZCOM} and Prop.~\ref{P:COMMUTATIONCURRENTDIVERGENCEFRAMEDECOMP})
yield that the factor $\partial^2 (\upmu \partial u)$ appearing in \eqref{E:MODELBADINTEGRAL}
is equal to the geometric quantity $\upmu \Singletan \mytr \upchi$,
where $\upchi$ is the symmetric type $\binom{0}{2}$ 
$\ell_{t,u}-$tangent\footnote{Note that the $\ell_{t,u}$ are one-dimensional curves and hence
for any $m$ and $n$, the space of all type $\binom{m}{n}$ $\ell_{t,u}-$tangent tensors is one-dimensional. 
Hence, the study of $\ell_{t,u}-$tangent tensorfields could be completely reduced to the study of scalar functions.
However, we do not carry out such a reduction in this article; we prefer to retain the tensorial
character of $\ell_{t,u}-$tangent tensorfields 
because that structure allows us to directly apply 
standard formulas and techniques from differential geometry.
} 
tensorfield defined by 
$\upchi_{\CoordAng \CoordAng} := g(\D_{\CoordAng} \Lunit, \CoordAng)$,
and $\mytr$ denotes the trace with respect to the Riemannian metric
$\gsphere$ induced on the $\ell_{t,u}$ by $g$.
To estimate $\mytr \upchi$, 
we rely on the well-known \emph{Raychaudhuri equation} from geometry, 
which yields the evolution equation
$\Lunit \mytr \upchi = - \Ric_{\Lunit \Lunit} + \cdots$, 
where $\Ric_{\Lunit \Lunit} := \Ric_{\alpha \beta} \Lunit^{\alpha} \Lunit^{\beta}$ 
is a component of the Ricci curvature tensor of $g(\Psi)$
and the terms $\cdots$ involve fewer derivatives.
The key point is that a careful decomposition 
(see Lemma~\ref{L:ALPHARENORMALIZED})
shows that for solutions to \eqref{E:GEOWAVE},
all top-order terms contain a perfect $\Lunit$ derivative:
$\upmu \Ric_{\Lunit \Lunit} =  \Lunit(-G_{\Lunit \Lunit} \Rad \Psi + \upmu \Singletan \Psi) + \cdots$,
where the factor $-G_{\Lunit \Lunit} \Rad \Psi$ is precisely depicted.
This remarkable structure was first\footnote{A related but simpler observation was made in \cite{dCsK1993}.} 
observed\footnote{Although the authors needed to exploit this structure to avoid losing
a derivative in their work \cite{sKiR2003},
they did not need to address the difficulty of obtaining estimates in regions
where $\upmu$ is near $0$.
}  
by Klainerman and Rodnianski 
in their proof of low regularity well-posedness
for quasilinear wave equations \cite{sKiR2003}
and was also used in \cites{dC2007,jS2014b,sMpY2014}. Combining, we find that
$
\Lunit 
\left\lbrace
	\upmu \mytr \upchi
	- 
	G_{\Lunit \Lunit} \Rad \Psi
	+
	\upmu \Singletan \Psi
\right\rbrace
= \cdots.
$
Taking one $\Singletan$ derivative
and setting 
$\mathsf{Modified} :=
	\upmu \Singletan \mytr \upchi
	- 
	G_{\Lunit \Lunit} \Rad \Singletan \Psi
	+
	\upmu \Singletan \Singletan \Psi$,
we find that
$\Lunit \mathsf{Modified} = l.o.t.$ as desired,
where $l.o.t.$ denotes terms with an allowable degree of
differentiability.

\begin{remark}[\textbf{The need for elliptic estimates in three or more spatial dimensions}]
	\label{R:ELLIPTIC}
	In $n$ spatial dimensions with $n \geq 3$,
	it is no longer possible to obtain an equation of the form
	$\Lunit \mathsf{Modified} = l.o.t$. The difficulty is that
	some third derivatives of $u$ still remain on the RHS:
	$\Lunit \mathsf{Modified} = \partial^2 (\upmu \partial u) + l.o.t$.
	However a careful decomposition of the remaining term 
	$\partial^2 (\upmu \partial u)$ on the RHS shows that
 $\partial^2 (\upmu \partial u) \sim \hat{\upchi} \cdot \Lie_{\Singletan} \hat{\upchi}$, 
	where $\Lie$ denotes Lie differentiation 
	and $\hat{\upchi}$ is the trace-free part of $\upchi$, which vanishes when $n=2$.
	To bound the top-order factor $\Lie_{\Singletan} \hat{\upchi}$ in $L^2$, one can derive elliptic estimates
	on the $n-1$ dimensional surfaces analogous to the $\ell_{t,u}$ in the 
	present article; see, for example, \cites{sKiR2003,dCsK1993,dC2007,jS2014b,sMpY2014}
	for more details.
\end{remark}

\subsubsection{Less degenerate energy estimates at the lower orders}
\label{SSS:LOWERORDERENERGYESTIMATES}
We now explain why the energy estimates 
\eqref{E:INTROONEBELOWTOPENERGY}-\eqref{E:INTROFIRSTNONDEGENERATEENERGY}
become two powers less degenerate relative to $\upmu_{\star}^{-1}$ 
at each level in the descent, which eventually brings us to the non-degenerate levels 
\eqref{E:INTROFIRSTNONDEGENERATEENERGY}-\eqref{E:INTROLOWESTENERGY}.
To illustrate the method, we now pretend that 
equation \eqref{E:VERYSCHEMATICCOMMUTED} represents
one level below top order (equivalently, that three derivatives of $\Psi$ in the norm
$\| \cdot \|_{L^2(\Sigma_t^u)}$ represents top order).
The main idea is to allow the loss of
one derivative in the factor $\partial^2 (\upmu \partial u) \sim \Singletan \mytr \upchi$
in \eqref{E:MODELBADINTEGRAL}; 
a loss of one derivative is permissible below top order. 

We refer to the just-below-top-order energy
that we are trying to estimate 
by $\mathbb{E}_{One-Below-Top}$.
In this case, we can use the non-degenerate low-level estimate
$\|\Rad \Psi\|_{L^{\infty}(\Sigma_t^u)} \lesssim 1$
and Cauchy-Schwarz to bound the error integral in \eqref{E:MODELBADINTEGRAL} by
\begin{align} \label{E:ONEBELOWTOP}
	\lesssim
	\int_{t'=0}^t
		\left\| 
			\Singletan \mytr \upchi
		\right\|_{L^2(\Sigma_{t'}^u)}
		\mathbb{E}_{One-Below-Top}^{1/2}(t',u)
	\, dt'.
\end{align}
The expression \eqref{E:ONEBELOWTOP} leads to a gain in powers of $\upmu_{\star}$
because of the following
critically important estimate (see \eqref{E:LOSSKEYMUINTEGRALBOUND}), 
which shows that integrating
in time produces the gain:
for constants $\Contwo > 1$, we have
\begin{align} \label{E:INTROLOSSKEYMUINTEGRALBOUND}
	\int_{t'=0}^t
		\frac{1}{\upmu_{\star}^{\Contwo}(t',u)}
	\, dt'
	\lesssim
	\upmu_{\star}^{1-\Contwo}(t,u).
\end{align}

The point is that there are two time integrations in \eqref{E:ONEBELOWTOP},
the obvious one, and the one that comes from the schematic relation
$
\left\| 
	\Lunit \Singletan \mytr \upchi
\right\|_{L^2(\Sigma_{t'}^u)}
\sim
\left\| 
	\Singletan \Ric_{\Lunit \Lunit}
\right\|_{L^2(\Sigma_{t'}^u)}
+ \cdots
\sim
\left\| 
	\Singletan \Singletan \Singletan \Psi
\right\|_{L^2(\Sigma_{t'}^u)}
+ \cdots
\sim \upmu_{\star}^{-1/2}(t',u) \mathbb{E}_{Top}^{1/2}(t',u)
+ \cdots
$,
where we have incurred the factor $\upmu_{\star}^{-1/2}$ in the last step
due to the fact that the energies
control ``geometric torus derivatives'' $\Singletan = \angdiff$ with a $\upmu^{1/2}$ weight
(see \eqref{E:INTROENERGY}).
By the already 
proven\footnote{In practice, we have to derive a Gronwall 
estimate for the top-order and just-below-top-order energies as a system, 
rather than treating the top energy completely separately.} 
bound $\mathbb{E}_{Top}^{1/2}(t',u) \lesssim \mathring{\upepsilon} \upmu_{\star}^{-5.9}$
(see \eqref{E:INTROTOPENERGY})
we can integrate the previous estimate in time
(see Lemma~\ref{L:L2NORMSOFTIMEINTEGRATEDFUNCTIONS})
to yield, via \eqref{E:INTROLOSSKEYMUINTEGRALBOUND},
the estimate
$
\left\| 
	\Singletan \mytr \upchi
\right\|_{L^2(\Sigma_{t'}^u)}
\sim \int_{s=0}^t 
				\left\| 
					\Lunit \Singletan \mytr \upchi
				\right\|_{L^2(\Sigma_s^u)}
			\, ds
		+ \cdots
	\sim 
			\mathring{\upepsilon}
			\int_{s=0}^{t'} 
				\upmu_{\star}^{-6.4}(s,u)
			\, ds
		+ \cdots
\lesssim 
\mathring{\upepsilon}
\upmu_{\star}^{-5.4}(t,u)
+ \cdots
$.
The outer time integration in \eqref{E:ONEBELOWTOP}
leads to the gain of another power of $\upmu_{\star}$,
which in total yields the a priori estimate\footnote{We have ignored some other error integrals which are slightly more
degenerate and only allow us to prove the slightly weaker estimate
$\mathbb{E}_{One-Below-Top}^{1/2}(t,u)
\lesssim \mathring{\upepsilon}
\upmu_{\star}^{-4.9}(t,u)$.}
$\mathbb{E}_{One-Below-Top}^{1/2}(t,u)
\lesssim \mathring{\upepsilon}
\upmu_{\star}^{-4.4}(t,u)
+ \cdots,
$
an improvement over the top-order degeneracy.
We can continue the descent in this fashion,
and when we reach the level \eqref{E:INTROLOWESTENERGY},
the following analog of \eqref{E:INTROLOSSKEYMUINTEGRALBOUND}
(proved below as \eqref{E:LESSSINGULARTERMSMPOINTNINEINTEGRALBOUND})
\emph{allows us to completely break the $\upmu_{\star}^{-1}$ degeneracy}:
\begin{align} \label{E:INTRONODEGENERACYKEYMUINTEGRALBOUND}
	\int_{t'=0}^t
		\frac{1}{\upmu_{\star}^{9/10}(t',u)}
	\, dt'
	\lesssim 1.
\end{align}
We conclude by remarking that the
proofs of \eqref{E:INTROLOSSKEYMUINTEGRALBOUND} and \eqref{E:INTRONODEGENERACYKEYMUINTEGRALBOUND} 
are based on knowing exactly how $\upmu_{\star}$ goes to $0$, that is, based on a
sharp version of the caricature estimate 
$\upmu_{\star}(t,u) \sim 1 - t \TranminusdatasizeWithFactor$;
see \eqref{E:MUSTARBOUNDS}.
In particular, it is very important that
$\upmu_{\star}$ goes to $0$ \emph{linearly} in time.

\subsubsection{The coercive spacetime integral}
\label{SSS:SIGNEDSPACETIME}
As we highlighted in Subsubsect.~\ref{SSS:INTROENERGYESTIMATES},
the energies \eqref{E:INTROENERGY} and null fluxes \eqref{E:INTROFLUX}
control geometric torus derivatives with $\upmu$ weights, which makes them
too weak to control certain error integrals involving torus derivatives
that lack $\upmu$ weights, at least in regions where $\upmu$ is small.
The saving grace is that as in \cites{dC2007,jS2014b,sMpY2014},
our energy estimates generate a spacetime integral
with a good sign. Under appropriate bootstrap assumptions,
the integral is strong in regions where $\upmu$ is small
and controls geometric torus derivatives without $\upmu$ weights.
For the $\Singletan$-commuted wave equation, this
integral takes the form
\begin{align} \label{E:INTROBONUSSPACETIME}
\coercivespacetime[\Singletan \Psi](t,u)
		& :=
	 	\frac{1}{2}
	 	\int_{\mathcal{M}_{t,u}}
			[\Lunit \upmu]_{-}
			|\angdiff \Singletan \Psi|^2
		\, d \vol,
\end{align}
where $[\Lunit \upmu]_{-} = \left|\min \lbrace \Lunit \upmu, 0 \rbrace \right|$.
The key estimate that makes \eqref{E:INTROBONUSSPACETIME} useful in
regions of small $\upmu$ is:
$\upmu(t,u,\vartheta) \leq \frac{1}{4}
\implies
\Lunit \upmu(t,u,\vartheta) \leq - \frac{1}{4} \TranminusdatasizeWithFactor$
(see \eqref{E:SMALLMUIMPLIESLMUISNEGATIVE}).
Here $\TranminusdatasizeWithFactor > 0$ is the data-dependent
parameter \eqref{E:INTROKEYBLOWUPTRANSVERALFACTOR}
that controls the blowup-time.
$\TranminusdatasizeWithFactor$ is large enough to be useful because
of our assumption that $\mathring{\upepsilon}$
is sufficiently small.
Note that the key estimate has a ``point of no return character'' 
in that once $\upmu$ becomes sufficiently small, it must continue to shrink along the integral
curves of $\Lunit$ to form a shock. The proof of the key estimate is non-trivial
and is part of the detailed analysis of $\upmu$ located in Sect.~\ref{S:SHARPESTIMATESFORUPMU}.

\subsection{Comparison with previous work}
\label{SS:PREVIOUSWORK}
\subsubsection{Blowup-results in one spatial dimension}
Under the assumption of plane symmetry,
the finite-time breakdown of solutions
to \eqref{E:GEOWAVE} or \eqref{E:NONGEOMETRICWAVE}
(for nonlinearities verifying the conditions described in Subsect.~\ref{SS:EQUATIONINRECTANGULARCOMPONENTS})
is well known and can be proved through 
the method of characteristics;
our analysis in Subsect.~\ref{SS:SPSW} was essentially a simple version of this method.
Readers may consult \cites{gHsKjSwW2016,jS2014b} 
for detailed examples derived with the help of sharp techniques 
paralleling the ones employed in the present article. 
There is a vast literature on the use of the method of 
characteristics to prove blowup for various nonlinear hyperbolic systems.
A far-from-exhaustive list of examples is:
the groundbreaking work of Riemann \cite{bR1860} mentioned at the beginning,
Lax's seminal finite-time breakdown results \cite{pL1973}
for scalar conservation laws
and his aforementioned application of the method of Riemann invariants 
to $2 \times 2$ genuinely nonlinear strictly hyperbolic systems \cite{pL1964},
Jeffrey's work \cite{aJ1965} on magnetoacoustics,
Jeffrey-Korobeinikov's work \cite{aJvK1969} on nonlinear electromagnetism,
Jeffrey-Teymur's work \cite{aJmT1974} on hyperelastic solids,
John's extension \cite{fJ1974} of Lax's work to systems in one spatial dimension
with more than two unknowns
(which required the development of new methods, in particularly identifying the important role played by \emph{simple waves}, since the method of Riemann invariants is
no longer applicable),
Liu's further refinement \cite{tL1979} of John's work,
John's work \cite{fJ1984} on spherically symmetric solutions to the equations of elasticity,
Klainerman-Majda's work \cite{sKaM1980} 
on nonlinear vibrating string equations,
Bloom's work \cite{fB1993} on nonlinear electrodynamics,
and Cheng-Young-Zhang's work \cite{gCrYqZ2013} on magnetohydrodynamics and related systems.
Roughly, the blowup in all of these works is proved
by finding a quantity $y(t)$ that verifies a 
Riccati-type equation $\dot{y}(t) = a(t) y^2(t) + \mbox{\upshape Error}$,
where $a(t)$ is non-integrable in time near $\infty$ and $\mbox{\upshape Error}$ is a small
error term that does not interfere with the blowup.
Recently, Christodoulou and Raoul Perez gave a new sharp proof \cite{dCdRP2016} 
of John's blowup-results \cite{fJ1974}
for genuinely nonlinear strictly hyperbolic quasilinear first-order systems
in one spatial dimension.
They showed that these systems can be treated with 
extensions of Christodoulou's framework \cite{dC2007},
which yields a sharp description of the blowup with upper and lower bounds
on the lifespan. Moreover, they applied their results
to prove shock formation in electromagnetic plane waves in a crystal. 

\subsubsection{Proofs of breakdown by a contradiction argument in more than one spatial dimension}
For nonlinear hyperbolic equations in more than one spatial dimension, 
many blowup-results have been proved by a contradiction argument
that bypasses the need to obtain a detailed description of the singularity.
For example, John gave a non-constructive proof \cite{fJ1981}  
showing that many wave equations in three spatial dimensions
with quadratic nonlinearities
exhibit finite-time blowup
for a large set\footnote{For some nonlinearities,
John's proof yields blowup for \emph{all} non-trivial, smooth, compactly supported data.} 
of smooth data. He did not need to impose any size restriction on the data for his proof to work,
but his proof did not provide any information about the blowup-time.
As a second example, we mention Sideris' well-known proof \cite{tS1985} 
of blowup for the compressible Euler equations in three spatial dimensions under a convexity
assumption on the equation of state and under signed integral conditions on the data. 
His proof was based on virial identity arguments that
yielded a manifestly non-negative weighted space-integrated quantity 
with a sufficiently negative time derivative,
which eventually leads to a contradiction
even if one assumes that the solution is otherwise smooth.
In particular, his proof gave an explicit upper bound on the solution's lifespan.
There are many similar results available which prove blowup
for various evolution equations via a virial identity argument. 
We do not aim to survey the extensive literature here,
but we do highlight the following examples:
semilinear Schr\"{o}dinger equations \cite{rG1977},
the relativistic Vlasov-Poisson equation \cite{rGjS1985},
and various semilinear wave and heat equations \cite{pKwS2007}.
We note that for semilinear Schr\"{o}dinger, wave, and related equations,
the state of the art knowledge 
of the blowup has advanced far beyond proof of blowup by contradiction;
see \cites{pR2013,yMfMpRjS2014} for surveys.

Though appealing in its shortness,
a serious limitation of the virial identity approach is that it 
relies specific algebraic structures of the equations that
are unstable under perturbations of the equations.
Another limitation is that it provides 
a lifespan upper bound that can be inaccurate; 
without additional information, one must concede that
the solution could in principle blow up much sooner by a different mechanism.
In contrast, our proof has many robust elements (see, however, Remark~\ref{R:SLIGHTEXTENSIONNOBADTERMS}),
and our work yields a sharp description
of the solution's lifespan and identifies the quantities that
blow up as well as the ones that remain regular.

\subsubsection{Detailed blowup-results in more than one spatial dimension}
Alinhac was the first \cites{sA1999a,sA2001a} 
to give a sharp description of singularity formation
in solutions to quasilinear wave equations
in more than one spatial dimension without symmetry assumptions.
He addressed a compactly supported small-data regime in which dispersive effects
are eventually overcome by sufficiently strong quadratic nonlinearities.
For convenience, even though these kinds of solutions eventually blow up,
we say that they belong to the ``small-data dispersive regime.''
Alinhac's results have been generalized to various equations by several authors;
see, for example, \cites{bDiWhY2013a,bDiWhY2015b,bDiWhY2015}.
In the case of three spatial dimensions (more precisely, the data are given on $\mathbb{R}^3$), 
Alinhac proved that whenever the nonlinearities in equation \eqref{E:NONGEOMETRICWAVE} fail
to satisfy Klainerman's null condition \cite{sk1984},
there exists a set of data of small size $\mathring{\upepsilon}$ 
(in a Sobolev norm) such that the solution
decays for a long time at the linear rate $t^{-1}$ before 
finally blowing up at the ``almost global existence'' time $\sim \mbox{exp}(c \mathring{\upepsilon}^{-1})$. 
More precisely, the singularity-forming quantities\footnote{In Alinhac's 
equations of type \eqref{E:NONGEOMETRICWAVE}, the second rectangular derivatives of
the solution blow up. In our work on equations of type \eqref{E:GEOWAVE}, the first rectangular derivatives blow up.} behave like 
$
\displaystyle
\frac{\mathring{\upepsilon}}{(1 + t)\left[1 + \mathcal{O}(\mathring{\upepsilon}) \ln(1 + t)\right]}
$,
where the $\mathcal{O}(\mathring{\upepsilon})$ term in the denominator
depends on the nonlinearities as well as the profile of the data
and the blowup (for some $t > 0$) occurs in regions where $\mathcal{O}(\mathring{\upepsilon}) < 0$.
%
Alinhac's data were posed in an annular region 
of $\mathbb{R}^3$, and he assumed that they verified
a non-degeneracy condition. 
His results showed that the almost global existence
lifespan lower bounds, 
obtained by John and Klainerman
\cites{sK1983, fJsK1984, sK1985} 
with the help of dispersive estimates that delay\footnote{By ``delay,'' we mean relative
to the case of one spatial dimension, where the lack of dispersion leads to blowup at time
$\mathcal{O}(\mathring{\upepsilon}^{-1})$.} 
the singularity formation,
are in fact saturated.
Moreover, his results confirmed John's conjecture \cite{fJ1987}
regarding the asymptotically correct description of 
the blowup-time in the limit $\mathring{\upepsilon} \downarrow 0$
for data verifying the non-degeneracy condition.

Christodoulou's remarkable work \cite{dC2007} yielded a sharp improvement (described below) of Alinhac's results
for a similar class of small compactly supported data given on $\mathbb{R}^3$, 
and he did not make any non-degeneracy assumption.
His main results applied to irrotational regions of solutions to the special relativistic Euler equations
in the small-data dispersive regime.
In such regions, the fluid equations reduce\footnote{Up to simple renormalizations outlined in Appendix~\ref{A:EXTENSIONTOEULER}.} 
to a special case of the wave equation \eqref{E:NONGEOMETRICWAVE}
in which additional structure is present. The non-relativistic Euler equations
were treated through the same approach in \cite{dCsM2014} and feature the same additional structure,
including that the irrotational fluid equations derive from a Lagrangian
(and thus can be written in Euler-Lagrange form)
and that solutions possess several conserved quantities associated to
various symmetries of the Lagrangian. These assumptions were used in the proofs, in particular
in exhibiting the good null structure\footnote{In particular, 
in Christodoulou's version of equation 
\eqref{E:PSIGEOMETRIC}, the RHS completely vanishes.
The vanishing occurs because he studies equations 
of the form \eqref{E:NONGEOMETRICWAVE} that derive from a Lagrangian.
} 
enjoyed by the equations.
The equations also had some additional structure due to the assumption that they model a physical fluid.
In addition to assuming that the data are of a small size $\mathring{\upepsilon}$ in a high Sobolev norm,
Christodoulou also made further assumptions on the data to ensure that a shock forms.
His sufficient conditions were phrased in terms of certain integrals of the data:
shocks form in the solution whenever the data integrals have the appropriate sign 
(determined by the nonlinearities) and are not too small in magnitude relative to 
$\mathring{\upepsilon}$.

Christodoulou's results were extended \cite{jS2014b}
to a larger class of equations and data by Speck
(see also the survey article \cite{gHsKjSwW2016}, 
joint with Holzegel, Klainerman, and Wong).
In particular, for data given on $\mathbb{R}^3$,
he proved a sharp small-data shock formation result
for equations \eqref{E:GEOWAVE} or \eqref{E:NONGEOMETRICWAVE}
whenever the null condition fails. That is, he showed that Christodoulou's sharp shock formation results
are not tied to the specific structure of the fluid equations
and that the additional structure present in those equations is not needed to close the proof.
Speck also showed that given any sufficiently regular non-trivial compactly supported initial data,
if they are rescaled by a small positive factor, then the solution forms a shock in finite time.
That is, \emph{all sufficiently regular data profiles lead to shock formation if they are suitably rescaled}.

Alinhac's and Christodoulou's approaches to proving shock formation share many common features.
For example, the main idea of Alinhac's proof was to resolve the singularity by
constructing an eikonal function $u$, as in Subsect.~\ref{SS:PROOFOVERVIEW}. 
Moreover, near the singularity, he changed variables to a new ``geometric'' coordinate system in which $u$ is one of the new coordinates. 
Relative to the geometric coordinates, he proved that the solution to \eqref{E:NONGEOMETRICWAVE}
remains regular all the way up to the point where the characteristics first intersect 
but that the change of variables map between the rectangular and geometric coordinates breaks down there.
Changing variables back to rectangular coordinates, he showed that the degeneracy implies that
$|\partial^2 \Phi|$ blows up in finite time precisely at the point where the characteristics intersect.
Alinhac also had to overcome the potential loss of derivatives that we described in Subsubsect.~\ref{SSS:INTROENERGYESTIMATES}
with the help of ``modified'' quantities. However, the methods he used did not immediately 
eliminate all of the derivative loss and thus differed in a fundamental way from Christodoulou's approach. Specifically,
to close his energy estimates, Alinhac employed a Nash-Moser iteration scheme.
His scheme featured a free boundary
due to the fact that the blowup-time for each iterate can be slightly different.
Although Alinhac gave a sharp description of the asymptotic behavior of the solution near the singularity,
his proof was not able to reveal information beyond the first blowup-point.
Moreover, in order for his proof to close, 
\emph{the constant-time hypersurface of first blowup was allowed to contain only one blowup-point}.
These fundamental technical limitations were tied to the presence of the free boundary in his Nash-Moser iteration scheme
and they are the reason that he had to make the non-degeneracy assumption on the data; 
see \cite{jS2014b} for additional discussion regarding his approach.

We now describe the most important difference between the approaches of Alinhac and Christodoulou.
The main advantage afforded by Christodoulou's framework,
as shown in \cites{dC2007,dCsM2014,jS2014b}, 
is that in the small-data dispersive regime,
there is a sharp criterion for blowup.
Specifically, the solution blows up at a given point $\iff$ $\upmu$ vanishes there.
In particular, in the small-data dispersive regime,
\emph{shocks are the only kinds of singularities that can form}.
Since the behavior of $\upmu$ is local in time and space, 
the vanishing of $\upmu$ at one point
does not preclude one from continuing the solution to a neighborhood of other nearby points where 
$\upmu > 0$. Moreover, $\lbrace \upmu = 0 \rbrace$ precisely
characterizes the singular portion of the boundary of the maximal development of the data,
that is, the portion of the boundary on which the solution blows up.
Thus, Christodoulou's framework is able to reveal detailed information about the structure of the maximal development
of the data, the shape of the various components of its boundary, and the
behavior of the solution along it. The same information can be extracted for
the solutions that we study here; see Remark~\ref{R:MAXIMALDEVELOPMENT}. 
The sharp description is an essential ingredient in setting up 
the problem of extending the Euler solution weakly beyond the first singularity. 
We note that an essential component of solving this problem is obtaining
information about the shock hypersurface across which 
discontinuities occur. The problem was recently solved in spherical symmetry \cite{dCaL2016},
while the non-symmetric problem remains open and is expected to be of immense difficulty.




\subsubsection{Differences between the proof of shock formation in the small-data dispersive regime and in the nearly plane symmetric regime}
\label{SSS:DIFFERENCESFROMDISPERSIVEREGIME}
As we mentioned near the beginning of Sect.~\ref{S:INTRO}, 
the most important new feature of the analysis in the 
nearly plane symmetric regime
is that we rely on a different mechanism to control the nonlinear error terms.
More precisely, since solutions do not decay in the nearly plane symmetric regime,
our approach is based on the propagation of the $\mathring{\upepsilon}-\mathring{\updelta}$
hierarchy described in Subsubsect.~\ref{SSS:INTROPOOFTHATMUVANISHES} (in other words, proving that our solution remain close to a simple outgoing wave),
rather than the smallness and dispersive decay estimates\footnote{We recall that in both regimes, 
the solution remains regular at the low derivative levels with respect to the geometric coordinates and
the blowup occurs in the partial derivatives of the solution with respect to the rectangular coordinates.}
used in the small-data dispersive regime \cites{dC2007,dCsM2014,jS2014b}.
We remark that there is a technical simplification in the nearly plane symmetric regime
that allows for a shorter proof compared to the small-data dispersive regime:
our propagation of the $\mathring{\upepsilon}-\mathring{\updelta}$ hierarchy 
does not involve weights in $t$ or the Euclidean radial coordinate $r$.



To propagate the $\mathring{\upepsilon}-\mathring{\updelta}$ hierarchy,
we must make some observations about various product/null structures in the equations that are
not needed for treating the small-data dispersive regime. 
Such structures are relevant both for obtaining suitable energy estimates up to top order
and for deriving non-degenerate $L^{\infty}$ estimates at the lower 
derivative levels. We now give one example of such a structure:
\begin{changemargin}{.25in}{.25in}
Repeatedly commuting the wave equation
$\upmu \square_{g(\Psi)} \Psi = 0$ 
up to top order
with $\mathcal{P}_u-$tangent vectorfields 
$\Singletan \in \Tanset = \lbrace \Lunit, \GeoAng \rbrace$
produces commutator error term products that are
quadratic and higher order in the derivatives of $\Psi$, $\upmu$, and $\Lunit^i$
\emph{with each product involving no more than one $\Rad$ derivative}.
\end{changemargin}
The above structure is 
a consequence of the schematic 
structures 
$[\Singletan_1,\Singletan_2] \sim \Singletan_3$
and
$[\Rad,\Singletan_1] \sim \Singletan_2$,
where $\Singletan_1$, $\Singletan_2$, and $\Singletan_3$ are arbitrary $\mathcal{P}_u-$tangent vectorfields.
These schematic commutator relations are easy to see relative to the geometric coordinates $(t,u,\vartheta)$.
To further explain these issues, we first note that 
$
\Rad
= \frac{\partial}{\partial u}
- 
\XiCoordComp \frac{\partial}{\partial \vartheta}
$,
where
$
\CoordAng
= \frac{\partial}{\partial \vartheta}
$
and $\XiCoordComp$ is a scalar function (see \eqref{E:RADSPLITINTOPARTTILAUANDXI}).
From these expressions, it easily follows that
for $Z_1, Z_2 \in \lbrace \Rad, P_1, P_2 \rbrace$,
the commutator $[Z_1,Z_2]$ belongs to
$\mbox{\upshape span} \lbrace \frac{\partial}{\partial t}, \frac{\partial}{\partial \vartheta} \rbrace$
and is therefore $\mathcal{P}_u^t-$tangent
(the key point is that the coefficient of 
$\frac{\partial}{\partial u}$ in the above expression for $\Rad$ is a constant!).
That is, we have shown that
$[\Singletan_1,\Singletan_2] \sim \Singletan_3$
and
$[\Rad,\Singletan_1] \sim \Singletan_2$.
Recalling the wave equation decomposition
\eqref{E:INTROWAVEEQUATIONFRAMEDECOMPOSED},
we easily obtain the structure 
for the commutators $[\upmu \square_{g(\Psi)},\Singletan]$
highlighted in the above indented sentence 
in the special case $\Singletan \in \lbrace \Lunit, \GeoAng \rbrace$
relevant for our energy estimates.\footnote{The detailed proof of the structure
of the commutators $[\upmu \square_{g(\Psi)},\Singletan]$
for $\Singletan \in \lbrace \Lunit, \GeoAng \rbrace$,
in the precise form that we need for our proof,
is based on straightforward but lengthy geometric computations 
carried out in
Lemma~\ref{L:BOXZCOM},
Prop.~\ref{P:COMMUTATIONCURRENTDIVERGENCEFRAMEDECOMP} with $Z \in \lbrace \Lunit, \GeoAng \rbrace$,
and Lemma~\ref{L:DEFORMATIONTENSORFRAMECOMPONENTS}.
}
The structure is a manifestation of the miraculous
null structure mentioned in the discussion surrounding
equation \eqref{E:INTROWAVEEQUATIONFRAMEDECOMPOSED},
and it allows us to derive energy estimates for the
$\Lunit$ and $\GeoAng$ derivatives of the solution up to top order
without having to derive energy estimates for its
high $\Rad$ derivatives. Put differently, 
there is a kind of decoupling
between energy estimates for the $\mathcal{P}_u-$tangential derivatives
and the $\mathcal{P}_u-$transversal derivatives.
Moreover, the structure has the following important consequence:
all energy estimate error integrands generated by commuting 
the wave equation with $\Lunit$ and $\GeoAng$
\emph{contain at most one $\mathring{\updelta}-$sized 
factor and thus are at least quadratically small in the 
quantities that are expected to be of size $\mathring{\upepsilon}$}.
This suggests that a Gronwall estimate will lead to the $C \mathring{\upepsilon}^2$ 
smallness of the energies for\footnote{As we explain in Subsect.~\ref{SS:NOTATION},
we use the convention that constants $C$ are allowed to depend on
$\mathring{\updelta}$
and
$\TranminusdatasizeWithFactor^{-1}$.} 
the relevant time scale $t < 2 \TranminusdatasizeWithFactor^{-1}$,
as described in Subsect.~\ref{SS:PROOFOVERVIEW}.
Indeed, modulo the many difficulties with high-order energy degeneracy with respect to $\upmu_{\star}^{-1}$ 
that we previously explained,
this is exactly what our energy estimate hierarchy \eqref{E:INTROTOPENERGY}-\eqref{E:INTROLOWESTENERGY} reveals.
This allows us to propagate the $\mathcal{O}(\mathring{\upepsilon}^2)$ smallness of
the energies of the $\mathcal{P}_u-$tangent
derivatives of $\Psi$ without having to bound the energies\footnote{We note, however, the following
non-obvious feature of our proof, described at the start of Sect.~\ref{S:LINFINITYESTIMATESFORHIGHERTRANSVERSAL}:
to close our energy estimates at any order,
we rely on the bound $\| \Rad^{\leq 3} \Psi \|_{L^{\infty}(\Sigma_t^u)} \lesssim \mathring{\updelta} \lesssim 1$,
which we obtain by commuting the wave equation up to two times with $\Rad$ and treating the wave equation 
as a transport equation up to derivative-losing terms.}
of the pure transversal derivatives such as $\Rad \Psi$, $\Rad \Rad \Psi$, etc.,
which can be of large size $\mathcal{O}(\mathring{\updelta}^2)$.

For illustration, we now give one example of how the $\mathcal{O}(\mathring{\upepsilon}^2)$ 
smallness of the energies is used in our proof. We recall that our bootstrap argument heavily relies on the
expectation (described just below equation \eqref{E:INTROMUDECAYESTIMATE})
that the first vanishing time of $\upmu$ 
(that is, the blowup-time of the first rectangular derivatives of $\Psi$) is 
$(1 + \mathcal{O}(\mathring{\upepsilon}))\TranminusdatasizeWithFactor^{-1}$.
To realize this expectation, we must show that
the $\Lunit \Psi-$involving products on RHS~\eqref{E:INTROMUEVOLUTION}
are of small size $\mathcal{O}(\mathring{\upepsilon})$ all the way up to the shock.
The desired smallness estimate
$\| \Lunit \Psi \|_{L^{\infty}(\Sigma_t^u)} \lesssim \mathring{\upepsilon}$
is a simple consequence
of the $\mathcal{O}(\mathring{\upepsilon}^2)$ smallness of the low-order energies,
a data smallness assumption,
and Sobolev embedding;
see Cor.~\ref{C:PSILINFTYINTERMSOFENERGIES} for a proof.

We now further explain how
the analysis of the small-data dispersive regime \cites{dC2007,dCsM2014,jS2014b}
is different than our analysis here.
In that regime, \emph{there is only one smallness parameter capturing the size of 
a full spanning set of directional derivatives of the solution at time $0$}, 
and the $\mathring{\upepsilon}^2$ smallness of all energies from level $0$ up to top order
can be propagated all the way up to the shock 
(modulo possible energy degeneracy relative to powers of $\upmu_{\star}^{-1}$ at the high orders).
Because all directional derivatives are controlled, 
there is no need to rely on the structure emphasized two paragraphs above, namely that
the energy estimates for the pure tangential derivatives (up to top order) 
effectively decouple from energy estimates for transversal derivatives.
The good null structure mentioned above does, however, play an important role in 
allowing one to control error terms and prove shock formation. 
The structure is used in a different way:
in place of the two-parameter $\mathring{\upepsilon}-\mathring{\updelta}$ hierarchy
exploited in the present article,
the error terms are controlled all the way up to the shock 
via a hierarchy of dispersive estimates.
More precisely, one relies on the fact that 
\emph{the transversal derivative of the solution decays in time at a non-integrable rate tied to the formation of a shock}, 
while the tangential derivatives decay at an integrable rate and generate only small error terms;
see the next paragraph for more details.
The availability of this decay hierarchy is intimately connected to the good null structure, 
and we explain it more detail two paragraphs below.

For the sake of comparison, we first provide some additional background on the behavior of solutions
in the small-data dispersive regime \cites{dC2007,dCsM2014,jS2014b}.
The data are compactly supported functions on $\mathbb{R}^3$
of small Sobolev\footnote{The work \cite{jS2014b} showed that for equations of type \eqref{E:GEOWAVE}, 
the proof closes for small data verifying
$(\mathring{\Psi},\mathring{\Psi}_0) \in H_{\Euct}^{25}(\Sigma_0) \times H_{\Euct}^{24}(\Sigma_0)$.
} 
size $\mathring{\upepsilon}$,
and the characteristics are outgoing null cones $\mathcal{C}_u$.
The $\mathcal{C}_u$, which are level sets of an eikonal function $u$,
are distorted versions of the Minkowskian cones $\lbrace t - r = \mbox{const} \rbrace$,
where $r$ is the standard radial coordinate on Minkowski spacetime.
The dispersive estimates take the following form:
\emph{relative to a suitable rescaled vectorfield frame} analogous to \eqref{E:INTROFRAME},
$\Psi$ and its $\mathcal{C}_u-$transversal derivative
decay like $\mathring{\upepsilon}(1 + t)^{-1}$
while its $\mathcal{C}_u-$tangential derivatives decay 
at the faster rate $\mathring{\upepsilon}(1 + t)^{-2}$.
Moreover, relative to the geometric coordinates,
related estimates hold for $\Psi$ at slightly higher derivative levels 
and for the low-order derivatives of $\upmu$
and the rectangular components $\Lunit^i$.
We now describe the mechanism for the vanishing of $\upmu$ 
(that is, for the formation of a shock)
in the small-data dispersive regime.
The most relevant estimate takes the form
$
\Lunit \upmu 
= \mathcal{O}(\mathring{\upepsilon})(1 + t)^{-1}
+ \mathcal{O}(\mathring{\upepsilon}) (1 + t)^{-2}
$
and is analogous to the estimate \eqref{E:INTROMUEVOLUTIONESTIAMTED}
in this paper.
The term $\mathcal{O}(\mathring{\upepsilon})(1 + t)^{-1}$ corresponds to the size of the $\mathcal{C}_u-$transversal derivative
of the solution, while the term $\mathcal{O}(\mathring{\upepsilon}) (1 + t)^{-2}$
is an error term that bounds the $\mathcal{C}_u-$tangential derivatives.
In view of the fact that 
$
\displaystyle
\Lunit = \frac{\partial}{\partial t}
$,
the small-data estimate $\upmu|_{t=0} \sim 1$,
and the observation that
$(1 + t)^{-1}$ is not integrable in $t$ 
while $(1 + t)^{-2}$ is,
we see that
$\upmu \sim 1 + \mathcal{O}(\mathring{\upepsilon}) \ln(1 + t)$.
Hence, $\upmu$ will vanish at a time
$ 
\displaystyle
\exp
\left(
	|\mathcal{O}(\mathring{\upepsilon})|^{-1}
\right)
$
for data such that the factor 
$\mathcal{O}(\mathring{\upepsilon})$
from the term
$
\mathcal{O}(\mathring{\upepsilon})(1 + t)^{-1}
$
is negative and sufficiently bounded from below in magnitude.\footnote{The precise behavior of the 
$
\mathcal{O}(\mathring{\upepsilon})(1 + t)^{-1}
$
term depends on the nonlinearities as well as the profile of the data
and is connected to Friedlander's radiation field;
see \cites{jS2014b,gHsKjSwW2016} for more details.}

The derivation of the above mentioned
\emph{directionally dependent} decay rates 
in the small-data dispersive regime
is based on a modified version of Klainerman's commuting vectorfield method \cite{sK1985},
the modification being that the vectorfields are dynamically adapted to the 
characteristics through an eikonal function,
much like the vectorfields $\Fullset$ that we use in the present article 
(as described at the start of Subsect.~\ref{SS:OVERVIEW}).
As we mentioned previously, the use of an eikonal function in the context
of deriving global estimates for quasilinear hyperbolic equations
originated in 
\cite{dCsK1993}. 
In the small-data dispersive regime,
one can exploit the decay properties mentioned above, 
the good null structure mentioned in
the discussion surrounding equation \eqref{E:INTROWAVEEQUATIONFRAMEDECOMPOSED},
and various structures present in the evolution equations for 
$\upmu$ and $\Lunit^i$
to show that the solution behaves,
relative to the rescaled frame, much like a
solution to a wave equation that verifies Klainerman's classic null condition.
In particular, upon commuting the wave equation
$\upmu \square_{g(\Psi)} \Psi = 0$ with an appropriate spanning commutation set,
one can show that the commutator error terms are
quadratic and higher-order products such that
\emph{each product contains no more than one slowly decaying factor 
corresponding to pure $\mathcal{C}_u-$transversal differentiations}.
This is an analog, for a full spanning set of commutation vectorfields, 
of the structure described in the second paragraph
of this subsubsection for the $\Lunit$ and $\GeoAng$ commutation vectorfields 
in our case.
Moreover, in the small-data dispersive regime,
relative to the rescaled frame, one can propagate the 
$\mathring{\upepsilon}$ smallness of the solution in various Sobolev norms
and prove \emph{conditional} global existence and decay-type estimates.
In particular, \emph{without any a priori restriction on $t$}
(such as the restriction $t < 2 \TranminusdatasizeWithFactor^{-1}$ made in our work here),
one can prove that the solution remains regular
relative to both the geometric coordinates 
and the rectangular coordinates
as long as $\upmu$ remains strictly positive.
We mention again that in contrast, 
in the nearly plane symmetric regime, 
there is no obvious structure in the equations hinting at
the validity of a conditional global existence-type result in which
the solution persists for all times as long as $\upmu$ remains strictly positive.
Rather, as we explained in Subsect.~\ref{SS:PROOFOVERVIEW},
we propagate the $\mathring{\upepsilon}-\mathring{\updelta}$ hierarchy only for 
times up to $2 \TranminusdatasizeWithFactor^{-1}$,
which is long enough for the shock to form.

\subsubsection{Blowup in a large-data regime featuring a one-parameter scaling of the data}
\label{SSS:LARGEDATASHOCKS}
Recently, Miao and Yu proved \cite{sMpY2014}
a related shock formation result 
for the wave equation 
$-\partial_t^2 \phi + [1 + (\partial_t \phi)^2] \Delta \phi = 0
$
in three spatial dimensions
with data that are compactly supported in an annular region of radius $\approx 1$ and thin width $\updelta$,
where $\updelta$ is a small positive parameter.
The data's amplitude and their functional dependence on a radial coordinate
are rescaled by powers of $\updelta$.
Consequently, the data and their derivatives
verify a hierarchy of estimates featuring various powers of $\updelta$.
For example, $\phi$ itself has small $L^{\infty}$ size $\updelta^{3/2}$, 
its rectangular derivatives $\partial_{\alpha} \phi$ have $L^{\infty}$ size $\delta^{1/2}$,
and a certain derivative of $\partial_{\alpha} \phi$ that is transversal to the characteristics
has \emph{large} $L^{\infty}$ size $\updelta^{-1/2}$.
Due to the largeness, the blowup of the second rectangular derivatives of $\phi$
happens within one unit of time.
The scaling of the data is closely related to the short-pulse ansatz
pioneered by Christodoulou in his aforementioned
proof of the formation of trapped surfaces
in solutions to the Einstein-vacuum equations \cite{dC2009}.
The main contribution of \cite{sMpY2014}  
was showing how to propagate the $\updelta$ hierarchy estimates
until the time of first shock formation. 
In the proof, dispersive effects are not relevant. Instead, the authors control 
nonlinear error terms by tracking the powers of $\updelta$
associated to each factor in the product.
Roughly, the error terms have a product structure, typically of the form
$\mbox{\upshape small} \cdot \mbox{\upshape large}$ 
(relative to powers of $\updelta$),
where the small factor often more than compensates for the large one.
That is, the authors show that 
the overall powers of $\updelta$
associated to the error term products are favorable 
in the sense that the smallness of $\updelta$ is sufficient for controlling them.
In this way, a class of large data solutions can be treated using techniques
borrowed from the usual small-data framework.

Our results are related to those of \cite{sMpY2014}
but are distinguished by our use of two size parameters
(the parameters $\mathring{\upepsilon}$ and $\mathring{\updelta}$
from Subsect.~\ref{SS:PROOFOVERVIEW} and Subsubsect.~\ref{SSS:INTROPOOFTHATMUVANISHES}),
which allows us to treat a set of initial conditions containing large data
and, unlike \cite{sMpY2014}, small data too. 
As we described above, a key aspect of our proof is that we can propagate the small size
$\mathring{\upepsilon}$ of the $\mathcal{P}_u-$tangent derivatives
long enough for the shock to form, 
even though the transversal derivatives can be of a relatively large size
$\mathring{\updelta}$.
To this end, we must exploit the good product/null structure in the equations,
as described in Subsubsect.~\ref{SSS:DIFFERENCESFROMDISPERSIVEREGIME},
in ways that go beyond the $\updelta$ scaling structures exploited in \cite{sMpY2014}.

\section{Geometric Setup}
\label{S:GEOMETRICSETUP}
In this section, we set up the geometric framework that we use for analyzing solutions. 
We note that most of the basic geometric insights are present in \cite{dC2007} 
and that the calculations in this section have analogs in \cite{dC2007}. 
For the reader's convenience, we re-derive 
the relevant results and adapt them in our setting.
Similar remarks apply throughout the article
(see Subsubsect.~\ref{SSS:DIFFERENCESFROMDISPERSIVEREGIME}
for an overview of the main new ideas of the present work).
We also note that for pedagogical reasons, 
there is some redundancy with Sect.~\ref{S:INTRO}.

\subsection{Notational conventions and shorthand notation}
\label{SS:NOTATION}
We start by summarizing some of our notational conventions;
the precise definitions of some of the concepts referred to
here are provided later in the article.

\begin{itemize}
	\item Lowercase Greek spacetime indices 
	$\alpha$, $\beta$, etc.\
	correspond to the rectangular spacetime coordinates 
	defined in Subsect.~\ref{SS:EQUATIONINRECTANGULARCOMPONENTS}
	and vary over $0,1,2$.
	Lowercase Latin spatial indices
	$a$,$b$, etc.\ 
	correspond to the rectangular spatial coordinates and vary over $1,2$.
	All lowercase Greek indices are lowered and raised with the spacetime metric
	$g$ and its inverse $g^{-1}$, and \emph{not with the Minkowski metric}.
\item We sometimes use $\cdot$ to denote the natural contraction between two tensors
		(and thus raising or lowering indices with a metric is not needed). 
		For example, if $\xi$ is a spacetime one-form and $V$ is a 
		spacetime vectorfield,
		then $\xi \cdot V := \xi_{\alpha} V^{\alpha}$.
\item If $\xi$ is a one-form and $V$ is a vectorfield, then
	$\xi_V := \xi_{\alpha} V^{\alpha}$. 
	Similarly, if $W$ is a vectorfield, then
	$W_V := W_{\alpha} V^{\alpha} = g(W,V)$.
	We use similar notation when contracting higher-order tensorfields
	against vectorfields. Similarly, if $\Gamma_{\alpha \kappa \beta}$
	are the rectangular Christoffel symbols \eqref{E:CHRISTOFEELRECT}, then
	$\Gamma_{UVW} := U^{\alpha} V^{\kappa} W^{\beta} \Gamma_{\alpha \kappa \beta}$.
\item If $\xi$ is an $\ell_{t,u}-$tangent one-form
	(as defined in Subsect.~\ref{SS:PROJECTIONTENSORFIELDANDPROJECTEDLIEDERIVATIVES}),
	then $\xi^{\#}$ denotes its $\gsphere-$dual vectorfield,
	where $\gsphere$ is the Riemannian metric induced on $\ell_{t,u}$ by $g$.
	Similarly, if $\xi$ is a symmetric type $\binom{0}{2}$ $\ell_{t,u}-$tangent tensor, 
	then $\xi^{\#}$ denotes the type $\binom{1}{1}$ $\ell_{t,u}-$tangent tensor formed by raising one index with $\ginversesphere$
	and $\xi^{\# \#}$ denotes the type $\binom{2}{0}$ $\ell_{t,u}-$tangent tensor formed by raising both indices with $\ginversesphere$.
\item If $\xi$ is an $\ell_{t,u}-$tangent tensor, then 
	the norm $|\xi|$ is defined relative to the Riemannian metric $\gsphere$, as in Def.~\ref{D:POINTWISENORM}.
\item Unless otherwise indicated, 
	all quantities in our estimates that are not explicitly under
	an integral are viewed as functions of 
	the geometric coordinates $(t,u,\vartheta)$
	of Def.~\ref{D:GEOMETRICCOORDINATES}.
	Unless otherwise indicated, quantities
	under integrals have the functional dependence 
	established below in
	Def.~\ref{D:NONDEGENERATEVOLUMEFORMS}.
\item If $Q_1$ and $Q_2$ are two operators, then
	$[Q_1,Q_2] = Q_1 Q_2 - Q_2 Q_1$ denotes their commutator.
\item $A \lesssim B$ means that there exists $C > 0$ such that $A \leq C B$.
\item $A = \mathcal{O}(B)$ means that $|A| \lesssim |B|$.
\item Constants such as $C$ and $c$ are free to vary from line to line.
	\textbf{Explicit and implicit constants are allowed to depend in an increasing, 
	continuous fashion on the data-size parameters 
	$\mathring{\updelta}$
	and $\TranminusdatasizeWithFactor^{-1}$
	from
	Subsect.~\ref{SS:DATAASSUMPTIONS}.
	However, the constants can be chosen to be 
	independent of the parameters $\mathring{\upepsilon}$ 
	and $\varepsilon$ whenever $\mathring{\upepsilon}$ and $\varepsilon$
	are sufficiently small relative to 
	$\mathring{\updelta}^{-1}$
	and $\TranminusdatasizeWithFactor$.
	}
\item $\lfloor \cdot \rfloor$
	and $\lceil \cdot \rceil$
	respectively denote the floor and ceiling functions. 
\end{itemize}

\subsection{The structure of the equation in rectangular components}
\label{SS:EQUATIONINRECTANGULARCOMPONENTS}
In this subsection, we formulate equation \eqref{E:GEOWAVE} in rectangular coordinates
and state our assumptions on the nonlinear terms.
We use $t = x^0 \in \mathbb{R}$ to denote the time coordinate
and $(x^1,x^2) \in \mathbb{R} \times \mathbb{T}$ 
to denote standard coordinates on $\Sigma$, where $x^2$ is locally defined.
The vectorfields $\partial_t$, $\partial_1$, $\partial_2$ are globally defined. 
We call $\lbrace x^{\alpha} \rbrace_{\alpha = 0,1,2}$ the \emph{rectangular coordinates}
because relative to them, the standard Minkowski metric on $\mathbb{R} \times \Sigma$
takes the form $m_{\mu \nu} = \mbox{diag}(-1,1,1)$.

We assume that relative to the rectangular coordinates,
\begin{align} \label{E:LITTLEGDECOMPOSED}
	g_{\mu \nu} 
	= g_{\mu \nu}(\Psi)
	& := m_{\mu \nu} 
		+ g_{\mu \nu}^{(Small)}(\Psi),
	&& (\mu, \nu = 0,1,2),
\end{align}
where $g_{\mu \nu}^{(Small)}(\Psi)$ is a given smooth function of $\Psi$ with
\begin{align} \label{E:METRICPERTURBATIONFUNCTION}
	g_{\mu \nu}^{(Small)}(0)
	& = 0.
\end{align}
Relative to the rectangular coordinates, \eqref{E:GEOWAVE} takes the form
\begin{align} \label{E:WAVEEQUATIONINCOORDINATES}
	(g^{-1})^{\alpha \beta} \partial_{\alpha} \partial_{\beta} \Psi
	- (g^{-1})^{\alpha \beta} (g^{-1})^{\kappa \lambda} \Gamma_{\alpha \kappa \beta} \partial_{\lambda} \Psi
	& = 0.
\end{align}
The $\Gamma_{\alpha \kappa \beta}$ are the lowered
Christoffel symbols\footnote{Our Christoffel symbol index conventions are such that for vectorfields $V$, we have
$\D_{\alpha} V^{\beta} = \partial_{\alpha} V^{\beta} + \Gamma_{\alpha \ \lambda}^{\ \beta} V^{\lambda}$,
where $\Gamma_{\alpha \ \lambda}^{\ \beta} := (g^{-1})^{\beta \kappa} \Gamma_{\alpha \kappa \lambda}$.} 
of $g$ relative to rectangular coordinates
and can be expressed as
\begin{align} \label{E:CHRISTRECT}
	\Gamma_{\alpha \kappa \beta}
	= \Gamma_{\alpha \kappa \beta}(\Psi,\partial \Psi)
	& := 
		\frac{1}{2}
		\left\lbrace
			\partial_{\alpha} g_{\kappa \beta}
			+ \partial_{\beta} g_{\alpha \kappa}
			- \partial_{\kappa} g_{\alpha \beta}
		\right\rbrace
			\\
	& = 
		\frac{1}{2}
		\left\lbrace
			G_{\kappa \beta} \partial_{\alpha} \Psi
			+ G_{\alpha \kappa} \partial_{\beta} \Psi
			- G_{\alpha \beta} \partial_{\kappa} \Psi
		\right\rbrace,
		\notag
\end{align}
where 
\begin{align} \label{E:BIGGDEF}
	G_{\alpha \beta}
	= G_{\alpha \beta}(\Psi)
	& := \frac{d}{d \Psi} g_{\alpha \beta}(\Psi).
\end{align}
For later use, we also define
\begin{align} \label{E:BIGGPRIMEDEF}
	G_{\alpha \beta}'
	= G_{\alpha \beta}'(\Psi)
	& := \frac{d^2}{d \Psi^2} g_{\alpha \beta}(\Psi).
\end{align}

We now describe our assumptions on the tensorfield $G_{\alpha \beta}(\Psi = 0)$, 
which can be viewed as a $3 \times 3$ matrix with \emph{constant entries} relative to rectangular coordinates.
We could prove the existence\footnote{The condition \eqref{E:NONVANISHINGNONLINEARCOEFFICIENT} would 
be sufficient for allowing us to prove the existence of stable large-data shock-forming solutions.
However, in order to handle the set of data (which includes some small data)
stated in Theorem~\ref{T:MAINTHEOREM},
we need the additional assumption \eqref{E:LFLAT}.
} 
of stable shock-forming solutions whenever
there exists a Minkowski-null vectorfield $\Lunit_{(Flat)}$
(that is, $m_{\alpha \beta}\Lunit_{(Flat)}^{\alpha} \Lunit_{(Flat)}^{\beta} = 0$)
such that
\begin{align} \label{E:NONVANISHINGNONLINEARCOEFFICIENT}
	G_{\alpha \beta}(\Psi = 0) \Lunit_{(Flat)}^{\alpha} \Lunit_{(Flat)}^{\beta} \neq 0.
\end{align}
The assumption \eqref{E:NONVANISHINGNONLINEARCOEFFICIENT} holds for most nonlinearities 
and is equivalent to the failure of Klainerman's classic null condition \cite{sk1984}.
We recall that the main results that we present in this article rely on the 
existence of a family of plane symmetric shock-forming solutions. 
The existence of the family is based on the following assumption:
there exists a vectorfield
$\Lunit_{(Flat)} \in \mbox{\upshape span} \lbrace \partial_t, \partial_1 \rbrace$
such that \eqref{E:NONVANISHINGNONLINEARCOEFFICIENT} holds.
We may then perform a Lorentz transformation on the $t,x^1$ coordinates if necessary
in order to put $\Lunit_{(Flat)}$ into the following form, which we assume throughout
the remainder of the article:
\begin{align} \label{E:LFLAT}
	\Lunit_{(Flat)} = \partial_t + \partial_1.
\end{align}
Note that under the above assumptions, 
LHS $\mbox{\eqref{E:NONVANISHINGNONLINEARCOEFFICIENT}}$ 
is equal to the non-zero constant
$G_{00}(\Psi = 0) + 2 G_{01}(\Psi = 0) + G_{11}(\Psi = 0)$.

\begin{remark}[\textbf{Genuinely nonlinear systems}]
\label{R:GENUINELYNONLINEAR}
	Our assumption that the vectorfield \eqref{E:LFLAT}
	verifies \eqref{E:NONVANISHINGNONLINEARCOEFFICIENT}
	is reminiscent of the
	well-known genuine nonlinearity condition for first-order strictly hyperbolic systems.
	In particular, for plane symmetric solutions with $\Psi$ sufficiently small, 
	the assumption ensures that there are quadratic 
	Riccati-type terms\footnote{The vectorfield frame that we construct in fact leads to the 
	cancellation of the Riccati-type terms; see the discussion just below \eqref{E:INTROWAVEEQUATIONFRAMEDECOMPOSED}.} 
	in the wave equation, which is the main mechanism
	driving the singularity formation in the $2 \times 2$
	genuinely nonlinear strictly hyperbolic systems studied by Lax
	\cite{pL1964}.

\end{remark}

By rescaling the metric by the scalar function $1/(g^{-1})^{00}(\Psi)$, 
we may assume without loss of generality\footnote{Technically, rescaling the metric 
introduces a semilinear term proportional to $(g^{-1})^{\alpha \beta}(\Psi) \partial_{\alpha} \Psi \partial_{\beta} \Psi$
in the covariant wave equation corresponding to the rescaled metric.
However, our proof will show that for the solutions under study, 
this term makes a negligible contribution to the dynamics because it has a special null structure
(it verifies the strong null condition mentioned in Remark~\ref{R:SLIGHTEXTENSIONNOBADTERMS})
that is visible relative to the frame \eqref{E:RESCALEDFRAME}.
Hence, we ignore it for simplicity.
\label{FN:GINVERSE00ISMINUSONE} } 
that
\begin{align} \label{E:GINVERSE00ISMINUSONE}
	(g^{-1})^{00}(\Psi) \equiv -1.
\end{align}
The assumption \eqref{E:GINVERSE00ISMINUSONE} simplifies many of our formulas.

\begin{remark}
	In total, our assumptions on the nonlinearities imply that the term 
	$\frac{1}{2} G_{\Lunit \Lunit} \Rad \Psi$
	on RHS~\eqref{E:UPMUFIRSTTRANSPORT},
	which lies at the heart of our analysis, is sufficiently
	strong to drive $\upmu$ to $0$ in the regime under study.
\end{remark}

\subsection{Basic constructions involving the eikonal function}
\label{SS:EIKONALFNCNBASICCONSTRUCTIONS}
As we described in Subsect.~\ref{SS:PROOFOVERVIEW}, 
our entire work is based on an eikonal function,
specifically, the solution to the hyperbolic initial value problem
\eqref{E:INTROEIKONAL}-\eqref{E:INTROEIKONALINITIALVALUE}.
We associate the following subsets of spacetime to $u$.
They were depicted in Figure~\ref{F:SOLIDREGION} on pg.~\pageref{F:SOLIDREGION}.

\begin{definition} [\textbf{Subsets of spacetime}]
\label{D:HYPERSURFACESANDCONICALREGIONS}
We define the following spacetime subsets:
\begin{subequations}
\begin{align}
	\Sigma_{t'} & := \lbrace (t,x^1,x^2) \in \mathbb{R} \times \mathbb{R} \times \mathbb{T}  
		\ | \ t = t' \rbrace, 
		\label{E:SIGMAT} \\
	\Sigma_{t'}^{u'} & := \lbrace (t,x^1,x^2) \in \mathbb{R} \times \mathbb{R} \times \mathbb{T} 
		 \ | \ t = t', \ 0 \leq u(t,x^1,x^2) \leq u' \rbrace, 
		\label{E:SIGMATU} 
		\\
	\mathcal{P}_{u'}^{t'} 
	& := 
		\lbrace (t,x^1,x^2) \in \mathbb{R} \times \mathbb{R} \times \mathbb{T} 
			\ | \ 0 \leq t \leq t', \ u(t,x^1,x^2) = u' 
		\rbrace, 
		\label{E:PUT} \\
	\ell_{t',u'} 
		&:= \mathcal{P}_{u'}^{t'} \cap \Sigma_{t'}^{u'}
		= \lbrace (t,x^1,x^2) \in \mathbb{R} \times \mathbb{R} \times \mathbb{T} 
			\ | \ t = t', \ u(t,x^1,x^2) = u' \rbrace, 
			\label{E:LTU} \\
	\mathcal{M}_{t',u'} & := \cup_{u \in [0,u']} \mathcal{P}_u^{t'} \cap \lbrace (t,x^1,x^2) 
		\in \mathbb{R} \times \mathbb{R} \times \mathbb{T}  \ | \ 0 \leq t < t' \rbrace.
		\label{E:MTUDEF}
\end{align}
\end{subequations}
\end{definition}

We refer to the $\Sigma_t$ and $\Sigma_t^u$ as ``constant time slices,'' 
the $\mathcal{P}_u^t$ as ``null hyperplanes,''
and the $\ell_{t,u}$ as ``curves.'' 
We sometimes use the notation $\mathcal{P}_u$ in place of $\mathcal{P}_u^t$ 
when we are not concerned with the truncation time $t$.
We restrict our attention to spacetime regions with 
$0 \leq u \leq U_0$, where 
we recall (see \eqref{E:FIXEDPARAMETER}) that 
$0 < U_0 \leq 1$ 
is a parameter,
fixed until Theorem~\ref{T:MAINTHEOREM}.

\begin{remark}
	\label{R:CONSTANTSDONOTDEPENDONU0}
	The constants in all of our estimates can be chosen to be independent of 
	$U_0 \in (0,1]$.
\end{remark}

We associate the following gradient vectorfield
to the eikonal function solution to \eqref{E:INTROEIKONAL}:
\begin{align} \label{E:LGEOEQUATION}
	\Lgeo^{\nu} & := - (g^{-1})^{\nu \alpha} \partial_{\alpha} u.
\end{align}
It is easy to see that $\Lgeo$ is future-directed\footnote{Here and throughout, 
a vectorfield $V$ is ``future-directed'' if its rectangular component $V^0$ is positive.} 
with
\begin{align}  \label{E:LGEOISNULL}
g(\Lgeo,\Lgeo) 
:= g_{\alpha \beta} \Lgeo^{\alpha} \Lgeo^{\beta}
= 0,
\end{align}
that is, $\Lgeo$ is $g-$null.
Moreover, we can differentiate the eikonal equation with
$\D^{\nu} := (g^{-1})^{\nu \alpha} \D_{\alpha}$
and use the torsion-free property of the connection $\D$ to deduce that
$
0
= 
(g^{-1})^{\alpha \beta} \D_{\alpha} u \D_{\beta} \D^{\nu} u
= - \D^{\alpha} u \D_{\alpha} \Lgeo^{\nu}
= \Lgeo^{\alpha} \D_{\alpha} \Lgeo^{\nu}
$. That is, $\Lgeo$ is geodesic:
\begin{align} \label{E:LGEOISGEODESIC}
	\D_{\Lgeo} \Lgeo & = 0.
\end{align}
In addition, since $\Lgeo$ is proportional to the metric dual of the one-form $d u$,
which is co-normal to the level sets $\mathcal{P}_u$ of the eikonal function,
it follows that $\Lgeo$ is $g-$orthogonal to $\mathcal{P}_u$.
Hence, the $\mathcal{P}_u$ have null normals. 
Such hypersurfaces are known as \emph{null hypersurfaces}.
Our analysis will show that the rectangular components of $\Lgeo$
blow up when the shock forms. 
In particular, 
as we described in Subsect.~\ref{SS:PROOFOVERVIEW}, 
the formation of a shock
is equivalent to the vanishing of the following quantity $\upmu$.

\begin{definition}[\textbf{Inverse foliation density}]
 \label{D:UPMUDEF}
Let $\Lgeo^0$ be the $0$ rectangular component of the vectorfield $\Lgeo$ defined in \eqref{E:LGEOEQUATION}.
We define the inverse foliation density $\upmu$ as follows:
\begin{align} \label{E:UPMUDEF}
	\upmu 
	& := 
		\frac{-1}{(g^{-1})^{\alpha \beta} \partial_{\alpha} t \partial_{\beta} u} 
	= \frac{-1}{(g^{-1})^{0 \alpha} \partial_{\alpha} u} 
	= \frac{1}{\Lgeo^0}.
\end{align}
\end{definition}
The quantity $1/\upmu$ measures the density of the level sets of $u$
relative to the constant-time hypersurfaces $\Sigma_t$. When $\upmu$
becomes $0$, the density becomes infinite and the level sets of $u$
intersect. We also note that the vanishing of $\upmu$
is equivalent to the blow up of $\D_V u$, where
$V := - (g^{-1})^{0 \alpha} \partial_{\alpha}$
is approximately equal to $\partial_t$ in the regime under study.

In our analysis, we work with a rescaled version of
$\Lgeo$ that we refer to as $\Lunit$.  
Our proof reveals that the rectangular components of $\Lunit$
remain near those of $\Lunit_{(Flat)}$ 
(which is defined in \eqref{E:LFLAT})
all the way up to the shock.
\begin{definition}[\textbf{Rescaled null vectorfield}]
	\label{D:LUNITDEF}
	We define the rescaled null (see \eqref{E:LGEOISNULL}) vectorfield $\Lunit$ as follows:
	\begin{align} \label{E:LUNITDEF}
		\Lunit
		& := \upmu \Lgeo.
	\end{align}
\end{definition}

\begin{definition}[\textbf{Geometric torus coordinate} $\vartheta$ \textbf{and the corresponding vectorfield} $\CoordAng$]
	\label{D:ANGULARCOORDINATE}
	Along $\Sigma_0^1$, we define $\vartheta(t=0,x^1,x^2) = x^2$.
 	We extend $\vartheta$ to regions of the form $\mathcal{M}_{t,u}$
 	by solving the transport equation $\Lunit \vartheta = 0$
 	with $\vartheta$ subject to the above 
 	initial conditions along $\Sigma_0^1$.
 
 	We define 
 	$\CoordAng 
 	=
 	\frac{\partial}{\partial \vartheta} 
 	:= \frac{\partial}{\partial \vartheta}|_{t,u}
 	$
 	to be the vectorfield corresponding to partial differentiation with respect
 	to $\vartheta$ at fixed $t$ and $u$.
\end{definition}

\begin{definition}[\textbf{Geometric coordinates}]
	\label{D:GEOMETRICCOORDINATES}
	We refer to $(t,u,\vartheta)$ as the geometric coordinates.
\end{definition}

\begin{remark}[\textbf{$C^1-$equivalent differential structures until shock formation}]
We often identify
spacetime regions of the form
$\mathcal{M}_{t,U_0}$ 
(see \eqref{E:MTUDEF})
with the region
$[0,t) \times [0,U_0] \times \mathbb{T}$
corresponding to the geometric coordinates.
This identification is 
justified by the fact that 
during the classical lifespan of the solutions under consideration,
the differential structure on 
$\mathcal{M}_{t,U_0}$ corresponding to the 
geometric coordinates is $C^1-$equivalent
to the differential structure on 
$\mathcal{M}_{t,U_0}$
corresponding to the rectangular coordinates.
The equivalence is captured by the fact that the change
of variables map $\Upsilon$ (see Def.~\ref{D:CHOVMAP})
from geometric to rectangular coordinates is differentiable 
with a differentiable inverse, until a shock forms;
see Lemma~\ref{L:CHOVREMAINSADIFFEOMORPHISM} 
and Theorem~\ref{T:MAINTHEOREM}.
However, at points where $\upmu$ vanishes,
the rectangular derivatives of $\Psi$ blow up
(see equation \eqref{E:BLOWUPPOINTINFINITE} and the discussion below it),
the inverse map $\Upsilon^{-1}$ becomes singular,
and the equivalence of the differential structures breaks down as well.
\end{remark}

\subsection{Important vectorfields, the rescaled frame, and the non-rescaled frame}
\label{SS:FRAMEANDRELATEDVECTORFIELDS}
In this subsection, we define some 
additional vectorfields that we use in our analysis
and exhibit their basic properties.

\begin{definition}[$\Radunit$, $\Rad$, \textbf{and} $\Timenormal$  
]
	\label{D:RADANDXIDEFS}
	We define $\Radunit$ to be the unique vectorfield
	that is $\Sigma_t-$tangent, $g-$orthogonal 
	to the $\ell_{t,u}$, and normalized by
	\begin{align} \label{E:GLUNITRADUNITISMINUSONE}
		g(\Lunit,\Radunit) = -1.
	\end{align}
	We define
	\begin{align} \label{E:RADDEF}
		\Rad := \upmu \Radunit.
	\end{align}
		We define 
	\begin{align} \label{E:TIMENORMAL}
		\Timenormal 
		& := \Lunit + \Radunit.
	\end{align}
\end{definition}

\begin{definition}[\textbf{Two frames}]
	\label{D:RESCALEDFRAME}
	We define, respectively, the rescaled frame and the non-rescaled frame as follows:
	\begin{subequations}
	\begin{align} \label{E:RESCALEDFRAME}
		& \lbrace \Lunit, \Rad, \CoordAng \rbrace,
		&& \mbox{Rescaled frame},	\\
		&\lbrace \Lunit, \Radunit, \CoordAng \rbrace,
		&& \mbox{Non-rescaled frame}.
			\label{E:UNITFRAME}
	\end{align}
	\end{subequations}
\end{definition}


\begin{lemma}[\textbf{Basic properties of} $\Radunit$, $\Rad$, $\Lunit$, \textbf{and} $\Timenormal$]
\label{L:BASICPROPERTIESOFFRAME}
The following identities hold:
\begin{subequations}
\begin{align} \label{E:LUNITOFUANDT}
	\Lunit u & = 0, \qquad \Lunit t = \Lunit^0 = 1,
		\\
	\Rad u & = 1,
	\qquad \Rad t = \Rad^0 = 0, \label{E:RADOFUANDT}
\end{align}
\end{subequations}

\begin{subequations}
\begin{align} \label{E:RADIALVECTORFIELDSLENGTHS}
	g(\Radunit,\Radunit)
	& = 1,
	\qquad
	g(\Rad,\Rad)
	= \upmu^2,
		\\
	g(\Lunit,\Radunit)
	& = -1,
	\qquad
	g(\Lunit,\Rad) = -\upmu.
	\label{E:LRADIALVECTORFIELDSNORMALIZATIONS}
\end{align}
\end{subequations}

Moreover, relative to the geometric coordinates, we have
\begin{align} \label{E:LISDDT}
	\Lunit = \frac{\partial}{\partial t}.
\end{align}

In addition, there exists an $\ell_{t,u}-$tangent vectorfield
$\NonRadialRad = \XiCoordComp \CoordAng$ 
(where $\XiCoordComp$ is a scalar function)
such that
\begin{align} \label{E:RADSPLITINTOPARTTILAUANDXI}
	\Rad 
	& = \frac{\partial}{\partial u} - \NonRadialRad
	= \frac{\partial}{\partial u} - \XiCoordComp \CoordAng.
\end{align}

The vectorfield $\Timenormal$ defined in \eqref{E:TIMENORMAL} is future-directed, $g-$orthogonal
to $\Sigma_t$ and is normalized by
\begin{align} \label{E:TIMENORMALUNITLENGTH}
	g(\Timenormal,\Timenormal) 
	& = - 1.
\end{align}
Moreover, relative to rectangular coordinates, we have
(for $\nu = 0,1,2$):
\begin{align} \label{E:TIMENORMALRECTANGULAR}
	\Timenormal^{\nu} = - (g^{-1})^{0 \nu}.
\end{align}

	Finally, the following identities hold relative to the rectangular coordinates 
	(for $\nu = 0,1,2$):
	\begin{align}  \label{E:DOWNSTAIRSUPSTAIRSSRADUNITPLUSLUNITISAFUNCTIONOFPSI} 
		\Radunit_{\nu} 
		& = - \Lunit_{\nu} - \delta_{\nu}^0,
		\qquad
		\Radunit^{\nu}
		= - \Lunit^{\nu}
			- (g^{-1})^{0\nu},
	\end{align}
	where $\delta_{\nu}^0$ is the standard Kronecker delta.
\end{lemma}

\begin{proof}
	We first prove \eqref{E:LUNITOFUANDT}.
	We begin by using \eqref{E:INTROEIKONAL}, \eqref{E:LGEOEQUATION}, and \eqref{E:LUNITDEF}
	to deduce that 
	$
	\Lunit u = \Lunit^{\alpha} \partial_{\alpha} u 
	= - \upmu (g^{-1})^{\alpha \beta} \partial_{\alpha} u \partial_{\beta} u 	= 0
	$ as desired. The fact that $\Lunit t = 1$ is a simple consequence of
	\eqref{E:LGEOEQUATION}, \eqref{E:UPMUDEF}, and \eqref{E:LUNITDEF}.

	We now prove \eqref{E:RADOFUANDT}. We begin by using 
	\eqref{E:LGEOEQUATION}, 
	\eqref{E:LUNITDEF},
	\eqref{E:GLUNITRADUNITISMINUSONE},
	and 
	\eqref{E:RADDEF}
	to deduce that
	$\Rad u = \upmu \Radunit^{\alpha} \partial_{\alpha} u = - \Radunit^{\alpha} \Lunit_{\alpha} = - g(\Radunit,\Lunit) = 1$.
	The fact that $\Rad t = 0$ is an immediate consequence of the fact that by construction, $\Rad$ is $\Sigma_t-$tangent.

	\eqref{E:RADSPLITINTOPARTTILAUANDXI} then follows easily
	from \eqref{E:RADOFUANDT} and the fact that 
	$\frac{\partial}{\partial u}$
	and $\CoordAng$ span the tangent space of $\Sigma_t$ at each point.

	\eqref{E:LRADIALVECTORFIELDSNORMALIZATIONS} is an easy consequence of
	\eqref{E:GLUNITRADUNITISMINUSONE} and \eqref{E:RADDEF}.

	To derive the properties of 
	$\Timenormal$, we consider the vectorfield 
	$V^{\nu} := - (g^{-1})^{0 \nu}$, which is $g-$dual to
	the one-form with rectangular components $- \delta_{\nu}^0$ and therefore 
	$g-$orthogonal to $\Sigma_t$. By \eqref{E:GINVERSE00ISMINUSONE},
	$g(V,V) = (g^{-1})^{\alpha \beta} \delta_{\alpha}^0 \delta_{\beta}^0 = -1$,
	so $V$ is future-directed, timelike, and unit-length. In particular,
	$V$ belongs to the $g-$orthogonal complement of $\ell_{t,u}$,
	a space spanned by $\lbrace \Lunit, \Radunit \rbrace$. Thus, there exist
	scalars $a,b$ such that $V = a \Lunit + b \Radunit$. 
	Since $V t = V^0 = 1 = \Lunit t = \Lunit^0$ and since $\Radunit t = \Radunit^0 = 0$,
	we find that $a = 1$, that is, that $V = \Lunit + b \Radunit$.
	Taking the inner product of this expression with $\Radunit$ and using \eqref{E:LRADIALVECTORFIELDSNORMALIZATIONS}
	together with the fact that
	$\Radunit$ is $\Sigma_t-$tangent (and hence $g-$orthogonal to $V$),
	we find that $0 = - 1 + b g(\Radunit,\Radunit)$.
	Similarly, using \eqref{E:LRADIALVECTORFIELDSNORMALIZATIONS}, 
	the fact that $\Lunit$ is null, 
	and the previous identity,
	we compute that $-1 = g(V,V) = - 2b + b^2 g(\Radunit, \Radunit) = - 2b + b = - b$.
	It follows that $V = \Lunit + \Radunit := \Timenormal$ and
	$g(\Radunit,\Radunit) = 1$. 
	We have thus obtained the properties of $\Timenormal$ and obtained \eqref{E:TIMENORMALUNITLENGTH}, 
	\eqref{E:TIMENORMALRECTANGULAR}, and the first identity in \eqref{E:RADIALVECTORFIELDSLENGTHS}.
	The second identity in \eqref{E:RADIALVECTORFIELDSLENGTHS}
	follows easily from the first one and definition \eqref{E:RADDEF}.
	\eqref{E:DOWNSTAIRSUPSTAIRSSRADUNITPLUSLUNITISAFUNCTIONOFPSI} 
	follows from the definition \eqref{E:TIMENORMAL} of $\Timenormal$
	and from lowering the indices in \eqref{E:TIMENORMALRECTANGULAR} with $g$.

	To obtain \eqref{E:LISDDT}, we simply use \eqref{E:LUNITOFUANDT} and the fact that
	by construction, we have $\Lunit \vartheta = 0$ (see Def.~\ref{D:ANGULARCOORDINATE}).

\end{proof}


\subsection{Projection tensorfields, \texorpdfstring{$G_{(Frame)}$}{frame components}, and projected Lie derivatives}
\label{SS:PROJECTIONTENSORFIELDANDPROJECTEDLIEDERIVATIVES}
Many of our constructions involve projections onto 
$\Sigma_t$ and $\ell_{t,u}$.
\begin{definition}[\textbf{Projection tensorfields}]
We define the $\Sigma_t$ projection tensorfield $\Sigmatproject$
and the $\ell_{t,u}$ projection tensorfield
$\Lineproject$ relative to rectangular coordinates as follows:
\begin{subequations}
\begin{align} 
	\Sigmatproject_{\nu}^{\ \mu} 
	&:=	\delta_{\nu}^{\ \mu}
			- \Timenormal_{\nu} \Timenormal^{\mu} 
		= \delta_{\nu}^{\ \mu}
			+ \delta_{\nu}^{\ 0} \Lunit^{\mu}
			+ \delta_{\nu}^{\ 0} \Radunit^{\mu},
			\label{E:SIGMATPROJECTION} \\
	\Lineproject_{\nu}^{\ \mu} 
	&:=	\delta_{\nu}^{\ \mu}
			+ \Radunit_{\nu} \Lunit^{\mu} 
			+ \Lunit_{\nu} (\Lunit^{\mu} + \Radunit^{\mu}) 
		= \delta_{\nu}^{\ \mu}
			- \delta_{\nu}^{\ 0} \Lunit^{\mu} 
			+  \Lunit_{\nu} \Radunit^{\mu}.
			\label{E:LINEPROJECTION}
	\end{align}
\end{subequations}
\end{definition}

\begin{definition}[\textbf{Projections of tensorfields}]
Given any spacetime tensorfield $\xi$,
we define its $\Sigma_t$ projection $\Sigmatproject \xi$
and its $\ell_{t,u}$ projection $\Lineproject \xi$
as follows:
\begin{subequations}
\begin{align} 
(\Sigmatproject \xi)_{\nu_1 \cdots \nu_n}^{\mu_1 \cdots \mu_m}
& :=
	\Sigmatproject_{\widetilde{\mu}_1}^{\ \mu_1} \cdots \Sigmatproject_{\widetilde{\mu}_m}^{\ \mu_m}
	\Sigmatproject_{\nu_1}^{\ \widetilde{\nu}_1} \cdots \Sigmatproject_{\nu_n}^{\ \widetilde{\nu}_n} 
	\xi_{\widetilde{\nu}_1 \cdots \widetilde{\nu}_n}^{\widetilde{\mu}_1 \cdots \widetilde{\mu}_m},
		\\
(\Lineproject \xi)_{\nu_1 \cdots \nu_n}^{\mu_1 \cdots \mu_m}
& := 
	\Lineproject_{\widetilde{\mu}_1}^{\ \mu_1} \cdots \Lineproject_{\widetilde{\mu}_m}^{\ \mu_m}
	\Lineproject_{\nu_1}^{\ \widetilde{\nu}_1} \cdots \Lineproject_{\nu_n}^{\ \widetilde{\nu}_n} 
	\xi_{\widetilde{\nu}_1 \cdots \widetilde{\nu}_n}^{\widetilde{\mu}_1 \cdots \widetilde{\mu}_m}.
	\label{E:STUPROJECTIONOFATENSOR}
\end{align}
\end{subequations}
\end{definition}
We say that a spacetime tensorfield $\xi$ is $\Sigma_t-$tangent 
(respectively $\ell_{t,u}-$tangent)
if $\Sigmatproject \xi = \xi$
(respectively if $\Lineproject \xi = \xi$).
Alternatively, we say that $\xi$ is a
$\Sigma_t$ tensor (respectively $\ell_{t,u}$ tensor).


\begin{definition}[\textbf{$\ell_{t,u}$ projection notation}]
	\label{D:STUSLASHPROJECTIONNOTATION}
	If $\xi$ is a spacetime tensor, then we define
	\begin{align} \label{E:TENSORSTUPROJECTED}
		\angxi := \Lineproject \xi.
	\end{align}

	If $\xi$ is a symmetric type $\binom{0}{2}$ spacetime tensor and $V$ is a spacetime
	vectorfield, then we define
	\begin{align} \label{E:TENSORVECTORANDSTUPROJECTED}
		\angxiarg{V} 
		& := \Lineproject (\xi_V),
	\end{align}
	where $\xi_V$ is the spacetime one-form with 
	rectangular components $\xi_{\alpha \nu} V^{\alpha}$, $(\nu = 0,1,2)$.
\end{definition}

We often refer to the following arrays of $\ell_{t,u}-$tangent tensorfields in our analysis.

\begin{definition}[\textbf{Components of $G$ and $G'$ relative to the non-rescaled frame}]
	\label{D:GFRAMEARRAYS}
	We define
	\[
		G_{(Frame)} := \left( G_{\Lunit \Lunit}, G_{\Lunit \Radunit}, G_{\Radunit \Radunit}, \angGarg{\Lunit}, \angGarg{\Radunit}, \angG \right)
	\]
	to be the array of components of the tensorfield \eqref{E:BIGGDEF} relative 
	to the non-rescaled frame \eqref{E:UNITFRAME}.
	Similarly, we define $G_{(Frame)}'$ to be the analogous array for the tensorfield \eqref{E:BIGGPRIMEDEF}.
\end{definition}

\begin{definition}[\textbf{Lie derivatives}]
\label{D:LIEDERIVATIVE}
If $V^{\mu}$ is a spacetime vectorfield and  
$\xi_{\nu_1 \cdots \nu_n}^{\mu_1 \cdots \mu_m}$
is a type $\binom{m}{n}$ spacetime tensorfield,
then relative to the arbitrary coordinates,\footnote{It is well-known that
RHS~\eqref{E:LIEDERIVATIVE} is coordinate invariant.}
the Lie derivative of $\xi$ with respect to $V$ is
the type $\binom{m}{n}$ spacetime tensorfield $\Lie_V \xi$ with the following components:
\begin{align} \label{E:LIEDERIVATIVE}
\Lie_V \xi_{\nu_1 \cdots \nu_n}^{\mu_1 \cdots \mu_m}
& := V^{\alpha} \partial_{\alpha} \xi_{\nu_1 \cdots \nu_n}^{\mu_1 \cdots \mu_m}
	- \sum_{a=1}^m \xi_{\nu_1 \cdots \nu_n}^{\mu_1 \cdots \mu_{a-1} \alpha \mu_{a+1} \cdots \mu_m} \partial_{\alpha} V^{\mu_a}
	+ \sum_{b=1}^n \xi_{\nu_1 \cdots \nu_{b-1} \alpha \nu_{b+1} \cdots \nu_n}^{\mu_1 \cdots \mu_m} \partial_{\nu_b} V^{\alpha}.
\end{align}

In addition, when
$V$ and $W$ are both vectorfields,
we often use the standard Lie bracket notation
$[V,W] := \Lie_V W$.
\end{definition}


It is a standard fact that Lie differentiation  
obeys the Leibniz rule as well as the Jacobi-type identity
\begin{align} \label{E:JACOBI}
	\Lie_V \Lie_W \xi
	-
	\Lie_W \Lie_V \xi
	& = \Lie_{[V,W]} \xi
		= \Lie_{\Lie_V W} \xi.
\end{align} 
Moreover, it is a standard fact based on the torsion-free property of $\D$ that
RHS~\eqref{E:LIEDERIVATIVE} is invariant upon
replacing all coordinate partial derivatives $\partial$ 
with covariant derivatives $\D$.

In our analysis, we will apply the Leibniz rule for Lie derivatives to 
contractions of tensor products of $\ell_{t,u}-$tensorfields.
Due in part to the special properties 
(such as \eqref{E:LIELANDLIERADOFELLTUTANGENTISELLTUTANGENT})
of the vectorfields that we use to differentiate,
the non-$\ell_{t,u}$ components of the differentiated factor in the products typically cancel.
This motivates the following definition.

\begin{definition}[$\ell_{t,u}$ \textbf{and} $\Sigma_t-$\textbf{projected Lie derivatives}]
\label{D:PROJECTEDLIE}
Given a tensorfield $\xi$
and a vectorfield $V$,
we define the $\Sigma_t-$projected Lie derivative
$\SigmatLie_V \xi$ of $\xi$
and the $\ell_{t,u}-$projected Lie derivative
$\angLie_V \xi$ of $\xi$ as follows:
\begin{align} 
	\SigmatLie_V \xi
	& := \Sigmatproject \Lie_V \xi,
		\qquad
	\angLie_V \xi
	:= \Lineproject \Lie_V \xi.
	\label{E:PROJECTIONS}
\end{align}
\end{definition}

\begin{definition}[\textbf{Geometric torus differential}]
	\label{D:ANGULARDIFFERENTIAL}
	If $f$ is a scalar function on $\ell_{t,u}$, 
	then $\angdiff f := \angD f = \Lineproject \D f$,
	where $\D f$ is the gradient one-form associated to $f$.
\end{definition}

The above definition avoids potentially confusing notation
such as $\angD \Lunit^i$ by replacing it with 
$\angdiff \Lunit^i$; the latter notation clarifies that
$\Lunit^i$ is to be viewed as a scalar rectangular component function.

\begin{lemma}[\textbf{Sometimes} $\ell_{t,u}$ \textbf{projection is redundant}]
\label{L:SOMETIMESPROJECTIONISNOTNECESSARY}
	Let $\xi$ be a type $\binom{0}{n}$ spacetime tensorfield. 
	Then $\SigmatLie_{\Timenormal} \xi = \SigmatLie_{\Timenormal} (\Sigmatproject \xi)$
	and $\angLie_{\Lunit} \xi = \angLie_{\Lunit} \angxi$.
\end{lemma}
\begin{proof}
	To prove $\SigmatLie_{\Timenormal} \xi = \SigmatLie_{\Timenormal} (\Sigmatproject \xi)$, 
	we will show that 
	$\Sigmatproject_{\mu}^{\ \alpha} \Lie_{\Timenormal} \Sigmatproject_{\alpha}^{\ \nu} = 0$.
	Once we have shown this, we combine this identity with the Leibniz rule to deduce the following identity,
	where the first term on the RHS is exact and the second one schematic:
	$\SigmatLie_{\Timenormal} (\Sigmatproject \xi) 
	= 
	\SigmatLie_{\Timenormal} \xi
	+ \Sigmatproject \cdot \Lie_{\Timenormal} \Sigmatproject \cdot \xi.
	$
	A careful analysis of the schematic term shows that it always contains 
	a factor of the form
	$\Sigmatproject_{\mu}^{\ \alpha} \SigmatLie_{\Timenormal} \Sigmatproject_{\alpha}^{\ \nu}$, 
	which vanishes. We have proved the desired result.

	We now show that $\Sigmatproject_{\mu}^{\ \alpha} \Lie_{\Timenormal} \Sigmatproject_{\alpha}^{\ \nu} = 0$.
	Actually, we prove a stronger result:
	$\Lie_{\Timenormal} \Sigmatproject_{\nu}^{\ \mu} = 0$.
	Since $\Lie_{\Timenormal} \delta_{\nu}^{\ \mu} = 0$
	and since $\Lie_{\Timenormal} \Timenormal^{\mu} = 0$,
	we see from \eqref{E:SIGMATPROJECTION} that it suffices to prove
	that $\Lie_{\Timenormal} \Timenormal_{\nu}  = 0$.
	The LHS of the previous identity is 
	equal to the one-form $(\Lie_{\Timenormal} g_{\nu \alpha}) \Timenormal^{\alpha}$.
	To show that it vanishes, we separately show that its 
	$\Timenormal$ and $\Sigma_t-$tangent components vanish.
	For the former, we use the identity $g(\Timenormal, \Timenormal) = -1$
	and the Leibniz rule for Lie derivatives
	to deduce the desired result $(\Lie_{\Timenormal} g_{\nu \alpha}) \Timenormal^{\alpha} \Timenormal^{\nu} = 0$.
	It remains only for us to show that
	$(\Lie_{\Timenormal} g)(V,\Timenormal) = 0$ for $\Sigma_t-$tangent vectorfields $V$.
	Using that $V t = 0$, we compute that
	$(\Lie_{\Timenormal} V) t
	 = \Timenormal(Vt) - V (\Timenormal t) = -V(1)= 0
	$.
	It follows that $\Lie_{\Timenormal} V$ is also $\Sigma_t-$tangent
	and hence $g(\Lie_{\Timenormal} V, \Timenormal) = 0$. By the Leibniz rule for Lie derivatives,
	we conclude that $0 = \Timenormal(g(V,\Timenormal)) = (\Lie_{\Timenormal} g)(V,\Timenormal)$ as desired.

	The proof that $\angLie_{\Lunit} \xi = \angLie_{\Lunit} \angxi$ is similar
	and reduces to showing that 
	$\Lineproject_{\nu}^{\ \alpha} \Lie_{\Lunit} \Lineproject_{\alpha}^{\ \mu} = 0$.
	From \eqref{E:LINEPROJECTION}, we see that it further reduces to showing that
	$\Lineproject_{\nu}^{\ \alpha} 
	\Lie_{\Lunit} 
	\left\lbrace
		\Radunit_{\alpha} \Lunit^{\mu} 
		+ \Lunit_{\alpha} (\Lunit^{\mu} + \Radunit^{\mu})
	\right\rbrace
	= 0
	$.
	Since $\Lineproject$ annihilates $\Lunit$ and $\Radunit$,
	we need only to confirm that 
	$\Lineproject_{\nu}^{\ \alpha} 
	\Lie_{\Lunit}(\Radunit_{\alpha} + \Lunit_{\alpha}) = 0
	$, which is equivalent to
	$\Lineproject_{\nu}^{\ \alpha} 
	\Lie_{\Lunit} \Timenormal_{\alpha} = 0
	$. Since $\CoordAng$ spans the tangent space of $\ell_{t,u}$,
	it suffices to show that $\CoordAng^{\alpha} \cdot \Lie_{\Lunit} \Timenormal_{\alpha} = 0$.
	This latter identity follows easily from differentiating 
	the identity $\CoordAng^{\alpha} \Timenormal_{\alpha} = 0$
	with $\Lie_{\Lunit}$ and using the identity $\Lie_{\Lunit} \CoordAng = 0$
	(since $\Lunit = \frac{\partial}{\partial t}$ and $\CoordAng = \frac{\partial}{\partial \vartheta}$).
\end{proof}

\subsection{First and second fundamental forms and covariant differential operators}

\begin{definition}[\textbf{First fundamental forms}] \label{D:FIRSTFUND}
	We define the first fundamental form $\gt$ of $\Sigma_t$ and the 
	first fundamental form $\gsphere$ of $\ell_{t,u}$ as follows:
	\begin{align}
		\gt
		:= \Sigmatproject g,
		\qquad
		\gsphere
		:= \Lineproject g.
		\label{E:GTANDGSPHERESPHEREDEF}
	\end{align}

	We define the corresponding inverse first fundamental forms by raising the 
	indices with $g^{-1}$:
	\begin{align}
		(\gt^{-1})^{\mu \nu}
		:= (g^{-1})^{\mu \alpha} (g^{-1})^{\nu \beta} \gt_{\alpha \beta},
		\qquad 
		(\gsphere^{-1})^{\mu \nu}
		:= (g^{-1})^{\mu \alpha} (g^{-1})^{\nu \beta} \gsphere_{\alpha \beta}.
		\label{E:GGTINVERSEANDGSPHEREINVERSEDEF}
	\end{align}
\end{definition}
Note that $\gt$ is the Riemannian metric on $\Sigma_t$ induced by $g$
and that $\gsphere$ is the Riemannian metric on $\ell_{t,u}$ induced by $g$.
Moreover, a straightforward calculation shows that
$(\gt^{-1})^{\mu \alpha} \gt_{\alpha \nu} = \Sigmatproject_{\nu}^{\ \mu}$
and $(\gsphere^{-1})^{\mu \alpha} \gsphere_{\alpha \nu} = \Lineproject_{\nu}^{\ \mu}$.

\begin{remark}
	\label{E:ELLTUTENSORSAREPURETRACE}
	Because the $\ell_{t,u}$ are one-dimensional manifolds, it follows that
	symmetric type $\binom{0}{2}$ $\ell_{t,u}-$tangent tensorfields $\xi$
	satisfy $\xi = (\mytr \xi) \gsphere$,
	where $\mytr \xi := \ginversesphere \cdot \xi$.
	This simple fact simplifies some of our formulas
	compared to the case of higher spatial dimensions.
	In the remainder of the article, we often use this fact without
	explicitly mentioning it. 
	Moreover, as we described in Remark~\ref{R:ELLIPTIC},
	this fact is the reason that we do not need to derive
	elliptic estimates in two spatial dimensions.
\end{remark}

\begin{definition}[\textbf{Differential operators associated to the metrics}] 
\label{D:CONNECTIONS}
	We use the following notation for various differential operators associated 
	to the spacetime metric $g$, the Minkowski metric $m$, and the Riemannian metric
	$\gsphere$ induced on the $\ell_{t,u}$.
	\begin{itemize}
		\item $\D$ denotes the Levi-Civita connection of the spacetime metric $g$.
		\item $\angD$ denotes the Levi-Civita connection of $\gsphere$.
		\item If $\xi$ is an $\ell_{t,u}-$tangent one-form,
			then $\angdiv \xi$ is the scalar-valued function
			$\angdiv \xi := \ginversesphere \cdot \angD \xi$.
		\item Similarly, if $V$ is an $\ell_{t,u}-$tangent vectorfield,
			then $\angdiv V := \ginversesphere \cdot \angD V_{\flat}$,
			where $V_{\flat}$ is the one-form $\gsphere-$dual to $V$.
		\item If $\xi$ is a symmetric type $\binom{0}{2}$ 
		 $\ell_{t,u}-$tangent tensorfield, then
		 $\angdiv \xi$ is the $\ell_{t,u}-$tangent 
		 one-form $\angdiv \xi := \ginversesphere \cdot \angD \xi$,
		 where the two contraction indices in $\angD \xi$
		 correspond to the operator $\angD$ and the first index of $\xi$.
		\end{itemize}
\end{definition}

\begin{definition}[\textbf{Covariant wave operators and Laplacians}] 
	\label{D:WAVEOPERATORSANDLAPLACIANS}
	We use the following standard notation.
	\begin{itemize}
		\item $\square_g := (g^{-1})^{\alpha \beta} \D_{\alpha \beta}^2$ denotes the covariant wave operator 
			corresponding to the spacetime metric $g$.
		\item $\angLap := \ginversesphere \cdot \angD^2$ denotes the covariant Laplacian 
			corresponding to $\gsphere$. 
	\end{itemize}
\end{definition}

\begin{definition}[\textbf{Second fundamental forms}]
We define the second fundamental form $k$ of $\Sigma_t$, 
which is a symmetric type $\binom{0}{2}$ $\Sigma_t-$tangent tensorfield,
by
\begin{align} \label{E:SECONDFUNDSIGMATDEF}
	k 
	&:= \frac{1}{2} \SigmatLie_{\Timenormal} \gt.
\end{align}

	We define the null second fundamental form $\upchi$ of $\ell_{t,u}$, 
	which is a symmetric type $\binom{0}{2}$ $\ell_{t,u}-$tangent tensorfield,
	by
\begin{align} \label{E:CHIDEF}
	\upchi
	& := \frac{1}{2} \angLie_{\Lunit} \gsphere.
\end{align}

\end{definition}

From Lemma~\ref{L:SOMETIMESPROJECTIONISNOTNECESSARY},
we see that the following alternate expressions hold:
\begin{align} \label{E:ALTERNATESECONDFUND}
	k 
	&= \frac{1}{2} \SigmatLie_{\Timenormal} g,
		\qquad
	\upchi
	= \frac{1}{2} \angLie_{\Lunit} g.
\end{align}

We now provide some identities that we use later.

\begin{lemma}[\textbf{Alternate expressions for the second fundamental forms}]
We have the following identities:
\begin{subequations}
\begin{align} 
	\upchi_{\CoordAng \CoordAng}
	& = g(\D_{\CoordAng} \Lunit, \CoordAng),
		\label{E:CHIUSEFULID} \\
	\angkdoublearg{\Radunit}{\CoordAng}
	& = g(\D_{\CoordAng} \Lunit, \Radunit).
	\label{E:ANGKRINTERMSOFCOVARIANTDERIVATIVES}
\end{align}
\end{subequations}
\end{lemma}
\begin{proof}
	We prove only \eqref{E:ANGKRINTERMSOFCOVARIANTDERIVATIVES}
	since the proof of \eqref{E:CHIUSEFULID} is similar.
	Using \eqref{E:ALTERNATESECONDFUND}, we
	compute that 
	$2 \angkdoublearg{\Radunit}{\CoordAng} 
	= (\SigmatLie_{\Timenormal} g)_{\Radunit \CoordAng}
	= (\Lie_{\Timenormal} g)_{\Radunit \CoordAng}
	= g(\D_{\Radunit} \Timenormal, \CoordAng)
		+ g(\D_{\CoordAng} \Timenormal, \Radunit)
	$.
	Since $g(\Radunit,\Radunit) = 1$ and $\Timenormal = \Lunit + \Radunit$,
	we see that 
	$g(\D_{\CoordAng} \Timenormal, \Radunit) 
	= g(\D_{\CoordAng} \Lunit, \Radunit) 
	$. Thus, to complete the proof, we need only to show that
	$g(\D_{\Radunit} \Timenormal, \CoordAng) = g(\D_{\CoordAng} \Lunit, \Radunit)$.
	To proceed, we note that since $g(\Timenormal,\Radunit) = 0$ and $g(\Radunit,\Radunit) = 1$, we have
	$g(\D_{\CoordAng} \Timenormal, \Radunit) = - g(\D_{\CoordAng} \Radunit, \Timenormal)
	= - g(\D_{\CoordAng} \Radunit, \Lunit)
	$. Then since $g(\Radunit,\Lunit) = -1$, we conclude that 
	$- g(\D_{\CoordAng} \Radunit, \Lunit) = g(\D_{\CoordAng} \Lunit, \Radunit)$ as desired.
\end{proof}

\subsection{Expressions for the metrics}
\label{SS:METRICEXPRESSIONS}
In this subsection, we decompose $g$ relative to
the non-rescaled frame and relative to the geometric coordinates.
We then provide expressions for various forms
relative to the geometric coordinates and
for the change of variables map from geometric to rectangular coordinates.

\begin{lemma}[\textbf{Expressions for $g$ and $g^{-1}$ in terms of the non-rescaled frame}]
\label{L:METRICDECOMPOSEDRELATIVETOTHEUNITFRAME}
We have the following identities:
\begin{subequations}
\begin{align}
	g_{\mu \nu} 
	& = - \Lunit_{\mu} \Lunit_{\nu}
			- (
					\Lunit_{\mu} \Radunit_{\nu} 
					+ 
					\Radunit_{\mu} \Lunit_{\nu}
				)
			+ \gsphere_{\mu \nu} 
			\label{E:METRICFRAMEDECOMPLUNITRADUNITFRAME},
			\\
	(g^{-1})^{\mu \nu} 
	& = 
			- \Lunit^{\mu} \Lunit^{\nu}
			- (
					\Lunit^{\mu} \Radunit^{\nu} 
					+ \Radunit^{\mu} \Lunit^{\nu}
				)
			+ (\ginversesphere)^{\mu \nu}.
			\label{E:GINVERSEFRAMEWITHRECTCOORDINATESFORGSPHEREINVERSE}
\end{align}
\end{subequations}

\end{lemma}

\begin{proof}
	It suffices to prove \eqref{E:METRICFRAMEDECOMPLUNITRADUNITFRAME} 
	since \eqref{E:GINVERSEFRAMEWITHRECTCOORDINATESFORGSPHEREINVERSE}
	then follows from raising the indices
	of \eqref{E:METRICFRAMEDECOMPLUNITRADUNITFRAME} 
	with $g^{-1}$.

	To verify the formula \eqref{E:METRICFRAMEDECOMPLUNITRADUNITFRAME}, 
	we contract each side against the rectangular coordinates of pairs of elements of
	the frame $\lbrace \Lunit, \Radunit, \CoordAng \rbrace$ and check that both sides agree.
	This of course requires that we know the inner products of all pairs of elements of the frame,
	some of which follow from the basic properties of the frame vectorfields,
	and some of which were established in Lemma~\ref{L:BASICPROPERTIESOFFRAME}.
	As an example, we note that contracting the LHS against $\Lunit^{\mu} \CoordAng^{\nu}$
	yields $g(\Lunit,\CoordAng) = 0$, while contracting the RHS yields
	$- g(\Lunit, \Lunit) g(\Lunit,\CoordAng)
	- g(\Lunit, \Lunit) g(\Radunit,\CoordAng)
	- g(\Radunit, \Lunit) g(\Lunit,\CoordAng)
	+ \gsphere(\Lunit,\CoordAng) = 0 + 0 + 0 + 0 = 0$
	as desired.
	As a second example, we note that
	contracting the LHS against $\Lunit^{\mu} \Radunit^{\nu}$
	yields $g(\Lunit,\Radunit) = -1$, 
	while contracting the RHS yields
	$- g(\Lunit, \Lunit) g(\Lunit,\Radunit)
	- g(\Lunit, \Lunit) g(\Radunit,\Radunit)
	- g(\Radunit, \Lunit) g(\Lunit,\Radunit)
	+ \gsphere(\Lunit,\Radunit) = 0 + 0 - 1+ 0 = -1$
	as desired.
\end{proof}

The following scalar function captures the $\ell_{t,u}$ part of $g$.

\begin{definition}[\textbf{The metric component} $\gtancomp$]
\label{D:METRICANGULARCOMPONENT}
We define the scalar function $\gtancomp > 0$ by
\begin{align} \label{E:METRICANGULARCOMPONENT}
	\gtancomp^2
	& :=  g(\CoordAng,\CoordAng) = \gsphere(\CoordAng,\CoordAng).
\end{align}
\end{definition}

It follows that relative to the geometric coordinates, we have
\begin{align} \label{E:GSPHEREINVERSERELATIVETOGEOMETRIC}
	\ginversesphere
	= \gtancomp^{-2} \CoordAng \otimes \CoordAng.
\end{align}

We now express $g$ relative to the geometric coordinates.

\begin{lemma}[\textbf{Expressions for $g$ and $g^{-1}$ in terms of the geometric coordinate frame}]
Relative to the geometric coordinate $(t,u,\vartheta)$, we have
\begin{align} \label{E:METRICRELATIVETOGEOMETRICCOORDINATES}
	g 
	& = - 2 \upmu dt \otimes du
		+ \upmu^2 du \otimes du
		+ \gtancomp^2(d \vartheta + \XiCoordComp du) \otimes (d \vartheta + \XiCoordComp du),
			\\
	g^{-1}
	& = - \frac{\partial}{\partial t} \otimes \frac{\partial}{\partial t}
		- \upmu^{-1} \frac{\partial}{\partial t} \otimes \frac{\partial}{\partial u}
		- \upmu^{-1} \frac{\partial}{\partial u} \otimes \frac{\partial}{\partial t}
		- \upmu^{-1} \XiCoordComp \frac{\partial}{\partial t} \otimes \CoordAng
		- \upmu^{-1} \XiCoordComp \CoordAng \otimes \frac{\partial}{\partial t}
		+ \gtancomp^{-2} \CoordAng \otimes \CoordAng.
			\label{E:INVERSEMETRICRELATIVETOGEOMETRICCOORDINATES}
\end{align}
The scalar functions
$\XiCoordComp$
and $\gtancomp$
from above
are defined respectively in
\eqref{E:RADSPLITINTOPARTTILAUANDXI}
and \eqref{E:METRICANGULARCOMPONENT}.
\end{lemma}
\begin{proof}
We recall that by Lemma~\ref{L:BASICPROPERTIESOFFRAME} and \eqref{E:GSPHEREINVERSERELATIVETOGEOMETRIC}, 
we have
$\Lunit = \frac{\partial}{\partial t}$,
and $\upmu \Radunit = \Rad = \frac{\partial}{\partial u} - \XiCoordComp \CoordAng$,
and $\ginversesphere = \gtancomp^{-2} \frac{\partial}{\partial \vartheta} \otimes \frac{\partial}{\partial \vartheta}$.
Moreover, by \eqref{E:GINVERSEFRAMEWITHRECTCOORDINATESFORGSPHEREINVERSE},
we have $g^{-1} = - \Lunit \otimes \Lunit - \Lunit \otimes \Radunit - \Radunit \otimes \Lunit + \ginversesphere$.
Combining these identities, we easily conclude \eqref{E:INVERSEMETRICRELATIVETOGEOMETRICCOORDINATES}.
\eqref{E:METRICRELATIVETOGEOMETRICCOORDINATES} then follows from \eqref{E:INVERSEMETRICRELATIVETOGEOMETRICCOORDINATES}
as a simple linear algebra exercise
(just compute the inverse of the $3 \times 3$ matrix corresponding to 
\eqref{E:INVERSEMETRICRELATIVETOGEOMETRICCOORDINATES}).
\end{proof}

We now provide expressions for the geometric volume form factors of $g$ and $\gt$.
\begin{corollary}[\textbf{The geometric volume form factors of} $g$ \textbf{and} $\gt$]
\label{C:SPACETIMEVOLUMEFORMWITHUPMU}
The following identity is verified by the spacetime metric $g$:
\begin{align} \label{E:SPACETIMEVOLUMEFORMWITHUPMU}
	|\mbox{\upshape{det}} g| 
	& = \upmu^2 \gtancomp^2,
\end{align}
where the determinant on the LHS is taken
relative to the geometric coordinates
$(t,u,\vartheta)$.

Furthermore, the following identity is verified by 
the first fundamental form $\gt$ of $\Sigma_t^{U_0}$: 
\begin{align} \label{E:SIGMATVOLUMEFORMWITHUPMU}
	\mbox{\upshape{det}} \gt|_{\Sigma_t^{U_0}}
	& = \upmu^2 \gtancomp^2,
\end{align}
where the determinant on the LHS is taken
relative to the geometric coordinates
$(u,\vartheta)$ induced on $\Sigma_t^{U_0}$.
\end{corollary}

\begin{proof}
	Equation \eqref{E:SPACETIMEVOLUMEFORMWITHUPMU}
	follows easily from computing the determinant of \eqref{E:METRICRELATIVETOGEOMETRICCOORDINATES}.

	Next, we note that 
	\eqref{E:METRICRELATIVETOGEOMETRICCOORDINATES} implies that
	$\gt = \upmu^2 du^2
	+ \gtancomp^2(d \vartheta + \XiCoordComp du) (d \vartheta + \XiCoordComp du)$.
	A simple calculation then yields \eqref{E:SIGMATVOLUMEFORMWITHUPMU}.
\end{proof}

\begin{definition}\label{D:CHOVMAP}
	We define $\Upsilon: [0,T) \times [0,U_0] \times \mathbb{T} \rightarrow \mathcal{M}_{T,U_0}$,
	$\Upsilon(t,u,\vartheta) := (t,x^1,x^2)$,
	to be the change of variables map from geometric to rectangular coordinates.
\end{definition}

\begin{lemma}[\textbf{Basic properties of the change of variables map}]
\label{L:CHOV}
	We have the following expression for the Jacobian of $\Upsilon$:
	\begin{align} \label{E:CHOV}
		\frac{\partial \Upsilon}{\partial (t,u,\vartheta)}
		& :=
		\frac{\partial (x^0,x^1,x^2)}{\partial (t,u,\vartheta)}
		= 
		\left(
		\begin{array}{ccc}
			1 & 0 & 0  \\
			\Lunit^1 & \Rad^1 + \NonRadialRad^1 & \CoordAng^1 \\
			\Lunit^2 & \Rad^2 + \NonRadialRad^2 & \CoordAng^2 \\
		\end{array}
		\right).
	\end{align}
	Moreover, the Jacobian determinant of $\Upsilon$ can be expressed as
\begin{align} \label{E:JACOBIAN}
	\mbox{\upshape{det}}
	\frac{\partial (x^0,x^1,x^2)}{\partial (t,u,\vartheta)}
	= \upmu (\mbox{\upshape{det}} \gt_{ij})^{-1/2} \gtancomp,
\end{align}
where $\gtancomp$ is the metric component from Def.~\ref{D:METRICANGULARCOMPONENT}
and $(\mbox{\upshape{det}} \gt_{ij})^{-1/2}$ is a smooth function of $\Psi$ in a neighborhood of $0$ 
with $(\mbox{\upshape{det}} \gt_{ij})^{-1/2}(\Psi=0) = 1$. In \eqref{E:JACOBIAN},
$\gt$ is viewed as the Riemannian metric on $\Sigma_t^{U_0}$ defined by \eqref{E:GTANDGSPHERESPHEREDEF}
and $\mbox{\upshape{det}} \gt_{ij}$ is the determinant of the corresponding $2 \times 2$ matrix
of components of $\gt$ relative to the rectangular spatial coordinates.
\end{lemma}
\begin{proof}
	Since $\frac{\partial}{\partial t} = \Lunit$, 
	first column of the matrix on RHS~\eqref{E:CHOV}
	is by definition $(\Lunit x^0, \Lunit x^1, \Lunit x^2)^{\top} = (1, \Lunit^1, \Lunit^2)^{\top}$,
	where $\top$ denotes the transpose operator
	and we have used \eqref{E:LUNITOFUANDT}.
	The second column is $(\frac{\partial}{\partial u} x^0, \frac{\partial}{\partial u} x^1, \frac{\partial}{\partial u} x^2)^{\top}$,
	and to obtain the form stated on RHS~\eqref{E:CHOV}, we use \eqref{E:RADSPLITINTOPARTTILAUANDXI}.
	The third column is 
	$(\CoordAng x^0, \CoordAng x^1, \CoordAng x^2)^{\top}$, and to obtain the stated form, we
	use the fact that $\CoordAng x^0 = \CoordAng t = 0$ (since $\CoordAng$ is $\Sigma_t-$tangent).

	To obtain \eqref{E:JACOBIAN}, we first observe that the determinant of the RHS is equal
	to the determinant of the $2 \times 2$ lower right block. Moreover, since
	$\NonRadialRad$ and $\CoordAng$ are parallel, we can assume that $\NonRadialRad \equiv 0$.
	Also recalling that $\Rad = \upmu \Radunit$, 
	we see that the determinant of interest
	is equal to
	$
	\upmu
	\mbox{\upshape{det}}
	N
	$,
	where
	$N :=
	\left(
	\begin{array}{cc}
			\Radunit^1 & \CoordAng^1 \\
			\Radunit^2 & \CoordAng^2. \\
		\end{array}
	\right)
	$.
	Next we consider the $2 \times 2$ matrix
	$
	M:=
	\left(
	\begin{array}{cc}
			g(\Radunit,\Radunit) &  g(\Radunit,\CoordAng)\\
			g(\CoordAng,\Radunit) & g(\CoordAng,\CoordAng) \\
		\end{array}
	\right)
	=
	\left(
	\begin{array}{cc}
			1 &  0 \\
			0 & \gtancomp^2 \\
		\end{array}
	\right)
	$.
	On the one hand, we clearly have $\mbox{\upshape{det}} M = \gtancomp^2$.
	On the other hand, we have the matrix identity
	$M = N^{\top} \cdot \gt \cdot N$ 
	(where $\gt$ is viewed as a $2 \times 2$ matrix expressed relative to the spatial rectangular coordinates),
	which implies that $\mbox{\upshape{det}} M  =  \mbox{\upshape{det}} \gt (\mbox{\upshape{det}} N)^2$.
	Combining these identities, we conclude \eqref{E:JACOBIAN}.
	Finally, we note that since $\gt_{ij} = \delta_{ij} + \smoothfunction(\Psi) \Psi$ with $\smoothfunction$ 
	a tensor depending smoothly on $\Psi$,
	we easily conclude that $(\mbox{\upshape{det}} \gt_{ij})^{-1/2}$ is a smooth function of $\Psi$ in a neighborhood of $0$ with 
	$(\mbox{\upshape{det}} \gt_{ij})^{-1/2}(\Psi=0) = 1$.
\end{proof}

\subsection{Commutation vectorfields}
\label{SS:COMMUTATIONVECTORFIELDS}
To obtain higher-order estimates for $\Psi$ and the eikonal function quantities along
$\ell_{t,u}$, we commute various evolution equations
with an $\ell_{t,u}-$tangent vectorfield.
A natural candidate commutator is the geometric coordinate partial derivative vectorfield
$\CoordAng$, which solves the transport equation
$\Lie_{\Lunit} \CoordAng = 0$.
In terms of the rectangular component functions, 
the transport equation reads
$\Lunit \CoordAng^i = \CoordAng \cdot \angdiff \Lunit^i$
and thus $\CoordAng^i$ is one degree less differentiable than 
$\Lunit^i$ in directions transversal to $\Lunit$. 
This loss of a derivative introduces 
technical complications into the analysis
that have no obvious resolution.
To circumvent this difficulty, we instead commute with the $\ell_{t,u}-$tangent vectorfield
$\GeoAng$, obtained by projecting a rectangular coordinate vectorfield 
$\GeoAng_{(Flat)}$ onto the $\ell_{t,u}$.
The identity \eqref{E:DOWNSTAIRSUPSTAIRSSRADUNITPLUSLUNITISAFUNCTIONOFPSI} and
Lemma~\ref{L:GEOANGDECOMPOSITION} below together show that
unlike $\CoordAng$,
the rectangular components $\GeoAng^i$
have the same degree of differentiability
as $\Psi$ and $\Lunit^i$. Another advantage of
using the commutator $\GeoAng$
is that its deformation tensor structure allows us
to derive our high-order
energy estimates without commuting
the wave equation with the transversal
vectorfield $\Rad$ at high orders
(see Def.~\ref{D:MAINCOERCIVEQUANT}
and Prop.~\ref{P:MAINAPRIORIENERGY}).
We note here that at first glance, the
top-order derivatives of 
the deformation tensor of $\GeoAng$ that appear
in the top-order wave equation energy estimates
seem to lose derivatives relative to $\Psi$.
However, we are able to overcome this difficulty
by working with modified 
quantities, which we construct in
Sect.~\ref{S:MODQUANTS}.

\begin{definition}[\textbf{The vectorfields} $\GeoAng_{(Flat)}$ \textbf{and} $\GeoAng$]
\label{D:ANGULARVECTORFIELDS}
We define the rectangular components of
the $\Sigma_t-$tangent vectorfields 
$\GeoAng_{(Flat)}$ and $\GeoAng$ 
as follows ($i=1,2$):
\begin{align}
	\GeoAng_{(Flat)}^i
	&: = \delta_2^i,
		\label{E:GEOANGEUCLIDEAN} \\
	\GeoAng^i
	& :=
	\Lineproject_a^{\ i} \GeoAng_{(Flat)}^a
		= \Lineproject_2^{\ i},
		\label{E:GEOANGDEF}
\end{align}
where $\Lineproject$ is the
$\ell_{t,u}$ projection tensorfield 
defined in \eqref{E:LINEPROJECTION}.
\end{definition}

To prove our main theorem, we commute the equations
with the elements of the following set of vectorfields.
\begin{definition}[\textbf{Commutation vectorfields}]
	\label{D:COMMUTATIONVECTORFIELDS}
	We define the commutation set $\Fullset$ as follows:
	\begin{align} \label{E:COMMUTATIONVECTORFIELDS}
		\Fullset
		:= \lbrace \Lunit, \Rad, \GeoAng \rbrace,
	\end{align}
	where $\Lunit$, $\Rad$, and $\GeoAng$ are respectively defined by
	\eqref{E:LUNITDEF}, \eqref{E:RADDEF}, and \eqref{E:GEOANGDEF}.

	We define the $\mathcal{P}_u-$tangent commutation set $\Tanset$ as follows:
	\begin{align} \label{E:TANGENTIALCOMMUTATIONVECTORFIELDS}
		\Tanset
		:= \lbrace \Lunit, \GeoAng \rbrace.
	\end{align}
\end{definition}

The rectangular spatial components of
$\Lunit$,
$\Radunit$,
and $\GeoAng$
deviate from their
flat values by a small amount
captured in the following definition.

\begin{definition}[\textbf{Perturbed part of various vectorfields}]
\label{D:PERTURBEDPART}
For $i=1,2$, we define the following scalar functions:
\begin{align} \label{E:PERTURBEDPART}
	\Lunit_{(Small)}^i
	& := \Lunit^i
		- \delta_1^i,
	\qquad
	\Radunit_{(Small)}^i
	:= \Radunit^i
		+ \delta_1^i,
	\qquad
	\GeoAng_{(Small)}^i
	:= \GeoAng^i - \delta_2^i.
\end{align}
The vectorfields 
$\Lunit$,
$\Radunit$,
and
$\GeoAng$
in \eqref{E:PERTURBEDPART}
are defined in
Defs.~\ref{D:LUNITDEF},
\ref{D:RADANDXIDEFS},
and \ref{D:ANGULARVECTORFIELDS}.

\end{definition}

\begin{remark}
	\label{R:RADUNITSMALLLUNITSMALLRELATION}
	From 
	\eqref{E:LITTLEGDECOMPOSED},
	\eqref{E:METRICPERTURBATIONFUNCTION},
	\eqref{E:DOWNSTAIRSUPSTAIRSSRADUNITPLUSLUNITISAFUNCTIONOFPSI}, 
	and \eqref{E:PERTURBEDPART}, 
	we have that 
	$\Radunit_{(Small)}^i = -\Lunit_{(Small)}^i - (g^{-1})^{0i}$,
	where $(g^{-1})^{0i}(\Psi = 0) = 0$.
	We will use this simple fact later on.
\end{remark}

In the next lemma, we characterize the
discrepancy between $\GeoAng_{(Flat)}$ and $\GeoAng$.

\begin{lemma}[\textbf{Decomposition of} $\GeoAng_{(Flat)}$]
\label{L:GEOANGDECOMPOSITION}
We can decompose $\GeoAng_{(Flat)}$ into an $\ell_{t,u}-$tangent vectorfield
and a vectorfield parallel to $\Radunit$ as follows:
since $\GeoAng$ is $\ell_{t,u}-$tangent, 
there exists a scalar function $\GeoAngFlatRadComponent$ such that
\begin{subequations}
\begin{align} \label{E:GEOANGINTERMSOFEUCLIDEANANGANDRADUNIT}
	\GeoAng_{(Flat)}^i
	& = \GeoAng^i 
		+ \GeoAngFlatRadComponent \Radunit^i,
			\\
	\GeoAng_{(Small)}^i
	& = - \GeoAngFlatRadComponent \Radunit^i.
	\label{E:GEOANGSMALLINTERMSOFRADUNIT}
\end{align}
\end{subequations}
Moreover, we have
\begin{align} \label{E:FLATYDERIVATIVERADIALCOMPONENT}
	\GeoAngFlatRadComponent 
	= g(\GeoAng_{(Flat)},\Radunit)
	= g_{ab} \GeoAng_{(Flat)}^a \Radunit^b
	= g_{2a} \Radunit^a
	= g_{21}^{(Small)} \Radunit^1 - g_{22} \Radunit_{(Small)}^2.
\end{align}
\end{lemma}
\begin{proof}
	The existence of the decomposition 
	\eqref{E:GEOANGINTERMSOFEUCLIDEANANGANDRADUNIT} 
	follows from the fact that 
	by construction,
	$\GeoAng_{(Flat)}$ and $\GeoAng$
	differ only by a vectorfield that is parallel to $\Radunit$
	(because the $\ell_{t,u}$ projection tensorfield $\Lineproject$ annihilates the 
	$\Radunit$ component of the $\Sigma_t-$tangent vectorfield $\GeoAng_{(Flat)}$
	while preserving its $\ell_{t,u}-$tangent component).

	The expression \eqref{E:GEOANGSMALLINTERMSOFRADUNIT}
	then follows from definition
	\eqref{E:PERTURBEDPART}
	and
	\eqref{E:GEOANGINTERMSOFEUCLIDEANANGANDRADUNIT}.

	To obtain \eqref{E:FLATYDERIVATIVERADIALCOMPONENT}, we
	contract \eqref{E:GEOANGINTERMSOFEUCLIDEANANGANDRADUNIT}
	against $\Radunit_i$ and use 
	\eqref{E:LITTLEGDECOMPOSED}-\eqref{E:METRICPERTURBATIONFUNCTION},
	the identities
	$\GeoAng^a \Radunit_a = 0$
	and
	$\Radunit^a \Radunit_a = 1$,
	and definition \eqref{E:PERTURBEDPART}
\end{proof}

\subsection{Deformation tensors and basic vectorfield commutator properties}
\label{SS:BASICVECTORFIELDCOMMUTATOR}
In this subsection, we start by recalling the standard definition of the deformation tensor of a vectorfield.
We then exhibit some basic properties enjoyed by the Lie derivatives
of various vectorfields.

\begin{definition}[\textbf{Deformation tensor of a vectorfield} $V$]
\label{D:DEFORMATIONTENSOR}
If $V$ is a spacetime vectorfield,
then its deformation tensor $\deform{V}$
(relative to the spacetime metric $g$)
is the symmetric type $\binom{0}{2}$ spacetime tensorfield
\begin{align} \label{E:DEFORMATIONTENSOR}
	\deformarg{V}{\alpha}{\beta}
	:= \Lie_V g_{\alpha \beta}
	= \D_{\alpha} V_{\beta} 
		+
		\D_{\beta} V_{\alpha},
\end{align}
where the last equality in \eqref{E:DEFORMATIONTENSOR} is a well-known consequence of 
\eqref{E:LIEDERIVATIVE} and the torsion-free property of the connection $\D$.

\end{definition}

\begin{lemma}[\textbf{Basic vectorfield commutator properties}]
\label{L:CONNECTIONBETWEENCOMMUTATORSANDDEFORMATIONTENSORS}
The vectorfields 
$[\Lunit, \Rad]$,
$[\Lunit, \GeoAng]$,
and
$[\Rad, \GeoAng] $
are $\ell_{t,u}-$ tangent,
and the following identities hold:
\begin{align}
	[\Lunit, \Rad] 
	& = \angdeformoneformupsharparg{\Rad}{\Lunit},
		\qquad
	[\Lunit, \GeoAng] 
	= \angdeformoneformupsharparg{\GeoAng}{\Lunit},
		\qquad
	[\Rad, \GeoAng] 
	= \angdeformoneformupsharparg{\GeoAng}{\Rad}.
		\label{E:CONNECTIONBETWEENCOMMUTATORSANDDEFORMATIONTENSORS}
\end{align}

Furthermore, if $Z \in \Fullset$, then
\begin{align} 
		\angLie_Z \gsphere 
		& = \angdeform{Z},
			\qquad
		\angLie_Z \ginversesphere
		= - \angdeform{Z}^{\# \#}.
			\label{E:CONNECTIONBETWEENANGLIEOFGSPHEREANDDEFORMATIONTENSORS}
	\end{align}

Finally, if $V$ is an $\ell_{t,u}-$tangent vectorfield, then
\begin{align} \label{E:LIELANDLIERADOFELLTUTANGENTISELLTUTANGENT}
	[\Lunit, V] \mbox{ and } [\Rad, V] \mbox{ are } \ell_{t,u}-\mbox{tangent}.
\end{align}
\end{lemma}

\begin{proof}
	We first prove \eqref{E:LIELANDLIERADOFELLTUTANGENTISELLTUTANGENT}. 
	We use the identities
	$\Lunit u = 0$ and $\Lunit t = 1$ to compute that
	$[\Lunit, V] t = \Lunit V t - V \Lunit t = \Lunit 0 - V 1 = 0$
	and $[\Lunit, V] u = \Lunit V u - V \Lunit u = \Lunit 0 - V 0 = 0$.
	Since $[\Lunit, V]$ annihilates $t$ and $u$, it must be $\ell_{t,u}-$tangent as desired.
	A similar argument based on the identities $\Rad u = 1$ and $\Rad t = 0$
	yields that $[\Rad,V]$ is $\ell_{t,u}-$tangent.

	We now prove \eqref{E:CONNECTIONBETWEENCOMMUTATORSANDDEFORMATIONTENSORS}.
	Using the arguments from the previous paragraph,
	we easily deduce that the left-hand and right-hand sides of the identities
	are vectorfields that annihilate the function $t$
	and are therefore $\Sigma_t-$tangent. 
	Hence, it suffices to show that the inner products of the two sides of \eqref{E:CONNECTIONBETWEENCOMMUTATORSANDDEFORMATIONTENSORS}
	with $\Rad$ are equal and that the same holds for inner products with $\GeoAng$.
	We give the details only in the case of the last identity 
	$[\Rad, \GeoAng] = \angdeformoneformupsharparg{\GeoAng}{\Rad}$
	since the other two can be proved similarly.
	First, we note that the inner product of $\Rad$ and $\angdeformoneformupsharparg{\GeoAng}{\Rad}$
	is trivially $0$.
	Moreover, since we showed in the first paragraph that $[\Rad, \GeoAng]$ is $\ell_{t,u}-$tangent,
	we conclude that $g([\Rad, \GeoAng],\Rad) = 0$ as desired.
	We now show that
	the inner products of $\GeoAng$ and the two sides of 
	the last identity in \eqref{E:CONNECTIONBETWEENCOMMUTATORSANDDEFORMATIONTENSORS}
	are equal. Using again the torsion-free property and the fact that $g(\Rad,\GeoAng) = 0$,
	we compute that
	$
	g([\Rad, \GeoAng],\GeoAng) 
	= g(\D_{\Rad} \GeoAng,\GeoAng) - g(\D_{\GeoAng} \Rad, \GeoAng)
	= g(\D_{\Rad} \GeoAng,\GeoAng) + g(\D_{\GeoAng} \GeoAng,\Rad)
	$.
	The RHS of this identity is equal to the inner product of the RHS of the last identity in
	\eqref{E:CONNECTIONBETWEENCOMMUTATORSANDDEFORMATIONTENSORS} with
	$\GeoAng$ as desired.

	To prove \eqref{E:CONNECTIONBETWEENANGLIEOFGSPHEREANDDEFORMATIONTENSORS}
	for $\angLie_Z \gsphere$,
	we apply 
	$\angLie_Z \gsphere$
	to the identity
	\eqref{E:METRICFRAMEDECOMPLUNITRADUNITFRAME}.
	The LHS of the resulting identity is $\angdeform{Z}$,
	while only the last term
	$\angLie_Z \gsphere$
	survives on the RHS since
	the $\ell_{t,u}-$projection 
	$\Lineproject$ annihilates the
	non-differentiated factors arising from the first
	three tensor products on RHS~\eqref{E:METRICFRAMEDECOMPLUNITRADUNITFRAME}.
	We have thus proved 
	\eqref{E:CONNECTIONBETWEENANGLIEOFGSPHEREANDDEFORMATIONTENSORS}
	for $\angLie_Z \gsphere$.
	The identity \eqref{E:CONNECTIONBETWEENANGLIEOFGSPHEREANDDEFORMATIONTENSORS} 
	for $\angLie_Z \ginversesphere$ is a simple consequence of the identity for
	$\angLie_Z \gsphere$, the identity
	$
	(\ginversesphere)^{\alpha \kappa} \gsphere_{\kappa \beta}
	= \Lineproject_{\beta}^{\ \alpha}
	$,
	the Leibniz rule,
	and the identity
	$
	(\ginversesphere)^{\alpha \kappa} \angLie_Z \Lineproject_{\kappa}^{\ \beta}
	= 0
	$,
	which we now prove.
	Since $(\ginversesphere)^{\alpha \beta} = \gtancomp^{-2} \CoordAng^{\alpha} \CoordAng^{\beta}$,
	the proof reduces to showing that
	$(\Lie_Z \Lineproject_i^{\ a}) \CoordAng^i = 0$. To this end,
	we differentiate the identity $\CoordAng = \Lineproject \cdot \CoordAng$ and use the Leibniz rule
	to deduce that
	$\Lie_Z \CoordAng = (\Lie_Z \Lineproject) \cdot \CoordAng + \Lineproject \cdot \Lie_Z \CoordAng$.
	Using \eqref{E:LIELANDLIERADOFELLTUTANGENTISELLTUTANGENT}, we see that
	$\Lie_Z \CoordAng = \Lineproject \cdot \Lie_Z \CoordAng$,
	which finishes the proof.

\end{proof}

\begin{lemma}[$\Lunit$, $\Rad$, $\GeoAng$ \textbf{commute with} $\angdiff$]
\label{L:LANDRADCOMMUTEWITHANGDIFF}
For scalar functions $f$ and $V \in \lbrace \Lunit, \Rad, \GeoAng \rbrace$,
we have
\begin{align} \label{E:ANGLIECOMMUTESWITHANGDIFF}
	\angLie_V \angdiff f
	& = \angdiff V f.
\end{align}
\end{lemma}

\begin{proof}
	We prove the identity only when $V = \Rad$ since the remaining identities
	can be proved similarly.
	To proceed, we contract \eqref{E:ANGLIECOMMUTESWITHANGDIFF} against $\GeoAng$
	and use the Leibniz rule on the LHS to
	find that the identity is equivalent to
	$\Rad \GeoAng f
	- (\Lie_{\Rad} \GeoAng) \cdot \angdiff f
	= \GeoAng \Rad f
	$.
	The previous identity is equivalent to
	$[\Rad,\GeoAng] = \angLie_{\Rad} \GeoAng$,
	which follows from \eqref{E:CONNECTIONBETWEENCOMMUTATORSANDDEFORMATIONTENSORS}.
\end{proof}

\subsection{The rectangular Christoffel symbols}
\label{SS:RECTANGULARCHRISTOFFELSYMBOLS}
In many of our subsequent calculations, we start by expressing quantities
in rectangular coordinates. The most important of these are the
Christoffel symbols.

\begin{lemma}[\textbf{Christoffel symbols of} $g$ \textbf{in rectangular coordinates}]
	\label{L:CHRISTOFEELRECT}
	Let 
	\[
	\Gamma_{\alpha \kappa \beta} 
		:=
		\frac{1}{2} 
		\left\lbrace
			\partial_{\alpha} g_{\kappa \beta} 
			+ 
			\partial_{\beta} g_{\alpha \kappa}
			-
			\partial_{\kappa} g_{\alpha \beta}
		\right\rbrace
	\]
	denote the lowered Christoffel symbols of $g$ relative to rectangular coordinates
	and recall that $G_{\alpha \beta}(\Psi) = \frac{d}{d \Psi} g_{\alpha \beta}(\Psi)$.
	Then we have
	\begin{align} \label{E:CHRISTOFEELRECT}
		\Gamma_{\alpha \kappa \beta}
		& = 
			\frac{1}{2}
			\left\lbrace
				 G_{\kappa \beta} \partial_{\alpha} \Psi
				+ 
				 G_{\alpha \kappa} \partial_{\beta} \Psi
				-
				 G_{\alpha \beta} \partial_{\kappa} \Psi
		\right\rbrace.
	\end{align}
\end{lemma}

\begin{proof}
	\eqref{E:CHRISTOFEELRECT} is a simple consequence of the chain rule.
\end{proof}

\subsection{Transport equations for the eikonal function quantities}
\label{SS:TRANSPORTEQUATIONSFOREIKONALFUNCTION}
We now use Lemma~\ref{L:CHRISTOFEELRECT} to derive evolution equations for $\upmu$ and the rectangular components 
$\Lunit_{(Small)}^i$, ($i=1,2$).

\begin{lemma} [\textbf{The transport equations verified by} $\upmu$ \textbf{and} $\Lunit^i$] 
\label{L:UPMUANDLUNITIFIRSTTRANSPORT}
The inverse foliation density $\upmu$ defined in \eqref{E:UPMUDEF} verifies the following transport equation:
\begin{align} \label{E:UPMUFIRSTTRANSPORT}
	\Lunit \upmu 
	& := \upomega
		=
		\frac{1}{2} G_{\Lunit \Lunit} \Rad \Psi
		- \frac{1}{2} \upmu G_{\Lunit \Lunit} \Lunit \Psi
		- \upmu G_{\Lunit \Radunit} \Lunit \Psi.
\end{align}

Moreover, the scalar-valued rectangular component functions $\Lunit_{(Small)}^i$,
($i=1,2$), defined in \eqref{E:PERTURBEDPART},
verify the following transport equation:
\begin{align}
	\Lunit \Lunit_{(Small)}^i
	&  = - \frac{1}{2} G_{\Lunit \Lunit} (\Lunit \Psi) \Lunit^i
			- \frac{1}{2} G_{\Lunit \Lunit} (\Lunit \Psi) (g^{-1})^{0i}
			- \angGmixedarg{\Lunit}{\#} \cdot (\angdiff x^i) (\Lunit \Psi) 
			+ \frac{1}{2} G_{\Lunit \Lunit} (\angdiffuparg{\#} \Psi) \cdot \angdiff x^i.
				\label{E:LLUNITI} 
\end{align}


\end{lemma}

\begin{proof}
	We first prove \eqref{E:UPMUFIRSTTRANSPORT}. 
	We start by writing the $0$ component of the geodesic equation 
	$\D_{\Lgeo} \Lgeo = 0$
	relative to rectangular coordinates with the help of \eqref{E:CHRISTOFEELRECT}:
	$
	\Lgeo \Lgeo^0
		= 
			(g^{-1})^{0 \kappa} 
			\left\lbrace
				 (1/2) G_{\Lgeo \Lgeo} \partial_{\kappa} \Psi
				 -
				 G_{\kappa \Lgeo} \Lgeo \Psi
			\right\rbrace
	$. Using this equation, the identity $(g^{-1})^{0 \kappa} 
	= - \Lunit^{\kappa} - \Radunit^{\kappa}
	$
	(see \eqref{E:GINVERSEFRAMEWITHRECTCOORDINATESFORGSPHEREINVERSE}
	and recall that 
	$\Lunit^0 =1$ and $\Radunit^0 = 0$),
	the relation $\upmu = 1/\Lgeo^0$, 
	and the definition $\Lunit = \upmu \Lgeo$,
	we conclude \eqref{E:UPMUFIRSTTRANSPORT} from straightforward computations.

	To prove \eqref{E:LLUNITI}, we use the definition
	$\Lunit = \upmu \Lgeo$ 
	and \eqref{E:PERTURBEDPART} 
	to write the spatial components of the 
	geodesic equation $\D_{\Lgeo} \Lgeo = 0$ relative to rectangular components as
	$\Lunit \Lunit_{(Small)}^i
	  = \Lunit \Lunit^i 
		= - (g^{-1})^{i \kappa} \Gamma_{\Lunit \kappa \Lunit}
			+ \upmu^{-1} (\Lunit \upmu) \Lunit^i.
	$
	Using 
	\eqref{E:GINVERSEFRAMEWITHRECTCOORDINATESFORGSPHEREINVERSE}
	and
	\eqref{E:CHRISTOFEELRECT},
	we compute that
	$
	\displaystyle
	- (g^{-1})^{i \kappa} \Gamma_{\Lunit \kappa \Lunit}
	= (\Lunit^i + \Radunit^i) \Gamma_{\Lunit \Lunit \Lunit}
		+ \Lunit^i \Gamma_{\Lunit \Radunit \Lunit}
		- (\ginversesphere)^{i \kappa} \Gamma_{\Lunit \kappa \Lunit}
	$.
	Using \eqref{E:DOWNSTAIRSUPSTAIRSSRADUNITPLUSLUNITISAFUNCTIONOFPSI}
	and \eqref{E:CHRISTOFEELRECT}, we express
	the RHS of the previous identity as
	$
	- (1/2)(g^{-1})^{0i} G_{\Lunit \Lunit} \Lunit \Psi
	+ \Lunit^i \left\lbrace G_{\Lunit \Radunit} \Lunit \Psi - (1/2) G_{\Lunit \Lunit} \Radunit \Psi \right\rbrace
	- \angGmixedarg{\Lunit}{i} \Lunit \Psi
	+ (1/2) G_{\Lunit \Lunit} \angdiffuparg{i} \Psi
	$.
	We then add this expression to the second product
	$\upmu^{-1}(\Lunit \upmu) \Lunit^i$
	in the formula for $\Lunit \Lunit_{(Small)}^i$ from above
	and use \eqref{E:UPMUFIRSTTRANSPORT} to substitute for
	$\Lunit \upmu$.
	We note in particular that the terms proportional to
	$
	\Lunit^i G_{\Lunit \Lunit} \Radunit \Psi
	$
	and 
	$
	\Lunit^i G_{\Lunit \Radunit} \Lunit \Psi
	$
	completely cancel.
	Also using the simple identities
	$\angGmixedarg{\Lunit}{i} = \angGmixedarg{\Lunit}{\#} \cdot \angdiff x^i$
	and $\angdiffuparg{i} \Psi = (\angdiffuparg{\#} \Psi) \cdot \angdiff x^i$,
	we conclude \eqref{E:LLUNITI}.

\end{proof}

\subsection{Connection coefficients of the rescaled frame}
\label{SS:CONNECTIONCOEFFIENTSOFRESCALEDFRAME}
We now derive expressions for the connection coefficients 
of the frame $\lbrace \Lunit, \Rad, \CoordAng \rbrace$ in terms
of $\Psi,\upmu,\Lunit^1,\Lunit^2$.
We also decompose some of the connection coefficients
into ``regular'' pieces and pieces that have a ``singular'' $\upmu^{-1}$ 
factor.

\begin{lemma}[\textbf{Connection coefficients of the rescaled frame $\lbrace \Lunit, \Rad, \CoordAng \rbrace$} 
		\textbf{and their decomposition into} $\upmu^{-1}-$\textbf{singular and} $\upmu^{-1}-$\textbf{regular pieces}]
	\label{L:CONNECTIONLRADFRAME}
	Let $\upzeta$ be the $\ell_{t,u}-$tangent one-form defined by (see the identity \eqref{E:ANGKRINTERMSOFCOVARIANTDERIVATIVES})
	\begin{align} \label{E:ZETADEF}
		\upzeta_{\CoordAng} 
		& := \angkdoublearg{\Radunit}{\CoordAng}
		= g(\D_{\CoordAng} \Lunit, \Radunit) 
		= \upmu^{-1} g(\D_{\CoordAng} \Lunit, \Rad).
	\end{align}
	Then the covariant derivatives of the rescaled frame vectorfields can be expressed as follows,
	where the tensorfields 
	$k$, 
	$\upchi$, 
	and $\upomega$ are defined in 
	\eqref{E:SECONDFUNDSIGMATDEF}, 
	\eqref{E:CHIDEF},
	and \eqref{E:UPMUFIRSTTRANSPORT}:
	\begin{subequations}
	\begin{align}
		\D_{\Lunit} \Lunit 
		& = \upmu^{-1} \upomega \Lunit, 
			\label{E:DLL} \\
		\D_{\Rad} \Lunit 
		& = - \upomega \Lunit 
			+ \upmu \upzeta^{\#}
			+ \angdiffuparg{\#} \upmu, 
			\label{E:DRADL} \\
		\D_{\CoordAng} \Lunit 
		& = - \upzeta_{\CoordAng} \Lunit
			+ \mytr \upchi \CoordAng, 
			\label{E:DAL} \\
		\D_{\Lunit} \Rad 
		& = - \upomega \Lunit 
			- \upmu \upzeta^{\#}, 
			\label{E:DLRAD} \\
		\D_{\Rad} \Rad 
		& = \upmu \upomega \Lunit
			+ \left\lbrace 
					\upmu^{-1} \Rad \upmu  
					+ \upomega \right
				\rbrace \Rad
			- \upmu \angdiffuparg{\#} \upmu, 
			\label{E:DRADRAD} \\
		\D_{\CoordAng} \Rad 
		& = \upmu \upzeta_{\CoordAng} \Lunit
			+ \upzeta_{\CoordAng} \Rad
			+ \upmu^{-1} (\angdiffarg{\CoordAng} \upmu) \Rad
			+ \upmu \mytr \angk \CoordAng
			- \upmu \mytr \upchi \CoordAng,
			\label{E:DARAD} \\
		\D_{\Lunit} \CoordAng 
		& = \D_{\CoordAng} \Lunit, 
			\label{E:DLA} \\
		\D_{\CoordAng} \CoordAng
		& = \angDarg{\CoordAng} \CoordAng
			+ \angkdoublearg{\CoordAng}{\CoordAng} \Lunit
			+ \upmu^{-1} \upchi_{\CoordAng \CoordAng} \Rad.
			\label{E:DAXBINTERMSOFANGDAXB}
	\end{align}
	\end{subequations}

	Furthermore, we can decompose the frame components of the 
	$\ell_{t,u}-$tangent
	tensorfields $\angk$ and $\upzeta$ into 
	$\upmu^{-1}-$singular and $\upmu^{-1}-$regular
	pieces as follows: 
	\begin{subequations}
	\begin{align}
		\upzeta & = \upmu^{-1} \upzeta^{(Trans-\Psi)}
		+ \upzeta^{(Tan-\Psi)},
		 \label{E:ZETADECOMPOSED} \\
	\angk 
	& = \upmu^{-1} \angk^{(Trans-\Psi)}
		+ \angk^{(Tan-\Psi)},
		\label{E:ANGKDECOMPOSED}
	\end{align}
\end{subequations}
where
\begin{subequations}
	\begin{align}
		\upzeta^{(Trans-\Psi)} 
		& :=
			- \frac{1}{2} \angGarg{\Lunit} \Rad \Psi,
			\label{E:ZETATRANSVERSAL} \\
		\angk^{(Trans-\Psi)} 
		& := \frac{1}{2} \angG \Rad \Psi,
			\label{E:KABTRANSVERSAL}
	\end{align}
\end{subequations}
and
\begin{subequations}
\begin{align}
	\upzeta^{(Tan-\Psi)}
	& := \frac{1}{2} \angGarg{\Radunit} \Lunit \Psi
			- \frac{1}{2} G_{\Lunit \Radunit} \angdiff \Psi
			- \frac{1}{2} G_{\Radunit \Radunit} \angdiff \Psi,
		\label{E:ZETAGOOD} \\
	\angk^{(Tan-\Psi)} 
	& := \frac{1}{2} \angG \Lunit \Psi
			- \frac{1}{2} \angGarg{\Lunit} \otimes \angdiff \Psi
			- \frac{1}{2} \angdiff \Psi \otimes \angGarg{\Lunit} 
			- \frac{1}{2} \angGarg{\Radunit} \otimes \angdiff \Psi
			- \frac{1}{2} \angdiff \Psi \otimes \angGarg{\Radunit}.
			\label{E:KABGOOD}
\end{align}
\end{subequations}
\end{lemma}

\begin{proof}
	The identity \eqref{E:DLL} follows easily from 
	the geodesic equation $\D_{\Lgeo} \Lgeo = 0$
	and the definition $\Lunit = \upmu \Lgeo$.

	To derive \eqref{E:DLRAD}, we expand
	$\D_{\Lunit} \Rad = a_{\Lunit} \Lunit + a_{\Rad} \Rad + a_{\CoordAng} \CoordAng$,
	where the $a$ are scalar functions.
	Taking the inner product of each side with $\Lunit$
	and using $g(\Lunit,\Rad) = - \upmu$,
	$g(\Lunit,\CoordAng) = 0$,
	and \eqref{E:DLL},
	we find that $- a_{\Rad} \upmu =  g(\D_{\Lunit} \Rad, \Lunit) = - \Lunit \upmu - g(\Rad, \D_{\Lunit} \Lunit) = 0$
	as desired.
	Taking the inner product of each side with $\Rad$
	and using in addition that $g(\Rad,\Rad) = \upmu^2$ and $g(\Rad,\CoordAng) = 0$,
	we find that
	$- \upmu a_{\Lunit} = g(\D_{\Lunit} \Rad, \Rad) = \upmu \Lunit \upmu$ as desired.
	Finally, taking the inner product of each side with $\CoordAng$
	and using in addition that $0 = [\Lunit,\CoordAng] = \D_{\Lunit} \CoordAng - \D_{\CoordAng} \Lunit$
	(where the second equality follows from the torsion-free property of $\D$),
	we find that
	$a_{\CoordAng} g(\CoordAng,\CoordAng) 
		= g(\D_{\Lunit} \Rad,\CoordAng) 
		= - g(\Rad,\D_{\Lunit} \CoordAng)
		= - g(\Rad,\D_{\CoordAng} \Lunit)
		= - \upmu \upzeta_{\CoordAng}
	$
	as desired. 
	A similar argument yields \eqref{E:DRADL};
	we omit the full details and instead only note
	that the argument relies in part on the identity
	$g(\D_{\Rad} \Lunit, \CoordAng) 
	= - g(\D_{\Rad} \CoordAng,\Lunit)
	= - g(\D_{\CoordAng} \Rad,\Lunit)
		- g([\Rad, \CoordAng],\Lunit)
		= - g(\D_{\CoordAng} \Rad,\Lunit)
	$.
	The second equality follows from the torsion-free property of $\D$, 
	while the last one follows from
	the fact that $[\Rad, \CoordAng]$ is $\ell_{t,u}-$tangent,
	which is a simple consequence of
	\eqref{E:RADSPLITINTOPARTTILAUANDXI}.
	A similar argument also yields
	\eqref{E:DRADRAD}; we omit the details.

	To derive \eqref{E:DAL},
	we expand
	$\D_{\CoordAng} \Lunit = a_{\Lunit} \Lunit + a_{\Rad} \Rad + a_{\CoordAng} \CoordAng$.
	Taking the inner product of each side with $\Lunit$
	and using the identities noted above as well as $g(\Lunit,\Lunit) = 0$,
	we find that $a_{\Rad} = 0$ as desired. 
	Similarly, taking the inner product of each side with $\Rad$,
	we find that 
	$- \upmu a_{\Lunit} 
	= g(\D_{\CoordAng} \Lunit, \Rad) 
	= \upmu \upzeta_{\CoordAng}
	$
	as desired.
	Similarly, taking the inner product of each side with $\CoordAng$
	and using \eqref{E:CHIUSEFULID},
	we find that 
	$a_{\CoordAng} g(\CoordAng,\CoordAng) 
	= g(\D_{\CoordAng} \Lunit, \CoordAng)
	= \upchi_{\CoordAng \CoordAng} 
	$,
	from which we easily conclude that $a_{\CoordAng} = \mytr \upchi$ as desired.

	\eqref{E:DLA} is a simple consequence of the identity $[\Lunit,\CoordAng] = 0$
	and the torsion-free property of $\D$.

	To prove \eqref{E:DAXBINTERMSOFANGDAXB},
	we expand
	$\D_{\CoordAng} \CoordAng = a_{\Lunit} \Lunit + a_{\Rad} \Rad + a_{\CoordAng} \CoordAng$.
	Taking the inner product of each side with $\Lunit$
	and using the identities noted above,
	we find that
	$- a_{\Rad} \upmu 
	=  g(\D_{\CoordAng} \CoordAng, \Lunit) 
	= - g(\D_{\CoordAng} \Lunit, \CoordAng)
	= - \mytr \upchi g(\CoordAng, \CoordAng)
	= - \upchi_{\CoordAng \CoordAng}
	$
	as desired.
	Taking the inner product of each side with $\Rad$,
 	and using 
 	\eqref{E:RADDEF},
 	\eqref{E:TIMENORMAL},
 	and \eqref{E:DAL},
	we find that 
	$- \upmu a_{\Lunit} 
	= g(\D_{\CoordAng} \CoordAng,\Rad)
		- \upmu \upchi_{\CoordAng \CoordAng}
	= - g(\D_{\CoordAng} \Rad, \CoordAng)
		- \upmu \mytr \upchi g(\CoordAng, \CoordAng)
	= - \upmu g(\D_{\CoordAng} \Radunit, \CoordAng)
		- \upmu \mytr \upchi g(\CoordAng, \CoordAng)
	= - \upmu g(\D_{\CoordAng} \Timenormal, \CoordAng)
	= -\upmu \angkdoublearg{\CoordAng}{\CoordAng}
	$
	as desired.
	Taking the inner product of each side with $\CoordAng$ and using
	that $\angD \xi = \Lineproject \D \xi$ for $\ell_{t,u}-$tangent tensorfields $\xi$,
	we find that 
	$a_{\CoordAng} g(\CoordAng,\CoordAng) 
	= g(\D_{\CoordAng} \CoordAng, \CoordAng)
	= g(\angD_{\CoordAng} \CoordAng, \CoordAng)
	$
	as desired.

	To prove \eqref{E:DARAD}, we expand
	$\D_{\CoordAng} \Rad = a_{\Lunit} \Lunit + a_{\Rad} \Rad + a_{\CoordAng} \CoordAng$.
	Taking the inner product of each side with $\Lunit$ and using 
	and using the identities noted above as well as \eqref{E:DAL},
	we find that 
	$- a_{\Rad} \upmu 
	=  g(\D_{\CoordAng} \Rad, \Lunit) 
	= - \CoordAng \upmu - g(\D_{\CoordAng} \Lunit, \Rad)
	= - \CoordAng \upmu - \upmu \upzeta_{\CoordAng}
	$
	as desired.
	Taking the inner product of each side with $\Rad$,
	we find that
	$- \upmu a_{\Lunit} 
	= g(\D_{\CoordAng} \Rad, \Rad) 
		-
		\upmu^2 a_{\Rad}
	= - \upmu^2 \upzeta_{\CoordAng}
	$ as desired.
	Finally, taking the inner product of each side with $\CoordAng$
	and using \eqref{E:DAXBINTERMSOFANGDAXB},
	we find that
	$a_{\CoordAng} g(\CoordAng,\CoordAng) 
	= g(\D_{\CoordAng} \Rad, \CoordAng) 
	= - g(\Rad,\D_{\CoordAng} \CoordAng) 
	=  \upmu \angkdoublearg{\CoordAng}{\CoordAng} 
		- \upmu \upchi_{\CoordAng \CoordAng}
	=  \upmu \mytr \angk g(\CoordAng,\CoordAng)
		- \upmu \mytr \upchi g(\CoordAng,\CoordAng)
	$
	as desired. 

	We now prove \eqref{E:ZETADECOMPOSED}, \eqref{E:ZETATRANSVERSAL}, and \eqref{E:ZETAGOOD}.
	Our proof relies on the identity
	\begin{align} \label{E:LIETIMENORMALGRADUNITCOORDANGCOMPONENT}
		(\Lie_{\Timenormal} g)_{\Radunit \CoordAng} 
		& = 
			G_{\Radunit \CoordAng} \Lunit \Psi
			- G_{\Lunit \CoordAng} \Radunit \Psi
			- G_{\Lunit \Radunit} \CoordAng \Psi
			- G_{\Radunit \Radunit} \CoordAng \Psi.
	\end{align}
	To prove \eqref{E:LIETIMENORMALGRADUNITCOORDANGCOMPONENT},
	we use
	\eqref{E:TIMENORMAL},
	\eqref{E:TIMENORMALRECTANGULAR},
	\eqref{E:METRICFRAMEDECOMPLUNITRADUNITFRAME},
	\eqref{E:GINVERSEFRAMEWITHRECTCOORDINATESFORGSPHEREINVERSE},
	the chain rule identity 
	$\Timenormal g_{\alpha \beta} = G_{\alpha \beta} \Timenormal \Psi$,
	and the inverse matrix differentiation identity 
	$- V (g^{-1})^{0 \alpha} = (g^{-1})^{0 \kappa} (g^{-1})^{\alpha \lambda} G_{\kappa \lambda} V \Psi$
	(valid for any vectorfield $V$)
	to deduce that
	\begin{align}
		(\Lie_{\Timenormal} g)_{\Radunit \CoordAng}
		& = G_{\Radunit \CoordAng} \Timenormal \Psi
			+ g_{\Radunit \alpha} \CoordAng \Timenormal^{\alpha} 
			+ g_{\alpha \CoordAng} \Radunit \Timenormal^{\alpha}
				\\
		& = G_{\Radunit \CoordAng} \Lunit \Psi
			+ G_{\Radunit \CoordAng} \Radunit \Psi
			+ g_{\Radunit \alpha} (g^{-1})^{0 \kappa} (g^{-1})^{\alpha \lambda} G_{\kappa \lambda} \CoordAng \Psi
			+ g_{\alpha \CoordAng} (g^{-1})^{0 \kappa} (g^{-1})^{\alpha \lambda} G_{\kappa \lambda} \Radunit \Psi
			\notag \\
		& = G_{\Radunit \CoordAng} \Lunit \Psi
			+ G_{\Radunit \CoordAng} \Radunit \Psi
			- G_{\Lunit \Radunit} \CoordAng \Psi
			- G_{\Radunit \Radunit} \CoordAng \Psi
			- G_{\Lunit \CoordAng} \Radunit \Psi
			- G_{\Radunit \CoordAng} \Radunit \Psi.
		\notag
	\end{align}
	We have thus proved \eqref{E:LIETIMENORMALGRADUNITCOORDANGCOMPONENT}.
	We now use
	\eqref{E:ALTERNATESECONDFUND},
	\eqref{E:ANGKRINTERMSOFCOVARIANTDERIVATIVES},
	\eqref{E:ZETADEF},
	and \eqref{E:LIETIMENORMALGRADUNITCOORDANGCOMPONENT}
	to compute that
	$\upzeta_{\CoordAng} = (1/2) \mbox{RHS~\eqref{E:LIETIMENORMALGRADUNITCOORDANGCOMPONENT}}$,
	which easily yields \eqref{E:ZETADECOMPOSED}, \eqref{E:ZETATRANSVERSAL}, and \eqref{E:ZETAGOOD}.

	The identities \eqref{E:ANGKDECOMPOSED}, \eqref{E:KABTRANSVERSAL}, and \eqref{E:KABGOOD}
	can be proved by employing a similar argument; we omit the details.
\end{proof}

We will use the next lemma when deriving $L^{\infty}$ estimates
for the transversal derivatives of the rectangular component functions 
$\Lunit_{(Small)}^i$.

\begin{lemma}[\textbf{Formula for} $\Rad \Lunit^i$]
\label{L:RADLUNITI}
We have the following identity for the scalar-valued functions $\Lunit_{(Small)}^i$, ($i=1,2$):
\begin{align}
	\Rad \Lunit_{(Small)}^i
	& = \left\lbrace
				- \frac{1}{2} G_{\Lunit \Lunit} \Rad \Psi
				+ \frac{1}{2} \upmu G_{\Lunit \Lunit} \Lunit \Psi
				+ \upmu G_{\Lunit \Radunit} \Lunit \Psi
				+ \frac{1}{2} \upmu G_{\Radunit \Radunit} \Lunit \Psi
			\right\rbrace
			\Lunit^i
			\label{E:RADLUNITI} 	\\
	& \ \ 
			+
			\left\lbrace
				- \frac{1}{2} G_{\Lunit \Lunit} \Rad \Psi
				+ \frac{1}{2} \upmu G_{\Lunit \Lunit} \Lunit \Psi
				+ \upmu G_{\Lunit \Radunit} \Lunit \Psi
				+ \frac{1}{2} \upmu G_{\Radunit \Radunit} \Lunit \Psi
			\right\rbrace
			(g^{-1})^{0i}
			\notag \\
	& \ \ 
		 - \left\lbrace
					\angGmixedarg{\Lunit}{\#} \Rad \Psi
					+ \frac{1}{2} \upmu G_{\Radunit \Radunit} \angdiffuparg{\#} \Psi
				\right\rbrace
				\cdot
				\angdiff x^i
			+ (\angdiffuparg{\#} \upmu) \cdot \angdiff x^i.
			\notag
\end{align}
\end{lemma}

\begin{proof}
	Throughout this proof, $\nabla$ denotes the Levi-Civita connection of 
	the background Minkowski metric $m_{\alpha \beta} = \mbox{\upshape diag}(-1,1,1)$. 
	Since $\Lunit^0 = 1$, we can view 
	$\nabla_{\Rad} \Lunit$ as a $\Sigma_t-$tangent vectorfield
	with rectangular spatial components $\Rad \Lunit^i$.
	Since $\Radunit$ and $\CoordAng$ span the tangent space of $\Sigma_t$ at each point,
	we can expand 
	$\nabla_{\Rad} \Lunit = a_{\Radunit} \Radunit + a_{\CoordAng} \CoordAng$,
	where $a_{\Radunit}$ and $a_{\CoordAng}$ are scalar functions.
	Taking the inner product of each side with $\Radunit$, we find that
	$a_{\Radunit} 
	= g(\nabla_{\Rad} \Lunit, \Radunit)
	= g(\D_{\Rad} \Lunit, \Radunit)
		- \Gamma_{\Rad \Radunit \Lunit}
	$,
	where
	$
	\Gamma_{\Rad \Radunit \Lunit}
	:= \Rad^{\alpha} \Radunit^{\kappa} \Lunit^{\beta} \Gamma_{\alpha \kappa \beta}
	$
	and $ \Gamma_{\alpha \kappa \beta}$ is given by \eqref{E:CHRISTOFEELRECT}.
	Using 
	\eqref{E:CHRISTOFEELRECT},
	\eqref{E:UPMUFIRSTTRANSPORT},
	and	\eqref{E:DRADL},
	we compute that
	$
	g(\D_{\Rad} \Lunit, \Radunit)
	= \upomega 
	= \frac{1}{2} G_{\Lunit \Lunit} \Rad \Psi
		- \frac{1}{2} \upmu G_{\Lunit \Lunit} \Lunit \Psi
		- \upmu G_{\Lunit \Radunit} \Lunit \Psi
	$
	and 
	$
	\Gamma_{\Rad \Radunit \Lunit}
	= \frac{1}{2} \upmu G_{\Radunit \Radunit} \Lunit \Psi
	$.
	Hence, 
	$a_{\Radunit} 
	= 
		\frac{1}{2} G_{\Lunit \Lunit} \Rad \Psi
		- \frac{1}{2} \upmu G_{\Lunit \Lunit} \Lunit \Psi
		- \upmu G_{\Lunit \Radunit} \Lunit \Psi
		- \frac{1}{2} \upmu G_{\Radunit \Radunit} \Lunit \Psi
	$.
	Similarly, we find that
	$a_{\CoordAng} 
		= g(\D_{\Rad} \Lunit, \CoordAng)
		- \Gamma_{\Rad \CoordAng \Lunit}
	$.
	Using 
	\eqref{E:DRADL},
	we compute that
	$g(\D_{\Rad} \Lunit, \CoordAng)
		= 
		\upmu \upzeta_{\CoordAng} 
		+ 
		\angdiff_{\CoordAng} \upmu
	$
	and 
	$
		\Gamma_{\Rad \CoordAng \Lunit}
		= 
		\frac{1}{2}
		\left\lbrace
			\upmu G_{\Radunit \CoordAng} \Lunit \Psi
			+ G_{\Lunit \CoordAng} \Rad \Psi
			- \upmu G_{\Lunit \Radunit} \angdiff_{\CoordAng} \Psi
		\right\rbrace
	$.
	Hence,
	using 
	\eqref{E:ZETADECOMPOSED}, 
	\eqref{E:ZETATRANSVERSAL}, 
	and \eqref{E:ZETAGOOD}
	to substitute for $\upmu \upzeta$,
	we deduce that
	$a_{\CoordAng}
	= - G_{\Lunit \CoordAng} \Rad \Psi
		- \frac{1}{2} \upmu  G_{\Radunit \Radunit} \angdiff_{\CoordAng} \Psi
		+ 
		\angdiff_{\CoordAng} \upmu
	$.
	Combining these identities,
	using the identity $\CoordAng^i = \angdiff_{\CoordAng} x^i$, 
	and using \eqref{E:DOWNSTAIRSUPSTAIRSSRADUNITPLUSLUNITISAFUNCTIONOFPSI}
	to replace
	$\Radunit^i$ with $-\Lunit^i - (g^{-1})^{0i}$,
	we conclude \eqref{E:RADLUNITI}.

\end{proof}

\subsection{Useful expressions for the null second fundamental form}

\begin{lemma}[\textbf{Identities involving} $\upchi$]
	We have the following identities:
	\begin{subequations}
	\begin{align} \label{E:CHIINTERMSOFOTHERVARIABLES}
	\upchi 
	& = g_{ab} (\angdiff \Lunit^a) \otimes \angdiff x^b
		+ \frac{1}{2} \angG \Lunit \Psi,
			\\
	\mytr \upchi 
	& = g_{ab} \ginversesphere \cdot \left\lbrace (\angdiff \Lunit^a) \otimes \angdiff x^b \right\rbrace
		+ \frac{1}{2} \ginversesphere \cdot \angG \Lunit \Psi,
			 \label{E:TRCHIINTERMSOFOTHERVARIABLES}
			\\
	\Lunit \ln \gtancomp
	& = \mytr \upchi,
	\label{E:LDERIVATIVEOFVOLUMEFORMFACTOR}
\end{align}
\end{subequations}
where 
$\upchi$ is the $\ell_{t,u}-$tangent tensorfield defined by \eqref{E:CHIDEF}
and $\gtancomp$ is the metric component from Def.~\ref{D:METRICANGULARCOMPONENT}.
\end{lemma}

\begin{proof}
	To prove \eqref{E:CHIINTERMSOFOTHERVARIABLES}, we use
	\eqref{E:CHIUSEFULID} and \eqref{E:CHRISTOFEELRECT} to compute, 
	relative to rectangular coordinates, 
	that
	$\upchi_{\CoordAng \CoordAng} 
	= g_{ab} (\CoordAng \Lunit^a) \CoordAng^b
	 + \Gamma_{\CoordAng \CoordAng \Lunit} 
	=  g_{ab} (\CoordAng \Lunit_{(Small)}^a) \CoordAng^b
		+ \frac{1}{2} G_{\CoordAng \CoordAng} \Lunit \Psi
	$.
	Noting that $\CoordAng x^b = \CoordAng \cdot \angdiff x^b$, 
	we easily conclude \eqref{E:CHIINTERMSOFOTHERVARIABLES}.
	To deduce \eqref{E:TRCHIINTERMSOFOTHERVARIABLES}, we simply take the
	$\gsphere-$trace of  \eqref{E:CHIINTERMSOFOTHERVARIABLES}.
	To prove \eqref{E:LDERIVATIVEOFVOLUMEFORMFACTOR}, we use the Leibniz rule,
	the fact that $[\Lunit,\CoordAng] = 0$, and the torsion-free property of $\D$
	to compute that
	$\Lunit (\gtancomp^2) 
	= \Lunit [g(\CoordAng,\CoordAng)] 
	= 2 g(\D_{\Lunit} \CoordAng,\CoordAng)
	= 2 g(\D_{\CoordAng} \Lunit,\CoordAng)
	= 2 \upchi_{\CoordAng \CoordAng}
	= 2 \mytr \upchi g(\CoordAng, \CoordAng)
	= 2 \mytr \upchi \gtancomp^2,
	$
	from which the desired identity easily follows.
\end{proof}

\subsection{Frame decomposition of the wave operator}
\label{SS:FRAMEDCOMPOFBOX}
In this subsection, we decompose $\upmu \square_{g(\Psi)}$ relative to the rescaled frame.
The factor of $\upmu$ is important for our decompositions.

\begin{proposition} [\textbf{Frame decomposition of $\upmu \square_{g(\Psi)} f$}]
	\label{P:GEOMETRICWAVEOPERATORFRAMEDECOMPOSED}
	Let $f$ be a scalar function.
	Then relative to the rescaled frame $\lbrace \Lunit, \Rad, \CoordAng \rbrace$,
	$\upmu \square_{g(\Psi)} f$ can be expressed in either of the following two forms:
	\begin{subequations}
	\begin{align} \label{E:LONOUTSIDEGEOMETRICWAVEOPERATORFRAMEDECOMPOSED}
		\upmu \square_{g(\Psi)} f 
			& = - \Lunit(\upmu \Lunit f + 2 \Rad f)
				+ \upmu \angLap f
				- \mytr \upchi \Rad f
				- \upmu \mytr \angk \Lunit f
				- 2 \upmu \upzeta^{\#} \cdot \angdiff f,
					\\
			& = - (\upmu \Lunit + 2 \Rad) (\Lunit f) 
				+ \upmu \angLap f
				- \mytr \upchi \Rad f
				- \upomega \Lunit f
				+ 2 \upmu \upzeta^{\#} \cdot \angdiff f
				+ 2 (\angdiffuparg{\#} \upmu) \cdot \angdiff f,
				\label{E:LONINSIDEGEOMETRICWAVEOPERATORFRAMEDECOMPOSED}
	\end{align}
	\end{subequations}
	where the $\ell_{t,u}-$tangent tensorfields
	$\upchi$,
	$\upzeta$,
	and
	$\angk$
	can be expressed via
	\eqref{E:CHIINTERMSOFOTHERVARIABLES},
	\eqref{E:ZETADECOMPOSED},
	and
	\eqref{E:ANGKDECOMPOSED}.
\end{proposition}

\begin{proof}
	To derive \eqref{E:LONOUTSIDEGEOMETRICWAVEOPERATORFRAMEDECOMPOSED}, 
	we first use \eqref{E:GINVERSEFRAMEWITHRECTCOORDINATESFORGSPHEREINVERSE}
	to decompose
	\begin{align} \label{E:FIRSTEXPRESSIONLONOUTSIDEGEOMETRICWAVEOPERATORFRAMEDECOMPOSED}
		\upmu \square_{g(\Psi)} f
		& = - \upmu \Lunit^{\alpha} \Lunit^{\beta} \D_{\alpha \beta}^2 f 
			- 2 \Lunit^{\alpha} \Rad^{\beta} \D_{\alpha \beta}^2 f 
			+ (\ginversesphere) \cdot \D^2 f
				\\
		& = -  \Lunit(\upmu \Lunit f + 2 \Rad f)
			+ (\ginversesphere) \cdot \D^2 f
			\notag \\
	& \ \
			+ \upmu (\D_{\Lunit} \Lunit^{\alpha}) \D_{\alpha} f
			+ 2 (\D_{\Lunit} \Rad^{\alpha}) \D_{\alpha} f
			+ (\Lunit \upmu) \Lunit f.
			\notag
	\end{align}
Next, we note that
$\CoordAng(\CoordAng f) 
= \CoordAng^{\alpha} \D_{\alpha}(\CoordAng^{\beta} \D_{\beta} f)
= \D_{\CoordAng \CoordAng}^2 f
	+ (\D_{\CoordAng} \CoordAng)^{\alpha} \D_{\alpha} f
= \angD_{\CoordAng \CoordAng}^2 f
	+ (\angD_{\CoordAng} \CoordAng) \cdot \angdiff f
$.
Hence, by \eqref{E:DAXBINTERMSOFANGDAXB}, we have
$
\D_{\CoordAng \CoordAng}^2 f
= \angD_{\CoordAng \CoordAng}^2 f
	- \angkdoublearg{\CoordAng}{\CoordAng} \Lunit f
	- \upmu^{-1} \upchi_{\CoordAng \CoordAng} \Rad f
$.
Consequently, 
$\upmu (\ginversesphere) \cdot \D^2 f 
= (1/g(\CoordAng,\CoordAng)) 
\left\lbrace 
	\upmu \angD_{\CoordAng \CoordAng}^2 f
	- \upmu \angkdoublearg{\CoordAng}{\CoordAng} \Lunit f
	- \upchi_{\CoordAng \CoordAng} \Rad f
\right\rbrace
= \upmu \angLap f
	- \upmu \mytr \angk \Lunit f
	- \mytr \upchi \Rad f
$.
We now substitute this identity into
RHS~\eqref{E:FIRSTEXPRESSIONLONOUTSIDEGEOMETRICWAVEOPERATORFRAMEDECOMPOSED}.
We also use Lemma~\ref{L:CONNECTIONLRADFRAME}
to substitute for
the terms 
$\upmu (\D_{\Lunit} \Lunit^{\alpha}) \D_{\alpha} f$
and $2 (\D_{\Lunit} \Rad^{\alpha}) \D_{\alpha} f$.
The identity \eqref{E:LONOUTSIDEGEOMETRICWAVEOPERATORFRAMEDECOMPOSED} then follows
from straightforward calculations.

The proof of \eqref{E:LONINSIDEGEOMETRICWAVEOPERATORFRAMEDECOMPOSED} is similar
and we omit the details.

\end{proof}



\subsection{Frame components of the deformation tensors of the commutation vectorfields}
\label{SS:FRAMECOMPONENTSOFDEFORMTEN}
In this subsection, we decompose the deformation tensors 
(see Def.~\ref{D:DEFORMATIONTENSOR})
of the commutation vectorfields
\eqref{E:COMMUTATIONVECTORFIELDS} relative to the rescaled frame. The exact structure
of a few of the terms, including the precise numerical constants, affects the
degree of degeneracy of our top-order energy estimates.

The main result of this subsection is Lemma~\ref{L:DEFORMATIONTENSORFRAMECOMPONENTS}.
We first provide a preliminary lemma in which 
we calculate certain covariant derivatives of the $\ell_{t,u}$ projection tensorfield $\Lineproject$.

\begin{lemma}[\textbf{Frame covariant derivatives of} $\Lineproject$]
\label{L:FRAMECOVARIANTDERIVATIVESOFSPHERICALPROJECTION}
Let $\Lineproject$ be the type $\binom{1}{1}$ $\ell_{t,u}$ projection tensorfield defined in \eqref{E:LINEPROJECTION}.
Then the following identities hold:
\begin{subequations}
\begin{align}
	(\D_{\Lunit} \Lineproject) \cdot \Radunit
	& = \upmu^{-1} \upzeta^{(Trans-\Psi) \# }
		+ \upzeta^{(Tan-\Psi) \#},
		\label{E:DLLINEPROJECTAPPLIEDTORADUNIT}
				\\
	(\D_{\Lunit} \Lineproject) \cdot \CoordAng
	& = - \upmu^{-1} \upzeta_{\CoordAng}^{(Trans-\Psi)} \Lunit
		- \upzeta_{\CoordAng}^{(Tan-\Psi)} \Lunit,
		\label{E:DLLINEPROJECTAPPLIEDTOCOORDANG}
				\\
	(\D_{\Rad} \Lineproject) \cdot \Radunit
	& = \angdiffuparg{\#} \upmu,
		\label{E:DRADLINEPROJECTAPPLIEDTORADUNIT} \\
	(\D_{\Rad} \Lineproject) \cdot \CoordAng
	& = \upzeta_{\CoordAng}^{(Trans-\Psi)} \Lunit
		+ \upmu \upzeta_{\CoordAng}^{(Tan-\Psi)} \Lunit
		+ \upzeta_{\CoordAng}^{(Trans-\Psi)} \Radunit
		+ \upmu \upzeta_{\CoordAng}^{(Tan-\Psi)} \Radunit
		+ (\angdiffarg{\CoordAng} \upmu) \Radunit,
		\label{E:DRADLINEPROJECTAPPLIEDTOA} \\
	(\D_{\CoordAng} \Lineproject) \cdot \Radunit
	& = \upchi_{\CoordAng}^{\ \#} 
		- \upmu^{-1} \angkmixedarg{\CoordAng}{(Trans-\Psi) \#} 
		- \angkmixedarg{\CoordAng}{(Tan-\Psi) \#},
		\label{E:DALINEPROJECTAPPLIEDTORADUNIT} \\
	(\D_{\CoordAng} \Lineproject) \cdot \CoordAng
	& = \upmu^{-1} \angktriplearg{\CoordAng}{\CoordAng}{(Trans-\Psi)} \Lunit
		+ \angktriplearg{\CoordAng}{\CoordAng}{(Tan-\Psi)} \Lunit
		+ \upchi_{\CoordAng \CoordAng} \Radunit.
		\label{E:DALINEPROJECTAPPLIEDTOCOORDANG}
\end{align}
\end{subequations}
	In the above expressions, the $\ell_{t,u}-$tangent tensorfields
	$\upchi$, 
	$\upzeta^{(Trans-\Psi)}$,
	$\angkuparg{(Trans-\Psi)}$, 
	$\upzeta^{(Tan-\Psi)}$,
	and
	$\angkuparg{(Tan-\Psi)}$
	are defined by
	\eqref{E:CHIDEF},
	\eqref{E:ZETATRANSVERSAL},
	\eqref{E:KABTRANSVERSAL},
	\eqref{E:ZETAGOOD},
	and \eqref{E:KABGOOD}.
\end{lemma}

\begin{proof}
The main idea of the proof is to use the decompositions provided by Lemma~\ref{L:CONNECTIONLRADFRAME}.
As examples, we prove \eqref{E:DLLINEPROJECTAPPLIEDTORADUNIT} and \eqref{E:DRADLINEPROJECTAPPLIEDTOA}. 
The remaining identities in the lemma
can be proved using similar arguments and we omit those details.
To prove \eqref{E:DLLINEPROJECTAPPLIEDTORADUNIT},
we differentiate the identity
$\Lineproject \cdot \Rad = 0$
and use the identity $\Rad = \upmu \Radunit$
to deduce that
$(\D_{\Lunit} \Lineproject) \cdot \Radunit 
= \upmu^{-1} (\D_{\Lunit} \Lineproject) \cdot \Rad
= - \upmu^{-1} \Lineproject \cdot \D_{\Lunit} \Rad$.
The desired identity \eqref{E:DLLINEPROJECTAPPLIEDTORADUNIT}
now follows easily from the previous identity, \eqref{E:DLRAD}, and \eqref{E:ZETADECOMPOSED}.

To prove \eqref{E:DRADLINEPROJECTAPPLIEDTOA},
we differentiate the identity $\Lineproject \cdot \CoordAng = \CoordAng$ to deduce 
$(\D_{\Rad} \Lineproject) \cdot \CoordAng = \D_{\Rad} \CoordAng - \Lineproject \cdot \D_{\Rad} \CoordAng$.
Since $\D_{\Rad} \CoordAng - \D_{\CoordAng} \Rad = [\Rad,\CoordAng]$ is $\ell_{t,u}-$tangent
(see Lemma~\ref{L:CONNECTIONBETWEENCOMMUTATORSANDDEFORMATIONTENSORS}),
it follows that 
$(\D_{\Rad} \Lineproject) \cdot \CoordAng = \D_{\CoordAng} \Rad - \Lineproject \cdot \D_{\CoordAng} \Rad$.
The desired identity \eqref{E:DRADLINEPROJECTAPPLIEDTOA} now
follows easily from the previous identity, 
\eqref{E:DARAD}, 
and
\eqref{E:ZETADECOMPOSED}.

\end{proof}

We now provide the main lemma of Subsect.~\ref{SS:FRAMECOMPONENTSOFDEFORMTEN}.

\begin{lemma}[\textbf{The frame components of} $\deform{Z}$]
\label{L:DEFORMATIONTENSORFRAMECOMPONENTS}
We have the following identities for the frame components of the 
deformation tensors (see Def.~\ref{D:DEFORMATIONTENSOR}) of the 
commutation vectorfields $Z \in \Fullset$ 
(see definition \eqref{E:COMMUTATIONVECTORFIELDS}):
\begin{subequations}
\begin{align}
	\deformarg{\Rad}{\Lunit}{\Lunit} & = 0, 
		\qquad
	\deformarg{\Rad}{\Rad}{\Radunit} 
		= 2 \Rad \upmu,
		\qquad
	\deformarg{\Rad}{\Lunit}{\Rad} 
		= - \Rad \upmu,
		\\
	\angdeformoneformarg{\Rad}{\Lunit}
	& = - \angdiff \upmu 
		- 2 \upzeta^{(Trans-\Psi)}
		- 2 \upmu \upzeta^{(Tan-\Psi)},
	\qquad
	\angdeformoneformarg{\Rad}{\Rad}
		= 0,
		\label{E:RADDEFORMSPHERERAD} \\
	\angdeform{\Rad}
	& = - 2 \upmu \upchi 
		+ 2 \angk^{(Trans-\Psi)}
		+ 2 \upmu \angk^{(Tan-\Psi)},
\end{align}
\end{subequations}

\begin{subequations}
\begin{align}
	\deformarg{\Lunit}{\Lunit}{\Lunit}
	& = 0, 
		\qquad
	\deformarg{\Lunit}{\Rad}{\Radunit} 
		= 2 \Lunit \upmu,
		\qquad
	\deformarg{\Lunit}{\Lunit}{\Rad} 
		= - \Lunit \upmu,
		\label{E:LUNITDEFORMSCALARS} \\
	\angdeformoneformarg{\Lunit}{\Lunit}
	& = 0, 
		\qquad
	\angdeformoneformarg{\Lunit}{\Rad}
		= \angdiff \upmu
	 		+ 2 \upzeta^{(Trans-\Psi)}
			+ 2 \upmu \upzeta^{(Tan-\Psi)}, 
	 		\label{E:LUNITDEFORMSPHERELUNITANDLUNITDEFORMSPHERERAD} \\
	\angdeform{\Lunit}
	& = 2 \upchi,
	\label{E:LUNITDEFORMSPHERE}
\end{align}
\end{subequations}

\begin{subequations}
\begin{align}
	\deformarg{\GeoAng}{\Lunit}{\Lunit}
	& = 0, 
		\qquad
	\deformarg{\GeoAng}{\Rad}{\Radunit} 
		= 2 \GeoAng \upmu,
		\qquad
	\deformarg{\GeoAng}{\Lunit}{\Rad} 
		= - \GeoAng \upmu,
		\label{E:GEOANGDEFORMSCALARS} \\
	\angdeformoneformarg{\GeoAng}{\Lunit}
	& = - \upchi \cdot \GeoAng
		+ \frac{1}{2} (\angG \cdot \GeoAng) \Lunit \Psi
		+ \GeoAngFlatRadComponent \angGarg{\Radunit} \Lunit \Psi
		+ \frac{1}{2} (\angGarg{\Lunit} \cdot \GeoAng) \angdiff \Psi
		- \GeoAngFlatRadComponent G_{\Lunit \Radunit} \angdiff \Psi
		- \frac{1}{2} \GeoAngFlatRadComponent G_{\Radunit \Radunit} \angdiff \Psi, 
		\label{E:GEOANGDEFORMSPHEREL} \\
	\angdeformoneformarg{\GeoAng}{\Rad}
	& = \upmu \upchi \cdot \GeoAng
			+ \GeoAngFlatRadComponent \angdiff \upmu
			+ \GeoAngFlatRadComponent \angGarg{\Radunit} \Rad \Psi
			- \frac{1}{2} \upmu \GeoAngFlatRadComponent G_{\Radunit \Radunit} \angdiff \Psi
			\label{E:GEOANGDEFORMSPHERERAD} \\
	& \ \
			- \frac{1}{2} \upmu (\angG \cdot \GeoAng) \Lunit \Psi
			+ \upmu (\angGarg{\Lunit} \cdot \GeoAng) \angdiff \Psi
			+ \upmu (\angGarg{\Radunit} \cdot \GeoAng) \angdiff \Psi, 
	 		\notag \\
	\angdeform{\GeoAng}
	& = 2 \GeoAngFlatRadComponent \upchi
		+ \frac{1}{2} (\angG \cdot \GeoAng) \otimes \angdiff \Psi
		+ \frac{1}{2} \angdiff \Psi \otimes (\angG \cdot \GeoAng) 
		- \GeoAngFlatRadComponent \angG \Lunit \Psi
			\label{E:GEOANGDEFORMSPHERE} \\
	& \ \
		+ \GeoAngFlatRadComponent \angGarg{\Lunit} \otimes \angdiff \Psi
		+ \GeoAngFlatRadComponent \angdiff \Psi \otimes \angGarg{\Lunit} 
		+ \GeoAngFlatRadComponent \angGarg{\Radunit} \otimes \angdiff \Psi
		+ \GeoAngFlatRadComponent \angdiff \Psi \otimes \angGarg{\Radunit}.
		\notag 
\end{align}
\end{subequations}
	The scalar function $\GeoAngFlatRadComponent$ from above
	is as in Lemma~\ref{L:GEOANGDECOMPOSITION}, while
	the $\ell_{t,u}-$tangent tensorfields
	$\upchi$,
	$\upzeta^{(Trans-\Psi)}$,
	$\angk^{(Trans-\Psi)}$,
	$\upzeta^{(Tan-\Psi)}$,
	and
	$\angk^{(Tan-\Psi)}$
from above 
are as in
\eqref{E:CHIDEF},
\eqref{E:ZETATRANSVERSAL},
\eqref{E:KABTRANSVERSAL},
\eqref{E:ZETAGOOD},
and \eqref{E:KABGOOD}.
\end{lemma}

\begin{proof}
	We give a detailed proof of the identities \eqref{E:GEOANGDEFORMSCALARS}-\eqref{E:GEOANGDEFORMSPHERE},
	some of which involve the observation of important cancellations.
	The proofs of the remaining identities do not involve such cancellations. Hence, they
	are easier to prove and we omit those details.

	First, we deduce
	$\deformarg{\GeoAng}{\Lunit}{\Lunit}
	= 2 g(\D_{\Lunit} \GeoAng, \Lunit)
	= - 2 g(\D_{\Lunit} \Lunit,\GeoAng)
	$,
	where to obtain the second equality,
	we differentiated the identity $g(\Lunit,\GeoAng) = 0$.
	Using \eqref{E:DLL}, we conclude that
	$g(\D_{\Lunit} \Lunit,\GeoAng) = 0$
	as desired.

	Next, we use similar reasoning to obtain
	$
		\deformarg{\GeoAng}{\Rad}{\Radunit}
		= 2 \upmu^{-1} g(\D_{\Rad} \GeoAng, \Rad)
		= - 2 \upmu^{-1} g(\D_{\Rad} \Rad, \GeoAng)
	$
	Using \eqref{E:DRADRAD},
	we conclude that the previous expression is equal to
	$2 g(\angdiff^{\#} \upmu, \GeoAng)
	= 2 \GeoAng \upmu
	$
	as desired.

	Next, we use similar reasoning to obtain
	$\deformarg{\GeoAng}{\Lunit}{\Rad}
	= g(\D_{\Lunit} \GeoAng, \Rad) 
		+
		g(\D_{\Rad} \GeoAng, \Lunit)
	= - g(\D_{\Lunit} \Rad, \GeoAng) 
		-
		g(\D_{\Rad} \Lunit, \GeoAng)
	$. 
	Using \eqref{E:DRADL} and \eqref{E:DLRAD}, 
	we conclude that the previous expression is equal to 
	$- g(\angdiff^{\#} \upmu,\GeoAng) = - \GeoAng \upmu$ as desired.

	Next, we use similar reasoning to obtain
	$g(\angdeformoneformupsharparg{\GeoAng}{\Lunit},\CoordAng)
	= g(\D_{\Lunit} \GeoAng, \CoordAng)
		+ g(\D_{\CoordAng} \GeoAng, \Lunit)
	= g(\D_{\Lunit} \GeoAng, \CoordAng)
		- g(\D_{\CoordAng} \Lunit, \GeoAng)
	$.
	By \eqref{E:DAL}, we have
	$- g(\D_{\CoordAng} \Lunit, \GeoAng) = - \upchi_{\GeoAng \CoordAng}$.
	From definition \eqref{E:GEOANGDEF}, we derive
	\begin{align} \label{E:DLGEOANGCOORDANGFIRSTCALC}
		g(\D_{\Lunit} \GeoAng, \CoordAng)
		= g((\D_{\Lunit} \Lineproject) \cdot \GeoAng_{(Flat)}, \CoordAng) 
			+ g(\D_{\Lunit} \GeoAng_{(Flat)}, \CoordAng).
	\end{align}
	Using 
	\eqref{E:GEOANGINTERMSOFEUCLIDEANANGANDRADUNIT},
	\eqref{E:DLLINEPROJECTAPPLIEDTORADUNIT},
	and 
	\eqref{E:DLLINEPROJECTAPPLIEDTOCOORDANG},
	we compute that
	\begin{align} \label{E:DLLINEPROJECTGEOANGFLATCOORDANGSINGULARTERM}
		g((\D_{\Lunit} \Lineproject) \cdot \GeoAng_{(Flat)}, \CoordAng)
		& = 
		\upmu^{-1} \GeoAngFlatRadComponent \upzeta_{\CoordAng}^{(Trans-\Psi)} 
		+
		\GeoAngFlatRadComponent \upzeta_{\CoordAng}^{(Tan-\Psi)}.
	\end{align}
	Next, using that $\Lunit \GeoAng_{(Flat)}^i = 0$
	and \eqref{E:CHRISTOFEELRECT}, we compute that
	\begin{align} \label{E:GOFDLUNITGEOANGFLATCOORDANG}
		g(\D_{\Lunit} \GeoAng_{(Flat)},\CoordAng)
		& 
		= \Gamma_{\Lunit \CoordAng \GeoAng_{(Flat)}}
		= \Gamma_{\Lunit \CoordAng \GeoAng}
			+ \GeoAngFlatRadComponent \Gamma_{\Lunit \CoordAng \Radunit}
			\\
		& = \frac{1}{2} \angGdoublearg{\Lunit}{\CoordAng} \GeoAng \Psi
			+ \frac{1}{2} \angGdoublearg{\CoordAng}{\GeoAng} \Lunit \Psi
			- \frac{1}{2} G_{\Lunit \Radunit} \CoordAng \Psi
				\notag \\
		& \ \
			+ \frac{1}{2} \upmu^{-1} \GeoAngFlatRadComponent \angGdoublearg{\Lunit}{\CoordAng} \Rad \Psi
			+ \frac{1}{2} \GeoAngFlatRadComponent \angGdoublearg{\Radunit}{\CoordAng} \Lunit \Psi
			- \frac{1}{2} \GeoAngFlatRadComponent G_{\Lunit \Radunit} \CoordAng \Psi.
			\notag
	\end{align}
	Combining
	the above calculations,
	noting that
	the term $\upmu^{-1} \GeoAngFlatRadComponent \upzeta_{\CoordAng}^{(Trans-\Psi)}$
	exactly cancels the dangerous term 
	$
		\frac{1}{2} \upmu^{-1} \GeoAngFlatRadComponent \angGdoublearg{\Lunit}{\CoordAng} \Rad \Psi
	$ 
	on RHS~\eqref{E:GOFDLUNITGEOANGFLATCOORDANG}
	(see \eqref{E:ZETATRANSVERSAL}),
	using \eqref{E:ZETAGOOD},
	and noting that 
	$
	\angGdoublearg{\Lunit}{\CoordAng} \GeoAng \Psi
	=
	(\angGarg{\Lunit} \cdot \GeoAng) \CoordAng \Psi
	$,
	we conclude \eqref{E:GEOANGDEFORMSPHEREL}.

	Next, we use similar reasoning to obtain
	$g(\angdeformoneformupsharparg{\GeoAng}{\Rad},\CoordAng)
	= g(\D_{\Rad} \GeoAng, \CoordAng)
		+ g(\D_{\CoordAng} \GeoAng, \Rad)
	= g(\D_{\Rad} \GeoAng, \CoordAng)
		- g(\D_{\CoordAng} \Rad, \GeoAng)
	$.
	By \eqref{E:DARAD}, we have
	$- g(\D_{\CoordAng} \Rad, \GeoAng) 
		= - \upmu \angkdoublearg{\GeoAng}{\CoordAng} 
			+ \upmu \upchi_{\CoordAng \GeoAng}$.
	From definition \eqref{E:GEOANGDEF}, we derive
	\begin{align} \label{E:DRADGEOANGCOORDANGFIRSTCALC}
		g(\D_{\Rad} \GeoAng, \CoordAng)
		= g((\D_{\Rad} \Lineproject) \cdot \GeoAng_{(Flat)}, \CoordAng) 
			+ g(\D_{\Rad} \GeoAng_{(Flat)}, \CoordAng).
	\end{align}
	Using 
	\eqref{E:GEOANGINTERMSOFEUCLIDEANANGANDRADUNIT},
	\eqref{E:DRADLINEPROJECTAPPLIEDTORADUNIT},
	and
	\eqref{E:DRADLINEPROJECTAPPLIEDTOA},
	we compute that
	\begin{align} \label{E:DRADLINEPROJECTGEOANGFLATCOORDANGSINGULARTERM}
		g((\D_{\Rad} \Lineproject) \cdot \GeoAng_{(Flat)}, \CoordAng)
		& = \GeoAngFlatRadComponent \CoordAng \upmu.
	\end{align}
	Next, using that $\Rad \GeoAng_{(Flat)}^i = 0$
	and \eqref{E:CHRISTOFEELRECT}, we compute that
	\begin{align} \label{E:GOFRADGEOANGFLATCOORDANG}
		g(\D_{\Rad} \GeoAng_{(Flat)},\CoordAng)
		& 
		= \Gamma_{\Rad \CoordAng \GeoAng_{(Flat)}}
		= \Gamma_{\Rad \CoordAng \GeoAng}
			+ \GeoAngFlatRadComponent \Gamma_{\Rad \CoordAng \Radunit}
			\\
		& = \frac{1}{2} \upmu \angGdoublearg{\Radunit}{\CoordAng} \GeoAng \Psi
			+ \frac{1}{2} \angGdoublearg{\GeoAng}{\CoordAng} \Rad \Psi
			- \frac{1}{2} \upmu (\angGarg{\Radunit} \cdot \GeoAng) \CoordAng \Psi
				\notag \\
		& \ \
			+ \GeoAngFlatRadComponent \angGdoublearg{\Radunit}{\CoordAng} \Rad \Psi
			- \frac{1}{2} \upmu \GeoAngFlatRadComponent G_{\Radunit \Radunit} \CoordAng \Psi.
			\notag
	\end{align}
	Combining the above calculations,
	using \eqref{E:ANGKDECOMPOSED}, \eqref{E:KABTRANSVERSAL}, and \eqref{E:KABGOOD}
	to substitute for $- \upmu \angkdoublearg{\GeoAng}{\CoordAng}$,
	and noting that 
	$
	\angGdoublearg{\Radunit}{\CoordAng} \GeoAng \Psi
	=
	(\angGarg{\Radunit} \cdot \GeoAng) \CoordAng \Psi
	$,
	we conclude \eqref{E:GEOANGDEFORMSPHERERAD}.

	Next, we note that 
	$\angdeformarg{\GeoAng}{\CoordAng}{\CoordAng}
	= 2 g(\D_{\CoordAng} \GeoAng, \CoordAng)
	$.
	From definition \eqref{E:GEOANGDEF}, we derive
	\begin{align} \label{E:DCORDANGGEOANGCOORDANGFIRSTCALC}
		g(\D_{\CoordAng} \GeoAng, \CoordAng)
		= g((\D_{\CoordAng} \Lineproject) \cdot \GeoAng_{(Flat)}, \CoordAng) 
			+ g(\D_{\CoordAng} \GeoAng_{(Flat)}, \CoordAng) .
	\end{align}
	Using 
	\eqref{E:GEOANGINTERMSOFEUCLIDEANANGANDRADUNIT},
	\eqref{E:DALINEPROJECTAPPLIEDTORADUNIT},
	and
	\eqref{E:DALINEPROJECTAPPLIEDTOCOORDANG},
	we compute that
	\begin{align} \label{E:DCOORDANGLINEPROJECTGEOANGFLATCOORDANGSINGULARTERM}
		g((\D_{\CoordAng} \Lineproject) \cdot \GeoAng_{(Flat)}, \CoordAng)
		& = \GeoAngFlatRadComponent \upchi_{\CoordAng \CoordAng} 
		- \upmu^{-1} \GeoAngFlatRadComponent \angkdoublearg{\CoordAng}{\CoordAng}^{(Trans-\Psi)} 
		- \GeoAngFlatRadComponent \angkdoublearg{\CoordAng}{\CoordAng}^{(Tan-\Psi) \#}.
	\end{align}
	Next, using that $\CoordAng \GeoAng_{(Flat)}^i = 0$
	and \eqref{E:CHRISTOFEELRECT}, we compute that
	\begin{align} \label{E:GOFCOORDANGGEOANGFLATCOORDANG}
		g(\D_{\CoordAng} \GeoAng_{(Flat)},\CoordAng)
		& 
		= \Gamma_{\CoordAng \CoordAng \GeoAng_{(Flat)}}
		= \Gamma_{\CoordAng \CoordAng \GeoAng}
			+ \GeoAngFlatRadComponent \Gamma_{\CoordAng \CoordAng \Radunit}
			\\
		& = \frac{1}{2} \angGdoublearg{\CoordAng}{\GeoAng} \CoordAng \Psi
			+ \frac{1}{2} \upmu^{-1} \GeoAngFlatRadComponent \angGdoublearg{\CoordAng}{\CoordAng} \Rad \Psi.
			\notag
	\end{align}
	Combining
	the above calculations,
	noting that the term 
	$- \upmu^{-1} \GeoAngFlatRadComponent \angkdoublearg{\CoordAng}{\CoordAng}^{(Trans-\Psi)}$
	on RHS~\eqref{E:DCOORDANGLINEPROJECTGEOANGFLATCOORDANGSINGULARTERM}
	exactly cancels the dangerous term 
	$\frac{1}{2} \upmu^{-1} \GeoAngFlatRadComponent \angGdoublearg{\CoordAng}{\CoordAng} \Rad \Psi$
	on RHS~\eqref{E:GOFCOORDANGGEOANGFLATCOORDANG}
	(see \eqref{E:KABTRANSVERSAL}),
	and using \eqref{E:KABGOOD},
	we conclude \eqref{E:GEOANGDEFORMSPHERE}.

	\end{proof}


\subsection{Arrays of fundamental unknowns}
\label{SS:ARRAYS}
Our goal in this subsection is to show that
many scalar functions and tensorfields
that we have introduced depend on just a handful of more fundamental functions and tensorfields.
This reduction highlights the structures that are relevant 
for deriving estimates, with the exception of the delicate top-order estimates
that are based on modified quantities (which we define in Sect.~\ref{S:MODQUANTS}).
The main result is Lemma~\ref{L:SCHEMATICDEPENDENCEOFMANYTENSORFIELDS}.
We start by introducing some convenient shorthand notation
that we use throughout the rest of the article.

\begin{definition}[\textbf{Shorthand notation for the unknowns}]
\label{D:ABBREIVATEDVARIABLES}
We define the following arrays $\GdVar$ and $\BadVar$ of scalar functions:
\begin{subequations}
\begin{align}
	\GdVar 
	& := \left(\Psi, \Lunit_{(Small)}^1, \Lunit_{(Small)}^2 \right),
			\label{E:GOODABBREIVATEDVARIABLES} \\
	\BadVar 
	& := \left(\Psi, \upmu - 1, \Lunit_{(Small)}^1, \Lunit_{(Small)}^2 \right).
	\label{E:BADABBREIVATEDVARIABLES}
\end{align}
\end{subequations}
\end{definition}

\begin{remark}[\textbf{Schematic functional dependence}]
\label{R:SCHEMATICTENSORFIELDPRODUCTS}
In the remainder of the article, we use the notation
$\smoothfunction(\xi_{(1)},\xi_{(2)},\cdots,\xi_{(m)})$ to schematically depict
an expression (often tensorial and involving contractions)
that depends smoothly on the $\ell_{t,u}-$tangent tensorfields $\xi_{(1)}, \xi_{(2)}, \cdots, \xi_{(m)}$.
Note that in general, $\smoothfunction(0) \neq 0$.
\end{remark}

\begin{lemma}[\textbf{Schematic structure of various tensorfields}]
	\label{L:SCHEMATICDEPENDENCEOFMANYTENSORFIELDS}
	We have the following schematic relations for scalar functions:
	\begin{subequations}
	\begin{align} \label{E:SCALARSDEPENDINGONGOODVARIABLES}
		g_{\alpha \beta},
		(g^{-1})^{\alpha \beta},
		(\gt^{-1})^{\alpha \beta},
		\gsphere_{\alpha \beta},
		(\ginversesphere)^{\alpha \beta},
		G_{\alpha \beta},
		G_{\alpha \beta}',
		\Lineproject_{\beta}^{\ \alpha},
		\Lunit^{\alpha}, 
		\Radunit^{\alpha},
		\GeoAng^{\alpha}
		& = \smoothfunction(\GdVar),
			\\
		G_{\Lunit \Lunit},
		G_{\Lunit \Radunit},
		G_{\Radunit \Radunit},
		G_{\Lunit \Lunit}',
		G_{\Lunit \Radunit}',
		G_{\Radunit \Radunit}'
		& = \smoothfunction(\GdVar),
			\label{E:GFRAMESCALARSDEPENDINGONGOODVARIABLES} \\
		g_{\alpha \beta}^{(Small)},
		\GeoAng_{(Small)}^{\alpha},
		\Radunit_{(Small)}^{\alpha},
		\GeoAngFlatRadComponent
		& = \smoothfunction(\GdVar) \GdVar,
			\label{E:LINEARLYSMALLSCALARSDEPENDINGONGOODVARIABLES}
			\\
		\Rad^{\alpha} 
		& = \smoothfunction(\BadVar).
			\label{E:SCALARSDEPENDINGONBADVARIABLES}
	\end{align}
	\end{subequations}

	Moreover, we have the following schematic relations for $\ell_{t,u}-$tangent tensorfields:
	\begin{subequations}
	\begin{align}
		\gsphere,
		\angGarg{\Lunit},
		\angGarg{\Radunit},
		\angG,
		\angGprimearg{\Lunit},
		\angGprimearg{\Radunit},
		\angGprime
		& = \smoothfunction(\GdVar,\angdiff x^1,\angdiff x^2),
			\label{E:TENSORSDEPENDINGONGOODVARIABLES} \\
	\GeoAng
	& = \smoothfunction(\GdVar,\ginversesphere,\angdiff x^1,\angdiff x^2),
			\label{E:TENSORSDEPENDINGONGOODVARIABLESANDGINVERSESPHERE}
			\\
	\upzeta^{(Tan-\Psi)},
	\angk^{(Tan-\Psi)} 
	& = \smoothfunction(\GdVar,\angdiff x^1,\angdiff x^2) \Singletan \Psi,
		\label{E:TENSORSDEPENDINGONGOODVARIABLESGOODPSIDERIVATIVES}
			\\
	\upzeta^{(Trans-\Psi)},
	\angk^{(Trans-\Psi)}
	& = \smoothfunction(\GdVar,\angdiff x^1,\angdiff x^2) \Rad \Psi,
		\label{E:TENSORSDEPENDINGONGOODVARIABLESBADDERIVATIVES} \\
	\upchi 
	& = \smoothfunction(\GdVar,\angdiff x^1,\angdiff x^2) \Singletan \GdVar,
		\label{E:TENSORSDEPENDINGONGOODVARIABLESGOODDERIVATIVES}
			\\
	\mytr \upchi 
	& = \smoothfunction(\GdVar,\ginversesphere,\angdiff x^1,\angdiff x^2) \Singletan \GdVar.
	\label{E:TENSORSDEPENDINGONGOODVARIABLESGOODDERIVATIVESANDGINVERSESPHERE}
  \end{align}
\end{subequations}

\end{lemma}

\begin{remark}[\textbf{Clarification regarding the dependence of $\smoothfunction$ on $\ginversesphere$}]
	On the RHS of
	\eqref{E:TENSORSDEPENDINGONGOODVARIABLESANDGINVERSESPHERE}
	and
	\eqref{E:TENSORSDEPENDINGONGOODVARIABLESGOODDERIVATIVESANDGINVERSESPHERE},
	we view $\ginversesphere$ as a type $\binom{2}{0}$ 
	$\ell_{t,u}-$tangent tensorfield. 
	In contrast, on LHS~\eqref{E:SCALARSDEPENDINGONGOODVARIABLES},
	we are viewing the \emph{rectangular components}
	$(\ginversesphere)^{\alpha \beta}$ to be scalar functions.
	Therefore, it is not redundant to include the dependence of
	$\smoothfunction$ on $\ginversesphere$ in the
	relations
	\eqref{E:TENSORSDEPENDINGONGOODVARIABLESANDGINVERSESPHERE}
	and
	\eqref{E:TENSORSDEPENDINGONGOODVARIABLESGOODDERIVATIVESANDGINVERSESPHERE}.
\end{remark}

\begin{proof}
	The relations 
	in \eqref{E:SCALARSDEPENDINGONGOODVARIABLES}-\eqref{E:SCALARSDEPENDINGONBADVARIABLES}
	all follow easily from the definitions of the quantities involved,
	so we prove only one representative relation.
	Specifically, to obtain the schematic form 
	of $\GeoAng_{(Small)}^{\alpha}$ in
	\eqref{E:LINEARLYSMALLSCALARSDEPENDINGONGOODVARIABLES},
	we use 
	\eqref{E:GEOANGSMALLINTERMSOFRADUNIT},
	\eqref{E:FLATYDERIVATIVERADIALCOMPONENT},
	Remark~\ref{R:RADUNITSMALLLUNITSMALLRELATION},
	and the fact that $\GeoAng_{(Small)}^0 = 0$.

	The relations in
	\eqref{E:TENSORSDEPENDINGONGOODVARIABLES}-\eqref{E:TENSORSDEPENDINGONGOODVARIABLESGOODDERIVATIVESANDGINVERSESPHERE}
	are also easy to derive from the definitions of the quantities involved and some simple observations.
	We give proofs of a few representative examples.
	To obtain \eqref{E:TENSORSDEPENDINGONGOODVARIABLESANDGINVERSESPHERE}, we
	let $\GeoAng_{\flat}$ be the $\gsphere$ dual of $\GeoAng$
	so that $\GeoAng = \ginversesphere \cdot \GeoAng_{\flat}$.
	It is easy to see that we have the following identity for $\ell_{t,u}-$tangent one-forms:
	$\GeoAng_{\flat} = \GeoAng_a \angdiff x^a$;
	it can be checked by contracting both sides against elements of $\lbrace \Lunit, \Radunit, \CoordAng \rbrace$.
	We now note that $\GeoAng_a = g_{ab} \GeoAng^b$. Thus, by \eqref{E:SCALARSDEPENDINGONGOODVARIABLES},
	we have $\GeoAng_a = \smoothfunction(\GdVar)$.
	Combining the above observations, we find that
	$\GeoAng_{\flat} = \smoothfunction(\GdVar,\angdiff x^1,\angdiff x^2)$, 
	from which the desired relation \eqref{E:TENSORSDEPENDINGONGOODVARIABLESANDGINVERSESPHERE} easily follows.
	To obtain 
	\eqref{E:TENSORSDEPENDINGONGOODVARIABLESGOODDERIVATIVES}-\eqref{E:TENSORSDEPENDINGONGOODVARIABLESGOODDERIVATIVESANDGINVERSESPHERE},
	we apply similar reasoning based on the identities 
	\eqref{E:CHIINTERMSOFOTHERVARIABLES}-\eqref{E:TRCHIINTERMSOFOTHERVARIABLES}.
\end{proof}

\section{Energy Identities and Basic Ingredients in the \texorpdfstring{$L^2$}{Square Integral} Analysis}
\label{S:ENERGY}
In this section, we establish the integral identities that we use in our $L^2$ analysis. 

\subsection{Fundamental Energy Identity}
\label{SS:FUNDAMENTALENERGYIDENTITY}
To derive energy estimates, we rely on the \emph{energy-momentum tensor} $\enmomtensor$,
which is the symmetric type $\binom{0}{2}$ tensor 
\begin{align} \label{E:ENERGYMOMENTUMTENSOR}
	\enmomtensor_{\mu \nu}
	=
	\enmomtensor_{\mu \nu}[\Psi]
	& := \D_{\mu} \Psi \D_{\nu} \Psi
	- \frac{1}{2} g_{\mu \nu} (g^{-1})^{\alpha \beta} \D_{\alpha} \Psi \D_{\beta} \Psi.
\end{align}

In the next lemma, we exhibit the basic divergence property of $\enmomtensor$;
we omit the proof, which is a simple calculation.

\begin{lemma}[\textbf{Basic divergence property of} $\enmomtensor$]
\label{L:DIVT}
	For solutions to $\upmu \square_g \Psi = \mathfrak{F}$, we have
	\begin{align} \label{E:DIVT}
		\upmu \D_{\alpha} \enmomtensor^{\alpha \nu}
		& = \mathfrak{F} \D^{\nu} \Psi.
	\end{align}
\end{lemma}
\hfill $\qed$

In the next lemma, we provide the components of
$\enmomtensor$ relative to the rescaled frame.
\begin{lemma} [\textbf{The frame components of} $\enmomtensor$]
\label{L:ENMOMEMFRAMECOMPONENTS}
The components of the energy-momentum tensor $\enmomtensor$
relative to the rescaled frame can be expressed as follows:
\begin{subequations}
	\begin{align}
		\enmomtensor_{\Lunit \Lunit}[\Psi] 
		& = (\Lunit \Psi)^2, 
			\qquad
		\enmomtensor_{\Lunit \Rad}[\Psi] 
			= - \frac{1}{2} \upmu (\Lunit \Psi)^2
				+ \frac{1}{2} \upmu |\angdiff \Psi|^2,
			\\
		\enmomtensor_{\Rad \Rad}[\Psi]
		& = \frac{1}{2} \upmu^2 (\Lunit \Psi)^2
			+ (\Rad \Psi)^2
			+ \upmu (\Lunit \Psi) \Rad \Psi
			- \frac{1}{2} \upmu^2 |\angdiff \Psi|^2, 
			\\
		\angENMOMarg{\Lunit}[\Psi] 
		& = (\Lunit \Psi) \angdiff \Psi, 
			\qquad
		\angENMOMarg{\Rad}[\Psi]  
		= (\Rad \Psi) \angdiff \Psi, 
			\\
		\angENMOM[\Psi]
		 & = 
			\frac{1}{2} 
			(\Lunit \Psi)^2
			\gsphere
			+ 
			\upmu^{-1}
			(\Lunit \Psi)
			(\Rad \Psi)
			\gsphere
			+ 
			\frac{1}{2} 
			|\angdiff \Psi|^2 
			\gsphere.
\end{align}
\end{subequations}
\end{lemma}

\begin{proof}
	The lemma is a simple consequence of the formula \eqref{E:ENERGYMOMENTUMTENSOR}
	and the frame decompositions of $g$ and $g^{-1}$
	provided by
	\eqref{E:METRICFRAMEDECOMPLUNITRADUNITFRAME}
	and \eqref{E:GINVERSEFRAMEWITHRECTCOORDINATESFORGSPHEREINVERSE}.
\end{proof}

We derive our energy estimates with the help of 
the following multiplier vectorfield.

\begin{definition}[\textbf{The timelike multiplier vectorfield} $\Mult$]
\label{D:DEFINITIONMULT} 
We define (see Footnote~\ref{FN:MULTNOTATION} on pg.~\pageref{FN:MULTNOTATION} regarding the notation)
\begin{align} \label{E:DEFINITIONMULT} 
		\Mult 
		& := (1 + 2 \upmu) \Lunit + 2 \Rad.
\end{align}
\end{definition}
A simple calculation yields that
$g(\Mult,\Mult) = - 4 \upmu (1 + \upmu)$. 
Thus, $\Mult$ is $g-$timelike whenever $\upmu > 0$.
This property is important because it leads
to coercive energy identities.

In the next lemma, we provide the frame components
of $\deform{\Mult}$. These are important for 
our energy estimates because $\deform{\Mult}$
appears in our fundamental energy-flux identity
(see Prop.~\ref{P:DIVTHMWITHCANCELLATIONS}).

\begin{lemma}[\textbf{The frame components of} $\deform{\Mult}$]
\label{L:MULTIPLIERDEFORM}
The components of the deformation tensor $\deform{\Mult}$
(see Def.~\ref{D:DEFORMATIONTENSOR})
of the multiplier vectorfield \eqref{E:DEFINITIONMULT} 
can be expressed as follows relative to the rescaled frame:
\begin{subequations}
\begin{align}
	\deformarg{\Mult}{\Lunit}{\Lunit} 
	& = 0, 
		\\
	\deformarg{\Mult}{\Lunit}{\Rad} 
	& = - \left\lbrace
				\Lunit \upmu 
		 		+ 4 \upmu \Lunit \upmu
		 		+ 2 \Rad \upmu
		 	\right\rbrace, 
		\\
	\deformarg{\Mult}{\Rad}{\Radunit}
	& = 2(1 + 2 \upmu) \Lunit \upmu,
		 \\
	\angdeformoneformarg{\Mult}{\Lunit} 
	& = - 2 \angdiff \upmu 
			- 4 
				\left\lbrace
					\upzeta^{(Trans-\Psi)}
					+
					\upmu 
					\upzeta^{(Tan-\Psi)}
				\right\rbrace,
		\\
	\angdeformoneformarg{\Mult}{\Rad} 
	& = \angdiff \upmu 
		+ 2 (1 + 2 \upmu)
				\left\lbrace
					\upzeta^{(Trans-\Psi)}
					+
					\upmu 
					\upzeta^{(Tan-\Psi)}
				\right\rbrace,
			\\
	\angdeform{\Mult}
	& = 2 \upchi
			+ 4 
				\left\lbrace
					\angk^{(Trans-\Psi)}
					+
					\upmu 
					\angk^{(Tan-\Psi)}
				\right\rbrace.
\end{align}
\end{subequations}
The $\ell_{t,u}$ 
tensorfields
$\upchi$,
$\upzeta^{(Trans-\Psi)}$,
$\angk^{(Trans-\Psi)}$,
$\upzeta^{(Tan-\Psi)}$,
and
$\angk^{(Tan-\Psi)}$
from above 
are as in
\eqref{E:CHIDEF},
\eqref{E:ZETATRANSVERSAL},
\eqref{E:KABTRANSVERSAL},
\eqref{E:ZETAGOOD},
and \eqref{E:KABGOOD}.
\end{lemma}
\begin{proof}
	The proof is similar to that of
	Lemma~\ref{L:DEFORMATIONTENSORFRAMECOMPONENTS}
	but is much simpler because cancellations do not play a role;
	we therefore omit the details.
\end{proof}

We define our geometric integrals in terms of length, area, and volume
forms that remain non-degenerate throughout the evolution,
all the way up to the shock.

\begin{definition}[\textbf{Non-degenerate forms and related integrals}]
	\label{D:NONDEGENERATEVOLUMEFORMS}
	We define the length form
	$d \spherevol$ on $\ell_{t,u}$,
	the area form $d \tvol$ on $\Sigma_t^u$,
	the area form $d \conevol$ on $\mathcal{P}_u^t$,
	and the volume form $d \vol$ on $\mathcal{M}_{t,u}$
	as follows (relative to the geometric coordinates):
	\begin{align} \label{E:RESCALEDVOLUMEFORMS}
			d \spherevol
			& = d \argspherevol{(t,u,\vartheta)}
			:= \gtancomp(t,u,\vartheta) \, d \vartheta,
				&&
			d \tvol
			=
			d \tvol(t,u',\vartheta)
			:= d \spherevol(t,u',\vartheta) du',
				\\
			d \conevol 
			& = d \conevol(t',u,\vartheta)
			:= d \spherevol(t',u,\vartheta) dt',
				&&
			d \vol 
			= d \vol(t',u',\vartheta')
			:= d \spherevol(t',u',\vartheta') du' dt',
				\notag
	\end{align}
	where $\gtancomp$ is the scalar function from Def.~\ref{D:METRICANGULARCOMPONENT}.

	If $f$ is a scalar function, then we define
	\begin{subequations}
	\begin{align}
	\int_{\ell_{t,u}}
			f
		\, d \spherevol
		& 
		:=
		\int_{\vartheta \in \mathbb{T}}
			f(t,u,\vartheta)
		\, \gtancomp(t,u,\vartheta) d \vartheta,
			\label{E:LINEINTEGRALDEF} \\
		\int_{\Sigma_t^u}
			f
		\, d \tvol
		& 
		:=
		\int_{u'=0}^u
		\int_{\vartheta \in \mathbb{T}}
			f(t,u',\vartheta)
		\, \gtancomp(t,u',\vartheta) d \vartheta du',
			\label{E:SIGMATUINTEGRALDEF} \\
		\int_{\mathcal{P}_u^t}
			f
		\, d \conevol
		& 
		:=
		\int_{t'=0}^t
		\int_{\vartheta \in \mathbb{T}}
			f(t',u,\vartheta)
		\, \gtancomp(t',u,\vartheta) d \vartheta dt',
			\label{E:PUTINTEGRALDEF} \\
		\int_{\mathcal{M}_{t,u}}
			f
		\, d \vol
		& 
		:=
		\int_{t'=0}^t
		\int_{u'=0}^u
		\int_{\vartheta \in \mathbb{T}}
			f(t',u',\vartheta)
		\, \gtancomp(t',u',\vartheta) d \vartheta du' dt'.
		\label{E:MTUTUINTEGRALDEF}
	\end{align}
	\end{subequations}
\end{definition}

\begin{remark}
	The canonical forms associated to
	$\gt$ and $g$ are respectively
	$\upmu d \tvol$ and $\upmu d \vol$.
\end{remark}

We now define energies and null fluxes,
which serve as building blocks for 
the quantities that we use in our $L^2$
analysis of solutions.

\begin{definition}[\textbf{Energy and null flux}]
\label{D:ENERGYFLUX}
In terms of the non-degenerate forms of Def.~\ref{D:NONDEGENERATEVOLUMEFORMS},
we define the energy functional $\enzero[\cdot]$ and null flux functional
$\flzero[\cdot]$ as follows:
\begin{align} \label{E:ENERGYFLUX}
	\enzero[\Psi](t,u)
		& := 
		\int_{\Sigma_t^u} 
			\upmu \enmomtensor_{\Timenormal \Mult}[\Psi]
		\, d \tvol,
	\qquad
	\flzero[\Psi](t,u)
	:= 
		\int_{\mathcal{P}_u^t} 
			\enmomtensor_{\Lunit \Mult}[\Psi]
		\, d \conevol,
\end{align}
where $\Timenormal$ and $\Mult$
are the vectorfields defined in
\eqref{E:TIMENORMAL}
and
\eqref{E:DEFINITIONMULT}.

\end{definition}

In the next lemma, we reveal the coercive nature
of $\enzero[\Psi]$ and $\flzero[\Psi]$.

\begin{lemma}[\textbf{Coercivity of the energy and null flux}]
\label{L:ORDERZEROCOERCIVENESS}
The energy and null flux from Def.~\ref{D:ENERGYFLUX} enjoy the following 
coerciveness properties:
\begin{subequations}
\begin{align} \label{E:ENERGYORDERZEROCOERCIVENESS}
	\enzero[\Psi](t,u)
		& = 
		\int_{\Sigma_t^u} 
			\frac{1}{2} (1 + 2 \upmu) \upmu (\Lunit \Psi)^2
			+ 2 \upmu (\Lunit \Psi) \Rad \Psi
			+ 2 (\Rad \Psi)^2
			+ \frac{1}{2} (1 + 2 \upmu)\upmu |\angdiff \Psi|^2
		\, d \tvol,
		\\
	\flzero[\Psi](t,u)
		& = 
		\int_{\mathcal{P}_u^t} 
			(1 + \upmu)(\Lunit \Psi)^2
			+ \upmu |\angdiff \Psi|^2
		\, d \conevol.
		\label{E:NULLFLUXENERGYORDERZEROCOERCIVENESS}
\end{align}
\end{subequations}
\end{lemma}

\begin{proof}
The lemma follows from the identities
\begin{align}
	\upmu \enmomtensor_{\Timenormal \Mult}
	& = \frac{1}{2} (1 + \upmu) \upmu (\Lunit \Psi)^2
			+ \frac{1}{2} (\upmu \Lunit \Psi + 2 \Rad \Psi)^2
			+ \frac{1}{2} (1 + 2 \upmu) \upmu |\angdiff \Psi|^2,
			\\
	\enmomtensor_{\Lunit \Mult}
	& =  (1 + \upmu)(\Lunit \Psi)^2
		+ \upmu |\angdiff \Psi|^2,
\end{align}
which are a simple consequence of Lemma~\ref{L:ENMOMEMFRAMECOMPONENTS}
and the identities 
\eqref{E:TIMENORMAL}
and \eqref{E:DEFINITIONMULT}.
\end{proof}

In the next proposition, we provide
the fundamental energy-flux
identities that hold for solutions
to the inhomogeneous wave equation
$\upmu \square_{g(\Psi)} \Psi = \waveinhom$.
The term $\waveinhom$ 
represents the error terms that arise upon
commuting the homogeneous equation \eqref{E:GEOWAVE} after it has been multiplied by the factor $\upmu$ 
(see Remark~\ref{R:MUFACTOR} below for an explanation of why we include the factor $\upmu$).
See Figure~\ref{F:SOLIDREGION} on pg.~\pageref{F:SOLIDREGION} for a picture of the 
spacetime region $\mathcal{M}_{t,u}$
on which we apply the divergence theorem
and the relevant boundary surfaces
$\Sigma_0^u$,
$\Sigma_t^u$,
$\mathcal{P}_0^t$,
and
$\mathcal{P}_u^t$.

\begin{remark}
	In Figure~\ref{F:SOLIDREGION}, 
	the (unlabeled) front and back boundaries should be identified;
	they represent
	the same ``periodic timelike surface'' $\lbrace \vartheta = \mbox{\upshape const} \rbrace$
	and thus they do not make a contribution to the energy identity 
	of Prop.~\ref{P:DIVTHMWITHCANCELLATIONS}.
\end{remark}


\begin{proposition}[\textbf{Fundamental energy-flux identity}]
	\label{P:DIVTHMWITHCANCELLATIONS}
		For solutions $\Psi$ to 
		\begin{align*}
			\upmu \square_{g(\Psi)} \Psi & = \waveinhom
		\end{align*}
		that vanish along the outer null hyperplane $\mathcal{P}_0$,
		we have the following identity involving
		the energy and flux from Def.~\ref{D:ENERGYFLUX}:
	\begin{align} \label{E:E0DIVID}
				&
				\enzero[\Psi](t,u)
				+
				\flzero[\Psi](t,u)
					\\
				& 
				= 
				\enzero[\Psi](0,u)
				- 
				\int_{\mathcal{M}_{t,u}}
					\left\lbrace
						(1 + 2 \upmu) (\Lunit \Psi)
						+ 
						2 \Rad \Psi 
					\right\rbrace
				\waveinhom 
				\, d \vol
				- 
				\frac{1}{2} 
				\int_{\mathcal{M}_{t,u}}
					\upmu \enmomtensor^{\alpha \beta}[\Psi] \deformarg{\Mult}{\alpha}{\beta}
				\, d \vol.
				\notag
			\end{align}
Furthermore, with $f_+: = \max \lbrace f,0 \rbrace$ and $f_- := \max \lbrace -f, 0 \rbrace$, we have
\begin{align}
\basicenergyerror{\Mult}[\Psi] 
	& := - \frac{1}{2} \upmu \enmomtensor^{\alpha \beta}[\Psi] \deformarg{\Mult}{\alpha}{\beta}
	 = - \frac{1}{2} \upmu |\angdiff \Psi|^2 \frac{[\Lunit \upmu]_-}{\upmu} 
	 	+ \sum_{i=1}^5 \basicenergyerrorarg{\Mult}{i}[\Psi],
		\label{E:MULTERRORINT} 
\end{align}
where
	\begin{subequations}
		\begin{align}
			\basicenergyerrorarg{\Mult}{1}[\Psi] 
				& := (\Lunit \Psi)^2 
						\left\lbrace
							- \frac{1}{2} \Lunit \upmu
							+ \Rad \upmu
							- \frac{1}{2} \upmu \mytr \upchi
							- \mytr \angk^{(Trans-\Psi)}
							- \upmu \mytr \angk^{(Tan-\Psi)}
						\right\rbrace,
						\label{E:MULTERRORINTEG1} \\
			\basicenergyerrorarg{\Mult}{2}[\Psi] 
			& := - (\Lunit \Psi) (\Rad \Psi)
					\left\lbrace
						 \mytr \upchi
						+ 2 \mytr \angk^{(Trans-\Psi)}
						+ 2 \upmu \mytr \angk^{(Tan-\Psi)}
					\right\rbrace,
						\label{E:MULTERRORINTEG2} \\
		  \basicenergyerrorarg{\Mult}{3}[\Psi] 
			& := 
				\upmu |\angdiff \Psi|^2
				\left\lbrace
					\frac{1}{2} \frac{[\Lunit \upmu]_+}{\upmu}
					+ \frac{\Rad \upmu}{\upmu}
					+ 2 \Lunit \upmu
					- \frac{1}{2} \mytr \upchi
					- \mytr \angk^{(Trans-\Psi)}
					- \upmu \mytr \angk^{(Tan-\Psi)}
				\right\rbrace,
				\label{E:MULTERRORINTEG3} \\
			\basicenergyerrorarg{\Mult}{4}[\Psi] 
			& := 	(\Lunit \Psi)(\angdiffuparg{\#} \Psi) 
					\cdot
					\left\lbrace
						(1 - 2 \upmu) \angdiff \upmu 
						+ 
						2 \upzeta^{(Trans-\Psi)}
						+
						2 \upmu \upzeta^{(Tan-\Psi)}
					\right\rbrace,
					\label{E:MULTERRORINTEG4} \\
			\basicenergyerrorarg{\Mult}{5}[\Psi] 
			& := - 2 (\Rad \Psi)(\angdiffuparg{\#} \Psi)
					\cdot
					 \left\lbrace
							\angdiff \upmu 
						 	+ 
						 	2 \upzeta^{(Trans-\Psi)}
							+
							2 \upmu \upzeta^{(Tan-\Psi)}
					 \right\rbrace.
					\label{E:MULTERRORINTEG5} 
		\end{align}
\end{subequations}
	The tensorfields
	$\upchi$, 
	$\upzeta^{(Trans-\Psi)}$,
	$\angk^{(Trans-\Psi)}$,
	$\upzeta^{(Tan-\Psi)}$,
	and
	$\angk^{(Tan-\Psi)}$
	from above 
	are as in
	\eqref{E:CHIDEF},
	\eqref{E:ZETATRANSVERSAL},
	\eqref{E:KABTRANSVERSAL},
	\eqref{E:ZETAGOOD},
	and \eqref{E:KABGOOD}.
\end{proposition}

\begin{proof}
	We define the vectorfield $J^{\alpha} := \enmomtensor^{\alpha \beta}[\Psi]\Mult_{\beta}$, 
	where
	$\Mult$ is defined in \eqref{E:DEFINITIONMULT}.
	We decompose $J = J^t \frac{\partial}{\partial t} + J^u \frac{\partial}{\partial u} + J^{\CoordAng} \CoordAng$,
	where $J^t$, $J^u$, $J^{\CoordAng}$ are scalar functions
	and we recall that $\CoordAng = \frac{\partial}{\partial \vartheta}$.
	We claim that
	\begin{align} \label{E:ENERGYESTIMATEMULTIPLIERUCOMPONENT}
		J^u
		& = - \upmu^{-1} \enmomtensor_{\Lunit \Mult}[\Psi],
			\\
		J^t
		& = \upmu J^u - J_{\Radunit}
			= - \enmomtensor_{\Lunit \Mult}[\Psi]
			- \enmomtensor_{\Radunit \Mult}[\Psi]
			= - \enmomtensor_{\Timenormal \Mult}[\Psi],
			\label{E:ENERGYESTIMATEMULTIPLIERTCOMPONENT}
	\end{align}
	where $\Timenormal$ is the vectorfield defined in \eqref{E:TIMENORMAL}.
	To derive 
	\eqref{E:ENERGYESTIMATEMULTIPLIERUCOMPONENT}, we take the inner product of the decomposition equation
	with $\Lunit$ and use \eqref{E:RADSPLITINTOPARTTILAUANDXI}
	to find that 
	$g(J,\Lunit) = J^u g(\Lunit, \frac{\partial}{\partial u}) = J^u g(\Lunit, \frac{\partial}{\partial u} - \XiCoordComp \CoordAng)
	= J^u g(\Lunit, \Rad) = - \upmu J^u$. 
	Since $g(J,\Lunit) = \enmomtensor_{\Lunit \Mult}[\Psi]$,
	we have obtained the desired identity \eqref{E:ENERGYESTIMATEMULTIPLIERUCOMPONENT}. 
	The proof of \eqref{E:ENERGYESTIMATEMULTIPLIERTCOMPONENT} is
	similar and we omit it.
	Next, we note the identity
	\begin{align} \label{E:FIRSTDIVERGENCEID}
		\int_{\mathcal{M}_{t,u}}
			\upmu \D_{\alpha} J^{\alpha}
		 \, d \vol
		 & = 
		 \int_{t'=0}^t
		 \int_{u'=0}^u
		 \int_{\vartheta \in \mathbb{T}}
		 	\frac{\partial}{\partial t}
		 	\left(
		 		\upmu \gtancomp J^t
		 	\right)
		 	+
		 	\frac{\partial}{\partial u}
		 	\left(
		 		\upmu \gtancomp J^u
		 	\right)
		 	+
		 	\frac{\partial}{\partial \vartheta}
		 	\left(
		 		\upmu \gtancomp J^{\CoordAng}
		 	\right)
		 \, dt' 
		 \, du'
		 \, d \vartheta.
	\end{align}
	\eqref{E:FIRSTDIVERGENCEID} follows from 
	the standard identity for the divergence of a vectorfield
	expressed relative to a coordinate frame 
	(in this case the geometric coordinates)
	and the formula \eqref{E:SPACETIMEVOLUMEFORMWITHUPMU}, which
	implies that $|\mbox{\upshape{det}} g|^{1/2} = \upmu \gtancomp$
	(where the determinant is taken relative to the geometric coordinates).
	Using Fubini's theorem, carrying out some integrations in \eqref{E:FIRSTDIVERGENCEID},
	and noting that the integral of 
	$\frac{\partial}{\partial \vartheta}
		 	\left(
		 		\upmu \gtancomp J^{\CoordAng}
	\right)$
	over $\mathbb{T}$ vanishes, we deduce 
	\begin{align} \label{E:FIRSTDIVERGENCEBOUNDARYERMS}
		 \mbox{RHS } \eqref{E:FIRSTDIVERGENCEID}
		 & = 
		 \int_{u'=0}^u
		 \int_{\vartheta \in \mathbb{T}}
		 	\left(
		 		\upmu \gtancomp J^t
		 	\right)(t,u',\vartheta)
		 	\, du'
		 \, d \vartheta
		 -
		 \int_{u'=0}^u
		 \int_{\vartheta \in \mathbb{T}}
		 	\left(
		 		\upmu \gtancomp J^t
		 	\right)(0,u',\vartheta)
		 	\, du'
		 \, d \vartheta
		 	\\
		& \ \
		+
		 \int_{t'=0}^t
		 \int_{\vartheta \in \mathbb{T}}
		 	\left(
		 		\upmu \gtancomp J^u
		 	\right)(t',u,\vartheta)
		 	\, dt'
		 \, d \vartheta
		 -
		 \int_{t'=0}^t
		 \int_{\vartheta \in \mathbb{T}}
		 	\left(
		 		\upmu \gtancomp J^u
		 	\right)(t',0,\vartheta)
		 	\, dt'
		 \, d \vartheta.
		 \notag
	\end{align}
	Inserting \eqref{E:ENERGYESTIMATEMULTIPLIERUCOMPONENT} and \eqref{E:ENERGYESTIMATEMULTIPLIERTCOMPONENT} into 
	\eqref{E:FIRSTDIVERGENCEBOUNDARYERMS}, we obtain all terms in
	\eqref{E:E0DIVID} except for the two $\mathcal{M}_{t,u}$ integrals on the RHS.
	The proof of \eqref{E:E0DIVID} will be complete once we show that the integrands
	under the $\mathcal{M}_{t,u}$ integrals
	sum to $\upmu \D_{\alpha}(\enmomtensor^{\alpha \beta}[\Psi]\Mult_{\beta})$.
	This fact follows from
	\eqref{E:DIVT} and the symmetry of $\enmomtensor$,
	which imply that
	$\upmu \D_{\alpha}(\enmomtensor^{\alpha \beta}[\Psi]\Mult_{\beta})
	=
	\frac{1}{2}
	\upmu \enmomtensor^{\alpha \beta} \deformarg{\Mult}{\alpha}{\beta}
	+
	\upmu (\Mult \Psi) \waveinhom
	$.

	It remains for us to derive \eqref{E:MULTERRORINT}.
	We first write
	$\enmomtensor^{\alpha \beta}[\Psi] \deformarg{\Mult}{\alpha}{\beta}
		= (g^{-1})^{\alpha \beta} (g^{-1})^{\kappa \lambda} \enmomtensor_{\alpha \kappa}[\Psi] \deformarg{\Mult}{\beta}{\lambda}
	$. We then decompose the two $g^{-1}$ factors relative to the frame $\lbrace \Lunit, \Radunit, \CoordAng \rbrace$
	with the formula \eqref{E:GINVERSEFRAMEWITHRECTCOORDINATESFORGSPHEREINVERSE}.
	Also using Lemmas~\ref{L:ENMOMEMFRAMECOMPONENTS} and \ref{L:MULTIPLIERDEFORM},
	we conclude \eqref{E:MULTERRORINT} from straightforward calculations.

\end{proof}

To close our top-order energy estimates,
we must perform some additional
integrations by parts, going
beyond those of Prop.~\ref{P:DIVTHMWITHCANCELLATIONS}.
We provide the required identities in the next lemma.

\begin{lemma}[\textbf{Identities connected to integration by parts}]
	\label{L:LDERIVATIVEOFLINEINTEGRAL}
	The following identities hold for scalar functions $f$:
	\begin{subequations}
	\begin{align} \label{E:LDERIVATIVEOFLINEINTEGRAL}
		\frac{\partial}{\partial t}
		\int_{\ell_{t,u}}
			f
		\, d \spherevol
		& 
		= 
		\int_{\ell_{t,u}}
			\Lunit f
			+ \mytr \upchi f
		\, d \spherevol,
			\\
		\frac{\partial}{\partial u}
		\int_{\ell_{t,u}}
			f
		\, d \spherevol
		& 
		= 
		\int_{\ell_{t,u}}
			\Rad f
			+ 
			\frac{1}{2} \mytr \angdeform{\Rad}
			f
		\, d \spherevol,
			\label{E:UDERIVATIVEOFLINEINTEGRAL} \\
		\frac{\partial}{\partial t}
		\int_{\Sigma_t^u}
			f
		\, d \tvol
		& 
		= 
		\int_{\Sigma_t^u}
			\Lunit f
			+ \mytr \upchi f
		\, d \tvol.
		\label{E:LDERIVATIVEOFSIGMATINTEGRAL}
	\end{align}
	\end{subequations}

	Moreover, we have the following integration by parts identities:
	\begin{align} \label{E:LINEIBP}
		\int_{\ell_{t,u}}
			(\GeoAng f_1) f_2
		\, d \spherevol
		& = - 
		\int_{\ell_{t,u}}
			f_1 (\GeoAng f_2)
		\, d \spherevol
		- \int_{\ell_{t,u}} 
				\overbrace{\angdiv Y}^{\frac{1}{2} \mytr \angdeform{Y}} f_1 f_2
			\, d \spherevol,
	\end{align}

	\begin{align}	\label{E:LUNITIBPIDENTITY}
		\int_{\mathcal{M}_{t,u}}
			(\Lunit f_1) f_2
		\, d \vol
		& 
		= 
		-
		\int_{\mathcal{M}_{t,u}}
			f_1 (\Lunit f_2)
		\, d \vol
		-
		\int_{\mathcal{M}_{t,u}}
			\mytr \upchi f_1 f_2
		\, d \vol
			\\
	& \ \ 
		+ 
		\int_{\Sigma_t^u}
			f_1 f_2
		\, d \tvol
		-
		\int_{\Sigma_0^u}
			f_1 f_2
		\, d \tvol.
		\notag
	\end{align}

	Finally, the following integration by parts identity holds for 
	scalar functions $\ThirdSmoothFunction$:
	\begin{align}	\label{E:LUNITANDANGULARIBPIDENTITY}
		& 
		\int_{\mathcal{M}_{t,u}}
			(1 + 2 \upmu) (\Rad \Psi) (\Lunit \Tanset^N \Psi) \GeoAng \ThirdSmoothFunction
		\, d \vol
			\\
	& = 
		\int_{\mathcal{M}_{t,u}}
			(1 + 2 \upmu) (\Rad \Psi) (\GeoAng \Tanset^N \Psi) \Lunit \ThirdSmoothFunction
		\, d \vol
		\notag	\\
	& \ \ 
		- \int_{\Sigma_t^u}
				(1 + 2 \upmu) (\Rad \Psi) (\GeoAng \Tanset^N \Psi) \ThirdSmoothFunction
			\, d \vol
		+ \int_{\Sigma_0^u}
				(1 + 2 \upmu) (\Rad \Psi) (\GeoAng \Tanset^N \Psi) \ThirdSmoothFunction
			\, d \vol
			\notag \\
	& \ \
		+
		\int_{\mathcal{M}_{t,u}}
			\mbox{\upshape Error}_1[\Tanset^N \Psi;\ThirdSmoothFunction]
		\, d \vol
		+ \int_{\Sigma_t^u}
				\mbox{\upshape Error}_2[\Tanset^N \Psi;\ThirdSmoothFunction]
			\, d \vol
		- \int_{\Sigma_0^u}
				\mbox{\upshape Error}_2[\Tanset^N \Psi;\ThirdSmoothFunction]
			\, d \vol,
			\notag
	\end{align}
	where
	\begin{subequations}
	\begin{align} \label{E:LUNITIBPSPACETIMEERROR}
		\mbox{\upshape Error}_1[\Tanset^N \Psi;\ThirdSmoothFunction] 
		& := 2 (\Lunit \upmu) (\Rad \Psi) (\GeoAng \Tanset^N \Psi) \ThirdSmoothFunction
			+ (1 + 2 \upmu) (\Lunit \Rad \Psi) (\GeoAng \Tanset^N \Psi) \ThirdSmoothFunction
				\\
		& \ \
			+ (1 + 2 \upmu) (\Rad \Psi) (\angdeformoneformupsharparg{\GeoAng}{\Lunit} \cdot \angdiff \Tanset^N \Psi) \ThirdSmoothFunction
			+  (1 + 2 \upmu) (\Rad \Psi) \mytr \upchi (\GeoAng \Tanset^N \Psi) \ThirdSmoothFunction
			\notag \\
		& \ \
			+ 2 (\GeoAng \upmu) (\Rad \Psi) (\Tanset^N \Psi) \Lunit \ThirdSmoothFunction
			+ (1 + 2 \upmu) (\GeoAng \Rad \Psi) (\Tanset^N \Psi) \Lunit \ThirdSmoothFunction
				\notag \\
		& \ \ + \frac{1}{2} (1 + 2 \upmu) (\Rad \Psi) \mytr \angdeform{Y} (\Tanset^N \Psi) \Lunit \ThirdSmoothFunction
				\notag \\
		& \ \ 
			+  2 (\Lunit \GeoAng \upmu) (\Rad \Psi) (\Tanset^N \Psi) \ThirdSmoothFunction
			+  2 (\GeoAng \upmu) (\Lunit \Rad \Psi) (\Tanset^N \Psi) \ThirdSmoothFunction
				\notag \\
		& \ \
			+  2 (\GeoAng \upmu) (\Rad \Psi) \mytr \upchi (\Tanset^N \Psi) \ThirdSmoothFunction
			+  (\Lunit \upmu) (\Rad \Psi) \mytr \angdeform{Y} (\Tanset^N \Psi) \ThirdSmoothFunction
				\notag \\
		& \ \
			+  (1 + 2 \upmu) (\Lunit \GeoAng \Rad \Psi) (\Tanset^N \Psi) \ThirdSmoothFunction
			+  (1 + 2 \upmu) (\Rad \Psi) (\GeoAng \mytr \upchi) (\Tanset^N \Psi) \ThirdSmoothFunction
			\notag \\
		& \ \
		   +  (1 + 2 \upmu) (\GeoAng \Rad \Psi) \mytr \upchi (\Tanset^N \Psi) \ThirdSmoothFunction
			+ \frac{1}{2} (1 + 2 \upmu) (\Lunit \Rad \Psi) \mytr \angdeform{Y} (\Tanset^N \Psi) \ThirdSmoothFunction 
			 \notag \\
		& \ \
			+ (1 + 2 \upmu) (\Rad \Psi) (\angdiv \angdeformoneformupsharparg{\GeoAng}{\Lunit}) (\Tanset^N \Psi) \ThirdSmoothFunction
			+ \frac{1}{2} (1 + 2 \upmu) (\Rad \Psi) \mytr \upchi \mytr \angdeform{Y} (\Tanset^N \Psi) \ThirdSmoothFunction,
				\notag
				\\
		\mbox{\upshape Error}_2[\Tanset^N \Psi;\ThirdSmoothFunction] 
		& := 
			- 2 (\GeoAng \upmu) (\Rad \Psi) (\Tanset^N \Psi) \ThirdSmoothFunction
			- (1 + 2 \upmu) (\GeoAng \Rad \Psi) (\Tanset^N \Psi) \ThirdSmoothFunction
				\label{E:LUNITIBPHYPERSURFACEERROR} \\
		& \ \
			- \frac{1}{2} (1 + 2 \upmu) (\Rad \Psi) \mytr \angdeform{Y} (\Tanset^N \Psi) \ThirdSmoothFunction.
			\notag	
	\end{align}
	\end{subequations}

\end{lemma}

\begin{proof}
	To prove \eqref{E:UDERIVATIVEOFLINEINTEGRAL},
	we first fix $t$ and construct
	a local coordinate $\widetilde{\vartheta}$ on $\Sigma_t^{U_0} = \cup_{u \in [0,U_0]} \ell_{t,u}$
	by setting $\widetilde{\vartheta} = \vartheta$ on $\ell_{t,0}$
	and then propagating $\widetilde{\vartheta}$ 
	by solving the transport equation $\Rad \widetilde{\vartheta} = 0$.
	Since $\Rad u = 1$, it follows 
	that relative to the coordinates $(u,\widetilde{\vartheta})$ on $\Sigma_t^{U_0}$, we have
	$[\Rad, \frac{\partial}{\partial \widetilde{\vartheta}}] = 0$
	and $\Rad = \frac{\partial}{\partial u}$.
	Relative to $(u,\widetilde{\vartheta})$ coordinates on $\Sigma_t^{U_0}$, we have
	$
	\ginversesphere = \widetilde{\gtancomp}^{-2} 
	\frac{\partial}{\partial \widetilde{\vartheta}} 
	\otimes
	\frac{\partial}{\partial \widetilde{\vartheta}}
	$
	and $d \spherevol = \widetilde{\gtancomp} d \, \widetilde{\vartheta}$,
	where 
	$\widetilde{\gtancomp}^2 
		=  g(\frac{\partial}{\partial \widetilde{\vartheta}},\frac{\partial}{\partial \widetilde{\vartheta}})$
	(see 
	\eqref{E:METRICANGULARCOMPONENT},
	\eqref{E:GSPHEREINVERSERELATIVETOGEOMETRIC},
	and
	\eqref{E:RESCALEDVOLUMEFORMS}).
	Differentiating with $\angLie_{\Rad}$ and using
	\eqref{E:CONNECTIONBETWEENANGLIEOFGSPHEREANDDEFORMATIONTENSORS},
	we find that
	$
	- \angdeform{\Rad}^{\# \#}
	=
	\angLie_{\Rad} \ginversesphere 
	= -2 \widetilde{\gtancomp}^{-3} \Rad \widetilde{\gtancomp}
		\frac{\partial}{\partial \widetilde{\vartheta}} 
		\otimes
		\frac{\partial}{\partial \widetilde{\vartheta}}
	$.
	Contracting against $\gsphere$, we obtain
	the identity $\Rad \ln \widetilde{\gtancomp} = (1/2) \mytr \angdeform{\Rad}$.
	We now express the integrand on LHS~\eqref{E:UDERIVATIVEOFLINEINTEGRAL} in $(u,\widetilde{\vartheta})$ coordinates,
	differentiate under the integral, and use that $\Rad = \frac{\partial}{\partial u}$ in these coordinates
	to obtain
	$
	\frac{\partial}{\partial u}
		\int_{\ell_{t,u}}
			f
		\, d \spherevol
		= 
		\int_{\ell_{t,u}}
			\left\lbrace
				\Rad f
				+ 
				(\Rad \ln \widetilde{\gtancomp}) f
			\right\rbrace
		\, \widetilde{\gtancomp} d \widetilde{\vartheta}
	$. 
	Using this identity and the previous expression for $\Rad \ln \widetilde{\gtancomp}$,
	we conclude \eqref{E:UDERIVATIVEOFLINEINTEGRAL}.

	The proof of \eqref{E:LDERIVATIVEOFLINEINTEGRAL} is similar
	and relies on the identity $\Lunit \ln \gtancomp = \mytr \upchi$
	(see \eqref{E:LDERIVATIVEOFVOLUMEFORMFACTOR});
	we omit the details.

	\eqref{E:LDERIVATIVEOFSIGMATINTEGRAL} follows from \eqref{E:LDERIVATIVEOFLINEINTEGRAL}
	and the fact that $\int_{\Sigma_t^u} \cdots \, d \tvol = \int_{u'=0}^u \int_{\ell_{t,u'}} \cdots \, d \spherevol du'$.

	\eqref{E:LINEIBP} follows easily from integrating the identity
	$
	(\GeoAng f_1) f_2
	+
	(\GeoAng f_2) f_1
	=
	\GeoAng(f_1 f_2) 
	= 
	\angdiv(f_1 f_2 \GeoAng) - f_1 f_2 \angdiv Y$ over $\ell_{t,u}
	$.

	\eqref{E:LUNITIBPIDENTITY} follows from integrating the identity
	\eqref{E:LDERIVATIVEOFSIGMATINTEGRAL} with $f = f_1 f_2$ with respect to time
	from time $0$ to time $t$.

	The identity \eqref{E:LUNITANDANGULARIBPIDENTITY} follows from
	a series of tedious but straightforward integrations by parts
	that we now describe. We first integrate by parts using \eqref{E:LUNITIBPIDENTITY}
	in order to move the $\Lunit$ operator off of $\Lunit \Tanset^N \Psi$.
	This procedure results in the presence of the integral
	$
		-
		\int_{\mathcal{M}_{t,u}}
			(1 + 2 \upmu) (\Rad \Psi) (\Tanset^N \Psi) \Lunit \GeoAng \ThirdSmoothFunction
	\, d \vol
	$
	(among others). 
	We then commute $\Lunit$ and $\GeoAng$ to obtain the identity
	$\Lunit \GeoAng \ThirdSmoothFunction 
	= \GeoAng \Lunit \ThirdSmoothFunction 
	+ \angdeformoneformupsharparg{\GeoAng}{\Lunit} \cdot \angdiff \ThirdSmoothFunction
	$ 
	(see \eqref{E:CONNECTIONBETWEENCOMMUTATORSANDDEFORMATIONTENSORS}),
	which we substitute into the previous integral.
	We then use \eqref{E:LINEIBP} to move all $\GeoAng$ derivatives off of all factors of
	$\ThirdSmoothFunction$ in all of the error integrals.
	Finally, we integrate by parts on the $\ell_{t,u}$ to move all
	$\angdiff$ derivatives off of all factors of
	$\ThirdSmoothFunction$ in all of the error integrals.
\end{proof}

\section{The Structure of the Terms in the Commuted Wave Equation}
\label{S:STRUCTUREOFTERMSINCOMMUTEDWAVE}
\setcounter{equation}{0}
To derive energy estimates for the higher derivatives of $\Psi$,
we commute the wave equation $\upmu \square_{g(\Psi)} \Psi = 0$
with vectorfields $Z \in \Fullset$.
In this section, 
we reveal the precise structure of the
commutator error terms $[Z, \upmu \square_{g(\Psi)}] \Psi$.
By precise structure, we mean that we decompose
all terms relative to the rescaled frame 
$\lbrace \Lunit, \Rad, \CoordAng \rbrace$
and keep track of the exact expressions including the
constant coefficients; some of these
constants affect the number of derivatives we need to
close our estimates.

\begin{remark} \label{R:MUFACTOR}
We have included the factor of $\upmu$ in front of the operator $\square_{g(\Psi)}$ in the wave equation
because it leads to important cancellations in the commutation identities.
\end{remark}

We start with the following standard commutation identity.

\begin{lemma}[\cite{dCsK1993}*{Lemma 7.1.3}]
	\label{L:BASICLIECOVARIANTZCOMMUTATION}
	If $\xi_{\alpha_1 \cdots \alpha_n}$ is any type $\binom{0}{n}$ spacetime tensorfield, 
	$V$ is any spacetime vectorfield, and $\deform{V}$ is its deformation tensor (see Def.~\ref{D:DEFORMATIONTENSOR}),
	then
	\begin{align} \label{E:BASICLIECOVARIANTZCOMMUTATION}
		\D_{\beta} \Lie_V \xi_{\alpha_1 \cdots \alpha_n}
		-
		\Lie_V  \D_{\beta} \xi_{\alpha_1 \cdots \alpha_n}
		& = 
			\frac{1}{2}
				\sum_{i=1}^n
				\left\lbrace
					\D_{\alpha_i} \deformmixedarg{V}{\beta}{\kappa} 
					+
					\D_{\beta} \deformmixedarg{V}{\alpha_i} {\kappa}
					-
					\D^{\kappa} \deformarg{V}{\alpha_i}{\beta} 
			\right\rbrace
			\xi_{\alpha_1 \cdots \alpha_{i-1} \kappa \alpha_{i+1} \cdots \alpha_n}.
	\end{align}
\end{lemma}

We now use Lemma~\ref{L:BASICLIECOVARIANTZCOMMUTATION} to derive an identity for
$[\upmu \square_g, Z]$.

\begin{lemma}[\textbf{Vectorfield-covariant wave operator commutation lemma}]
\label{L:BOXZCOM}
For the vectorfields $Z \in \Fullset$ (see Def.~\ref{D:COMMUTATIONVECTORFIELDS}),
which by Lemma~\ref{L:DEFORMATIONTENSORFRAMECOMPONENTS} satisfy 
$\deformarg{Z}{\Lunit}{\Lunit} = 0$
and $\deformarg{Z}{\Lunit}{\Rad} = -Z \upmu$,
we have the following commutation identity:
\begin{align} \label{E:BOXZCOM}
		\upmu \square_{g(\Psi)} (Z \Psi)
		& =   \upmu 
					\D_{\alpha} 
					\left\lbrace
					\deformuparg{Z}{\alpha}{\beta} \D_{\beta} \Psi
						- \frac{1}{2} \myspacetimetr \deform{Z} \D^{\alpha} \Psi 
					\right\rbrace
				+ Z (\upmu \square_{g(\Psi)} \Psi)
				+ \frac{1}{2} \mytr \angdeform{Z} (\upmu \square_{g(\Psi)} \Psi),
\end{align}
where $\myspacetimetr \deform{Z} := (g^{-1})^{\alpha \beta} \deformarg{Z}{\alpha}{\beta}$.
\end{lemma}

\begin{proof}
	We begin by applying $Z$ to $\upmu \square_{g(\Psi)} \Psi$ and using the Leibniz rule for Lie
	derivatives and the Lie derivative identity 
	$(\Lie_Z g^{-1})^{\alpha \beta} = - \deformuparg{Z}{\alpha}{\beta}$
	to obtain the identity
	$Z(\upmu \square_{g(\Psi)} \Psi) 
	= (Z \upmu) \square_{g(\Psi)} \Psi
		- \upmu\deformuparg{Z}{\alpha}{\beta} \D_{\alpha} \D_{\beta} \Psi
		+ \upmu (g^{-1})^{\alpha \beta} \Lie_Z \D_{\alpha} \D_{\beta} \Psi
	$.
	Applying \eqref{E:BASICLIECOVARIANTZCOMMUTATION} with $\xi := \D \Psi$, we obtain
  $
		\D_{\alpha} \D_{\beta} Z \Psi
		= \Lie_Z \D_{\alpha} \D_{\beta} \Psi
			+ \frac{1}{2}
				\left\lbrace
					\D_{\alpha} \deformarg{Z}{\lambda}{\beta} \D^{\lambda} \Psi
					+
					\D_{\beta} \deformarg{Z}{\alpha}{\lambda} \D^{\lambda} \Psi
					-
					\D_{\lambda} \deformarg{Z}{\alpha}{\beta} \D^{\lambda} \Psi
				\right\rbrace
	$.
	Contracting the previous identity against $(g^{-1})^{\alpha \beta}$ 
	and using the Leibniz rule for Lie
	derivatives and the aforementioned identity
	$(\Lie_Z g^{-1})^{\alpha \beta} = - \deformuparg{Z}{\alpha}{\beta}$,
	we find that
	\begin{align} \label{E:BOXOFZPSIFIRSTEXPRESSION}
		\square_{g(\Psi)} Z \Psi
		& = Z(\square_{g(\Psi)} \Psi)
			+ 
			\deformuparg{Z}{\alpha}{\beta} \D_{\alpha} \D_{\beta} \Psi
			+
			(\D_{\alpha} \deformuparg{Z}{\alpha}{\lambda}) \D_{\lambda} \Psi
			- \frac{1}{2} 
				 \D_{\lambda}
				 \left\lbrace
					(g^{-1})^{\alpha \beta}
					\deformarg{Z}{\alpha}{\beta} 
				\right\rbrace
				\D^{\lambda} \Psi
					\\
			& = Z(\square_{g(\Psi)} \Psi)
			+ 
			\D_{\alpha}
			\left\lbrace
				\deformuparg{Z}{\alpha}{\beta}  \D_{\beta} \Psi
				-
				\frac{1}{2}
				\myspacetimetr \deform{Z} 
				\D^{\alpha} \Psi
			\right\rbrace
			+
			\frac{1}{2} 
			\myspacetimetr \deform{Z} 
			\square_{g(\Psi)} \Psi.
			\notag
	\end{align}
	The identity \eqref{E:BOXZCOM} now follows easily from
	\eqref{E:BOXOFZPSIFIRSTEXPRESSION},
	the identity
	$\frac{1}{2} \myspacetimetr \deform{Z}
	= 
	- \frac{1}{2} \deformarg{Z}{\Lunit}{\Lunit}
	- \upmu^{-1} \deformarg{Z}{\Lunit}{\Rad}
	+ \frac{1}{2} \mytr \angdeform{Z}
	$
	(see \eqref{E:GINVERSEFRAMEWITHRECTCOORDINATESFORGSPHEREINVERSE}), 
	and the assumed properties
	$\deformarg{Z}{\Lunit}{\Lunit} = 0$ and $\deformarg{Z}{\Lunit}{\Rad} = -Z \upmu$.

\end{proof}

In the next proposition, we decompose the first term
on RHS~\eqref{E:BOXZCOM} relative to the rescaled frame.
Our proof of the proposition relies on the following lemma.

\begin{lemma}[\textbf{Spacetime divergence in terms of derivatives of frame components}]
	\label{L:DIVERGENCEFRAME}
	Let $\mathscr{J}$ be a spacetime vectorfield.
	Let 
	$\upmu \mathscr{J} 
	= 
	- \upmu \mathscr{J}_{\Lunit} \Lunit 
	- \mathscr{J}_{\Rad} \Lunit 
	- \mathscr{J}_{\Lunit} \Rad 
	+ \upmu \angJ$
	be its decomposition relative to the rescaled frame, where
	$\mathscr{J}_{\Lunit} = \mathscr{J}^{\alpha} \Lunit_{\alpha}$,
	$\mathscr{J}_{\Rad} = \mathscr{J}^{\alpha} \Rad_{\alpha}$,
	and $\angJ = \Lineproject \mathscr{J}$. Then
	\begin{align} \label{E:DIVERGENCEFRAME}
	\upmu \D_{\alpha} \mathscr{J}^{\alpha} 
	& = - \Lunit (\upmu \mathscr{J}_{\Lunit})
		- \Lunit (\mathscr{J}_{\Rad})
		- \Rad (\mathscr{J}_{\Lunit})
		+ \angdiv (\upmu \angJ)
		- \upmu \mytr \angk \mathscr{J}_{\Lunit}
		- \mytr \upchi \mathscr{J}_{\Rad},
\end{align}
where the $\ell_{t,u}-$tangent tensorfields
$\upchi$ and
$\angk$ can be expressed via 
\eqref{E:CHIINTERMSOFOTHERVARIABLES}
and
\eqref{E:ANGKDECOMPOSED}.
\end{lemma}

\begin{proof}
Using \eqref{E:GINVERSEFRAMEWITHRECTCOORDINATESFORGSPHEREINVERSE}, 
we find that
\begin{align} \label{E:PRELIMINARDIVERGENCEFRAME}
\upmu \D_{\alpha} \mathscr{J}^{\alpha} 
& = \upmu (g^{-1})^{\alpha \beta} \D_{\alpha} \mathscr{J}_{\beta}
	\\
& = 	- \Lunit(\upmu \mathscr{J}_{\Lunit})
			- \Lunit(\mathscr{J}_{\Rad}) 
			- \Rad(\mathscr{J}_{\Lunit}) 
			+ (\ginversesphere)^{\alpha \beta} \D_{\alpha} (\upmu \mathscr{J}_{\beta})
				\notag \\
& \ \ + (\Lunit \upmu) \mathscr{J}_{\Lunit}
			+ \upmu (\D_{\Lunit} \Lunit^{\alpha}) \mathscr{J}_{\alpha}
			+ (\D_{\Lunit} \Rad^{\alpha}) \mathscr{J}_{\alpha}
			+ (\D_{\Rad} \Lunit^{\alpha}) \mathscr{J}_{\alpha}
			- \angJ \cdot \angdiff \upmu.
			\notag
\end{align}
Next, we use 
Lemma~\ref{L:CONNECTIONLRADFRAME}
to substitute for 
$\D_{\Lunit} \Lunit$,
$\D_{\Lunit} \Rad$,
and $\D_{\Rad} \Lunit$;
we find that all terms on the last line of
\eqref{E:PRELIMINARDIVERGENCEFRAME} cancel.
We then use the rescaled frame decomposition formula to express
$
(\ginversesphere)^{\alpha \beta} \D_{\alpha} (\upmu \mathscr{J}_{\beta})
= \angdiv (\upmu \angJ)
- \upmu \mathscr{J}_{\Lunit} (\ginversesphere)^{\alpha \beta} \D_{\alpha} \Lunit_{\beta}
- \mathscr{J}_{\Rad} (\ginversesphere)^{\alpha \beta} \D_{\alpha} \Lunit_{\beta}
- \mathscr{J}_{\Lunit} (\ginversesphere)^{\alpha \beta} \D_{\alpha} \Rad_{\beta} 
$
and then Lemma~\ref{L:CONNECTIONLRADFRAME} to deduce the following
identities, which we substitute into the previous equation:
$(\ginversesphere)^{\alpha \beta} \D_{\alpha} \Lunit_{\beta} = \mytr \upchi$
and $(\ginversesphere)^{\alpha \beta} \D_{\alpha} \Rad_{\beta} = 
\upmu \mytr \angk - \upmu \mytr \upchi$.
Straightforward calculations then lead to \eqref{E:DIVERGENCEFRAME}.
\end{proof}

We now decompose the term
$\upmu 
					\D_{\alpha} 
					\left\lbrace\deformuparg{Z}{\alpha}{\beta} \D_{\beta} \Psi
						- \frac{1}{2} \myspacetimetr \deform{Z} \D^{\alpha} \Psi 
					\right\rbrace
$
from RHS~\eqref{E:BOXZCOM} relative to the rescaled frame.

\begin{proposition} [\textbf{Frame decomposition of the divergence of the key inhomogeneous term}]
\label{P:COMMUTATIONCURRENTDIVERGENCEFRAMEDECOMP}
	For vectorfields $Z \in \Fullset$,
	which have $\deformarg{Z}{\Lunit}{\Lunit} = 0$
	and $\deformarg{Z}{\Lunit}{\Rad} = -Z \upmu$,
	we have the following identity for the
	first term on RHS~\eqref{E:BOXZCOM}:
	\begin{align} \label{E:DIVCOMMUTATIONCURRENTDECOMPOSITION}
				\upmu 
					\D_{\alpha} 
					\left\lbrace
						\deformuparg{Z}{\alpha}{\beta} \D_{\beta} \Psi
						- \frac{1}{2} \myspacetimetr \deform{Z} \D^{\alpha} \Psi 
					\right\rbrace
		& = \mathscr{K}_{(\pi-Danger)}^{(Z)}[\Psi]
			\\
		& \ \
			+ \mathscr{K}_{(\pi-Cancel-1)}^{(Z)}[\Psi]
			+ \mathscr{K}_{(\pi-Cancel-2)}^{(Z)}[\Psi]
			\notag	\\
		& \ \ 
			+ \mathscr{K}_{(\pi-Less \ Dangerous)}^{(Z)}[\Psi]
			+ \mathscr{K}_{(\pi-Good)}^{(Z)}[\Psi]
				\notag \\
		& \	\
			+ \mathscr{K}_{(\Psi)}^{(Z)}[\Psi]
			+ \mathscr{K}_{(Low)}^{(Z)}[\Psi],
			\notag
	\end{align}
	where
	\begin{subequations}
		\begin{align}
			\mathscr{K}_{(\pi-Danger)}^{(Z)}[\Psi]
			& := - (\angdiv \angdeformoneformupsharparg{Z}{\Lunit}) \Rad \Psi,
				\label{E:DIVCURRENTTRANSVERSAL}
				\\
			\mathscr{K}_{(\pi-Cancel-1)}^{(Z)}[\Psi]
			& := \left\lbrace 
						\frac{1}{2} \Rad \mytr  \angdeform{Z}
						- \angdiv \angdeformoneformupsharparg{Z}{\Rad}
						- \upmu \angdiv \angdeformoneformupsharparg{Z}{\Lunit}
					\right\rbrace 
						\Lunit \Psi,
					\label{E:DIVCURRENTCANEL1} \\
			\mathscr{K}_{(\pi-Cancel-2)}^{(Z)}[\Psi]
			& :=
				\left\lbrace
					- \angLie_{\Rad} \angdeformoneformupsharparg{Z}{\Lunit}
					+ \angdiffuparg{\#} \deformarg{Z}{\Lunit}{\Rad}
				\right\rbrace 
				\cdot
				\angdiff \Psi,
				\label{E:DIVCURRENTCANEL2} \\
			\mathscr{K}_{(\pi-Less \ Dangerous)}^{(Z)}[\Psi]
			& := \frac{1}{2} \upmu (\angdiffuparg{\#} \mytr \angdeform{Z}) \cdot \angdiff \Psi, 
				\label{E:DIVCURRENTELLIPTIC} \\
			\mathscr{K}_{(\pi-Good)}^{(Z)}[\Psi] 
			& := \frac{1}{2} \upmu (\Lunit \mytr \angdeform{Z}) \Lunit \Psi
				+ (\Lunit \deformarg{Z}{\Lunit}{\Rad}) \Lunit \Psi
				+ (\Lunit \deformarg{Z}{\Rad}{\Radunit}) \Lunit \Psi
				\label{E:DIVCURRENTGOOD} \\
			& \ \ + \frac{1}{2} (\Lunit \mytr \angdeform{Z}) \Rad \Psi
				- \upmu (\angLie_{\Lunit} \angdeformoneformupsharparg{Z}{\Lunit}) \cdot \angdiff \Psi
				- (\angLie_{\Lunit} \angdeformoneformupsharparg{Z}{\Rad}) \cdot \angdiff \Psi,
				\notag 
	\end{align}
	\end{subequations}
	\begin{align}
		\mathscr{K}_{(\Psi)}^{(Z)}[\Psi]
			& := \left\lbrace
							\frac{1}{2}
							\upmu \mytr \angdeform{Z}
							+ \deformarg{Z}{\Lunit}{\Rad}
							+ \deformarg{Z}{\Rad}{\Radunit}
						\right\rbrace
						\Lunit^2 \Psi  
					\label{E:DIVCURRENTPSI} \\
		& \ \
				+ \mytr \angdeform{Z}
					 \Lunit \Rad \Psi
					\notag \\
			& \ \ - 2 \upmu \angdeformoneformupsharparg{Z}{\Lunit} \cdot \angdiff \Lunit \Psi
					- 2 \angdeformoneformupsharparg{Z}{\Rad} \cdot \angdiff \Lunit \Psi
					- 2 \angdeformoneformupsharparg{Z}{\Lunit} \cdot \angdiff \Rad \Psi
				\notag \\
			& \ \ + \deformarg{Z}{\Lunit}{\Rad} \angLap \Psi
				+ \frac{1}{2} \upmu \mytr \angdeform{Z} \angLap \Psi,
				\notag 
	\end{align}
	and
	\begin{align}
		\mathscr{K}_{(Low)}^{(Z)}[\Psi] 
		& := \left\lbrace
				 	\frac{1}{2} (\Lunit \upmu) \mytr \angdeform{Z}
				 	+ \frac{1}{2} \upmu \mytr \angk \mytr \angdeform{Z}
				 	+ \mytr \upchi \deformarg{Z}{\Lunit}{\Rad}
				 	+ \mytr \upchi \deformarg{Z}{\Rad}{\Radunit}
				 	- \angdeformoneformupsharparg{Z}{\Lunit} \cdot \angdiff \upmu
				 \right\rbrace
				 \Lunit \Psi
				\label{E:DIVCURRENTLOW}  \\
		& \ \ 
				+
				\frac{1}{2} \mytr \upchi \mytr \angdeform{Z} \Rad \Psi
				\notag \\
		& \ \  
				+ \left\lbrace
						- (\Lunit \upmu) \angdeformoneformupsharparg{Z}{\Lunit}
						- \upmu \mytr \angk \angdeformoneformupsharparg{Z}{\Lunit}
						- \mytr \upchi \angdeformoneformupsharparg{Z}{\Rad}
						+ \mytr \angdeform{Z} \angdiffuparg{\#} \upmu
						+  \mytr \upchi \upmu \upzeta^{\#}
					\right\rbrace
					\cdot
					\angdiff \Psi.
				\notag
	\end{align}
	In the above expressions, the $\ell_{t,u}-$tangent tensorfields
	$\upchi$, 
	$\upzeta$,
	and
	$\angk$, 
	are as in
	\eqref{E:CHIINTERMSOFOTHERVARIABLES},
	\eqref{E:ZETADECOMPOSED},
	and 
	\eqref{E:ANGKDECOMPOSED}.
\end{proposition}

\begin{proof}
	We define $\mathscr{J}$
	to be the spacetime vectorfield whose
	divergence is taken on LHS~\eqref{E:DIVCOMMUTATIONCURRENTDECOMPOSITION}:
	$\mathscr{J}^{\alpha}
	:= \deformuparg{Z}{\alpha}{\beta} \D_{\beta} \Psi
			- 
			\frac{1}{2} (g^{-1})^{\kappa \lambda}\deformarg{Z}{\kappa}{\lambda} \D^{\alpha} \Psi 
	$.
	With the help of Lemma~\ref{L:METRICDECOMPOSEDRELATIVETOTHEUNITFRAME}, 
	we compute that
	\begin{align} \label{E:COMMUTATIONJLFRAMEDECOMPOSED}
		\mathscr{J}_{\Lunit}
		& = - \frac{1}{2} \mytr \angdeform{Z} \Lunit \Psi
			+ \angdeformoneformupsharparg{Z}{\Lunit} \cdot \angdiff \Psi, 
				\\
		\mathscr{J}_{\Rad}
		& = - \deformarg{Z}{\Lunit}{\Rad} \Lunit \Psi 
			- \deformarg{Z}{\Rad}{\Radunit} \Lunit \Psi 
			+ \angdeformoneformupsharparg{Z}{\Rad} \cdot \angdiff \Psi
			- \frac{1}{2} \mytr \angdeform{Z}\Rad \Psi,
			\label{E:COMMUTATIONJRADFRAMEDECOMPOSED}\\
		\upmu \angJ
		& = - \upmu \angdeformoneformupsharparg{Z}{\Lunit} \Lunit \Psi
			- \angdeformoneformupsharparg{Z}{\Rad} \Lunit \Psi
			- \angdeformoneformupsharparg{Z}{\Lunit} \Rad \Psi
			+ \deformarg{Z}{\Lunit}{\Rad} \angdiffuparg{\#} \Psi
			+  \frac{1}{2} \upmu \mytr \angdeform{Z} \angdiffuparg{\#} \Psi.
			\label{E:COMMUTATIONJANGFRAMEDECOMPOSED}
	\end{align}
	The proposition then follows from 
	the divergence formula \eqref{E:DIVERGENCEFRAME}
	and tedious but straightforward calculations.
	We remark that in our calculations, 
	we use the identity 
	$[\Lunit,\Rad] 
	= 
	- \angdiff^{\#} \upmu 
		- 2 \upmu \zeta^{\#}$
	(see \eqref{E:CONNECTIONBETWEENCOMMUTATORSANDDEFORMATIONTENSORS} 
	and \eqref{E:RADDEFORMSPHERERAD})
	to replace the term 
	$\frac{1}{2} \mytr \angdeform{Z} \Rad \Lunit \Psi$
	arising from \eqref{E:COMMUTATIONJLFRAMEDECOMPOSED}
	with 
	$\frac{1}{2} \mytr \angdeform{Z} \Lunit \Rad \Psi
		+ \frac{1}{2} \mytr \angdeform{Z} (\angdiff^{\#} \upmu) \cdot \angdiff \Psi
		+ \mytr \angdeform{Z} \upmu \zeta^{\#} \cdot \angdiff \Psi
	$.
\end{proof}

\section{Differential Operator Commutation Identities}
\label{S:DIFFERENTIALOPERATORSANDCOMMUTATIONIDENTITIES}
\setcounter{equation}{0}
In this section, 
we provide a collection of
commutation identities that we use when commuting the equations.
The precise numerical constants and the structure
of tensor contractions in these identities is
not important for our estimates.
Thus, we present some of the identities in schematic form.

\begin{definition}[\textbf{Notation for repeated differentiation}]
\label{D:REPEATEDDIFFERENTIATIONSHORTHAND}
We recall the commutation sets $\Fullset$ and $\Tanset$ from Def.~\ref{D:COMMUTATIONVECTORFIELDS}.
We label the three vectorfields in 
$\Fullset$ as follows: $Z_{(1)} = \Lunit, Z_{(2)} = \GeoAng, Z_{(3)} = \Rad$.
Note that $\Tanset = \lbrace Z_{(1)}, Z_{(2)} \rbrace$.
We define the following vectorfield operators:
\begin{itemize}
	\item If $\vec{I} = (\iota_1, \iota_2, \cdots, \iota_N)$ is a multi-index
		of order $|\vec{I}| := N$
		with $\iota_1, \iota_2, \cdots, \iota_N \in \lbrace 1,2,3 \rbrace$,
		then $\Fullset^{\vec{I}} := Z_{(\iota_1)} Z_{(\iota_2)} \cdots Z_{(\iota_N)}$ 
		denotes the corresponding $N^{th}$ order differential operator.
		We write $\Fullset^N$ rather than $\Fullset^{\vec{I}}$
		when we are not concerned with the structure of $\vec{I}$.
	\item Similarly, $\angLie_{\Fullset}^{\vec{I}} 
	:= \angLie_{Z_{(\iota_1)}} \angLie_{Z_{(\iota_2})} \cdots \angLie_{Z_{(\iota_N})}$
		denotes an $N^{th}$ order $\ell_{t,u}-$projected Lie derivative operator
		(see Def.~\ref{D:PROJECTEDLIE}),
		and we write $\angLie_{\Fullset}^N$
		when we are not concerned with the structure of $\vec{I}$.
	\item If $\vec{I} = (\iota_1, \iota_2, \cdots, \iota_N)$,
		then 
		$\vec{I}_1 + \vec{I}_2 = \vec{I}$ 
		means that
		$\vec{I}_1 = (\iota_{k_1}, \iota_{k_2}, \cdots, \iota_{k_m})$
		and
		$\vec{I}_2 = (\iota_{k_{m+1}}, \iota_{k_{m+2}}, \cdots, \iota_{k_N})$,
		where $1 \leq m \leq N$ and
		$k_1, k_2, \cdots, k_N$ is a permutation of 
		$1,2,\cdots,N$. 
	\item Sums such as $\vec{I}_1 + \vec{I}_2 + \cdots + \vec{I}_M = \vec{I}$
		have an analogous meaning.
	\item $\mathcal{P}_u-$tangent operators such as 
		$\Tanset^{\vec{I}}$ are defined analogously,  
		except in this case we clearly have
		$\iota_1, \iota_2, \cdots, \iota_N \in \lbrace 1,2 \rbrace$.
\end{itemize}
\end{definition}

\begin{remark}[\textbf{Schematic depiction of the structure of} $\Fullset^{\vec{I}}$ \textbf{and} $\Tanset^{\vec{I}}$]
In deriving our estimates,
we often need only partial information about the structure 
of the operators $\Fullset^{\vec{I}}$ and $\Tanset^{\vec{I}}$.
Thus, in Subsect.~\ref{SS:STRINGSOFCOMMUTATIONVECTORFIELDS},
we introduce additional shorthand notation that captures the information that we need.
\end{remark}

\begin{lemma} [\textbf{Preliminary identities for commuting $Z \in \Fullset$ with $\angD$}]
\label{L:COMMUTINGVEDCTORFIELDSWITHANGD}
For each $\Fullset-$multi-index $\vec{I}$ and integer $n \geq 1$,
there exist constants
$C_{\vec{I}_1,\vec{I}_2,\cdots,\vec{I}_{M+1},n}$
such that the following commutator identity holds
for all type $\binom{0}{n}$ $\ell_{t,u}-$tangent tensorfields $\xi$:
\begin{align} \label{E:ANGDANGLIEZELLTUTENSORFIELDCOMMUTATOR}
	[\angD, \angLie_{\Fullset}^{\vec{I}}] \xi
	& = 	\sum_{M=1}^{|\vec{I}|}
				\mathop{\sum_{\vec{I}_1 + \cdots + \vec{I}_{M+1} = \vec{I}}}_{|\vec{I}_a| \geq 1 \mbox{\upshape \ for } 1 \leq a \leq M} 
				C_{\vec{I}_1,\vec{I}_2,\cdots,\vec{I}_{M+1},n}
				(\ginversesphere)^M
				\underbrace{
				(\angLie_{\Fullset}^{\vec{I}_1} \gsphere)
				\cdots
				(\angLie_{\Fullset}^{\vec{I}_{M-1}} \gsphere)
				}_{\mbox{absent when $M=1$}}
				(\angD \angLie_{\Fullset}^{\vec{I}_M} \gsphere)
				(\angLie_{\Fullset}^{\vec{I}_{M+1}} \xi).
\end{align}
Moreover, with $\angdiv$ denoting the torus divergence operator from Def.~\ref{D:CONNECTIONS},
for each $\Fullset-$multi-index $\vec{I}$,
there exist constants
$C_{\vec{I}_1,\vec{I}_2,\cdots,\vec{I}_{M+1},i_1,i_2}$
such that the following commutator identity holds
for all symmetric type $\binom{0}{2}$ $\ell_{t,u}-$tangent tensorfields $\xi$:
\begin{align}
	[\angdiv, \angLie_{\Fullset}^{\vec{I}}] \xi
	& = 	\sum_{i_1 + i_2 = 1}
				\sum_{M=1}^{|\vec{I}|}
				\mathop{\sum_{\vec{I}_1 + \cdots + \vec{I}_{M+1} = \vec{I}}}_{|\vec{I}_a| \geq 1 \mbox{\upshape \ for } 1 \leq a \leq M} 
					\label{E:ANGDIVANGLIEZELLTUTENSORFIELDCOMMUTATOR} \\
	& \ \ \ \ \ \
				 C_{\vec{I}_1,\vec{I}_2,\cdots,\vec{I}_{M+1},i_1,i_2}
				(\ginversesphere)^{M+1}
				\underbrace{
				(\angLie_{\Fullset}^{\vec{I}_1} \gsphere)
				\cdots
				(\angLie_{\Fullset}^{\vec{I}_{M-1}} \gsphere)
				}_{\mbox{absent when $i_1=M=1$}}
				(\angD^{i_1} \angLie_{\Fullset}^{\vec{I}_M} \gsphere)
				(\angD^{i_2} \angLie_{\Fullset}^{\vec{I}_{M+1}} \xi).
				\notag
\end{align}

Finally, for each $\Fullset-$multi-index $\vec{I}$
and each commutation vectorfield $Z \in \Fullset$,
there exist constants
$C_{\vec{I}_1,\vec{I}_2,\cdots,\vec{I}_{M+1}}$
and
$C_{\vec{I}_1,\vec{I}_2,\cdots,\vec{I}_{M+1},i_1,i_2}$
such that the following commutator identity holds
for all scalar-valued functions $f$:
\begin{subequations}
\begin{align} \label{E:COMMUTINGANGDSQUAREDANDLIEZ}
	[\angD^2, \angLie_{\Fullset}^{\vec{I}}] f
	& = 	\sum_{M=1}^{|\vec{I}|}
				\mathop{\sum_{\vec{I}_1 + \cdots + \vec{I}_{M+1} = \vec{I}}}_{|\vec{I}_a| \geq 1 \mbox{\upshape \ for } 1 \leq a \leq M} 
				C_{\vec{I}_1,\vec{I}_2,\cdots,\vec{I}_{M+1}}
				(\ginversesphere)^M
				\underbrace{
				(\angLie_{\Fullset}^{\vec{I}_1} \gsphere)
				\cdots
				(\angLie_{\Fullset}^{\vec{I}_{M-1}} \gsphere)
				}_{\mbox{absent when $M=1$}}
				(\angD \angLie_{\Fullset}^{\vec{I}_M} \gsphere)
				(\angdiff \Fullset^{\vec{I}_{M+1}} f),
				\\
	[\angLap, \Fullset^{\vec{I}}] f
	& = 	\sum_{i_1 + i_2 = 1}
				\sum_{M=1}^{|\vec{I}|}
				\mathop{\sum_{\vec{I}_1 + \cdots + \vec{I}_{M+1} = \vec{I}}}_{|\vec{I}_a| \geq 1 \mbox{\upshape \ for } 1 \leq a \leq M} 
					\label{E:COMMUTINGANGANGLAPANDLIEZ}
					\\
	& \ \ \ \ \ \
				C_{\vec{I}_1,\vec{I}_2,\cdots,\vec{I}_{M+1},i_1,i_2}
				(\ginversesphere)^{M+1}
				\underbrace{
				(\angLie_{\Fullset}^{\vec{I}_1} \gsphere)
				\cdots
				(\angLie_{\Fullset}^{\vec{I}_{M-1}} \gsphere)
				}_{\mbox{absent when $i_1=M=1$}}
				(\angD^{i_1} \angLie_{\Fullset}^{\vec{I}_M} \gsphere)
				(\angD^{i_2+1} \Fullset^{\vec{I}_{M+1}} f).
				\notag
\end{align}
\end{subequations}
In equations \eqref{E:ANGDANGLIEZELLTUTENSORFIELDCOMMUTATOR}-\eqref{E:COMMUTINGANGANGLAPANDLIEZ}, 
we have omitted all tensorial contractions in order to condense the presentation.
\end{lemma}

\begin{proof}
	We claim that for $Z \in \Fullset$, we have the schematic identity
	$[\angD, \angLie_Z] \xi \sim (\angD \angLie_Z \gsphere)^{\#} \cdot \xi$,
	correct up to constants.
	From this identity
	and the fact that $\angLie_Z \ginversesphere = - (\angLie_Z \gsphere)^{\# \#}$
	(see \eqref{E:CONNECTIONBETWEENANGLIEOFGSPHEREANDDEFORMATIONTENSORS}),
	a straightforward argument involving induction in $|\vec{I}|$, omitted here, yields \eqref{E:ANGDANGLIEZELLTUTENSORFIELDCOMMUTATOR}.
	We now prove the claim in the case that $\xi$ is an $\ell_{t,u}-$tangent one-form. The case of higher-order tensorfields
	then follows easily from the Leibniz rules for $\angD$ and $\angLie_Z$ 
	and we omit those details.
	Moreover, it is easy to reduce the proof
	to the case $\xi_{\CoordAng} = 1$ (that is, $\xi = \angdiff \vartheta$);
	we thus assume for the remainder of the proof that $\xi = \angdiff \vartheta$.
	We now recall that $\gtancomp^2 = \gsphere(\CoordAng,\CoordAng)$ 
	(see \eqref{E:METRICANGULARCOMPONENT}).
	Note that $\angD$ is entirely determined by the formula
	$\angD_{\CoordAng} \CoordAng = \gtancomp^{-1} (\CoordAng \gtancomp) \CoordAng$.
	We first address the case $Z = \Lunit$.
	By Lemma~\ref{L:LANDRADCOMMUTEWITHANGDIFF} and the fact that $\Lunit \vartheta = 0$, 
	we have
	$\angD \angLie_{\Lunit} \xi 
	= \angD^2 (\Lunit \vartheta)
	= 0
	$.
	Next, we compute compute that
	$\angD_{\CoordAng \CoordAng}^2 \vartheta
	:=
	(\angD^2 \vartheta) \cdot \CoordAng \otimes \CoordAng
	=
	\CoordAng \CoordAng \vartheta
	- (\angD_{\CoordAng} \CoordAng) \cdot \angdiff \vartheta 
	=
	- \gtancomp^{-1} \CoordAng \gtancomp
	$.
	Also using $[\Lunit, \CoordAng] = 0$, 
	we compute that
	$\angLie_{\Lunit} \angD_{\CoordAng \CoordAng}^2 \vartheta
	:=
	(\angLie_{\Lunit} \angD^2 \vartheta) \cdot \CoordAng \otimes \CoordAng 
	= - \gtancomp^{-1} \CoordAng \Lunit \gtancomp
		+ \gtancomp^{-2} (\Lunit \gtancomp) \CoordAng \gtancomp
	$.
	Similarly, we compute that
	$(\angLie_{\Lunit} \gsphere)_{\CoordAng \CoordAng} 
		= 2 \gtancomp \Lunit \gtancomp
	$,
	$(\angD \angLie_{\Lunit} \gsphere)_{\CoordAng \CoordAng \CoordAng}
	=	2 \gtancomp \CoordAng \Lunit \gtancomp
		- 2 (\CoordAng \gtancomp) \Lunit \gtancomp
	$,
	and
	$(\angD \angLie_{\Lunit} \gsphere)^{\#} 
	\cdot (\xi \otimes \CoordAng \otimes \CoordAng)
	= 2 \gtancomp^{-1} \Lunit \CoordAng \gtancomp
		- 2 \gtancomp^{-2} (\CoordAng \gtancomp) \Lunit \gtancomp
	$.
	Combining the above computations, we conclude that
	$[\angD, \angLie_{\Lunit}] \xi
	= \frac{1}{2} (\angD \angLie_{\Lunit} \gsphere)^{\#} \cdot \xi
	$
	as desired.
	To treat the case $Z = \Rad$, we first fix $t$ and construct
	a local coordinate $\widetilde{\vartheta}$ on $\Sigma_t^{U_0}$
	such that $\Rad = \frac{\partial}{\partial u}$,
	as in the proof of Lemma~\ref{L:LDERIVATIVEOFLINEINTEGRAL}.
	The proof then mirrors the proof in the case $Z = \Lunit$. 
	We now treat the case $Z = \GeoAng$, which is 
	$\ell_{t,u}-$intrinsic. The result follows from Lemma~\ref{L:BASICLIECOVARIANTZCOMMUTATION},
	which for the $\ell_{t,u}-$tangent vectorfield $\GeoAng$ applies with 
	$\D$ replaced by $\angD$ and $g$ replaced by $\gsphere$.
	We have thus proved the claim, which completes the proof of \eqref{E:ANGDANGLIEZELLTUTENSORFIELDCOMMUTATOR}.

	\eqref{E:ANGDIVANGLIEZELLTUTENSORFIELDCOMMUTATOR} then follows
	as a straightforward consequence of \eqref{E:ANGDANGLIEZELLTUTENSORFIELDCOMMUTATOR},
	the fact that $\angdiv \xi = \ginversesphere \cdot \angD \xi$
	for symmetric type $\binom{0}{2}$ $\ell_{t,u}-$tangent tensorfields $\xi$,
	and the aforementioned identity $\angLie_Z \ginversesphere = - (\angLie_Z \gsphere)^{\# \#}$.

	To prove \eqref{E:COMMUTINGANGDSQUAREDANDLIEZ},
	we first use Lemma~\ref{L:LANDRADCOMMUTEWITHANGDIFF} to deduce that
	$
	\angD^2 \angLie_Z f 
	= \angD (\angLie_Z \angdiff f)
	= \Lie_Z \angD^2 f
		+  [\angD, \angLie_Z] \cdot \angdiff f
	$.
	The identity \eqref{E:COMMUTINGANGDSQUAREDANDLIEZ} now follows from
	\eqref{E:ANGDANGLIEZELLTUTENSORFIELDCOMMUTATOR}
	with $\xi := \angdiff f$
	and a straightforward argument involving induction in $|\vec{I}|$
	and Lemma~\ref{L:LANDRADCOMMUTEWITHANGDIFF};
	we omit the details.

	\eqref{E:COMMUTINGANGANGLAPANDLIEZ} follows easily 
	from \eqref{E:COMMUTINGANGDSQUAREDANDLIEZ}, 
	the identity $\angLap = \ginversesphere \cdot \angD^2$,
	and the aforementioned identity
	$\angLie_Z \ginversesphere 
		= 
	- (\angLie_Z \gsphere)^{\# \#}
	$.
\end{proof}

\begin{lemma}[\textbf{Preliminary Lie derivative commutation identities}]
	%

	Let $\vec{I} = (\iota_1,\iota_2,\cdots, \iota_N)$ be an $N^{th}-$order $\Fullset$ multi-index,
	let $f$ be a function, and let
	$\xi$ be a type $\binom{m}{n}$ $\ell_{t,u}-$tangent tensorfield with $m + n \geq 1$.
	Let $i_1,i_2,\cdots,i_N$ be any permutation of
	$1,2,\cdots,N$ and let
	$\vec{I}' = (\iota_{i_1},\iota_{i_2},\cdots, \iota_{i_N})$.
	Then there exist constants
	$
	C_{\vec{I}_1,\vec{I}_2,\iota_{k_1},\iota_{k_2}}
	$
	such that
	\begin{subequations}
	\begin{align} \label{E:ALGEBRAICSUBTRACTINGTWOPERMUTEDORDERNVECTORFIELDS}
		\left\lbrace
			\Fullset^{\vec{I}}
			- \Fullset^{\vec{I}'}
		\right\rbrace
		f
		& = 	\mathop{\sum_{\vec{I}_1 + \vec{I}_2 + \iota_{k_1} + \iota_{k_2} = \vec{I}}}_{
					Z_{(\iota_{k_1})} \in \lbrace \Lunit, \Rad \rbrace, 
					\ Z_{(\iota_{k_2})} \in \lbrace \Rad, \GeoAng \rbrace, 
					\ Z_{(\iota_{k_1})} \neq Z_{(\iota_{k_2})}}
					C_{\vec{I}_1,\vec{I}_2,\iota_{k_1},\iota_{k_2}}
					\angLie_{\Fullset}^{\vec{I}_1} \angdeformoneformupsharparg{Z_{(\iota_{k_2})}}{Z_{(\iota_{k_1})}} 
					\cdot 
					\angdiff \Fullset^{\vec{I}_2} f,
					\\
		\left\lbrace
			\angLie_{\Fullset}^{\vec{I}}
			- \angLie_{\Fullset}^{\vec{I}'}
		\right\rbrace
		\xi
		& = 
			 	\mathop{\sum_{\vec{I}_1 + \vec{I}_2 + \iota_{k_1} + \iota_{k_2} = \vec{I}}}_{
			 		Z_{(\iota_{k_1})} \in \lbrace \Lunit, \Rad \rbrace, 
					\ Z_{(\iota_{k_2})} \in \lbrace \Rad, \GeoAng \rbrace, 
					\ Z_{(\iota_{k_1})} \neq Z_{(\iota_{k_2})}
			 	}
				C_{\vec{I}_1,\vec{I}_2,\iota_{k_1},\iota_{k_2}}
				\angLie_{\angLie_{\Fullset}^{\vec{I}_1} \angdeformoneformupsharparg{Z_{(\iota_{k_2})}}{Z_{(\iota_{k_1})}}} 
				\angLie_{\Fullset}^{\vec{I}_2} \xi.
				\label{E:TENSORFIELDACTINGALGEBRAICSUBTRACTINGTWOPERMUTEDORDERNVECTORFIELDS}
	\end{align}
	\end{subequations}
	In \eqref{E:ALGEBRAICSUBTRACTINGTWOPERMUTEDORDERNVECTORFIELDS}-\eqref{E:TENSORFIELDACTINGALGEBRAICSUBTRACTINGTWOPERMUTEDORDERNVECTORFIELDS},
	$
	\vec{I}_1 + \vec{I}_2 + \iota_{k_1} + \iota_{k_2} = \vec{I}
	$
	means that 
	$\vec{I}_1 = (\iota_{k_3}, \iota_{k_4}, \cdots, \iota_{k_m})$,
	and
	$\vec{I}_2 = (\iota_{k_{m+1}}, \iota_{k_{m+2}}, \cdots, \iota_{k_N})$,
	where
	$k_1, k_2, \cdots, k_N$ is a permutation of 
	$1,2,\cdots,N$. In particular, 
	$|\vec{I}_1| + |\vec{I}_2| = N-2$.
	\end{lemma}

\begin{proof}
	The identities 
	\eqref{E:ALGEBRAICSUBTRACTINGTWOPERMUTEDORDERNVECTORFIELDS}
	and
	\eqref{E:TENSORFIELDACTINGALGEBRAICSUBTRACTINGTWOPERMUTEDORDERNVECTORFIELDS}
	are straightforward to verify using Lemmas~\ref{L:CONNECTIONBETWEENCOMMUTATORSANDDEFORMATIONTENSORS}
	and \ref{L:LANDRADCOMMUTEWITHANGDIFF},
	the facts that 
	$ \angdeformoneformupsharparg{W}{Z}
		= 
		- \angdeformoneformupsharparg{Z}{W}
	$ 
	for $W \in \lbrace \Rad, \GeoAng \rbrace$ 
	and $Z \in \lbrace \Lunit, \Rad \rbrace$ 
	(see Lemma~\ref{L:DEFORMATIONTENSORFRAMECOMPONENTS}),
	and the Lie derivative commutation property \eqref{E:JACOBI}.
\end{proof}

\section{Modified Quantities Needed for Top-Order Estimates}
\label{S:MODQUANTS}
As we explained in Subsubsect.~\ref{SSS:INTROENERGYESTIMATES},
in order to close our top-order energy estimates without
incurring derivative loss, we must work with modified quantities.
The modified quantities allow us to control the top-order derivatives of
$\mytr \upchi$. In this section, we define these ``fully modified quantities''
and derive transport equations for them.
At the top order, we also need ``partially modified quantities,''
which are similar but serve a different purpose: 
they enable us to avoid the appearance of certain $\Sigma_t^u$ 
error integrals in the energy estimates that are too large to be controlled
all the way up to the shock. These error integrals arise when we integrate by parts with 
respect to $\Lunit$ using the identity \eqref{E:LUNITANDANGULARIBPIDENTITY}.

\subsection{Curvature tensors and the key Ricci component identity}
\label{SS:CURVATURETENSORSRICCIID}
The calculations related to the modified quantities 
are involved. A convenient way to organize them is
to rely on the curvature tensors of $g$.

\begin{definition}[\textbf{Curvature tensors} \textbf{of} $g$]
\label{D:SPACETIMECURVATURE}
The Riemann curvature tensor $\Cur_{\alpha \beta \kappa \lambda}$ of the spacetime metric $g$ is
the type $\binom{0}{4}$ spacetime tensorfield defined by
\begin{align} \label{E:SPACETIMERIEMANN}
	g(\D_{UV}^2 W - \D_{VU}^2 W,Z)
	& = - \Cur(U,V,W,Z),
\end{align}
where $U$, $V$, $W$, and $Z$ are arbitrary spacetime vectors.
In \eqref{E:SPACETIMERIEMANN}, 
$\D_{UV}^2 W := U^{\alpha} V^{\beta} \D_{\alpha} \D_{\beta} W$.

The Ricci curvature tensor $\Ric_{\alpha \beta}$ of $g$ is the following type 
$\binom{0}{2}$ tensorfield:
\begin{align} \label{E:RICCIDEF}
	\Ric_{\alpha \beta}
	& := (g^{-1})^{\kappa \lambda} \Cur_{\alpha \kappa \beta \lambda}.
\end{align}
\end{definition}

We now provide the lemma that forms the crux
of the construction of the modified quantities.

\begin{lemma}[\textbf{The key identity verified by} $\upmu \Ric_{\Lunit \Lunit}$]
\label{L:ALPHARENORMALIZED}
Assume that $\square_{g(\Psi)} \Psi = 0$.
Then the following identity holds for the Ricci curvature
component $\Ric_{\Lunit \Lunit} := \Ric_{\alpha \beta} \Lunit^{\alpha} \Lunit^{\beta}$:
\begin{align} \label{E:ALPHARENORMALIZED}
	\upmu \Ric_{\Lunit \Lunit}
	& = \Lunit
			\left\lbrace
				- G_{\Lunit \Lunit} \Rad \Psi
				- \frac{1}{2} \upmu \mytr \angG \Lunit \Psi
				- \frac{1}{2} \upmu G_{\Lunit \Lunit} \Lunit \Psi 
				+ \upmu \angGnospacemixedarg{\Lunit}{\#} \cdot \angdiff \Psi
			\right\rbrace
		+ \mathfrak{A},
\end{align}
where $\mathfrak{A}$ has the following schematic structure,
where $Z \in \mathscr{Z}$ and $P$ is $\mathcal{P}_u^t-$tangent:
\begin{align} \label{E:RENORMALIZEDALPHAINHOMOGENEOUSTERM}
	\mathfrak{A}
	& = \smoothfunction(\BadVar,\ginversesphere,\angdiff x^1,\angdiff x^2,Z \Psi)
			\Singletan \Psi.
\end{align}

Furthermore, without assuming $\square_{g(\Psi)} \Psi = 0$,
we have
\begin{align} \label{E:ALPHAPARTIALRENORMALIZED}
	\Ric_{\Lunit \Lunit}
	& = \frac{(\Lunit \upmu)}{\upmu} \mytr \upchi
		+ \Lunit
			\left\lbrace
				- \frac{1}{2} \mytr \angG \Lunit \Psi
				- \frac{1}{2} G_{\Lunit \Lunit} \Lunit \Psi 
				+ \angGnospacemixedarg{\Lunit}{\#} \cdot \angdiff \Psi
			\right\rbrace
		- \frac{1}{2} G_{\Lunit \Lunit} \angLap \Psi
		+ \mathfrak{B},
\end{align}
where $\mathfrak{B}$ has the following schematic structure:
\begin{align} \label{E:PARTIALRENORMALIZEDALPHAINHOMOGENEOUSTERM}
	\mathfrak{B}
	& = 	\smoothfunction(\GdVar,\ginversesphere,\angdiff x^1,\angdiff x^2)
				(\Singletan \Psi) 
				\Singletan \GdVar.
\end{align}

\end{lemma}

\begin{proof}[Sketch of a proof]
	We sketch the proof instead of providing complete details since 
	the computations are lengthy and since the identities
	follow from inserting the schematic relations
	provided by Lemma~\ref{L:SCHEMATICDEPENDENCEOFMANYTENSORFIELDS}
	into the identities derived in \cite{jS2014b}*{Corollary 11.1.13}.
	We start by sketching the proof of \eqref{E:ALPHAPARTIALRENORMALIZED}.
	First, we note that straightforward but tedious computations 
	imply that relative to the \textbf{rectangular} coordinate system, 
	the components of $\Cur$ can be expressed as
	\begin{align}
		\Cur_{\mu \nu \alpha \beta}
			& = \frac{1}{2} \Big\lbrace
					G_{\beta \mu} \D_{\alpha \nu}^2 \Psi
				+ G_{\alpha \nu} \D_{\beta \mu}^2 \Psi
				- G_{\beta \nu} \D_{\alpha \mu}^2 \Psi
				- G_{\alpha \mu} \D_{\beta \nu}^2 \Psi
				\Big \rbrace 
				\label{E:RECTANGULRCURVATURECOMPONENTS} \\
		& \ \ 
				+ 
				\frac{1}{4} 
				G_{\mu \alpha}
				G_{\nu \beta} 
				(g^{-1})^{\kappa \lambda} 
				(\partial_{\kappa} \Psi)
				(\partial_{\lambda} \Psi)
				- 
				\frac{1}{4} 
				G_{\mu \beta} 
				G_{\nu \alpha} 
				(g^{-1})^{\kappa \lambda}
				(\partial_{\kappa} \Psi)
				(\partial_{\lambda} \Psi)
			\notag \\
		& \ \ 
				+ 
				\frac{1}{4} 
				(g^{-1})^{\kappa \lambda} 
				[G_{\lambda (\nu} \partial_{\alpha)} \Psi]
				[G_{\kappa (\mu} \partial_{\beta)} \Psi]
				- 
				\frac{1}{4} 
				(g^{-1})^{\kappa \lambda} 
				[G_{\lambda (\mu} \partial_{\alpha)} \Psi]
				[G_{\kappa (\nu} \partial_{\beta)} \Psi]
			\notag
				\\
			& \ \ + \frac{1}{2} \Big\lbrace
					G_{\beta \mu}' (\partial_{\alpha} \Psi) (\partial_{\nu} \Psi)
				+ G_{\alpha \nu}' (\partial_{\beta} \Psi) (\partial_{\mu} \Psi)
				- G_{\beta \nu}' (\partial_{\alpha} \Psi) (\partial_{\mu} \Psi)
				- G_{\alpha \mu}' (\partial_{\beta} \Psi) (\partial_{\nu} \Psi)
				\Big \rbrace,
				\notag 
\end{align}
where 
$\D^2 \Psi$ (which is a symmetric type $\binom{0}{2}$ tensorfield) denotes the second covariant derivative of $\Psi$.
We then contract both sides of \eqref{E:RECTANGULRCURVATURECOMPONENTS}
against $\Lunit^{\mu} \Lunit^{\alpha} (\ginversesphere)^{\nu \kappa}$,
which yields the desired term $\Ric_{\Lunit \Lunit}$ on the LHS.
Finally, we use Lemmas~\ref{L:CONNECTIONLRADFRAME} and \ref{L:SCHEMATICDEPENDENCEOFMANYTENSORFIELDS}
to express the RHS of the contracted identity 
in the form written on RHS~\eqref{E:ALPHAPARTIALRENORMALIZED}.

The main idea behind the proof of \eqref{E:ALPHARENORMALIZED}
is to multiply both sides of \eqref{E:ALPHAPARTIALRENORMALIZED} by $\upmu$
and use the wave equation in the form 
	\eqref{E:LONOUTSIDEGEOMETRICWAVEOPERATORFRAMEDECOMPOSED};
	the wave equation allows us
	to replace $\upmu$ times the term
	$
	\displaystyle
	- \frac{1}{2} G_{\Lunit \Lunit} \angLap \Psi
	$
	from \eqref{E:ALPHAPARTIALRENORMALIZED} with 
	$
	\displaystyle
	- 
	\frac{1}{2}
	\Lunit
	\left\lbrace 
		\upmu G_{\Lunit \Lunit} \Lunit \Psi + 2 G_{\Lunit \Lunit} \Rad \Psi) 
	\right\rbrace
	$
	up to error terms involving an admissible number of derivatives.
	This completes our proof sketch of the lemma.
\end{proof}

\subsection{The definitions of the modified quantities and their transport equations}
\label{SS:MODIFIEDQUANTITIES}
We now define the modified quantities.

\begin{definition}[\textbf{Modified versions of the pure} {$\mathcal{P}_u$-}\textbf{tangent derivatives of} $\mytr \upchi$]
\label{D:TRANSPORTRENORMALIZEDTRCHIJUNK}
Let $\Tanset^N$ be an $N^{th}$ order pure $\mathcal{P}_u-$tangent commutation vectorfield operator.
We define the fully modified quantity $\upchifullmodarg{\Tanset^N}$ as follows:
\begin{subequations}
\begin{align}
	\upchifullmodarg{\Tanset^N}
	& := \upmu \Tanset^N \mytr \upchi 
			 + \Tanset^N \upchifullmodinhom,
		\label{E:TRANSPORTRENORMALIZEDTRCHIJUNK} 
			\\
	\upchifullmodinhom
	& := - G_{\Lunit \Lunit} \Rad \Psi
				- \frac{1}{2} \upmu \mytr \angG \Lunit \Psi
				- \frac{1}{2} \upmu G_{\Lunit \Lunit} \Lunit \Psi 
				+ \upmu \angGnospacemixedarg{\Lunit}{\#} \cdot \angdiff \Psi.
			\label{E:LOWESTORDERTRANSPORTRENORMALIZEDTRCHIJUNKDISCREPANCY}
\end{align}
\end{subequations}

We define the partially modified quantity $\upchipartialmodarg{\Tanset^N}$ as follows:
\begin{subequations}
\begin{align}
	\upchipartialmodarg{\Tanset^N}
	& := \Tanset^N \mytr \upchi 
		+ \upchipartialmodinhomarg{\Tanset^N},
		\label{E:TRANSPORTPARTIALRENORMALIZEDTRCHIJUNK} \\
	\upchipartialmodinhomarg{\Tanset^N}
	& := - \frac{1}{2} \mytr \angG \Lunit \Tanset^N \Psi
			- \frac{1}{2} G_{\Lunit \Lunit} \Lunit \Tanset^N \Psi
			+ \angGmixedarg{\Lunit}{\#} \cdot \angdiff \Tanset^N \Psi.
			\label{E:TRANSPORTPARTIALRENORMALIZEDTRCHIJUNKDISCREPANCY}
\end{align}
\end{subequations}
We also define the following ``$0^{th}-$order'' version of \eqref{E:TRANSPORTPARTIALRENORMALIZEDTRCHIJUNKDISCREPANCY}:
\begin{align} \label{E:LOWESTORDERTRANSPORTPARTIALRENORMALIZEDTRCHIJUNKDISCREPANCY}
	\upchipartialmodinhom
	& := - \frac{1}{2} \mytr \angG \Lunit \Psi
			- \frac{1}{2} G_{\Lunit \Lunit} \Lunit \Psi
			+ \angGmixedarg{\Lunit}{\#} \cdot \angdiff \Psi.
\end{align}

\end{definition}

We now derive the transport equation
verified by the fully modified quantities.

\begin{proposition}[\textbf{The transport equation for the fully modified version of} $\Tanset^N \mytr \upchi$]
\label{P:TOPORDERTRCHIJUNKRENORMALIZEDTRANSPORT}
Assume that $\square_{g(\Psi)} \Psi = 0$.
Let $\Tanset^N$ be 
an $N^{th}-$order $\mathcal{P}_u-$tangent commutation vectorfield operator,
and let $\upchifullmodarg{\Tanset^N}$ and $\upchifullmodinhom$ be the corresponding quantities 
defined in \eqref{E:TRANSPORTRENORMALIZEDTRCHIJUNK} and \eqref{E:LOWESTORDERTRANSPORTRENORMALIZEDTRCHIJUNKDISCREPANCY}.
Then $\upchifullmodarg{\Tanset^N}$ verifies the following transport equation:
\begin{align} \label{E:TOPORDERTRCHIJUNKRENORMALIZEDTRANSPORT}
\Lunit \upchifullmodarg{\Tanset^N}
	- \left(
			2 \frac{\Lunit \upmu}{\upmu}
			- 2 \mytr \upchi
		\right)
		\upchifullmodarg{\Tanset^N}
	& = 
	\upmu [\Lunit, \Tanset^N] \mytr \upchi
		- 
				\left(
					2 \frac{\Lunit \upmu}{\upmu} 
				\right)
				\Tanset^N \upchifullmodinhom
		+ 2 \mytr \upchi \Tanset^N \upchifullmodinhom
			\\
	& \ \ 
			+ [\Lunit, \Tanset^N] \upchifullmodinhom
			+ [\upmu, \Tanset^N] \Lunit \mytr \upchi
			+ [\Tanset^N,\Lunit \upmu] \mytr \upchi
			\notag \\
	& \ \ - \left\lbrace
					\Tanset^N \left(\upmu (\mytr \upchi)^2 \right)
					- 2 \upmu \mytr \upchi \Tanset^N \mytr \upchi
				\right\rbrace 
			- \Tanset^N \mathfrak{A},
			\notag
\end{align}
where $\mathfrak{A}$ is the term on RHS~\eqref{E:ALPHARENORMALIZED}.
\end{proposition}

\begin{proof}
	From \eqref{E:CHIUSEFULID},
	the identity $[\Lunit,\CoordAng] = 0$,
	the torsion-free property of $\D$,
	and Def.~\ref{D:SPACETIMECURVATURE},
	we deduce
	$
	\displaystyle
		\angLie_{\Lunit} \upchi_{\CoordAng \CoordAng}
		= g(\D_{\CoordAng} \D_{\Lunit} \Lunit, \CoordAng)
			+ g(\D_{\CoordAng} \Lunit, \D_{\CoordAng} \Lunit)
			- \Cur_{\Lunit \CoordAng \Lunit \CoordAng}
	$.
	Viewing both sides to be a type $\binom{0}{2}$ $\ell_{t,u}-$tangent tensorfield,
	we take the $\gsphere$ trace, 
	use $\angLie_{\Lunit} \ginversesphere 
	= - 2 \upchi^{\# \#}$
	(which follows from definition \eqref{E:CHIDEF} and \eqref{E:CONNECTIONBETWEENANGLIEOFGSPHEREANDDEFORMATIONTENSORS}),
	use Lemma~\ref{L:CONNECTIONLRADFRAME},
	and carry out straightforward calculations to derive
	\begin{align} \label{E:TRCHIEVOLUTION}
		\upmu \Lunit \mytr \upchi
		& = (\Lunit \upmu) \mytr \upchi
			- \upmu (\mytr \upchi)^2
			- \upmu \Ric_{\Lunit \Lunit}.
	\end{align}
	Then from \eqref{E:ALPHARENORMALIZED} and \eqref{E:TRCHIEVOLUTION},
	we find that
	\begin{align} \label{E:LOWESTORDERTRCHIJUNKRENORMALIZEDTRANSPORT}
		\Lunit
		\left\lbrace
			\upmu \mytr \upchi 
			 + 
			 \upchifullmodinhom
		\right\rbrace
		& = 
			2 (\Lunit \upmu) \mytr \upchi
			- \upmu (\mytr \upchi)^2
			- \mathfrak{A}.
	\end{align}
	Applying $\Tanset^N$ to \eqref{E:LOWESTORDERTRCHIJUNKRENORMALIZEDTRANSPORT} 
	and performing straightforward commutations, we find that
	\begin{align} \label{E:NOTFINALFORMTOPORDERTRCHIJUNKRENORMALIZEDTRANSPORT}
		\Lunit
		\left\lbrace
			\upmu \Tanset^N \mytr \upchi 
			 + 
			 \Tanset^N \upchifullmodinhom
		\right\rbrace
		& = 
			2 (\Lunit \upmu) \Tanset^N \mytr \upchi
			- 2 \upmu \mytr \upchi \Tanset^N \mytr \upchi
			+ \upmu [\Lunit, \Tanset^N] \mytr \upchi
				\\
		& \ \ 
			+ [\Lunit, \Tanset^N] \upchifullmodinhom
			+ [\upmu, \Tanset^N] \Lunit \mytr \upchi
			+ [\Tanset^N, \Lunit \upmu] \mytr \upchi
			\notag \\
		& \ \
			- \left\lbrace 
					\Tanset^N \left(\upmu (\mytr \upchi)^2\right)
					- 2 \upmu \mytr \upchi \Tanset^N \mytr \upchi
				\right\rbrace
			- \Tanset^N \mathfrak{A}.
			\notag
	\end{align}
	The identity \eqref{E:TOPORDERTRCHIJUNKRENORMALIZEDTRANSPORT}
	now follows easily from \eqref{E:NOTFINALFORMTOPORDERTRCHIJUNKRENORMALIZEDTRANSPORT}
	and the definition of $\upchifullmodarg{\Tanset^N}$.
\end{proof}

We now derive the transport equation
verified by the partially modified quantities.

\begin{proposition}[\textbf{The transport equation for the partially modified version of} $\Tanset^{N-1} \mytr \upchi$]
\label{P:COMMUTEDTRCHIJUNKFIRSTPARTIALRENORMALIZEDTRANSPORTEQUATION}
Let $\Tanset^{N-1}$ be 
an $(N-1)^{st}$ order pure $\mathcal{P}_u-$tangent commutation vectorfield operator,
and let $\upchipartialmodarg{\Tanset^{N-1}}$ 
be the corresponding partially modified quantity defined in \eqref{E:TRANSPORTPARTIALRENORMALIZEDTRCHIJUNK}.
Then $\upchipartialmodarg{\Tanset^{N-1}}$ verifies the following transport equation:
\begin{align} \label{E:COMMUTEDTRCHIJUNKFIRSTPARTIALRENORMALIZEDTRANSPORTEQUATION}
	\Lunit \upchipartialmodarg{\Tanset^{N-1}}
	& = \frac{1}{2} G_{\Lunit \Lunit} \angLap \Tanset^{N-1} \Psi
			+ \upchipartialmodsourcearg{\Tanset^{N-1}},
\end{align}
where the inhomogeneous term $\upchipartialmodsourcearg{\Tanset^{N-1}}$ is given by
\begin{align}  \label{E:TRCHIJUNKCOMMUTEDTRANSPORTEQNPARTIALRENORMALIZATIONINHOMOGENEOUSTERM}
\upchipartialmodsourcearg{\Tanset^{N-1}}
	& = - \Tanset^{N-1} \mathfrak{B}
			- \Tanset^{N-1} (\mytr \upchi)^2
				\\
	& \ \ 
			+ \frac{1}{2} [\Tanset^{N-1}, G_{\Lunit \Lunit}] \angLap \Psi
			+ \frac{1}{2} G_{\Lunit \Lunit} [\Tanset^{N-1}, \angLap] \Psi
			+ [\Lunit, \Tanset^{N-1}] \mytr \upchi
			  \notag \\
	& \ \ 
				+ 
				[\Lunit, \Tanset^{N-1}] \upchipartialmodinhom
				+ 
				\Lunit
				\left\lbrace
					\upchipartialmodinhomarg{\Tanset^{N-1}}
					- \Tanset^{N-1} \upchipartialmodinhom
				\right\rbrace,
				\notag
\end{align}
$\mathfrak{B}$ is defined in \eqref{E:PARTIALRENORMALIZEDALPHAINHOMOGENEOUSTERM},
$	\upchipartialmodinhomarg{\Tanset^{N-1}}$ is defined in 
\eqref{E:TRANSPORTPARTIALRENORMALIZEDTRCHIJUNKDISCREPANCY},
and $\upchipartialmodinhom$
is defined in \eqref{E:LOWESTORDERTRANSPORTPARTIALRENORMALIZEDTRCHIJUNKDISCREPANCY}.
\end{proposition}

\begin{proof}
From \eqref{E:TRCHIEVOLUTION}, we find that
\begin{align} \label{E:SECONDTRCHIEVOLUTION}
		\Lunit \mytr \upchi
		& = \frac{\Lunit \upmu}{\upmu} \mytr \upchi
			- (\mytr \upchi)^2
			-\Ric_{\Lunit \Lunit}.
\end{align}
Then from \eqref{E:ALPHAPARTIALRENORMALIZED} and \eqref{E:SECONDTRCHIEVOLUTION},
we deduce that
\begin{align} \label{E:LOWESTORDERTRCHIJUNKPARTIALRENORMALIZEDTRANSPORT}
		\Lunit
		\left\lbrace
			\mytr \upchi 
			 + 
			 \upchipartialmodinhom
		\right\rbrace
		& = 
			\frac{1}{2} G_{\Lunit \Lunit} \angLap \Psi
			- (\mytr \upchi)^2
			- \mathfrak{B}.
\end{align}
We note in particular that 
the dangerous product
$
\displaystyle
\frac{\Lunit \upmu}{\upmu} \mytr \upchi
$
from
\eqref{E:ALPHAPARTIALRENORMALIZED}
and
\eqref{E:SECONDTRCHIEVOLUTION}
cancels from \eqref{E:LOWESTORDERTRCHIJUNKPARTIALRENORMALIZEDTRANSPORT}.
The desired identity \eqref{E:TRCHIJUNKCOMMUTEDTRANSPORTEQNPARTIALRENORMALIZATIONINHOMOGENEOUSTERM}
now follows from 
applying $\Tanset^{N-1}$ to
\eqref{E:LOWESTORDERTRCHIJUNKPARTIALRENORMALIZEDTRANSPORT}
and carrying out straightforward operator commutations.
\end{proof}

\section{Norms, Initial Data, Bootstrap Assumptions, and Smallness Assumptions}
\label{S:NORMSANDBOOTSTRAP}
In this section, we first introduce the pointwise norms that we use
to control solutions. We then describe our assumptions on the size of the initial
data. Finally, we state bootstrap assumptions that we use
throughout most of the rest of the paper to derive estimates.

\subsection{Norms}
\label{SS:NORMS}
In our analysis, we primarily estimate scalar functions and
$\ell_{t,u}-$tangent tensorfields. We always use the metric
$\gsphere$ when taking the pointwise norm of 
$\ell_{t,u}-$tangent tensorfields, a concept which we
make precise in the next definition.

\begin{definition}[\textbf{Pointwise norms}]
	\label{D:POINTWISENORM}
	If $\xi_{\nu_1 \cdots \nu_n}^{\mu_1 \cdots \mu_m}$ 
	is a type $\binom{m}{n}$ $\ell_{t,u}$ tensor,
	then we define the norm $|\xi| \geq 0$ by
	\begin{align} \label{E:POINTWISENORM}
		|\xi|^2
		:= 
		\gsphere_{\mu_1 \widetilde{\mu}_1} \cdots \gsphere_{\mu_m \widetilde{\mu}_m}
		(\ginversesphere)^{\nu_1 \widetilde{\nu}_1} \cdots (\ginversesphere)^{\nu_n \widetilde{\nu}_n}
		\xi_{\nu_1 \cdots \nu_n}^{\mu_1 \cdots \mu_m}
		\xi_{\widetilde{\nu}_1 \cdots \widetilde{\nu}_n}^{\widetilde{\mu}_1 \cdots \widetilde{\mu}_m}.
	\end{align}
\end{definition}

Our analysis relies on the following $L^2$ and $L^{\infty}$ norms.

\begin{definition}[$L^2$ \textbf{and} $L^{\infty}$ \textbf{norms}]
In terms of the non-degenerate forms of Def.~\ref{D:NONDEGENERATEVOLUMEFORMS},
we define the following norms for 
$\ell_{t,u}-$tangent tensorfields:
\label{D:SOBOLEVNORMS}
	\begin{subequations}
	\begin{align}  \label{E:L2NORMS}
			\left\|
				\xi
			\right\|_{L^2(\ell_{t,u})}^2
			& :=
			\int_{\ell_{t,u}}
				|\xi|^2
			\, d \spherevol,
				\qquad
			\left\|
				\xi
			\right\|_{L^2(\Sigma_t^u)}^2
			:=
			\int_{\Sigma_t^u}
				|\xi|^2
			\, d \tvol,
				\\
			\left\|
				\xi
			\right\|_{L^2(\mathcal{P}_u^t)}^2
			& :=
			\int_{\mathcal{P}_u^t}
				|\xi|^2
			\, d \conevol,
			\notag
	\end{align}

	\begin{align} 
			\left\|
				\xi
			\right\|_{L^{\infty}(\ell_{t,u})}
			& :=
				\mbox{ess sup}_{\vartheta \in \mathbb{T}}
				|\xi|(t,u,\vartheta),
			\qquad
			\left\|
				\xi
			\right\|_{L^{\infty}(\Sigma_t^u)}
			:=
			\mbox{ess sup}_{(u',\vartheta) \in [0,u] \times \mathbb{T}}
				|\xi|(t,u',\vartheta),
			\label{E:LINFTYNORMS}
				\\
			\left\|
				\xi
			\right\|_{L^{\infty}(\mathcal{P}_u^t)}
			& :=
			\mbox{ess sup}_{(t',\vartheta) \in [0,t] \times \mathbb{T}}
				|\xi|(t',u,\vartheta).
			\notag
	\end{align}
	\end{subequations}
\end{definition}

\begin{remark}[\textbf{Subset norms}]
	\label{R:SUBSETNORMS}
	In our analysis below, we occasionally use norms 
	$\| \cdot \|_{L^2(\Omega)}$
	and
	$\| \cdot \|_{L^{\infty}(\Omega)}$,
	where $\Omega$ is a subset of $\Sigma_t^u$.
	These norms are defined by replacing 
	$\Sigma_t^u$ with $\Omega$ in 
	\eqref{E:L2NORMS} and \eqref{E:LINFTYNORMS}.
\end{remark}

\subsection{Strings of commutation vectorfields and vectorfield seminorms}
\label{SS:STRINGSOFCOMMUTATIONVECTORFIELDS}
The following shorthand notation captures the relevant structure
of our vectorfield operators and allows us to depict estimates schematically.

\begin{remark}
	Some operators in Def.~\ref{D:VECTORFIELDOPERATORS} are decorated with a $*$.
	These operators involve $\mathcal{P}_u-$tangent differentiations that often
	lead to a gain in smallness in the estimates.
	More precisely, the operators $\Tanset_*^N$ always lead to a gain in smallness while
	the operators $\Fullset_*^{N;M}$ lead to a gain in smallness except 
	perhaps when they are applied to $\upmu$
	(because $\Lunit \upmu$ and its $\Rad$ derivatives are not small).
\end{remark}

\begin{definition}[\textbf{Strings of commutation vectorfields and vectorfield seminorms}] \label{D:VECTORFIELDOPERATORS}
	\ \\
	\begin{itemize}
		\item $\Fullset^{N;M} f$ 
			denotes an arbitrary string of $N$ commutation
			vectorfields in $\Fullset$ (see \eqref{E:COMMUTATIONVECTORFIELDS})
			applied to $f$, where the string contains \emph{at most} $M$ factors of the $\mathcal{P}_u^t-$transversal
			vectorfield $\Rad$. 
		\item $\Tanset^N f$
			denotes an arbitrary string of $N$ commutation
			vectorfields in $\Tanset$ (see \eqref{E:TANGENTIALCOMMUTATIONVECTORFIELDS})
			applied to $f$.
		\item 
			For $N \geq 1$,
			$\Fullset_*^{N;M} f$
			denotes an arbitrary string of $N$ commutation
			vectorfields in $\Fullset$ 
			applied to $f$, where the string contains \emph{at least} one $\mathcal{P}_u^t-$tangent factor 
			and \emph{at most} $M$ factors of $\Rad$.
			We also set  $\Fullset_*^{0;0} f := f$.
		\item For $N \geq 1$,
					$\Tanset_*^N f$ 
					denotes an arbitrary string of $N$ commutation
					vectorfields in $\Tanset$ 
					applied to $f$, where the string contains
					\emph{at least one factor} of $\GeoAng$ or \emph{at least two factors} of $\Lunit$.
		\item For $\ell_{t,u}-$tangent tensorfields $\xi$, 
					we similarly define strings of $\ell_{t,u}-$projected Lie derivatives 
					such as $\angLie_{\Fullset}^{N;M} \xi$.
	\end{itemize}

	We also define pointwise seminorms constructed out of sums of the above strings of vectorfields:
	\begin{itemize}
		\item $|\Fullset^{N;M} f|$ 
		simply denotes the magnitude of one of the $\Fullset^{N;M} f$ as defined above
		(there is no summation).
	\item $|\Fullset^{\leq N;M} f|$ is the \emph{sum} over all terms of the form $|\Fullset^{N';M} f|$
			with $N' \leq N$ and $\Fullset^{N';M} f$ as defined above.
			When $N=M=1$, we sometimes write $|\Fullset^{\leq 1} f|$ instead of $|\Fullset^{\leq 1;1} f|$.
		\item $|\Fullset^{[1,N];M} f|$ is the sum over all terms of the form $|\Fullset^{N';M} f|$
			with $1 \leq N' \leq N$ and $\Fullset^{N';M} f$ as defined above.
		\item Sums such as 
			$|\Tanset_*^{[1,N]} f|$,
			$|\angLie_{\Fullset}^{\leq N;M} \xi|$,
			$|\GeoAng^{\leq 1} f|$,
			$|\Rad^{[1,N]} f|$,
			etc.,
			are defined analogously.
			For example, 
			$|\Rad^{[1,N]} f| 
			= |\Rad f| 
			+ |\Rad \Rad f| 
			+ \cdots 
			+ |\overbrace{\Rad \Rad \cdots \Rad}^{N \mbox{ \upshape copies}} f|
			$.
			We write $|\Tanset_* f|$ instead of $|\Tanset_*^{[1,1]} f|$.
		\end{itemize}

\end{definition}

\subsection{Assumptions on the initial data and the behavior of quantities along \texorpdfstring{$\Sigma_0$}{the initial data hypersurface}}
\label{SS:DATAASSUMPTIONS}
In this subsection, we introduce our Sobolev norm assumptions on the data,
which involve several size parameters.
We then derive identities and estimates
for various quantities on $\Sigma_0$.
In Subsect.~\ref{SS:SMALLNESSASSUMPTIONS}, we describe 
our assumptions on the size parameters.

We first recall that
$
(\Psi|_{\Sigma_0},\partial_t \Psi|_{\Sigma_0})
:=
(\mathring{\Psi},\mathring{\Psi}_0)
$
and that we assume
$(\mathring{\Psi},\mathring{\Psi}_0) \in H_{\Euct}^{19}(\Sigma_0^1) \times H_{\Euct}^{18}(\Sigma_0^1)$.
We assume that the data verify the following size estimates
(see Subsect.~\ref{SS:STRINGSOFCOMMUTATIONVECTORFIELDS} regarding the vectorfield operator notation):
\begin{align} \label{E:PSIDATAASSUMPTIONS}
	\left\|
		\Fullset_*^{\leq 17;3} \Psi
	\right\|_{L^{\infty}(\Sigma_0^1)},
		\,
	\left\|
		\Fullset_*^{\leq 19;3} \Psi
	\right\|_{L^2(\Sigma_0^1)}
	& \leq \mathring{\upepsilon},
		\qquad
	\left\|
		\Rad^{[1,3]} \Psi
	\right\|_{L^{\infty}(\Sigma_0^1)}
	:= \mathring{\updelta}
		> 0.
\end{align}

\begin{remark}[\textbf{Non-optimal regularity assumptions involving higher transversal derivatives}]
\label{R:REGULARITYINVOLVINGHIGHERTRANSVERSAL}
The data assumptions \eqref{E:PSIDATAASSUMPTIONS} involving two or more
transversal derivatives of $\Psi$ are not optimal relative to our proof;
we have stated our assumptions in the form \eqref{E:PSIDATAASSUMPTIONS} only 
for convenience.
For example, in the parts of our proof that involve three transversal derivatives of
$\Psi$
(see Prop.~\ref{P:IMPROVEMENTOFHIGHERTRANSVERSALBOOTSTRAP}), 
we use only the data assumptions
$
\left\|
	\Fullset_*^{\leq 4;2} \Psi
\right\|_{L^{\infty}(\Sigma_0^1)}
\leq \mathring{\upepsilon}
$,
$
\left\|
	\Lunit \Rad \Rad \Rad \Psi
\right\|_{L^{\infty}(\Sigma_0^1)}
\leq \mathring{\upepsilon}
$,
and
$
\left\|
	\Rad \Rad \Rad \Psi
\right\|_{L^{\infty}(\Sigma_0^1)}
\leq \mathring{\updelta}
$.
Similar remarks apply to the initial regularity of the eikonal function quantities, 
which we exhibit in Lemma~\ref{L:BEHAVIOROFEIKONALFUNCTIONQUANTITIESALONGSIGMA0}.
\end{remark}

In the next definition, we introduce the data-dependent number 
$\TranminusdatasizeWithFactor$, which is of crucial importance.
Our main theorem shows that 
for $\mathring{\upepsilon}$ sufficiently small,
the time of first shock formation is 
$(1 + \mathcal{O}(\mathring{\upepsilon}))\TranminusdatasizeWithFactor^{-1}$.

\begin{definition}[\textbf{The quantity that controls the blowup-time}]
	\label{D:CRITICALBLOWUPTIMEFACTOR}
	We define
	\begin{align} \label{E:CRITICALBLOWUPTIMEFACTOR}
		\TranminusdatasizeWithFactor
		& := \frac{1}{2} \sup_{\Sigma_0^1} \left[G_{\Lunit \Lunit} \Rad \Psi \right]_-.
	\end{align}
\end{definition}

\begin{remark}
Our proof of shock formation relies on the assumptions that $\TranminusdatasizeWithFactor > 0$
and that $\TranminusdatasizeWithFactor$ is not too small compared to other 
quantities; see Subsect.~\ref{SS:SMALLNESSASSUMPTIONS}.
These assumptions are tied to the fact that we are studying 
blowup for perturbations of right-moving simple plane symmetric waves,
as is depicted in Figure~\ref{F:FRAME}.
\end{remark}

To prove our main theorem,
we make assumptions on the
relative sizes of the above
parameters; see Subsect.~\ref{SS:SMALLNESSASSUMPTIONS}.

Our next goal is to derive estimates for various quantities along
$\Sigma_0^1$ that hold whenever the data verify
\eqref{E:PSIDATAASSUMPTIONS} and $\mathring{\upepsilon}$
is sufficiently small. We start by providing two lemmas that
yield some identities that are relevant for that analysis.

\begin{lemma}[\textbf{Identities involving} $\partial_i u$]
	The following identity holds:
	\begin{align} \label{E:ALTERNATEEIKONAL}
		(\gt^{-1})^{ab} \partial_a u \partial_b u 
		 = \upmu^{-2}.
	\end{align}
	In \eqref{E:ALTERNATEEIKONAL},
	$\gt^{-1}$ is the inverse of the 
	Riemannian metric $\gt$ on $\Sigma_t^{U_0}$ defined by \eqref{E:GTANDGSPHERESPHEREDEF}.

	Furthermore, the rectangular spatial derivatives 
	of $u$ verify (for $i=1,2$):
	\begin{align} \label{E:PARTIALIUINTERMSOFRADUNIT}
		\upmu \partial_i u 
		& = \Radunit_i.
	\end{align}
\end{lemma}

\begin{proof}
	Since $- \upmu \partial_{\alpha} u = \Lunit_{\alpha}$ 
	(see \eqref{E:LGEOEQUATION} and \eqref{E:LUNITDEF}),
	\eqref{E:PARTIALIUINTERMSOFRADUNIT} follows from
	the first identity in \eqref{E:DOWNSTAIRSUPSTAIRSSRADUNITPLUSLUNITISAFUNCTIONOFPSI}.
	To deduce \eqref{E:ALTERNATEEIKONAL}, we 
	use \eqref{E:PARTIALIUINTERMSOFRADUNIT}
	to substitute for $\partial_a u$ and $\partial_b u$ on LHS~\eqref{E:PARTIALIUINTERMSOFRADUNIT} 
	and use the normalization condition
	$(\gt^{-1})^{ab} \Radunit_a \Radunit_b = g(\Radunit,\Radunit) = 1$
	(see \eqref{E:RADIALVECTORFIELDSLENGTHS}).
\end{proof}

\begin{lemma}[\textbf{Algebraic identities along} $\Sigma_0$]
	\label{L:ALGEBRAICIDALONGSIGMA0}
The following identities hold along $\Sigma_0$ (for $i=1,2$):
\begin{align}
	\upmu
	& = \frac{1}{\sqrt{(\gtinverse)^{11}}},
	\qquad
	\Lunit_{(Small)}^i
	 = 	\frac{(\gtinverse)^{i1}}{\sqrt{(\gtinverse)^{11}}} 
	 		- \delta^{i1}
			- (g^{-1})^{0i},
	\qquad
	\NonRadialRad^i
	 =  \frac{(\gtinverse)^{i1}}{(\gtinverse)^{11}}
	 	- \delta^{i1},
		\label{E:INITIALRELATIONS}
\end{align}
where 
$\gt$ is viewed as the $2 \times 2$ matrix of rectangular spatial components 
of the Riemannian metric on $\Sigma_0$ defined by \eqref{E:GTANDGSPHERESPHEREDEF},
$\gt^{-1}$ is the corresponding inverse matrix,
and $\NonRadialRad$ is the $\ell_{t,u}-$tangent vectorfield from \eqref{E:RADSPLITINTOPARTTILAUANDXI}.
\end{lemma}

\begin{proof}
	The identity for $\upmu$ is a simple consequence of
	\eqref{E:ALTERNATEEIKONAL}
	and the fact that by construction, $u|_{\Sigma_0} = 1 - x^1$
	(see \eqref{E:INTROEIKONALINITIALVALUE}).
	Next, using in addition \eqref{E:PARTIALIUINTERMSOFRADUNIT}, 
	we deduce that 
	$\Radunit^i 
	= - \upmu (\gtinverse)^{ia} \delta_a^1 
	= - \upmu (\gtinverse)^{i1}
	= - \frac{(\gtinverse)^{i1}}{\sqrt{(\gtinverse)^{11}}}
	$ along $\Sigma_0$.
	Also using that $\Lunit^i = \Timenormal^i - \Radunit^i$
	(see \eqref{E:TIMENORMAL}),
	$\Timenormal^i = - (g^{-1})^{0i}$ (see \eqref{E:TIMENORMALRECTANGULAR}),
	and $\Lunit_{(Small)}^i = \Lunit^i - \delta_1^i$
	(see \eqref{E:PERTURBEDPART}), 
	we easily conclude the desired identity for $\Lunit_{(Small)}^i$.
	Next, we recall that by construction, $\vartheta|_{\Sigma_0} = x^2$
	(see Def.~\ref{D:ANGULARCOORDINATE}).
	Hence, $\frac{\partial}{\partial u} = - \partial_1$ and $\CoordAng = \partial_2$
	along $\Sigma_0$.
	Also using that $\Rad = \frac{\partial}{\partial u} - \NonRadialRad$ (see \eqref{E:RADSPLITINTOPARTTILAUANDXI}),
	\eqref{E:PARTIALIUINTERMSOFRADUNIT}, and $\Rad = \upmu \Radunit$,
	we conclude that
	$\NonRadialRad^i 
	= - \delta^{i1} - \Rad^i
	= - \delta^{i1} 
	+ \upmu \frac{(\gtinverse)^{i1}}{\sqrt{(\gtinverse)^{11}}}
	= - \delta^{i1} 
		+ \frac{(\gtinverse)^{i1}}{(\gtinverse)^{11}}
	$
	as desired.
\end{proof}

In the next lemma, we provide estimates verified by the
eikonal function quantities $\upmu$ and $\Lunit_{(Small)}^i$ along $\Sigma_0^1$. The estimates
are a consequence of the assumptions 
\eqref{E:PSIDATAASSUMPTIONS} on the initial data of $\Psi$
as well as the evolution equations verified
by $\upmu$ and $\Lunit_{(Small)}^i$.

\begin{lemma}[\textbf{Behavior of the eikonal function quantities along} $\Sigma_0^1$]
\label{L:BEHAVIOROFEIKONALFUNCTIONQUANTITIESALONGSIGMA0}
For initial data verifying \eqref{E:PSIDATAASSUMPTIONS},
the following $L^2$ and $L^{\infty}$ estimates hold along $\Sigma_0^1$
whenever $\mathring{\upepsilon}$ is sufficiently small,
where the implicit constants are allowed to depend on 
$\mathring{\updelta}$:
\begin{align} \label{E:LUNITIDATAL2CONSEQUENCES}
	\left\|
		\Fullset_*^{\leq 19;3} \Lunit_{(Small)}^i
	\right\|_{L^2(\Sigma_0^1)}
	& \lesssim \mathring{\upepsilon},
		\qquad
	\left\|
		\Rad^{[1,3]} \Lunit_{(Small)}^i
	\right\|_{L^2(\Sigma_0^1)}
	\lesssim 1,
\end{align}

\begin{subequations}
\begin{align}  \label{E:UPMUDATATANGENTIALL2CONSEQUENCES}
	\left\|
		\upmu - 1
	\right\|_{L^2(\Sigma_0^1)},
		\,
	\left\|
		\Tanset_*^{[1,19]} \upmu
	\right\|_{L^2(\Sigma_0^1)}
	& \lesssim \mathring{\upepsilon},
		\\
	\left\|
		\Lunit \Rad^{[0,2]} \upmu 
	\right\|_{L^2(\Sigma_0^1)},
		\,
	\left\|
		\Rad^{[0,2]} \Lunit \upmu 
	\right\|_{L^2(\Sigma_0^1)},
		\,
	\left\|
		\Rad \Lunit \Rad \upmu 
	\right\|_{L^2(\Sigma_0^1)},
		\,
	\left\|
		\Rad^{[1,2]} \upmu 
	\right\|_{L^2(\Sigma_0^1)}
	& \lesssim 1,
	\label{E:UPMUDATARADIALL2CONSEQUENCES}
\end{align}
\end{subequations}

\begin{align} \label{E:LUNITIDATALINFTYCONSEQUENCES}
	\left\|
		\Fullset_*^{\leq 17;2} \Lunit_{(Small)}^i
	\right\|_{L^{\infty}(\Sigma_0^1)}
	& \lesssim \mathring{\upepsilon},
		\qquad
	\left\|
		\Rad^{[1,2]} \Lunit_{(Small)}^i
	\right\|_{L^{\infty}(\Sigma_0^1)}
	\lesssim 1,
\end{align}

\begin{subequations}
\begin{align}  \label{E:UPMUDATATANGENTIALLINFINITYCONSEQUENCES}
	\left\|
		\upmu - 1
	\right\|_{L^{\infty}(\Sigma_0^1)},
		\,
	\left\|
		\Tanset_*^{[1,17]} \upmu
	\right\|_{L^{\infty}(\Sigma_0^1)}
	& \lesssim \mathring{\upepsilon},
		\\
	\left\|
		\Lunit \Rad^{[0,2]} \upmu 
	\right\|_{L^{\infty}(\Sigma_0^1)},
		\,
	\left\|
		\Rad^{[0,2]} \Lunit \upmu 
	\right\|_{L^{\infty}(\Sigma_0^1)},
		\,
	\left\|
		\Rad \Lunit \Rad \upmu 
	\right\|_{L^{\infty}(\Sigma_0^1)},
		\,
	\left\|
		\Rad^{[1,2]} \upmu 
	\right\|_{L^{\infty}(\Sigma_0^1)}
	& \lesssim 1.
	\label{E:UPMUDATARADIALLINFINITYCONSEQUENCES}
\end{align}
\end{subequations}

\end{lemma}

\begin{proof}[Sketch of proof]
	We only sketch the proofs of the estimates because they have a lengthy component 
	and because we give complete details of related but more complicated 
	estimates in our proofs of
	Propositions~\ref{P:IMPROVEMENTOFAUX} and \ref{P:IMPROVEMENTOFHIGHERTRANSVERSALBOOTSTRAP}
	below. 
	An easy part of the proof is deriving estimates
	involving derivatives with respect to 
	the $\Sigma_0-$tangent vectorfields
	$\Rad$ 
	and
	$\GeoAng$;
	we can use the identities \eqref{E:INITIALRELATIONS} 
	to express $(\upmu - 1)|_{\Sigma_0} $ and $\Lunit_{(Small)}^i|_{\Sigma_0}$
	in the form $\smoothfunction(\Psi)\Psi$ with $\smoothfunction$ smooth. 
	We can then repeatedly differentiate $\smoothfunction(\Psi)\Psi$ with respect to
	$\Rad, \GeoAng$ and use the assumptions
	\eqref{E:PSIDATAASSUMPTIONS} and the standard Sobolev calculus
	to obtain the desired $L^2$ and $L^{\infty}$ estimates.
	Under the umbrella of the Sobolev calculus, we include the embedding estimate
	$\| f \|_{L^{\infty}(\ell_{0,u})} \lesssim \| \GeoAng f \|_{L^2(\ell_{0,u})} + \| f \|_{L^2(\ell_{0,u})}$
	(valid with a uniform implicit constant for $u \in [0,U_0]$).
	This Sobolev embedding estimate is easy to derive because
	along $\ell_{0,u}$, we have
	$\CoordAng = \partial_2$,
	$\GeoAng = (1 + \mathcal{O}(\mathring{\upepsilon})) \partial_2$,
	and
	$d \spherevol = (1 + \mathcal{O}(\mathring{\upepsilon})) \, d \vartheta$.
	Thus, the embedding result follows from the standard one on the 
	torus $\mathbb{T}$ equipped with the standard Euclidean metric.

	Another easy part of the proof is 
	obtaining the desired estimates for $\Lunit^{\leq 19} \upmu$, $\Lunit^{\leq 19} \Lunit_{(Small)}^i$,
	and their derivatives up to top order with
	respect to $\Rad$ and $\GeoAng$
	(where all $\Rad$ and $\GeoAng$ derivatives occur after the $\Lunit$ differentiations).
	To derive them, 
	we can repeatedly use the evolution equations
	\eqref{E:UPMUFIRSTTRANSPORT}
	and
	\eqref{E:LLUNITI}
	to substitute for $\Lunit \upmu$ 
	and $\Lunit \Lunit_{(Small)}^i$
	and argue as in the previous paragraph.
	In obtaining these estimates, 
	the main point 
	(which is easy to see with the help of
	\eqref{E:UPMUFIRSTTRANSPORT}
	and
	\eqref{E:LLUNITI})
	is that the quantities 
	that we claim are $\lesssim \mathring{\upepsilon}$
	contain at least one small factor
	$\Fullset_*^{\leq 19;3} \Psi$,
	which by \eqref{E:PSIDATAASSUMPTIONS} 
	yields the desired smallness factor
	$\mathring{\upepsilon}$.

	The lengthy part of the proof is deriving estimates for
	the derivatives of $\upmu$ and $\Lunit_{(Small)}^i$
	that involve both $\Lunit$ and the $\Sigma_0-$tangent operators 
	$\lbrace \Rad, \GeoAng \rbrace$,
	where the $\Sigma_0-$tangent differentiation acts before $\Lunit$ does.
	A model term is $\Lunit \Rad \Lunit_{(Small)}^i$.
	The main idea of the argument is to first write 
	$\Lunit \Rad \Lunit_{(Small)}^i = \Rad \Lunit \Lunit_{(Small)}^i + [\Lunit,\Rad] \Lunit_{(Small)}^i$.
	The advantage of this decomposition is that the arguments given in the previous paragraph 
	imply that 
	$\left\| \Rad \Lunit \Lunit_{(Small)}^i \right\|_{L^2(\Sigma_0^1)} \lesssim \mathring{\upepsilon}$. 
	Thus, the main step remaining is to establish
	commutation estimates showing
	that, roughly speaking, the commutator operators
	$[\Lunit, \Rad]$, 
	$[\Lunit, \GeoAng]$, 
	and
	$[\Rad, \GeoAng]$
	lead to products involving at least one (good)
	$\ell_{t,u}-$tangent differentiation, which provides the smallness factor
	$\mathring{\upepsilon}$ in the relevant expressions.
	The identities in \eqref{E:CONNECTIONBETWEENCOMMUTATORSANDDEFORMATIONTENSORS}
	feature $\ell_{t,u}-$tangent right-hand sides and thus
	imply the availability of the desired structure after one commutation.
	To derive estimates up to top order, 
	we must also show that a suitable version of this structure
	survives under higher-order 
	differentiations and commutations. 
	We establish the necessary commutation estimates in
	Lemmas~\ref{L:COMMUTATORESTIMATES},
	\ref{L:TRANSVERALTANGENTIALCOMMUTATOR},
	and \ref{L:HIGHERTRANSVERALTANGENTIALCOMMUTATOR}
	below under bootstrap assumptions that are consistent
	with the evolution of the solution. 
	In the inequalities stated in those lemmas, 
	the bootstrap assumptions are used 
	to gain a factor of $\varepsilon^{1/2}$
	in various quadratic terms that appear on the right-hand sides,
	where $\varepsilon$ is a small bootstrap parameter.
	Along $\Sigma_0$, 
	the lemmas can be established 
	without the bootstrap assumptions by using the same arguments given in their proofs.
	In fact, one can ignore the availability of the smallness factor.
	After establishing the commutation Lemmas along $\Sigma_0$
	we can derive the desired estimates
	using induction in the number of derivatives.
	For example, to deduce that
	$\left\| \Lunit \Lunit \Rad \upmu \right\|_{L^2(\Sigma_0^1)} \lesssim \mathring{\upepsilon}$,
	we write 
	$\Lunit \Lunit \Rad \upmu 
	= \Rad \Lunit \Lunit \upmu
	+ [\Lunit \Lunit, \Rad] \upmu
	$. The argument sketched in the previous paragraph implies that
	$\left\| \Rad \Lunit \Lunit \upmu \right\|_{L^2(\Sigma_0^1)} \lesssim \mathring{\upepsilon}$.
	To conclude that
	$\left\| [\Lunit \Lunit, \Rad] \upmu \right\|_{L^2(\Sigma_0^1)} \lesssim \mathring{\upepsilon}$,
	we can use the commutator estimate \eqref{E:ONERADIALTANGENTIALFUNCTIONCOMMUTATORESTIMATE}
	to derive the pointwise bound 
	$\left| [\Lunit \Lunit, \Rad] \upmu \right| 
	\lesssim 
	\left|
		\GeoAng \Fullset^{\leq 1;1} \upmu
	\right|
	+
	\left|
		\Tanset_*^{[1,2]} \BadVar
	\right|
	+
	\left|
		\Fullset_*^{\leq 2;1} \GdVar
	\right|
	+
	\left|
		\Tanset^{\leq 2} \GdVar
	\right|
	$.
	The RHS involves only up-to-order $2$ derivatives of quantities
	that, by induction, 
	would have been shown to be bounded in the norm
	$\left\| \cdot \right\|_{L^2(\Sigma_0^1)}$ 
	by
	$\lesssim \mathring{\upepsilon}$.
\end{proof}

\subsection{\texorpdfstring{$\Tboot$}{The bootstrap time}, the positivity of 
\texorpdfstring{$\upmu$}{the inverse foliation density}, and the diffeomorphism property of 
\texorpdfstring{$\Upsilon$}{the change of variables map}}
\label{SS:SIZEOFTBOOT}
We now state some basic bootstrap assumptions.
We start by fixing a real number $\Tboot$ with
\begin{align} \label{E:TBOOTBOUNDS}
	0 < \Tboot \leq 2 \TranminusdatasizeWithFactor^{-1}.
\end{align}

We assume that on the spacetime domain $\mathcal{M}_{\Tboot,U_0}$
(see \eqref{E:MTUDEF}), we have
\begin{align} \label{E:BOOTSTRAPMUPOSITIVITY} \tag{$\mathbf{BA} \upmu > 0$}
	\upmu > 0.
\end{align}
Inequality \eqref{E:BOOTSTRAPMUPOSITIVITY} implies that no shocks are present in
$\mathcal{M}_{\Tboot,U_0}$.

We also assume that
\begin{align} \label{E:BOOTSTRAPCHOVISDIFFEO}
	& \mbox{The change of variables map $\Upsilon$ from Def.~\ref{D:CHOVMAP}
	is a $C^1$ diffeomorphism from} \\
	& [0,\Tboot) \times [0,U_0] \times \mathbb{T}
	\mbox{ onto its image}.
	\notag
\end{align}

\subsection{Fundamental \texorpdfstring{$L^{\infty}$}{essential sup-norm} bootstrap assumptions}
\label{SS:PSIBOOTSTRAP}
 Our fundamental bootstrap assumptions for $\Psi$ are that the following inequalities hold on $\mathcal{M}_{\Tboot,U_0}$
 (see Subsect.~\ref{SS:STRINGSOFCOMMUTATIONVECTORFIELDS} regarding the vectorfield operator notation):
\begin{align} \label{E:PSIFUNDAMENTALC0BOUNDBOOTSTRAP} \tag{$\mathbf{BA}\Psi$}
	\left\| 
		\Tanset^{\leq 11} \Psi 
	\right\|_{L^{\infty}(\Sigma_t^u)}
	& \leq \varepsilon,
\end{align}
where $\varepsilon$ is a small positive bootstrap parameter whose smallness 
we describe in Subsect.~\ref{SS:SMALLNESSASSUMPTIONS}.

\subsection{Auxiliary  \texorpdfstring{$L^{\infty}$}{essential sup-norm} bootstrap assumptions}
\label{SS:AUXILIARYBOOTSTRAP}
In deriving pointwise estimates, we find it convenient to make the
following auxiliary bootstrap assumptions.
In Prop.~\ref{P:IMPROVEMENTOFAUX}, we will derive strict improvements
of these assumptions.

\noindent \underline{\textbf{Auxiliary bootstrap assumptions for small quantities}.}
We assume that the following inequalities hold on $\mathcal{M}_{\Tboot,U_0}$:
\begin{align} \label{E:PSIAUXLINFINITYBOOTSTRAP} \tag{$\mathbf{AUX1}\Psi$}
	\left\| 
		\Fullset_*^{\leq 10;1} \Psi 
	\right\|_{L^{\infty}(\Sigma_t^u)}
	& \leq \varepsilon^{1/2},
\end{align}

\begin{align}
	\left\| 
		\Lunit \Tanset^{[1,9]} \upmu
	\right\|_{L^{\infty}(\Sigma_t^u)},
		\,
	\left\| 
		\Tanset_*^{[1,9]} \upmu
	\right\|_{L^{\infty}(\Sigma_t^u)}
	& \leq \varepsilon^{1/2},
		\label{E:UPMUBOOT}  \tag{$\mathbf{AUX1}\upmu$} 
		\\
	\left\| 
		\Fullset_*^{\leq 9;1} \Lunit_{(Small)}^i 
	\right\|_{L^{\infty}(\Sigma_t^u)}
	& \leq \varepsilon^{1/2},
		\label{E:FRAMECOMPONENTSBOOT} \tag{$\mathbf{AUX1}\Lunit_{(Small)}$} 
		\\
	\left\| 
		\angLie_{\Fullset}^{\leq 8;1} \upchi
	\right\|_{L^{\infty}(\Sigma_t^u)}
	& \leq \varepsilon^{1/2}.
	\label{E:CHIBOOT} \tag{$\mathbf{AUX1}\upchi$}
\end{align}

\noindent \underline{\textbf{Auxiliary bootstrap assumptions for quantities that are allowed to be large}.}
\begin{align}
	\left\| 
		\Rad \Psi 
	\right\|_{L^{\infty}(\Sigma_t^u)}
	& \leq 
	\left\| 
		\Rad \Psi 
	\right\|_{L^{\infty}(\Sigma_0^u)}
	+ \varepsilon^{1/2},
	&& 
		\label{E:PSITRANSVERSALLINFINITYBOUNDBOOTSTRAP} \tag{$\mathbf{AUX2}\Psi$}
\end{align}

\begin{align}
	\left\| 
		\Lunit \upmu
	\right\|_{L^{\infty}(\Sigma_t^u)}
	& \leq 
		\frac{1}{2}
		\left\| 
			G_{\Lunit \Lunit} \Rad \Psi
		\right\|_{L^{\infty}(\Sigma_0^u)}
		+ \varepsilon^{1/2},
		 \label{E:LUNITUPMUBOOT}  \tag{$\mathbf{AUX2}\upmu$}  \\
		\left\| 
			\upmu
		\right\|_{L^{\infty}(\Sigma_t^u)}
		& \leq
		1
		+ 
		2 \TranminusdatasizeWithFactor^{-1} 
		\left\| 
				G_{\Lunit \Lunit} \Rad \Psi
		\right\|_{L^{\infty}(\Sigma_0^u)}
		+ \varepsilon^{1/2},
			\label{E:UPMUTRANSVERSALBOOT}  \tag{$\mathbf{AUX3}\upmu$} \\
		\left\| 
			\Rad \Lunit_{(Small)}^i
		\right\|_{L^{\infty}(\Sigma_t^u)}
		& \leq
	 	\left\| 
			\Rad \Lunit_{(Small)}^i
		\right\|_{L^{\infty}(\Sigma_0^u)}
		+ \varepsilon^{1/2}.
		\label{E:LUNITITRANSVERSALBOOT}  \tag{$\mathbf{AUX2}\Lunit_{(Small)}$} 
\end{align}

\subsection{Smallness assumptions}
\label{SS:SMALLNESSASSUMPTIONS}
For the remainder of the article, 
when we say that ``$A$ is small relative to $B$''
we mean that $B > 0$ and that
there exists a continuous increasing function 
$f :(0,\infty) \rightarrow (0,\infty)$ 
such that 
$
\displaystyle
A \leq f(B)
$.
In principle, the functions $f$ could always be chosen to be 
polynomials with positive coefficients or exponential functions.\footnote{The exponential functions appear, for example, in our energy estimates, 
during our Gronwall argument; see the proof of Prop.~\ref{P:MAINAPRIORIENERGY}
given in Subsect.~\ref{SS:PROOFOFPROPMAINAPRIORIENERGY}.} 
However, to avoid lengthening the paper, we typically do not 
specify the form of $f$.

Throughout the rest of the paper, we make the following
relative smallness assumptions. We
continually adjust the required smallness
in order to close our estimates.
\begin{itemize}
	\item $\varepsilon$ 
		is small relative to $\mathring{\updelta}^{-1}$,
		where $\mathring{\updelta}$ is the data-size parameter 
		from \eqref{E:PSIDATAASSUMPTIONS}.
	\item $\varepsilon$ is small relative to 
		the data-size parameter $\TranminusdatasizeWithFactor$ 
		from \eqref{E:CRITICALBLOWUPTIMEFACTOR}.
\end{itemize}
The first assumption will allow us to control error terms that,
roughly speaking, are of size $\varepsilon \mathring{\updelta}^k$ 
for some integer $k \geq 0$. The second assumption 
is relevant because the expected blowup-time is 
approximately $\TranminusdatasizeWithFactor^{-1}$,
and the assumption will allow us to show that various
error products featuring a small factor $\varepsilon$
remain small for $t < 2 \TranminusdatasizeWithFactor^{-1}$, 
which is plenty of time for us to show that a shock forms.

Finally, we assume that
\begin{align} \label{E:DATAEPSILONVSBOOTSTRAPEPSILON}
	\mathring{\upepsilon} 
	& \leq \varepsilon,
\end{align}
where $\mathring{\upepsilon}$ is the data smallness parameter from \eqref{E:PSIDATAASSUMPTIONS}.

\begin{remark}
	$\mathring{\updelta}$ and $\TranminusdatasizeWithFactor$
	do not have to be small.
\end{remark}

\begin{remark}[\textbf{The existence of data verifying the size assumptions}]
\label{R:EXISTENCEOFDATA}
We now sketch why there exists an open set of data
$(\Psi|_{\Sigma_0},\partial_t \Psi|_{\Sigma_0}) = (\mathring{\Psi},\mathring{\Psi}_0)$ 
that are compactly supported in $\Sigma_0^1$ and that satisfy the above size assumptions
involving $\mathring{\upepsilon}$, $\mathring{\updelta}^{-1}$, and $\TranminusdatasizeWithFactor$.
It is enough to show that there exist plane symmetric data $(\mathring{\Psi},\mathring{\Psi}_0)$ 
(that is, data depending only on $x^1$) because the size assumptions are stable under  
Sobolev-class perturbations (without symmetry), where the 
relevant Sobolev space is
$H_{\Euct}^{19}(\Sigma_0^1) \times H_{\Euct}^{18}(\Sigma_0^1)$.
Note that in plane symmetry, 
$\CoordAng$ and $\GeoAng$ are proportional to $\partial_2$
and all $\partial_2$ derivatives of all scalar functions defined throughout the article vanish.
Moreover, even though $\Lunit$ and $\Rad$ do not necessarily commute,
$[\Lunit, \Rad]$ is $\ell_{t,u}-$tangent 
(see Lemma~\ref{L:CONNECTIONBETWEENCOMMUTATORSANDDEFORMATIONTENSORS})
and therefore proportional to $\partial_2$.
Thus, when $[\Lunit, \Rad]$ acts as a differential operator on a scalar function,
it annihilates it.

To see that suitable plane symmetric data exist,
we first note that since 
Remark~\ref{R:RADUNITSMALLLUNITSMALLRELATION}
and \eqref{E:INITIALRELATIONS} imply that
$\Rad = - (1 + \mathcal{O}(\Psi)) \partial_1$ along $\Sigma_0$
and $\Lunit_{(Small)}^i = \mathcal{O}(\Psi)$
along $\Sigma_0^1$,
a simple argument relative to rectangular coordinates (omitted here) yields 
that it is possible to find smooth plane symmetric data such that
$\| \mathring{\Psi} \|_{L^{\infty}(\Sigma_0^1)}$ is as small as we want relative to
$\frac{1}{2} \sup_{\Sigma_0^1} \left[G_{\Lunit \Lunit} \Rad \Psi \right]_-$
(see definition \eqref{E:CRITICALBLOWUPTIMEFACTOR})
and relative to $1/\| \Rad^{[1,3]} \mathring{\Psi} \|_{L^{\infty}(\Sigma_0^1)}$
(see \eqref{E:PSIDATAASSUMPTIONS});
for example, one can consider functions $\mathring{\Psi}$ that have a small amplitude
but with $|\partial_1 \mathring{\Psi}|$ relatively large in some very small sub-interval of $[0,1]$
(that is, with a short-but-steep peak).
If $\partial_1 \mathring{\Psi}$ has the correct sign in the sub-interval,
this will produce the desired relative largeness of
$\frac{1}{2} \sup_{\Sigma_0^1} \left[G_{\Lunit \Lunit} \Rad \Psi \right]_-$.
Moreover, since
$\Lunit \Psi|_{\Sigma_0^1} 
= \partial_t \Psi|_{\Sigma_0^1} 
+ \Lunit^1 \partial_1 \Psi|_{\Sigma_0^1}
= \mathring{\Psi}_0
+ \Lunit^1 \partial_1 \mathring{\Psi}
$,
we can choose $\mathring{\Psi}_0$ in terms of $\mathring{\Psi}$
so that 
$
\| \Rad^{[1,3]} \Lunit \Psi \|_{L^{\infty}(\Sigma_0^1)}
$
is as small as we want
(we could even make $\Lunit \Psi|_{\Sigma_0^1} \equiv 0$
by setting $\mathring{\Psi}_0 := - \Lunit^1 \partial_1 \mathring{\Psi}$). 
Moreover, from the above remarks, we conclude 
that the same smallness holds for all permutations of the operators
$\Rad^{[1,3]} \Lunit$ acting on $\Psi$ along $\Sigma_0^1$.
To obtain the desired smallness of 
$
\| \Rad^K \Lunit^J \Psi \|_{L^{\infty}(\Sigma_0^1)}
$
for $1 \leq J$, $K \leq 3$, and $J + K \leq 17$
and
$
\| \Rad^K \Lunit^J \Psi \|_{L^2(\Sigma_0^1)}
$
for $1 \leq J$, $K \leq 3$, and $J + K \leq 19$,
we can inductively use the evolution equations
\eqref{E:LONINSIDEGEOMETRICWAVEOPERATORFRAMEDECOMPOSED},
\eqref{E:UPMUFIRSTTRANSPORT},
and
\eqref{E:LLUNITI}
and the relations \eqref{E:INITIALRELATIONS},
much like we described in the proof sketch of 
Lemma~\ref{L:BEHAVIOROFEIKONALFUNCTIONQUANTITIESALONGSIGMA0}
(note that we must simultaneously derive estimates for the derivatives of 
$\upmu$ and $\Lunit_{(Small)}^i$ in order to obtain the smallness estimates for $\Psi$).
See also Subsect.~\ref{SS:REMARKSONSMALLNESS} for
a discussion, based on the method of Riemann invariants,
of the existence of plane symmetric data verifying the desired
size assumptions in the case of the irrotational relativistic Euler equations.
\end{remark}

\section{Preliminary pointwise estimates}
\label{S:PRELIMINARYPOINTWISE}
In this section, we use the assumptions on the data and the
bootstrap assumptions from Sect.~\ref{S:NORMSANDBOOTSTRAP}
to derive pointwise estimates for the simplest error terms that appear in the commuted wave equation.
The arguments are tedious but not too difficult.
In Sect.~\ref{S:LINFINITYESTIMATESFORHIGHERTRANSVERSAL}, we 
derive related estimates involving higher transversal derivatives.
In Sects.~\ref{S:SHARPESTIMATESFORUPMU},
and \ref{S:POINTWISEESTIMATESFORWAVEEQUATIONERRORTERMS},
we use the preliminary estimates to derive related but more difficult estimates.

In the remainder of the article, 
we schematically express many of our inequalities
by stating them in terms of the arrays $\GdVar$ and $\BadVar$ from Def.~\ref{D:ABBREIVATEDVARIABLES}.
We also remind the reader that we often use the abbreviations
introduced in Subsect.~\ref{SS:STRINGSOFCOMMUTATIONVECTORFIELDS} 
to schematically indicate the structure of various derivative operators.


\subsection{Differential operator comparison estimates}
\label{SS:DIFFOPCOMPARISONESTIMATES}
We start by establishing comparison estimates for various differential operators.

\begin{lemma}[\textbf{The norm of $\ell_{t,u}-$tangent tensors can be measured via $\GeoAng$ contractions}]
	\label{L:TENSORSIZECONTROLLEDBYYCONTRACTIONS}
	Let $\xi_{\alpha_1 \cdots \alpha_n}$ be a type $\binom{0}{n}$ $\ell_{t,u}-$tangent tensor with $n \geq 1$.
	Under the data-size and bootstrap assumptions 
	of Subsects.~\ref{SS:SIZEOFTBOOT}-\ref{SS:AUXILIARYBOOTSTRAP}
	and the smallness assumptions of Subsect.~\ref{SS:SMALLNESSASSUMPTIONS}, 
	we have
	\begin{align} \label{E:TENSORSIZECONTROLLEDBYYCONTRACTIONS}
		|\xi| = 
			\left\lbrace 
				1 + \mathcal{O}(\varepsilon^{1/2})
			\right\rbrace
			|\xi_{\GeoAng \GeoAng \cdots \GeoAng}|.
	\end{align}
	The same result holds if 
	$|\xi_{\GeoAng \GeoAng \cdots \GeoAng}|$
	is replaced with 
	$|\xi_{\GeoAng \cdot}|$, 
	$|\xi_{\GeoAng \GeoAng \cdot}|$,
	etc., where $\xi_{\GeoAng \cdot}$
	is the type $\binom{0}{n-1}$ tensor with components
  $\GeoAng^{\alpha_1} \xi_{\alpha_1 \alpha_2 \cdots \alpha_n}$,
  and similarly for $\xi_{\GeoAng \GeoAng \cdot}$, etc.
\end{lemma}
\begin{proof}
	\eqref{E:TENSORSIZECONTROLLEDBYYCONTRACTIONS} is easy to derive 
	relative to rectangular coordinates by using 
	the decomposition $(\ginversesphere)^{ij} = \frac{1}{|\GeoAng|^2} \GeoAng^i \GeoAng^j$
	and the estimate $|\GeoAng| = 1 + \mathcal{O}(\varepsilon^{1/2})$.
	This latter estimate follows from
	the identity 
	$|\GeoAng|^2 
	= g_{ab} \GeoAng^a \GeoAng^b
	= (\delta_{ab} + g_{ab}^{(Small)}) (\delta_2^a + \GeoAng_{(Small)}^a)(\delta_2^b + \GeoAng_{(Small)}^b)$,
	the fact that $g_{ab}^{(Small)} = \smoothfunction(\GdVar)\GdVar$ with $\smoothfunction$ smooth
	and similarly for $\GeoAng_{(Small)}^a$
	(see Lemma~\ref{L:SCHEMATICDEPENDENCEOFMANYTENSORFIELDS}),
	and the bootstrap assumptions.
\end{proof}

\begin{lemma}[\textbf{Controlling $\angD$ derivatives in terms of $\GeoAng$ derivatives}]
\label{L:ANGDERIVATIVESINTERMSOFTANGENTIALCOMMUTATOR}
	Let $f$ be a scalar function on $\ell_{t,u}$.
	Under the data-size and bootstrap assumptions 
	of Subsects.~\ref{SS:SIZEOFTBOOT}-\ref{SS:AUXILIARYBOOTSTRAP}
	and the smallness assumptions of Subsect.~\ref{SS:SMALLNESSASSUMPTIONS}, 
	the following comparison estimates hold
	on $\mathcal{M}_{\Tboot,U_0}$:
\begin{align} \label{E:ANGDERIVATIVESINTERMSOFTANGENTIALCOMMUTATOR}
		|\angdiff f|
		& \leq (1 + C \varepsilon^{1/2})\left| \GeoAng f \right|,
		\qquad
		|\angD^2 f|
		\leq (1 + C \varepsilon^{1/2})\left| \angdiff(\GeoAng f) \right|
			+ C \varepsilon |\angdiff f|.
\end{align}
\end{lemma}

\begin{proof}
	The first inequality in \eqref{E:ANGDERIVATIVESINTERMSOFTANGENTIALCOMMUTATOR}
	follows directly from Lemma~\ref{L:TENSORSIZECONTROLLEDBYYCONTRACTIONS}.
	To prove the second, we first use
	Lemma~\ref{L:TENSORSIZECONTROLLEDBYYCONTRACTIONS},
	the identity 
	$\angD_{\GeoAng \GeoAng}^2 f = \GeoAng \cdot \angdiff (\GeoAng f) - \angD_{\GeoAng} \GeoAng \cdot \angdiff f$,
	and the estimate $|\GeoAng| = 1 + \mathcal{O}(\varepsilon^{1/2})$ noted in the proof of
	Lemma~\ref{L:TENSORSIZECONTROLLEDBYYCONTRACTIONS}
	to deduce that 
	\begin{align} \label{E:ANGDSQUAREFUNCTIONFIRSTBOUNDINTERMSOFGEOANG}
		|\angD^2 f| 
		& \leq (1 + C \varepsilon^{1/2})|\angD_{\GeoAng \GeoAng}^2 f|
			\leq (1 + C \varepsilon^{1/2})|\angdiff(\GeoAng f)|
			+ |\angD_{\GeoAng} \GeoAng||\angdiff f|.
	\end{align}
	Next, we use Lemma~\ref{L:TENSORSIZECONTROLLEDBYYCONTRACTIONS} 
	and the identity 
	$	\angdeformarg{\GeoAng}{\GeoAng}{\GeoAng} 
		= \angD_{\GeoAng} (\gsphere(\GeoAng, \GeoAng))
		= \GeoAng (g_{ab} \GeoAng^a \GeoAng^b)
	$
	to deduce that
	\begin{align} \label{E:ANGDGEOANGOFGEOANGINTERMSOFGEOANGDEFORMSPHERE}
		\left|
			\angD_{\GeoAng} \GeoAng
		\right|
		& \lesssim
			\left|
				g(\angD_{\GeoAng} \GeoAng,\GeoAng)
			\right|
			\lesssim
			\left|
				\angdeformarg{\GeoAng}{\GeoAng}{\GeoAng}
			\right|
			\lesssim
			\left|
				\GeoAng(g_{ab} \GeoAng^a \GeoAng^b)
			\right|.
	\end{align}
	Since Lemma~\ref{L:SCHEMATICDEPENDENCEOFMANYTENSORFIELDS} implies that
	$g_{ab} \GeoAng^a \GeoAng^b = \smoothfunction(\GdVar)$ with $\smoothfunction$ smooth,
	the bootstrap assumptions yield that RHS~\eqref{E:ANGDGEOANGOFGEOANGINTERMSOFGEOANGDEFORMSPHERE}
	is $\lesssim |\GeoAng \GdVar| \lesssim \varepsilon^{1/2}$.
	The desired inequality now follows from this estimate,
	\eqref{E:ANGDSQUAREFUNCTIONFIRSTBOUNDINTERMSOFGEOANG},
	and
	\eqref{E:ANGDGEOANGOFGEOANGINTERMSOFGEOANGDEFORMSPHERE}.
\end{proof}

\begin{lemma} [\textbf{Controlling $\angLie_V$ and $\angD$ derivatives in terms of $\angLie_{\GeoAng}$ derivatives}]
	\label{L:ANGLIEPXIINTERMSOFANGLIEGEOANGXI}
	Let $\xi_{\alpha_1 \cdots \alpha_n}$ be a type $\binom{0}{n}$ $\ell_{t,u}-$tangent tensor with $n \geq 1$
	and let $V$ be an $\ell_{t,u}-$tangent vectorfield.
	Under the data-size and bootstrap assumptions 
	of Subsects.~\ref{SS:SIZEOFTBOOT}-\ref{SS:AUXILIARYBOOTSTRAP}
	and the smallness assumptions of Subsect.~\ref{SS:SMALLNESSASSUMPTIONS}, 
	the following comparison estimates hold
	on $\mathcal{M}_{\Tboot,U_0}$:
	\begin{align} \label{E:ANGLIEPXIINTERMSOFANGLIEGEOANGXI}
		\left| 
			\angLie_V \xi
		\right|
		& \lesssim 
			|V| 
			 \left|
			 	\angLie_{\GeoAng} \xi
			 \right|
			+ |\xi|
				\left|
					\angLie_{\GeoAng} V
				\right|
			+ \varepsilon^{1/2} |\xi| |V|,
				\\
		|\angD \xi|
		& \lesssim
			|\angLie_{\GeoAng} \xi|
			+ \varepsilon^{1/2} |\xi|.
		\label{E:ANGDPIINTERMSOFLIEDERIVATIVES}
	\end{align}

\end{lemma}

\begin{proof}
	To prove \eqref{E:ANGLIEPXIINTERMSOFANGLIEGEOANGXI},
	we first note 
	the schematic identity
	$\angLie_V \xi 
		= \angD_V \xi + \sum \xi \cdot \angD V
	$,
	which follows from applying $\Lineproject$ to both sides of \eqref{E:LIEDERIVATIVE},
	recalling that RHS~\eqref{E:LIEDERIVATIVE} is invariant upon
	replacing all coordinate partial derivatives $\partial$ 
	with covariant derivatives $\D$, 
	and recalling that $\angD = \Lineproject \D$
	when acting on $\ell_{t,u}-$tangent tensorfields $\xi$.
	Also using Lemma~\ref{L:TENSORSIZECONTROLLEDBYYCONTRACTIONS}, we find that
	\begin{align} \label{E:FIRSTBOUNDANGLIEPXIINTERMSOFANGD}
		|\angLie_V \xi| \lesssim |V||\angD_{\GeoAng} \xi| + |\xi| |\angD_{\GeoAng} V|.
	\end{align}
	Next, we note that the torsion-free property of $\angD$ implies that
	$\angD_{\GeoAng} V 
		= \angLie_{\GeoAng} V 
		+ \angD_V \GeoAng
	$.
	Hence, using Lemma~\ref{L:TENSORSIZECONTROLLEDBYYCONTRACTIONS},
	\eqref{E:ANGDGEOANGOFGEOANGINTERMSOFGEOANGDEFORMSPHERE},
	and the estimate $|\angD_{\GeoAng} \GeoAng| \lesssim \varepsilon^{1/2}$
	shown in the proof of Lemma~\ref{L:ANGDERIVATIVESINTERMSOFTANGENTIALCOMMUTATOR},
	we find that
	\begin{align}  \label{E:ANGDGEOANGXINTERMSOFLIEDERIVATIVESBOUND}
		|\angD_{\GeoAng} V|
		& \lesssim 
			|\angLie_{\GeoAng} V|
			+ |V||\angD \GeoAng|
			\lesssim 
			|\angLie_{\GeoAng} V|
			+ |V||\angD_{\GeoAng} \GeoAng|
			\lesssim 
			|\angLie_{\GeoAng} V|
			+ \varepsilon^{1/2} |V|.
	\end{align}
	Similarly, we have
	\begin{align} \label{E:ANGDGEOANGXIINTERMSOFLIEDERIVATIVESBOUND}
		|\angD_{\GeoAng} \xi|
		& \lesssim |\angLie_{\GeoAng} \xi|
			+ \varepsilon^{1/2} |\xi|.
	\end{align}
	The desired estimate \eqref{E:ANGLIEPXIINTERMSOFANGLIEGEOANGXI}
	now follows from
	\eqref{E:FIRSTBOUNDANGLIEPXIINTERMSOFANGD},
	\eqref{E:ANGDGEOANGXINTERMSOFLIEDERIVATIVESBOUND},
	and
	\eqref{E:ANGDGEOANGXIINTERMSOFLIEDERIVATIVESBOUND}.

	The estimate \eqref{E:ANGDPIINTERMSOFLIEDERIVATIVES}
	follows from applying Lemma~\ref{L:TENSORSIZECONTROLLEDBYYCONTRACTIONS} to $\angD \xi$ and
	using \eqref{E:ANGDGEOANGXIINTERMSOFLIEDERIVATIVESBOUND}.
\end{proof}

\subsection{Basic facts and estimates that we use silently}
\label{SS:OFTENUSEDESTIMATES}
For the reader's convenience, we present here some basic facts and estimates
that we silently use throughout the rest of the paper when deriving
estimates.

\begin{enumerate}
	\item All quantities that we estimate can be controlled in terms of the small quantities
		$\BadVar = \lbrace \Psi, \upmu - 1, \Lunit_{(Small)}^1, \Lunit_{(Small)}^2 \rbrace$
		and their derivatives (where the $\Rad$ derivatives do not have to be small,
		nor does $\Lunit \upmu$).
	\item We typically use the Leibniz rule for 
		the operators $\angLie_Z$ and $\angD$ when deriving
		pointwise estimates for the $\angLie_Z$ and $\angD$ 
		derivatives of tensor products of the schematic form 
		$\prod_{i=1}^m v_i$, where the $v_i$ are scalar functions or
		$\ell_{t,u}-$tangent tensors. Our derivative counts are such that
		all $v_i$ except at most one are uniformly bounded in $L^{\infty}$
		on $\mathcal{M}_{\Tboot,U_0}$.
		Thus, our pointwise estimates often explicitly feature 
		(on the right-hand sides)
		only the factor with the most derivatives on it,
		multiplied by a constant that uniformly bounds the other factors.
		In some estimates, the right-hand sides also gain a smallness factor,
		such as $\varepsilon^{1/2}$,
		generated by the remaining $v_i's$.
	\item The operators $\angLie_{\Fullset}^N$ commute through
		$\angdiff$, as shown by Lemma~\ref{L:LANDRADCOMMUTEWITHANGDIFF}.
	\item As differential operators acting on scalar functions, we have
		$\GeoAng 
		= \left(1 + \mathcal{O}(\GdVar)\right) \angdiff
		= (1 + \mathcal{O}(\varepsilon^{1/2})) \angdiff
		$,
		a fact which follows from 
		the proof of Lemma~\ref{L:ANGDERIVATIVESINTERMSOFTANGENTIALCOMMUTATOR},
		\eqref{E:GEOANGPOINTWISE}, and the bootstrap assumptions.
		Hence, for scalar functions $f$, we sometimes schematically depict $\angdiff f$ 
		as 
		$\left(1 + \mathcal{O}(\GdVar)\right) \Singletan f$
		or
		$\Singletan f$
		when the factor 
		$1 + \mathcal{O}(\GdVar)$ is not important.
		Similarly, the proofs of
		Lemmas~\ref{L:ANGDERIVATIVESINTERMSOFTANGENTIALCOMMUTATOR} 
		and \ref{L:ANGLIEPXIINTERMSOFANGLIEGEOANGXI} show that
		we can depict $\angLap f$ by 
		$
		\smoothfunction(\Tanset^{\leq 1} \GdVar,\ginversesphere)
		\Tanset_*^{[1,2]} f
		$
		(or $\Tanset_*^{[1,2]} f$
		when the factor $\smoothfunction(\Tanset^{\leq 1} \GdVar,\ginversesphere)$
		is not important)
		and, 
		for type $\binom{0}{n}$ $\ell_{t,u}-$tangent tensorfields $\xi$,
		$\angD \xi$ by 
		$
		\smoothfunction(\Tanset^{\leq 1} \GdVar,\ginversesphere)
		\angLie_{\Tanset}^{\leq 1} \xi
		$
		(or $\angLie_{\Tanset}^{\leq 1} \xi$
		when the factor $\smoothfunction(\Tanset^{\leq 1} \GdVar,\ginversesphere)$
		is not important).
	\item We remind the reader that all constants are allowed to depend on 
		the data-size parameters
		$\mathring{\updelta}$
		and 
		$\TranminusdatasizeWithFactor^{-1}$.
\end{enumerate}

\subsection{Pointwise estimates for the rectangular coordinates and the rectangular components of some vectorfields}

\begin{lemma}[\textbf{Pointwise estimates for} $x^i$ \textbf{and the rectangular components of several vectorfields}]
	\label{L:POINTWISEFORRECTANGULARCOMPONENTSOFVECTORFIELDS}
	Assume that $N \leq 18$
	and
	$V \in \lbrace \Lunit, \Radunit, \GeoAng \rbrace$.
	Let $x^i = x^i(t,u,\vartheta)$ denote the rectangular coordinate function
	and let $\mathring{x}^i = \mathring{x}^i(u,\vartheta) := x^i(0,u,\vartheta)$.
	Under the data-size and bootstrap assumptions 
	of Subsects.~\ref{SS:SIZEOFTBOOT}-\ref{SS:AUXILIARYBOOTSTRAP}
	and the smallness assumptions of Subsect.~\ref{SS:SMALLNESSASSUMPTIONS}, 
	the following pointwise estimates hold
	on $\mathcal{M}_{\Tboot,U_0}$,
	for $i = 1,2$
	(see Subsect.~\ref{SS:STRINGSOFCOMMUTATIONVECTORFIELDS} regarding the vectorfield operator notation):
	\begin{subequations}
	\begin{align}
		\left| 
			V^i
		\right|
		& \lesssim 
			1
			+
			\left|
				\GdVar
			\right|,
		\label{E:ORDERZEROVECTORFIELDSRECTCOMPPOINTWISE} \\		 
		\left| 
			\Tanset^{[1,N]} V^i
		\right|
		& \lesssim 
			\left|
				\Tanset^{\leq N} \GdVar
			\right|,
		\label{E:TANGENTIALCOMMMUTEDVECTORFIELDSRECTCOMPPOINTWISE} \\
		\left| 
			\Fullset_*^{[1,N];1} V^i
		\right|
		& \lesssim 
			\left|
				\Fullset_*^{\leq N;1} \GdVar
			\right|,
		\label{E:ONERADIALNOTPURERADIALCOMMMUTEDVECTORFIELDSRECTCOMPPOINTWISE} \\
		\left| 
			\Fullset^{[1,N];1} V^i
		\right|
		& \lesssim 
			\left|
				\Fullset^{\leq N;1} \GdVar
			\right|,
		\label{E:ONERADIALCOMMMUTEDVECTORFIELDSRECTCOMPPOINTWISE} \\
		\left| 
			\Rad^i
		\right|
		& \lesssim
			1 + |\BadVar|,
			\label{E:RADRECTCOMPPOINTWISE} \\
		\left| 
			\Tanset^{[1,N]} \Rad^i
		\right|
		& \lesssim
			\left|
				\Tanset^{\leq N}
				\BadVar
			\right|,
			\label{E:TANGENTIALCOMMUTEDRADIPOINTWISE} 
				\\
		\left| 
			\Fullset_*^{[1,N];1} \Rad^i
		\right|
		& \lesssim
			\left|
				\Fullset_*^{\leq N;1} \BadVar
			\right|,
			\label{E:ONERADIALNOTPURERADIALCOMMUTEDRADIPOINTWISE} 
				\\
		\left| 
			\Fullset^{[1,N];1} \Rad^i
		\right|
		& \lesssim
			\left|
				\Fullset^{\leq N;1} \BadVar
			\right|,
			\label{E:ONERADIALCOMMUTEDRADIPOINTWISE} 
	\end{align}
	\end{subequations}

	\begin{subequations}
		\begin{align} \label{E:XIPOINTWISE}
			\left|
				x^i - \mathring{x}^i
			\right|
			& \lesssim 1,
				\\
			\left|
				\angdiff x^i
			\right|
			& \lesssim 
				1 
				+
				\left| 
					\GdVar 
				\right|,
				 \label{E:ANGDIFFXI} \\
			\left|
				\angdiff \Tanset^{[1,N]} x^i
			\right|
			& \lesssim 
				\left| 
					\Tanset^{\leq N} \GdVar 
				\right|,
					\label{E:ANGDIFFXIPURETANGENTIALDIFFERENTIATED} \\
			\left|
				\angdiff \Fullset_*^{[1,N];1} x^i
			\right|
			& \lesssim 
				\left| 
					\Fullset_*^{\leq N;1} \GdVar 
				\right|
				+
				\left| 
					\Tanset_*^{[1,N]} \BadVar
				\right|,
					\label{E:ANGDIFFXIONERADIALNOTPURERADIALDIFFERENTIATED} \\
			\left|
				\angdiff \Fullset^{[1,N];1} x^i
			\right|
			& \lesssim 
				\left| 
					\Tanset^{\leq N} \GdVar
				\right|
				+
				\left| 
					\Tanset_*^{[1,N]} \BadVar
				\right|,
					\label{E:ANGDIFFXIONERADIALDIFFERENTIATED} 
		\end{align} 
	\end{subequations}

	\begin{subequations}
	\begin{align} \label{E:PURETANGENTIALGEOANGSMALLIPOINTWISE}
	\left|
		\Tanset^N \GeoAng_{(Small)}^i
	\right|
	& \lesssim 
		\left|
			\Tanset^{\leq N} \GdVar
		\right|,
			\\
	\left|
		\Fullset_*^{N;1} \GeoAng_{(Small)}^i
	\right|
	& \lesssim 
		\left|
			\Fullset_*^{\leq N;1} \GdVar
		\right|,
		\label{E:ONERADIALNOTPURERADIALGEOANGSMALLIPOINTWISE}
			\\
		\left|
			\Fullset^{N;1} \GeoAng_{(Small)}^i
		\right|
	& \lesssim 
		\left|
			\Fullset^{\leq N;1} \GdVar
		\right|.
			\label{E:ONERADIALGEOANGSMALLIPOINTWISE} 
	\end{align}
	\end{subequations}
	In the case $i=2$ at fixed $u,\vartheta$,
	LHS~\eqref{E:XIPOINTWISE} is to be interpreted as
	the Euclidean distance traveled
	by the point $x^2$
	in the flat universal covering space 
	$\mathbb{R}$ of $\mathbb{T}$
	along the corresponding integral curve of $\Lunit$
	over the time interval $[0,t]$.
\end{lemma}

\begin{proof}
See Subsect.~\ref{SS:OFTENUSEDESTIMATES} for some comments on the analysis.
Lemma~\ref{L:SCHEMATICDEPENDENCEOFMANYTENSORFIELDS} implies
that for $V \in \lbrace \Lunit, \Radunit, \GeoAng \rbrace$,
the component $V^i = V x^i$ verifies $V^i = \smoothfunction(\GdVar)$
with $\smoothfunction$ smooth.
Similarly, $\GeoAng_{(Small)}^i$ verifies
$\GeoAng_{(Small)}^i = \smoothfunction(\GdVar)\GdVar$
with $\smoothfunction$ smooth
and
$\Rad x^i = \Rad^i$ verifies $\Rad^i = \smoothfunction(\BadVar)$
with $\smoothfunction$ smooth.
The estimates of the lemma therefore follow easily
from the bootstrap assumptions,
except for the estimates \eqref{E:XIPOINTWISE}-\eqref{E:ANGDIFFXIONERADIALDIFFERENTIATED}.
To obtain \eqref{E:XIPOINTWISE}, 
we first argue as above to deduce $|\Lunit x^i| = |\Lunit^i| = |\smoothfunction(\GdVar)| \lesssim 1$.
Since $\Lunit = \frac{\partial}{\partial t}$, 
we may integrate along the integral curves of $\Lunit$ starting from $t = 0$
and use the previous estimate to conclude \eqref{E:XIPOINTWISE}.
To derive \eqref{E:ANGDIFFXI},
we use \eqref{E:ANGDERIVATIVESINTERMSOFTANGENTIALCOMMUTATOR}
with $f=x^i$ to deduce 
$|\angdiff x^i| 
\lesssim |\GeoAng x^i| 
= |\GeoAng^i| 
= |\smoothfunction(\GdVar)| 
\lesssim 
1 
+
|\GdVar|
$
as desired.
The proofs of 
\eqref{E:ANGDIFFXIPURETANGENTIALDIFFERENTIATED}-\eqref{E:ANGDIFFXIONERADIALDIFFERENTIATED}
are similar, but we also use
Lemma~\ref{L:LANDRADCOMMUTEWITHANGDIFF}
to commute vectorfields under $\angdiff$.
\end{proof}



\subsection{Pointwise estimates for various \texorpdfstring{$\ell_{t,u}-$}{}tensorfields}

\begin{lemma}[\textbf{Crude pointwise estimates for the Lie derivatives of $\gsphere$ and $\ginversesphere$}]
\label{L:POINTWISEESTIMATESFORGSPHEREANDITSDERIVATIVES}
	Assume that $N \leq 18$.
	Under the data-size and bootstrap assumptions 
	of Subsects.~\ref{SS:SIZEOFTBOOT}-\ref{SS:AUXILIARYBOOTSTRAP}
	and the smallness assumptions of Subsect.~\ref{SS:SMALLNESSASSUMPTIONS}, 
	the following pointwise estimates hold
	on $\mathcal{M}_{\Tboot,U_0}$
	(see Subsect.~\ref{SS:STRINGSOFCOMMUTATIONVECTORFIELDS} regarding the vectorfield operator notation):
	\begin{subequations}
	\begin{align} \label{E:POINTWISEESTIMATESFORGSPHEREANDITSTANGENTIALDERIVATIVES}
		\left|
			\angLie_{\Tanset}^{N+1} \gsphere
		\right|,
			\,
		\left|
			\angLie_{\Tanset}^{N+1} \ginversesphere
		\right|
			\,
		\left|
			\angLie_{\Tanset}^N \upchi
		\right|,
			\,
		\left|
			\Tanset^N \mytr \upchi
		\right|
		& \lesssim 
			\left| 
				\Tanset^{\leq N+1} \GdVar
			\right|,
				\\
		\left|
			\angLie_{\Fullset_*}^{N+1;1} \gsphere
		\right|,
			\,
		\left|
			\angLie_{\Fullset_*}^{N+1;1} \ginversesphere
		\right|,
			\,
		\left|
			\angLie_{\Fullset}^{N;1} \upchi
		\right|,
			\,
		\left|
			\Fullset^{N;1} \mytr \upchi
		\right|
		& \lesssim 
			\left|
				\Fullset_*^{\leq N+1;1} \GdVar
			\right|
			+
			\left|
				\Tanset_*^{[1, N+1]} \BadVar
			\right|,
			\label{E:ONERADIALNOTPURERADIALFORGSPHERE}
				\\
		\left|
			\angLie_{\Fullset}^{N+1;1} \gsphere
		\right|,
			\,
		\left|
			\angLie_{\Fullset}^{N+1;1} \ginversesphere
		\right|
		& \lesssim 
			\left|
				\Fullset^{\leq N+1;1} \GdVar
			\right|
			+
			\left| 
				\Tanset_*^{[1, N+1]} \BadVar
			\right|.
			\label{E:ONERADIALFORGSPHERE}
\end{align}
\end{subequations}
\end{lemma}

\begin{proof}
		See Subsect.~\ref{SS:OFTENUSEDESTIMATES} for some comments on the analysis.
		By Lemma~\ref{L:SCHEMATICDEPENDENCEOFMANYTENSORFIELDS}, we have
		$\gsphere = \smoothfunction(\GdVar,\angdiff x^1,\angdiff x^2)$.
		The desired estimates
		for $\angLie_{\Tanset}^{N+1} \gsphere$
		thus follow from
		Lemma~\ref{L:POINTWISEFORRECTANGULARCOMPONENTSOFVECTORFIELDS}
		and the bootstrap assumptions.
		The desired estimates for 
		$\angLie_{\Tanset}^{N+1} \ginversesphere$
		then follow from repeated use of the second
		identity in \eqref{E:CONNECTIONBETWEENANGLIEOFGSPHEREANDDEFORMATIONTENSORS}
		and the estimates for $\angLie_{\Tanset}^{N+1} \gsphere$.
		The estimates for $\angLie_{\Tanset}^N \upchi$ 
		and
		$\Tanset^N \mytr \upchi$
		follow from the estimates for $\angLie_{\Tanset}^{N+1} \gsphere$
		and $\angLie_{\Tanset}^{N+1} \ginversesphere$
		since $\upchi \sim \angLie_{\Singletan} \gsphere$
		(see \eqref{E:CHIDEF})
		and $\mytr \upchi \sim \ginversesphere \cdot \angLie_{\Singletan} \gsphere$.
\end{proof}

\begin{lemma}[\textbf{Pointwise estimates for the Lie derivatives of $\GeoAng$ and some deformation tensor components}]
\label{L:MOREPRECISEANGSPHERELESTIMATES}
Assume that $N \leq 18$.
Under the data-size and bootstrap assumptions 
of Subsects.~\ref{SS:SIZEOFTBOOT}-\ref{SS:AUXILIARYBOOTSTRAP}
and the smallness assumptions of Subsect.~\ref{SS:SMALLNESSASSUMPTIONS}, 
the following pointwise estimates hold
on $\mathcal{M}_{\Tboot,U_0}$
(see Subsect.~\ref{SS:STRINGSOFCOMMUTATIONVECTORFIELDS} regarding the vectorfield operator notation):
\begin{subequations}
		\begin{align} \label{E:GEOANGPOINTWISE}
		\left|
			\GeoAng
		\right|
		& \leq 1
			+
			C
			\left|
				\GdVar
			\right|,
				\\
		\left|
			\angLie_{\Tanset}^{[1,N]} \GeoAng
		\right|
		& \leq 
			C
			\left|
				\Tanset^{\leq N} \GdVar
			\right|,
				\label{E:TANGENTIALDIFFERENTIATEDGEOANGPOINTWISE} 
				\\
	\left|
			\angLie_{\Fullset_*}^{[1,N];1} \GeoAng
		\right|
		& \leq
			C
			\left|
				\Fullset_*^{\leq N;1} \GdVar
			\right|
			+ 
			C
			\left|
				\Tanset_*^{[1,N]} \BadVar
			\right|,
				\label{E:ONERADIALANDTANGENTIALDIFFERENTIATEDGEOANGPOINTWISE}
				\\
		\left|
			\angLie_{\Fullset}^{[1,N];1} \GeoAng
		\right|
		& \leq
			C
			\left|
				\Fullset^{\leq N;1} \GdVar
			\right|
			+ 
			C
			\left|
				\Tanset_*^{[1,N]} \BadVar
			\right|,
				\label{E:ONERADIALANDDIFFERENTIATEDGEOANGPOINTWISE}
	\end{align}
	\end{subequations}

\begin{subequations}
\begin{align}
	\left|
		\angLie_{\Tanset}^N \angdeformoneformarg{\GeoAng}{\Lunit}
	\right|,
		\,
	\left|
		\angLie_{\Tanset}^N \angdeformoneformupsharparg{\GeoAng}{\Lunit}
	\right|
	& \lesssim 
		\left|
			\Tanset^{\leq N+1} \GdVar
		\right|,
				\label{E:TANGENTDIFFERNTIATEDGEOANGDEFORMSPHERELSHARPPOINTWISE} \\
	\left|
		\angLie_{\Fullset}^{N;1} \angdeformoneformarg{\GeoAng}{\Lunit}
	\right|,
		\,
	\left|
		\angLie_{\Fullset}^{N;1} \angdeformoneformupsharparg{\GeoAng}{\Lunit}
	\right|
	& \lesssim 
		\left|
			\Fullset_*^{\leq N+1;1} \GdVar
		\right|
		+
		\varepsilon^{1/2}
		\left|
			\Fullset^{\leq N;1} \GdVar
		\right|
		+
		\varepsilon^{1/2}
		\left|
			\Tanset_*^{[1,N]} \BadVar
		\right|,
				\label{E:ONERADIALTANGENTDIFFERNTIATEDGEOANGDEFORMSPHERELSHARPPOINTWISE} 
	\end{align}
\end{subequations}

\begin{align}
	\left|
		\angLie_{\Tanset}^N \angdeformoneformarg{\GeoAng}{\Rad}
	\right|,
		\,
	\left|
		\angLie_{\Tanset}^N \angdeformoneformupsharparg{\GeoAng}{\Rad}
	\right|
	& \lesssim 
		\left|
			\Fullset_*^{\leq N+1;1} \Psi
		\right|
		+ 
		\left|
			\Tanset^{\leq N+1} \GdVar
		\right|
		+ 
		\left|
			\Tanset_*^{[1,N]} \BadVar
		\right|,
				\label{E:TANGENTDIFFERNTIATEDGEOANGDEFORMSPHERERADPOINTWISE}
\end{align}

\begin{subequations}
\begin{align}
	\left|
		\angdeformoneformarg{\Rad}{\Lunit}
	\right|,
		\,
	\left|
		\angdeformoneformupsharparg{\Rad}{\Lunit}
	\right|
	& \lesssim 
		\left|
			\Fullset^{\leq 1} \Psi
		\right|
		+ 
		\left|
			\Tanset_* \upmu
		\right|,
			\label{E:RADDEFORMSPHERELPOINTWISE}
			\\
	\left|
		\angLie_{\Tanset}^{[1,N]} \angdeformoneformarg{\Rad}{\Lunit}
	\right|,
		\,
	\left|
		\angLie_{\Tanset}^{[1,N]} \angdeformoneformupsharparg{\Rad}{\Lunit}
	\right|
	& \lesssim 
		\left|
			\Fullset_*^{\leq N+1;1} \Psi
		\right|
		+ 
		\left|
			\Tanset^{\leq N+1} \GdVar
		\right|
		+ 
		\left|
			\Tanset_*^{[1, N+1]} \BadVar
		\right|,
			\label{E:TANGENTDIFFERNTIATEDRADDEFORMSPHERELPOINTWISE}
	\end{align}
\end{subequations}

\begin{subequations}
\begin{align} \label{E:TANGENTDIFFERNTIATEDANGDEFORMTANGENTPOINTWISE}
	\left|
		\angLie_{\Tanset}^N \angdeform{\Lunit}
	\right|,
		\,
	\left|
		\angLie_{\Tanset}^N \angdeform{\GeoAng}
	\right|
	& \lesssim 
		\left|
			\Tanset^{\leq N+1} \GdVar
		\right|,
			\\
	\left|
		\angLie_{\Fullset}^{N;1} \angdeform{\Lunit}
	\right|,
		\,
	\left|
		\angLie_{\Fullset}^{N;1} \angdeform{\GeoAng}
	\right|
	& \lesssim 
	  \left|
			\Fullset_*^{\leq N+1;1} \GdVar
		\right|
		+
		\varepsilon^{1/2}
		\left|
			\Fullset^{\leq N;1} \GdVar
		\right|
		+
		\varepsilon^{1/2}
		\left|
			\Tanset_*^{[1,N]} \BadVar
		\right|.
		\label{E:ONERADIALDIFFERNTIATEDANGDEFORMTANGENTPOINTWISE}
\end{align}
\end{subequations}

\end{lemma}

\begin{proof}
	See Subsect.~\ref{SS:OFTENUSEDESTIMATES} for some comments on the analysis.
	To prove \eqref{E:TANGENTDIFFERNTIATEDGEOANGDEFORMSPHERELSHARPPOINTWISE}
	for $\angdeformoneformupsharparg{\GeoAng}{\Lunit}$,
	we first note that by 
	Lemma~\ref{L:SCHEMATICDEPENDENCEOFMANYTENSORFIELDS}
	and
	\eqref{E:GEOANGDEFORMSPHEREL},
	we have
	$\angdeformoneformupsharparg{\GeoAng}{\Lunit}
	=
	\smoothfunction(\GdVar,\ginversesphere,\angdiff x^1,\angdiff x^2) \Singletan \GdVar
	$.
	We now apply $\angLie_{\Tanset}^N$ to the previous relation.
	We bound the derivatives of $\ginversesphere$ and $\angdiff x$
	with Lemmas~\ref{L:POINTWISEFORRECTANGULARCOMPONENTSOFVECTORFIELDS}
	and \ref{L:POINTWISEESTIMATESFORGSPHEREANDITSDERIVATIVES}.
	Also using the bootstrap assumptions, we conclude the desired result.
	A similar argument yields the same estimate for
	$\angdeformoneformarg{\GeoAng}{\Lunit}$.
	The proof of \eqref{E:ONERADIALTANGENTDIFFERNTIATEDGEOANGDEFORMSPHERELSHARPPOINTWISE} 
	is similar and we omit the details.

	Since Lemma~\ref{L:SCHEMATICDEPENDENCEOFMANYTENSORFIELDS} implies that
	$\GeoAng = \smoothfunction(\GdVar,\ginversesphere,\angdiff x^1,\angdiff x^2)$,
	similar reasoning yields 
	\eqref{E:ONERADIALANDTANGENTIALDIFFERENTIATEDGEOANGPOINTWISE}-\eqref{E:ONERADIALANDDIFFERENTIATEDGEOANGPOINTWISE}.

	Inequality \eqref{E:GEOANGPOINTWISE}
	follows from the slightly more
	precise arguments already given in the proof of Lemma~\ref{L:TENSORSIZECONTROLLEDBYYCONTRACTIONS}.

	The proof of \eqref{E:TANGENTDIFFERNTIATEDGEOANGDEFORMSPHERERADPOINTWISE}
	is similar and is based on the observation that by
	Lemma~\ref{L:SCHEMATICDEPENDENCEOFMANYTENSORFIELDS}
	and
	\eqref{E:GEOANGDEFORMSPHERERAD},
	we have
	\[
	\angdeformoneformupsharparg{\GeoAng}{\Rad}
	=
	\smoothfunction(\BadVar,\ginversesphere,\angdiff x^1,\angdiff x^2) \Singletan \GdVar
	+ \smoothfunction(\GdVar,\ginversesphere,\angdiff x^1,\angdiff x^2,\Rad \Psi) \GdVar
	+ \smoothfunction(\GdVar,\ginversesphere) \angdiff \upmu,
	\]
	and a similar schematic relation holds for 
	$\angdeformoneformarg{\GeoAng}{\Rad}$.

	The proofs of \eqref{E:RADDEFORMSPHERELPOINTWISE}-\eqref{E:TANGENTDIFFERNTIATEDRADDEFORMSPHERELPOINTWISE}
	are similar and are based on the observation that by
	Lemma~\ref{L:SCHEMATICDEPENDENCEOFMANYTENSORFIELDS},
	\eqref{E:ZETADECOMPOSED},
	and 
	\eqref{E:RADDEFORMSPHERERAD},
	we have
	\[
	\angdeformoneformupsharparg{\Rad}{\Lunit}
	=
	\smoothfunction(\GdVar,\ginversesphere,\angdiff x^1,\angdiff x^2) \Rad \Psi
	+
	\smoothfunction(\BadVar,\ginversesphere,\angdiff x^1,\angdiff x^2) \Singletan \Psi
	+ \ginversesphere \angdiff \upmu,
	\]
	and a similar schematic relation holds for 
	$\angdeformoneformarg{\Rad}{\Lunit}$.

	The proofs of 
	\eqref{E:TANGENTDIFFERNTIATEDANGDEFORMTANGENTPOINTWISE}-\eqref{E:ONERADIALDIFFERNTIATEDANGDEFORMTANGENTPOINTWISE} 
	are similar and are based on the fact that by
	Lemma~\ref{L:SCHEMATICDEPENDENCEOFMANYTENSORFIELDS},
	\eqref{E:LUNITDEFORMSPHERE},
	and
	\eqref{E:GEOANGDEFORMSPHERE},
	we have
	$\angdeform{\Lunit},
	\angdeform{\GeoAng}
	= \smoothfunction(\GdVar,\ginversesphere,\angdiff x^1,\angdiff x^2) \Singletan \GdVar
	$.

\end{proof}

\subsection{Commutator estimates}
\label{SS:COMMUTATORESTIMATES}
In this subsection, we establish some commutator estimates.

\begin{lemma}[\textbf{Pure $\mathcal{P}_u-$tangent commutator estimates}]
	\label{L:COMMUTATORESTIMATES}
	Assume that $1 \leq N \leq 18$.
	Let $\vec{I}$ be an order $|\vec{I}| = N + 1$
	multi-index for the set $\Tanset$ of 
	$\mathcal{P}_u-$tangent commutation vectorfields
	(see Def.~\ref{D:REPEATEDDIFFERENTIATIONSHORTHAND}),
	and let $\vec{I}'$ be any permutation of $\vec{I}$.
	Let $f$ be a scalar function, and let
	$\xi$ be an $\ell_{t,u}-$tangent one-form or a type $\binom{0}{2}$ $\ell_{t,u}-$tangent tensorfield.
	Under the data-size and bootstrap assumptions 
	of Subsects.~\ref{SS:SIZEOFTBOOT}-\ref{SS:AUXILIARYBOOTSTRAP}
	and the smallness assumptions of Subsect.~\ref{SS:SMALLNESSASSUMPTIONS}, 
	the following commutator estimates hold
	on $\mathcal{M}_{\Tboot,U_0}$
	(see Subsect.~\ref{SS:STRINGSOFCOMMUTATIONVECTORFIELDS} regarding the vectorfield operator notation):
	\begin{subequations}
	\begin{align}
		\left|
			\Tanset^{\vec{I}} f
			-
			\Tanset^{\vec{I}'} f
		\right|
		& \lesssim
			\varepsilon^{1/2}
			\left|
				\Tanset_*^{[1,N]} f
			\right|
			+
			\left|
				\Tanset_*^{[1,\lfloor N/2 \rfloor]} f
			\right|
			\left|
				\Tanset^{\leq N} \GdVar
			\right|.
			\label{E:PURETANGENTIALFUNCTIONCOMMUTATORESTIMATE}  
		\end{align}
		\end{subequations}

		Moreover, if $1 \leq N \leq 17$
		and $\vec{I}$ is as above,
		then the following commutator estimates hold:
		\begin{subequations}
		\begin{align}
		\left|
			[\angD^2, \Tanset^N] f
		\right|
		& \lesssim
			\varepsilon^{1/2}
			\left|
				\Tanset_*^{[1,N]} f
			\right|
			+ 
			\left|
				\Tanset_*^{[1,\lceil N/2 \rceil]} f
			\right|
			\left|
				\Tanset^{\leq N+1} \GdVar
			\right|,
				\label{E:ANGDSQUAREDPURETANGENTIALFUNCTIONCOMMUTATOR} \\
		\left|
			[\angLap, \Tanset^N] f
		\right|
		& \lesssim
			\varepsilon^{1/2}
			\left|
				\Tanset_*^{[1, N+1]} f
			\right|
			+ 
			\left|
				\Tanset_*^{[1,\lceil N/2 \rceil]} f
			\right|
			\left|
				\Tanset^{\leq N+1} \GdVar
			\right|,
			\label{E:ANGLAPPURETANGENTIALFUNCTIONCOMMUTATOR}
		\end{align}
		\end{subequations}

		\begin{subequations}
		\begin{align}
		\left|
			\angLie_{\Tanset}^{\vec{I}} \xi
			-
			\angLie_{\Tanset}^{\vec{I}'} \xi
		\right|
		& \lesssim
			\varepsilon^{1/2}
			\left|
				\angLie_{\Tanset}^{[1,N]} \xi
			\right|
			+
			\left|
				\angLie_{\Tanset}^{\leq \lfloor N/2 \rfloor} \xi
			\right|
			\left|
				\Tanset^{\leq N+1} \GdVar
			\right|,
			\label{E:PURETANGENTIALTENSORFIELDCOMMUTATORESTIMATE}
		\\
		\left|
			[\angD, \angLie_{\Tanset}^N] \xi
		\right|
		& \lesssim
			\underbrace{
			\varepsilon^{1/2}
			\left|
				\angLie_{\Tanset}^{\leq [1,N-1]} \xi
			\right|
			}_{\mbox{Absent if $N=1$}}
			+
			\left|
				\angLie_{\Tanset}^{\leq \lfloor N/2 \rfloor} \xi
			\right|
			\left|
				\Tanset^{\leq N+1} \GdVar
			\right|,
				\label{E:ANGDANGLIETANGENTIALTENSORFIELDCOMMUTATORESTIMATE} \\
		\left|
			[\angdiv, \angLie_{\Tanset}^N] \xi
		\right|
		& \lesssim
			\varepsilon^{1/2}
			\left|
				\angLie_{\Tanset}^{[1,N]} \xi
			\right|
			+
			\left|
				\angLie_{\Tanset}^{\leq \lfloor N/2 \rfloor} \xi
			\right|
			\left|
				\Tanset^{\leq N+1} \GdVar
			\right|.
			\label{E:ANGDIVANGLIETANGENTIALTENSORFIELDCOMMUTATORESTIMATE}
		\end{align}
		\end{subequations}

		Finally, 
		if $1 \leq N \leq 17$,
		then we have the following alternate version of
		\eqref{E:ANGDSQUAREDPURETANGENTIALFUNCTIONCOMMUTATOR}:
		\begin{align}	 \label{E:ALTERNATEANGDSQUAREDPURETANGENTIALFUNCTIONCOMMUTATOR}
		\left|
			[\angD^2, \Tanset^N] f
		\right|
		& \lesssim
			\left|
				\Tanset^{\leq \lceil N/2 \rceil +1} \GdVar
			\right|
			\left|
				\Tanset_*^{[1,N]} f
			\right|
			+ 
			\left|
				\Tanset_*^{[1,\lceil N/2 \rceil]} f
			\right|
			\left|
				\Tanset^{\leq N+1} \GdVar
			\right|.
	\end{align}

\end{lemma}

\begin{proof}
	See Subsect.~\ref{SS:OFTENUSEDESTIMATES} for some comments on the analysis.
	We first prove \eqref{E:PURETANGENTIALFUNCTIONCOMMUTATORESTIMATE}. 
	Using \eqref{E:ALGEBRAICSUBTRACTINGTWOPERMUTEDORDERNVECTORFIELDS}
	and Lemma~\ref{L:CONNECTIONBETWEENCOMMUTATORSANDDEFORMATIONTENSORS}, 
	we see that it suffices to bound
	\begin{align} \label{E:FIRSTESTIMATEPURETANGENTIALFUNCTIONCOMMUTATORESTIMATE}
		\sum_{N_1 + N_2 \leq N-1}
					\left|
						\angLie_{\Tanset}^{N_1} \angdeformoneformupsharparg{\GeoAng}{\Lunit}
					\right|
					\left|
						\GeoAng \Tanset^{N_2} f
					\right|.
	\end{align}
	The desired bound 
	of \eqref{E:FIRSTESTIMATEPURETANGENTIALFUNCTIONCOMMUTATORESTIMATE} by $\lesssim$ RHS 
	\eqref{E:PURETANGENTIALFUNCTIONCOMMUTATORESTIMATE}
	now follows easily from \eqref{E:TANGENTDIFFERNTIATEDGEOANGDEFORMSPHERELSHARPPOINTWISE} 
	and the bootstrap assumptions.

	The proof of \eqref{E:PURETANGENTIALTENSORFIELDCOMMUTATORESTIMATE}
	is similar
	but relies on \eqref{E:TENSORFIELDACTINGALGEBRAICSUBTRACTINGTWOPERMUTEDORDERNVECTORFIELDS}
	in place of
	\eqref{E:ALGEBRAICSUBTRACTINGTWOPERMUTEDORDERNVECTORFIELDS}
	and also
	\eqref{E:ANGLIEPXIINTERMSOFANGLIEGEOANGXI}
	with $V:= \angdeformoneformupsharparg{Z_{\iota_{k_2}}}{Z_{\iota_{k_1}}}$
	(to handle the first Lie derivative operator on RHS~\eqref{E:TENSORFIELDACTINGALGEBRAICSUBTRACTINGTWOPERMUTEDORDERNVECTORFIELDS}).

	The proofs of 
	\eqref{E:ANGDSQUAREDPURETANGENTIALFUNCTIONCOMMUTATOR},
	\eqref{E:ANGLAPPURETANGENTIALFUNCTIONCOMMUTATOR},
	and \eqref{E:ALTERNATEANGDSQUAREDPURETANGENTIALFUNCTIONCOMMUTATOR}
	are similar
	and are based on the commutation identities 
	\eqref{E:COMMUTINGANGDSQUAREDANDLIEZ}-\eqref{E:COMMUTINGANGANGLAPANDLIEZ},
	the identity \eqref{E:CONNECTIONBETWEENANGLIEOFGSPHEREANDDEFORMATIONTENSORS},
	and the estimate \eqref{E:TANGENTDIFFERNTIATEDANGDEFORMTANGENTPOINTWISE}.

	The proofs of 
	\eqref{E:ANGDANGLIETANGENTIALTENSORFIELDCOMMUTATORESTIMATE}-\eqref{E:ANGDIVANGLIETANGENTIALTENSORFIELDCOMMUTATORESTIMATE}
	are similar and are based on the commutation identities
	\eqref{E:ANGDANGLIEZELLTUTENSORFIELDCOMMUTATOR}-\eqref{E:ANGDIVANGLIEZELLTUTENSORFIELDCOMMUTATOR}
\end{proof}

\begin{lemma}[\textbf{Mixed $\mathcal{P}_u-$transversal-tangent commutator estimates}]
		\label{L:TRANSVERALTANGENTIALCOMMUTATOR}
		Assume that $1 \leq N \leq 18$.
		Let $\Fullset^{\vec{I}}$ be a 
		$\Fullset-$multi-indexed operator containing
		\textbf{exactly one} $\Rad$ factor, and assume that
		$|\vec{I}| := N+1$. 
		Let $\vec{I}'$ be any permutation of $\vec{I}$.
		Let $f$ be a scalar function.
		Under the data-size and bootstrap assumptions 
		of Subsects.~\ref{SS:SIZEOFTBOOT}-\ref{SS:AUXILIARYBOOTSTRAP}
		and the smallness assumptions of Subsect.~\ref{SS:SMALLNESSASSUMPTIONS}, 
		the following commutator estimates hold
		on $\mathcal{M}_{\Tboot,U_0}$
		(see Subsect.~\ref{SS:STRINGSOFCOMMUTATIONVECTORFIELDS} regarding the vectorfield operator notation):
		\begin{align}
		\left|
			\Fullset^{\vec{I}} f
			-
			\Fullset^{\vec{I}'} f
		\right|
		& \lesssim
			\left|
				\Tanset_*^{[1,N]} f
			\right|
			+
			\varepsilon^{1/2}
			\left|
				\GeoAng \Fullset^{\leq N-1;1} f
			\right|
				\label{E:ONERADIALTANGENTIALFUNCTIONCOMMUTATORESTIMATE} \\
	& \ \
			+
			\left|
				\Tanset_*^{[1,\lfloor N/2 \rfloor]} f
			\right|
			\left|
				\myarray
					[\Tanset_*^{[1,N]} \BadVar]
					{\Fullset_*^{\leq N;1} \GdVar}
			\right|
			+
			\left|
				\GeoAng \Fullset^{\leq \lfloor N/2 \rfloor - 1;1} f
			\right|
			\left|
				\Tanset^{\leq N} \GdVar
			\right|.
			\notag 
		\end{align}

		Moreover, if $1 \leq N \leq 17$, then the following estimates hold:
		\begin{subequations}
		\begin{align}
		\left|
			[\angD^2, \Fullset^{N;1}] f
		\right|
		& \lesssim
			\left|
				\Fullset_*^{\leq N;1} f
			\right|
				\label{E:ANGDSQUAREDONERADIALTANGENTIALFUNCTIONCOMMUTATOR} 
				\\
		& \ \
			+
			\left|
				\Tanset^{\leq \lceil N/2 \rceil} f
			\right|
			\left|
				\myarray
					[\Tanset_*^{[1, N+1]} \BadVar]
					{\Fullset_*^{\leq N+1;1} \GdVar}
			\right|
			+
			\left|
				\Fullset_*^{\leq \lceil N/2 \rceil} f 
			\right|
			\left|
				\Tanset^{\leq N+1} \GdVar
			\right|,
			\notag \\
		\left|
			[\angLap, \Fullset^{N;1}] f
		\right|
		& \lesssim
			\left|
				\Fullset_*^{\leq N+1;1} f
			\right|
				\label{E:ANGLAPONERADIALTANGENTIALFUNCTIONCOMMUTATOR} 
				\\
		& \ \
			+
			\left|
				\Tanset^{\leq \lceil N/2 \rceil} f
			\right|
			\left|
				\myarray
					[\Tanset_*^{[1, N+1]} \BadVar]
					{\Fullset_*^{\leq N+1;1} \GdVar}
			\right|
			+
			\left|
				\Fullset_*^{\leq \lceil N/2 \rceil} f
			\right|
			\left|
				\Tanset^{\leq N+1} \GdVar
			\right|.
			\notag 
		\end{align}
		\end{subequations}
\end{lemma}

\begin{proof}
	See Subsect.~\ref{SS:OFTENUSEDESTIMATES} for some comments on the analysis.
	The proof is similar to that of Lemma~\ref{L:COMMUTATORESTIMATES}, 
	so we only sketch it by highlighting the
	few differences worth mentioning.
	To illustrate the differences, 
	we prove \eqref{E:ONERADIALTANGENTIALFUNCTIONCOMMUTATORESTIMATE} in detail.
	To proceed, we argue as in the proof of \eqref{E:FIRSTESTIMATEPURETANGENTIALFUNCTIONCOMMUTATORESTIMATE} 
	and use that precisely one factor of $\Fullset^{\vec{I}}$ is equal to $\Rad$,
	thereby deducing that 
	\begin{align} \label{E:FIRSTESTIMATEONERADIALTANGENTIALFUNCTIONCOMMUTATORESTIMATE}
			\mbox{LHS } \eqref{E:ONERADIALTANGENTIALFUNCTIONCOMMUTATORESTIMATE}
			& \lesssim
			\sum_{N_1 + N_2 \leq N-1}
			\left|
				\angLie_{\Tanset}^{N_1} \angdeformoneformupsharparg{\Rad}{\Lunit}
			\right|
			\left|
				\GeoAng \Tanset^{N_2} f
			\right|
			+
			\sum_{N_1 + N_2 \leq N-1}
			\left|
				\angLie_{\Tanset}^{N_1} \angdeformoneformupsharparg{\GeoAng}{\Rad}
			\right|
			\left|
				\GeoAng \Tanset^{N_2} f
			\right|
				\\
			&
			\ \ 
			+
			\sum_{N_1 + N_2 \leq N-1}
			\left|
				\angLie_{\Fullset}^{N_1;1} \angdeformoneformupsharparg{\GeoAng}{\Lunit}
			\right|
			\left|
				\GeoAng \Tanset^{N_2} f
			\right|
			+
			\sum_{N_1 + N_2 \leq N-1}
			\left|
				\angLie_{\Tanset}^{N_1} \angdeformoneformupsharparg{\GeoAng}{\Lunit}
			\right|
			\left|
				\GeoAng \Fullset^{N_2;1} f
			\right|.
			\notag
	\end{align}
	The key point in \eqref{E:FIRSTESTIMATEONERADIALTANGENTIALFUNCTIONCOMMUTATORESTIMATE}
	is that all $\ell_{t,u}-$projected Lie derivatives that fall on 
	$\angdeformoneformupsharparg{\Rad}{\Lunit}$ or $\angdeformoneformupsharparg{\GeoAng}{\Rad}$
	are with respect to vectorfields in $\Tanset$.
	The desired bound \eqref{E:ONERADIALTANGENTIALFUNCTIONCOMMUTATORESTIMATE} 
	now follows easily from the estimates
	\eqref{E:TANGENTDIFFERNTIATEDGEOANGDEFORMSPHERELSHARPPOINTWISE}-\eqref{E:TANGENTDIFFERNTIATEDRADDEFORMSPHERELPOINTWISE}
	and the bootstrap assumptions. Note that the first term on RHS~\eqref{E:RADDEFORMSPHERELPOINTWISE} is
	not necessarily small and hence, in contrast to \eqref{E:PURETANGENTIALFUNCTIONCOMMUTATORESTIMATE}, 
	we do not gain a smallness factor of
	$\varepsilon^{1/2}$ in front of the first term on RHS~\eqref{E:ONERADIALTANGENTIALFUNCTIONCOMMUTATORESTIMATE}.

	The remaining estimates stated in Lemma~\ref{L:TRANSVERALTANGENTIALCOMMUTATOR} 
	can be proved by making similar modifications
	to our proof of Lemma~\ref{L:COMMUTATORESTIMATES}
	and employing the estimates of Lemma~\ref{L:MOREPRECISEANGSPHERELESTIMATES}.

\end{proof}

\begin{corollary}
	\label{C:TANGENTIALDERIVATIVESOFANGLAPPSIPOINTWISE}
	Assume that $1 \leq N \leq 18$. Under the assumptions of 
	Lemma~\ref{L:TRANSVERALTANGENTIALCOMMUTATOR}, 
	the following pointwise estimates hold
	on $\mathcal{M}_{\Tboot,U_0}$:
	\begin{align} \label{E:TANGENTIALDERIVATIVESOFANGLAPPSIPOINTWISE}
		\left|
			\Tanset^{N-1} \angLap \Psi
		\right|
		& \lesssim 
			\left|
				\Tanset^{\leq N+1} \Psi
			\right|
			+
			\left|
				\Tanset^{\leq N} \GdVar
			\right|.
	\end{align}
\end{corollary}

\begin{proof}
	See Subsect.~\ref{SS:OFTENUSEDESTIMATES} for some comments on the analysis.
	Writing 
	$
	\Tanset^{N-1} \angLap \Psi 
	= 
	\angLap \Tanset^{N-1} \Psi 
	+ [\Tanset^{N-1}, \angLap] \Psi
	$,
	we see that the corollary is a simple consequence of
	\eqref{E:ANGDERIVATIVESINTERMSOFTANGENTIALCOMMUTATOR},
	\eqref{E:ANGLAPPURETANGENTIALFUNCTIONCOMMUTATOR} with $f = \Psi$,
	and the bootstrap assumptions.
\end{proof}

%
%
%
%


%

\subsection{Transport inequalities and improvements of the auxiliary bootstrap assumptions}
\label{SS:IMPROVEMENTOFAUX}
In the next proposition, 
we use the previous estimates to derive transport inequalities for the eikonal function quantities and improvements
of the auxiliary bootstrap assumptions. 
The transport inequalities form the starting point for our
derivation of $L^2$ estimates for the below-top-order derivatives of the eikonal function quantities
(see Subsect.~\ref{SS:L2FOREIKONALNOMODIFIED}).
In proving the proposition, we must in particular propagate the smallness of the
$\mathring{\upepsilon}-$sized quantities even though
some terms in the evolution equations involve $\mathring{\updelta}-$sized
quantities, which are allowed to be large.
To this end, we must find and exploit effective partial decoupling 
between various quantities, which is present because of the special
structure of the evolution equations relative to the geometric coordinates
and because of the good properties of the commutation vectorfield
sets $\Fullset$ and $\Tanset$.

\begin{proposition}[\textbf{Transport inequalities and improvements of the auxiliary bootstrap assumptions}] 
\label{P:IMPROVEMENTOFAUX}
Under the data-size and bootstrap assumptions 
of Subsects.~\ref{SS:SIZEOFTBOOT}-\ref{SS:AUXILIARYBOOTSTRAP}
and the smallness assumptions of Subsect.~\ref{SS:SMALLNESSASSUMPTIONS}, 
the following estimates hold
on $\mathcal{M}_{\Tboot,U_0}$
(see Subsect.~\ref{SS:STRINGSOFCOMMUTATIONVECTORFIELDS} regarding the vectorfield operator notation):

\medskip
\noindent \underline{\textbf{Transport inequalities for the eikonal function quantities}.}

\medskip


\noindent \textbf{$\bullet$Transport inequalities for} $\upmu$.
	The following pointwise estimate holds:
	\begin{subequations}
	\begin{align} 
		\left|
			\Lunit \upmu
		\right|
		& \lesssim 
			\left|
				\Fullset^{\leq 1} \Psi
			\right|.
			\label{E:LUNITUPMUPOINTWISE} 
		\end{align}

	Moreover, for $1 \leq N \leq 18$, the following estimates hold:
	\begin{align} \label{E:PURETANGENTIALLUNITUPMUCOMMUTEDESTIMATE}
		\left|
			\Lunit \Tanset^N \upmu
		\right|,
			\,
		\left|
			\Tanset^N \Lunit \upmu
		\right|
		& \lesssim 
			\left|
				\Fullset_*^{\leq N+1;1} \Psi
			\right|
			+
			\left|
				\Tanset^{\leq N} \GdVar
			\right|
			+
			\varepsilon
			\left|
				\Tanset_*^{[1,N]} \BadVar
			\right|.
	\end{align}
	\end{subequations}

	\noindent \textbf{$\bullet$Transport inequalities for} $\Lunit_{(Small)}^i$ and $\mytr \upchi$.
	For $N \leq 18$, the following estimates hold:
	\begin{subequations}
	\begin{align}
		\left|
			\myarray
				[\Lunit \Tanset^N \Lunit_{(Small)}^i]
				{\Lunit \Tanset^{N-1} \mytr \upchi}
		\right|,
			\,
		\left|
			\myarray
				[\Tanset^N \Lunit \Lunit_{(Small)}^i]
				{\Tanset^{N-1} \Lunit \mytr \upchi}
		\right|
		& \lesssim 
			\left|
				\Tanset^{\leq N+1} \Psi
			\right|
		+ \varepsilon
			\left|
				\Tanset^{\leq N} \GdVar
			\right|,
				\label{E:LUNITTANGENTDIFFERENTIATEDLUNITSMALLIMPROVEDPOINTWISE} \\
		\left|
			\myarray
				[\Lunit \Fullset^{N;1} \Lunit_{(Small)}^i]
				{\Lunit \Fullset^{N-1;1} \mytr \upchi}
		\right|,
			\,
		\left|
			\myarray
				[\Fullset^{N;1} \Lunit \Lunit_{(Small)}^i]
				{\Fullset^{N-1;1} \Lunit \mytr \upchi}
		\right|
		& \lesssim 
		\left|
			\Fullset_*^{\leq N+1;1} \Psi
		\right|
		+ 
		\left|
			\myarray[\varepsilon \Tanset_*^{[1,N]} \BadVar]
				{\Fullset_*^{\leq N;1} \GdVar}
		\right|.
			\label{E:LUNITONERADIALTANGENTDIFFERENTIATEDLUNITSMALLIMPROVEDPOINTWISE} 
\end{align}
\end{subequations}

\medskip
\noindent \underline{$L^{\infty}$ \textbf{estimates for} $\Psi$ \textbf{and the eikonal function quantities}.}

\medskip

\noindent \textbf{$\bullet$$L^{\infty}$ estimates involving at most one transversal derivative of $\Psi$}. 
The following estimates hold:
\begin{subequations}
\begin{align} \label{E:PSITRANSVERSALLINFINITYBOUNDBOOTSTRAPIMPROVED}
	\left\| 
		\Rad \Psi 
	\right\|_{L^{\infty}(\Sigma_t^u)}
	& \leq 
	\left\| 
		\Rad \Psi 
	\right\|_{L^{\infty}(\Sigma_0^u)}
	+ C \varepsilon,
		\\
		\left\| 
			\Fullset_*^{\leq 10;1} \Psi
		\right\|_{L^{\infty}(\Sigma_t^u)}
		& \leq C \varepsilon.
		\label{E:PSIMIXEDTRANSVERSALTANGENTBOOTSTRAPIMPROVED}
\end{align}
\end{subequations}

\noindent \textbf{$\bullet$$L^{\infty}$ estimates for $\upmu$}. 
	The following estimates hold:
	\begin{subequations}
	\begin{align} \label{E:LUNITUPMULINFINITY}
		\left\| 
			\Lunit \upmu
		\right\|_{L^{\infty}(\Sigma_t^u)}
		& 
		= 
		\frac{1}{2}
		\left\| 
			G_{\Lunit \Lunit} \Rad \Psi
		\right\|_{L^{\infty}(\Sigma_0^u)}
		+ \mathcal{O}(\varepsilon),
			\\
		\left\|
			\Lunit \Tanset^{[1,9]} \upmu
		\right\|_{L^{\infty}(\Sigma_t^u)},
			\,
		\left\| 
			\Tanset_*^{[1,9]} \upmu
		\right\|_{L^{\infty}(\Sigma_t^u)}
		& \leq
			C \varepsilon,
			\label{E:LUNITAPPLIEDTOTANGENTIALUPMUANDTANSETSTARLINFTY}
	\end{align}
	\end{subequations}

	\begin{subequations}
	\begin{align} \label{E:UPMULINFTY}
	\left\| 
			\upmu - 1
		\right\|_{L^{\infty}(\Sigma_t^u)}
		& \leq
	 	2 \TranminusdatasizeWithFactor^{-1} 
		\left\| 
			G_{\Lunit \Lunit} \Rad \Psi
		\right\|_{L^{\infty}(\Sigma_0^u)}
		+ C \varepsilon.
	\end{align}
	\end{subequations}

\noindent \textbf{$\bullet$$L^{\infty}$ estimates for $\Lunit_{(Small)}^i$ and $\upchi$}.
The following estimates hold:
\begin{subequations}
\begin{align}  \label{E:PURETANGENTIALLUNITAPPLIEDTOLISMALLANDLISMALLINFTYESTIMATE}
	\left\|
		\Lunit \Tanset^{\leq 10} \Lunit_{(Small)}^i
	\right\|_{L^{\infty}(\Sigma_t^u)},
		\,
	\left\|
		\Tanset^{\leq 10} \Lunit_{(Small)}^i
	\right\|_{L^{\infty}(\Sigma_t^u)}
	& \leq C \varepsilon,
		\\
	\left\|
		\Lunit \Fullset^{\leq 9;1} \Lunit_{(Small)}^i
	\right\|_{L^{\infty}(\Sigma_t^u)},
		\,
	\left\|
		\Fullset_*^{\leq 9;1} \Lunit_{(Small)}^i
	\right\|_{L^{\infty}(\Sigma_t^u)}
	& \leq C \varepsilon,
		\label{E:LUNITAPPLIEDTOLISMALLANDLISMALLINFTYESTIMATE} \\
	\left\|
		\Rad \Lunit_{(Small)}^i
	\right\|_{L^{\infty}(\Sigma_t^u)}
	& \leq
	\left\| 
		\Rad \Lunit_{(Small)}^i
	\right\|_{L^{\infty}(\Sigma_0^u)}
	+  C \varepsilon,
	\label{E:LISMALLLONERADIALINFINITYESTIMATE}
\end{align}
\end{subequations}


\begin{align}  \label{E:PURETANGENTIALCHICOMMUTEDLINFINITY}
		\left\|
			\angLie_{\Tanset}^{\leq 9} \upchi
		\right\|_{L^{\infty}(\Sigma_t^u)},
			\,
		\left\|
			\angLie_{\Tanset}^{\leq 9} \upchi^{\#}
		\right\|_{L^{\infty}(\Sigma_t^u)},
			\,
		\left\|
			\Tanset^{\leq 9} \mytr \upchi
		\right\|_{L^{\infty}(\Sigma_t^u)}
		& \leq C \varepsilon,
			\\
		\left\|
			\angLie_{\Fullset}^{\leq 8;1} \upchi
		\right\|_{L^{\infty}(\Sigma_t^u)},
			\,
		\left\|
			\angLie_{\Fullset}^{\leq 8;1} \upchi^{\#}
		\right\|_{L^{\infty}(\Sigma_t^u)},
			\,
		\left\|
			\Fullset^{\leq 8;1} \mytr \upchi
		\right\|_{L^{\infty}(\Sigma_t^u)}
		& \leq C \varepsilon.
		\label{E:ONERADIALCHICOMMUTEDLINFINITY}
\end{align}
\end{proposition}

\begin{remark}[\textbf{The auxiliary bootstrap assumptions of Subsect.~\ref{SS:AUXILIARYBOOTSTRAP}
  are now redundant}]
  \label{R:AUXILIARYBOOTSTRAPASSUMPTIONS}
	Since Prop.~\ref{P:IMPROVEMENTOFAUX} in particular
	provides an improvement of the auxiliary bootstrap assumptions of
	Subsect.~\ref{SS:AUXILIARYBOOTSTRAP},
	we do not bother to include those bootstrap assumptions
	in the hypotheses of any 
	of the lemmas or propositions proved in the remainder of the article.
\end{remark}

\begin{proof}[Proof of Prop.~\ref{P:IMPROVEMENTOFAUX}]
	See Subsect.~\ref{SS:OFTENUSEDESTIMATES} for some comments on the analysis.
	We must derive the estimates in a viable order.
	Throughout this proof, we use
	the estimates of Lemma~\ref{L:BEHAVIOROFEIKONALFUNCTIONQUANTITIESALONGSIGMA0}
	and the assumption \eqref{E:DATAEPSILONVSBOOTSTRAPEPSILON}
	without explicitly mentioning them each time. We refer to these as ``conditions on the data.''
	Similarly, when we say that we use the ``bootstrap assumptions,''
	we mean the assumptions on $\mathcal{M}_{\Tboot,U_0}$ stated in 
	Subsects.~\ref{SS:SIZEOFTBOOT}-\ref{SS:AUXILIARYBOOTSTRAP}.

	\medskip
	\noindent \textbf{Proof of the estimates \eqref{E:LUNITTANGENTDIFFERENTIATEDLUNITSMALLIMPROVEDPOINTWISE} for 
	$\Lunit \Tanset^N \Lunit_{(Small)}^i$ and $\Tanset^N \Lunit \Lunit_{(Small)}^i$ and 
	\eqref{E:PURETANGENTIALLUNITAPPLIEDTOLISMALLANDLISMALLINFTYESTIMATE}:}
	We prove the estimate
	\eqref{E:LUNITTANGENTDIFFERENTIATEDLUNITSMALLIMPROVEDPOINTWISE} for $\Lunit \Tanset^N \Lunit_{(Small)}^i$
	and omit the proof for $\Tanset^N \Lunit \Lunit_{(Small)}^i$; the proof of the latter estimate is similar but simpler because it 
	involves fewer commutation estimates.
	To derive the desired bound, 
	in place of \eqref{E:LUNITTANGENTDIFFERENTIATEDLUNITSMALLIMPROVEDPOINTWISE},
	we first show that 
	\begin{align}  \label{E:PROOFLUNITTANGENTDIFFERENTIATEDLUNITSMALLIMPROVEDPOINTWISE}
		\left|
			\Lunit \Tanset^N \Lunit_{(Small)}^i
		\right|
		& \lesssim 
			\left|
				\Tanset^{\leq N+1} \Psi
			\right|
		+ \varepsilon^{1/2}
			\left|
				\Tanset^{\leq N} \GdVar
			\right|.
	\end{align}
	The factors of $\varepsilon^{1/2}$ in \eqref{E:PROOFLUNITTANGENTDIFFERENTIATEDLUNITSMALLIMPROVEDPOINTWISE}
	arise from the auxiliary bootstrap assumptions of Subsect.~\ref{SS:AUXILIARYBOOTSTRAP}.
	At the end of the proof, we will have shown that the auxiliary bootstrap assumptions
	have been improved in that they hold with $C \varepsilon$ 
	in place of $\varepsilon^{1/2}$. Using this improvement, we easily conclude 
	\eqref{E:LUNITTANGENTDIFFERENTIATEDLUNITSMALLIMPROVEDPOINTWISE} for 
	$\Lunit \Tanset^N \Lunit_{(Small)}^i$
	by repeating the proof of
	\eqref{E:PROOFLUNITTANGENTDIFFERENTIATEDLUNITSMALLIMPROVEDPOINTWISE} 
	with $C \varepsilon$ in place of the factor $\varepsilon^{1/2}$.
	To prove \eqref{E:PROOFLUNITTANGENTDIFFERENTIATEDLUNITSMALLIMPROVEDPOINTWISE},
	we commute equation \eqref{E:LLUNITI} with
	$\Tanset^N$ and use Lemma~\ref{L:SCHEMATICDEPENDENCEOFMANYTENSORFIELDS}
	to derive the schematic equation
	\begin{align} \label{E:LLUNITISCHEMATICCOMMUTED}
		\Lunit \Tanset^N \Lunit_{(Small)}^i
		& = [\Lunit, \Tanset^N] \Lunit_{(Small)}^i
			+
			\Tanset^N 
			\left\lbrace
				\smoothfunction(\GdVar,\ginversesphere,\angdiff x^1,\angdiff x^2) \Singletan \Psi
			\right\rbrace.
	\end{align}
	To bound the second term on RHS~\eqref{E:LLUNITISCHEMATICCOMMUTED}
	by RHS~\eqref{E:PROOFLUNITTANGENTDIFFERENTIATEDLUNITSMALLIMPROVEDPOINTWISE}, 
	we use
	Lemmas~\ref{L:POINTWISEFORRECTANGULARCOMPONENTSOFVECTORFIELDS}
	and
	\ref{L:POINTWISEESTIMATESFORGSPHEREANDITSDERIVATIVES}
	and the bootstrap assumptions.
	To bound the remaining term $\left|[\Lunit, \Tanset^N] \Lunit_{(Small)}^i \right|$, 
	we use the commutator estimate \eqref{E:PURETANGENTIALFUNCTIONCOMMUTATORESTIMATE} with 
	$f = L_{(Small)}^i$ and the bootstrap assumptions. We have thus proved the desired bound
	\eqref{E:PROOFLUNITTANGENTDIFFERENTIATEDLUNITSMALLIMPROVEDPOINTWISE}.
	Next, to derive 
	the estimate
	$
	\left\|
		\Lunit \Tanset^{\leq 10} \Lunit_{(Small)}^i
	\right\|_{L^{\infty}(\Sigma_t^u)}
	\lesssim \varepsilon
	$
	stated in \eqref{E:PURETANGENTIALLUNITAPPLIEDTOLISMALLANDLISMALLINFTYESTIMATE},
	we use \eqref{E:PROOFLUNITTANGENTDIFFERENTIATEDLUNITSMALLIMPROVEDPOINTWISE}
	and the bootstrap assumptions.
	To obtain the estimates
	$
	\left\|
		\Tanset^{\leq 10} \Lunit_{(Small)}^i
	\right\|_{L^{\infty}(\Sigma_t^u)}
	\lesssim \varepsilon
	$
	stated in \eqref{E:PURETANGENTIALLUNITAPPLIEDTOLISMALLANDLISMALLINFTYESTIMATE},
	we first use the fundamental theorem of calculus to write
	\begin{align} \label{E:INTEGRATINGALONGINTEGRALCURVES}
		\Tanset^{\leq 10} \Lunit_{(Small)}^i(t,u,\vartheta)
		= \Tanset^{\leq 10} \Lunit_{(Small)}^i(0,u,\vartheta)
			+ 
			\int_{s=0}^t
				\Lunit \Tanset^{\leq 10} \Lunit_{(Small)}^i(s,u,\vartheta)
			\, ds.
	\end{align}
	We then use the conditions on the data to bound 
	$
	\left|
		\Tanset^{\leq 10} \Lunit_{(Small)}^i(0,u,\vartheta)
	\right|
	\lesssim \varepsilon
	$
	and the inequality
	$
	\left\|
		\Lunit \Tanset^{\leq 10} \Lunit_{(Small)}^i
	\right\|_{L^{\infty}(\Sigma_t^u)}
	\lesssim \varepsilon
	$
	to bound the time integral
	on RHS~\eqref{E:INTEGRATINGALONGINTEGRALCURVES}
	by $\lesssim \Tboot \varepsilon \lesssim \varepsilon$,
	which in total yields the desired result.

	\medskip
	\noindent \textbf{Proof of \eqref{E:LUNITUPMUPOINTWISE}:}
	We first use equation \eqref{E:UPMUFIRSTTRANSPORT}
	and Lemma~\ref{L:SCHEMATICDEPENDENCEOFMANYTENSORFIELDS}
	to deduce 
	$
	\Lunit \upmu 
	= \smoothfunction(\BadVar) \Singletan \Psi
		+ \smoothfunction(\GdVar) \Rad \Psi
	$.
	The desired estimate \eqref{E:LUNITUPMUPOINTWISE}
	now follows easily from the previous expression 
	and the bootstrap assumptions.

	\medskip
	\noindent \textbf{Proof of \eqref{E:PSITRANSVERSALLINFINITYBOUNDBOOTSTRAPIMPROVED} and
	\eqref{E:PSIMIXEDTRANSVERSALTANGENTBOOTSTRAPIMPROVED}:}
	We first note that in proving \eqref{E:PSIMIXEDTRANSVERSALTANGENTBOOTSTRAPIMPROVED}, 
	we may assume that the operator $\Fullset_*^{\leq 10;1}$ contains the factor $\Rad$ since
	otherwise the estimate is implied by the bootstrap assumption \eqref{E:PSIFUNDAMENTALC0BOUNDBOOTSTRAP}.
	To proceed, 
	we use equation \eqref{E:LONOUTSIDEGEOMETRICWAVEOPERATORFRAMEDECOMPOSED},
	Lemma~\ref{L:SCHEMATICDEPENDENCEOFMANYTENSORFIELDS},
	and the aforementioned fact
	$
	\Lunit \upmu 
	= \smoothfunction(\BadVar) \Singletan \Psi
		+ \smoothfunction(\GdVar) \Rad \Psi
	$
	to rewrite the wave equation as
	\begin{align} \label{E:WAVEEQUATIONTRANSPORTINTERPRETATION}
		\Lunit \Rad \Psi 
		& 
		= \smoothfunction(\BadVar) \angLap \Psi 
			+ \smoothfunction(\BadVar,\ginversesphere,\angdiff x^1,\angdiff x^2,\Singletan \Psi, \Rad \Psi) 
				\Singletan \Singletan \Psi
			+ \smoothfunction(\BadVar,\ginversesphere,\angdiff x^1,\angdiff x^2,\Singletan \Psi, \Rad \Psi) 
				\Singletan \GdVar.
	\end{align}
	Commuting \eqref{E:WAVEEQUATIONTRANSPORTINTERPRETATION}
	with $\Tanset^N$, $(0 \leq N \leq 9)$,
	and using Lemmas~\ref{L:POINTWISEFORRECTANGULARCOMPONENTSOFVECTORFIELDS}
	and
	\ref{L:POINTWISEESTIMATESFORGSPHEREANDITSDERIVATIVES}
	and the bootstrap assumptions, 
	we find that
	\begin{align} \label{E:PROOFWAVEEQNTRANSPORTTANGENTIALCOMMUTED}
		\left|
			\Lunit \Tanset^N \Rad \Psi
		\right|
		& \leq
		\left|
			\Lunit \Tanset^N \Rad \Psi
			-
			\Tanset^N \Lunit \Rad \Psi
		\right| 
		+
		\left|
			\Tanset^N \Lunit \Rad \Psi
		\right|
		\\
		& \lesssim
				\left|
					\Lunit \Tanset^N \Rad \Psi
					-  \Tanset^N \Lunit \Rad \Psi
				\right|
				+ 
				\left|
					[\angLap, \Tanset^N] \Psi
				\right|
				\notag	\\
		& \ \
				+
				\left|
					\Tanset^{\leq N+2;1} \Psi
				\right|
				+
				\varepsilon^{1/2}
				\left|
					\Fullset_*^{\leq N+1;1} \Psi
				\right|
				+
				\left|
					\Tanset^{\leq N+1} \GdVar
				\right|
				+
				\varepsilon^{1/2}
				\left|
					\Tanset_*^{[1,N]} \BadVar
				\right|.
				\notag
	\end{align}
	Using in addition the commutator estimates \eqref{E:ONERADIALTANGENTIALFUNCTIONCOMMUTATORESTIMATE}
	and \eqref{E:ANGLAPPURETANGENTIALFUNCTIONCOMMUTATOR} with $f = \Psi$,
	we bound the two commutator terms on the second line of RHS~\eqref{E:PROOFWAVEEQNTRANSPORTTANGENTIALCOMMUTED}
	by $\lesssim$ the terms on the last line of RHS~\eqref{E:PROOFWAVEEQNTRANSPORTTANGENTIALCOMMUTED}. 
	The bootstrap assumptions imply 
	that most terms on the last line of \eqref{E:PROOFWAVEEQNTRANSPORTTANGENTIALCOMMUTED}
	are $\lesssim \varepsilon$. The exceptional terms (that is, the ones not included in ``most terms'') are
	$\Tanset^{\leq 10} \Lunit_{(Small)}^i$ for $i = 1,2$,
	but we have already shown that these terms are bounded
	in the norm $\| \cdot \|_{L^{\infty}(\Sigma_t^u)}$ by $\lesssim \varepsilon$.
	In total, we find that
	$\left\|
		\Lunit \Tanset^{\leq 9} \Rad \Psi
	\right\|_{L^{\infty}(\Sigma_t^u)} \lesssim \varepsilon$. 
	Integrating
	along the integral curves of
	$\Lunit$ as in \eqref{E:INTEGRATINGALONGINTEGRALCURVES} and using the conditions on the data, 
	we conclude 
	\eqref{E:PSITRANSVERSALLINFINITYBOUNDBOOTSTRAPIMPROVED}
	and also the estimate
	$\left\|
		\Tanset^{[1,9]} \Rad \Psi
	\right\|_{L^{\infty}(\Sigma_t^u)}
	\lesssim \varepsilon
	$.
	Moreover, we 
	use the commutator estimate \eqref{E:ONERADIALTANGENTIALFUNCTIONCOMMUTATORESTIMATE}
	with $f = \Psi$ 
	and the bootstrap assumptions
	to commute the factor of $\Rad$ in $\Fullset_*^{\leq 10;1}$
	so that it hits $\Psi$ first, 
	thereby concluding that
	$
	\left\|
		\Fullset_*^{\leq 10;1} \Psi
	\right\|_{L^{\infty}(\Sigma_t^u)}
	\lesssim
	\left\|
		\Tanset^{[1,9]} \Rad \Psi
	\right\|_{L^{\infty}(\Sigma_t^u)}
	+
	\varepsilon
	\lesssim \varepsilon
	$.
	We have thus proved the desired bound \eqref{E:PSIMIXEDTRANSVERSALTANGENTBOOTSTRAPIMPROVED}.

	\medskip
	\noindent \textbf{Proof of \eqref{E:LUNITUPMULINFINITY} and \eqref{E:UPMULINFTY}:}
	To derive \eqref{E:LUNITUPMULINFINITY},
	we first use
	equation \eqref{E:UPMUFIRSTTRANSPORT}
	and Lemma~\ref{L:SCHEMATICDEPENDENCEOFMANYTENSORFIELDS} to write
	$ \Lunit \upmu
		=
		\frac{1}{2} G_{\Lunit \Lunit} \Rad \Psi
		+
		\smoothfunction(\BadVar) \Singletan \Psi
	$.
	From the previous expression 
	and the bootstrap assumptions, 
	we deduce that
	$
	\left\| 
			\Lunit \upmu
		\right\|_{L^{\infty}(\Sigma_t^u)}
		= 
		\frac{1}{2}
		\left\| 
			G_{\Lunit \Lunit} \Rad \Psi
		\right\|_{L^{\infty}(\Sigma_t^u)}
		+ \mathcal{O}(\varepsilon)
	$.
	Next, we use Lemma~\ref{L:SCHEMATICDEPENDENCEOFMANYTENSORFIELDS} to deduce that
	$G_{\Lunit \Lunit} \Rad \Psi = \smoothfunction(\GdVar) \Rad \Psi$.
	Applying $\Lunit$ to the previous expression and using 
	the bounds
	$\left\| 
		\Lunit \Lunit_{(Small)}^i
	\right\|_{L^{\infty}(\Sigma_t^u)}
	$,
	$
	\left\| 
		\Lunit \Rad \Psi
	\right\|_{L^{\infty}(\Sigma_t^u)}
	$
	$
	\lesssim \varepsilon
	$
	proven above and
	the bootstrap assumptions,
	we find that
	$
	\left\| 
		\Lunit(G_{\Lunit \Lunit} \Rad \Psi)
	\right\|_{L^{\infty}(\Sigma_t^u)}
	$
	$
	\lesssim \varepsilon
	$.
	Integrating 
	along the integral curves of $\Lunit$
	as in \eqref{E:INTEGRATINGALONGINTEGRALCURVES}
	and using the previous inequality,
	we find that
	$
	\left\| 
		G_{\Lunit \Lunit} \Rad \Psi
	\right\|_{L^{\infty}(\Sigma_t^u)}
	= 
	\left\| 
		G_{\Lunit \Lunit} \Rad \Psi
	\right\|_{L^{\infty}(\Sigma_0^u)}
	+ \mathcal{O}(\varepsilon)
	$. Inserting this estimate into the first
	estimate of this paragraph, we conclude \eqref{E:LUNITUPMULINFINITY}.
	The estimate \eqref{E:UPMULINFTY}
	then follows from integrating 
	along the integral curves of $\Lunit$ 
	as in \eqref{E:INTEGRATINGALONGINTEGRALCURVES} and using 
	\eqref{E:LUNITUPMULINFINITY},
	the conditions on the data,
	and the assumption $\Tboot \leq 2 \TranminusdatasizeWithFactor^{-1}$.

	\noindent \textbf{Proof of \eqref{E:PURETANGENTIALLUNITUPMUCOMMUTEDESTIMATE} and 
	\eqref{E:LUNITAPPLIEDTOTANGENTIALUPMUANDTANSETSTARLINFTY}:}
	We now prove \eqref{E:PURETANGENTIALLUNITUPMUCOMMUTEDESTIMATE}
	for $\Lunit \Tanset^N \upmu$. The proof for $\Tanset^N \Lunit \upmu$
	is similar but simpler because it involves
	fewer commutation estimates; we omit these details.
	To proceed, in place of \eqref{E:PURETANGENTIALLUNITUPMUCOMMUTEDESTIMATE},
	we first prove that
	\begin{align}  \label{E:PROOFPURETANGENTIALLUNITUPMUCOMMUTEDESTIMATE}
		\left|
			\Lunit \Tanset^N \upmu
		\right|
		& \lesssim 
			\left|
				\Fullset_*^{\leq N+1;1} \Psi
			\right|
		+ \left|
				\Tanset^{\leq N} \GdVar
			\right|
		+ \varepsilon^{1/2}
			\left|
				\Tanset_*^{[1,N]} \BadVar
			\right|.
	\end{align}
	As we described above,
	at the end of the proof, we will have shown that the auxiliary bootstrap assumptions
	have been improved in that they hold with $C \varepsilon$ 
	in place of $\varepsilon^{1/2}$ and this improvement 
	implies that \eqref{E:PROOFPURETANGENTIALLUNITUPMUCOMMUTEDESTIMATE} 
	holds with $C \varepsilon$ in place of $\varepsilon^{1/2}$ as desired.
	To prove \eqref{E:PROOFPURETANGENTIALLUNITUPMUCOMMUTEDESTIMATE}, 
	we commute the equation
	$
	\Lunit \upmu 
	= \smoothfunction(\BadVar) \Singletan \Psi
		+ \smoothfunction(\GdVar) \Rad \Psi
	$
	(see equation \eqref{E:UPMUFIRSTTRANSPORT} and Lemma~\ref{L:SCHEMATICDEPENDENCEOFMANYTENSORFIELDS})
	with $\Tanset^N$ to deduce the schematic identity
	\begin{align} \label{E:LUPMUTANGENTIALLYSCHEMATICCOMMUTED}
		\Lunit \Tanset^N \upmu
		& = [\Lunit, \Tanset^N] \upmu
			+
			\Tanset^N 
			\left\lbrace
				\smoothfunction(\GdVar) \Rad \Psi
				+ \smoothfunction(\BadVar) \Singletan \Psi
			\right\rbrace.
	\end{align}
	To bound the second term on RHS~\eqref{E:LUPMUTANGENTIALLYSCHEMATICCOMMUTED} by
	RHS~\eqref{E:PROOFPURETANGENTIALLUNITUPMUCOMMUTEDESTIMATE}, 
	we use the already proven bound \eqref{E:PURETANGENTIALLUNITAPPLIEDTOLISMALLANDLISMALLINFTYESTIMATE}
	and the bootstrap assumptions.
	To bound the term $[\Lunit, \Tanset^N] \upmu$ on RHS~\eqref{E:LUPMUTANGENTIALLYSCHEMATICCOMMUTED}
	by RHS~\eqref{E:PROOFPURETANGENTIALLUNITUPMUCOMMUTEDESTIMATE},
	we use the commutator estimate \eqref{E:PURETANGENTIALFUNCTIONCOMMUTATORESTIMATE} with 
	$f = \upmu$
	and the bootstrap assumptions.
	To prove the estimate
	$
	\left\|
		\Lunit \Tanset^{[1,9]} \upmu
	\right\|_{L^{\infty}(\Sigma_t^u)}
	\lesssim \varepsilon
	$
	stated in \eqref{E:LUNITAPPLIEDTOTANGENTIALUPMUANDTANSETSTARLINFTY},
	we use \eqref{E:PROOFPURETANGENTIALLUNITUPMUCOMMUTEDESTIMATE},
	the already proven bounds 
	\eqref{E:PURETANGENTIALLUNITAPPLIEDTOLISMALLANDLISMALLINFTYESTIMATE}
	and \eqref{E:PSIMIXEDTRANSVERSALTANGENTBOOTSTRAPIMPROVED},
	and the bootstrap assumptions.
	The estimate
	\eqref{E:LUNITAPPLIEDTOTANGENTIALUPMUANDTANSETSTARLINFTY}
	for
	$
	\left\|
		\Tanset_*^{[1,9]} \upmu
	\right\|_{L^{\infty}(\Sigma_t^u)}
	$
	then follows from integrating 
	along the integral curves of
	$\Lunit$ 
	as in \eqref{E:INTEGRATINGALONGINTEGRALCURVES}
	and using the estimate
	$
	\left\|
		\Lunit \Tanset^{[1,9]} \upmu
	\right\|_{L^{\infty}(\Sigma_t^u)}
	\lesssim \varepsilon
	$
	and the conditions on the data.

	\noindent \textbf{Proof of \eqref{E:LUNITONERADIALTANGENTDIFFERENTIATEDLUNITSMALLIMPROVEDPOINTWISE}
	for $\Lunit \Fullset^{N;1} \Lunit_{(Small)}^i$ and $\Fullset^{N;1} \Lunit \Lunit_{(Small)}^i$
	and \eqref{E:LUNITAPPLIEDTOLISMALLANDLISMALLINFTYESTIMATE}:}
	We now prove \eqref{E:LUNITONERADIALTANGENTDIFFERENTIATEDLUNITSMALLIMPROVEDPOINTWISE} 
	for $\Lunit \Fullset^{N;1} \Lunit_{(Small)}^i$. 
	The proof of \eqref{E:LUNITONERADIALTANGENTDIFFERENTIATEDLUNITSMALLIMPROVEDPOINTWISE} 
	for $\Fullset^{N;1} \Lunit \Lunit_{(Small)}^i$ is similar but simpler
	because it involves fewer commutation estimates; we omit these details.
	We may assume that $\Fullset^{N;1}$ contains a factor $\Rad$ since otherwise the desired estimate is implied by
	\eqref{E:LUNITTANGENTDIFFERENTIATEDLUNITSMALLIMPROVEDPOINTWISE}.
	The proof is similar to the proof of 
	\eqref{E:PROOFLUNITTANGENTDIFFERENTIATEDLUNITSMALLIMPROVEDPOINTWISE},
	the new feature being that we need to exploit the already proven estimates
	$\left\|
		\Tanset^{\leq 10} \Lunit_{(Small)}^i
	\right\|_{L^{\infty}(\Sigma_t^u)}
	\lesssim \varepsilon
	$
	and
	$\left\| 
		\Fullset_*^{\leq 10;1} \Psi
	\right\|_{L^{\infty}(\Sigma_t^u)}
	\lesssim \varepsilon
	$.
	To proceed, we note that \eqref{E:LLUNITISCHEMATICCOMMUTED} holds
	with $\Fullset^{N;1}$ in place of $\Tanset^N$ on both sides
	and that the non-commutator term is easy to bound by using arguments similar to the ones we used in
	proving \eqref{E:LUNITTANGENTDIFFERENTIATEDLUNITSMALLIMPROVEDPOINTWISE}.
	It remains for us to bound the commutator term
	$\left|[\Lunit, \Fullset^{N;1}] \Lunit_{(Small)}^i \right|$
	by $\lesssim \mbox{RHS~\eqref{E:LUNITONERADIALTANGENTDIFFERENTIATEDLUNITSMALLIMPROVEDPOINTWISE}}$.
	This estimate follows from the commutator estimate \eqref{E:ONERADIALTANGENTIALFUNCTIONCOMMUTATORESTIMATE}
	with $f = L_{(Small)}^i$,
	the already proven estimate for 
	$
	\Tanset^{\leq 10} \Lunit_{(Small)}^i
	$ 
	mentioned above
	(to bound the factor 
	$
	\left|
			\Tanset_*^{[1,\lfloor N/2 \rfloor]} f
	\right|
	$
	from the second line of RHS~\eqref{E:ONERADIALTANGENTIALFUNCTIONCOMMUTATORESTIMATE}
	by $\lesssim \varepsilon$),
	and the bootstrap assumptions.
	We have thus obtained the desired estimate \eqref{E:LUNITONERADIALTANGENTDIFFERENTIATEDLUNITSMALLIMPROVEDPOINTWISE} 
	for $\Lunit \Fullset^{N;1} \Lunit_{(Small)}^i$.
	Next, from from the estimate \eqref{E:LUNITONERADIALTANGENTDIFFERENTIATEDLUNITSMALLIMPROVEDPOINTWISE} 
	for $\Lunit \Fullset^{N;1} \Lunit_{(Small)}^i$,
	the already proven estimates 
	$\left\|
		\Tanset^{\leq 10} \Lunit_{(Small)}^i
	\right\|_{L^{\infty}(\Sigma_t^u)}
	\lesssim \varepsilon
	$
	and 
	$\left\| 
		\Fullset_*^{\leq 10;1} \Psi
	\right\|_{L^{\infty}(\Sigma_t^u)}
	\lesssim \varepsilon
	$,
	and the bootstrap assumptions,
	we find that
	$\left\|
		\Lunit \Fullset^{\leq 9;1} \Lunit_{(Small)}^i
	\right\|_{L^{\infty}(\Sigma_t^u)} \lesssim \varepsilon$.
	This completes the proof of \eqref{E:LUNITAPPLIEDTOLISMALLANDLISMALLINFTYESTIMATE} for the first term 
	$
	\left\|
		\Lunit \Fullset^{\leq 9;1} \Lunit_{(Small)}^i
	\right\|_{L^{\infty}(\Sigma_t^u)}
	$
	on the LHS.
	Integrating 
	along the integral curves of $\Lunit$ as in \eqref{E:INTEGRATINGALONGINTEGRALCURVES}
	and using the estimate \eqref{E:LUNITAPPLIEDTOLISMALLANDLISMALLINFTYESTIMATE} for 
	$
	\left\|
		\Lunit \Fullset^{\leq 9;1} \Lunit_{(Small)}^i
	\right\|_{L^{\infty}(\Sigma_t^u)}
	$
	as well as the conditions on the data,
	we conclude the estimate \eqref{E:LUNITAPPLIEDTOLISMALLANDLISMALLINFTYESTIMATE} for
	$
	\left\|
		\Fullset_*^{\leq 9;1} \Lunit_{(Small)}^i
	\right\|_{L^{\infty}(\Sigma_t^u)}
	$
	as well as \eqref{E:LISMALLLONERADIALINFINITYESTIMATE}.

	\noindent \textbf{Proof of \eqref{E:PURETANGENTIALCHICOMMUTEDLINFINITY} and \eqref{E:ONERADIALCHICOMMUTEDLINFINITY}:}
	These two estimates follow from Lemma~\ref{L:POINTWISEESTIMATESFORGSPHEREANDITSDERIVATIVES},
	the already proven estimates
	$
	\left\| 
		\Fullset_*^{\leq 10;1} \Psi
	\right\|_{L^{\infty}(\Sigma_t^u)},
	$
	$
	\left\|
		\Tanset^{\leq 10} \Lunit_{(Small)}^i
	\right\|_{L^{\infty}(\Sigma_t^u)},
	$
	$
	\left\|
		\Fullset_*^{\leq 9;1} \Lunit_{(Small)}^i
	\right\|_{L^{\infty}(\Sigma_t^u)}
	\lesssim \varepsilon
	$
	and the bootstrap assumptions.

	\noindent \textbf{Proof of the estimate \eqref{E:LUNITTANGENTDIFFERENTIATEDLUNITSMALLIMPROVEDPOINTWISE} 
	for $\Tanset^{N-1} \Lunit \mytr \upchi$:}
	We first take the $\gsphere-$trace of equation \eqref{E:CHIINTERMSOFOTHERVARIABLES},
	apply $\Lunit$, 
	and use the schematic identity $\angLie_{\Lunit} \ginversesphere = (\ginversesphere)^{-2} \upchi
	=\smoothfunction(\GdVar,\ginversesphere,\angdiff x^1,\angdiff x^2) \Singletan \GdVar$
	to deduce that
	$\Lunit \mytr \upchi 
	= 
	\smoothfunction(\GdVar,\ginversesphere,\angdiff x^1,\angdiff x^2) \Singletan \Lunit \GdVar
	+ l.o.t.$,
	where
	$l.o.t.
	:=
	\smoothfunction(\Tanset^{\leq 1} \GdVar, \ginversesphere,\angdiff x^1,\angdiff x^2) \Singletan \GdVar
	$
	$
	+
	\smoothfunction(\GdVar, \angLie_{\Tanset}^{\leq 1} \ginversesphere,\angdiff x^1,\angdiff x^2) \Singletan \GdVar
	$
	$
	+
	\smoothfunction(\GdVar, \ginversesphere,\angdiff \Tanset^{\leq 1} x^1, \angdiff \Tanset^{\leq 1} x^2) 
	\Singletan \GdVar
	$.
	We now apply $\Tanset^{N-1}$ to this identity 
	and use Lemmas~\ref{L:POINTWISEFORRECTANGULARCOMPONENTSOFVECTORFIELDS}
	and 
	\ref{L:POINTWISEESTIMATESFORGSPHEREANDITSDERIVATIVES}
	and the already proven estimates
	$
	\left\| 
		\Fullset_*^{\leq 10;1} \Psi
	\right\|_{L^{\infty}(\Sigma_t^u)}
	\lesssim \varepsilon
	$
	and
	$
	\left\|
		\Fullset_*^{\leq 9;1} \Lunit_{(Small)}^i
	\right\|_{L^{\infty}(\Sigma_t^u)}
	\lesssim \varepsilon
	$,
	which implies that
	$
	\left|
		\Tanset^{N-1} \Lunit \mytr \upchi
	\right|
	\lesssim
	\sum_{i=1}^2
	\left|
			\Tanset^{N+1} \Lunit_{(Small)}^i
	\right|
	+
	\mbox{{\upshape RHS}~\eqref{E:LUNITTANGENTDIFFERENTIATEDLUNITSMALLIMPROVEDPOINTWISE}}
	$,
	where $\Tanset^{N+1}$ contains a factor of $\Lunit$.
	We may commute the factor of $\Lunit$ to the front
	using the commutator estimate \eqref{E:ONERADIALTANGENTIALFUNCTIONCOMMUTATORESTIMATE}
	with $f = L_{(Small)}^i$,
	the already proven estimate for 
	$
	\Tanset^{\leq 10} \Lunit_{(Small)}^i
	$ 
	mentioned above
	(to bound the factor 
	$
	\left|
			\Tanset_*^{[1,\lfloor N/2 \rfloor]} f
	\right|
	$
	from the second line of RHS~\eqref{E:ONERADIALTANGENTIALFUNCTIONCOMMUTATORESTIMATE}
	by $\lesssim \varepsilon$),
	and the bootstrap assumptions,
	which yields
	$
	\left|
		\Tanset^{N+1} \Lunit_{(Small)}^i
	\right|
	\lesssim
	\left|
		\Lunit \Tanset^N \Lunit_{(Small)}^i
	\right|
	+
	\mbox{{\upshape RHS}~\eqref{E:LUNITTANGENTDIFFERENTIATEDLUNITSMALLIMPROVEDPOINTWISE}}
	$.
	Moreover, we have already shown that
	$
	\left|
		\Lunit \Tanset^N \Lunit_{(Small)}^i
	\right|
	\lesssim
	\mbox{{\upshape RHS}~\eqref{E:LUNITTANGENTDIFFERENTIATEDLUNITSMALLIMPROVEDPOINTWISE}}
	$.
	We have thus proved
	the estimate for 
	$\left|
		\Tanset^{N-1} \Lunit \mytr \upchi
	\right|
	$
	stated in \eqref{E:LUNITTANGENTDIFFERENTIATEDLUNITSMALLIMPROVEDPOINTWISE}.
	To obtain the same estimate for
	$
	\left|
		\Lunit \Tanset^{N-1} \mytr \upchi
	\right|
	$,
	we use the commutator estimate 
	\eqref{E:PURETANGENTIALFUNCTIONCOMMUTATORESTIMATE} with $f = \mytr \upchi$,
	\eqref{E:POINTWISEESTIMATESFORGSPHEREANDITSTANGENTIALDERIVATIVES},
	the already proven estimates
	$
	\left\| 
		\Fullset_*^{\leq 10;1} \Psi
	\right\|_{L^{\infty}(\Sigma_t^u)}
	\lesssim \varepsilon
	$
	and
	$
	\left\|
		\Fullset_*^{\leq 9;1} \Lunit_{(Small)}^i
	\right\|_{L^{\infty}(\Sigma_t^u)}
	\lesssim \varepsilon
	$
	to deduce that
	$\left|
		\Lunit \Tanset^{N-1} \mytr \upchi
	\right|
	\lesssim 
	\left|
		\Tanset^{N-1} \Lunit \mytr \upchi
	\right|
	+
	\varepsilon
	\left|
		\Tanset^{\leq N} \GdVar
	\right|
	$.
	The desired bound \eqref{E:LUNITTANGENTDIFFERENTIATEDLUNITSMALLIMPROVEDPOINTWISE} for
	$
	\left|
		\Lunit \Tanset^{N-1} \mytr \upchi
	\right|
	$
	now follows from this estimate and the
	one we established for
	$
	\left|
		\Tanset^{N-1} \Lunit \mytr \upchi
	\right|
	$
	just above.

	\medskip
	\noindent \textbf{Proof sketch of the estimate \eqref{E:LUNITONERADIALTANGENTDIFFERENTIATEDLUNITSMALLIMPROVEDPOINTWISE}
	for 
	$\Lunit \Fullset^{N-1;1} \mytr \upchi$
	and $\Fullset^{N-1;1} \Lunit \mytr \upchi$:}
	The proof is much like the proof
	of the estimates for 
	$\Tanset^{N-1} \Lunit \mytr \upchi$
	and $\Lunit \Tanset^{N-1} \mytr \upchi$
	given in the previous paragraph.
	The only notable change is that we must use the commutator estimate \eqref{E:ONERADIALTANGENTIALFUNCTIONCOMMUTATORESTIMATE}
	with $f = \mytr \upchi$
	(in place of the one \eqref{E:PURETANGENTIALFUNCTIONCOMMUTATORESTIMATE} used in the previous paragraph)
	in order to obtain the estimate for 
	$\Lunit \Fullset^{N-1;1} \mytr \upchi$
	from the one for $\Fullset^{N-1;1} \Lunit \mytr \upchi$.
	We remark that all factors leading to the gain of the factor
	$\varepsilon$ on RHS~\eqref{E:LUNITONERADIALTANGENTDIFFERENTIATEDLUNITSMALLIMPROVEDPOINTWISE}
	have already been bounded in $\| \cdot \|_{L^{\infty}}$ by $\lesssim \varepsilon$.


\end{proof}

The following corollary is an immediate consequence of the fact that we have improved the auxiliary
bootstrap assumptions by showing that they hold with $\varepsilon^{1/2}$ replaced by $C \varepsilon$.

\begin{corollary}[\textbf{$\varepsilon^{1/2}$ can be replaced by $C \varepsilon$}]
	\label{C:SQRTEPSILONTOCEPSILON}
	All prior inequalities whose right-hand sides feature an explicit factor of $\varepsilon^{1/2}$
	remain true with $\varepsilon^{1/2}$ replaced by $C \varepsilon$.
\end{corollary}

\section{\texorpdfstring{$L^{\infty}$}{Essential Sup-Norm} Estimates Involving Higher Transversal Derivatives}
\label{S:LINFINITYESTIMATESFORHIGHERTRANSVERSAL}
Our energy estimates are difficult to derive
when $\upmu$ is small because
some products in the energy identities contain
the dangerous factor $1/\upmu$.
In order to control the degeneracy, 
we rely on the estimate $\| \Rad \Rad \upmu \|_{L^{\infty}(\Sigma_t^u)} \lesssim 1$.
In particular, we use this estimate
in proving inequality \eqref{E:UNIFORMBOUNDFORMRADMUOVERMU}
(see the estimate \eqref{E:RADMUOVERMUALGEBRAICBOUND}),
\emph{which is essential for showing that the low-order energies
do not blow up as $\upmu \to 0$.}
We derive the bound
$\| \Rad \Rad \upmu \|_{L^{\infty}(\Sigma_t^u)} \lesssim 1$
by commuting the evolution equation \eqref{E:UPMUFIRSTTRANSPORT} for $\upmu$
with up to two factors of $\Rad$.
Since RHS~\eqref{E:UPMUFIRSTTRANSPORT} depends on 
$\Rad \Psi$
and 
$\Lunit_{(Small)}^i$,
in order to derive the desired bound,
we must obtain estimates for
$\| \Rad \Rad \Rad \Psi \|_{L^{\infty}(\Sigma_t^u)}$,
$\| \Rad \Rad \Lunit_{(Small)}^i \|_{L^{\infty}(\Sigma_t^u)}$,
etc. We provide the necessary estimates in Sect.~\ref{S:LINFINITYESTIMATESFORHIGHERTRANSVERSAL}.
The main result is Prop.~\ref{P:IMPROVEMENTOFHIGHERTRANSVERSALBOOTSTRAP}.

\subsection{Auxiliary bootstrap assumptions}
\label{SS:BOOTSTRAPFORHIGHERTRANSVERSAL}
To facilitate the analysis, we introduce the following 
auxiliary bootstrap assumptions.
In Prop.~\ref{P:IMPROVEMENTOFHIGHERTRANSVERSALBOOTSTRAP}, we derive strict improvements of the assumptions based on our assumptions
\eqref{E:PSIDATAASSUMPTIONS} on the data.

\noindent \underline{\textbf{Auxiliary bootstrap assumptions for small quantities}.}
We assume that the following inequalities hold on $\mathcal{M}_{\Tboot,U_0}$
(see Subsect.~\ref{SS:STRINGSOFCOMMUTATIONVECTORFIELDS} regarding the vectorfield operator notation):
\begin{align}
	\left\| 
		\Fullset_*^{\leq 4;2} \Psi 
	\right\|_{L^{\infty}(\Sigma_t^u)}
	& \leq \varepsilon^{1/2},
		\label{E:HIGHERTRANSVERSALPSIMIXEDFUNDAMENTALLINFINITYBOUNDBOOTSTRAP} 
		\tag{$\mathbf{BA'}1\Psi$}
\end{align}

\begin{subequations}
\begin{align} \label{E:UPMUONERADIALNOTPURERADIALBOOTSTRAP} \tag{$\mathbf{BA'}1\upmu$}
	\left\| 
		\Rad \GeoAng \upmu
	\right\|_{L^{\infty}(\Sigma_t^u)},
		\,
	\left\| 
		\Rad \Lunit \Lunit \upmu
	\right\|_{L^{\infty}(\Sigma_t^u)},
		\,
	\left\| 
		\Rad \GeoAng \GeoAng \upmu
	\right\|_{L^{\infty}(\Sigma_t^u)},
		\,
	\left\| 
		\Rad \Lunit \GeoAng \upmu
	\right\|_{L^{\infty}(\Sigma_t^u)}
	& \leq \varepsilon^{1/2},
\end{align}
\end{subequations}
and
\begin{align} \label{E:PERMUTEDUPMUONERADIALNOTPURERADIALBOOTSTRAP} \tag{$\mathbf{BA''}1\upmu$}
	\eqref{E:UPMUONERADIALNOTPURERADIALBOOTSTRAP} \mbox{ also holds for all permutations of the vectorfield operators on LHS } \eqref{E:UPMUONERADIALNOTPURERADIALBOOTSTRAP},
\end{align}

\begin{align} \label{E:UPMUTWORADIALNOTPURERADIALBOOTSTRAP} \tag{$\mathbf{BA'}1\Lunit_{(Small)}$}
	\left\|
		\Fullset_*^{\leq 3;2}
		\Lunit_{(Small)}^i
	\right\|_{L^{\infty}(\Sigma_t^u)}
	& \leq \varepsilon^{1/2}.
\end{align}

\noindent \underline{\textbf{Auxiliary bootstrap assumptions for quantities that are allowed to be large}.}
We assume that the following inequalities hold on $\mathcal{M}_{\Tboot,U_0}$:
\begin{align}
	\left\| 
		\Rad^M \Psi 
	\right\|_{L^{\infty}(\Sigma_t^u)}
	& \leq 
	\left\| 
		\Rad^M \Psi 
	\right\|_{L^{\infty}(\Sigma_0^u)}
	+ \varepsilon^{1/2},
	&& (2 \leq M \leq 3),
		\label{E:HIGHERPSITRANSVERSALFUNDAMENTALC0BOUNDBOOTSTRAP} 
		\tag{$\mathbf{BA'}2\Psi$}
\end{align}

\begin{align}
	\left\| 
		\Lunit \Rad^M \upmu
	\right\|_{L^{\infty}(\Sigma_t^u)}
	& \leq 
		\frac{1}{2}
		\left\| 
			\Rad^M
			\left\lbrace
				G_{\Lunit \Lunit} \Rad \Psi
			\right\rbrace
		\right\|_{L^{\infty}(\Sigma_0^u)}
		+ \varepsilon^{1/2},
		&& (1 \leq M \leq 2),
			\label{E:HIGHERLUNITUPMUBOOT}  \tag{$\mathbf{BA'}2\upmu$}  \\
		\left\| 
			\Rad^M \upmu
		\right\|_{L^{\infty}(\Sigma_t^u)}
		& \leq
	 	\left\| 
				\Rad^M \upmu
		\right\|_{L^{\infty}(\Sigma_0^u)}
		+ 
		2 \TranminusdatasizeWithFactor^{-1} 
		\left\| 
			\Rad^M
			\left\lbrace
				G_{\Lunit \Lunit} \Rad \Psi
			\right\rbrace
		\right\|_{L^{\infty}(\Sigma_0^u)}
		+ \varepsilon^{1/2},
		&& (1 \leq M \leq 2),
			\label{E:HIGHERUPMUTRANSVERSALBOOT} 
			\tag{$\mathbf{BA'}3\upmu$} \\
		\left\| 
			\Rad \Rad \Lunit_{(Small)}^i
		\right\|_{L^{\infty}(\Sigma_t^u)}
		& \leq
	 	\left\| 
				\Rad \Rad \Lunit_{(Small)}^i
		\right\|_{L^{\infty}(\Sigma_0^u)}
		+ \varepsilon^{1/2}.
		 &&
			\label{E:HIGHERLUNITITRANSVERSALBOOT}  
			\tag{$\mathbf{BA'}2\Lunit_{(Small)}$} 
\end{align}

\subsection{Commutator estimates involving two transversal derivatives}
\label{SS:TWORADDERIVATIVESCOMMUTATORESTIMATES}
In this subsection, we provide some basic commutation estimates that
complement those of Subsect.~\ref{SS:COMMUTATORESTIMATES}.

\begin{lemma}[\textbf{Mixed $\mathcal{P}_u-$transversal-tangent commutator estimates involving two $\Rad$ derivatives}]
		\label{L:HIGHERTRANSVERALTANGENTIALCOMMUTATOR}
		Let $\Fullset^{\vec{I}}$ be a 
		$\Fullset-$multi-indexed operator containing
		exactly two $\Rad$ factors, and assume that
		$3 \leq |\vec{I}| := N+1 \leq 4$. 
		Let $\vec{I}'$ be any permutation of $\vec{I}$.
		Under the data-size and bootstrap assumptions
		of Subsects.~\ref{SS:SIZEOFTBOOT}-\ref{SS:PSIBOOTSTRAP} and Subsect.~\ref{SS:BOOTSTRAPFORHIGHERTRANSVERSAL}
		and the smallness assumptions of Subsect.~\ref{SS:SMALLNESSASSUMPTIONS}, 
		the following commutator 
		estimates hold for functions $f$ on $\mathcal{M}_{\Tboot,U_0}$
		(see Subsect.~\ref{SS:STRINGSOFCOMMUTATIONVECTORFIELDS} regarding the vectorfield operator notation):
		\begin{align}
		\left|
			\Fullset^{\vec{I}} f
			-
			\Fullset^{\vec{I}'} f
		\right|
		& \lesssim
			\left|
				\GeoAng \Fullset^{\leq N-1;1} f
			\right|
			+
			\underbrace{
			\varepsilon
			\left|
				\GeoAng \Fullset^{\leq N-1;2} f
			\right|}_{\mbox{\upshape Absent if $N=2$}}.
				\label{E:TWORADIALTANGENTIALFUNCTIONCOMMUTATORESTIMATE} 
		\end{align}

		Moreover, we have
		\begin{align}
		\left|
			[\angLap, \Rad \Rad] f
		\right|
		& \lesssim
			\left|
				\GeoAng \Fullset^{\leq 2;1} f
			\right|.
				\label{E:ANGLAPTWORADIALTANGENTIALFUNCTIONCOMMUTATOR} 
		\end{align}

\end{lemma}

\begin{proof}
	See Subsect.~\ref{SS:OFTENUSEDESTIMATES} for some comments on the analysis.
	To prove \eqref{E:TWORADIALTANGENTIALFUNCTIONCOMMUTATORESTIMATE},
	we split the terms on RHS~\eqref{E:ALGEBRAICSUBTRACTINGTWOPERMUTEDORDERNVECTORFIELDS}
	into the case where at most one $\Rad$ derivative falls
	on $f$ and the case where both $\Rad$ derivatives fall on $f$. In the former case,
	at most two derivatives fall on the deformation tensors, while in the latter case,
	at most one derivative (which must be $\mathcal{P}_u-$tangent) falls on them. 
	We thus find that
	\begin{align} \label{E:FIRSTESTIMATETWORADIALTANGENTIALFUNCTIONCOMMUTATORESTIMATE}
			\mbox{LHS } \eqref{E:TWORADIALTANGENTIALFUNCTIONCOMMUTATORESTIMATE} 
			& \lesssim
			\left\lbrace
				\left|
					\angLie_{\Fullset}^{\leq 2;2} \angdeformoneformupsharparg{\GeoAng}{\Lunit}
				\right|
				+
				\left|
					\angLie_{\Fullset}^{\leq 2;1} \angdeformoneformupsharparg{\Rad}{\Lunit}
				\right|
				+
				\left|
					\angLie_{\Fullset}^{\leq 2;1} \angdeformoneformupsharparg{\GeoAng}{\Rad}
				\right|
			\right\rbrace
			\left|
				\GeoAng \Fullset^{\leq N-1;1} f
			\right|
			 \\
			&
			\ \ 
			+
			\left\lbrace
				\left|
					\angLie_{\Tanset}^{\leq 1} \angdeformoneformupsharparg{\GeoAng}{\Lunit}
				\right|
				+
				\left|
					\angLie_{\Tanset}^{\leq 1} \angdeformoneformupsharparg{\GeoAng}{\Rad}
				\right|
			\right\rbrace
			\left|
				\GeoAng \Fullset^{\leq N-1;2} f
			\right|.
			\notag
	\end{align}
	From the identities
	\eqref{E:RADDEFORMSPHERERAD},
	\eqref{E:GEOANGDEFORMSPHEREL},
	and
	\eqref{E:GEOANGDEFORMSPHERERAD}
	and Lemma~\ref{L:SCHEMATICDEPENDENCEOFMANYTENSORFIELDS},
	we deduce that the terms in braces on 
	the first line of RHS~\eqref{E:FIRSTESTIMATETWORADIALTANGENTIALFUNCTIONCOMMUTATORESTIMATE} are 
	$
	\lesssim
	\left|
		\angLie_{\Fullset}^{\leq 2;2} \smoothfunction(\Tanset^{\leq 1} \GdVar,\ginversesphere,\angdiff x^1,\angdiff x^2) 
	\right|
	+ 
	\left|
		\angLie_{\Fullset}^{\leq 2;1} \smoothfunction(\Tanset^{\leq 1} \BadVar,\ginversesphere,\angdiff x^1,\angdiff x^2,\Fullset^{\leq 1} \Psi)
	\right|
	$.
	We now show that both terms from the previous inequality are $\lesssim 1$,
	which yields the desired bound
	(this is a simple estimate, 
	where the main point that requires demonstration is that all terms in the braces are sufficiently regular
	such that we have control of their relevant derivatives in $L^{\infty}$).
	To handle the first term from the previous inequality,
	we first commute $\angLie_{\Fullset}^{\leq 2;2}$ under $\angdiff$ and use that
	$Z x^i = Z^i = \smoothfunction(\BadVar)$ for $Z \in \Fullset$
	to bound factors involving the derivatives of $\angdiff x$ by
	$\lesssim \left| \Fullset^{\leq 2;1} \BadVar \right|$.
	Moreover, using \eqref{E:CONNECTIONBETWEENANGLIEOFGSPHEREANDDEFORMATIONTENSORS} 
	and the fact that $\gsphere = \smoothfunction(\GdVar,\angdiff x^1,\angdiff x^2)$,
	we deduce that
	$
	\left|\angLie_{\Fullset}^{\leq 2;2}  \ginversesphere \right|
	\lesssim 
	\left| \Fullset^{\leq 2;2} \GdVar \right|
	+
	\left| \Fullset^{\leq 2;1} \BadVar \right|
	$.
	We thus find that
	$
	\left|
		\angLie_{\Fullset}^{\leq 2;2} \smoothfunction(\Tanset^{\leq 1} \GdVar,\ginversesphere,\angdiff x^1,\angdiff x^2) 
	\right|
	\lesssim 
	\left| \Fullset^{\leq 3;2} \GdVar \right|
	+
	\left| \Fullset^{\leq 2;1} \BadVar \right|
	$.
	From the $L^{\infty}$ estimates of Prop.~\ref{P:IMPROVEMENTOFAUX}
	and the bootstrap assumptions of Subsect.~\ref{SS:BOOTSTRAPFORHIGHERTRANSVERSAL},
	we deduce that the RHS of the previous inequality is $\lesssim 1$ as desired.
	Similar reasoning yields that
	$
	\left|
		\angLie_{\Fullset}^{\leq 2;1} \smoothfunction(\Tanset^{\leq 1} \BadVar,\ginversesphere,\angdiff x^1,\angdiff x^2,\Fullset^{\leq 1} \Psi)
	\right|
	\lesssim 1
	$,
	which completes the proof of the bound for the terms in braces on the first 
	line of RHS~\eqref{E:FIRSTESTIMATETWORADIALTANGENTIALFUNCTIONCOMMUTATORESTIMATE}.
	To handle the terms in braces on the second line of RHS~\eqref{E:FIRSTESTIMATETWORADIALTANGENTIALFUNCTIONCOMMUTATORESTIMATE},
	we use Lemma~\ref{L:MOREPRECISEANGSPHERELESTIMATES}
	and the $L^{\infty}$ estimates of Prop.~\ref{P:IMPROVEMENTOFAUX}
	to bound them by $\lesssim \varepsilon$.
	We have thus proved \eqref{E:TWORADIALTANGENTIALFUNCTIONCOMMUTATORESTIMATE}.

	The proof of \eqref{E:ANGLAPTWORADIALTANGENTIALFUNCTIONCOMMUTATOR}
	is similar and relies on the commutation identity 
	\eqref{E:COMMUTINGANGANGLAPANDLIEZ}
	and the estimates of Lemma~\ref{L:POINTWISEESTIMATESFORGSPHEREANDITSDERIVATIVES};
	we omit the details.

\end{proof}

\subsection{The main estimates involving higher-order transversal derivatives}
\label{SS:HIGHERORDERTRANSVERALMAINESTIMATES}
In the next proposition, we provide the main estimates of Sect.~\ref{S:LINFINITYESTIMATESFORHIGHERTRANSVERSAL}.
In particular, the proposition yields strict improvements of the bootstrap assumptions of Subsect.~\ref{SS:BOOTSTRAPFORHIGHERTRANSVERSAL}.

\begin{proposition}[\textbf{$L^{\infty}$ estimates involving higher-order transversal derivatives}]
	\label{P:IMPROVEMENTOFHIGHERTRANSVERSALBOOTSTRAP}
 		Under the data-size and bootstrap assumptions 
		of Subsects.~\ref{SS:SIZEOFTBOOT}-\ref{SS:PSIBOOTSTRAP} and Subsect.~\ref{SS:BOOTSTRAPFORHIGHERTRANSVERSAL}
		and the smallness assumptions of Subsect.~\ref{SS:SMALLNESSASSUMPTIONS}, 
		the following estimates hold
		on $\mathcal{M}_{\Tboot,U_0}$
		(see Subsect.~\ref{SS:STRINGSOFCOMMUTATIONVECTORFIELDS} regarding the vectorfield operator notation):

	\noindent \underline{\textbf{$L^{\infty}$ estimates involving two or three transversal derivatives of $\Psi$}.}
	\begin{subequations}
	\begin{align} 
	\left\| 
		\Lunit \Fullset_*^{\leq 4;2} \Psi
	\right\|_{L^{\infty}(\Sigma_t^u)}
	& \leq C \varepsilon,
		\label{E:IMPROVEDTRANSVERALESTIMATESFORLUNITTWORADBUTNOTPURERAD} 
		\\
	\left\| 
		\Fullset_*^{\leq 4;2} \Psi
	\right\|_{L^{\infty}(\Sigma_t^u)}
	& \leq C \varepsilon,
		\label{E:IMPROVEDTRANSVERALESTIMATESFORTWORADBUTNOTPURERAD} 
		\\
	\left\| 
		\Rad \Rad \Psi
	\right\|_{L^{\infty}(\Sigma_t^u)}
	& \leq 
		\left\| 
			\Rad \Rad \Psi
		\right\|_{L^{\infty}(\Sigma_0^u)}
		+
		C \varepsilon,
	\label{E:IMPROVEDTRANSVERALESTIMATESFORTWORAD}
		\\
	\left\| 
		\Lunit \Rad \Rad \Rad \Psi
	\right\|_{L^{\infty}(\Sigma_t^u)}
	& \leq 
		C \varepsilon,
	\label{E:IMPROVEDTRANSVERALESTIMATESFORLUNITTHREERAD}
		\\
	\left\| 
		\Rad \Rad \Rad \Psi
	\right\|_{L^{\infty}(\Sigma_t^u)}
	& \leq 
		\left\| 
			\Rad \Rad \Rad \Psi
		\right\|_{L^{\infty}(\Sigma_0^u)}
		+
		C \varepsilon.
	\label{E:IMPROVEDTRANSVERALESTIMATESFORTHREERAD}
	\end{align}
	\end{subequations}

	\noindent \underline{\textbf{$L^{\infty}$ estimates involving one or two transversal derivatives of $\upmu$}.}
	\begin{subequations}
	\begin{align}
		\left\|
			\Lunit \Rad \upmu
		\right\|_{L^{\infty}(\Sigma_t^u)}
		& \leq
			\frac{1}{2}
			\left\|
				\Rad
				\left(
					G_{\Lunit \Lunit} \Rad \Psi
				\right)
			\right\|_{L^{\infty}(\Sigma_0^u)}
			+
			C \varepsilon,
				\label{E:LUNITRADUPMULINFTY}
					\\
		\left\|
			\Rad \upmu
		\right\|_{L^{\infty}(\Sigma_t^u)}
		& \leq
			\left\|
				\Rad \upmu
			\right\|_{L^{\infty}(\Sigma_0^u)}
			+
			\TranminusdatasizeWithFactor^{-1}
			\left\|
				\Rad
				\left(
					G_{\Lunit \Lunit} \Rad \Psi
				\right)
			\right\|_{L^{\infty}(\Sigma_0^u)}
			+
			C \varepsilon,
				\label{E:RADUPMULINFTY}
\end{align}

\begin{align}
	\left\| 
		\Lunit \Rad \GeoAng \upmu
	\right\|_{L^{\infty}(\Sigma_t^u)},
		\,
	\left\| 
		\Lunit \Rad \Lunit \Lunit \upmu
	\right\|_{L^{\infty}(\Sigma_t^u)},
		\,
	\left\| 
		\Lunit \Rad \GeoAng \GeoAng \upmu
	\right\|_{L^{\infty}(\Sigma_t^u)},
		\,
	\left\| 
		\Lunit \Rad \Lunit \GeoAng \upmu
	\right\|_{L^{\infty}(\Sigma_t^u)}
	& \leq C \varepsilon,
				\label{E:LUNITRADTANGENTIALUPMULINFTY}
				\\
	\left\| 
		\Rad \GeoAng \upmu
	\right\|_{L^{\infty}(\Sigma_t^u)},
		\,
	\left\| 
		\Rad \Lunit \Lunit \upmu
	\right\|_{L^{\infty}(\Sigma_t^u)},
		\,
	\left\| 
		\Rad \GeoAng \GeoAng \upmu
	\right\|_{L^{\infty}(\Sigma_t^u)},
		\,
	\left\| 
		\Rad \Lunit \GeoAng \upmu
	\right\|_{L^{\infty}(\Sigma_t^u)}
	& \leq C \varepsilon,
				\label{E:RADTANGENTIALUPMULINFTY}
\end{align}

\begin{align} \label{E:PERMUTEDRADTANGENTIALUPMULINFTY}
	\eqref{E:LUNITRADTANGENTIALUPMULINFTY}-\eqref{E:RADTANGENTIALUPMULINFTY} \mbox{ also hold for all permutations of the vectorfield 
	operators on the LHS },
\end{align}

\begin{align}
		\left\|
			\Lunit \Rad \Rad \upmu
		\right\|_{L^{\infty}(\Sigma_t^u)}
		& \leq
			\frac{1}{2}
			\left\|
				\Rad
				\Rad
				\left(
					G_{\Lunit \Lunit} \Rad \Psi
				\right)
			\right\|_{L^{\infty}(\Sigma_0^u)}
			+
			C \varepsilon
		\label{E:LUNITRADRADUPMULINFTY},
			\\
		\left\|
			\Rad \Rad \upmu
		\right\|_{L^{\infty}(\Sigma_t^u)}
		& \leq
			\left\|
				\Rad \Rad \upmu
			\right\|_{L^{\infty}(\Sigma_0^u)}
			+
			\TranminusdatasizeWithFactor^{-1}
			\left\|
				\Rad
				\Rad
				\left(
					G_{\Lunit \Lunit} \Rad \Psi
				\right)
			\right\|_{L^{\infty}(\Sigma_0^u)}
			+
			C \varepsilon.
		\label{E:RADRADUPMULINFTY}
	\end{align}
	\end{subequations}

	\noindent \underline{\textbf{$L^{\infty}$ estimates involving one or two transversal derivatives of $\Lunit_{(Small)}^i$}.}
	\begin{subequations}
	\begin{align} 
		\left\| 
		\Fullset_*^{\leq 3;2} \Lunit_{(Small)}^i
	\right\|_{L^{\infty}(\Sigma_t^u)}
	& \leq C \varepsilon,
		\label{E:TANGENTIALANDTWORADAPPLIEDTOLUNITILINFINITY}
		\\
	\left\| 
		\Rad \Rad \Lunit_{(Small)}^i
	\right\|_{L^{\infty}(\Sigma_t^u)}
	& \leq 
	\left\| 
		\Rad \Rad \Lunit_{(Small)}^i
	\right\|_{L^{\infty}(\Sigma_0^u)}
	+
	C \varepsilon.
	\label{E:TWORADAPPLIEDTOLUNITILINFINITY}
	\end{align}
	\end{subequations}

	\noindent \underline{\textbf{Sharp pointwise estimates involving the critical factor $G_{\Lunit \Lunit}$}.}
	Moreover, if $0 \leq M \leq 2$
	and $0 \leq s \leq t < \Tboot$, 
	then we have the following estimates:
	\begin{align}
		\left|
			\Rad^M G_{\Lunit \Lunit}(t,u,\vartheta)
			-
			\Rad^M G_{\Lunit \Lunit}(s,u,\vartheta)
		\right|
		& \leq C \varepsilon (t - s),
			\label{E:RADDERIVATIVESOFGLLDIFFERENCEBOUND} \\
		\left|
			\Rad^M 
			\left\lbrace
				G_{\Lunit \Lunit}
				\Rad \Psi
			\right\rbrace
			(t,u,\vartheta)
			-
			\Rad^M 
			\left\lbrace
				G_{\Lunit \Lunit}
				\Rad \Psi
			\right\rbrace
			(s,u,\vartheta)
		\right|
		& \leq C \varepsilon (t - s).
		\label{E:RADDERIVATIVESOFGLLRADPSIDIFFERENCEBOUND}
	\end{align}

	Furthermore, with $\Lunit_{(Flat)} = \partial_t + \partial_1$,
	we have
	\begin{align} \label{E:LUNITUPMUDOESNOTDEVIATEMUCHFROMTHEDATA}
		\Lunit \upmu(t,u,\vartheta)
		& = \frac{1}{2} G_{\Lunit_{(Flat)} \Lunit_{(Flat)}}(\Psi = 0) \Rad \Psi(t,u,\vartheta)
			+ \mathcal{O}(\varepsilon),
	\end{align}
	where 
	$G_{\Lunit_{(Flat)} \Lunit_{(Flat)}}(\Psi = 0)$ is a non-zero constant.
\end{proposition}

\begin{remark}[\textbf{The auxiliary bootstrap assumptions of Subsect.~\ref{SS:BOOTSTRAPFORHIGHERTRANSVERSAL}
  are now redundant}]
  \label{R:HIGHERTRANSVERSALREDUNDANTAUXILIARYBOOTSTRAPASSUMPTIONS}
	Since Prop.~\ref{P:IMPROVEMENTOFHIGHERTRANSVERSALBOOTSTRAP} in particular
	provides an improvement of the auxiliary bootstrap assumptions of
	Subsect.~\ref{SS:BOOTSTRAPFORHIGHERTRANSVERSAL},
	we do not bother to include those bootstrap assumptions
	in the hypotheses of any 
	of the lemmas or propositions proved in the remainder of the article.
\end{remark}

\begin{proof}[Proof of Prop.~\ref{P:IMPROVEMENTOFHIGHERTRANSVERSALBOOTSTRAP}]
	See Subsect.~\ref{SS:OFTENUSEDESTIMATES} for some comments on the analysis.
	We must derive the estimates in a viable order.
	Throughout this proof, 
	we use the estimates of Lemma~\ref{L:BEHAVIOROFEIKONALFUNCTIONQUANTITIESALONGSIGMA0}
	and the assumption \eqref{E:DATAEPSILONVSBOOTSTRAPEPSILON}
	without explicitly mentioning them each time. We refer to these as ``conditions on the data.''
	Similarly, when we say that we use ``the bootstrap assumptions,''
	we mean the assumptions 
	stated in 
	Subsect.~\ref{SS:BOOTSTRAPFORHIGHERTRANSVERSAL}.

	\medskip
	\noindent \textbf{Proof of 
	\eqref{E:IMPROVEDTRANSVERALESTIMATESFORLUNITTWORADBUTNOTPURERAD}-\eqref{E:IMPROVEDTRANSVERALESTIMATESFORTWORAD}:}
	We may assume that the operator $\Fullset_*^{\leq 4;2}$ contains two factors of $\Rad$,
	since otherwise the desired estimates are implied by \eqref{E:PSIMIXEDTRANSVERSALTANGENTBOOTSTRAPIMPROVED}.
	To proceed, we commute the wave equation \eqref{E:WAVEEQUATIONTRANSPORTINTERPRETATION}
	with $\Rad \Tanset^M$, ($0 \leq M \leq 2$), 
	and use Lemmas~\ref{L:POINTWISEFORRECTANGULARCOMPONENTSOFVECTORFIELDS}
	and
	\ref{L:POINTWISEESTIMATESFORGSPHEREANDITSDERIVATIVES},
	the $L^{\infty}$ estimates of Prop.~\ref{P:IMPROVEMENTOFAUX},
	and the bootstrap assumptions to deduce that
	\begin{align} \label{E:WAVEEQNONCETRANSVERSALLYCOMMUTEDTRANSPORTPOINTOFVIEW}
		\left|
			\Lunit \Rad \Tanset^M \Rad \Psi
		\right|
		& \lesssim
				\left|
					\Fullset_*^{\leq M+3;1} \Psi
				\right|
				+
				\left|
					\Fullset_*^{\leq M+2;1} \GdVar
				\right|
			+ \left|
					[\angLap, \Rad \Tanset^M] \Psi
				\right|
			+ 
		\left|
			[\Lunit, \Rad \Tanset^M] \Rad \Psi
		\right|.
	\end{align}
	The $L^{\infty}$ estimates of Prop.~\ref{P:IMPROVEMENTOFAUX}
	imply that the first two terms on RHS~\eqref{E:WAVEEQNONCETRANSVERSALLYCOMMUTEDTRANSPORTPOINTOFVIEW} are
	$\lesssim \varepsilon$.
	Moreover, using in addition the commutator estimate \eqref{E:ANGLAPONERADIALTANGENTIALFUNCTIONCOMMUTATOR} with $f = \Psi$,
	we see that 
	$\left|
		[\angLap, \Rad \Tanset^M] \Psi
	\right|$
	is $\lesssim$ the first term on RHS~\eqref{E:WAVEEQNONCETRANSVERSALLYCOMMUTEDTRANSPORTPOINTOFVIEW}
	and hence $\lesssim \varepsilon$ as well.
	To bound $[\Lunit , \Rad \Tanset^M] \Rad \Psi$,
	we use the commutator estimate \eqref{E:ONERADIALTANGENTIALFUNCTIONCOMMUTATORESTIMATE}
	with $f = \Rad \Psi$,
	Cor.~\ref{C:SQRTEPSILONTOCEPSILON},
	and the $L^{\infty}$ estimates of Prop.~\ref{P:IMPROVEMENTOFAUX}
	to deduce
	$
	\left|
		[\Lunit, \Rad \Tanset^M] \Rad \Psi
	\right|
	\lesssim
	\left|
		\Fullset_*^{\leq M+2;1} \Psi
	\right|
	+
	\varepsilon
	\left|
		\Fullset_*^{\leq M+2;2} \Psi
	\right|
	$.
	The estimates of Prop.~\ref{P:IMPROVEMENTOFAUX} 
	and the bootstrap assumptions
	imply that 
	$
	\left|
		\Fullset_*^{\leq M+2;1} \Psi
	\right|
	\lesssim 
	\varepsilon
	$, 
	while the bootstrap assumptions imply that
	$
	\varepsilon
	\left|
		\Fullset_*^{\leq M+2;2} \Psi
	\right|
	\lesssim \varepsilon
	$.
	Combining these estimates, we deduce that
	$\left| \Lunit \Rad \Tanset^M \Rad \Psi \right| \lesssim \varepsilon$.
	Integrating along the integral curves of $\Lunit$ as in \eqref{E:INTEGRATINGALONGINTEGRALCURVES}
	and using the previous estimate, 
	we find that
	$\left\| \Rad \Tanset^M \Rad \Psi \right\|_{L^{\infty}(\Sigma_t^u)}
	\leq 
	\left\| \Rad \Tanset^M \Rad \Psi \right\|_{L^{\infty}(\Sigma_0^u)}
	+ C \varepsilon
	$.
	Using the previous estimate
	and the conditions on the data,
	and using
	\eqref{E:TWORADIALTANGENTIALFUNCTIONCOMMUTATORESTIMATE} 
	with $f = \Psi$,
	the $L^{\infty}$ estimates of Prop.~\ref{P:IMPROVEMENTOFAUX},
	and the bootstrap assumptions to reorder the factors in the operator 
	$\Rad \Tanset^M \Rad$ as desired (up to error terms bounded in
	$\| \cdot \|_{L^{\infty}(\Sigma_t^u)}$ by $\lesssim \varepsilon$),
	we conclude the desired estimates
	\eqref{E:IMPROVEDTRANSVERALESTIMATESFORTWORADBUTNOTPURERAD} 
	and \eqref{E:IMPROVEDTRANSVERALESTIMATESFORTWORAD}.
	Finally, we similarly reorder the factors in 
	$
		\Lunit \Rad \Tanset^M \Rad \Psi
	$
	and 
	use the estimates
	$\left| \Lunit \Rad \Tanset^M \Rad \Psi \right| \lesssim \varepsilon$
	and \eqref{E:IMPROVEDTRANSVERALESTIMATESFORTWORADBUTNOTPURERAD}
	to obtain \eqref{E:IMPROVEDTRANSVERALESTIMATESFORLUNITTWORADBUTNOTPURERAD}.

	\medskip
	\noindent \textbf{Proof of 
	\eqref{E:RADDERIVATIVESOFGLLDIFFERENCEBOUND}-\eqref{E:RADDERIVATIVESOFGLLRADPSIDIFFERENCEBOUND}
	in the cases $0 \leq M \leq 1$:}
	It suffices to prove 
	\begin{align} \label{E:LDERIVATIVEOFRADIALDERIVATIVESOFCRITICALFACTOR}
		\left|
			\Lunit \Rad^M G_{\Lunit \Lunit}
		\right|,
			\,
		\left|
			\Lunit \Rad^M 
			\left\lbrace
				G_{\Lunit \Lunit}
				\Rad \Psi
			\right\rbrace
		\right|
		& \lesssim \varepsilon,
	\end{align}
	for once we have shown \eqref{E:LDERIVATIVEOFRADIALDERIVATIVESOFCRITICALFACTOR},
	we can obtain the desired estimates by integrating
	along the integral curves of $\Lunit$ from time $s$ to $t$
	(in analogy with \eqref{E:INTEGRATINGALONGINTEGRALCURVES})
	and using the estimates \eqref{E:LDERIVATIVEOFRADIALDERIVATIVESOFCRITICALFACTOR}.
	To proceed, we first use Lemma~\ref{L:SCHEMATICDEPENDENCEOFMANYTENSORFIELDS} to deduce that
	$G_{\Lunit \Lunit} = \smoothfunction(\GdVar)$ 
	and $G_{\Lunit \Lunit} \Rad \Psi = \smoothfunction(\GdVar) \Rad \Psi$.
	Hence, to obtain \eqref{E:LDERIVATIVEOFRADIALDERIVATIVESOFCRITICALFACTOR}
	when $M=0$, we differentiate these two identities with $\Lunit$ and 
	use the $L^{\infty}$ estimates of Prop.~\ref{P:IMPROVEMENTOFAUX}
	and the bootstrap assumptions.
	The proof is similar in the case $M=1$,
	but we must also use the estimate
	$\left\|\Lunit \Rad \Rad \Psi \right\|_{L^{\infty}(\Sigma_t^u)} \lesssim \varepsilon$,
	which is a consequence of the previously established estimate 
	\eqref{E:IMPROVEDTRANSVERALESTIMATESFORTWORADBUTNOTPURERAD}.

	\medskip
	\noindent \textbf{Proof of \eqref{E:LUNITUPMUDOESNOTDEVIATEMUCHFROMTHEDATA}:}
	We first use \eqref{E:UPMUFIRSTTRANSPORT},
	the fact that 
	$G_{\Lunit \Lunit}, G_{\Lunit \Radunit} = \smoothfunction(\GdVar)$
	(see Lemma~\ref{L:SCHEMATICDEPENDENCEOFMANYTENSORFIELDS}),
	and the $L^{\infty}$ estimates of Prop.~\ref{P:IMPROVEMENTOFAUX}
	to deduce that 
	$\Lunit \upmu(t,u,\vartheta)
		= 
		\frac{1}{2}
	 [G_{\Lunit \Lunit} \Rad \Psi](t,u,\vartheta)
	 + \mathcal{O}(\varepsilon)
	$.
	Since 
	$\Lunit^0 = \Lunit_{(Flat)}^0 = 1$,
	$\Lunit^i = \Lunit_{(Flat)}^i + \Lunit_{(Small)}^i$,
	and $G_{\alpha \beta} = G_{\alpha \beta}(\Psi = 0) + \mathcal{O}(\Psi)$,
	we can use the $L^{\infty}$ estimates of Prop.~\ref{P:IMPROVEMENTOFAUX}
	to deduce that 
	$G_{\Lunit \Lunit}(t,u,\vartheta) = G_{\Lunit_{(Flat)} \Lunit_{(Flat)}}(\Psi = 0) + \mathcal{O}(\varepsilon)$.
	Combining this estimate with the previous one
	and using \eqref{E:PSITRANSVERSALLINFINITYBOUNDBOOTSTRAPIMPROVED},
	we conclude \eqref{E:LUNITUPMUDOESNOTDEVIATEMUCHFROMTHEDATA}.

	\medskip
	\noindent \textbf{Proof of \eqref{E:LUNITRADUPMULINFTY}-\eqref{E:PERMUTEDRADTANGENTIALUPMULINFTY}:}
	Let $1 \leq K \leq 3$ be an integer.
	We commute equation \eqref{E:UPMUFIRSTTRANSPORT}
	with $\Fullset^{K;1}$ and
	use the aforementioned relations $G_{\Lunit \Lunit}, G_{\Lunit \Radunit} = \smoothfunction(\GdVar)$,
	the $L^{\infty}$ estimates of Prop.~\ref{P:IMPROVEMENTOFAUX},
	and the bootstrap assumptions to deduce
	\begin{align} \label{E:UPMUEVOLUTIONEQUATIONRADANDTANGENTIALCOMMUTEDFIRSTBOUND}
		\left|
			\Lunit \Fullset^{K;1} \upmu
		\right|
		& \leq 
			\frac{1}{2}
			\left|
				\Fullset^{K;1}
				\left\lbrace
					G_{\Lunit \Lunit}
					\Rad \Psi
				\right\rbrace
			\right|
			+
			\left|
				\Fullset_*^{\leq K+1;1} \Psi
			\right|
			+
			\left|
				[\Lunit, \Fullset^{K;1}] \upmu
			\right|.
		\end{align}
		We now show that the last two terms on RHS~\eqref{E:UPMUEVOLUTIONEQUATIONRADANDTANGENTIALCOMMUTEDFIRSTBOUND}
		are $\lesssim \varepsilon$.
		We already proved 
		$
		\left|
			\Fullset_*^{\leq K+1;1} \Psi
		\right|
		\lesssim \varepsilon
		$
		in Prop.~\ref{P:IMPROVEMENTOFAUX}.
		To bound $[\Lunit, \Fullset^{K;1}] \upmu$,
		we use the commutator estimate \eqref{E:ONERADIALTANGENTIALFUNCTIONCOMMUTATORESTIMATE}
		with $f = \upmu$,
		the $L^{\infty}$ estimates of Prop.~\ref{P:IMPROVEMENTOFAUX},
		and Cor.~\ref{C:SQRTEPSILONTOCEPSILON}
		to deduce that
		$
		\left|
			[\Lunit, \Fullset^{K;1}] \upmu
		\right|
		\lesssim 
		\left|
			\Tanset_*^{[1, K]} \upmu
		\right|
		+
		\varepsilon
		\left|
			\GeoAng \Fullset^{\leq K-1;1} \upmu
		\right|
		$.
		The $L^{\infty}$ estimates of Prop.~\ref{P:IMPROVEMENTOFAUX} imply that 
		$
		\left|
			\Tanset_*^{[1,K]} \upmu
		\right|
		\lesssim \varepsilon
		$,
		while the bootstrap assumptions
		imply that
		$
		\varepsilon
		\left|
			\GeoAng \Fullset^{\leq K-1;1} \upmu
		\right|
		\lesssim \varepsilon
		$
		as well. We have thus shown that
		\begin{align} \label{E:UPMUEVOLUTIONEQUATIONRADANDTANGENTIALCOMMUTEDSECONDBOUND}
		\left\|
			\Lunit \Fullset^{K;1} \upmu
		\right\|_{L^{\infty}(\Sigma_t^u)}
		& \leq 
			\frac{1}{2}
			\left\|
				\Fullset^{K;1}
				\left\lbrace
					G_{\Lunit \Lunit}
					\Rad \Psi
				\right\rbrace
			\right\|_{L^{\infty}(\Sigma_t^u)}
			+ 
			C \varepsilon.
		\end{align}
		We split the remainder of the proof into two cases, starting with the case $\Fullset^{K;1} = \Rad$.
		Using the bound \eqref{E:RADDERIVATIVESOFGLLRADPSIDIFFERENCEBOUND} 
		with $s=0$ and $M=1$ (established above),
		we can replace the norm 
		$\| \cdot \|_{L^{\infty}(\Sigma_t^u)}$
		on RHS~\eqref{E:UPMUEVOLUTIONEQUATIONRADANDTANGENTIALCOMMUTEDSECONDBOUND}
		with the norm $\| \cdot \|_{L^{\infty}(\Sigma_0^u)}$
		plus an error term that is 
		bounded in the norm $\| \cdot \|_{L^{\infty}(\Sigma_t^u)}$
		by $\leq C \varepsilon$,
		which yields \eqref{E:LUNITRADUPMULINFTY}.
		Integrating along the integral curves of $\Lunit$ as in \eqref{E:INTEGRATINGALONGINTEGRALCURVES},
		using the resulting estimate for 
		$
		\left\|
			\Lunit \Fullset^{K;1} \upmu
		\right\|_{L^{\infty}(\Sigma_t^u)}
		$,
		and using the assumption
		$\Tboot \leq 2 \TranminusdatasizeWithFactor^{-1}$,
		we conclude \eqref{E:RADUPMULINFTY}.
		In the remaining case, $\Fullset^{K;1}$ is not the operator $\Rad$.
		That is, $K > 1$ and $\Fullset^{K;1}$ must contain a $\mathcal{P}_u-$tangent factor,
		which is equivalent to $\Fullset^{K;1} = \Fullset_*^{K;1}$.
		Recalling that
		$G_{\Lunit \Lunit} \Rad \Psi = \smoothfunction(\GdVar) \Rad \Psi$
		and using the estimates of Prop.~\ref{P:IMPROVEMENTOFAUX},
		the bootstrap assumptions,
		and
		\eqref{E:IMPROVEDTRANSVERALESTIMATESFORTWORADBUTNOTPURERAD},
		we find that
		$
		\left\|
			\Fullset_*^{K;1}
		\left\lbrace
					G_{\Lunit \Lunit}
					\Rad \Psi
				\right\rbrace
			\right\|_{L^{\infty}(\Sigma_t^u)}
			\lesssim \varepsilon
		$.
		Thus, in this case, we have shown that 
		$
		\mbox{RHS~\eqref{E:UPMUEVOLUTIONEQUATIONRADANDTANGENTIALCOMMUTEDSECONDBOUND}}
		\lesssim
		\varepsilon
		$
		as desired.
		Integrating 
		along the integral curves of $\Lunit$ as in \eqref{E:INTEGRATINGALONGINTEGRALCURVES} and using the estimate
		$
		\left\|
			\Lunit \Fullset_*^{K;1} \upmu
		\right\|_{L^{\infty}(\Sigma_t^u)}
		\lesssim 
		\varepsilon
		$
		just obtained, we conclude that
		$
		\left\|
			\Fullset_*^{K;1} \upmu
		\right\|_{L^{\infty}(\Sigma_t^u)}
		\leq
		\left\|
			\Fullset_*^{K;1} \upmu
		\right\|_{L^{\infty}(\Sigma_0^u)}
		+ C \varepsilon
		$.
		All bounds in \eqref{E:LUNITRADTANGENTIALUPMULINFTY}-\eqref{E:PERMUTEDRADTANGENTIALUPMULINFTY}
		now follow from the previous estimate and the conditions on the data
		except for the estimate \eqref{E:PERMUTEDRADTANGENTIALUPMULINFTY} concerning
		the permutations of the vectorfields in \eqref{E:LUNITRADTANGENTIALUPMULINFTY}.
		To obtain the remaining estimate \eqref{E:PERMUTEDRADTANGENTIALUPMULINFTY},
		we use the commutation estimate
		\eqref{E:ONERADIALTANGENTIALFUNCTIONCOMMUTATORESTIMATE}
		with $f = \upmu$,
		the $L^{\infty}$ estimates of Prop.~\ref{P:IMPROVEMENTOFAUX},
		the estimate \eqref{E:RADTANGENTIALUPMULINFTY},
		and the bootstrap assumptions.

\medskip
	\noindent \textbf{Proof of \eqref{E:TANGENTIALANDTWORADAPPLIEDTOLUNITILINFINITY} and \eqref{E:TWORADAPPLIEDTOLUNITILINFINITY}:}
We may assume that the operator $\Fullset_*^{\leq 3;2}$ in \eqref{E:TANGENTIALANDTWORADAPPLIEDTOLUNITILINFINITY}
contains two factors of $\Rad$ since otherwise the desired estimate is implied by \eqref{E:LUNITAPPLIEDTOLISMALLANDLISMALLINFTYESTIMATE}.
To proceed, we express \eqref{E:RADLUNITI} in the schematic form
$
\Rad \Lunit_{(Small)}^i 
= 
\smoothfunction(\GdVar,\ginversesphere,\angdiff x^1,\angdiff x^2) \Rad \Psi 
+ 
\smoothfunction(\BadVar,\ginversesphere,\angdiff x^1,\angdiff x^2) \Singletan \Psi
+ \smoothfunction(\ginversesphere,\angdiff x^1,\angdiff x^2) \angdiff \upmu
$.
We now apply $\Singletan \Rad$ to this identity, where $\Singletan \in \Tanset$.
Using Lemmas~\ref{L:POINTWISEFORRECTANGULARCOMPONENTSOFVECTORFIELDS}
and
\ref{L:POINTWISEESTIMATESFORGSPHEREANDITSDERIVATIVES},
the $L^{\infty}$ estimates of Prop.~\ref{P:IMPROVEMENTOFAUX},
the already proven estimates 
\eqref{E:IMPROVEDTRANSVERALESTIMATESFORTWORADBUTNOTPURERAD},
\eqref{E:RADTANGENTIALUPMULINFTY},
and
\eqref{E:PERMUTEDRADTANGENTIALUPMULINFTY},
and the bootstrap assumptions, 
we deduce that
\begin{align} \label{E:TANGENTIALRADRADLUNITIFIRSTBOUND}
	\left|
		\Singletan \Rad \Rad \Lunit_{(Small)}^i
	\right|
	& \lesssim
		\left|
			\Fullset_*^{\leq 3;2} \Psi
		\right|
		+
		\left|
			\Fullset_*^{\leq 3;1} \GdVar
		\right|
		+
		\left|
			\GeoAng \Fullset_*^{\leq 2;1} \upmu
		\right|
		+
		\left|
			\Tanset_*^{[1,2]} \upmu
		\right|
		\lesssim \varepsilon.
\end{align}
Also using the commutator estimate \eqref{E:TWORADIALTANGENTIALFUNCTIONCOMMUTATORESTIMATE} 
with $f = \Lunit_{(Small)}^i$ to reorder the factors of the operator $\Singletan \Rad \Rad$
as desired up to error terms bounded in the norm 
$\| \cdot \|_{L^{\infty}(\Sigma_t^u)}$ by $\lesssim \varepsilon$,
we conclude \eqref{E:TANGENTIALANDTWORADAPPLIEDTOLUNITILINFINITY}.
Moreover, a special case of \eqref{E:TANGENTIALANDTWORADAPPLIEDTOLUNITILINFINITY}
is the bound $\left|\Lunit \Rad \Rad \Lunit_{(Small)}^i \right| \lesssim \varepsilon$.
Integrating along the integral curves of $\Lunit$ 
as in \eqref{E:INTEGRATINGALONGINTEGRALCURVES}
and using the previous estimate,
we conclude \eqref{E:TWORADAPPLIEDTOLUNITILINFINITY}.

\medskip
\noindent \textbf{Proof of \eqref{E:IMPROVEDTRANSVERALESTIMATESFORLUNITTHREERAD} and \eqref{E:IMPROVEDTRANSVERALESTIMATESFORTHREERAD}:}
	We commute equation \eqref{E:WAVEEQUATIONTRANSPORTINTERPRETATION} with $\Rad \Rad$
	and argue as in the proof of \eqref{E:WAVEEQNONCETRANSVERSALLYCOMMUTEDTRANSPORTPOINTOFVIEW} 
	to deduce that
	\begin{align} \label{E:WAVEEQNTWICETRANSVERSALLYCOMMUTEDTRANSPORTPOINTOFVIEW}
		\left|
			\Lunit \Rad \Rad \Rad \Psi
		\right|
		& \lesssim
				\left|
					\Fullset_*^{\leq 4;2} \Psi
				\right|
				+
				\left|
					\Fullset_*^{\leq 3;2} \GdVar
				\right|
			+ \left|
					[\angLap, \Rad \Rad] \Psi
				\right|
			+ 
		\left|
			\Lunit \Rad \Rad \Rad \Psi
			- \Rad \Rad \Lunit \Rad \Psi
		\right|.
	\end{align}
	We clarify that the proof of \eqref{E:WAVEEQNTWICETRANSVERSALLYCOMMUTEDTRANSPORTPOINTOFVIEW}
	requires the bounds 
	$\left|
		\angLie_{\Rad} \angLie_{\Rad} \ginversesphere
	\right|,
	\left|
		\angLie_{\Rad} \angLie_{\Rad} \angdiff x
	\right|
	\lesssim 1
	$, which we obtained in the proof of
	Lemma~\ref{L:HIGHERTRANSVERALTANGENTIALCOMMUTATOR}.
	Next, we note that the already proven estimates
	\eqref{E:IMPROVEDTRANSVERALESTIMATESFORTWORADBUTNOTPURERAD}
	and \eqref{E:TANGENTIALANDTWORADAPPLIEDTOLUNITILINFINITY}
	imply that 
	$
	\left|
		\Fullset_*^{\leq 4;2} \Psi
	\right|,
		\,
	\left|
		\Fullset_*^{\leq 3;2} \GdVar
	\right|
	\lesssim \varepsilon
	$.
	Next, we use \eqref{E:ANGLAPTWORADIALTANGENTIALFUNCTIONCOMMUTATOR}
	with $f = \Psi$
	to bound the commutator term
	$\left|
		[\angLap, \Rad \Rad] \Psi
	\right|$
	by $\lesssim$ the first term on RHS~\eqref{E:WAVEEQNTWICETRANSVERSALLYCOMMUTEDTRANSPORTPOINTOFVIEW}
	(and hence it is $\lesssim \varepsilon$ too).
	Next, we use
 	\eqref{E:TWORADIALTANGENTIALFUNCTIONCOMMUTATORESTIMATE} 
 	with $f=\Rad \Psi$ and $N=2$ to deduce that
 	$
 	\left|
		\Lunit \Rad \Rad \Rad \Psi
		- \Rad \Rad \Lunit \Rad \Psi
	\right|
	\lesssim
	\left|
		\Fullset_*^{\leq 4;2} \Psi
	\right|
	$.
 	As we have mentioned, we already have shown
 	that 
 	$
 	\left|
		\Fullset_*^{\leq 4;2} \Psi
	\right|
	\lesssim \varepsilon
 	$.
 Combining these estimates,
 we deduce that
	$
	\left|
		\Lunit \Rad \Rad \Rad \Psi
	\right|
	\lesssim \varepsilon
	$,
	which implies \eqref{E:IMPROVEDTRANSVERALESTIMATESFORLUNITTHREERAD}.
	Integrating along the integral curves of $\Lunit$ as in \eqref{E:INTEGRATINGALONGINTEGRALCURVES}
	and using the previous estimate,
	we conclude the desired estimate \eqref{E:IMPROVEDTRANSVERALESTIMATESFORTHREERAD}.

	\medskip
\noindent \textbf{Proof of \eqref{E:RADDERIVATIVESOFGLLDIFFERENCEBOUND}-\eqref{E:RADDERIVATIVESOFGLLRADPSIDIFFERENCEBOUND}
	in the case $M=2$:}
	The proof is very similar to the proof given above in the cases
	$M=0,1$, so we only highlight the main new ingredients needed in the case $M=2$:
	we must use the estimates
	$
	\left|
		\Lunit \Rad \Rad \Rad \Psi
	\right|
	\lesssim \varepsilon
	$
	and
	$\left|
		\Lunit \Rad \Rad \Lunit_{(Small)}^i 
	\right| \lesssim \varepsilon
	$
	established in 
	\eqref{E:IMPROVEDTRANSVERALESTIMATESFORLUNITTHREERAD}
	and 
	\eqref{E:TANGENTIALANDTWORADAPPLIEDTOLUNITILINFINITY}
	in order to deduce
	\eqref{E:LDERIVATIVEOFRADIALDERIVATIVESOFCRITICALFACTOR}
	in the case $M=2$.

\medskip
	\noindent \textbf{Proof of \eqref{E:LUNITRADRADUPMULINFTY}-\eqref{E:RADRADUPMULINFTY}:}
	We commute equation \eqref{E:UPMUFIRSTTRANSPORT} with $\Rad \Rad$ 
	and argue as in the proof of \eqref{E:UPMUEVOLUTIONEQUATIONRADANDTANGENTIALCOMMUTEDFIRSTBOUND}
	to obtain
		\begin{align} \label{E:UPMUEVOLUTIONEQUATIONRADDOUBLECOMMUTEDFIRSTBOUND}
		\left|
			\Lunit \Rad \Rad \upmu
		\right|
		& \leq 
			\frac{1}{2}
			\left|
				\Rad \Rad
				\left\lbrace
					G_{\Lunit \Lunit}
					\Rad \Psi
				\right\rbrace
			\right|
			+
			\left|
				\Fullset_*^{\leq 3;2} \Psi
			\right|
			+
			\left|
				\Lunit \Rad \Rad \upmu
				- 
				\Rad \Rad \Lunit \upmu
			\right|.
		\end{align}
		Using the commutator estimate
		\eqref{E:TWORADIALTANGENTIALFUNCTIONCOMMUTATORESTIMATE}
		with $f = \upmu$,
		the $L^{\infty}$ estimates of Prop.~\ref{P:IMPROVEMENTOFAUX}, 
		and the already proven bound \eqref{E:RADTANGENTIALUPMULINFTY},
		we deduce that 
		$\left|
				\Lunit \Rad \Rad \upmu
				- 
				\Rad \Rad \Lunit \upmu
		\right|
		\lesssim
			\left|
				\GeoAng \Fullset^{\leq 1} \upmu
			\right|
		\lesssim 
		\varepsilon
		$.
		Next, we use
		\eqref{E:IMPROVEDTRANSVERALESTIMATESFORTWORADBUTNOTPURERAD}
		to deduce that 
		$
		\left|
			\Fullset_*^{\leq 3;2} \Psi
		\right|
		\lesssim \varepsilon
		$.
		Thus, we have shown that the last two terms 
		on RHS~\eqref{E:UPMUEVOLUTIONEQUATIONRADDOUBLECOMMUTEDFIRSTBOUND} are 
		$\lesssim \varepsilon$. 
		The remainder of the proof of \eqref{E:LUNITRADRADUPMULINFTY}-\eqref{E:RADRADUPMULINFTY}
		now proceeds as in the proof
		of \eqref{E:LUNITRADUPMULINFTY}-\eqref{E:RADUPMULINFTY},
		thanks to the availability of the already proven estimates
		\eqref{E:RADDERIVATIVESOFGLLDIFFERENCEBOUND}-\eqref{E:RADDERIVATIVESOFGLLRADPSIDIFFERENCEBOUND}
		in the case $M=2$.

\end{proof}

\section{Sharp Estimates for \texorpdfstring{$\upmu$}{the Inverse Foliation Density}}
\label{S:SHARPESTIMATESFORUPMU}
In this section, we derive sharp pointwise estimates for $\upmu$
and its derivatives that are far more detailed than those of
Sects.~\ref{S:PRELIMINARYPOINTWISE} and \ref{S:LINFINITYESTIMATESFORHIGHERTRANSVERSAL}.
We use these estimates in Sect.~\ref{S:ENERGYESTIMATES}
when we derive a priori energy estimates. To close the energy estimates,
we must have precise knowledge of how $\upmu$ vanishes, which is
the main information derived in Sect.~\ref{S:SHARPESTIMATESFORUPMU}.

Many results derived in this section are based on a posteriori estimates
in which the behavior of a quantity at times $0 \leq s \leq t$ is
tied to the behavior of other quantities at the ``late time'' $t$, where $t < \Tboot$.
For this reason, some of our analysis refers to quantities 
$q = q(s,u,\vartheta;t)$ that are functions of the geometric coordinates
$(s,u,\vartheta)$ and the ``late time parameter'' $t$.
When we state and derive estimates for such quantities, 
$s$ is the ``moving'' time variable verifying $0 \leq s \leq t$.

\subsection{Auxiliary quantities for analyzing \texorpdfstring{$\upmu$}{the inverse foliation density} and first estimates}
\label{SS:AUXMUQUANTITIES}
We start by defining some quantities that play a role in our analysis of $\upmu$.

\begin{definition}[\textbf{Auxiliary quantities used to analyze $\upmu$}]
	\label{D:AUXQUANTITIES}
	We define the following quantities, 
	where we assume that $0 \leq s \leq t$
	for those quantities that depend on both $s$ and $t:$
	\begin{subequations}
	\begin{align}
	M(s,u,\vartheta;t) 
	& := \int_{s'=s}^{s'=t} 
					\left\lbrace
						\Lunit \upmu(t,u,\vartheta) - \Lunit \upmu(s',u,\vartheta) 
					\right\rbrace
				\, ds',
		\label{E:BIGMDEF} \\
	\mathring{\upmu}(u,\vartheta)
	& := \upmu(s=0,u,\vartheta),
		\\
	\widetilde{M}(s,u,\vartheta;t)
	& := \frac{M(s,u,\vartheta;t)}{\mathring{\upmu}(u,\vartheta) - M(0,u,\vartheta;t)},
		\label{E:WIDETILDEBIGMDEF} \\
	\upmu_{(Approx)}(s,u,\vartheta;t)
		& := 1
			+  \frac{\Lunit \upmu(t,u,\vartheta)}{\mathring{\upmu}(u,\vartheta) - M(0,u,\vartheta;t)}
				 s
			+ \widetilde{M}(s,u,\vartheta;t).
		\label{E:MUAPPROXDEF}
	\end{align}
	\end{subequations}
\end{definition}

The following quantity captures the worst-case smallness
of $\upmu$ along $\Sigma_t^u$. We use it to capture
the degeneracy of our high-order energy estimates.

\begin{definition}[\textbf{Definition of} $\upmu_{\star}$]
	\label{D:MUSTARDEF}
	\begin{align} \label{E:MUSTARDEF}
		\upmu_{\star}(t,u)
		& := \min \lbrace 1, \min_{\Sigma_t^u} \upmu \rbrace.
	\end{align}
\end{definition}

\begin{remark} \label{R:REDUNDANTMIN}
	It is redundant to take the min with $1$ in \eqref{E:MUSTARDEF} because $\upmu \equiv 1$ along $\mathcal{P}_0$;
	we have done this only to emphasize that $\upmu_{\star}(t,u) \leq 1$.
\end{remark}

We now provide some basic estimates for the auxiliary quantities.
\begin{lemma}[\textbf{First estimates for the auxiliary quantities}]
\label{L:FIRSTESTIMATESFORAUXILIARYUPMUQUANTITIES}
Under the data-size and bootstrap assumptions 
of Subsects.~\ref{SS:SIZEOFTBOOT}-\ref{SS:PSIBOOTSTRAP}
and the smallness assumptions of Subsect.~\ref{SS:SMALLNESSASSUMPTIONS}, 
the following  
estimates hold for $(t,u,\vartheta) \in [0,\Tboot) \times [0,U_0] \times \mathbb{T}$
and $0 \leq s \leq t$:
\begin{align}
	\mathring{\upmu}(u,\vartheta)
	& = 1 + \mathcal{O}(\varepsilon),
		\label{E:MUINITIALDATAESTIMATE}
		\\
	\mathring{\upmu}(u,\vartheta)
	& = 1 + M(0,u,\vartheta;t) + \mathcal{O}(\varepsilon).
		\label{E:MUAMPLITUDENEARONE}
\end{align}
In addition, the following pointwise estimates hold:
\begin{align} \label{E:LUNITUPMUATTIMETMINUSLUNITUPMUATTIMESPOINTWISEESTIMATE}
	\left|
		\Lunit \upmu(t,u,\vartheta) 
		- 
		\Lunit \upmu(s,u,\vartheta)
	\right|
	& \lesssim \varepsilon(t-s),
	\\
	|M(s,u,\vartheta;t)|, |\widetilde{M}(s,u,\vartheta;t)|
	& \lesssim 
		\varepsilon (t - s)^2,
		\label{E:BIGMEST} 
\end{align}
\begin{align} \label{E:MUAPPROXMISLIKEMU}
	\upmu(s,u,\vartheta)
	& = (1 + \mathcal{O}(\varepsilon)) 
	\upmu_{(Approx)}(s,u,\vartheta;t).
\end{align}
\end{lemma}

\begin{proof}
	\eqref{E:MUINITIALDATAESTIMATE} follows from \eqref{E:UPMUDATATANGENTIALLINFINITYCONSEQUENCES} 
	and \eqref{E:DATAEPSILONVSBOOTSTRAPEPSILON}.
	To prove \eqref{E:LUNITUPMUATTIMETMINUSLUNITUPMUATTIMESPOINTWISEESTIMATE}, 
	we note that \eqref{E:LUNITAPPLIEDTOTANGENTIALUPMUANDTANSETSTARLINFTY}
	implies that $\left| \Lunit \Lunit \upmu \right| \lesssim \varepsilon$.
	Integrating this estimate along the integral curves of $\Lunit$
	from time $s$ to time $t$, we conclude \eqref{E:LUNITUPMUATTIMETMINUSLUNITUPMUATTIMESPOINTWISEESTIMATE}.
	The estimate \eqref{E:MUAMPLITUDENEARONE}
	and the estimate \eqref{E:BIGMEST} for $M$
	then follow from definition \eqref{E:BIGMDEF} 
	and the estimate \eqref{E:LUNITUPMUATTIMETMINUSLUNITUPMUATTIMESPOINTWISEESTIMATE}.
	\eqref{E:BIGMEST} for $\widetilde{M}$
	follows from definition \eqref{E:WIDETILDEBIGMDEF}, 
	the estimate \eqref{E:BIGMEST} for $M$,
	and \eqref{E:MUAMPLITUDENEARONE}.
To prove \eqref{E:MUAPPROXMISLIKEMU}, we first note the following identity, which is
a straightforward consequence of Def.~\ref{D:AUXQUANTITIES}:
\begin{align} \label{E:MUSPLIT}
	\upmu(s,u,\vartheta) 
	& = \left\lbrace
				\mathring{\upmu}(u,\vartheta) - M(0,u,\vartheta;t)
			\right\rbrace
			\upmu_{(Approx)}(s,u,\vartheta;t).
\end{align}
The desired estimate \eqref{E:MUAPPROXMISLIKEMU} now follows from
\eqref{E:MUSPLIT}	and \eqref{E:MUAMPLITUDENEARONE}.

\end{proof}


To derive sharp estimates for $\upmu$, it is convenient to
distinguish between regions where
$\upmu$ is shrinking and regions where it is not.
This motivates the following definition.

\begin{definition}[\textbf{Regions of distinct $\upmu$ behavior}]
\label{D:REGIONSOFDISTINCTUPMUBEHAVIOR}
For each 
$t \in [0,\Tboot)$,
$s \in [0,t]$, 
and $u \in [0,U_0]$, 
we partition 
\begin{subequations}
\begin{align}
	[0,u] \times \mathbb{T} 
	& = \Vplus{t}{u} \cup \Vminus{t}{u},
		\label{E:OUINTERVALCROSSS2SPLIT} \\
	\Sigma_s^u
	& = \Sigmaplus{s}{t}{u} \cup \Sigmaminus{s}{t}{u},
	\label{E:SIGMASSPLIT}
\end{align}
\end{subequations}
where
\begin{subequations}
\begin{align}
	\Vplus{t}{u}
	& := 
	\left\lbrace
		(u',\vartheta) \in [0,u] \times \mathbb{T} \ | \  
			\frac{\Lunit \upmu(t,u',\vartheta)}{\mathring{\upmu}(u',\vartheta) - M(0,u',\vartheta;t)}
		\geq 0
	\right\rbrace,
		\label{E:ANGLESANDUWITHNONDECAYINUPMUGBEHAVIOR} \\
	\Vminus{t}{u}
	& := 
	\left\lbrace
		(u',\vartheta) \in [0,u] \times \mathbb{T} \ | \ 
			\frac{\Lunit \upmu(t,u',\vartheta)}{\mathring{\upmu}(u',\vartheta) - M(0,u',\vartheta;t)} < 0
	\right\rbrace,
		\label{E:ANGLESANDUWITHDECAYINUPMUGBEHAVIOR} \\
	\Sigmaplus{s}{t}{u}
	& := 
	\left\lbrace
		(s,u',\vartheta) \in \Sigma_s^u \ | \ (u',\vartheta) \in \Vplus{t}{u}
	\right\rbrace,
		\label{E:SIGMAPLUS} \\
	\Sigmaminus{s}{t}{u}
	& := 
	\left\lbrace
		(s,u',\vartheta) \in \Sigma_s^u \ | \ (u',\vartheta) \in \Vminus{t}{u}
	\right\rbrace.
	\label{E:SIGMAMINUS}
\end{align}
\end{subequations}

\end{definition}

\begin{remark}[\textbf{Positive denominators}]
\label{R:POSITIVEDENOMINATORS}
The estimate \eqref{E:MUAMPLITUDENEARONE} implies that the denominator 
$\mathring{\upmu}(u',\vartheta) - M(0,u',\vartheta;t)$
in \eqref{E:ANGLESANDUWITHNONDECAYINUPMUGBEHAVIOR}-\eqref{E:ANGLESANDUWITHDECAYINUPMUGBEHAVIOR}
remains strictly positive all the way up to the shock in the
solution regime under consideration.
We include the denominator in the definitions
\eqref{E:ANGLESANDUWITHNONDECAYINUPMUGBEHAVIOR}-\eqref{E:ANGLESANDUWITHDECAYINUPMUGBEHAVIOR}
because it helps 
to clarify the connection
between the sets
$\Vplus{t}{u},
\Vminus{t}{u}$
and the parameter
$\LateTimeLUnitMu$
defined in \eqref{E:CRUCIALLATETIMEDERIVATIVEDEF}.
\end{remark}

\subsection{Sharp pointwise estimates for \texorpdfstring{$\upmu$}{the inverse foliation density} and its derivatives}
\label{SS:SHARPPOINTWISEESTIMATESFORUPMU}
In the next proposition, we derive 
sharp pointwise estimates for $\upmu$ and its derivatives.

\begin{proposition}[\textbf{Sharp pointwise estimates for $\upmu$, $\Lunit \upmu$, and $\Rad \upmu$}]
\label{P:SHARPMU} 
Under the data-size and bootstrap assumptions 
of Subsects.~\ref{SS:SIZEOFTBOOT}-\ref{SS:PSIBOOTSTRAP}
and the smallness assumptions of Subsect.~\ref{SS:SMALLNESSASSUMPTIONS}, 
the following  
estimates hold for $(t,u,\vartheta) \in [0,\Tboot) \times [0,U_0] \times \mathbb{T}$
and $0 \leq s \leq t$.
\medskip

\noindent \underline{\textbf{Upper bound for $\displaystyle \frac{[\Lunit \upmu]_+}{\upmu}$}.}
\begin{align} \label{E:POSITIVEPARTOFLMUOVERMUISBOUNDED}
	\left\|
		\frac{[\Lunit \upmu]_+}{\upmu}
	\right\|_{L^{\infty}(\Sigma_s^u)}
	& \leq C.
\end{align}

\medskip

\noindent \underline{\textbf{Small $\upmu$ implies $\Lunit \upmu$ is negative}.}
\begin{align} \label{E:SMALLMUIMPLIESLMUISNEGATIVE}
	\upmu(s,u,\vartheta) \leq \frac{1}{4}
	\implies
	\Lunit \upmu(s,u,\vartheta) \leq - \frac{1}{4} \TranminusdatasizeWithFactor.
\end{align}

\medskip

\noindent \underline{\textbf{Upper bound for $\displaystyle \frac{[\Rad \upmu]_+}{\upmu}$}.}
\begin{align} \label{E:UNIFORMBOUNDFORMRADMUOVERMU}
	\left\|
		\frac{[\Rad \upmu]_+}{\upmu}
	\right\|_{L^{\infty}(\Sigma_s^u)}
	& \leq 
		\frac{C}{\sqrt{\Tboot - s}}.
\end{align}

\medskip


\medskip

\noindent \underline{\textbf{Sharp spatially uniform estimates}.}
Consider a time interval $s \in [0,t]$ and define
the ($t,u-$dependent) constant $\LateTimeLUnitMu$ by
\begin{align} \label{E:CRUCIALLATETIMEDERIVATIVEDEF}
	\LateTimeLUnitMu 
	& := 
			- \min_{(u',\vartheta) \in [0,u] \times \mathbb{T}} 
			\frac{\Lunit \upmu(t,u',\vartheta)}{\mathring{\upmu}(u',\vartheta) - M(0,u',\vartheta;t)}.
\end{align}
Note that $\LateTimeLUnitMu \geq 0$ in view of the estimate \eqref{E:MUAMPLITUDENEARONE}
and the fact that $\Lunit \upmu$ vanishes along the flat null hyperplane $\mathcal{P}_0$.
Then
\begin{subequations}
\begin{align}
	\upmu_{\star}(s,u)
	& = \left\lbrace
				1 + \mathcal{O}(\varepsilon)
			\right\rbrace
			\left\lbrace 
				1 - \LateTimeLUnitMu s 
			\right\rbrace,
	\label{E:MUSTARBOUNDS}  
		\\
	\left\| 
		[\Lunit \upmu]_- 
	\right\|_{L^{\infty}(\Sigma_s^u)}
	& =
	\begin{cases}
		\left\lbrace
			1 + \mathcal{O}(\varepsilon^{1/2})
		\right\rbrace
		\LateTimeLUnitMu,
		& \mbox{if } \LateTimeLUnitMu \geq \sqrt{\varepsilon},
		\\
		\mathcal{O}(\varepsilon^{1/2}),
		& \mbox{if } \LateTimeLUnitMu \leq \sqrt{\varepsilon}.
	\end{cases}
	\label{E:LUNITUPMUMINUSBOUND}
\end{align}
\end{subequations}

We also have
\begin{subequations}
\begin{align} \label{E:UNOTNECESSARILYEQUALTOONECRUCIALLATETIMEDERIVATIVECOMPAREDTODATAPARAMETER}
	\LateTimeLUnitMu 
	& \leq \left\lbrace
				1 + \mathcal{O}(\varepsilon)
			\right\rbrace
			\TranminusdatasizeWithFactor.
\end{align}

Moreover, when $u = 1$, we have
\begin{align} \label{E:CRUCIALLATETIMEDERIVATIVECOMPAREDTODATAPARAMETER}
	\LateTimeLUnitMu 
	& = \left\lbrace
				1 + \mathcal{O}(\varepsilon)
			\right\rbrace
			\TranminusdatasizeWithFactor.
\end{align}
\end{subequations}

\medskip

\noindent \underline{\textbf{Sharp estimates when $(u',\vartheta) \in \Vplus{t}{u}$}.}
We recall that the set $\Vplus{t}{u}$ is defined in \eqref{E:ANGLESANDUWITHNONDECAYINUPMUGBEHAVIOR}.

If $0 \leq s_1 \leq s_2 \leq t$, then the following estimate holds:
\begin{align} \label{E:LOCALIZEDMUCANTGROWTOOFAST}
	\sup_{(u',\vartheta) \in \Vplus{t}{u}}
	\frac{\upmu(s_2,u',\vartheta)}{\upmu(s_1,u',\vartheta)}
	& \leq C.
\end{align}

In addition, if $s \in [0,t]$ and $\Sigmaplus{s}{t}{u}$ is as defined in \eqref{E:SIGMAPLUS}, then we have
\begin{align}  \label{E:KEYMUNOTDECAYBOUND}
		\inf_{\Sigmaplus{s}{t}{u}} \upmu 
	& \geq 1 - C \varepsilon.
\end{align}

In addition, if $s \in [0,t]$ and $\Sigmaplus{s}{t}{u}$ is as defined in \eqref{E:SIGMAPLUS}, then we have
\begin{align} 
	\left\| \frac{[\Lunit \upmu]_-}{\upmu} \right\|_{L^{\infty}(\Sigmaplus{s}{t}{u})}
	& \leq C \varepsilon.
		\label{E:KEYMUNOTDECAYINGMINUSPARTLMUOVERMUBOUND}
\end{align}

\medskip

\noindent \underline{\textbf{Sharp estimates when $(u',\vartheta) \in \Vminus{t}{u}$}.}
Assume that the set $\Vminus{t}{u}$ defined in \eqref{E:ANGLESANDUWITHDECAYINUPMUGBEHAVIOR} 
is non-empty, and consider a time interval $s \in [0,t]$. Let $\LateTimeLUnitMu > 0$
be as in \eqref{E:CRUCIALLATETIMEDERIVATIVEDEF}.
Then the following estimate holds:
\begin{align} \label{E:LOCALIZEDMUMUSTSHRINK}
	\mathop{\sup_{0 \leq s_1 \leq s_2 \leq t}}_{(u',\vartheta) \in \Vminus{t}{u}}
	\frac{\upmu(s_2,u',\vartheta)}{\upmu(s_1,u',\vartheta)}
	& \leq 1 + C \varepsilon.
\end{align}
Furthermore, if $s \in [0,t]$ and $\Sigmaminus{s}{t}{u}$
is as defined in \eqref{E:SIGMAMINUS}, then the following estimate holds:
\begin{align} \label{E:LMUPLUSNEGLIGIBLEINSIGMAMINUS}
	\left\| 
		[\Lunit \upmu]_+ 
	\right\|_{L^{\infty}(\Sigmaminus{s}{t}{u})}
	& \leq C \varepsilon.
\end{align}
Finally, there exists a constant $C > 0$ 
such that if $0 \leq s \leq t$, then
\begin{align}		\label{E:HYPERSURFACELARGETIMEHARDCASEOMEGAMINUSBOUND}
	\left\| 
		[\Lunit \upmu]_- 
	\right\|_{L^{\infty}(\Sigmaminus{s}{t}{u})}
	& \leq 
		\begin{cases}
		\left\lbrace
			1 + C \varepsilon^{1/2}
		\right\rbrace
		\LateTimeLUnitMu,
		& \mbox{if } \LateTimeLUnitMu \geq \sqrt{\varepsilon},
		\\
		C \varepsilon^{1/2},
		& \mbox{if } \LateTimeLUnitMu \leq \sqrt{\varepsilon}.
	\end{cases}
\end{align}

\noindent \underline{\textbf{Approximate time-monotonicity of $\upmu_{\star}^{-1}(s,u)$}.}
There exists a constant $C > 0$ such that if 
$0 \leq s_1 \leq s_2 \leq t$, then
\begin{align} \label{E:MUSTARINVERSEMUSTGROWUPTOACONSTANT}
	\upmu_{\star}^{-1}(s_1,u) & \leq (1 + C \varepsilon) \upmu_{\star}^{-1}(s_2,u).
\end{align}


\end{proposition}

\begin{proof}
See Subsect.~\ref{SS:OFTENUSEDESTIMATES} for some comments on the analysis.

\medskip
\noindent \textbf{Proof of \eqref{E:POSITIVEPARTOFLMUOVERMUISBOUNDED}:}
Clearly it suffices for us to prove that for
$(s,u,\vartheta) \in [0,t] \times[0,U_0] \times \mathbb{T}$, 
we have
$
\displaystyle
[\Lunit \upmu(s,u,\vartheta)]_+/\upmu(s,u,\vartheta) 
\leq C
$.
We may assume that $[\Lunit \upmu(s,u,\vartheta)]_+ > 0$ since otherwise the desired estimate is trivial.
Then by \eqref{E:LUNITUPMUATTIMETMINUSLUNITUPMUATTIMESPOINTWISEESTIMATE},
for $0 \leq s' \leq s \leq t < \Tboot \leq 2 \TranminusdatasizeWithFactor^{-1}$,
we have that 
$\Lunit \upmu(s',u,\vartheta) 
\geq 
\Lunit \upmu(s,u,\vartheta)
- C \varepsilon(s-s') 
\geq 
- C \varepsilon
$.
Integrating this estimate with respect to $s'$ starting from $s'=0$
and using \eqref{E:MUINITIALDATAESTIMATE},
we find that 
$\upmu(s,u,\vartheta) \geq 1 - C \varepsilon s \geq 1 - C \varepsilon$
and thus $1/\upmu(s,u,\vartheta) \leq 1 + C \varepsilon$.
Also using the bound 
$\left|
	\Lunit \upmu(s,u,\vartheta)
\right|
\leq C
$
(that is, \eqref{E:LUNITUPMULINFINITY}),
we conclude the desired estimate.

\medskip

\noindent \textbf{Proof of \eqref{E:SMALLMUIMPLIESLMUISNEGATIVE}:}
By \eqref{E:LUNITUPMUATTIMETMINUSLUNITUPMUATTIMESPOINTWISEESTIMATE},
for $0 \leq s \leq t < \Tboot \leq 2 \TranminusdatasizeWithFactor^{-1}$,
we have that 
$
\Lunit \upmu(s,u,\vartheta) = \Lunit \upmu(0,u,\vartheta) + \mathcal{O}(\varepsilon)
$.
Integrating this estimate with respect to $s$ starting from $s=0$
and using \eqref{E:MUINITIALDATAESTIMATE},
we find that 
$\upmu(s,u,\vartheta) 
	= 1 
	+ 
	\mathcal{O}(\varepsilon) 
	+
	s \Lunit \upmu(0,u,\vartheta) 
$.
It follows that 
whenever $\upmu(s,u,\vartheta) < 1/4$, we have
$
\Lunit \upmu(0,u,\vartheta) 
< 
- \frac{1}{2} \TranminusdatasizeWithFactor(3/4 + \mathcal{O}(\varepsilon))
= - \frac{3}{8} \TranminusdatasizeWithFactor + \mathcal{O}(\varepsilon)
$.
Again using \eqref{E:LUNITUPMUATTIMETMINUSLUNITUPMUATTIMESPOINTWISEESTIMATE}
to deduce that $\Lunit \upmu(s,u,\vartheta) = \Lunit \upmu(0,u,\vartheta) + \mathcal{O}(\varepsilon)$,
we arrive at the desired estimate \eqref{E:SMALLMUIMPLIESLMUISNEGATIVE}.

\medskip

\noindent \textbf{Proof of \eqref{E:UNOTNECESSARILYEQUALTOONECRUCIALLATETIMEDERIVATIVECOMPAREDTODATAPARAMETER} and \eqref{E:CRUCIALLATETIMEDERIVATIVECOMPAREDTODATAPARAMETER}:}
We prove only \eqref{E:CRUCIALLATETIMEDERIVATIVECOMPAREDTODATAPARAMETER}
since \eqref{E:UNOTNECESSARILYEQUALTOONECRUCIALLATETIMEDERIVATIVECOMPAREDTODATAPARAMETER}
follows from nearly identical arguments. 
Above we showed that
for $0 \leq t < \Tboot \leq 2 \TranminusdatasizeWithFactor^{-1}$,
we have 
$
\Lunit \upmu(t,u,\vartheta) = \Lunit \upmu(0,u,\vartheta) + \mathcal{O}(\varepsilon)
$.
Moreover, equation \eqref{E:UPMUFIRSTTRANSPORT}
and Lemma~\ref{L:SCHEMATICDEPENDENCEOFMANYTENSORFIELDS} imply that
$\Lunit \upmu = \frac{1}{2} G_{\Lunit \Lunit} \Rad \Psi + \smoothfunction(\BadVar) \Singletan \Psi$.
From this relation and the $L^{\infty}$ estimates of Prop.~\ref{P:IMPROVEMENTOFAUX},
we deduce that
$\Lunit \upmu(0,u,\vartheta) = \frac{1}{2} [G_{\Lunit \Lunit} \Rad \Psi](0,u,\vartheta) + \mathcal{O}(\varepsilon)$.
In addition, from 
\eqref{E:MUAMPLITUDENEARONE}, we deduce that
$\mathring{\upmu}(u,\vartheta) - M(0,u,\vartheta;t) = 1 + \mathcal{O}(\varepsilon)$.
Combining these estimates and appealing to definitions 
\eqref{E:CRITICALBLOWUPTIMEFACTOR}
and
\eqref{E:CRUCIALLATETIMEDERIVATIVEDEF},
we conclude \eqref{E:CRUCIALLATETIMEDERIVATIVECOMPAREDTODATAPARAMETER}.

\medskip
\noindent \textbf{Proof of \eqref{E:MUSTARBOUNDS} and \eqref{E:MUSTARINVERSEMUSTGROWUPTOACONSTANT}:}
We first prove \eqref{E:MUSTARBOUNDS}.
We start by establishing the following preliminary estimate
for the crucial quantity $\LateTimeLUnitMu = \LateTimeLUnitMu(t,u)$ 
(see \eqref{E:CRUCIALLATETIMEDERIVATIVEDEF}):
\begin{align} \label{E:LATETIMELMUTIMESTISLESSTHANONE}
	t \LateTimeLUnitMu
	< 1.
\end{align} 
Using
\eqref{E:MUAPPROXDEF},
\eqref{E:MUSPLIT},
\eqref{E:MUAMPLITUDENEARONE},
and
\eqref{E:BIGMEST},
we deduce that the following estimate holds
for $(s,u',\vartheta) \in [0,t] \times [0,u] \times \mathbb{T}$:
\begin{align} \label{E:MUFIRSTLOWERBOUND}
	\upmu(s,u',\vartheta)
	& =
		(1 + \mathcal{O}(\varepsilon))
		\left\lbrace
			1 
			+
			\frac{\Lunit \upmu(t,u',\vartheta)}{\mathring{\upmu}(u',\vartheta) - M(0,u',\vartheta;t)} s
			+ 
			\mathcal{O}(\varepsilon) (t-s)^2
		\right\rbrace.
\end{align}
Setting $s=t$ in equation \eqref{E:MUFIRSTLOWERBOUND},
taking the min of both sides over $(u',\vartheta) \in [0,u] \times \mathbb{T}$,
and appealing to definitions
\eqref{E:MUSTARDEF} and 
\eqref{E:CRUCIALLATETIMEDERIVATIVEDEF},
we deduce that
$\upmu_{\star}(t,u)
= (1 + \mathcal{O}(\varepsilon))(1-\LateTimeLUnitMu t)
$.
Since $\upmu_{\star}(t,u) > 0$ by \eqref{E:BOOTSTRAPMUPOSITIVITY},
we conclude \eqref{E:LATETIMELMUTIMESTISLESSTHANONE}.

Having established the preliminary estimate, 
we now 
take the min of both sides over $(u',\vartheta) \in [0,u] \times \mathbb{T}$,
and appeal to definitions
\eqref{E:MUSTARDEF} and 
\eqref{E:CRUCIALLATETIMEDERIVATIVEDEF}
%
to obtain:
\begin{align} \label{E:HARDERCASEMUFIRSTLOWERBOUND}
	\min_{(u',\vartheta) \in [0,u] \times \mathbb{T}} \upmu(s,u',\vartheta)
	& = 
		(1 + \mathcal{O}(\varepsilon))
		\left\lbrace
			1 
			- \LateTimeLUnitMu s
			+ \mathcal{O}(\varepsilon) (t-s)^2
		\right\rbrace.
\end{align}
We will show that the terms in braces on RHS~\eqref{E:HARDERCASEMUFIRSTLOWERBOUND} verify
\begin{align} \label{E:MUSECONDLOWERBOUND}
	1 
	- \LateTimeLUnitMu s
	+ \mathcal{O}(\varepsilon) (t-s)^2
	& 
	=
	(1 + \smoothfunction(s,u;t))
	\left\lbrace
		1 - \LateTimeLUnitMu s
	\right\rbrace,
\end{align}
where
\begin{align} \label{E:AMPLITUDEDEVIATIONFUNCTIONMUSECONDLOWERBOUND}
	\smoothfunction(s,u;t)
	& = \mathcal{O}(\varepsilon).
\end{align}
The desired estimate \eqref{E:MUSTARBOUNDS}
then follows easily from 
\eqref{E:HARDERCASEMUFIRSTLOWERBOUND}-\eqref{E:AMPLITUDEDEVIATIONFUNCTIONMUSECONDLOWERBOUND}
and 
definition \eqref{E:MUSTARDEF}.
To prove \eqref{E:AMPLITUDEDEVIATIONFUNCTIONMUSECONDLOWERBOUND}, 
we first use \eqref{E:MUSECONDLOWERBOUND} to solve for $\smoothfunction(s,u;t)$: 
\begin{align} \label{E:AMPLITUDEDEVIATIONFUNCTIONEXPRESSION}
	\smoothfunction(s,u;t)
	=
	\frac{\mathcal{O}(\varepsilon) (t-s)^2
				}
				{
				1 - \LateTimeLUnitMu s
				}
	=
	\frac{\mathcal{O}(\varepsilon) (t-s)^2
				}
				{
				1 - \LateTimeLUnitMu t
				+ 
				\LateTimeLUnitMu (t-s)
				}.
\end{align}
We start by considering the case $\LateTimeLUnitMu \leq (1/4) \TranminusdatasizeWithFactor$.
Since $0 \leq s \leq t < \Tboot \leq 2 \TranminusdatasizeWithFactor^{-1}$,
the denominator in the middle expression in \eqref{E:AMPLITUDEDEVIATIONFUNCTIONEXPRESSION}
is $\geq 1/2$, and the desired estimate
\eqref{E:AMPLITUDEDEVIATIONFUNCTIONMUSECONDLOWERBOUND}
follows easily whenever
$\varepsilon$ is sufficiently small.
In remaining case, we have
$\LateTimeLUnitMu > (1/4) \TranminusdatasizeWithFactor$.
Using \eqref{E:LATETIMELMUTIMESTISLESSTHANONE},
we deduce that 
RHS~\eqref{E:AMPLITUDEDEVIATIONFUNCTIONEXPRESSION}
$
\leq \frac{1}{\LateTimeLUnitMu} \mathcal{O}(\varepsilon) (t-s)
\leq C \varepsilon \TranminusdatasizeWithFactor^{-2} 
\lesssim \varepsilon
$
as desired.

Inequality \eqref{E:MUSTARINVERSEMUSTGROWUPTOACONSTANT} then follows as a simple consequence of \eqref{E:MUSTARBOUNDS}.

\medskip

\noindent \textbf{Proof of \eqref{E:LUNITUPMUMINUSBOUND} and \eqref{E:HYPERSURFACELARGETIMEHARDCASEOMEGAMINUSBOUND}:}
To prove \eqref{E:LUNITUPMUMINUSBOUND},
we first use \eqref{E:LUNITUPMUATTIMETMINUSLUNITUPMUATTIMESPOINTWISEESTIMATE}
to deduce that for $0 \leq s \leq t < \Tboot \leq 2 \TranminusdatasizeWithFactor^{-1}$ 
and $(u',\vartheta) \in [0,u] \times \mathbb{T}$,
we have $\Lunit \upmu(s,u',\vartheta) = \Lunit \upmu(t,u',\vartheta) + \mathcal{O}(\varepsilon)$.
Appealing to definition \eqref{E:CRUCIALLATETIMEDERIVATIVEDEF} and
using the estimate \eqref{E:MUAMPLITUDENEARONE},
we find that
$
\displaystyle
\left\| 
	[\Lunit \upmu]_- 
\right\|_{L^{\infty}(\Sigma_s^u)}
=
	\LateTimeLUnitMu
	+ 
	\mathcal{O}(\varepsilon)
$.
If
$
\displaystyle
\sqrt{\varepsilon}
\leq \LateTimeLUnitMu
$,
we see that as long as $\varepsilon$ is sufficiently small, 
we have the desired bound
$
\displaystyle
	\LateTimeLUnitMu
	+ 
	\mathcal{O}(\varepsilon)
	=
	(1 + \mathcal{O}(\varepsilon^{1/2}))
	\LateTimeLUnitMu
$.
On the other hand, if
$
\displaystyle
\LateTimeLUnitMu
\leq 
\sqrt{\varepsilon}
$,
then similar reasoning yields that
$
\displaystyle
\left\| 
	[\Lunit \upmu]_- 
\right\|_{L^{\infty}(\Sigma_s^u)}
=
	\LateTimeLUnitMu
	+ 
	\mathcal{O}(\varepsilon)
= \mathcal{O}(\sqrt{\varepsilon})
$
as desired. We have thus proved \eqref{E:LUNITUPMUMINUSBOUND}.

The proof of \eqref{E:HYPERSURFACELARGETIMEHARDCASEOMEGAMINUSBOUND} is similar
and we omit the details.

\medskip

\noindent \textbf{Proof of \eqref{E:UNIFORMBOUNDFORMRADMUOVERMU}:}
We fix times $s$ and $t$ with $0 \leq s \leq t < \Tboot$
and a point $p \in \Sigma_s^u$ with geometric coordinates 
$(s,\widetilde{u},\widetilde{\vartheta})$.
Let $\iota : [0,u] \rightarrow \Sigma_s^u$ be the integral curve
of $\Rad$ that passes through $p$ and that is parametrized by the values $u'$ of the eikonal function.
We set $F(u') := \upmu \circ \iota(u')$ and 
$
\displaystyle
\dot{F}(u') := \frac{d}{d u'} F(u') = (\Rad \upmu)\circ \iota(u')
$.
We must bound 
$
\displaystyle
\frac{[\Rad \upmu]_+}{\upmu}|_p
= \frac{[\dot{F}(\widetilde{u})]_+}{F(\widetilde{u})}
$.
We split the proof into three cases that exhaust all possibilities.
In the first case, we assume that $\dot{F}(\widetilde{u}) = \Rad \upmu|_p \leq 0$. Then 
$[\dot{F}(\widetilde{u})]_+/F(\widetilde{u}) = 0 \leq \mbox{RHS~\eqref{E:UNIFORMBOUNDFORMRADMUOVERMU}}$ as desired.
In the second case, we assume that $\dot{F}(u') \geq 0$ for $u' \in [0,\widetilde{u}]$. 
Then since $F(0) = 1$ (because the solution is trivial in the exterior of the flat null hyperplane $\mathcal{P}_0$),
we have that $F(u') \geq F(0)= 1$
for $u' \in [0,\widetilde{u}]$. Also using the bounds 
$\| \upmu\|_{L^{\infty}(\Sigma_s^u)}, \| \Rad \upmu \|_{L^{\infty}(\Sigma_s^u)} \leq C$ 
(that is, \eqref{E:UPMULINFTY} and \eqref{E:RADUPMULINFTY}),
we deduce that
$
	[\dot{F}(\widetilde{u})]_+/F(\widetilde{u}) \leq C \leq \mbox{RHS~\eqref{E:UNIFORMBOUNDFORMRADMUOVERMU}}
$
as desired. In the final case, we have
$\dot{F}(\widetilde{u}) = \Rad \upmu|_p > 0$ and there exists a largest number $u_* \in (0,\widetilde{u})$
such that $\dot{F}(u_*) = 0$ (and hence $\dot{F}(u') > 0$ for $u' \in (u_*,\widetilde{u}]$).
We will use the following estimate
for $\upmu_{(Min)}(s,u') 
:= \min_{(u'',\vartheta) \in [0,u'] \times \mathbb{T}} \upmu(s,u'',\vartheta)$, 
which holds for all $0 \leq s \leq t < \Tboot$ and $u' \in [0,u]$:
\begin{align} \label{E:UPMUMINLOWERBOUND}
	\upmu_{(Min)}(s,u')
	& \geq 
		\max 
		\left\lbrace
			(1 - C \varepsilon)
			\LateTimeLUnitMu 
			(t-s),
			(1 - C \varepsilon) 
			(1 - \LateTimeLUnitMu s)
		\right\rbrace,
\end{align}
where $\LateTimeLUnitMu = \LateTimeLUnitMu(t,u)$ is defined in
\eqref{E:CRUCIALLATETIMEDERIVATIVEDEF}.
We prove \eqref{E:UPMUMINLOWERBOUND} below in the last paragraph of the proof.

To proceed, we set $H := \sup_{\mathcal{M}_{\Tboot,U_0}} \Rad \Rad \upmu$.
By the mean value theorem, we have
$H > 0$ (since $\dot{F}(\widetilde{u}) > 0$ and $\dot{F}(u_*) = 0$).
Moreover, by \eqref{E:RADRADUPMULINFTY}, we have $H \leq C$.
In the next paragraph, 
we will use the mean value theorem to prove that \emph{at the point $p$ of interest},
we have
\begin{align} \label{E:FDIFFERENCELOWERBOUND}
	\upmu 
	- 
	\upmu_{(Min)}
	& \geq 
	\frac{1}{4} \frac{[\Rad \upmu]_+^2}{H}.
\end{align}
Rearranging \eqref{E:FDIFFERENCELOWERBOUND}, we find that
$[\Rad \upmu]_+ \leq 2 H^{1/2} \sqrt{\upmu - \upmu_{(Min)}}$
and thus the following bound holds at $p$:
\begin{align} \label{E:RADMUOVERMUALGEBRAICBOUND}
	\frac{[\Rad \upmu]_+}{\upmu}
	& \leq 
		2 H^{1/2} \frac{\sqrt{\upmu - \upmu_{(Min)}}}{\upmu}.
\end{align}
We now view RHS~\eqref{E:RADMUOVERMUALGEBRAICBOUND} as a function
of the real variable $\upmu$ (with all other parameters fixed)
on the domain $[\upmu_{(Min)},\infty)$. A simple calculus exercise 
yields that RHS~\eqref{E:RADMUOVERMUALGEBRAICBOUND}
$\leq H^{1/2}/\sqrt{\upmu_{(Min)}}$. 
Combining this estimate with \eqref{E:UPMUMINLOWERBOUND}
and using the aforementioned bound $H \leq C$,
we deduce that
$
\displaystyle
\frac{[\Rad \upmu]_+}{\upmu}|_p
\leq C 
	\min 
	\left\lbrace
		1/(\LateTimeLUnitMu^{1/2}(t-s)^{1/2}),
		1/(1 - \LateTimeLUnitMu s)^{1/2}
	\right\rbrace
$.
If $\LateTimeLUnitMu \leq (1/4) \TranminusdatasizeWithFactor$,
then $1 - \LateTimeLUnitMu s \geq 1 - (1/4) \TranminusdatasizeWithFactor \Tboot \geq 1/2$,
and the
desired bound 
$
\displaystyle
\frac{[\Rad \upmu]_+}{\upmu}|_p
\leq C 
\leq C/\Tboot^{1/2}
\leq C/(\Tboot-s)^{1/2}
\leq \mbox{RHS~\eqref{E:UNIFORMBOUNDFORMRADMUOVERMU}}
$
follows easily from the second term in the $\min$.
If $\LateTimeLUnitMu \geq (1/4) \TranminusdatasizeWithFactor$, 
then $1/\LateTimeLUnitMu \leq C$, 
and using the first term in the $\min$, we
deduce that 
$
\displaystyle
\frac{[\Rad \upmu]_+}{\upmu}|_p
\leq C/(t-s)^{1/2}
$.
Since this estimate holds for all $t < \Tboot$ with a uniform constant $C$,
we conclude \eqref{E:UNIFORMBOUNDFORMRADMUOVERMU} in this case.

To prove the bound \eqref{E:FDIFFERENCELOWERBOUND} used above, 
we set $u_1 := \widetilde{u} - (1/2)\dot{F}(\widetilde{u})/H$ and use the mean value theorem to deduce that
for $u' \in [u_1, \widetilde{u}]$, we have
$\dot{F}(\widetilde{u}) - \dot{F}(u') 
	\leq H (\widetilde{u} - u')
	\leq (1/2) \dot{F}(\widetilde{u})
$.
Thus, we have
\begin{align} \label{E:FPRIMELOWERBOUND}
	\dot{F}(u') 
	& \geq \frac{1}{2} \dot{F}(\widetilde{u}),
	& \mbox{for } u' 
	&\in [u_1,\widetilde{u}].
\end{align}
Again using the mean value theorem and also \eqref{E:FPRIMELOWERBOUND},
we deduce that
$
F(\widetilde{u}) 
- 
F(u_1)
\geq 
(1/2) \dot{F}(\widetilde{u}) (\widetilde{u}-u_1)
= (1/4) \dot{F}^2(\widetilde{u})/H
$. Noting that 
the definition of $\upmu_{(Min)}$ implies that
$F(u_1) \geq \upmu_{(Min)}(s,\widetilde{u})$,
we conclude the desired estimate \eqref{E:FDIFFERENCELOWERBOUND}.

It remains for us to prove \eqref{E:UPMUMINLOWERBOUND}. 
Reasoning as in the proof of \eqref{E:MUFIRSTLOWERBOUND}-\eqref{E:AMPLITUDEDEVIATIONFUNCTIONMUSECONDLOWERBOUND}
and using \eqref{E:LATETIMELMUTIMESTISLESSTHANONE},
we find that for $0 \leq s \leq t < \Tboot$ and $u' \in [0,u]$,
we have
$\upmu_{(Min)}(s,u') 
\geq 
(1 - C \varepsilon)
\left\lbrace
	1 - \LateTimeLUnitMu s
\right\rbrace
\geq
(1 - C \varepsilon)
\LateTimeLUnitMu (t-s)
$.
From these two inequalities, 
we conclude the desired bound \eqref{E:UPMUMINLOWERBOUND}.

%

\medskip

\noindent {\textbf{Proof of} \eqref{E:LOCALIZEDMUMUSTSHRINK}:}
A straightforward modification of the proof of \eqref{E:MUSTARBOUNDS},
based on equation 
\eqref{E:MUFIRSTLOWERBOUND}
and on replacing $\LateTimeLUnitMu$ in 
\eqref{E:HARDERCASEMUFIRSTLOWERBOUND}-\eqref{E:MUSECONDLOWERBOUND} 
with $\Lunit \upmu(t,u',\vartheta)$
(without taking the $\min$ on the LHS of the analog of \eqref{E:HARDERCASEMUFIRSTLOWERBOUND}),
yields that for $0 \leq s \leq t < \Tboot$ 
and $(u',\vartheta) \in \Vminus{t}{u}$,
we have
$
\displaystyle
\upmu(s,u',\vartheta) 
= \left\lbrace
				1 + \mathcal{O}(\varepsilon)
			\right\rbrace
			\left\lbrace 
				1 - \left|\Lunit \upmu(t,u',\vartheta)\right| s 
			\right\rbrace$.
The estimate \eqref{E:LOCALIZEDMUMUSTSHRINK} 
then follows as a simple consequence.

\medskip

\noindent {\textbf{Proof of}
\eqref{E:LOCALIZEDMUCANTGROWTOOFAST}, \eqref{E:KEYMUNOTDECAYBOUND}, \textbf{and} \eqref{E:KEYMUNOTDECAYINGMINUSPARTLMUOVERMUBOUND}:}
By \eqref{E:LUNITUPMUATTIMETMINUSLUNITUPMUATTIMESPOINTWISEESTIMATE},
if $(u',\vartheta) \in \Vplus{t}{u}$ and
$0 \leq s \leq t < \Tboot$,
then
$[\Lunit \upmu]_-(s,u,\vartheta) \leq C \varepsilon$
and
$\Lunit \upmu(s,u,\vartheta) \geq - C \varepsilon$.
Integrating the latter estimate with respect to $s$ from $0$ to $t$
and using \eqref{E:MUINITIALDATAESTIMATE},
we find that if $0 \leq s < \Tboot$ 
and $(u',\vartheta) \in \Vplus{t}{u}$,
then
$\upmu(s,u',\vartheta) 
\geq 1 - C \varepsilon(1 + s)
\geq 1 - C \varepsilon
$.
Moreover, from \eqref{E:UPMULINFTY}, 
we have the crude bound $\upmu(s,u',\vartheta) \leq C$.
The desired bounds
\eqref{E:LOCALIZEDMUCANTGROWTOOFAST}, \eqref{E:KEYMUNOTDECAYBOUND}, and \eqref{E:KEYMUNOTDECAYINGMINUSPARTLMUOVERMUBOUND}
now readily follow
from these estimates.

\medskip

\noindent {\textbf{Proof of} \eqref{E:LMUPLUSNEGLIGIBLEINSIGMAMINUS}:}
By \eqref{E:LUNITUPMUATTIMETMINUSLUNITUPMUATTIMESPOINTWISEESTIMATE},
if $(u',\vartheta) \in \Vminus{t}{u}$ and
$0 \leq s \leq t < \Tboot$,
then
$[\Lunit \upmu]_+(s,u',\vartheta)
=  [\Lunit \upmu]_+(t,u',\vartheta) + \mathcal{O}(\varepsilon) = \mathcal{O}(\varepsilon)$.
The desired bound \eqref{E:LMUPLUSNEGLIGIBLEINSIGMAMINUS} thus follows.

\end{proof}

\subsection{Sharp time-integral estimates involving \texorpdfstring{$\upmu$}{the inverse foliation density}}
\label{SS:SHARPTIMEINTEGRALESTIMATES}
In Prop.~\ref{P:MUINVERSEINTEGRALESTIMATES}, 
we use the sharp pointwise estimates of Prop.~\ref{P:SHARPMU} 
to derive sharp estimates for time integrals involving powers of $\upmu_{\star}^{-1}$.
The time-integral estimates are a primary ingredient 
in the Gronwall-type argument that we use
to derive a priori energy estimates (see Prop.~\ref{P:MAINAPRIORIENERGY}),
which are degenerate with respect to powers of $\upmu_{\star}^{-1}$ at the high orders
(see inequality \eqref{E:MULOSSMAINAPRIORIENERGYESTIMATES}).
The estimates of Prop.~\ref{P:MUINVERSEINTEGRALESTIMATES}
directly influence the degree of $\upmu_{\star}^{-1}-$degeneracy
found in our high-order energy estimates.

\begin{proposition}[\textbf{Fundamental estimates for time integrals involving $\upmu^{-1}$}] 
\label{P:MUINVERSEINTEGRALESTIMATES}
	Let $\upmu_{\star}(t,u)$ be as defined in \eqref{E:MUSTARDEF}.
	Let
	\begin{align*}
		\Contwo > 1
	\end{align*}
	be a real number. 
	Under the data-size and bootstrap assumptions 
	of Subsects.~\ref{SS:SIZEOFTBOOT}-\ref{SS:PSIBOOTSTRAP}
	and the smallness assumptions of Subsect.~\ref{SS:SMALLNESSASSUMPTIONS}, 
	the following  
	estimates hold for $(t,u) \in [0,\Tboot) \times [0,U_0]$:

	\medskip

	\noindent \underline{\textbf{Estimates relevant for borderline top-order spacetime integrals}.}
	There exists a constant $C > 0$ 
	such that if $\Contwo \sqrt{\varepsilon} \leq 1$,
	then
	\begin{align} \label{E:KEYMUTOAPOWERINTEGRALBOUND}
		\int_{s=0}^t 
			\frac{\left\| [\Lunit \upmu]_- \right\|_{L^{\infty}(\Sigma_s^u)}} 
					 {\upmu_{\star}^{\Contwo}(s,u)}
		\, ds 
		& \leq \frac{1 + C \sqrt{\varepsilon}}{\Contwo-1} 
			\upmu_{\star}^{1-\Contwo}(t,u).
	\end{align}

	\noindent \underline{\textbf{Estimates relevant for borderline top-order hypersurface integrals}.}
	There exists a constant $C > 0$ such that
	\begin{align} \label{E:KEYHYPERSURFACEMUTOAPOWERINTEGRALBOUND}
		\left\| 
			\Lunit \upmu 
		\right\|_{L^{\infty}(\Sigmaminus{t}{t}{u})} 
		\int_{s=0}^t 
			\frac{1} 
				{\upmu_{\star}^{\Contwo}(s,u)}
			\, ds 
		& \leq \frac{1 + C \sqrt{\varepsilon}}{\Contwo-1} 
			\upmu_{\star}^{1-\Contwo}(t,u).
	\end{align}

	\medskip

	\noindent \underline{\textbf{Estimates relevant for less dangerous top-order spacetime integrals}.}
	There exists a constant $C > 0$ 
	such that if $\Contwo \sqrt{\varepsilon} \leq 1$,
	then
	\begin{align} \label{E:LOSSKEYMUINTEGRALBOUND}
		\int_{s=0}^t \frac{1} 
			{\upmu_{\star}^{\Contwo}(s,u)}
		\, ds 
		& \leq C \left\lbrace 2^{\Contwo} + \frac{1}{\Contwo-1} \right\rbrace \upmu_{\star}^{1-\Contwo}(t,u).
	\end{align} 

	\medskip

	\noindent \underline{\textbf{Estimates for integrals that lead to only $\ln \upmu_{\star}^{-1}$ degeneracy}.} 
	There exists a constant $C > 0$ such that
	\begin{align} \label{E:KEYMUINVERSEINTEGRALBOUND}
		\int_{s=0}^t 
			\frac{\left\| [\Lunit \upmu]_- \right\|_{L^{\infty}(\Sigma_s^u)}} 
					 {\upmu_{\star}(s,u)}
		\, ds 
		& \leq (1 + C \sqrt{\varepsilon}) \ln \upmu_{\star}^{-1}(t,u) + C \sqrt{\varepsilon}.
	\end{align}
	In addition, there exists a constant $C > 0$ such that
	\begin{align} \label{E:LOGLOSSMUINVERSEINTEGRALBOUND}
		\int_{s=0}^t 
			\frac{1}{\upmu_{\star}(s,u)}
		\, ds 
		& \leq  C \left\lbrace \ln \upmu_{\star}^{-1}(t,u) + 1 \right\rbrace.
	\end{align}

	\medskip

	\noindent \underline{\textbf{Estimates for integrals that break the $\upmu_{\star}^{-1}$ degeneracy}.} 
	There exists a constant $C > 0$ such that
	\begin{align} \label{E:LESSSINGULARTERMSMPOINTNINEINTEGRALBOUND}
		\int_{s=0}^t 
			\frac{1} 
			{\upmu_{\star}^{9/10}(s,u)}
		\, ds 
		& \leq C.
	\end{align}

\end{proposition}

\begin{proof}
\noindent {\textbf{Proof of} \eqref{E:KEYMUTOAPOWERINTEGRALBOUND}, \eqref{E:KEYHYPERSURFACEMUTOAPOWERINTEGRALBOUND},
and \eqref{E:KEYMUINVERSEINTEGRALBOUND}:}
To prove \eqref{E:KEYMUTOAPOWERINTEGRALBOUND}, we 
first consider the case $\LateTimeLUnitMu \geq \sqrt{\varepsilon}$
in \eqref{E:LUNITUPMUMINUSBOUND}.
Using \eqref{E:MUSTARBOUNDS} and \eqref{E:LUNITUPMUMINUSBOUND}, we deduce that
\begin{align} \label{E:PROOFKEYMUTOAPOWERINTEGRALBOUND}
		\int_{s=0}^t 
			\frac{\left\| [\Lunit \upmu]_- \right\|_{L^{\infty}(\Sigma_s^u)}} 
					 {\upmu_{\star}^{\Contwo}(s,u)}
		\, ds 
		& = (1 + \mathcal{O}(\varepsilon^{1/2}))
			\int_{s=0}^t 
				\frac{\LateTimeLUnitMu}{(1 - \LateTimeLUnitMu s)^{\Contwo}}
			\, ds 
				\\
		& \leq \frac{1 + \mathcal{O}(\varepsilon^{1/2})}{\Contwo - 1}
				\frac{1}{(1 - \LateTimeLUnitMu t)^{\Contwo-1}}
			= 	\frac{1 + \mathcal{O}(\varepsilon^{1/2})}{\Contwo - 1}
					\upmu_{\star}^{1-\Contwo}(t,u)
					\notag
	\end{align}
	as desired. We now consider the case
	$\LateTimeLUnitMu \leq \sqrt{\varepsilon}$
	in \eqref{E:LUNITUPMUMINUSBOUND}.
	Using 
	\eqref{E:MUSTARBOUNDS}
	and
	\eqref{E:LUNITUPMUMINUSBOUND}
	and
	the fact that $0 \leq s \leq t < \Tboot \leq 2 \TranminusdatasizeWithFactor^{-1}$,
	we see that for $\varepsilon$ sufficiently small relative to 
	$\TranminusdatasizeWithFactor$, we have
	\begin{align} \label{E:SECONDCASEPROOFKEYMUTOAPOWERINTEGRALBOUND}
		\int_{s=0}^t 
			\frac{\left\| [\Lunit \upmu]_- \right\|_{L^{\infty}(\Sigma_s^u)}} 
					 {\upmu_{\star}^{\Contwo}(s,u)}
		\, ds 
		& \leq
			C \varepsilon^{1/2}
			\int_{s=0}^t 
				\frac{1}{(1 - \LateTimeLUnitMu s)^{\Contwo}}
			\, ds 
				\\
		& \leq
			C \varepsilon^{1/2}
			\frac{1}{(1 - \LateTimeLUnitMu t)^{\Contwo-1}}
			\leq \frac{1}{\Contwo - 1}
			\upmu_{\star}^{1-\Contwo}(t,u)
					\notag
	\end{align}
	as desired.
	We have thus proved \eqref{E:KEYMUTOAPOWERINTEGRALBOUND}.
	The estimate \eqref{E:KEYMUINVERSEINTEGRALBOUND}
	can be proved in a similar fashion; we omit the details.

	Inequality \eqref{E:KEYHYPERSURFACEMUTOAPOWERINTEGRALBOUND} 
	can be proved in a similar fashion
	with the help of the estimate \eqref{E:HYPERSURFACELARGETIMEHARDCASEOMEGAMINUSBOUND};
	we omit the details.

	\medskip

\noindent {\textbf{Proof of} \eqref{E:LOSSKEYMUINTEGRALBOUND},
\eqref{E:LOGLOSSMUINVERSEINTEGRALBOUND}, and \eqref{E:LESSSINGULARTERMSMPOINTNINEINTEGRALBOUND}:}
	To prove \eqref{E:LOSSKEYMUINTEGRALBOUND}, 
	we first use \eqref{E:MUSTARBOUNDS} to deduce
	\begin{align} \label{E:PROOFLOSSKEYMUINTEGRALBOUND}
		\int_{s=0}^t \frac{1} 
			{\upmu_{\star}^{\Contwo}(s,u)}
		\, ds 
		& \leq 
		C
		\int_{s=0}^t 
			\frac{1}{(1 - \LateTimeLUnitMu s)^{\Contwo}}
		\, ds.
	\end{align}
	We first assume that $\LateTimeLUnitMu \leq (1/4) \TranminusdatasizeWithFactor$. 
	Then since $0 \leq t < \Tboot < 2 \TranminusdatasizeWithFactor^{-1}$,
	we see from \eqref{E:MUSTARBOUNDS} that $\upmu_{\star}(s,u) \geq (1/4)$ for $0 \leq s \leq t$
	and that RHS~\eqref{E:PROOFLOSSKEYMUINTEGRALBOUND} 
	$
	\leq C 2^{\Contwo} t 
	\leq C 2^{\Contwo} \TranminusdatasizeWithFactor^{-1} 
	\leq C 2^{\Contwo} 
	\leq C 2^{\Contwo} \upmu_{\star}^{1-\Contwo}(t,u)$
	as desired.
	In the remaining case, we have $\LateTimeLUnitMu > (1/4) \TranminusdatasizeWithFactor$,
	and we can use \eqref{E:MUSTARBOUNDS} 
	and the estimate $1/\LateTimeLUnitMu \leq C$
	to bound RHS~\eqref{E:PROOFLOSSKEYMUINTEGRALBOUND} by
	\begin{align} \label{E:SECONDPROORELATIONFLOSSKEYMUINTEGRALBOUND}
		&
		\leq
		\frac{C}{\LateTimeLUnitMu}
		\frac{1}{\Contwo - 1}
		\frac{1}{(1 - \LateTimeLUnitMu t)^{\Contwo-1}}
		\leq 
		\frac{C}{\Contwo - 1}
		\upmu_{\star}^{1-\Contwo}(t,u)
	\end{align}
	as desired.

	Inequalities 
	\eqref{E:LOGLOSSMUINVERSEINTEGRALBOUND}
	and
	\eqref{E:LESSSINGULARTERMSMPOINTNINEINTEGRALBOUND} can be proved in a similar
	fashion; we omit the details, aside from remarking that the last step of the proof of 
	\eqref{E:LESSSINGULARTERMSMPOINTNINEINTEGRALBOUND}
	relies on the trivial estimate $(1 - \LateTimeLUnitMu t)^{1/10} \leq 1$.

\end{proof}

\section{Pointwise estimates for the error integrands}
\label{S:POINTWISEESTIMATESFORWAVEEQUATIONERRORTERMS}
\setcounter{equation}{0}
Recall that if $\Tanset^N$ is an $N^{th}-$order $\mathcal{P}_u-$tangent vectorfield operator, 
then $\Tanset^N \Psi$ solves an inhomogeneous wave equation of the form
$\upmu \square_g (\Tanset^N \Psi) = \mathfrak{F}$.
In this section, 
we start by identifying
the difficult error terms in $\mathfrak{F}$; see 
Subsect.~\ref{SS:IDENTIFICATIONOFKEYDIFFICULTERRORTERMFACTORS}.
The difficult terms are products that contain a factor involving certain top derivatives
of the eikonal function; we have to work hard to 
control these products in the energy estimates. 
The remaining terms have a structure that we call
``$Harmless^{\leq N}$'' 
(see Def.~\ref{D:HARMLESSTERMS})
and are easy to control.
Next, we derive pointwise estimates
for the difficult products. After deriving some preliminary estimates, we 
provide the main result in this direction in Prop.~\ref{P:KEYPOINTWISEESTIMATE}.
Finally, in Subsect.~\ref{SS:POINTWISEESTIMATESFORTHEMULTIPLIERVECTORFIELDERRORTERMS},
we derive pointwise estimates for the 
error terms $\sum_{i=1}^5 \basicenergyerrorarg{\Mult}{i}$ on RHS~\eqref{E:MULTERRORINT},
which are generated by the deformation tensor of the multiplier vectorfield $\Mult$.

The following definition encapsulates error term factors
that are easy to bound in the energy estimates.
Most factors that arise in our analysis are of this form.

\begin{definition}[\textbf{Harmless terms}]
	\label{D:HARMLESSTERMS}
	A $Harmless^{\leq N}$ term is any term such that 
	under the data-size and bootstrap assumptions 
	of Subsects.~\ref{SS:SIZEOFTBOOT}-\ref{SS:PSIBOOTSTRAP}
	and the smallness assumptions of Subsect.~\ref{SS:SMALLNESSASSUMPTIONS}, 
	the following bound holds on
	$\mathcal{M}_{\Tboot,U_0}$
	(see Subsect.~\ref{SS:STRINGSOFCOMMUTATIONVECTORFIELDS} regarding the vectorfield operator notation):
	\begin{align} \label{E:HARMESSTERMPOINTWISEESTIMATE}
		\left| 
			Harmless^{\leq N}
		\right|
		& \lesssim 
			\left|
				\Fullset_*^{\leq N+1;1} \Psi
			\right|
			+
			\left|
				\Fullset_*^{\leq N;1} \GdVar
			\right|
			+
			\left|
				\Tanset_*^{[1,N]} \BadVar
			\right|.
	\end{align}

\end{definition}

In the next lemma, we provide 
$L^{\infty}$ estimates for $Harmless^{\leq 9}$ terms.

\begin{lemma}[$L^{\infty}$ \textbf{estimate for} $Harmless^{\leq 9}$ \textbf{terms}]
	\label{L:LINFINITYESTIMATEFORHARMLESSLEQ9}
	Under the data-size and bootstrap assumptions 
	of Subsects.~\ref{SS:SIZEOFTBOOT}-\ref{SS:PSIBOOTSTRAP}
	and the smallness assumptions of Subsect.~\ref{SS:SMALLNESSASSUMPTIONS}, 
	the following pointwise estimates hold for the terms 
	$Harmless^{\leq 9}$ from Def.~\ref{D:HARMLESSTERMS}
	on $\mathcal{M}_{\Tboot,U_0}$:
	\begin{align} \label{E:LINFINITYESTIMATEFORHARMLESSLEQ9}
		\left|
			Harmless^{\leq 9}
		\right|
		\lesssim \varepsilon.
	\end{align}
\end{lemma}

\begin{proof}
	The estimate \eqref{E:LINFINITYESTIMATEFORHARMLESSLEQ9} follows directly from
	the $L^{\infty}$ estimates of Prop.~\ref{P:IMPROVEMENTOFAUX}.
\end{proof}

\subsection{Identification of the key difficult error term factors}
\label{SS:IDENTIFICATIONOFKEYDIFFICULTERRORTERMFACTORS}
In the next proposition, we identify the products that are difficult to control
in the energy estimates.

\begin{proposition}[\textbf{Identification of the key difficult error term factors}]
\label{P:IDOFKEYDIFFICULTENREGYERRORTERMS}
For $1 \leq N \leq 18$ in \eqref{E:LISTHEFIRSTCOMMUTATORIMPORTANTTERMS} 
and $N \leq 18$ in \eqref{E:GEOANGANGISTHEFIRSTCOMMUTATORIMPORTANTTERMS},
we have the following estimates:
\begin{subequations}
\begin{align}
	\upmu \square_g (\GeoAng^{N-1} \Lunit \Psi)
	& = (\angdiffuparg{\#} \Psi) \cdot (\upmu \angdiff \GeoAng^{N-1} \mytr \upchi)
			+ Harmless^{\leq N},
			\label{E:LISTHEFIRSTCOMMUTATORIMPORTANTTERMS} \\
	\upmu \square_g (\GeoAng^N \Psi)
	& = (\Rad \Psi) \GeoAng^N \mytr \upchi
			+ \GeoAngFlatRadComponent (\angdiffuparg{\#} \Psi) \cdot (\upmu \angdiff \GeoAng^{N-1} \mytr \upchi)
				 + Harmless^{\leq N}.
				 \label{E:GEOANGANGISTHEFIRSTCOMMUTATORIMPORTANTTERMS}
\end{align}
Furthermore, if $2 \leq N \leq 18$ and $\Tanset^N$ is any $N^{th}$ order 
$\mathcal{P}_u-$tangent operator except for $\GeoAng^{N-1} \Lunit$ or $\GeoAng^N$,
then
\begin{align} \label{E:HARMLESSORDERNCOMMUTATORS}
	\upmu \square_g (\Tanset^N \Psi)
	& = Harmless^{\leq N}.
\end{align}
\end{subequations}
\end{proposition}

We provide the proof of Prop.~\ref{P:IDOFKEYDIFFICULTENREGYERRORTERMS} in 
Subsect.~\ref{SS:PROOFOFPROPIDOFKEYDIFFICULTENREGYERRORTERMS}.
In Subsects.~\ref{SS:PRELIMINARYCOMMUTATIONLEMMAS}-\ref{SS:IMPORTANTDEFORMATIONTENSORTERMS}, 
we establish some preliminary identities and estimates.
As in Sect.~\ref{S:MODQUANTS}, 
we use the Riemann curvature tensor of $g$
to aid our calculations.

\subsection{Preliminary lemmas connected to commutation} 
\label{SS:PRELIMINARYCOMMUTATIONLEMMAS}
We start with a lemma that provides an identity for the
curvature component $\mytr \Cur_{\Rad \cdot \Lunit \cdot}$.
It is an analog of Lemma~\ref{L:ALPHARENORMALIZED}.

\begin{lemma}[\textbf{An expression for} $\mytr \Cur_{\Rad \cdot \Lunit \cdot}$]
\label{L:GAMMACURVATURECOMPONENTEXPRESSION}
Let $\Cur_{\alpha \beta \kappa \lambda}$ be the Riemann curvature tensor from Def.~\ref{D:SPACETIMECURVATURE}.
Then the curvature component 
$\mytr \Cur_{\Rad \cdot \Lunit \cdot} 
:= 
(\ginversesphere) \cdot \Cur_{\Rad \cdot \Lunit \cdot}
$
can be expressed as follows, where all terms are exact except for $\smoothfunction(\cdots)$:
\begin{align}  \label{E:GAMMACURVATURECOMPONENTEXPRESSION}
	\mytr \Cur_{\Rad \cdot \Lunit \cdot}
	& = \frac{1}{2} 
			\left\lbrace
				\upmu \angGnospacemixedarg{\Radunit}{\#} \cdot \angdiff \Lunit \Psi
				+ \angGnospacemixedarg{\Lunit}{\#} \cdot \angdiff \Rad \Psi
					- \mytr \angG \Lunit \Rad \Psi
					- \upmu G_{\Lunit \Radunit} \angLap \Psi
			\right\rbrace
				\\
& \ \ + \frac{1}{4} 
				\upmu^{-1}
				\angGnospacemixedarg{\Lunit}{\#} \cdot \angGarg{\Lunit}
				(\Rad \Psi)^2
		- \frac{1}{2} 
				\upmu^{-1}
				\angGnospacemixedarg{\Lunit}{\#}
				\cdot
				(\angdiff \upmu) 
				\Rad \Psi
		\notag \\
& \ \ + \frac{1}{2} 
				\left\lbrace
					G_{\Lunit \Radunit} \mytr \upchi \Rad \Psi 
					+ \upmu \mytr \upchi \angGnospacemixedarg{\Lunit}{\#} \cdot \angdiff \Psi
					- \upmu \mytr \upchi  \angGnospacemixedarg{\Radunit}{\#} \cdot \angdiff \Psi
				\right\rbrace
			\notag \\
& \ \   +
				\smoothfunction(\BadVar,\ginversesphere,\angdiff x^1,\angdiff x^2,Z \Psi) \Singletan \Psi.
			\notag
\end{align}

\end{lemma}

\begin{proof}
	Lemma~\ref{L:GAMMACURVATURECOMPONENTEXPRESSION} follows from inserting the schematic relations
	provided by Lemma~\ref{L:SCHEMATICDEPENDENCEOFMANYTENSORFIELDS}
	into the identities derived in \cite{jS2014b}*{Lemma 15.1.3}.
	We therefore do not provide a detailed proof here.
	We remark that the main ideas behind the proof
	are the same as those of Lemma~\ref{L:ALPHARENORMALIZED}.
	In particular, the main idea is to contract the curvature
	tensor $\Cur_{\mu \nu \alpha \beta}$, given by \eqref{E:RECTANGULRCURVATURECOMPONENTS},
	against $\Rad^{\mu} \Lunit^{\alpha} (\ginversesphere)^{\nu \kappa}$
	and to use Lemmas~\ref{L:CONNECTIONLRADFRAME} and \ref{L:SCHEMATICDEPENDENCEOFMANYTENSORFIELDS}
	to express the RHS of the contracted identity 
	in the form written on RHS~\eqref{E:GAMMACURVATURECOMPONENTEXPRESSION}.
\end{proof}

We now use Lemma~\ref{L:GAMMACURVATURECOMPONENTEXPRESSION}
to derive an identity for $\Rad \mytr \upchi$.

\begin{lemma}[\textbf{An expression for} $\Rad \mytr \upchi$ \textbf{in terms of other variables}]
\label{L:RADDIRECTIONCHITRANSPORT}
$\Rad \mytr \upchi$ can be expressed as follows,
where the term $\angLap \upmu$ on RHS~\eqref{E:RADDIRECTIONCHITRANSPORT} 
is exactly depicted and the terms $\smoothfunction(\cdots)$ 
are schematically depicted:
\begin{align}  \label{E:RADDIRECTIONCHITRANSPORT}
	\Rad \mytr \upchi
	& = \angLap \upmu
			+ \smoothfunction(\GdVar,\ginversesphere,\angdiff x^1,\angdiff x^2) \Singletan \Rad \Psi
			+ \smoothfunction(\BadVar,\ginversesphere, \angdiff x^1, \angdiff x^2) \Singletan \Singletan \Psi
		\\
	& \ \ + 
		\smoothfunction(\BadVar,\Rad \Psi, \Singletan \BadVar,\ginversesphere,\angdiff x^1,\angdiff x^2) 
		\cdot (\Singletan \GdVar, \angD^2 x^1,\angD^2 x^2).
		\notag
\end{align}

\end{lemma}

\begin{proof}
	We start with the following analog of \eqref{E:TRCHIEVOLUTION}.
	The proof is similar but is slightly more computationally intensive
	due in part to the fact that $[\Rad,\CoordAng] \neq 0$ 
	(see \cite{jS2014b}*{Lemma 15.1.4} for detailed computations):
\begin{align} \label{E:RADTRCHIFIRSTEXPRESSION}
	\Rad \mytr \upchi
	& = \angLap \upmu
		+ \mytr \angD (\upmu \upzeta)
		- \mytr \Cur_{\Rad \cdot \Lunit \cdot}
			\\
	& \ \
		+ \upmu \upzeta^{\#} \cdot \upzeta
		+ \upzeta^{\#} \cdot \angdiff \upmu
		- \mytr \angdeform{\Rad} \mytr \upchi
		- (\Lunit \upmu) \mytr \upchi
		- \upmu (\mytr \upchi)^2
		+ \mytr \upchi (\upmu \mytr \angk).
		\notag
\end{align}
We now substitute RHS~\eqref{E:GAMMACURVATURECOMPONENTEXPRESSION}
for the third term on RHS~\eqref{E:RADTRCHIFIRSTEXPRESSION}
and use
\eqref{E:TRCHIINTERMSOFOTHERVARIABLES},
\eqref{E:UPMUFIRSTTRANSPORT},
\eqref{E:ZETADECOMPOSED},
and \eqref{E:ANGKDECOMPOSED}
to substitute for
$\mytr \upchi$,
$\Lunit \upmu$,
$\upzeta$,
and $\angk$
on RHS~\eqref{E:RADTRCHIFIRSTEXPRESSION}.
The only important observation is that
the two $\upmu^{-1}-$singular products
on the second line of RHS~\eqref{E:GAMMACURVATURECOMPONENTEXPRESSION}
(generated by the term $-	\mytr \Cur_{\Rad \cdot \Lunit \cdot}$ on RHS~\eqref{E:RADTRCHIFIRSTEXPRESSION})
are, in view of the expression \eqref{E:ZETATRANSVERSAL}
for $\upzeta^{(Trans-\Psi)\#}$, exactly canceled by the corresponding terms 
$\upmu^{-1} \upzeta^{(Trans-\Psi)\#} \cdot \upzeta^{(Trans-\Psi)}$
and $\upmu^{-1} \upzeta^{(Trans-\Psi)\#} \cdot \angdiff \upmu$
(generated, 
in view of equation \eqref{E:ZETADECOMPOSED},
by the products 
$\upmu \upzeta^{\#} \cdot \upzeta$ 
and
$\upzeta^{\#} \cdot \angdiff \upmu$
on RHS~\eqref{E:RADTRCHIFIRSTEXPRESSION}).
Also using Lemma~\ref{L:SCHEMATICDEPENDENCEOFMANYTENSORFIELDS},
we arrive at the desired expression \eqref{E:RADDIRECTIONCHITRANSPORT}.
\end{proof}

We now derive higher-order analogs of Lemma~\ref{L:RADDIRECTIONCHITRANSPORT}.

\begin{lemma}[\textbf{Identification of the only non-}$Harmless^{\leq N}$ \textbf{term in} $\Tanset^{N-1} \Rad \mytr \upchi$]
	\label{L:TANGENTIALCOMMUTEDRADTRCHIANGLAPUPMUCOMPARISON}
	Assume that $1 \leq N \leq 18$.
	Under the data-size and bootstrap assumptions 
	of Subsects.~\ref{SS:SIZEOFTBOOT}-\ref{SS:PSIBOOTSTRAP}
	and the smallness assumptions of Subsect.~\ref{SS:SMALLNESSASSUMPTIONS}, 
	the following estimate holds on
	$\mathcal{M}_{\Tboot,U_0}$
	(see Subsect.~\ref{SS:STRINGSOFCOMMUTATIONVECTORFIELDS} regarding the vectorfield operator notation):
	\begin{align} \label{E:TANGENTIALCOMMUTEDRADTRCHIANGLAPUPMUCOMPARISON}
		\left|
		\Tanset^{N-1} \Rad \mytr \upchi
		- \angLap \Tanset^{N-1} \upmu
		\right|
		& \lesssim
			\left|
				\Fullset_*^{\leq N+1;1} \Psi
			\right|
			+
			\left|
				\Tanset^{\leq N} \GdVar
			\right|
			+
			\varepsilon
			\left|
				\Tanset_*^{[1,N]} \BadVar
			\right|.
		\end{align}
\end{lemma}

\begin{proof}
	We start by applying $\Tanset^{N-1}$ 
	to \eqref{E:RADDIRECTIONCHITRANSPORT}.
	We then decompose the first term as 
	$\Tanset^{N-1} \angLap \upmu 
	= \angLap \Tanset^{N-1} \upmu
		+ [\Tanset^{N-1}, \angLap] \upmu
	$.
	We put the principal term $\angLap \Tanset^{N-1} \upmu$
	on LHS~\eqref{E:TANGENTIALCOMMUTEDRADTRCHIANGLAPUPMUCOMPARISON},
	while to bound 
	$
	\left| 
		[\Tanset^{N-1}, \angLap] \upmu
	\right|
	$
	by $\lesssim$ RHS~\eqref{E:TANGENTIALCOMMUTEDRADTRCHIANGLAPUPMUCOMPARISON}, we use
	the commutator estimate \eqref{E:ANGLAPPURETANGENTIALFUNCTIONCOMMUTATOR}
	with $N-1$ in the role of $N$ and $f = \upmu$,
	the $L^{\infty}$ estimates of Prop.~\ref{P:IMPROVEMENTOFAUX},
	and Cor.~\ref{C:SQRTEPSILONTOCEPSILON}.
	To deduce that the $\Tanset^{N-1}$ derivative of 
	the remaining terms on RHS~\eqref{E:RADDIRECTIONCHITRANSPORT},
	with the exception of terms involving $\angD^2 x^1$
	and $\angD^2 x^2$,
	are bounded in magnitude by $\lesssim$ RHS~\eqref{E:TANGENTIALCOMMUTEDRADTRCHIANGLAPUPMUCOMPARISON},
	we use
	Lemmas~\ref{L:POINTWISEFORRECTANGULARCOMPONENTSOFVECTORFIELDS}
	and
	\ref{L:POINTWISEESTIMATESFORGSPHEREANDITSDERIVATIVES}
	and the $L^{\infty}$ estimates of Prop.~\ref{P:IMPROVEMENTOFAUX}.
	We now bound the $\Tanset^{N-1}$ derivative
	of the terms on RHS~\eqref{E:RADDIRECTIONCHITRANSPORT}
	involving derivatives of $\angD^2 x^1$ 
	and $\angD^2 x^2$.
	We first show that
	$\left| \angD^2 x^i \right|
	\lesssim
	\left|
		\Tanset^{\leq 1} \GdVar
	\right|
	$.
	To this end, we note that inequality \eqref{E:ANGDSQUAREFUNCTIONFIRSTBOUNDINTERMSOFGEOANG}
	and the argument given just below it imply that
	$
	\left| \angD^2 x^i \right| 
	\lesssim 
	\left| \angdiff \GeoAng x^i \right|
	+ \left|
			\angdiff x^i
		\right|
		|\GeoAng \GdVar|
	$. The desired bound
	now follows from the previous inequality, 
	Lemma~\ref{L:POINTWISEFORRECTANGULARCOMPONENTSOFVECTORFIELDS},
	and the $L^{\infty}$ estimates of Prop.~\ref{P:IMPROVEMENTOFAUX}.
	We now show that for $1 \leq M \leq N-1$, we have
	$\left| \angLie_{\Tanset}^M \angD^2 x^i \right|
	\lesssim
	\left|
		\Tanset^{\leq M+1} \GdVar
	\right|
	$.
	To this end, we first decompose
	$\angLie_{\Tanset}^M \angD^2 x^i
	=	\angD^2 \Tanset^M x^i
		+
		[\angLie_{\Tanset}^M , \angD^2] x^i
	$.
	Next, using \eqref{E:ANGDERIVATIVESINTERMSOFTANGENTIALCOMMUTATOR}, 
	we deduce 
	$\left| \angD^2 \Tanset^M x^i \right| 
	\lesssim \left| \angdiff \Tanset^{[1,M+1]} x^i \right|
	$.
	Using Lemma~\ref{L:POINTWISEFORRECTANGULARCOMPONENTSOFVECTORFIELDS},
	we bound the RHS of the previous inequality
	by 
	$
	\lesssim
	\left|
		\Tanset^{\leq M+1} \GdVar
	\right|
	$ 
	as desired.
	To deduce that
	$
	\left|
		[\angLie_{\Tanset}^M , \angD^2] x^i
	\right|
	\lesssim
	\left|
		\Tanset^{\leq M+1} \GdVar
	\right|
	$,
	we also use the commutator estimate \eqref{E:ALTERNATEANGDSQUAREDPURETANGENTIALFUNCTIONCOMMUTATOR} with $f=x^i$
	and the $L^{\infty}$ estimates of Prop.~\ref{P:IMPROVEMENTOFAUX}.
	We have thus shown that for $0 \leq M \leq N-1$, 
	we have
	$\left| \angLie_{\Tanset}^M \angD^2 x^i \right|
	\lesssim
	\left|
		\Tanset^{\leq M+1} \GdVar
	\right|
	$
	Combining the previous inequality with the $L^{\infty}$ estimates of Prop.~\ref{P:IMPROVEMENTOFAUX},
	we find that the $\Tanset^{N-1}$ derivative
	of the products on RHS~\eqref{E:RADDIRECTIONCHITRANSPORT}
	involving a factor $\angD^2 x^i$
	are bounded in magnitude by $\lesssim$ RHS~\eqref{E:TANGENTIALCOMMUTEDRADTRCHIANGLAPUPMUCOMPARISON}.
	We have thus proved the lemma.

\end{proof}

\subsection{The important terms in the derivatives of \texorpdfstring{$\deform{\Lunit}$ and $\deform{\GeoAng}$}{the deformation tensors of the commutation vectorfields}}
\label{SS:IMPORTANTDEFORMATIONTENSORTERMS}

The most difficult terms in our energy estimates depend on the 
``top-order non$-\Lunit-$involving derivatives''
of the eikonal function quantities, which appear in
some frame components of the top derivatives of the deformation tensors
$\deform{\Lunit}$ and $\deform{\GeoAng}$.
In the next lemma, we identify the difficult terms.

\begin{lemma} [\textbf{Identification of the important top-order terms in} $\deform{\Lunit}$ \textbf{and} $\deform{\GeoAng}$]
\label{L:IMPORTANTDEFTENSORTERMS}
Assume that $1 \leq N \leq 18$.
Under the data-size and bootstrap assumptions 
of Subsects.~\ref{SS:SIZEOFTBOOT}-\ref{SS:PSIBOOTSTRAP}
and the smallness assumptions of Subsect.~\ref{SS:SMALLNESSASSUMPTIONS}, 
the following estimates hold on
$\mathcal{M}_{\Tboot,U_0}$
(see Subsect.~\ref{SS:STRINGSOFCOMMUTATIONVECTORFIELDS} regarding the vectorfield operator notation):

\underline{\textbf{Important top-order terms in} $\deform{\Lunit}$.}
We have
	\begin{subequations}
	\begin{align} \label{E:LDEFIMPORTANTRADTRACEANGTERMS}
		\left|
			\Tanset^{N-1} \Rad \mytr  \angdeform{\Lunit}
			- 2 \angLap \Tanset^{N-1} \upmu 
		\right|
		& \lesssim
			\left|
				\Fullset_*^{\leq N+1;1} \Psi
			\right|
			+
			\left|
				\Fullset_*^{\leq N;1} \GdVar
			\right|
			+
			\varepsilon
			\left|
				\Tanset_*^{[1,N]} \BadVar
			\right|,
	\end{align}

	\begin{align}  \label{E:LDEFIMPORTANTLIERADSPHERERADANDAGNDIVSPHERETERMS}
		&
		\left| 
			\angLie_{\Tanset}^{N-1} \angdiffuparg{\#} \deformarg{\Lunit}{\Lunit}{\Rad}
		\right|,
		\,
		\left|
				\Tanset^{N-1} \angdiv \angdeformoneformupsharparg{\Lunit}{\Rad}
				- \angLap \Tanset^{N-1} \upmu
		\right|,
		\,
		\left|
				\angLie_{\Tanset}^{N-1} \angdiffuparg{\#} \mytr \angdeform{\Lunit}
				-  
				2 \angdiffuparg{\#} \Tanset^{N-1} \mytr \upchi
		\right|
			\\
		& \lesssim
			\left|
				\Fullset_*^{\leq N+1;1} \Psi
			\right|
			+
			\left|
				\Fullset_*^{\leq N;1} \GdVar
			\right|
			+
			\varepsilon
			\left|
				\Tanset_*^{[1,N]} \BadVar
			\right|.
			\notag
	\end{align}
	\end{subequations}

	\underline{\textbf{Important top-order terms in} $\deform{\GeoAng}$.}
	We have
	\begin{subequations}
	\begin{align}  \label{E:GEOANGDEFIMPORTANTLIERADSPHERELANDRADTRACESPHERETERMS}
		& \left|
			 \angLie_{\Tanset}^{N-1}	\angLie_{\Rad} \angdeformoneformupsharparg{\GeoAng}{\Lunit}
				+ (\angLap \Tanset^{N-1} \upmu) \GeoAng 
			\right|,
			\,
		\left|
			\Tanset^{N-1} \Rad \mytr  \angdeform{\GeoAng}
			- 2 \GeoAngFlatRadComponent \angLap \Tanset^{N-1} \upmu
		\right|
			\\
		& \lesssim
			\left|
				\Fullset_*^{\leq N+1;1} \Psi
			\right|
			+
			\left|
				\Fullset_*^{\leq N;1} \GdVar
			\right|
			+
			\varepsilon
			\left|
				\Tanset_*^{[1,N]} \BadVar
			\right|,
			\notag
	\end{align}

	\begin{align} \label{E:GEOANGDEFIMPORTANTANGDIVSPHERELANDANGDIFFPILRADTERMS}
		& 
		\left|
			\Tanset^{N-1} \angdiv \angdeformoneformupsharparg{\GeoAng}{\Lunit}
				+ \GeoAng \Tanset^{N-1} \mytr \upchi
		\right|,
 			\,
		\left|
			\Tanset^{N-1} \angdiv \angdeformoneformupsharparg{\GeoAng}{\Rad}
			- \left\lbrace
					\upmu \GeoAng \Tanset^{N-1} \mytr \upchi
					+ \GeoAngFlatRadComponent \angLap \Tanset^{N-1} \upmu
				\right\rbrace
		\right|,
			\\
		& \left| 
				\angLie_{\Tanset}^{N-1} \angdiffuparg{\#} \deformarg{\GeoAng}{\Lunit}{\Rad}
				+ (\angLap \Tanset^{N-1} \upmu) \GeoAng
			\right|,
			\,
		\left|
			\angLie_{\Tanset}^{N-1} \angdiffuparg{\#} \mytr \angdeform{\GeoAng}
			- 2 \GeoAngFlatRadComponent \angdiffuparg{\#} \Tanset^{N-1} \mytr \upchi
		\right|
			\notag \\
		& \lesssim
			\left|
				\Fullset_*^{\leq N+1;1} \Psi
			\right|
			+
			\left|
				\Fullset_*^{\leq N;1} \GdVar
			\right|
			+
			\varepsilon
			\left|
				\Tanset_*^{[1,N]} \BadVar
			\right|.
			\notag
	\end{align}
	\end{subequations}
	Above, within a given inequality, the symbol $\Tanset^{N-1}$ on the LHS always denotes the same
	order $N-1$ $\mathcal{P}_u-$tangent vectorfield operator.
\end{lemma}

\begin{proof}
See Subsect.~\ref{SS:OFTENUSEDESTIMATES} for some comments on the analysis.
We first prove the estimate \eqref{E:GEOANGDEFIMPORTANTLIERADSPHERELANDRADTRACESPHERETERMS} 
for
$
\left|
	\angLie_{\Tanset}^{N-1}	\angLie_{\Rad} \angdeformoneformupsharparg{\GeoAng}{\Lunit}
	+ 
	(\angLap \Tanset^{N-1} \upmu) \GeoAng 
\right|
$.
We apply $\angLie_{\Tanset}^{N-1} \angLie_{\Rad}$ 
to the $\gsphere-$dual of \eqref{E:GEOANGDEFORMSPHEREL}.
Note that $\upchi = \mytr \upchi \gsphere$.
The principal top-order term is generated when all derivatives fall
on $\mytr \upchi$ in the product $-\mytr \upchi \GeoAng$.
The top-order product of interest is therefore
$-(\Tanset^{N-1} \Rad \mytr \upchi) \GeoAng$.
Using
\eqref{E:TANGENTIALCOMMUTEDRADTRCHIANGLAPUPMUCOMPARISON},
Lemma~\ref{L:MOREPRECISEANGSPHERELESTIMATES},
and the $L^{\infty}$ estimates of Prop.~\ref{P:IMPROVEMENTOFAUX},
we see that this top-order product
is equal to 
$- (\angLap \Tanset^{N-1} \upmu) \GeoAng$
(whose negative is found on LHS~\eqref{E:GEOANGDEFIMPORTANTLIERADSPHERELANDRADTRACESPHERETERMS})
plus an error term that is
$\lesssim$ RHS~\eqref{E:GEOANGDEFIMPORTANTLIERADSPHERELANDRADTRACESPHERETERMS}
as desired.
Also using Lemma~\ref{L:POINTWISEESTIMATESFORGSPHEREANDITSDERIVATIVES},
we deduce that the remaining terms 
in the Leibniz expansion of 
$- \angLie_{\Tanset}^{N-1} \angLie_{\Rad} (\mytr \upchi \GeoAng)$
are $\lesssim$ RHS~\eqref{E:GEOANGDEFIMPORTANTLIERADSPHERELANDRADTRACESPHERETERMS}.
Using Lemma~\ref{L:SCHEMATICDEPENDENCEOFMANYTENSORFIELDS}, 
we see that the remaining terms on the $\gsphere-$dual 
of RHS~\eqref{E:GEOANGDEFORMSPHEREL} are of the form
$\smoothfunction(\GdVar,\ginversesphere,\angdiff x^1,\angdiff x^2) \Singletan \Psi$.
Hence, their 
$\angLie_{\Tanset}^{N-1} \angLie_{\Rad}
= \angLie_{\Fullset}^{N;1}
$ 
derivatives can be bounded by 
$\lesssim$ RHS~\eqref{E:GEOANGDEFIMPORTANTLIERADSPHERELANDRADTRACESPHERETERMS}
via the estimates of 
Lemmas~\ref{L:POINTWISEFORRECTANGULARCOMPONENTSOFVECTORFIELDS}
and \ref{L:POINTWISEESTIMATESFORGSPHEREANDITSDERIVATIVES}
and the $L^{\infty}$ estimates of Prop.~\ref{P:IMPROVEMENTOFAUX}.
We have thus proved the desired estimate.
The proof of \eqref{E:GEOANGDEFIMPORTANTLIERADSPHERELANDRADTRACESPHERETERMS}
for 
$
\left|
	\Tanset^{N-1} \Rad \mytr \angdeform{\GeoAng}
	- 
	2 \GeoAngFlatRadComponent \angLap \Tanset^{N-1} \upmu
\right|
$
follows similarly from the identity \eqref{E:GEOANGDEFORMSPHERE} for $\angdeform{\GeoAng}$,
and the proof of 
\eqref{E:LDEFIMPORTANTRADTRACEANGTERMS}
follows similarly from the identity
\eqref{E:LUNITDEFORMSPHERE} for $\angdeform{\Lunit}$;
we omit the details. 

To prove \eqref{E:GEOANGDEFIMPORTANTANGDIVSPHERELANDANGDIFFPILRADTERMS} for
$
\left|
	\Tanset^{N-1} \angdiv \angdeformoneformupsharparg{\GeoAng}{\Lunit}
	+ 
	\GeoAng \Tanset^{N-1} \mytr \upchi
\right|
$,
we first use the commutator estimate \eqref{E:ANGDIVANGLIETANGENTIALTENSORFIELDCOMMUTATORESTIMATE}
with $\xi = \angdeformoneformarg{\GeoAng}{\Lunit}$ and $N-1$ in the role of $N$,
the estimate \eqref{E:TANGENTDIFFERNTIATEDGEOANGDEFORMSPHERELSHARPPOINTWISE},
and the $L^{\infty}$ estimates of Prop.~\ref{P:IMPROVEMENTOFAUX}
to commute the operator
$\angLie_{\Tanset}^{N-1}$ through $\angdiv$
in the term
$\Tanset^{N-1} \angdiv \angdeformoneformupsharparg{\GeoAng}{\Lunit}
= \angLie_{\Tanset}^{N-1} \angdiv \angdeformoneformarg{\GeoAng}{\Lunit} 
$,
thus obtaining that 
$\Tanset^{N-1} \angdiv \angdeformoneformupsharparg{\GeoAng}{\Lunit} 
= \angdiv \angLie_{\Tanset}^{N-1} \angdeformoneformarg{\GeoAng}{\Lunit}
$
up to error terms that are bounded in magnitude by
$\lesssim \mbox{{\upshape RHS}~\eqref{E:GEOANGDEFIMPORTANTANGDIVSPHERELANDANGDIFFPILRADTERMS}}$.
It remains for us to analyze the terms
that arise from applying 
$\angdiv \angLie_{\Tanset}^{N-1} = \ginversesphere \cdot \angD \angLie_{\Tanset}^{N-1}$ to
RHS~\eqref{E:GEOANGDEFORMSPHEREL}.
As in the previous paragraph, the principal top-order term is generated when all derivatives fall
on $\mytr \upchi$ in the product 
$
- \upchi \cdot \GeoAng
= -\mytr \upchi \GeoAng_{\flat} 
= -\mytr \upchi \gsphere \cdot \GeoAng$.
The top-order product of interest
(whose negative is found on LHS~\eqref{E:GEOANGDEFIMPORTANTANGDIVSPHERELANDANGDIFFPILRADTERMS}), 
is therefore
$- (\angdiff \Tanset^{N-1} \mytr \upchi) \cdot \GeoAng 
= 
- \GeoAng \Tanset^{N-1} \mytr \upchi$.
Moreover, 
we see that the remaining terms in the Leibniz expansion of
$- \angdiv \angLie_{\Tanset}^{N-1} (\mytr \upchi \GeoAng_{\flat})$ are 
$
\displaystyle
\lesssim
\mathop{\sum_{N_1 + N_2 + N_3 \leq N}}_{N_1 \leq N-1}
\left|
	\Tanset^{N_1} \mytr \upchi
\right|
\left|
	\angLie_{\Tanset}^{N_2} 
	\gsphere
\right|
\left|
	\angLie_{\Tanset}^{N_3} 
	\GeoAng
\right|
$.
We now bound these terms by
$\lesssim$ RHS~\eqref{E:GEOANGDEFIMPORTANTANGDIVSPHERELANDANGDIFFPILRADTERMS}
via the estimates of 
Lemmas~\ref{L:POINTWISEESTIMATESFORGSPHEREANDITSDERIVATIVES}
and~\ref{L:MOREPRECISEANGSPHERELESTIMATES}
and the $L^{\infty}$ estimates of Prop.~\ref{P:IMPROVEMENTOFAUX}.
Using Lemma~\ref{L:SCHEMATICDEPENDENCEOFMANYTENSORFIELDS}, we see that 
the remaining terms 
on RHS~\eqref{E:GEOANGDEFORMSPHEREL} are of the form
$\smoothfunction(\GdVar,\angdiff x^1,\angdiff x^2) \Singletan \Psi$.
To bound their 
$
 \angdiv \angLie_{\Tanset}^{N-1}
$
derivatives
by $\lesssim$ RHS~\eqref{E:GEOANGDEFIMPORTANTANGDIVSPHERELANDANGDIFFPILRADTERMS},
we use the same arguments given in the previous paragraph with 
the minor new addition that we also need the bound
$\left| \angD^2 \Tanset^{N-1} x^i \right| 
\lesssim
	\left|
		\Tanset^{\leq N} \GdVar
	\right|
$,
which we obtained in the proof of Lemma~\ref{L:TANGENTIALCOMMUTEDRADTRCHIANGLAPUPMUCOMPARISON}. 
We have thus obtained the desired result.
The proof of \eqref{E:GEOANGDEFIMPORTANTANGDIVSPHERELANDANGDIFFPILRADTERMS} 
for
$
\left|
			\Tanset^{N-1} \angdiv \angdeformoneformupsharparg{\GeoAng}{\Rad}
			- \left\lbrace
					\upmu \GeoAng \Tanset^{N-1} \mytr \upchi
					+ \GeoAngFlatRadComponent \angLap \Tanset^{N-1} \upmu
				\right\rbrace
		\right|
$
follows similarly from the identity
\eqref{E:GEOANGDEFORMSPHERERAD}
for $\angdeformoneformarg{\GeoAng}{\Rad}$,
the fact that 
$\GeoAngFlatRadComponent = \smoothfunction(\GdVar) \GdVar$
(see Lemma~\ref{L:SCHEMATICDEPENDENCEOFMANYTENSORFIELDS}),
the estimate \eqref{E:TANGENTDIFFERNTIATEDGEOANGDEFORMSPHERERADPOINTWISE},
and Cor.~\ref{C:SQRTEPSILONTOCEPSILON};
we omit the details, noting only that
Cor.~\ref{C:SQRTEPSILONTOCEPSILON} allows 
us to replace the factor $\varepsilon^{1/2}$
on RHS~\eqref{E:ANGDIVANGLIETANGENTIALTENSORFIELDCOMMUTATORESTIMATE}
with $C \varepsilon$.
The proof of \eqref{E:LDEFIMPORTANTLIERADSPHERERADANDAGNDIVSPHERETERMS}
for
$
\left|
	\Tanset^{N-1} \angdiv \angdeformoneformupsharparg{\Lunit}{\Rad}
	- \angLap \Tanset^{N-1} \upmu
\right|
$
follows similarly from the identity
\eqref{E:LUNITDEFORMSPHERELUNITANDLUNITDEFORMSPHERERAD}
for $\angdeformoneformarg{\Lunit}{\Rad}$,
the estimates 
\eqref{E:RADDEFORMSPHERELPOINTWISE}-\eqref{E:TANGENTDIFFERNTIATEDRADDEFORMSPHERELPOINTWISE},
and Cor.~\ref{C:SQRTEPSILONTOCEPSILON};
we omit the details, noting only that
Cor.~\ref{C:SQRTEPSILONTOCEPSILON} allows 
us to replace the factor $\varepsilon^{1/2}$
on RHS~\eqref{E:ANGDIVANGLIETANGENTIALTENSORFIELDCOMMUTATORESTIMATE}
with $C \varepsilon$.
The proof of \eqref{E:LDEFIMPORTANTLIERADSPHERERADANDAGNDIVSPHERETERMS}
for
$
\left|
	\angLie_{\Tanset}^{N-1} \angdiffuparg{\#} \deformarg{\Lunit}{\Lunit}{\Rad}
\right|
$
is based on the identity \eqref{E:LUNITDEFORMSCALARS}
for $\deformarg{\Lunit}{\Lunit}{\Rad}$
and the estimate \eqref{E:PURETANGENTIALLUNITUPMUCOMMUTEDESTIMATE}
and is similar but simpler; we omit the details.
The proof of \eqref{E:LDEFIMPORTANTLIERADSPHERERADANDAGNDIVSPHERETERMS}
for
$
\left|
	\angLie_{\Tanset}^{N-1} \angdiffuparg{\#} \mytr \angdeform{\Lunit}
	-  
	2 \angdiffuparg{\#} \Tanset^{N-1} \mytr \upchi
\right|
$
is similar and is based on the identity \eqref{E:LUNITDEFORMSPHERE}
for $\angdeform{\Lunit}$; we omit the details.
The proof of \eqref{E:GEOANGDEFIMPORTANTANGDIVSPHERELANDANGDIFFPILRADTERMS} 
for
$
\left| 
	\angLie_{\Tanset}^{N-1} \angdiffuparg{\#} \deformarg{\GeoAng}{\Lunit}{\Rad}
	+ 
	(\angLap \Tanset^{N-1} \upmu) \GeoAng
\right|
$
follows similarly from the identity \eqref{E:GEOANGDEFORMSCALARS} for 
$\deformarg{\GeoAng}{\Lunit}{\Rad}$
and the trivial identity
$- \GeoAng \upmu = - \GeoAng \cdot \angdiff \upmu$
relevant for the term on the RHS of the identity; we omit the details.
The proof of \eqref{E:GEOANGDEFIMPORTANTANGDIVSPHERELANDANGDIFFPILRADTERMS} 
for
$
\left|
			\angLie_{\Tanset}^{N-1} \angdiffuparg{\#} \mytr \angdeform{\GeoAng}
			- 2 \GeoAngFlatRadComponent \angdiffuparg{\#} \Tanset^{N-1} \mytr \upchi
		\right|
$
follows similarly from the identity \eqref{E:GEOANGDEFORMSPHERE} for 
$\angdeform{\GeoAng}$;
we omit the details.
\end{proof}

The top-order derivatives of $\deform{\Lunit}$ and $\deform{\GeoAng}$
involving at least one $\Lunit$ differentiation
and the below-top-order derivatives of
$\deform{\Lunit}$ and $\deform{\GeoAng}$
lead to negligible error terms in the energy estimates.
In the next lemma, we derive the relevant pointwise estimates
that will allow us to establish this fact.

\begin{lemma}[\textbf{Pointwise estimates for the negligible derivatives of $\deform{\Lunit}$ and $\deform{\GeoAng}$}]
		\label{L:TANGENTIALDERIVATIVESOFDEFORMATIONINVOLVINGONELUNITDERIVATIVE}
	 Assume that $1 \leq N \leq 18$ and let
	 	$\Singletan \in \Tanset = \lbrace \Lunit, \GeoAng \rbrace$.
	 	Under the data-size and bootstrap assumptions 
		of Subsects.~\ref{SS:SIZEOFTBOOT}-\ref{SS:PSIBOOTSTRAP}
		and the smallness assumptions of Subsect.~\ref{SS:SMALLNESSASSUMPTIONS}, 
		the following estimates hold on
		$\mathcal{M}_{\Tboot,U_0}$:
		\begin{align} \label{E:TANGENTIALDERIVATIVESOFDEFORMATIONINVOLVINGONELUNITDERIVATIVE}
		&
		\left|
			\Tanset^{\leq N-1} \Lunit \mytr \angdeform{\Singletan}
		\right|,
			\,
		\left|
			\Tanset^{\leq N-1} \Lunit \deformarg{\Singletan}{\Lunit}{\Rad}
		\right|,
			\,
		\left|
			\Tanset^{\leq N-1} \Lunit \deformarg{\Singletan}{\Rad}{\Radunit}
		\right|,
			\,
		\left|
			\angLie_{\Tanset}^{\leq N-1} \angLie_{\Lunit} \angdeformoneformupsharparg{\Singletan}{\Lunit}
		\right|,
			\,
		\left|
			\angLie_{\Tanset}^{\leq N-1} \angLie_{\Lunit} \angdeformoneformupsharparg{\Singletan}{\Rad}
		\right|
			\\
		& \lesssim
			\left|
				\Fullset_*^{\leq N+1;1} \Psi
			\right|
			+
			\left|
				\Tanset^{\leq N} \GdVar
			\right|
			+
			\left|
				\Tanset_*^{[1,N]} \BadVar
			\right|
			\notag.
	\end{align}

	Moreover, the following below-top-order estimates hold:
	\begin{align} \label{E:TANGENTIALDERIVATIVESOFDEFORMATION}
		&
		\left|
			\Tanset^{\leq N-1} \mytr \angdeform{\Singletan}
		\right|,
			\,
		\left|
			\Tanset^{\leq N-1} \deformarg{\Singletan}{\Lunit}{\Rad}
		\right|,
			\,
		\left|
			\Tanset^{\leq N-1} \deformarg{\Singletan}{\Rad}{\Radunit}
		\right|,
			\,
		\left|
			\angLie_{\Tanset}^{\leq N-1} \angdeformoneformupsharparg{\Singletan}{\Lunit}
		\right|,
			\,
		\left|
			\angLie_{\Tanset}^{\leq N-1} \angdeformoneformupsharparg{\Singletan}{\Rad}
		\right|
			\\
		& \lesssim
			\left|
				\Fullset_*^{\leq N;1} \Psi
			\right|
			+
			\left|
				\Tanset^{\leq N} \GdVar
			\right|
			+
			\left|
				\Tanset_*^{[1,N]} \BadVar
			\right|
			+
			1.
			\notag
	\end{align}

\end{lemma}

\begin{proof}
	See Subsect.~\ref{SS:OFTENUSEDESTIMATES} for some comments on the analysis.
	We first prove \eqref{E:TANGENTIALDERIVATIVESOFDEFORMATIONINVOLVINGONELUNITDERIVATIVE}.
	From Prop.~\ref{L:DEFORMATIONTENSORFRAMECOMPONENTS},
	equation \eqref{E:UPMUFIRSTTRANSPORT},
	and Lemma~\ref{L:SCHEMATICDEPENDENCEOFMANYTENSORFIELDS},
	we see that the deformation tensor components 
	$\mytr \angdeform{\Singletan}$,
	$\deformarg{\Singletan}{\Lunit}{\Rad}$,
	$\cdots$,
	$\angdeformoneformupsharparg{\Singletan}{\Rad}$
	on LHS~\eqref{E:TANGENTIALDERIVATIVESOFDEFORMATIONINVOLVINGONELUNITDERIVATIVE}
	are schematically of the form
	\[
	\smoothfunction(\BadVar,\ginversesphere,\angdiff x^1,\angdiff x^2) \Singletan \Psi
	+ \smoothfunction(\GdVar,\ginversesphere,\angdiff x^1,\angdiff x^2) 
		\Rad \Psi 
	+ \smoothfunction(\GdVar,\ginversesphere,\angdiff x^1,\angdiff x^2) 
		\mytr \upchi
	+ \smoothfunction(\GdVar,\ginversesphere,\angdiff x^1,\angdiff x^2) 
		\angdiff \upmu.
	\]
	We now apply
	$\angLie_{\Tanset}^{\leq N-1} \angLie_{\Lunit}$.
	If all derivatives fall on 
	$\mytr \upchi$, then we use
	\eqref{E:LUNITTANGENTDIFFERENTIATEDLUNITSMALLIMPROVEDPOINTWISE}
	to bound $\Tanset^{\leq N-1} \Lunit \mytr \upchi$ and
	we bound the remaining factors multiplying $\Tanset^{\leq N-1} \Lunit \mytr \upchi$ by $\lesssim 1$ via 
	Lemmas~\ref{L:POINTWISEFORRECTANGULARCOMPONENTSOFVECTORFIELDS}
	and \ref{L:POINTWISEESTIMATESFORGSPHEREANDITSDERIVATIVES}
	and the $L^{\infty}$ estimates of Prop.~\ref{P:IMPROVEMENTOFAUX}.
	Similarly, if all derivatives fall on
	$\angdiff \upmu$, 
	we bound
	$
		\angdiff \Tanset^{\leq N-1} \Lunit \upmu
	$
	with \eqref{E:PURETANGENTIALLUNITUPMUCOMMUTEDESTIMATE}
	and we bound the remaining factors multiplying 
	$
		\angdiff \Tanset^{\leq N-1} \Lunit \upmu
	$
	by $\lesssim 1$.
	If most (but not all) derivatives fall on
	$\mytr \upchi$ or $\angdiff \upmu$, then we bound all terms
	using the above arguments
	and also \eqref{E:POINTWISEESTIMATESFORGSPHEREANDITSTANGENTIALDERIVATIVES}.
	If most derivatives fall on 
	$\Singletan \Psi$ or $\Rad \Psi$,  
	then we bound these factors by the first term on 
	RHS~\eqref{E:TANGENTIALDERIVATIVESOFDEFORMATIONINVOLVINGONELUNITDERIVATIVE} 
	and use the above arguments
	to bound the remaining factors by $\lesssim 1$.

	The proof of \eqref{E:TANGENTIALDERIVATIVESOFDEFORMATION}
	is similar but simpler and we therefore omit the details.

\end{proof}

\subsection{Proof of Prop.~\ref{P:IDOFKEYDIFFICULTENREGYERRORTERMS}}
\label{SS:PROOFOFPROPIDOFKEYDIFFICULTENREGYERRORTERMS}
We now use the previous results to establish Prop.~\ref{P:IDOFKEYDIFFICULTENREGYERRORTERMS}.
See Subsect.~\ref{SS:OFTENUSEDESTIMATES} for some comments on the analysis.
Throughout the proof, we silently use the definition of $Harmless^{\leq N}$ terms from Def.~\ref{D:HARMLESSTERMS}
and the estimates of Lemma~\ref{L:LINFINITYESTIMATEFORHARMLESSLEQ9}.
We give a detailed proof of \eqref{E:GEOANGANGISTHEFIRSTCOMMUTATORIMPORTANTTERMS} and then at the end,
we sketch the minor changes needed to prove \eqref{E:LISTHEFIRSTCOMMUTATORIMPORTANTTERMS}.
To condense the notation, we define	the following commutation vectorfield 
${\Jcurrent{Z}^{\alpha}[\Psi]}$, which is just alternate notation for the term in braces on RHS~\eqref{E:BOXZCOM}:
\begin{align} \label{E:COMMCURRENTZ}
{\Jcurrent{Z}^{\alpha}[\Psi]}
			:= \deformuparg{Z}{\alpha}{\beta} \D_{\beta} \Psi 
			- \frac{1}{2} \myspacetimetr \deform{Z} \D^{\alpha} \Psi.
\end{align}
Iterating \eqref{E:BOXZCOM}, 
using $\square_{g(\Psi)} \Psi = 0$, 
and using the estimate
$
\left\|
	\Tanset^{\leq 9} \mytr \angdeform{\GeoAng}
\right\|_{L^{\infty}(\Sigma_t^u)} \lesssim \varepsilon
$
(which follows from \eqref{E:CONNECTIONBETWEENANGLIEOFGSPHEREANDDEFORMATIONTENSORS}, 
\eqref{E:POINTWISEESTIMATESFORGSPHEREANDITSTANGENTIALDERIVATIVES},
and the $L^{\infty}$ estimates of Prop.~\ref{P:IMPROVEMENTOFAUX}),
we find that
\begin{align} \label{E:COMMUTEDWAVEGEOANGFIRSTMAINTERMPLUSERRORTERM}
		\upmu \square_{g(\Psi)} (\GeoAng^N \Psi)
		&  =  \GeoAng^{N-1}
					\left(
					\upmu 
					\D_{\alpha} 
					\Jcurrent{\GeoAng}^{\alpha}[\Psi]
					\right)
				+ \mbox{Error},
\end{align}
where
\begin{align} \label{E:GEOANGLOWERORDERRORTERMESTIMATE}
	\left|
		\mbox{Error}
	\right|
	& \lesssim 
		\mathop{\sum_{N_1 + N_2 + N_3 \leq N-1}}_{N_1, N_2 \leq N-2}
		\left(
			1 +
			\left|
				\GeoAng^{N_1} \mytr \angdeform{\GeoAng}
			\right|
		\right)
		\left|
			\GeoAng^{N_2}
			\left(
				\upmu 
				\D_{\alpha} 
				\Jcurrent{\GeoAng}^{\alpha}[\GeoAng^{N_3} \Psi]
			\right)
		\right|.
\end{align}

We first analyze the main term 
$
\displaystyle
					\GeoAng^{N-1}
					\left(
					\upmu 
					\D_{\alpha} 
					\Jcurrent{\GeoAng}^{\alpha}[\Psi]
					\right)
$
on RHS~\eqref{E:COMMUTEDWAVEGEOANGFIRSTMAINTERMPLUSERRORTERM},
which contains all of the top-order derivatives of the eikonal
function quantities. By \eqref{E:DIVCOMMUTATIONCURRENTDECOMPOSITION}, we may equivalently 
analyze 
$
\GeoAng^{N-1} \mathscr{K}_{(\pi-Danger)}^{(\GeoAng)}[\Psi]
+ \GeoAng^{N-1} \mathscr{K}_{(\pi-Cancel-1)}^{(\GeoAng)}[\Psi]
+ \cdots
+ \GeoAng^{N-1} \mathscr{K}_{(Low)}^{(\GeoAng)}[\Psi].
$
We argue one term at a time.

\medskip
\noindent
\underline{\textbf{Analysis of} $\GeoAng^{N-1} \mathscr{K}_{(\pi-Danger)}^{(\GeoAng)}[\Psi]$.}
By \eqref{E:DIVCURRENTTRANSVERSAL}, we have
\begin{align} \label{E:GEOANGDANGERTERMMAINPOINTWISEESTIMATE}
	\GeoAng^{N-1} \mathscr{K}_{(\pi-Danger)}^{(\GeoAng)}[\Psi]
	& = - \sum_{N_1 + N_2 = N-1}
				(\GeoAng^{N_1} \angdiv \angdeformoneformupsharparg{\GeoAng}{\Lunit}) 
				\GeoAng^{N_2} \Rad \Psi.
\end{align}
We first consider the case $N_1 = N-1$.
Inequality \eqref{E:PSITRANSVERSALLINFINITYBOUNDBOOTSTRAPIMPROVED}
and the first inequality in \eqref{E:GEOANGDEFIMPORTANTANGDIVSPHERELANDANGDIFFPILRADTERMS}
yield that 
$
- (\GeoAng^{N-1} \angdiv \angdeformoneformupsharparg{\GeoAng}{\Lunit}) \Rad \Psi
= (\GeoAng^N \mytr \upchi) \Rad \Psi
+ Harmless^{\leq N}
$,
which in particular yields the desired first product on 
RHS~\eqref{E:GEOANGANGISTHEFIRSTCOMMUTATORIMPORTANTTERMS}.
To show that the remaining summands are $Harmless^{\leq N}$, 
we again use the first inequality in 
\eqref{E:GEOANGDEFIMPORTANTANGDIVSPHERELANDANGDIFFPILRADTERMS}
(now with $N_1+1$ in the role of $N$)
to deduce that
$
(\GeoAng^{N_1} \angdiv \angdeformoneformupsharparg{\GeoAng}{\Lunit}) \GeoAng^{N_2} \Rad \Psi
= (\GeoAng^{N_1+1} \mytr \upchi) \GeoAng^{N_2} \Rad \Psi 
	+ (Harmless^{\leq N_1}) \GeoAng^{N_2} \Rad \Psi
$.
Since $N_1 \leq N - 2$, 
\eqref{E:POINTWISEESTIMATESFORGSPHEREANDITSTANGENTIALDERIVATIVES}
implies that
$
\GeoAng^{N_1+1} \mytr \upchi 
= Harmless^{\leq N_1 + 2}
\leq Harmless^{\leq N}
$.
From these estimates and the $L^{\infty}$ estimates of Prop.~\ref{P:IMPROVEMENTOFAUX}, 
we easily conclude that 
$
(\GeoAng^{N_1} \angdiv \angdeformoneformupsharparg{\GeoAng}{\Lunit}) \GeoAng^{N_2} \Rad \Psi 
= Harmless^{\leq N}
$
as desired.

\medskip
\noindent
\underline{\textbf{Analysis of} $\GeoAng^{N-1} \mathscr{K}_{(\pi-Cancel-1)}^{(\GeoAng)}[\Psi]$.}
From \eqref{E:DIVCURRENTCANEL1}, we have
\begin{align} \label{E:GEOANGCANCELTERM1TERMMAINPOINTWISEESTIMATE}
	& \GeoAng^{N-1} \mathscr{K}_{(\pi-Cancel-1)}^{(\GeoAng)}[\Psi]
		\\
	& =  \sum_{N_1 + N_2 = N-1}
			\left\lbrace 
				\frac{1}{2} \GeoAng^{N_1} \Rad \mytr \angdeform{\GeoAng}
				- \GeoAng^{N_1} \angdiv \angdeformoneformupsharparg{\GeoAng}{\Rad}
				- \GeoAng^{N_1} (\upmu \angdiv \angdeformoneformupsharparg{\GeoAng}{\Lunit})
			\right\rbrace 
			\GeoAng^{N_2} \Lunit \Psi.
			\notag
\end{align}
We first consider the case in which $N_1 = N-1$ in all three terms in braces on RHS
\eqref{E:GEOANGCANCELTERM1TERMMAINPOINTWISEESTIMATE} and all derivatives fall
on the deformation tensor components.
From the second inequality in \eqref{E:GEOANGDEFIMPORTANTLIERADSPHERELANDRADTRACESPHERETERMS}
and the first and second inequalities in
\eqref{E:GEOANGDEFIMPORTANTANGDIVSPHERELANDANGDIFFPILRADTERMS},
we see that the main top-order eikonal function terms
$
	2 \GeoAngFlatRadComponent \angLap \GeoAng^{N-1} \upmu
$
and 
$
	\upmu \GeoAng^N \mytr \upchi
$
completely cancel from the terms in braces,
leaving only products of the form 
$Harmless^{\leq N} \times \Lunit \Psi = Harmless^{\leq N}$.
When $N_1 = N-1$, there are also terms in which at least one derivative 
falls on the factor $\upmu$ in the product
$
\left\lbrace
	\GeoAng^{N-1} (\upmu \angdiv \angdeformoneformupsharparg{\GeoAng}{\Lunit})
\right\rbrace
\Lunit \Psi
$
on RHS~\eqref{E:GEOANGCANCELTERM1TERMMAINPOINTWISEESTIMATE}.
We will show that these terms $= Harmless^{\leq N}$.
We first bound them by
$
\lesssim
	\mathop{\sum_{N_1 + N_2 \leq N-1}}_{N_2 \leq N-2} 
		\left|
			\GeoAng^{N_1} \upmu
		\right| 
		\left|
			\GeoAng^{N_2} \angdiv \angdeformoneformupsharparg{\GeoAng}{\Lunit} 
		\right|
		\left|
			\Lunit \Psi
		\right|
$.
Again using the first inequality in
\eqref{E:GEOANGDEFIMPORTANTANGDIVSPHERELANDANGDIFFPILRADTERMS}
(now with $N_2+1$ in the role of $N$)
to control
$
\GeoAng^{N_2} \angdiv \angdeformoneformupsharparg{\GeoAng}{\Lunit} 
$,
we bound the RHS of the previous inequality by
\[
\displaystyle
\lesssim
	\mathop{\sum_{N_1 + N_2 \leq N-1}}_{N_2 \leq N-2} 
	Harmless^{\leq N_1} 
	\left|
		\GeoAng^{N_2+1} \mytr \upchi
		+ Harmless^{\leq N_2}
	\right|
	\left|
		\Lunit \Psi
	\right|.
\]
Since $N_2 \leq N-2$,
the arguments given in our analysis of 
$\GeoAng^{N-1} \mathscr{K}_{(\pi-Danger)}^{(\GeoAng)}[\Psi]$
yield that
$
\GeoAng^{N_2+1} \mytr \upchi
= Harmless^{\leq N}
$.
From these estimates and the $L^{\infty}$ estimates of Prop.~\ref{P:IMPROVEMENTOFAUX}, 
we easily conclude that 
the terms under consideration $= Harmless^{\leq N}$ as desired.
We now consider the remaining cases, 
in which $N_1 \leq N-2$ in all three terms in braces on RHS
\eqref{E:GEOANGCANCELTERM1TERMMAINPOINTWISEESTIMATE}.
Again using second inequality in \eqref{E:GEOANGDEFIMPORTANTLIERADSPHERELANDRADTRACESPHERETERMS}
and the first and second inequalities in
\eqref{E:GEOANGDEFIMPORTANTANGDIVSPHERELANDANGDIFFPILRADTERMS}
(now with $N_1+1$ in the role of $N$)
and the arguments given above,
we deduce that the products under consideration 
$
= Harmless^{\leq N}
$
plus error products generated by the terms on 
LHS~\eqref{E:GEOANGDEFIMPORTANTLIERADSPHERELANDRADTRACESPHERETERMS} and
LHS~\eqref{E:GEOANGDEFIMPORTANTANGDIVSPHERELANDANGDIFFPILRADTERMS}.
The error products are in magnitude
\[
\lesssim
	\mathop{\sum_{N_1 + N_2 + N_3 \leq N-1}}_{N_1 \leq N-2} 
	\left|
		\GeoAng^{N_1+1} \mytr \upchi
	\right|
	\left|
		\GeoAng^{N_2} \upmu
	\right| 
	\left|
		\GeoAng^{N_3} \Lunit \Psi
	\right|
	+
	\mathop{\sum_{N_1 + N_2 \leq N-1}}_{N_1 \leq N-2} 
	\left|
		\GeoAngFlatRadComponent \angLap \GeoAng^{N_1} \upmu
	\right|
	\left|
		\GeoAng^{N_2} \Lunit \Psi
	\right|.
\]
Since $N_1 \leq N-2$, the arguments given above 
and the fact that $\GeoAngFlatRadComponent = \smoothfunction(\GdVar) \GdVar$ (see Lemma~\ref{L:SCHEMATICDEPENDENCEOFMANYTENSORFIELDS})
combine to yield that the RHS of the previous expression $= Harmless^{\leq N}$ as desired.

\medskip
\noindent
\underline{\textbf{Analysis of} $\GeoAng^{N-1} \mathscr{K}_{(\pi-Cancel-2)}^{(\GeoAng)}[\Psi]$.}
From \eqref{E:DIVCURRENTCANEL2}, we have
\begin{align} \label{E:GEOANGCANCELTERM2TERMMAINPOINTWISEESTIMATE}
	\GeoAng^{N-1} \mathscr{K}_{(\pi-Cancel-2)}^{(\GeoAng)}[\Psi]
	& = \sum_{N_1 + N_2 = N-1}
			\left\lbrace
					- \angLie_{\GeoAng}^{N_1} \angLie_{\Rad} \angdeformoneformupsharparg{\GeoAng}{\Lunit}
					+ \angLie_{\GeoAng}^{N_1} \angdiffuparg{\#} \deformarg{\GeoAng}{\Lunit}{\Rad}
				\right\rbrace 
				\cdot
				\angdiff \GeoAng^{N_2} \Psi.
\end{align}
In the case $N_1 = N-1$, we use
the first inequality in \eqref{E:GEOANGDEFIMPORTANTLIERADSPHERELANDRADTRACESPHERETERMS},
and the third inequality in
\eqref{E:GEOANGDEFIMPORTANTANGDIVSPHERELANDANGDIFFPILRADTERMS}
to deduce
that the top-order eikonal function terms
$
(\angLap \Tanset^{N-1} \upmu) \GeoAng
$
completely cancel from the terms in braces
on RHS~\eqref{E:GEOANGCANCELTERM2TERMMAINPOINTWISEESTIMATE}.
The remaining analysis now parallels our analysis of
$\GeoAng^{N-1} \mathscr{K}_{(\pi-Cancel-1)}^{(\GeoAng)}[\Psi]$,
with the minor addition that we must also use
the estimate \eqref{E:GEOANGPOINTWISE} to bound the factors of $|\GeoAng|$ that arise.
We thus conclude 
that RHS~\eqref{E:GEOANGCANCELTERM2TERMMAINPOINTWISEESTIMATE} 
$= Harmless^{\leq N}$ as desired.

\medskip
\noindent
\underline{\textbf{Analysis of} $\GeoAng^{N-1} \mathscr{K}_{(\pi-Less \ Dangerous)}^{(\GeoAng)}[\Psi]$.}
From \eqref{E:DIVCURRENTELLIPTIC}, we have
\begin{align} \label{E:GEOANGELLIPTICTERMMAINPOINTWISEESTIMATE}
	\GeoAng^{N-1} \mathscr{K}_{(\pi-Less \ Dangerous)}^{(\GeoAng)}[\Psi]
	& = 	\sum_{N_1 + N_2 = N-1}
				\frac{1}{2} 
				\left\lbrace 
					\angLie_{\GeoAng}^{N_1} (\upmu \angdiffuparg{\#} \mytr \angdeform{\GeoAng})  
				\right\rbrace
				\cdot
				\angdiff \GeoAng^{N_2} \Psi.
\end{align}
We first consider the case $N_1 = N-1$ on RHS~\eqref{E:GEOANGELLIPTICTERMMAINPOINTWISEESTIMATE}
and all derivatives fall on $\mytr \angdeform{\GeoAng}$.
Using the fourth inequality in \eqref{E:GEOANGDEFIMPORTANTANGDIVSPHERELANDANGDIFFPILRADTERMS},
we see that
$
\frac{1}{2} \upmu (\angdiffuparg{\#} \GeoAng^{N-1} \mytr \angdeform{\GeoAng}) \cdot \angdiff \Psi
= \GeoAngFlatRadComponent \upmu (\angdiffuparg{\#} \Tanset^{N-1} \mytr \upchi) \cdot \angdiff \Psi
+ Harmless^{\leq N}
$,
which in particular yields the desired second product on 
RHS~\eqref{E:GEOANGANGISTHEFIRSTCOMMUTATORIMPORTANTTERMS}.
All remaining terms on
RHS~\eqref{E:GEOANGCANCELTERM2TERMMAINPOINTWISEESTIMATE} 
have $\leq N-2$ derivatives falling on 
$\angdiffuparg{\#} \mytr \angdeform{\GeoAng}$,
and the arguments given in 
our analysis of $\mathscr{K}_{(\pi-Cancel-1)}^{(\GeoAng)}[\Psi]$
yield that the corresponding products
$= Harmless^{\leq N}$ as desired.

\medskip
\noindent
\underline{\textbf{Analysis of} $\GeoAng^{N-1} \mathscr{K}_{(\pi-Good)}^{(\GeoAng)}[\Psi]$.}
From \eqref{E:DIVCURRENTGOOD}, we have
\begin{align} \label{E:GEOANGGOODTERMMAINPOINTWISEESTIMATE}
	\GeoAng^{N-1} \mathscr{K}_{(\pi-Good)}^{(\GeoAng)}[\Psi]
	& = \sum_{N_1 + N_2 = N-1}
			\left\lbrace
				\frac{1}{2} \GeoAng^{N_1} (\upmu \Lunit \mytr \angdeform{\Rad}) 
				+ \GeoAng^{N_1} \Lunit \deformarg{\Rad}{\Lunit}{\Rad}
				+ \GeoAng^{N_1} \Lunit \deformarg{\Rad}{\Rad}{\Radunit}
			\right\rbrace
			\GeoAng^{N_2} \Lunit \Psi
			\\
	& \ \ + \frac{1}{2} \sum_{N_1 + N_2 = N-1} (\GeoAng^{N_1} \Lunit \mytr \angdeform{Z}) \GeoAng^{N_2} \Rad \Psi
				\notag \\
	& \ \	- 
				\sum_{N_1 + N_2 = N-1}
					\left\lbrace
						\angLie_{\GeoAng}^{N_1} (\upmu \angLie_{\Lunit} \angdeformoneformupsharparg{Z}{\Lunit})
						+ 
						(\angLie_{\GeoAng}^{N_1} \angLie_{\Lunit} \angdeformoneformupsharparg{Z}{\Rad})
					\right\rbrace 
					\cdot \angdiff \GeoAng^{N_2} \Psi.
				\notag 
\end{align}
We claim that all terms on RHS~\eqref{E:GEOANGGOODTERMMAINPOINTWISEESTIMATE} 
$= Harmless^{\leq N}$ without the need to observe any cancellations.
The main point is that all deformation tensor components are hit
with an $\Lunit$ derivative and hence can be bounded with the estimate
\eqref{E:TANGENTIALDERIVATIVESOFDEFORMATIONINVOLVINGONELUNITDERIVATIVE}.
Otherwise, the analysis is essentially the same as our analysis of
$\GeoAng^{N-1} \mathscr{K}_{(\pi-Cancel-1)}^{(\GeoAng)}[\Psi]$.

\medskip
\noindent
\underline{\textbf{Analysis of} $\GeoAng^{N-1} \mathscr{K}_{(\Psi)}^{(\GeoAng)}[\Psi]$.}
We will show that these terms $= Harmless^{\leq N}$.
The terms in $\mathscr{K}_{(\Psi)}^{(\GeoAng)}[\Psi]$ (see \eqref{E:DIVCURRENTPSI}) 
are of the form
$
\smoothfunction(\BadVar) \pi P Z \Psi
+
\smoothfunction(\BadVar) \angLap \Psi
$
where $\Singletan \in \lbrace \angdiff \rbrace \cup \Tanset$, 
$Z \in \Fullset$,
and
$
\pi
\in
\left\lbrace
	\mytr \angdeform{\GeoAng},
	\deformarg{\GeoAng}{\Lunit}{\Rad},
	\deformarg{\GeoAng}{\Rad}{\Radunit},
	\angdeformoneformupsharparg{\GeoAng}{\Lunit},
	\angdeformoneformupsharparg{\GeoAng}{\Rad}
\right\rbrace
$.
We therefore conclude that
$
\left|
	\GeoAng^{N-1} \mathscr{K}_{(\Psi)}^{(\GeoAng)}[\Psi]
\right|
\lesssim
\left|
	\Fullset_*^{\leq N+1;1} \Psi
\right|
+
\left|
	\Tanset^{\leq N} \GdVar
\right|
+
\left|
	\Tanset_*^{[1,N]} \BadVar
\right|
= Harmless^{\leq N}
$
by using 
\eqref{E:TANGENTIALDERIVATIVESOFANGLAPPSIPOINTWISE},
\eqref{E:TANGENTIALDERIVATIVESOFDEFORMATION},
and the $L^{\infty}$ estimates of Prop.~\ref{P:IMPROVEMENTOFAUX}

\medskip
\noindent
\underline{\textbf{Analysis of} $\GeoAng^{N-1} \mathscr{K}_{(Low)}^{(\GeoAng)}[\Psi]$.}
We will show that these terms $= Harmless^{\leq N}$.
Using Lemma~\ref{L:SCHEMATICDEPENDENCEOFMANYTENSORFIELDS}, 
we see that (see \eqref{E:DIVCURRENTLOW}) 
$
\mathscr{K}_{(Low)}^{(\GeoAng)}[\Psi]
=
\smoothfunction(\Tanset^{\leq 1} \BadVar,\ginversesphere,\angdiff x^1,\angdiff x^2,\Rad \Psi) \pi P \GdVar
$
where 
$\pi$ and $\Singletan$ are as in the previous paragraph.
Hence, we conclude that $\GeoAng^{N-1} \mathscr{K}_{(Low)}^{(\GeoAng)}[\Psi] = Harmless^{\leq N}$
by using the same arguments as in the previous paragraph together with 
Lemmas~\ref{L:POINTWISEFORRECTANGULARCOMPONENTSOFVECTORFIELDS}
and
~\ref{L:POINTWISEESTIMATESFORGSPHEREANDITSDERIVATIVES}
(to bound the derivatives of $\ginversesphere$ and $\angdiff x$).

\medskip
Summing the above estimates and recalling the splitting \eqref{E:DIVCOMMUTATIONCURRENTDECOMPOSITION},
we conclude that the main term
$
\displaystyle
\GeoAng^{N-1}
					\left(
					\upmu 
					\D_{\alpha} 
					\Jcurrent{\GeoAng}^{\alpha}[\Psi]
					\right)
$
on RHS~\eqref{E:COMMUTEDWAVEGEOANGFIRSTMAINTERMPLUSERRORTERM}
is equal to RHS~\eqref{E:GEOANGANGISTHEFIRSTCOMMUTATORIMPORTANTTERMS} as desired.

To complete the proof of \eqref{E:GEOANGANGISTHEFIRSTCOMMUTATORIMPORTANTTERMS},
it remains only for us to show that
RHS~\eqref{E:GEOANGLOWERORDERRORTERMESTIMATE} $\lesssim Harmless^{\leq N}$.
The main point is that $N_2 \leq N-2$ in these terms and hence 
they do not involve the top-order derivatives
of $\upmu$ or $\Lunit_{(Small)}^1$, $\Lunit_{(Small)}^2$.
We first consider the case $N_1 \leq 9$ in inequality \eqref{E:GEOANGLOWERORDERRORTERMESTIMATE}.
Using the bound
$
\left\|
	\Tanset^{\leq 9} \mytr \deform{\GeoAng}
\right\|_{L^{\infty}(\Sigma_t^u)} \lesssim \varepsilon
$
mentioned just below \eqref{E:COMMCURRENTZ},
we see that when $N_1 \leq 9$, 
it suffices to bound
\begin{align}
\mathop{\sum_{N_2 + N_3 \leq N-1}}_{N_2 \leq N-2}
		\left|
			\GeoAng^{N_2}
			\left(
				\upmu 
				\D_{\alpha} 
				\Jcurrent{\GeoAng}^{\alpha}[\GeoAng^{N_3} \Psi]
			\right)
		\right|.
\end{align}
That is, we must bound the terms
$
\GeoAng^{N_2} \mathscr{K}_{(\pi-Danger)}^{(\GeoAng)}[\GeoAng^{N_3} \Psi]
+ \cdots
+ \GeoAng^{N_2} \mathscr{K}_{(Low)}^{(\GeoAng)}[\GeoAng^{N_3} \Psi]
$.
To this end, we repeat the proofs of the 
above estimates for
$
\GeoAng^{N-1} \mathscr{K}_{(\pi-Danger)}^{(\GeoAng)}[\Psi]
+ \cdots
+ \GeoAng^{N-1} \mathscr{K}_{(Low)}^{(\GeoAng)}[\Psi]
$
but with $N_2$ in place of $N-1$ and
$\GeoAng^{N_3} \Psi$ in place of
the explicitly written $\Psi$ terms,
the key point being that $N_2 \leq N-2$.
The same arguments immediately yield that all terms 
$= Harmless^{\leq  N_2 + N_3 + 1} \leq Harmless^{\leq N}$
except for products of the form
$
(\Rad \GeoAng^{N_3} \Psi) \GeoAng^{N_2+1} \mytr \upchi
$
and
$
\GeoAngFlatRadComponent (\angdiffuparg{\#} \GeoAng^{N_3} \Psi) \cdot (\upmu \angdiff \GeoAng^{N_2} \mytr \upchi)
$
corresponding to the 
two explicitly written products on RHS
\eqref{E:GEOANGANGISTHEFIRSTCOMMUTATORIMPORTANTTERMS}.
Since 
$N_2 \leq N-2$,
we can bound these two products by
$
\leq
Harmless^{\leq \max \lbrace N_3, N_2 + 2 \rbrace}
\leq 
Harmless^{\leq N}
$
with the help of
the relation $\GeoAngFlatRadComponent = \smoothfunction(\GdVar) \GdVar$ (see Lemma~\ref{L:SCHEMATICDEPENDENCEOFMANYTENSORFIELDS}),
\eqref{E:POINTWISEESTIMATESFORGSPHEREANDITSTANGENTIALDERIVATIVES},
and the $L^{\infty}$ estimates of Prop.~\ref{P:IMPROVEMENTOFAUX}.

To complete the proof of the desired bound for RHS~\eqref{E:GEOANGLOWERORDERRORTERMESTIMATE},
we must handle the case $N_1 \geq 10$ on RHS~\eqref{E:GEOANGLOWERORDERRORTERMESTIMATE}
(and thus $N_2 + N_3 \leq 7$).
The arguments given in the previous paragraph yield that
$
\displaystyle
\GeoAng^{N_2}
			\left(
				\upmu 
				\D_{\alpha} 
				\Jcurrent{\GeoAng}^{\alpha}[\GeoAng^{N_3} \Psi]
			\right)
= Harmless^{\leq \max \lbrace N_2 + N_3 + 1, N_2 + 2 \rbrace}
\leq Harmless^{\leq 9}
$.
Using this estimate and Lemma~\ref{L:LINFINITYESTIMATEFORHARMLESSLEQ9}, we deduce that
$
\displaystyle
\left|
	\GeoAng^{N_2}
			\left(
				\upmu 
				\D_{\alpha} 
				\Jcurrent{\GeoAng}^{\alpha}[\GeoAng^{N_3} \Psi]
			\right)
\right|
\lesssim \varepsilon
$.
From this estimate, we deduce that
the terms on RHS~\eqref{E:GEOANGLOWERORDERRORTERMESTIMATE}
with $N_1 \geq 10$ are
$
\displaystyle
	\lesssim
	\left|
		\GeoAng^{\leq N-1} \mytr \angdeform{\GeoAng}
	\right|
$.
Finally, 
\eqref{E:CONNECTIONBETWEENANGLIEOFGSPHEREANDDEFORMATIONTENSORS}, 
\eqref{E:POINTWISEESTIMATESFORGSPHEREANDITSTANGENTIALDERIVATIVES},
and the $L^{\infty}$ estimates of Prop.~\ref{P:IMPROVEMENTOFAUX})
together yield that
$
\displaystyle
	\left|
		\GeoAng^{\leq N-1} \mytr \angdeform{\GeoAng}
	\right|
\lesssim 
	\left|
		\Tanset^{\leq N} \GdVar
	\right|
= Harmless^{\leq N}
$ as desired.
We have thus proved \eqref{E:GEOANGANGISTHEFIRSTCOMMUTATORIMPORTANTTERMS}.

The proof of \eqref{E:LISTHEFIRSTCOMMUTATORIMPORTANTTERMS} is essentially the same
with a few minor differences that we now mention.
The term $\GeoAng^{N-1} \mathscr{K}_{(\pi-Danger)}^{(\Lunit)}[\Psi]$
(see \eqref{E:DIVCURRENTTRANSVERSAL}) is actually trivial in this case because 
$\angdeformoneformupsharparg{\Lunit}{\Lunit} = 0$
(see \eqref{E:LUNITDEFORMSPHERELUNITANDLUNITDEFORMSPHERERAD}).
We again observe cancellation of the top-order eikonal function
quantities in
$\GeoAng^{N-1} \mathscr{K}_{(\pi-Cancel-1)}^{(\Lunit)}[\Psi]$
(see \eqref{E:DIVCURRENTCANEL1})
up to $Harmless^{\leq N}$ errors.
Specifically, with the help of
\eqref{E:LDEFIMPORTANTRADTRACEANGTERMS}-\eqref{E:LDEFIMPORTANTLIERADSPHERERADANDAGNDIVSPHERETERMS}
and the fact that $\angdeformoneformupsharparg{\Lunit}{\Lunit} = 0$,
we observe cancellation of
$\angLap \GeoAng^{N-1} \upmu$.
In contrast, 
without the need to observe any cancellations,
all terms in $\GeoAng^{N-1} \mathscr{K}_{(\pi-Cancel-2)}^{(\Lunit)}[\Psi]$
(see \eqref{E:DIVCURRENTCANEL2})
$= Harmless^{\leq N}$,
thanks to the estimate \eqref{E:LDEFIMPORTANTLIERADSPHERERADANDAGNDIVSPHERETERMS} for
$\angLie_{\Tanset}^{N-1} \angdiffuparg{\#} \deformarg{\Lunit}{\Lunit}{\Rad}$
and the fact that $\angdeformoneformupsharparg{\Lunit}{\Lunit} = 0$.
The main term on RHS~\eqref{E:LISTHEFIRSTCOMMUTATORIMPORTANTTERMS}
comes from the case when all $N-1$ derivatives in the term
$\GeoAng^{N-1} \mathscr{K}_{(\pi-Less \ Dangerous)}^{(\Lunit)}[\Psi]$
fall on the factor $\mytr \angdeform{\Lunit}$ from
\eqref{E:DIVCURRENTELLIPTIC},
where we substitute RHS~\eqref{E:LUNITDEFORMSPHERE}
for $\angdeform{\Lunit}$.
All other terms on RHS
\eqref{E:LISTHEFIRSTCOMMUTATORIMPORTANTTERMS}
are $Harmless^{\leq N}$, as in the proof of
\eqref{E:GEOANGANGISTHEFIRSTCOMMUTATORIMPORTANTTERMS}.

To prove \eqref{E:HARMLESSORDERNCOMMUTATORS}, we first note that
$\Tanset^N$ must be either of the form
$\Tanset^{N-1} \Lunit$
or $\Tanset^{N-1} \GeoAng$,
where $\Tanset^{N-1}$ \emph{contains a factor of} $\Lunit$.
In the former case, by using essentially the same 
arguments we used in the proof of \eqref{E:LISTHEFIRSTCOMMUTATORIMPORTANTTERMS}
but with $\Tanset^{N-1} \Lunit$ in the role of $\GeoAng^{N-1} \Lunit$,
we deduce that
\begin{align} \label{E:TWOLSWITHLFIRSTBOXPSICOMMUTATION}
	\square_g (\Tanset^{N-1} \Lunit \Psi)
	& = (\angdiffuparg{\#} \Psi) \cdot (\upmu \angdiff \Tanset^{N-1} \mytr \upchi)
			+ Harmless^{\leq N}.
\end{align}
The $L^{\infty}$ estimates of Prop.~\ref{P:IMPROVEMENTOFAUX}
imply that the first product on RHS~\eqref{E:TWOLSWITHLFIRSTBOXPSICOMMUTATION} is
$
\lesssim
	\left|
		\Tanset^N \mytr \upchi
	\right|
$,
where $\Tanset^N$ contains a factor of $\Lunit$.
We now commute the operator $\Lunit$ to the front
and use
\eqref{E:PURETANGENTIALFUNCTIONCOMMUTATORESTIMATE}
with $f = \mytr \upchi$,
\eqref{E:POINTWISEESTIMATESFORGSPHEREANDITSTANGENTIALDERIVATIVES},
and the $L^{\infty}$ estimates of Prop.~\ref{P:IMPROVEMENTOFAUX}
to deduce that the commutator error terms
$= Harmless^{\leq N}$.
We then use \eqref{E:LUNITTANGENTDIFFERENTIATEDLUNITSMALLIMPROVEDPOINTWISE}
to bound the non-commutator term as follows:
$
\displaystyle
\left|
	\Lunit \Tanset^{N-1} \mytr \upchi
\right|
\lesssim
\left|
	\Tanset^{\leq N+1} \Psi
\right|
+ 
\left|
	\Tanset^{\leq N} \GdVar
\right|
=
Harmless^{\leq N}
$.
We have thus proved \eqref{E:HARMLESSORDERNCOMMUTATORS} in this case.

In the remaining case of \eqref{E:HARMLESSORDERNCOMMUTATORS},
in which $\Tanset^N$ is of the form
$\Tanset^{N-1} \GeoAng$
and $\Tanset^{N-1}$ \emph{contains a factor of} $\Lunit$,
we use essentially the same 
arguments we used in the proof of \eqref{E:GEOANGANGISTHEFIRSTCOMMUTATORIMPORTANTTERMS}
but with $\Tanset^{N-1} \GeoAng$ in the role of $\GeoAng^N$
to deduce
\begin{align} \label{E:GEOANGFIRSTWITHAFACTOROFLBOXPSICOMMUTATION}
	\square_g (\Tanset^{N-1} \GeoAng \Psi)
	& = (\Rad \Psi) \Tanset^{N-1} \GeoAng \mytr \upchi
			+ \GeoAngFlatRadComponent (\angdiffuparg{\#} \Psi) \cdot (\upmu \angdiff \Tanset^{N-1} \mytr \upchi)
			+ Harmless^{\leq N},
\end{align}
where the operators $\Tanset^{N-1}$ in \eqref{E:GEOANGFIRSTWITHAFACTOROFLBOXPSICOMMUTATION} contain a factor of $\Lunit$.
The $L^{\infty}$ estimates of Prop.~\ref{P:IMPROVEMENTOFAUX} 
imply that the first product on RHS~\eqref{E:GEOANGFIRSTWITHAFACTOROFLBOXPSICOMMUTATION} is
$
\lesssim
	\left|
		\Tanset^N \mytr \upchi
	\right|
$,
where $\Tanset^N$ contains a factor of $\Lunit$.
Therefore, the arguments from the previous paragraph yield that this term
$= Harmless^{\leq N}$ as desired.
Also using that $\GeoAngFlatRadComponent = \smoothfunction(\GdVar) \GdVar$
(see \eqref{E:LINEARLYSMALLSCALARSDEPENDINGONGOODVARIABLES}),
we find that
the second product on RHS~\eqref{E:GEOANGFIRSTWITHAFACTOROFLBOXPSICOMMUTATION} is
$
\lesssim
	\left|
		\Tanset^N \mytr \upchi
	\right|
$,
where $\Tanset^N$ contains a factor of $\Lunit$.
The arguments from the previous paragraph yield that this term
$= Harmless^{\leq N}$ as desired.
We have thus proved \eqref{E:HARMLESSORDERNCOMMUTATORS}
and completed the proof of Prop.~\ref{P:IDOFKEYDIFFICULTENREGYERRORTERMS}.

\subsection{Pointwise estimates for the fully modified quantities}
\label{SS:POINTWISEFULLYMODIFIED}
In this subsection, we obtain pointwise estimates 
for the most difficult product we encounter in our energy
estimates: $(\Rad \Psi) \GeoAng^N \mytr \upchi$.
The main result is Prop.~\ref{P:KEYPOINTWISEESTIMATE}.
The proof of the proposition relies 
on pointwise estimates for the fully modified quantities,
which we first derive. We start with a simple lemma in which we obtain pointwise estimates
for some of the inhomogeneous terms in the transport equations
verified by the fully modified and partially modified quantities.

\begin{lemma}[\textbf{Pointwise estimates for} $\Tanset^N \upchifullmodinhom$ \textbf{and} $\Tanset^N \upchipartialmodinhom$]
	\label{L:CHIPARTIALMODSOURCETERMPOINTWISE}
	Assume that $N \leq 18$. Let 
	$\upchifullmodinhom$ be the quantity defined in 
	\eqref{E:LOWESTORDERTRANSPORTRENORMALIZEDTRCHIJUNKDISCREPANCY},
	let $\upchipartialmodinhom$
	be the quantity defined in \eqref{E:LOWESTORDERTRANSPORTPARTIALRENORMALIZEDTRCHIJUNKDISCREPANCY},
	let $\upchipartialmodinhomarg{\GeoAng^{N-1}}$
	be the quantity from
	\eqref{E:TRANSPORTPARTIALRENORMALIZEDTRCHIJUNKDISCREPANCY}
	(with $\GeoAng^{N-1}$ in the role of $\Tanset^N$),
	and let
	$
	\upchipartialmodsourcearg{\Tanset^{N-1}}
	$
	be the quantity defined in
	\eqref{E:TRCHIJUNKCOMMUTEDTRANSPORTEQNPARTIALRENORMALIZATIONINHOMOGENEOUSTERM}.
	Under the data-size and bootstrap assumptions 
	of Subsects.~\ref{SS:SIZEOFTBOOT}-\ref{SS:PSIBOOTSTRAP}
	and the smallness assumptions of Subsect.~\ref{SS:SMALLNESSASSUMPTIONS}, 
	the following pointwise estimates hold
	on $\mathcal{M}_{\Tboot,U_0}$: 
\begin{subequations}
	\begin{align} 
		\left|
			\GeoAng^N \upchifullmodinhom
			+ 
			G_{\Lunit \Lunit}
			\Rad \GeoAng^N \Psi
		\right|
		& 
		\lesssim
		\upmu
		\left|
			\Tanset^{\leq N+1} \Psi
		\right|
		+
		\left|
			\Fullset_*^{\leq N;1} \Psi
		\right|
		+
		\left|
			\Tanset^{\leq N} \GdVar
		\right|
		+
		\left|
			\Tanset_*^{[1,N]} \BadVar
		\right|,
			\label{E:TOPDERIVATIVESOFXPOINTWISEBOUND} \\
		\left|
			\Tanset^{\leq N} \upchifullmodinhom
		\right|
		&
		\lesssim
		\left|
			\Fullset_*^{\leq N+1;1} \Psi
		\right|
		+
		\left|
			\Tanset^{\leq N} \GdVar
		\right|
		+ 
		\left|
			\Tanset_*^{[1,N]} \BadVar
		\right|,
			\label{E:CRUDEXPOINTWISEBOUND} \\
		\left|
			\Tanset^N \upchipartialmodinhom
		\right|
		& \lesssim
			\left|
				\Tanset^{\leq N+1} \Psi
			\right|
			+
			\left|
				\Tanset^{\leq N} \GdVar
			\right|,
				\label{E:POINTWISELOWESTORDERTRANSPORTPARTIALRENORMALIZEDTRCHIJUNKDISCREPANCY} \\
		\left|
			\GeoAng \upchipartialmodinhomarg{\GeoAng^{N-1}}
		\right|
		& 
		\lesssim
			\left| \Tanset^{\leq N+1} \Psi \right|,
			\label{E:HARMLESSNATUREOFPARTIALLYMODIFIEDDISCREPANCY}
			\\
		\left|
			\upchipartialmodsourcearg{\Tanset^{N-1}}
		\right|
		& \lesssim 
			\varepsilon
			\left|
				\Tanset^{\leq N} \GdVar
			\right|.
			\label{E:CHIPARTIALMODSOURCETERMPOINTWISE}
	\end{align}
	\end{subequations}

\end{lemma}

\begin{proof}
	See Subsect.~\ref{SS:OFTENUSEDESTIMATES} for some comments on the analysis.
	Throughout this proof, we silently use the $L^{\infty}$ estimates 
	of Prop.~\ref{P:IMPROVEMENTOFAUX}.

To prove \eqref{E:TOPDERIVATIVESOFXPOINTWISEBOUND},
we first use \eqref{E:LOWESTORDERTRANSPORTRENORMALIZEDTRCHIJUNKDISCREPANCY}
and Lemma~\ref{L:SCHEMATICDEPENDENCEOFMANYTENSORFIELDS}
to obtain
$
\upchifullmodinhom 
	= 
	-
	 G_{\Lunit \Lunit} \Rad \Psi
	+
	\upmu \smoothfunction(\GdVar,\ginversesphere,\angdiff x^1,\angdiff x^2) \Singletan \Psi
$.
We now apply $\GeoAng^N$ to this identity 
and bring the top-order term 
$G_{\Lunit \Lunit} \Rad \GeoAng^N \Psi$
over to the left 
(as indicated on LHS~\eqref{E:TOPDERIVATIVESOFXPOINTWISEBOUND}), 
which leaves the commutator terms
$[G_{\Lunit \Lunit},\GeoAng^N] \Rad \Psi$
and
$G_{\Lunit \Lunit} [\Rad,\GeoAng^N] \Psi$
on the RHS.
To bound 
$
\left|
\GeoAng^N 
\left\lbrace 
	\upmu \smoothfunction(\GdVar,\ginversesphere,\angdiff x^1,\angdiff x^2) \Singletan \Psi
\right\rbrace
\right|
$
by $\leq$ RHS~\eqref{E:TOPDERIVATIVESOFXPOINTWISEBOUND}, 
we use
Lemmas~\ref{L:POINTWISEFORRECTANGULARCOMPONENTSOFVECTORFIELDS}
and \ref{L:POINTWISEESTIMATESFORGSPHEREANDITSDERIVATIVES}.
Note that we have paid special attention to terms in which
all derivatives $\GeoAng^N$ fall on $\Singletan \Psi$;
these terms are bounded by 
the first term on RHS~\eqref{E:TOPDERIVATIVESOFXPOINTWISEBOUND}.
To bound 
$
\left|
	[G_{\Lunit \Lunit},\GeoAng^N] \Rad \Psi
\right|
$
by $\leq$ RHS~\eqref{E:TOPDERIVATIVESOFXPOINTWISEBOUND},
we use the fact that 
$G_{\Lunit \Lunit} = \smoothfunction(\GdVar)$ (see Lemma~\ref{L:SCHEMATICDEPENDENCEOFMANYTENSORFIELDS}).
To bound
$
\left|
	G_{\Lunit \Lunit} [\Rad,\GeoAng^N] \Psi
\right|
$
by $\leq$ RHS~\eqref{E:TOPDERIVATIVESOFXPOINTWISEBOUND}, 
we also use the commutator estimate \eqref{E:ONERADIALTANGENTIALFUNCTIONCOMMUTATORESTIMATE}
with $f = \Psi$. The proof of \eqref{E:CRUDEXPOINTWISEBOUND}
is similar but simpler and we omit the details.
The same is true for the proof of
\eqref{E:POINTWISELOWESTORDERTRANSPORTPARTIALRENORMALIZEDTRCHIJUNKDISCREPANCY}
since by Lemma~\ref{L:SCHEMATICDEPENDENCEOFMANYTENSORFIELDS},
we have
$
\upchipartialmodinhom 
= 
\smoothfunction(\GdVar,\ginversesphere,\angdiff x^1,\angdiff x^2)
\Singletan \Psi
$.

	To derive \eqref{E:HARMLESSNATUREOFPARTIALLYMODIFIEDDISCREPANCY}, 
	we first use \eqref{E:TRANSPORTPARTIALRENORMALIZEDTRCHIJUNKDISCREPANCY} 
	and Lemma~\ref{L:SCHEMATICDEPENDENCEOFMANYTENSORFIELDS}
	to deduce that
	$
	\GeoAng \upchipartialmodinhomarg{\GeoAng^{N-1}} 
	= 
	\GeoAng 
	\left\lbrace 
		\smoothfunction(\GdVar,\ginversesphere,\angdiff x^1,\angdiff x^2) 
		\Tanset^N \Psi 
	\right\rbrace$.
	The estimate \eqref{E:HARMLESSNATUREOFPARTIALLYMODIFIEDDISCREPANCY} now follows easily 
	from the previous expression and
	Lemmas~\ref{L:POINTWISEFORRECTANGULARCOMPONENTSOFVECTORFIELDS}
  and \ref{L:POINTWISEESTIMATESFORGSPHEREANDITSDERIVATIVES}.

	We now prove \eqref{E:CHIPARTIALMODSOURCETERMPOINTWISE}.
	We bound term
	$
	\displaystyle
	\Tanset^{N-1} (\mytr \upchi)^2
	$
	from RHS~\eqref{E:TRCHIJUNKCOMMUTEDTRANSPORTEQNPARTIALRENORMALIZATIONINHOMOGENEOUSTERM}
	by $\leq$ RHS~\eqref{E:CHIPARTIALMODSOURCETERMPOINTWISE}
	with the help of inequality \eqref{E:POINTWISEESTIMATESFORGSPHEREANDITSTANGENTIALDERIVATIVES}.
	We bound the term
	$
	\displaystyle
	[G_{\Lunit \Lunit}, \Tanset^{N-1}] \angLap \Psi
	$
	using the aforementioned relation $G_{\Lunit \Lunit} = \smoothfunction(\GdVar)$ and 
	Cor.~\ref{C:TANGENTIALDERIVATIVESOFANGLAPPSIPOINTWISE}.
	To bound the term
	$
	\displaystyle
	G_{\Lunit \Lunit} [\angLap, \Tanset^{N-1}] \Psi
	$,
	we also use
	the commutator estimate \eqref{E:ANGLAPPURETANGENTIALFUNCTIONCOMMUTATOR}
	with $\Psi$ in the role of $f$
	and Cor.~\ref{C:SQRTEPSILONTOCEPSILON}.
	We bound the term
	$
	\displaystyle
	[\Lunit, \Tanset^{N-1}] \mytr \upchi
	$
	with the help of the commutator estimate
	\eqref{E:PURETANGENTIALFUNCTIONCOMMUTATORESTIMATE}
	with $\mytr \upchi$ in the role of $f$
	and inequality \eqref{E:POINTWISEESTIMATESFORGSPHEREANDITSTANGENTIALDERIVATIVES}.
	We bound 
	$
	\displaystyle
	[\Lunit, \Tanset^{N-1}] \upchipartialmodinhom
	$
	with the help of the commutator estimate \eqref{E:PURETANGENTIALFUNCTIONCOMMUTATORESTIMATE} with $f =\upchipartialmodinhom$
	and \eqref{E:POINTWISELOWESTORDERTRANSPORTPARTIALRENORMALIZEDTRCHIJUNKDISCREPANCY}.
	To bound
	$
	\displaystyle
	\Lunit
				\left\lbrace
					\upchipartialmodinhomarg{\Tanset^{N-1}}
					- \Tanset^{N-1} \upchipartialmodinhom
				\right\rbrace
	$,
	we first note that
	\eqref{E:TRANSPORTPARTIALRENORMALIZEDTRCHIJUNKDISCREPANCY},
	\eqref{E:LOWESTORDERTRANSPORTPARTIALRENORMALIZEDTRCHIJUNKDISCREPANCY},
	and the Leibniz rule imply that the magnitude of this term is
	$
	\displaystyle
	\lesssim
	\mathop{\sum_{N_1 + N_2 \leq N}}_{N_1 \geq 1}
	\left|
		\angLie_{\Tanset}^{N_1} G_{(Frame)}^{\#}
	\right|
	\left|
		\Tanset^{N_2+1} \Psi
	\right|
	+
	\left|
		G_{(Frame)}^{\#}
	\right|
	\left|
		[\Lunit,\Lunit \Tanset^{N-1}] \Psi
	\right|
	$.
	Since Lemma~\ref{L:SCHEMATICDEPENDENCEOFMANYTENSORFIELDS} implies that 
	$G_{(Frame)}^{\#} = \smoothfunction(\GdVar,\ginversesphere,\angdiff x^1,\angdiff x^2)$,
	the desired bound for the sum follows from
	Lemmas~\ref{L:POINTWISEFORRECTANGULARCOMPONENTSOFVECTORFIELDS}
  and \ref{L:POINTWISEESTIMATESFORGSPHEREANDITSDERIVATIVES}.
	To bound the term
	$\left|
		G_{(Frame)}^{\#}
	\right|
	\left|
		[\Lunit,\Lunit \Tanset^{N-1}] \Psi
	\right|
	$
	by $\leq$ RHS~\eqref{E:CHIPARTIALMODSOURCETERMPOINTWISE}, 
	we also use the commutator estimate \eqref{E:PURETANGENTIALFUNCTIONCOMMUTATORESTIMATE}
	with $f = \Psi$
	and Cor.~\ref{C:SQRTEPSILONTOCEPSILON}.
	This completes the proof of \eqref{E:CHIPARTIALMODSOURCETERMPOINTWISE} 
	and finishes the proof of the lemma.

\end{proof}

Recall that the fully modified quantities $\upchifullmodarg{\Tanset^N}$
verify the transport equation \eqref{E:TOPORDERTRCHIJUNKRENORMALIZEDTRANSPORT}.
In the next lemma, we integrate this transport equation and derive pointwise estimates
for $\upchifullmodarg{\Tanset^N}$. The lemma is a preliminary ingredient
in the proof of Prop.~\ref{P:KEYPOINTWISEESTIMATE}.

\begin{lemma}[\textbf{Estimates for solutions to the transport equation verified by} $\upchifullmodarg{\Tanset^N}$]
	\label{L:RENORMALIZEDTOPORDERTRCHIJUNKTRANSPORTINVERTED}
	Assume that $1 \leq N \leq 18$ and 
	let $\upchifullmodarg{\Tanset^N}$ 
	and
	$\upchifullmodinhom$
	be as in Prop.~\ref{P:TOPORDERTRCHIJUNKRENORMALIZEDTRANSPORT}.
	Assume that $\Tanset^N = \GeoAng^N$.
	Under the data-size and bootstrap assumptions 
	of Subsects.~\ref{SS:SIZEOFTBOOT}-\ref{SS:PSIBOOTSTRAP}
	and the smallness assumptions of Subsect.~\ref{SS:SMALLNESSASSUMPTIONS}, 
	the following pointwise estimate holds on
	$\mathcal{M}_{\Tboot,U_0}$:
	\begin{align} \label{E:RENORMALIZEDTOPORDERTRCHIJUNKTRANSPORTINVERTED}
		\left|\upchifullmodarg{\GeoAng^N} \right|(t,u,\vartheta)
		& \leq 
			C
			\left| \upchifullmodarg{\GeoAng^N } \right|(0,u,\vartheta)
				\\
		& \ \ + 	\boxed{2} (1 + C \varepsilon)
							\int_{s=0}^t 
								\frac{[\Lunit \upmu(s,u,\vartheta)]_-}{\upmu(s,u,\vartheta)}
								\left| \GeoAng^N  \upchifullmodinhom \right|(s,u,\vartheta)
							\, ds
					\notag \\
	& \ \ + 	C
						\int_{s=0}^t 
							\left\lbrace
								\left|
									\Fullset_*^{\leq N+1;1} \Psi
								\right|
								+
								\left|
									\Tanset^{\leq N} \GdVar
								\right|
								+
								\left|
									\Tanset_*^{[1,N]} \BadVar
								\right|
							\right\rbrace
							(s,u,\vartheta)
						\, ds.
						\notag 
\end{align}

\end{lemma}

\begin{proof}
		To prove \eqref{E:RENORMALIZEDTOPORDERTRCHIJUNKTRANSPORTINVERTED}, 
		we set $\Tanset^N = \GeoAng^N$
		in equation
		\eqref{E:TOPORDERTRCHIJUNKRENORMALIZEDTRANSPORT}
		and, in this part of the proof,
		we view the terms in the equation as functions of $(s,u,\vartheta)$.
		Corresponding to the factor
		$
		\displaystyle
		\left(
			- 2 \frac{\Lunit \upmu}{\upmu}
			+ 2 \mytr \upchi
		\right)
		$
		on the left-hand side, 
		we define the integrating factor
		\begin{align}
			\iota(s,u,\vartheta)
				:= 
				\exp
				\left\lbrace
				\int_{t'=0}^s
					\left(
						- 2 \frac{\Lunit \upmu}{\upmu}(t',u,\vartheta)
						+ 2 \mytr \upchi
					\right)
					(t',u,\vartheta)
				\, dt'
				\right\rbrace.
		\end{align}
		We then rewrite
		\eqref{E:TOPORDERTRCHIJUNKRENORMALIZEDTRANSPORT}
		as 
		$
		\displaystyle
		\Lunit
		(\iota \upchifullmodarg{\Tanset^N})
		= \iota \times \mbox{RHS } \eqref{E:TOPORDERTRCHIJUNKRENORMALIZEDTRANSPORT}
		$
		and integrate this equation with respect to $s$
		from $s=0$ to $s=t$.
		Using the estimate \eqref{E:PURETANGENTIALCHICOMMUTEDLINFINITY}
		for $\mytr \upchi$, we find that
		\begin{align} \label{E:KEYINTEGRATINFACTORESTIMATE}
			\iota(s,u,\vartheta)
			& = (1 + \mathcal{O}(\varepsilon))
					\frac{\upmu^2(0,u,\vartheta)}{\upmu^2(s,u,\vartheta)}.
		\end{align}
		From Def.~\ref{D:REGIONSOFDISTINCTUPMUBEHAVIOR}
		and the estimates
		\eqref{E:LOCALIZEDMUCANTGROWTOOFAST}
		and
		\eqref{E:LOCALIZEDMUMUSTSHRINK},
		we find that
		\begin{align} \label{E:CRUDEMUOVERMUBOUND}
			\sup_{0 \leq s' \leq t}
				\frac{\upmu(t,u,\vartheta)}{\upmu(s',u,\vartheta)}
				& \leq C.
		\end{align}
		From \eqref{E:KEYINTEGRATINFACTORESTIMATE}
		and \eqref{E:CRUDEMUOVERMUBOUND},
		it is straightforward to see that the desired bound
		\eqref{E:RENORMALIZEDTOPORDERTRCHIJUNKTRANSPORTINVERTED}
		follows once we establish the following bounds for terms
		generated by the terms on RHS~\eqref{E:TOPORDERTRCHIJUNKRENORMALIZEDTRANSPORT}
		(recall that $\Tanset^N = \GeoAng^N$):
		\begin{align}
			& \left|
					\upmu [\Lunit, \GeoAng^N] \mytr \upchi
				\right|
			(s,u,\vartheta)
				\label{E:SECONDRERNORMALIZEDTRANSPORTEQNINHOMOGENEOUSTERMBOUND}  \\
			& \leq
				C \varepsilon
			\left|
				\upchifullmodarg{\GeoAng^N} 
			\right|
			(s,u,\vartheta)
			+
			C \varepsilon
			\left|
				\Fullset_*^{\leq N+1;1} \Psi
			\right|
			(s,u,\vartheta)
			+
			C \varepsilon
			\left|
				\Tanset^{\leq N} \GdVar
			\right|
			(s,u,\vartheta)
			+
			C 
			\varepsilon
			\left|
				\Tanset_*^{[1,N]} \BadVar
			\right|
			(s,u,\vartheta),
				\notag \\
			& 2 
			\left(
				\frac{\upmu(t,u,\vartheta)}{\upmu(s,u,\vartheta)}
			\right)^2 
			\left|
				\frac{\Lunit \upmu(s,u,\vartheta)}{\upmu(s,u,\vartheta)}
			\right|
			\left|
				\GeoAng^N \upchifullmodinhom
			\right|
			(s,u,\vartheta)
				\label{E:FIRSTRERNORMALIZEDTRANSPORTEQNINHOMOGENEOUSTERMBOUND} \\
			& \leq
				2 (1 + C \varepsilon)
				\left|
					\frac{[\Lunit \upmu]_-(s,u,\vartheta)}{\upmu(s,u,\vartheta)}
				\right|
				\left|
					\GeoAng^N \upchifullmodinhom
				\right|
				(s,u,\vartheta)
				\notag \\
			& \ \
				+ 
				C
				\left|
					\Fullset_*^{\leq N+1;1} \Psi
				\right|
				(s,u,\vartheta)
				+
				C
				\left|
					\Tanset^{\leq N} \GdVar
				\right|
				(s,u,\vartheta)
				+
				C
				\left|
					\Tanset_*^{[1,N]} \BadVar
	  		\right|
	  		(s,u,\vartheta),
			 \notag	\\
			& \mbox{all remaining terms on RHS } \eqref{E:TOPORDERTRCHIJUNKRENORMALIZEDTRANSPORT}
				\mbox{ are in magnitude}
				\label{E:THIRDRERNORMALIZEDTRANSPORTEQNINHOMOGENEOUSTERMBOUND} \\
			& \leq 
			C 
			\left|
				\Fullset_*^{\leq N+1;1} \Psi
			\right|
			(s,u,\vartheta)
			+
			C
			\left|
				\Tanset^{\leq N} \GdVar
			\right|
			(s,u,\vartheta)
			+
			C 
			\left|
				\Tanset_*^{[1,N]} \BadVar
			\right|
			(s,u,\vartheta).
			\notag
		\end{align}
		We note that in deriving \eqref{E:RENORMALIZEDTOPORDERTRCHIJUNKTRANSPORTINVERTED},
		the product 
		$ C 
			\varepsilon
			\iota 
			\left|
				\upchifullmodarg{\GeoAng^N} 
			\right|
		$
		arising from the first term on RHS~\eqref{E:SECONDRERNORMALIZEDTRANSPORTEQNINHOMOGENEOUSTERMBOUND}
		needs to be treated with Gronwall's inequality.
		However, due to the small factor $\varepsilon$,
		this product has only the negligible effect of 
		contributing to the factors of
		$C \varepsilon$
		on RHS~\eqref{E:RENORMALIZEDTOPORDERTRCHIJUNKTRANSPORTINVERTED}.

		To derive \eqref{E:FIRSTRERNORMALIZEDTRANSPORTEQNINHOMOGENEOUSTERMBOUND}, 
		we first note the trivial bound
		$
		\displaystyle
		\left|
			\frac{\Lunit \upmu(s,u,\vartheta)}{\upmu(s,u,\vartheta)}
		\right|
		\leq
		\left|
			\frac{[\Lunit \upmu]_-(s,u,\vartheta)}{\upmu(s,u,\vartheta)}
		\right|
		+
		\left|
			\frac{[\Lunit \upmu]_+(s,u,\vartheta)}{\upmu(s,u,\vartheta)}
		\right|
		$.
		To bound the terms on LHS~\eqref{E:FIRSTRERNORMALIZEDTRANSPORTEQNINHOMOGENEOUSTERMBOUND} arising from the factor
		$
		\displaystyle
		\left|
			\frac{[\Lunit \upmu]_+(s,u,\vartheta)}{\upmu(s,u,\vartheta)}
		\right|
		$
		by $\leq$ the terms on the last line of RHS~\eqref{E:FIRSTRERNORMALIZEDTRANSPORTEQNINHOMOGENEOUSTERMBOUND},
		we use \eqref{E:POSITIVEPARTOFLMUOVERMUISBOUNDED},
		\eqref{E:CRUDEMUOVERMUBOUND},
		and the estimate \eqref{E:CRUDEXPOINTWISEBOUND}.
To bound the terms on LHS~\eqref{E:FIRSTRERNORMALIZEDTRANSPORTEQNINHOMOGENEOUSTERMBOUND} arising from the factor
$
\displaystyle
\left|
	\frac{[\Lunit \upmu]_-(s,u,\vartheta)}{\upmu(s,u,\vartheta)}
\right|
$,
we consider the partitions from Def.~\ref{D:REGIONSOFDISTINCTUPMUBEHAVIOR}.
When $(u,\vartheta) \in \Vplus{t}{u}$,
we use the bounds 
\eqref{E:KEYMUNOTDECAYINGMINUSPARTLMUOVERMUBOUND}
and
\eqref{E:CRUDEMUOVERMUBOUND}
to deduce that
$
\displaystyle
\left(
	\frac{\upmu(t,u,\vartheta)}{\upmu(s,u,\vartheta)}
\right)^2
\left|
	\frac{[\Lunit \upmu]_-(s,u,\vartheta)}{\upmu(s,u,\vartheta)}
\right|
\leq C \varepsilon
$.
Combining this bound with \eqref{E:CRUDEXPOINTWISEBOUND}, 
we easily conclude that the terms of interest
are $\leq$ the terms on the last line of
RHS~\eqref{E:FIRSTRERNORMALIZEDTRANSPORTEQNINHOMOGENEOUSTERMBOUND}.
Finally, when $(u,\vartheta) \in \Vminus{t}{u}$,
we use \eqref{E:LOCALIZEDMUMUSTSHRINK}
to deduce that
\[
2
\left(
	\frac{\upmu(t,u,\vartheta)}{\upmu(s,u,\vartheta)}
\right)^2
\left|
	\frac{[\Lunit \upmu]_-(s,u,\vartheta)}{\upmu(s,u,\vartheta)}
\right|
\leq 2(1 + C \varepsilon) 
\left|
	\frac{[\Lunit \upmu]_-(s,u,\vartheta)}{\upmu(s,u,\vartheta)}
\right|.
\]
Thus, we conclude that the terms under consideration are
$\leq$ the terms on the first line of
RHS~\eqref{E:FIRSTRERNORMALIZEDTRANSPORTEQNINHOMOGENEOUSTERMBOUND}
as desired.

		To deduce \eqref{E:SECONDRERNORMALIZEDTRANSPORTEQNINHOMOGENEOUSTERMBOUND},
		we use definition \eqref{E:TRANSPORTRENORMALIZEDTRCHIJUNK},
		use the commutator estimate
		\eqref{E:PURETANGENTIALFUNCTIONCOMMUTATORESTIMATE} with $f = \mytr \upchi$,
		the estimates \eqref{E:PURETANGENTIALCHICOMMUTEDLINFINITY},
		\eqref{E:POINTWISEESTIMATESFORGSPHEREANDITSTANGENTIALDERIVATIVES},
		and \eqref{E:UPMULINFTY},
		and Cor.~\ref{C:SQRTEPSILONTOCEPSILON}
		to deduce that
		\begin{align} \label{E:ANNOYINGTRANSPORTINHOMOGENEOUSCOMMUTATORTERM}
		\upmu [\Lunit, \GeoAng^N] \mytr \upchi
		&
		\lesssim
			\varepsilon
			\upmu 
			\left|
				\GeoAng^N \mytr \upchi
			\right|
			+ 
			\varepsilon
			\left|
				\Tanset^{\leq N} \GdVar
			\right|
				\\
		& \lesssim
			\varepsilon
			\left|\upchifullmodarg{\GeoAng^N} \right|
			+ 
			\varepsilon
			\left|
				\Tanset^{\leq N} \GdVar
			\right|
			+
			\varepsilon
			\left|
				\GeoAng^N \upchifullmodinhom
			\right|.
			\notag
\end{align}
To bound the last term on
RHS~\eqref{E:ANNOYINGTRANSPORTINHOMOGENEOUSCOMMUTATORTERM},
we simply quote \eqref{E:CRUDEXPOINTWISEBOUND}.
We have thus proved the desired estimate \eqref{E:SECONDRERNORMALIZEDTRANSPORTEQNINHOMOGENEOUSTERMBOUND}.

We now prove \eqref{E:THIRDRERNORMALIZEDTRANSPORTEQNINHOMOGENEOUSTERMBOUND}.
To bound the term
$
\displaystyle
2 \mytr \upchi \GeoAng^N \upchifullmodinhom
$
from RHS~\eqref{E:TOPORDERTRCHIJUNKRENORMALIZEDTRANSPORT},
we use the estimate \eqref{E:PURETANGENTIALCHICOMMUTEDLINFINITY}
for $\mytr \upchi$ and \eqref{E:CRUDEXPOINTWISEBOUND}.

To bound the term
$[\Lunit, \GeoAng^N] \upchifullmodinhom$
from RHS~\eqref{E:TOPORDERTRCHIJUNKRENORMALIZEDTRANSPORT},
we first note that
\eqref{E:CRUDEXPOINTWISEBOUND}
and the $L^{\infty}$ estimates of Prop.~\ref{P:IMPROVEMENTOFAUX} yield that
$
\displaystyle
\left\|
	\Tanset^{\leq 9} \upchifullmodinhom
\right\|_{L^{\infty}(\Sigma_t^u)}
\leq C \varepsilon
$.
Hence, using the commutator estimate \eqref{L:TRANSVERALTANGENTIALCOMMUTATOR}
with $f = \upchifullmodinhom$ 
and \eqref{E:CRUDEXPOINTWISEBOUND},
we conclude that
$
\displaystyle
\left|
	[\Lunit, \GeoAng^N] \upchifullmodinhom
\right|
\lesssim 
\left|
	\Fullset_*^{\leq N+1;1} \Psi
\right|
+
\left|
	\Tanset^{\leq N} \GdVar
\right|
+
\left|
	\Tanset_*^{[1,N]} \BadVar
\right|
$
as desired.

To bound the product
$[\Tanset^N, \Lunit \upmu] \mytr \upchi$
from RHS~\eqref{E:TOPORDERTRCHIJUNKRENORMALIZEDTRANSPORT},
we first note that its magnitude is
$
\displaystyle
\lesssim
\left\|
	\Tanset^{[1,9]} \Lunit \upmu
\right\|_{L^{\infty}(\Sigma_t^u)}
\left|
	\Tanset^{\leq N-1} \mytr \upchi
\right|
+ 
\left|
	\Tanset^{\leq N} \Lunit \upmu
\right|
\left\|
	\Tanset^{\leq 8} \mytr \upchi
\right\|_{L^{\infty}(\Sigma_t^u)}
$.
We now bound the first product in the previous inequality 
by $\lesssim$ RHS~\eqref{E:THIRDRERNORMALIZEDTRANSPORTEQNINHOMOGENEOUSTERMBOUND}
with the help of
\eqref{E:POINTWISEESTIMATESFORGSPHEREANDITSTANGENTIALDERIVATIVES}
and \eqref{E:LUNITAPPLIEDTOTANGENTIALUPMUANDTANSETSTARLINFTY}
and the second by $\lesssim$ RHS~\eqref{E:THIRDRERNORMALIZEDTRANSPORTEQNINHOMOGENEOUSTERMBOUND}
with the help of
\eqref{E:PURETANGENTIALLUNITUPMUCOMMUTEDESTIMATE}
and
\eqref{E:PURETANGENTIALCHICOMMUTEDLINFINITY}.
A similar argument 
that takes into account the estimates
\eqref{E:LUNITTANGENTDIFFERENTIATEDLUNITSMALLIMPROVEDPOINTWISE},
\eqref{E:LUNITUPMULINFINITY},
\eqref{E:LUNITAPPLIEDTOTANGENTIALUPMUANDTANSETSTARLINFTY}
and
\eqref{E:UPMULINFTY}
yields the same bound for the term
$[\upmu, \Tanset^N] \Lunit \mytr \upchi$
from RHS~\eqref{E:TOPORDERTRCHIJUNKRENORMALIZEDTRANSPORT}.
A similar argument yields the same bound
for the term
$\left\lbrace
					\Tanset^N \left(\upmu (\mytr \upchi)^2 \right)
					- 2 \upmu \mytr \upchi \Tanset^N \mytr \upchi
\right\rbrace 
$
from RHS~\eqref{E:TOPORDERTRCHIJUNKRENORMALIZEDTRANSPORT},
the key point being that the top-order term
$\Tanset^N \mytr \upchi$ cancels from this difference.
To bound the last term
$\Tanset^N \mathfrak{A}$
from RHS~\eqref{E:TOPORDERTRCHIJUNKRENORMALIZEDTRANSPORT}
in magnitude 
by $\lesssim$ RHS~\eqref{E:THIRDRERNORMALIZEDTRANSPORTEQNINHOMOGENEOUSTERMBOUND},
we apply $\Tanset^N$ to both sides of
\eqref{E:RENORMALIZEDALPHAINHOMOGENEOUSTERM}.
The desired bound now follows from 
Lemmas~\ref{L:POINTWISEFORRECTANGULARCOMPONENTSOFVECTORFIELDS}
and \ref{L:POINTWISEESTIMATESFORGSPHEREANDITSDERIVATIVES}
and the $L^{\infty}$ estimates of Prop.~\ref{P:IMPROVEMENTOFAUX}.
We have thus established
\eqref{E:RENORMALIZEDTOPORDERTRCHIJUNKTRANSPORTINVERTED}.
\end{proof}

We now use Lemmas~\ref{L:RENORMALIZEDTOPORDERTRCHIJUNKTRANSPORTINVERTED} and \ref{L:CHIPARTIALMODSOURCETERMPOINTWISE}
to establish the proposition.

\begin{remark}[\textbf{Boxed constants affect high-order energy blowup-rates}]
	\label{R:BOXEDCONSTANTS}
The ``boxed constants'' such as the $\boxed{2}$
appearing on the RHS of inequality \eqref{E:KEYPOINTWISEESTIMATE}
and
the $\boxed{8.1}$
appearing on the RHS of 
inequality \eqref{E:TOPORDERTANGENTIALENERGYINTEGRALINEQUALITIES}
are important because they affect the blowup-rate of our high-order energy estimates
with respect to powers of $\upmu_{\star}^{-1}$. 
\end{remark}

\begin{proposition}[\textbf{The key pointwise estimate for} $(\Rad \Psi) \GeoAng^N \mytr \upchi$]
	\label{P:KEYPOINTWISEESTIMATE}
	Assume that $1 \leq N \leq 18$.
	Under the assumptions of Lemma~\ref{L:RENORMALIZEDTOPORDERTRCHIJUNKTRANSPORTINVERTED},
	the following pointwise estimate holds on $\mathcal{M}_{\Tboot,U_0}$:
	\begin{align} \label{E:KEYPOINTWISEESTIMATE}
		\left|
			(\Rad \Psi) \GeoAng^N \mytr \upchi
		\right|(t,u,\vartheta)
		& \leq
			\boxed{2} 
			\left\|
				\frac{[\Lunit \upmu]_-}{\upmu}
			\right\|_{L^{\infty}(\Sigma_t^u)}
			\left| 
				\Rad \GeoAng^N \Psi 
			\right|(t,u,\vartheta)
				\\
		 &  \ \ + 
						\boxed{4} (1 + C \varepsilon)
						\frac
						{
							\left\|
								[\Lunit \upmu]_-
							\right\|_{L^{\infty}(\Sigma_t^u)}}
						{\upmu_{\star}(t,u)}
					\int_{t'=0}^t 
						\frac
						{
							\left\|
								[\Lunit \upmu]_-
							\right\|_{L^{\infty}(\Sigma_{t'}^u)}}
						{\upmu_{\star}(t',u)}
						\left| 
							\Rad \GeoAng^N \Psi 
						\right|(t',u,\vartheta)
				\, dt'
				\notag \\
	& \ \	+ \mbox{\upshape Error},
		\notag
	\end{align}
	where
	\begin{align}  \label{E:ERRORTERMKEYPOINTWISEESTIMATE}
		\left|
			\mbox{\upshape Error}
		\right|(t,u,\vartheta)
		& \lesssim
			\frac{1}{\upmu_{\star}(t,u)}
			\left|\upchifullmodarg{\GeoAng^N} \right|(0,u,\vartheta)
			+ 
			\left|
				\Fullset_*^{\leq N+1;1} \Psi
			\right|(t,u,\vartheta)
				\\
			& 
			\ \
			+
			\frac{1}{\upmu_{\star}(t,u)}
			\left|
				\Fullset_*^{\leq N;1} \Psi
			\right|(t,u,\vartheta)
			\notag \\
		& \ \
			+
			\frac{1}{\upmu_{\star}(t,u)}
			\left|
				\Tanset^{\leq N} \GdVar
			\right|(t,u,\vartheta)
			+
			\frac{1}{\upmu_{\star}(t,u)}
			\left|
				\Tanset_*^{[1,N]} \BadVar
			\right|(t,u,\vartheta)
				\notag \\
		& \ \ 
			  + 
				\varepsilon
				\frac{1}{\upmu_{\star}(t,u)}
				\int_{t'=0}^t
					\frac{1}{\upmu_{\star}(t',u)}
					\left|
						\Fullset_*^{\leq N+1;1} \Psi
					\right|
					(t',u,\vartheta)
				\, dt'
				\notag 
				\\
		& \ \ 
			  + 
				\frac{1}{\upmu_{\star}(t,u)}
				\int_{t'=0}^t
					\left|
						\Fullset_*^{\leq N+1;1} \Psi
					\right|
					(t',u,\vartheta)
				\, dt'
				\notag 
				\\
		& \ \
				+ 
				\frac{1}{\upmu_{\star}(t,u)}
				\int_{t'=0}^t
					\frac{1}{\upmu_{\star}(t',u)}
					\left\lbrace
						\left|
							\Fullset_*^{\leq N;1} \Psi
						\right|
						+
						\left|
							\Tanset^{\leq N} \GdVar
						\right|
						+
					\left|
						\Tanset_*^{[1,N]} \BadVar
					\right|
					\right\rbrace
					(t',u,\vartheta)
				\, dt'.
				\notag
	\end{align}

	Furthermore, we have the following less precise pointwise estimate:
	\begin{align} \label{E:LESSPRECISEKEYPOINTWISEESTIMATE}
		& \left|
			\upmu \GeoAng^N \mytr \upchi
		\right|
		(t,u,\vartheta)
			\\
		& \lesssim
		\left|\upchifullmodarg{\GeoAng^N} \right|(0,u,\vartheta)
		+
		\upmu
		\left|
			\Tanset^{N+1} \Psi
		\right|(t,u,\vartheta)
		+
		\left|
			\Rad \Tanset^N \Psi
		\right|(t,u,\vartheta)
		+
		\left|
			\Fullset_*^{\leq N;1} \Psi
		\right|(t,u,\vartheta)
			\notag \\
	& \ \ 
		+
		\left|
			\Tanset^{\leq N} \GdVar
		\right|(t,u,\vartheta)
		+
		\left|
			\Tanset_*^{[1,N]} \BadVar
		\right|(t,u,\vartheta)
		\notag	\\
		& \ \
				+
				\int_{t'=0}^t 
						\frac{1}{\upmu_{\star}(t',u)}
						\left| 
							\Rad \Tanset^N \Psi 
						\right|
						(t',u,\vartheta)
				\, dt'
				+ 
				\int_{t'=0}^t 
					\left|
						\Fullset_*^{\leq N+1;1} \Psi
					\right|
					(t',u,\vartheta)
				\, dt'
							\notag \\
		& \ \   + 
						\int_{t'=0}^t 
							\frac{1}{\upmu_{\star}(t',u)}
							\left\lbrace
							\left|
								\Fullset_*^{\leq N;1} \Psi
							\right|
							+
							\left|
								\Tanset^{\leq N} \GdVar
							\right|
							+
							\left|
								\Tanset_*^{[1,N]} \BadVar
							\right|
							\right\rbrace
							(t',u,\vartheta)
						\, dt'.
						\notag
		\notag
	\end{align}

\end{proposition}

\begin{proof}
	We first prove \eqref{E:KEYPOINTWISEESTIMATE}-\eqref{E:ERRORTERMKEYPOINTWISEESTIMATE}.
	Using
	 \eqref{E:TRANSPORTRENORMALIZEDTRCHIJUNK}
	 and
	 \eqref{E:LOWESTORDERTRANSPORTRENORMALIZEDTRCHIJUNKDISCREPANCY},
	we split
	\begin{align} \label{E:DIFFICULTTERMFIRSTSPLITTING}
		(\Rad \Psi) \GeoAng^N \mytr \upchi
		& = 
		- \frac{\Rad \Psi}{\upmu} \GeoAng^N \upchifullmodinhom
		+
		\frac{\Rad \Psi}{\upmu} \upchifullmodarg{\GeoAng^N}.
	\end{align}
	We now bound the first product on RHS~\eqref{E:DIFFICULTTERMFIRSTSPLITTING}.
	Using \eqref{E:UPMUFIRSTTRANSPORT}
	and \eqref{E:TOPDERIVATIVESOFXPOINTWISEBOUND},
	we deduce that
	\begin{align} \label{E:RADPSIOVERUPMUTIMESTOPORDERXPOINTWISEBOUND}
		\left|
			\frac{\Rad \Psi}{\upmu} \GeoAng^N \upchifullmodinhom
		\right|
		& \leq
			2 
			\left|
				\frac{\Lunit \upmu}{\upmu}
			\right|
			\left|
				\Rad \GeoAng^N \Psi
			\right|
			+ C
				\left|
					\frac{1}{2} G_{\Lunit \Lunit} \Lunit \Psi
					+
					G_{\Lunit \Radunit} \Lunit \Psi
				\right|
				\left|
					\Rad \GeoAng^N \Psi
				\right|
					\\
				& \ \
		+ C
						\left|
							\Rad \Psi
						\right|
						\left|
							\Tanset^{\leq N+1} \Psi
						\right|
		+
		C
		\left|
			\frac{\Rad \Psi}{\upmu}
		\right|
		\left\lbrace
			\left|
				\Fullset_*^{\leq N;1} \Psi
			\right|
			+
			\left|
				\Tanset^{\leq N} \GdVar
			\right|
			+
			\left|
				\Tanset_*^{[1,N]} \BadVar
			\right|
		\right\rbrace.
	\notag
	\end{align}
	To handle the first product on 
	RHS~\eqref{E:RADPSIOVERUPMUTIMESTOPORDERXPOINTWISEBOUND},
	we use \eqref{E:POSITIVEPARTOFLMUOVERMUISBOUNDED}
	to deduce that
	$
	\displaystyle
	2 
	\left|
		\frac{\Lunit \upmu}{\upmu}
	\right|
	\leq
	\displaystyle
	2 
	\left|
		\frac{[\Lunit \upmu]_-}{\upmu}
	\right|
	+
	2 
	\left|
		\frac{[\Lunit \upmu]_+}{\upmu}
	\right|
	\leq
	2 
	\left\|
		\frac{[\Lunit \upmu]_-}{\upmu}
	\right\|_{L^{\infty}(\Sigma_t^u)}
	+ C
	$,
	which easily leads to the product under consideration
	being bounded by $\leq$ the sum of the first term on RHS~\eqref{E:KEYPOINTWISEESTIMATE}
	and the second term on RHS~\eqref{E:ERRORTERMKEYPOINTWISEESTIMATE}.
	To handle the second product on 
	RHS~\eqref{E:RADPSIOVERUPMUTIMESTOPORDERXPOINTWISEBOUND},
	we first note that by Lemma~\ref{L:SCHEMATICDEPENDENCEOFMANYTENSORFIELDS}, we have
	$G_{\Lunit \Lunit}, G_{\Lunit \Radunit}
	= \smoothfunction(\GdVar)
	$. 
	The $L^{\infty}$ estimates of Prop.~\ref{P:IMPROVEMENTOFAUX}
	thus yield that
	$
	\displaystyle
	\left|
		\frac{1}{2} G_{\Lunit \Lunit} \Lunit \Psi
		+
		G_{\Lunit \Radunit} \Lunit \Psi
	\right|
	\leq C \varepsilon
	$,
 	from which we easily deduce that
 	the second product on RHS~\eqref{E:RADPSIOVERUPMUTIMESTOPORDERXPOINTWISEBOUND}
	is 
	$\leq$ the term 
	$
	\left|
		\mbox{\upshape Error}
	\right|
	$
	from \eqref{E:ERRORTERMKEYPOINTWISEESTIMATE}.
	Moreover, it is easy to deduce that all products on 
	the second line of RHS~\eqref{E:RADPSIOVERUPMUTIMESTOPORDERXPOINTWISEBOUND} are  
	$
	\leq
	\left|
		\mbox{\upshape Error}
	\right|
	$,
	thanks to the estimate
	$
	\left\|
		\Rad \Psi
	\right\|_{L^{\infty}(\Sigma_t^u)}
	\leq C
	$
	(that is, \eqref{E:PSITRANSVERSALLINFINITYBOUNDBOOTSTRAPIMPROVED}).

	We now bound the second product 
	$
	\displaystyle
	\frac{\Rad \Psi}{\upmu} \upchifullmodarg{\GeoAng^N}
	$
	on RHS~\eqref{E:DIFFICULTTERMFIRSTSPLITTING}.
	We start by multiplying both sides of
	\eqref{E:RENORMALIZEDTOPORDERTRCHIJUNKTRANSPORTINVERTED}
	by 
	$
	\displaystyle
	\frac{\Rad \Psi}{\upmu}.
	$
	We first address
	the product of
	$
	\displaystyle
	\frac{\Rad \Psi}{\upmu}
	$
	and the second product 
	$
	\boxed{2} (1 + C \varepsilon)
	\cdots
	$
	on RHS~\eqref{E:RENORMALIZEDTOPORDERTRCHIJUNKTRANSPORTINVERTED}.
	We now use inequality \eqref{E:TOPDERIVATIVESOFXPOINTWISEBOUND}
	to substitute for the term
	$
	\left| \GeoAng^N \upchifullmodinhom \right|
	$
	appearing in the integrand.
	The easy terms to bound are those
	that arise from RHS~\eqref{E:TOPDERIVATIVESOFXPOINTWISEBOUND};
	using the bound
	$
	\left\|
		\Rad \Psi
	\right\|_{L^{\infty}(\Sigma_t^u)}
	\leq C
	$
	mentioned above
	and the bound
	$
	\left\|
		\Lunit \upmu
	\right\|_{L^{\infty}(\Sigma_t^u)}
	\leq C
	$
	(that is, \eqref{E:LUNITUPMULINFINITY}),
	it is easy to
	see that 
	their contribution
	to the product of
	$
	\displaystyle
	\frac{\Rad \Psi}{\upmu}
	$
	and the first product 
	$
	\boxed{2} (1 + C \varepsilon)
	\cdots
	$
	is
	$\leq$ the term 
	$
	\left|
		\mbox{\upshape Error}
	\right|
	$
	from \eqref{E:ERRORTERMKEYPOINTWISEESTIMATE}.
	It remains for us bound
	the error term generated
	by the main part of
	the integrand factor
	$
	\left| \GeoAng^N \upchifullmodinhom \right|
	$,
	which is given by the term $G_{\Lunit \Lunit} \Rad \GeoAng^N \Psi$ from LHS
	\eqref{E:TOPDERIVATIVESOFXPOINTWISEBOUND}.
	Specifically, we must bound
	\begin{align} \label{E:MAINDIFFICULTTOPORDERTERMPOINTWISEBOUND}
	\boxed{2} (1 + C \varepsilon)
	\left\lbrace
		\frac{\Rad \Psi}{\upmu}
	\right\rbrace(t,u,\vartheta)
							\int_{s=0}^t 
								\frac{[\Lunit \upmu]_-(s,u,\vartheta)}{\upmu(s,u,\vartheta)}
								\left| 
									G_{\Lunit \Lunit} \Rad \GeoAng^N \Psi 
								\right|
								(s,u,\vartheta)
							\, ds.
	\end{align}
	We use \eqref{E:RADDERIVATIVESOFGLLDIFFERENCEBOUND}
	to replace the factor
	$G_{\Lunit \Lunit}(s,u,\vartheta)$
	with $G_{\Lunit \Lunit}(t,u,\vartheta)$
	up to the error factor $C \varepsilon$.
	We then pull
	$G_{\Lunit \Lunit}(t,u,\vartheta)$
	out of the $ds$ integral,
	multiply it against
	$
	\displaystyle
	\left\lbrace
		\frac{\Rad \Psi}{\upmu}
	\right\rbrace(t,u,\vartheta)
	$,
	and use the arguments used to deduce
	\eqref{E:RADPSIOVERUPMUTIMESTOPORDERXPOINTWISEBOUND}
	as well as those given just below it
	to deduce that
	$
	\displaystyle
	\left\lbrace
		G_{\Lunit \Lunit} \frac{\Rad \Psi}{\upmu}
	\right\rbrace(t,u,\vartheta)
	= \boxed{2} \frac{[\Lunit \upmu]_-(t,u,\vartheta)}{\upmu(t,u,\vartheta)}
		+ \mathcal{O}(1)
	$.
	The portion of \eqref{E:MAINDIFFICULTTOPORDERTERMPOINTWISEBOUND}
	corresponding to the factor
	$
	\displaystyle
	\boxed{2} \frac{[\Lunit \upmu]_-(t,u,\vartheta)}{\upmu(t,u,\vartheta)}
	$
	is clearly $\leq$ the
	product 
	$\boxed{4} (1 + C \varepsilon) \cdots$
	on the second line of RHS~\eqref{E:KEYPOINTWISEESTIMATE}.
	Moreover, 
	using the bound
	$
	\left\|
		\Lunit \upmu
	\right\|_{L^{\infty}(\Sigma_t^u)}
	\leq C
	$
	mentioned above
	to bound the 
	$[\Lunit \upmu]_-(s,u,\vartheta)$ integrand factor
	and also using the bound
	$
	\left\|
		\Rad \Psi
	\right\|_{L^{\infty}(\Sigma_t^u)}
	\leq C
	$
	mentioned above,
	we find that the portion of \eqref{E:MAINDIFFICULTTOPORDERTERMPOINTWISEBOUND}
	corresponding to the factor
	$
	\mathcal{O}(1)
	$
	from above is
	$\leq$ 
	the last product on RHS~\eqref{E:ERRORTERMKEYPOINTWISEESTIMATE}
	(here we are using the simple fact that 
	the factor 
	$
	\displaystyle
	\frac{1}{\upmu_{\star}(t,u)}
	$
	on the outside of the integral in the last product is 
	$\geq 1$).
	Similarly, 
	the error factor $C \varepsilon$ mentioned above, 
	generated by replacing
	$G_{\Lunit \Lunit}(s,u,\vartheta)$
	with 
	$G_{\Lunit \Lunit}(t,u,\vartheta)$,
	leads to a term that is
	$\lesssim$
	the term on RHS~\eqref{E:ERRORTERMKEYPOINTWISEESTIMATE}
	featuring the coefficient $\varepsilon$.

	To complete the proof,
	it remains for us to bound the magnitude of the product of
	$
	\displaystyle
	\frac{\Rad \Psi}{\upmu}
	$
	and the term
	$
	C \left| \upchifullmodarg{\GeoAng^N } \right|(0,u,\vartheta)
	$
	on RHS~\eqref{E:RENORMALIZEDTOPORDERTRCHIJUNKTRANSPORTINVERTED}
	by 
	$\leq
	\left|
		\mbox{\upshape Error}
	\right|
	$
	and the magnitude of the product of
	$
	\displaystyle
	\frac{\Rad \Psi}{\upmu}
	$
	and the last time integral
	$
	C \int_{s=0}^t \cdots
	$
	on RHS~\eqref{E:RENORMALIZEDTOPORDERTRCHIJUNKTRANSPORTINVERTED}
	by 
	$\leq
	\left|
		\mbox{\upshape Error}
	\right|
	$,
	where 
	$
	\left|
		\mbox{\upshape Error}
	\right|
	$
	verifies \eqref{E:ERRORTERMKEYPOINTWISEESTIMATE}
	in each case.
	The desired estimates follow easily from the bound
	$
	\left\|
		\Rad \Psi
	\right\|_{L^{\infty}(\Sigma_t^u)}
	\leq C
	$
	mentioned above.
	We have thus proved \eqref{E:KEYPOINTWISEESTIMATE}-\eqref{E:ERRORTERMKEYPOINTWISEESTIMATE}.

	The proof of \eqref{E:LESSPRECISEKEYPOINTWISEESTIMATE} is similar
	but simpler so we omit the details. The main simplifications are the presence
	of an additional power of $\upmu$ on LHS~\eqref{E:LESSPRECISEKEYPOINTWISEESTIMATE}
	and that we no longer have to observe the special structure
	that led to the factor $\Lunit \upmu$ on RHS 
	\eqref{E:RADPSIOVERUPMUTIMESTOPORDERXPOINTWISEBOUND}.
\end{proof}

\subsection{Pointwise estimates for the error terms generated by the multiplier vectorfield}
\label{SS:POINTWISEESTIMATESFORTHEMULTIPLIERVECTORFIELDERRORTERMS}
In this subsection, we derive simple pointwise
estimates for the energy estimate error terms
generated by the deformation tensor of the multiplier vectorfield $\Mult$.

\begin{lemma}[\textbf{Pointwise bounds for the error terms generated by the deformation tensor of} $\Mult$]
\label{L:MULTIPLIERVECTORFIEDERRORTERMPOINTWISEBOUND}
	Consider the multiplier vectorfield error terms 
	$\basicenergyerrorarg{\Mult}{1}[\Psi],\cdots,\basicenergyerrorarg{\Mult}{5}[\Psi]$
	defined in \eqref{E:MULTERRORINTEG1}-\eqref{E:MULTERRORINTEG5}.
	Let $\varsigma > 0$ be a real number.
	Under the data-size and bootstrap assumptions 
	of Subsects.~\ref{SS:SIZEOFTBOOT}-\ref{SS:PSIBOOTSTRAP}
	and the smallness assumptions of Subsect.~\ref{SS:SMALLNESSASSUMPTIONS}, 
	the following pointwise inequality holds
	on $\mathcal{M}_{\Tboot,U_0}$ (without any absolute value taken on the left),
	where the implicit constants are independent of $\varsigma$:
	\begin{align} \label{E:MULTIPLIERVECTORFIEDERRORTERMPOINTWISEBOUND}
		\sum_{i=1}^5 \basicenergyerrorarg{\Mult}{i}[\Psi]
		& \lesssim
			(1 + \varsigma^{-1})(\Lunit \Psi)^2
			+ (1 + \varsigma^{-1}) (\Rad \Psi)^2
			+ \upmu |\angdiff \Psi|^2
			+ \varsigma \TranminusdatasizeWithFactor |\angdiff \Psi|^2
				\\
		& \ \
				+ \frac{1}{\sqrt{\Tboot - t}} 
				\upmu |\angdiff \Psi|^2.
				\notag
\end{align}
\end{lemma}

\begin{proof}
	See Subsect.~\ref{SS:OFTENUSEDESTIMATES} for some comments on the analysis.
	Using 
	the schematic relations
	\eqref{E:TENSORSDEPENDINGONGOODVARIABLESGOODPSIDERIVATIVES}
	and
	\eqref{E:TENSORSDEPENDINGONGOODVARIABLESBADDERIVATIVES},
	the estimate \eqref{E:ANGDIFFXI},
	and the $L^{\infty}$ estimates
	of Props.~\ref{P:IMPROVEMENTOFAUX} and \ref{P:IMPROVEMENTOFHIGHERTRANSVERSALBOOTSTRAP},
	it is straightforward to verify that the terms in braces on 
	the RHS of \eqref{E:MULTERRORINTEG1}, \eqref{E:MULTERRORINTEG2},
	\eqref{E:MULTERRORINTEG4}, and \eqref{E:MULTERRORINTEG5}
	are bounded in magnitude by $\lesssim 1$.
	It follows that 
	for $i=1,2,4,5$,
	$\left| 
		\basicenergyerrorarg{\Mult}{i}[\Psi]
	\right|
	$
	is
	$\lesssim$ the sum of the terms on the first line of
	RHS~\eqref{E:MULTIPLIERVECTORFIEDERRORTERMPOINTWISEBOUND}.
	The factors of $\varsigma$ and $\TranminusdatasizeWithFactor$
	appear on RHS~\eqref{E:MULTIPLIERVECTORFIEDERRORTERMPOINTWISEBOUND} because
	we use Young's inequality to bound
	$\basicenergyerrorarg{\Mult}{4}[\Psi] \lesssim |\Lunit \Psi||\angdiff \Psi| 
	\leq \varsigma^{-1} \TranminusdatasizeWithFactor^{-1}(\Lunit \Psi)^2 
	+ \varsigma \TranminusdatasizeWithFactor |\angdiff \Psi|^2
	\leq C \varsigma^{-1}  (\Lunit \Psi)^2 + C \varsigma \TranminusdatasizeWithFactor |\angdiff \Psi|^2
	$.
	Similar remarks apply to $\basicenergyerrorarg{\Mult}{5}[\Psi]$.

	To bound the term 
	$\basicenergyerrorarg{\Mult}{3}[\Psi]$,
	we also need to use the estimates
	\eqref{E:POSITIVEPARTOFLMUOVERMUISBOUNDED}
	and 
	\eqref{E:UNIFORMBOUNDFORMRADMUOVERMU},
	which allow us to bound the
	first two terms in braces on RHS~\eqref{E:MULTERRORINTEG3}.
	Note that since no absolute value 
	is taken on LHS~\eqref{E:MULTIPLIERVECTORFIEDERRORTERMPOINTWISEBOUND},
	we are free to replace the factor 
	$(\Rad \upmu)/\upmu$ 
	from RHS~\eqref{E:MULTERRORINTEG3}
	with the factor $[\Rad \upmu]_+/\upmu$
	(which we bounded in \eqref{E:UNIFORMBOUNDFORMRADMUOVERMU}).
\end{proof}

\subsection{Pointwise estimates for the partially modified quantities}
\label{SS:POINTWISEFORPARTIALLYMODIFIED}
Recall that the partially modified quantity $\upchipartialmodarg{\GeoAng^{N-1}}$
verifies the transport equation \eqref{E:COMMUTEDTRCHIJUNKFIRSTPARTIALRENORMALIZEDTRANSPORTEQUATION}.
In this subsection, we use the transport equation to derive pointwise estimates
for $\upchipartialmodarg{\GeoAng^{N-1}}$ and its $\Lunit$ derivative.

\begin{lemma}[\textbf{Pointwise estimates for the partially modified quantities}]
	\label{L:SHARPPOINTWISEPARTIALLYMODIFIED}
	Assume that $1 \leq N \leq 18$ and let $\upchipartialmodarg{\GeoAng^{N-1}}$
	be the partially modified quantity defined by \eqref{E:TRANSPORTPARTIALRENORMALIZEDTRCHIJUNK}.
	Under the data-size and bootstrap assumptions 
	of Subsects.~\ref{SS:SIZEOFTBOOT}-\ref{SS:PSIBOOTSTRAP}
	and the smallness assumptions of Subsect.~\ref{SS:SMALLNESSASSUMPTIONS}, 
	the following pointwise estimates hold
	on $\mathcal{M}_{\Tboot,U_0}$: 
	\begin{subequations}
	\begin{align} 
		\left|
			\Lunit \upchipartialmodarg{\GeoAng^{N-1}}
		\right|
		& 
		\leq
		\frac{1}{2} 
		\left| 
			G_{\Lunit \Lunit} 
		\right|
			\left|
					\angLap \GeoAng^{N-1} \Psi
			\right|
				 +
			C \varepsilon 
			\left|
				\Tanset^{\leq N+1} \Psi
			\right|
			+
			C \varepsilon
			\left|
				\Tanset^{\leq N} \GdVar
			\right|,
				\label{E:SHARPPOINTWISELUNITPARTIALLYMODIFIED} \\
		\left|
			\upchipartialmodarg{\GeoAng^{N-1}}
		\right|
		(t,u,\vartheta)
		& 
		\leq
		\left|
			\upchipartialmodarg{\GeoAng^{N-1}}
		\right|
		(0,u,\vartheta)
		+
		\frac{1}{2} 
		\left| 
			G_{\Lunit \Lunit} 
		\right|
		(t,u,\vartheta)
			\int_{t'=0}^t
				\left|
					\angLap \GeoAng^{N-1} \Psi
				\right|
				(t',u,\vartheta)
			\, dt'
			\label{E:SHARPPOINTWISEPARTIALLYMODIFIED} \\
		& \ \
		+
		C \varepsilon 
			\int_{t'=0}^t
			\left\lbrace
				\left|
					\Tanset^{\leq N+1} \Psi
				\right|
				+
				\left|
					\Tanset^{\leq N} \GdVar
				\right|
			\right\rbrace
			(t',u,\vartheta)
			\, dt'.
			\notag
	\end{align}
	\end{subequations}
\end{lemma}

\begin{proof}
	To prove \eqref{E:SHARPPOINTWISELUNITPARTIALLYMODIFIED},
	we must bound the terms on RHS~\eqref{E:COMMUTEDTRCHIJUNKFIRSTPARTIALRENORMALIZEDTRANSPORTEQUATION},
	where $\GeoAng^{N-1}$ is in the role of $\Tanset^{N-1}$.
	Clearly the first term on RHS~\eqref{E:SHARPPOINTWISELUNITPARTIALLYMODIFIED}
	arises from the first term on RHS~\eqref{E:COMMUTEDTRCHIJUNKFIRSTPARTIALRENORMALIZEDTRANSPORTEQUATION}.
	To bound the terms on RHS~\eqref{E:TRCHIJUNKCOMMUTEDTRANSPORTEQNPARTIALRENORMALIZATIONINHOMOGENEOUSTERM},
	we simply quote \eqref{E:CHIPARTIALMODSOURCETERMPOINTWISE}.

	To derive \eqref{E:SHARPPOINTWISEPARTIALLYMODIFIED}, we integrate
	\eqref{E:SHARPPOINTWISELUNITPARTIALLYMODIFIED} along the integral curves of 
	$\Lunit$. The only subtle point is that we 
	bound the time integral of the first term on RHS~\eqref{E:SHARPPOINTWISELUNITPARTIALLYMODIFIED}
	as follows
	by using \eqref{E:RADDERIVATIVESOFGLLDIFFERENCEBOUND} 
	with $M=0$ and $s=t'$:
	$
	 	\int_{t'=0}^t
	 		\left\lbrace
	 		\left| G_{\Lunit \Lunit} \right|
				\left|
						\angLap \GeoAng^{N-1} \Psi
				\right|
			\right\rbrace
			(t',u,\vartheta)
		\, dt'
		\leq
		\left| G_{\Lunit \Lunit} \right|(t,u,\vartheta)
		\int_{t'=0}^t
	 		\left|
				\angLap \GeoAng^{N-1} \Psi
			\right|
			(t',u,\vartheta)
		\, dt'
		+ 
		C \varepsilon
		\int_{t'=0}^t
	 		\left|
				\Tanset^{\leq N+1} \Psi
			\right|
			(t',u,\vartheta)
		\, dt'
	$.
\end{proof}

\section{Sobolev embedding and estimates for the change of variables map}
\label{S:SOBOLEVEMBEDDING}
\setcounter{equation}{0}
In this section, 
we provide some simple Sobolev embedding
estimates adapted to the $\ell_{t,u}$. We use them  
in the proof of our main theorem, 
after deriving energy estimates, 
in order to recover the fundamental $L^{\infty}$ bootstrap assumptions
\eqref{E:PSIFUNDAMENTALC0BOUNDBOOTSTRAP} for $\Psi$.
We also derive a basic regularity estimate for the change
of variables map $\Upsilon$ 
from Def.~\ref{D:CHOVMAP}.

\subsection{Estimates for some \texorpdfstring{$\ell_{t,u}-$tangent}{tangent} vectorfields}
We start with the following preliminary lemma.

\begin{lemma}[\textbf{Comparison of} $\GeoAng$, $\NonRadialRad$, \textbf{and} $\CoordAng$]
\label{L:COORDCOMPESTIMATES}
	Recall that $\GeoAng$ is the commutation vectorfield \eqref{E:GEOANGDEF}
	and that $\CoordAng$ is the geometric torus coordinate partial derivative vectorfield.
	There exists a scalar function $\GeoAngCoordComp$ such that
	\begin{align} \label{E:GEOANGCOORD}
		\GeoAng = \GeoAngCoordComp \CoordAng.
	\end{align}   

	Moreover, 
	under the data-size and bootstrap assumptions 
	of Subsects.~\ref{SS:SIZEOFTBOOT}-\ref{SS:PSIBOOTSTRAP}
	and the smallness assumptions of Subsect.~\ref{SS:SMALLNESSASSUMPTIONS}, 
	the following pointwise estimates hold
	on $\mathcal{M}_{\Tboot,U_0}$: 
	\begin{align} \label{E:COORDCOMPESTIMATES}
		\left|
			\GeoAngCoordComp - 1
		\right|
		& \lesssim \varepsilon,
		\qquad
		\left|
			\Lunit \GeoAngCoordComp
		\right|,
			\,
		\left|
			\GeoAng \GeoAngCoordComp
		\right|
		\lesssim 
			\varepsilon,
		\qquad 
		\left|
			\Rad \GeoAngCoordComp
		\right|
		\lesssim 1.
	\end{align}

	Similarly, the following estimate holds for the scalar-valued function $\XiCoordComp$
	from \eqref{E:RADSPLITINTOPARTTILAUANDXI}:
	\begin{align} \label{E:XICOORDCOMPESTIMATES}
		\left|
			\Fullset^{\leq 1} \XiCoordComp
		\right|
		\lesssim 1.
	\end{align}

\end{lemma}
\begin{proof}
	The existence of $\GeoAngCoordComp$ 
	is a trivial consequence of the
	fact that $\GeoAng$ is $\ell_{t,u}-$tangent.

	To prove \eqref{E:COORDCOMPESTIMATES},
	we first note the data estimates
	$\left\|
		\GeoAngCoordComp - 1
	\right\|_{L^{\infty}(\ell_{0,u})},
		\,
	\left\|
		\GeoAng \GeoAngCoordComp 
	\right\|_{L^{\infty}(\ell_{0,u})}
	\lesssim \varepsilon
	$
	and
	$
	\left\|
		\Rad \GeoAngCoordComp 
	\right\|_{L^{\infty}(\ell_{0,u})}
	\lesssim 1
	$.
	These data estimates are a simple consequence of the fact that
	$\CoordAng = \partial_2$ when $t = 0$,
	equations 
	\eqref{E:LITTLEGDECOMPOSED}-\eqref{E:METRICPERTURBATIONFUNCTION},
	\eqref{E:DOWNSTAIRSUPSTAIRSSRADUNITPLUSLUNITISAFUNCTIONOFPSI},
	\eqref{E:GEOANGINTERMSOFEUCLIDEANANGANDRADUNIT},
	and
	\eqref{E:FLATYDERIVATIVERADIALCOMPONENT},
	Remark~\ref{R:RADUNITSMALLLUNITSMALLRELATION},
	the data assumptions \eqref{E:PSIDATAASSUMPTIONS},
	and Lemma~\ref{L:BEHAVIOROFEIKONALFUNCTIONQUANTITIESALONGSIGMA0}.
	A similar argument that also relies on the last
	identity of \eqref{E:INITIALRELATIONS}
	yields that
	$\left\|
		\XiCoordComp
	\right\|_{L^{\infty}(\ell_{0,u})},
		\,
	\left\|
		\GeoAng \XiCoordComp
	\right\|_{L^{\infty}(\ell_{0,u})}
	\lesssim \varepsilon
	$
	and
	$
	\left\|
		\Rad \XiCoordComp
	\right\|_{L^{\infty}(\ell_{0,u})}
	\lesssim 1
	$.

	See Subsect.~\ref{SS:OFTENUSEDESTIMATES} for some comments on the analysis.
	Next, we use
	\eqref{E:DOWNSTAIRSUPSTAIRSSRADUNITPLUSLUNITISAFUNCTIONOFPSI},
	\eqref{E:GEOANGINTERMSOFEUCLIDEANANGANDRADUNIT},
	\eqref{E:FLATYDERIVATIVERADIALCOMPONENT},
	\eqref{E:CHIINTERMSOFOTHERVARIABLES},
	\eqref{E:GEOANGDEFORMSPHEREL},
	and Lemma~\ref{L:SCHEMATICDEPENDENCEOFMANYTENSORFIELDS}
	to deduce that
	$\GeoAngCoordComp$ satisfies the evolution equation
	\begin{align} \label{E:COORDCOMPEVOLUTIONEQN}
		\Lunit \ln \GeoAngCoordComp
		& = \frac{\angdeformoneformarg{\GeoAng}{\Lunit} \cdot \GeoAng}{|\GeoAng|^2}
		= \smoothfunction(\GdVar, \ginversesphere, \angdiff x^1, \angdiff x^2) \Singletan \GdVar.
	\end{align}
	In deriving \eqref{E:COORDCOMPEVOLUTIONEQN}, 
	we used
	the identities
	$\angLie_{\Lunit} \CoordAng = 0$,
	$\angLie_{\Lunit} \GeoAng = [\Lunit, \GeoAng] = \angdeformoneformupsharparg{\GeoAng}{\Lunit}$
	and $\angLie_{\GeoAng} \gsphere = \angdeform{\GeoAng}$
	(see Lemma~\ref{L:CONNECTIONBETWEENCOMMUTATORSANDDEFORMATIONTENSORS}).
	Hence, from 
	Lemmas~\ref{L:POINTWISEFORRECTANGULARCOMPONENTSOFVECTORFIELDS}
	and
	\ref{L:POINTWISEESTIMATESFORGSPHEREANDITSDERIVATIVES}
	and the $L^{\infty}$ estimates of Prop.~\ref{P:IMPROVEMENTOFAUX},
	we deduce
	$
	\displaystyle
	\left|
		\Lunit \ln \GeoAngCoordComp
	\right|
	\lesssim \varepsilon
	$.
	Integrating along the integral curves of
	$\Lunit$ as in \eqref{E:INTEGRATINGALONGINTEGRALCURVES}
	and using the data estimates and the previous estimate,
	we conclude the desired estimates \eqref{E:COORDCOMPESTIMATES}
	for $\GeoAngCoordComp$ and $\Lunit \GeoAngCoordComp$.

	To derive the estimate for $\GeoAng \GeoAngCoordComp$, we commute \eqref{E:COORDCOMPEVOLUTIONEQN}
	with $\GeoAng$ to obtain
	\begin{align}
	\Lunit \GeoAng \ln \GeoAngCoordComp 
		& = \angdeformoneformupsharparg{\GeoAng}{\Lunit} \cdot \angdiff \ln \GeoAngCoordComp
		+ \GeoAng
		\left\lbrace
			\smoothfunction(\GdVar, \ginversesphere, \angdiff x^1, \angdiff x^2) \Singletan \GdVar
		\right\rbrace
			\label{E:LAPPLIEDTOGEOANGLNGOEANGCOORD} \\
	& = 
		\smoothfunction(\GdVar, \ginversesphere, \angdiff x^1, \angdiff x^2) (\Singletan \GdVar) \angdiff \ln \GeoAngCoordComp
		+ 
		\GeoAng
		\left\lbrace
			\smoothfunction(\GdVar, \ginversesphere, \angdiff x^1, \angdiff x^2) \Singletan \GdVar
		\right\rbrace.
		\notag
	\end{align}
	Using the same estimates as before, we deduce from \eqref{E:LAPPLIEDTOGEOANGLNGOEANGCOORD} that
	$\left|
		\Lunit \GeoAng \ln \GeoAngCoordComp 
	\right|
	\lesssim \varepsilon \left| \GeoAng \ln \GeoAngCoordComp  \right|
	+ \varepsilon
	$. 
	Hence, integrating along the integral curves of $\Lunit$ as before
	and using the data estimates and Gronwall's inequality, we
	conclude that $\left| \GeoAng \ln \GeoAngCoordComp  \right| \lesssim \varepsilon$.
	Combining this estimate with $\GeoAngCoordComp = 1 + \mathcal{O}(\varepsilon)$, we 
	conclude the desired estimate \eqref{E:COORDCOMPESTIMATES} 
	for $\GeoAng \GeoAngCoordComp$.

	To derive the estimate for $\Rad \GeoAngCoordComp$, we commute \eqref{E:COORDCOMPEVOLUTIONEQN}
	with $\Rad$ and use above reasoning as well as the formula \eqref{E:RADDEFORMSPHERERAD} 
	to obtain
	\begin{align}
	\Lunit \Rad \ln \GeoAngCoordComp 
	& = \angdeformoneformupsharparg{\Rad}{\Lunit} \cdot \angdiff \ln \GeoAngCoordComp
		+ \Rad
			\left\lbrace
				\smoothfunction(\GdVar, \ginversesphere, \angdiff x^1, \angdiff x^2) \Singletan \GdVar
			\right\rbrace
			\label{E:LAPPLIEDTORADLNGOEANGCOORD} \\
	& = 
		\smoothfunction(\Fullset^{\leq 1} \GdVar, \Tanset^{\leq 1} \BadVar, \ginversesphere, \angdiff x^1, \angdiff x^2) 
			\angdiff \ln \GeoAngCoordComp
		+ 
		\Rad
		\left\lbrace
			\smoothfunction(\GdVar, \ginversesphere, \angdiff x^1, \angdiff x^2) \Singletan \GdVar
		\right\rbrace.
			\notag
	\end{align}
	Using the same estimates as before
	and the already proven estimates for $\GeoAngCoordComp$ and $\GeoAng \GeoAngCoordComp$ 
	(which imply that $\left|\angdiff \ln \GeoAngCoordComp \right| \lesssim \varepsilon$),
	we deduce from \eqref{E:LAPPLIEDTORADLNGOEANGCOORD} that
	$\left|
		\Lunit \Rad \ln \GeoAngCoordComp 
	\right|
	\lesssim \varepsilon
	$. 
	Hence, integrating along the integral curves of $\Lunit$ as before
	and using the data estimates, 
	we conclude that $\left| \Rad \ln \GeoAngCoordComp  \right| \lesssim 1$.
	Combining this estimate with $\GeoAngCoordComp = 1 + \mathcal{O}(\varepsilon)$, we 
	conclude the desired estimate \eqref{E:COORDCOMPESTIMATES} 
	for $\Rad \GeoAngCoordComp$.

	Next, we use \eqref{E:RADSPLITINTOPARTTILAUANDXI},
	\eqref{E:CONNECTIONBETWEENCOMMUTATORSANDDEFORMATIONTENSORS},
	and the fact that $[\Lunit,\frac{\partial}{\partial u}] = 0$
	to deduce that
	$\angLie_{\Lunit} \NonRadialRad 
	= - \angLie_{\Lunit} \Rad
	= - \angdeformoneformarg{\Rad}{\Lunit}
	$.
	Combining this identity
	with
	\eqref{E:RADDEFORMSPHERERAD}
	and
	\eqref{E:GEOANGCOORD}
	and arguing as in the previous paragraph,
	we derive the evolution equation
	\begin{align} \label{E:XICOORDCOMPEVOLUTIONEQN}
		\Lunit \XiCoordComp
		& =
			-
			\frac{1}{g(\CoordAng,\CoordAng)}
			\angdeformoneformarg{\Rad}{\Lunit}
			\cdot \CoordAng
			=
			\frac{\GeoAngCoordComp }{|\GeoAng|^2} 
			\angdeformoneformarg{\Rad}{\Lunit}
			\cdot \GeoAng
			= \GeoAngCoordComp 
				\smoothfunction(\Fullset^{\leq 1} \GdVar, \Tanset^{\leq 1} \BadVar, \ginversesphere, \angdiff x^1, \angdiff x^2).
	\end{align}
	Using the same estimates as in the previous paragraph and \eqref{E:COORDCOMPESTIMATES},
	we find that 
	$
	\left|
		\Lunit \XiCoordComp
	\right|
	\lesssim 1
	$
	as desired.
	Moreover, integrating along the integral curves of $\Lunit$ as before
	and using the data estimates, we
	conclude that $\left| \XiCoordComp  \right| \lesssim 1$ as desired.
	It remains for us to derive the desired estimates for
	$\GeoAng \XiCoordComp$
	and $\Rad \XiCoordComp$.
	To this end, we commute \eqref{E:XICOORDCOMPEVOLUTIONEQN}
	with $\GeoAng$ and $\Rad$ and use
	the same arguments as in the previous two paragraphs
	as well as the $L^{\infty}$ estimates of Prop.~\ref{P:IMPROVEMENTOFHIGHERTRANSVERSALBOOTSTRAP}
	and \eqref{E:COORDCOMPESTIMATES}
	to deduce that
	$
	\left|
		\Lunit \GeoAng \XiCoordComp
	\right|,
	\left|
		\Lunit \Rad \XiCoordComp
	\right|
	\lesssim 1
	$.
	Integrating along the integral curves of $\Lunit$ as before
	and using the data estimates, we conclude the desired bounds
	$
	\left|
		\GeoAng \XiCoordComp
	\right|,
	\left|
		\Rad \XiCoordComp
	\right|
	\lesssim 1
	$.

\end{proof}

\subsection{Comparison estimates for length forms \texorpdfstring{on $\ell_{t,u}$}{}}
\label{SS:LENGTHFORMESTIMATE}
Before proving our Sobolev embedding result,
we first establish a comparison result for the length forms
$d \spherevol$ and $d \vartheta$ on $\ell_{t,u}$.
We start with a preliminary lemma in which we derive simple
pointwise estimates for the metric component $\gtancomp$
defined in \eqref{E:METRICANGULARCOMPONENT}.

\begin{lemma}[\textbf{Pointwise estimates for} $\gtancomp$]
\label{L:POINTWISEESTIMATEFORGTANCOMP}
Let $\gtancomp$ be the metric component from
Def.~\ref{D:METRICANGULARCOMPONENT}.
Under the data-size and bootstrap assumptions 
of Subsects.~\ref{SS:SIZEOFTBOOT}-\ref{SS:PSIBOOTSTRAP}
and the smallness assumptions of Subsect.~\ref{SS:SMALLNESSASSUMPTIONS}, 
the following estimate holds
on $\mathcal{M}_{\Tboot,U_0}$:
\begin{align} \label{E:POINTWISEESTIMATEFORGTANCOMP}
	\gtancomp = 1 + \mathcal{O}(\varepsilon).
\end{align}
\end{lemma}
\begin{proof}
	From \eqref{E:LDERIVATIVEOFVOLUMEFORMFACTOR} and the estimate
	\eqref{E:PURETANGENTIALCHICOMMUTEDLINFINITY}, we deduce
	that $\Lunit \ln \gtancomp = \mathcal{O}(\varepsilon)$.
	Integrating the previous estimate
	along the integral curves of
	$\Lunit$ as in \eqref{E:INTEGRATINGALONGINTEGRALCURVES},
	we find that $\ln \gtancomp(t,u,\vartheta) = \ln \gtancomp(0,u,\vartheta) + \mathcal{O}(\varepsilon)$.
	To complete the proof, we need only to show that
	$\gtancomp(0,u,\vartheta) = 1 + \mathcal{O}(\varepsilon)$.
	To this end, we note that by construction of the geometric coordinates, 
	at $t=0$ we have $u = 1 - x^1$ and $\vartheta = x^2$, which implies that
	$\CoordAng = \partial_2$. Therefore,
	$\gtancomp^2|_{t=0} = g(\CoordAng,\CoordAng)|_{t=0} = g_{22}|_{t=0}$.
	Using \eqref{E:LITTLEGDECOMPOSED}-\eqref{E:METRICPERTURBATIONFUNCTION}
	and the bootstrap assumptions \eqref{E:PSIFUNDAMENTALC0BOUNDBOOTSTRAP},
	we conclude that $g_{22} = 1 + \mathcal{O}(\Psi) = 1 + \mathcal{O}(\varepsilon)$,
	from which the desired estimate $\gtancomp(0,u,\vartheta) = 1 + \mathcal{O}(\varepsilon)$ easily follows.
\end{proof}

\begin{lemma}[\textbf{Comparison of the forms} $d \spherevol$ \textbf{and} $d \vartheta$]
\label{L:LINEVOLUMEFORMCOMPARISON}
Let $p=p(\vartheta)$ be a non-negative function of $\vartheta$. 
Under the data-size and bootstrap assumptions 
of Subsects.~\ref{SS:SIZEOFTBOOT}-\ref{SS:PSIBOOTSTRAP}
and the smallness assumptions of Subsect.~\ref{SS:SMALLNESSASSUMPTIONS}, 
the following estimates hold
for $(t,u) \in [0,\Tboot) \times [0,U_0]$:
\begin{align} \label{E:LINEVOLUMEFORMCOMPARISON}
	(1 - C \varepsilon) \int_{\vartheta \in \mathbb{T}} p(\vartheta) d \vartheta
	\leq
	\int_{\ell_{t,u}} p(\vartheta) d \argspherevol{(t,u,\vartheta)}
	& \leq (1 + C \varepsilon) \int_{\vartheta \in \mathbb{T}} p(\vartheta) d \vartheta,
\end{align}
where $d \vartheta$ denotes the standard integration measure on $\mathbb{T}$.

Furthermore, let $p=p(u',\vartheta)$ be a non-negative function of $(u',\vartheta) \in [0,u] \times \mathbb{T}$
that \textbf{does not depend on $t$}.
Then for $s, t \in [0,\Tboot)$ and $u \in [0,U_0]$, 
we have:
\begin{align} \label{E:SIGMATVOLUMEFORMCOMPARISON}
	(1 - C \varepsilon)
	\int_{\Sigma_s^u} p \, d \tvol
	\leq
	\int_{\Sigma_t^u} p \, d \tvol
	& \leq 
	(1 + C \varepsilon)
	\int_{\Sigma_s^u} p \, d \tvol.
\end{align}

Finally, we have
\begin{align}\label{E:L2NORMOFCONSTANT}
	\left\|
		1
	\right\|_{L^2(\Sigma_t^u)}
	& \leq C.
\end{align}

\end{lemma}

\begin{proof}
	From \eqref{E:RESCALEDVOLUMEFORMS} and inequality \eqref{E:POINTWISEESTIMATEFORGTANCOMP},
	we deduce that 
	$d \spherevol = (1 + \mathcal{O}(\varepsilon)) \, d \vartheta$,
	which yields \eqref{E:LINEVOLUMEFORMCOMPARISON}.
	%
	\eqref{E:SIGMATVOLUMEFORMCOMPARISON} then follows as a simple consequence of
	\eqref{E:LINEVOLUMEFORMCOMPARISON} and the fact that along $\Sigma_t^u$, we have
	$d \tvol = d \argspherevol{(t,u,\vartheta)} du'$.
	Finally, to derive \eqref{E:L2NORMOFCONSTANT}, we use 
	\eqref{E:LINEVOLUMEFORMCOMPARISON} 
	and \eqref{E:SIGMATVOLUMEFORMCOMPARISON} to deduce that
	$
	\left\|
		1
	\right\|_{L^2(\Sigma_t^u)}^2
	\leq C
	\int_{u'=0}^u
		\int_{\vartheta \in \mathbb{T}}
			1
		\, d \vartheta
	\, du'
	\leq C
	$
	as desired.
\end{proof}

\subsection{Sobolev embedding \texorpdfstring{along $\ell_{t,u}$}{}}
\label{SS:SOBOLEVEMBEDDING}
We now state and prove our main Sobolev embedding result of interest.

\begin{lemma}[\textbf{Sobolev embedding along} $\ell_{t,u}$]
	\label{L:SOBOLEV}
Under the data-size and bootstrap assumptions 
of Subsects.~\ref{SS:SIZEOFTBOOT}-\ref{SS:PSIBOOTSTRAP}
and the smallness assumptions of Subsect.~\ref{SS:SMALLNESSASSUMPTIONS}, 
the following estimate holds for 
scalar-valued functions $f$ defined on $\ell_{t,u}$
for $(t,u) \in [0,\Tboot) \times [0,U_0]$:
\begin{align} \label{E:SOBOLEV}
		\left\|
			f
		\right\|_{L^{\infty}(\ell_{t,u})}
		& \leq 
			C
		\left\|
			\GeoAng^{\leq 1} f
		\right\|_{L^2(\ell_{t,u})}.
	\end{align}
\end{lemma}

\begin{proof}
	Standard Sobolev embedding yields that
	$\left\|
			f
	\right\|_{L^{\infty}(\mathbb{T})}
	\leq 
	C
	\left\|
		\CoordAng^{\leq 1} f
	\right\|_{L^2(\mathbb{T})}$,
	where the integration measure defining $\| \cdot \|_{L^2(\mathbb{T})}$ is $d \vartheta$.
	Thus, in view of Lemma~\ref{L:LINEVOLUMEFORMCOMPARISON},
	the desired estimate \eqref{E:SOBOLEV} follows from 
	\eqref{E:GEOANGCOORD}-\eqref{E:COORDCOMPESTIMATES}.
\end{proof}

\subsection{Basic estimates connected to the change of variables map}
\label{SS:BASICESTIMATESFORTHECHANGEOFVARIABLESMAP}

\begin{lemma}[\textbf{Basic estimates for the rectangular components} $\CoordAng^i$ \textbf{and} $\NonRadialRad^i$]
\label{L:XIIPOINTWISEBOUNDS}
Under the data-size and bootstrap assumptions 
of Subsects.~\ref{SS:SIZEOFTBOOT}-\ref{SS:PSIBOOTSTRAP}
and the smallness assumptions of Subsect.~\ref{SS:SMALLNESSASSUMPTIONS}, 
the following estimates hold
on $\mathcal{M}_{\Tboot,U_0}$ (for $i=1,2$): 
\begin{align} 
		\left\|
			\Fullset^{\leq 1} \CoordAng^i
		\right\|_{L^{\infty}(\Sigma_t^u)} 
		& \lesssim 1,
			\label{E:COORDANGIPOINTWISE} 
			\\
		\left\|
			\Fullset^{\leq 1} \NonRadialRad^i
		\right\|_{L^{\infty}(\Sigma_t^u)}
		& \lesssim 1,
		\label{E:XIIPOINTWISEBOUNDS}
	\end{align}
	where $\NonRadialRad$ is the $\ell_{t,u}-$tangent vectorfield from \eqref{E:RADSPLITINTOPARTTILAUANDXI}.
\end{lemma}

\begin{proof}
	The estimate \eqref{E:COORDANGIPOINTWISE}
	follows from 
	\eqref{E:GEOANGCOORD},
	\eqref{E:COORDCOMPESTIMATES},
	and the bounds 
	$
	\displaystyle
	\left\|
			\Fullset^{\leq 1} \GeoAng^i
	\right\|_{L^{\infty}(\Sigma_t^u)} 
	\lesssim 
	1 
	+
	\left\|
			\Fullset^{\leq 1} \GeoAng_{(Small)}^i
	\right\|_{L^{\infty}(\Sigma_t^u)}
	\lesssim 1
	$,
	which follow from
	\eqref{E:PERTURBEDPART},
	\eqref{E:ONERADIALGEOANGSMALLIPOINTWISE},
	and the $L^{\infty}$ estimates of Prop.~\ref{P:IMPROVEMENTOFAUX}.
	Similarly, to prove \eqref{E:XIIPOINTWISEBOUNDS},
	we use \eqref{E:RADSPLITINTOPARTTILAUANDXI} to express
	$\NonRadialRad^i = \XiCoordComp \CoordAng^i$.
	The desired bounds then follow from
	\eqref{E:XICOORDCOMPESTIMATES}
	and
	\eqref{E:COORDANGIPOINTWISE}.

\end{proof}

\begin{lemma}[\textbf{Uniform} $C^{1,1}$ \textbf{bounds for} $\Upsilon$]
	\label{L:CHOVMAPREGULARITY}
	Under the data-size and bootstrap assumptions 
	of Subsects.~\ref{SS:SIZEOFTBOOT}-\ref{SS:PSIBOOTSTRAP}
	and the smallness assumptions of Subsect.~\ref{SS:SMALLNESSASSUMPTIONS}, 
	the change of variables map $\Upsilon$ from Def.~\ref{D:CHOVMAP} 
	is a $C^{1,1}$ function of the geometric coordinates\footnote{The notation ``$C^{1,1}$'' means that
	the up-to-first order geometric coordinate partial derivatives of the $\Upsilon^{\alpha}$
	are Lipschitz continuous.} 
	that verifies the following estimates on $\mathcal{M}_{\Tboot,U_0}$:
	\begin{subequations}
	\begin{align} 
		\sum_{i_1 + i_2 + i_3 \leq 2} 
		\sum_{\alpha = 0}^2
			&
			\left\|
				\left(\frac{\partial}{\partial t}\right)^{i_1} 
				\left(\frac{\partial}{\partial u}\right)^{i_2}
				\left(\frac{\partial}{\partial \vartheta}\right)^{i_3}
				\Upsilon^{\alpha}
			\right\|_{L^{\infty}(\Sigma_t^u)}
			\leq C,
				\label{E:CHOVMAPSECONDDERIVATIVESINLINFINITY} \\
		\sum_{i_1 + i_2 + i_3 \leq 1} 
		\sum_{\alpha = 0}^2
			&
			\left|
				\left(\frac{\partial}{\partial t}\right)^{i_1} 
				\left(\frac{\partial}{\partial u}\right)^{i_2}
				\left(\frac{\partial}{\partial \vartheta} \right)^{i_3}
				\Upsilon^{\alpha}(t_2,u_2,\vartheta_2)
				-
				\left(\frac{\partial}{\partial t}\right)^{i_1} 
				\left(\frac{\partial}{\partial u}\right)^{i_2}
				\left(\frac{\partial}{\partial \vartheta}\right)^{i_3}
				\Upsilon^{\alpha}(t_1,u_1,\vartheta_1)
			\right|
				\label{E:CHOVMAPREGULARITY} \\
			& \leq
				C
				\left\lbrace
					|t_2 - t_1|
					+
					|u_2 - u_1|
					+
					|\vartheta_2 - \vartheta_1|
				\right\rbrace.
				\notag
	\end{align}
	\end{subequations}
	On RHS~\eqref{E:CHOVMAPREGULARITY}, 
	$|\vartheta_2 - \vartheta_1|$ denotes the flat distance
	between $\vartheta_2$ and $\vartheta_1$ on $\mathbb{T}$.
\end{lemma}
\begin{proof}
	Recall that $\Upsilon^{\alpha} = x^{\alpha}$ (we view these quantities as a function of the geometric coordinates).
	It is a standard embedding result relative to geometric coordinates (Morrey's inequality) 
	that \eqref{E:CHOVMAPREGULARITY} follows once we prove \eqref{E:CHOVMAPSECONDDERIVATIVESINLINFINITY}.
	Clearly $t=x^0$ is uniformly bounded in the norm 
		$\| \cdot \|_{L^{\infty}(\Sigma_t^u)}$, while the $x^i$ were bounded in \eqref{E:XIPOINTWISE}.
		The first derivatives of the $x^{\alpha}$ are the terms on RHS~\eqref{E:CHOV}.
		They were bounded in the norm $\| \cdot \|_{L^{\infty}(\Sigma_t^u)}$ in Lemmas~\ref{L:POINTWISEFORRECTANGULARCOMPONENTSOFVECTORFIELDS}
		and \ref{L:XIIPOINTWISEBOUNDS}. To bound the second derivatives of
		the $x^{\alpha}$, we first note that
		$\frac{\partial}{\partial t} = \Lunit$,
		$\frac{\partial}{\partial \vartheta} = \CoordAng = (1 + \mathcal{O}(\varepsilon))\GeoAng$ (see Lemma~\ref{L:COORDCOMPESTIMATES}),
		and $\frac{\partial}{\partial u} = \Rad + \mathcal{O}(1) \GeoAng$
		(see \eqref{E:RADSPLITINTOPARTTILAUANDXI} and Lemma~\ref{L:XIIPOINTWISEBOUNDS}).
		Hence, it suffices to bound the norm $\| \cdot \|_{L^{\infty}(\Sigma_t^u)}$ of the 
		$\Lunit$, $\GeoAng$, and $\Rad$ derivatives of the 
		scalar functions on RHS~\eqref{E:CHOV}.
		The desired bounds were derived in Lemmas~\ref{L:POINTWISEFORRECTANGULARCOMPONENTSOFVECTORFIELDS}
		and \ref{L:XIIPOINTWISEBOUNDS}. We have thus proved \eqref{E:CHOVMAPSECONDDERIVATIVESINLINFINITY}.
\end{proof}

\section{The fundamental \texorpdfstring{$L^2-$}{square integral-}controlling quantities}
\label{S:FUNDAMENTALL2CONTROLLINGQUANTITIES}
\setcounter{equation}{0}
In this section, 
we define the controlling quantities that we use in our $L^2$
analysis of solutions and exhibit their coercivity.

\begin{definition}[\textbf{The main coercive quantities used for controlling the solution and its derivatives in} $L^2$]
\label{D:MAINCOERCIVEQUANT}
In terms of the energy-flux quantities of Def.~\ref{D:ENERGYFLUX}, we define
\begin{subequations}
\begin{align}
	\totTanmax{N}(t,u)
	& := \max_{|\vec{I}| = N}
		\sup_{(t',u') \in [0,t] \times [0,u]} 
		\left\lbrace
			\enzero[\Tanset^{\vec{I}} \Psi](t',u')
			+ \flzero[\Tanset^{\vec{I}} \Psi](t',u')
		\right\rbrace,
			\label{E:Q0TANNDEF} 
			\\
	\totTanmax{[1,N]}(t,u)
	& := \max_{1 \leq M \leq N} \totTanmax{M}(t,u).
		\label{E:MAXEDQ0TANLEQNDEF} 
\end{align}
\end{subequations}

%

\end{definition}

We use the following coercive spacetime integrals 
to control non$-\upmu-$weighted error integrals involving geometric torus derivatives.

\begin{definition}[\textbf{Key coercive spacetime integrals}]
\label{D:COERCIVEINTEGRAL}
	We associate the following integrals to $\Psi$,
	where $[\Lunit \upmu]_- = |\Lunit \upmu|$ when $\Lunit \upmu < 0$
	and $[\Lunit \upmu]_- = 0$ when $\Lunit \upmu \geq 0$:
	\begin{subequations}
	\begin{align} \label{E:COERCIVESPACETIMEDEF} 
		\coercivespacetime[\Psi](t,u)
		& :=
	 	\frac{1}{2}
	 	\int_{\mathcal{M}_{t,u}}
			[\Lunit \upmu]_-
			|\angdiff \Psi|^2
		\, d \vol, 
				\\
		\coerciveTanspacetimemax{N}(t,u) 
		& := \max_{|\vec{I}| = N} \coercivespacetime[\Tanset^{\vec{I}} \Psi](t,u),
			\\
		\coerciveTanspacetimemax{[1,N]}(t,u) 
		& := \max_{1 \leq M \leq N} \coerciveTanspacetimemax{M}(t,u).
		\label{E:MAXEDCOERCIVESPACETIMEDEF}
	\end{align}
	\end{subequations}
\end{definition}

\begin{remark}[\textbf{We derive energy estimates only for the $P-$commuted wave equation with $P \in \Tanset$}]
We stress that definitions
\eqref{E:MAXEDQ0TANLEQNDEF}
and \eqref{E:MAXEDCOERCIVESPACETIMEDEF}
provide $L^2-$type quantities
that correspond to commuting the wave 
equation \emph{only with the elements of
the set $\Tanset$, which are
$\mathcal{P}_u-$tangent.
} 
As we described in Subsubsect.~\ref{SSS:DIFFERENCESFROMDISPERSIVEREGIME},
we rely on the special null structure of the equations
and the special properties of the vectorfields in $\Tanset$
to close our energy estimates without
deriving energy estimates for the
$\Rad-$commuted wave equation.
\end{remark}

In the next lemma, we quantify the coercive nature of the
spacetime integrals from Def.~\ref{D:COERCIVEINTEGRAL}.

\begin{lemma}[\textbf{Strength of the coercive spacetime integral}]
\label{L:KEYSPACETIMECOERCIVITY}
Let $\mathbf{1}_{\lbrace \upmu \leq 1/4 \rbrace}$
denote the characteristic function of the spacetime
subset 
$
\displaystyle
\lbrace (t,u,\vartheta) \in [0,\infty) \times [0,1] \times \mathbb{T}
	\ | \ 
\upmu(t,u,\vartheta) \leq 1/4 \rbrace
$.
Under the data-size and bootstrap assumptions 
of Subsects.~\ref{SS:SIZEOFTBOOT}-\ref{SS:PSIBOOTSTRAP}
and the smallness assumptions of Subsect.~\ref{SS:SMALLNESSASSUMPTIONS}, 
the following lower bound holds
for $(t,u) \in [0,\Tboot) \times [0,U_0]$:
\begin{align} \label{E:KEYSPACETIMECOERCIVITY}
		\coercivespacetime[\Psi](t,u) 
		& \geq 
		\frac{1}{8}
		\TranminusdatasizeWithFactor
		\int_{\mathcal{M}_{t,u}}
			\mathbf{1}_{\lbrace \upmu \leq 1/4 \rbrace}
			\left|
				\angdiff \Psi
			\right|^2
		\, d \vol.
	\end{align}
\end{lemma}

\begin{proof}
	The lemma follows easily from 
	definition \eqref{E:COERCIVESPACETIMEDEF}
	and the estimate \eqref{E:SMALLMUIMPLIESLMUISNEGATIVE}.
\end{proof}

We now provide a simple technical lemma, based on Minkowski's integral inequality,
that we will use throughout our $L^2$ analysis.

\begin{lemma}[\textbf{Estimate for the norm} $\| \cdot \|_{L^2(\Sigma_t^u)}$ \textbf{of time-integrated functions}] 
\label{L:L2NORMSOFTIMEINTEGRATEDFUNCTIONS}
Let $f$ be a scalar function on $\mathcal{M}_{\Tboot,U_0}$ and let
\begin{align} \label{E:BIGFISTIMEINETGRALOFLITTLEF}
	F(t,u,\vartheta) := \int_{t'=0}^t f(t',u,\vartheta) \, dt'.
\end{align}
Under the data-size and bootstrap assumptions 
of Subsects.~\ref{SS:SIZEOFTBOOT}-\ref{SS:PSIBOOTSTRAP}
and the smallness assumptions of Subsect.~\ref{SS:SMALLNESSASSUMPTIONS}, 
the following estimate holds
for $(t,u) \in [0,\Tboot) \times [0,U_0]$:
\begin{align} \label{E:L2NORMSOFTIMEINTEGRATEDFUNCTIONS}
	\| F \|_{L^2(\Sigma_t^u)} 
	& \leq (1 + C \varepsilon) 
		\int_{t'=0}^t 
			\| f \|_{L^2(\Sigma_{t'}^u)}
		\, dt'.
\end{align}
\end{lemma}

\begin{proof}
	Recall that 
	$\| F \|_{L^2(\Sigma_t^u)}
	: = 
				\left\lbrace
				\int_{u'=0}^u
					\int_{\ell_{t,u'}}
						F^2(t,u',\vartheta) 
					\, d \spherevol
				\, du'
				\right\rbrace^{1/2}
		$.
	Using the estimate \eqref{E:LINEVOLUMEFORMCOMPARISON},
	we may replace $d \spherevol$ in the previous formula 
	with the standard integration measure
	$d \vartheta$ on the torus $\mathbb{T}$
	up to an overall multiplicative error factor of $1 + \mathcal{O}(\varepsilon)$.
	The desired estimate \eqref{E:L2NORMSOFTIMEINTEGRATEDFUNCTIONS} follows from this estimate and
	from applying Minkowski's inequality for integrals to equation \eqref{E:BIGFISTIMEINETGRALOFLITTLEF}.

\end{proof}

In the next lemma, we quantify the coercive nature of the 
controlling quantities from Def.~\ref{D:MAINCOERCIVEQUANT}.

\begin{remark}
	\label{R:SHARPCONSTANTSMATTER}
	The sharp constants 
	$1$ and $\frac{1}{2}$
	in front
	of the quantities
	$
	\left\|
					\Rad \Tanset^{[1,N]} \Psi
				\right\|_{L^2(\Sigma_t^u)}^2
	$
	and 
	$
	\frac{1}{2}
				\left\|
					\sqrt{\upmu} \angdiff \Tanset^{[1,N]} \Psi
				\right\|_{L^2(\Sigma_t^u)}^2
	$
	in the estimate \eqref{E:COERCIVENESSOFCONTROLLING}
	influence the blowup-rate of our top-order energy estimates.
	In turn, this affects the number of derivatives that we need to close
	our estimates.
\end{remark}

\begin{lemma}[\textbf{The coercivity of} $\totTanmax{[1,N]}$]
	\label{L:COERCIVENESSOFCONTROLLING}
	Let $1 \leq M \leq N \leq 18$, and let
	$\Tanset^M$ be an $M^{th}-$order $\mathcal{P}_u$-tangent vectorfield operator.
	Under the assumptions of Lemma~\ref{L:KEYSPACETIMECOERCIVITY},
	the following lower bounds hold
	for $(t,u) \in [0,\Tboot) \times [0,U_0]$:
	\begin{align} \label{E:COERCIVENESSOFCONTROLLING}
			\totTanmax{[1,N]}(t,u)
			\geq 
			\max
			\Big\lbrace
				&
				\frac{1}{2}
				\left\|
					\sqrt{\upmu} \Lunit \Tanset^M \Psi
				\right\|_{L^2(\Sigma_t^u)}^2,
					\,
				\left\|
					\Rad \Tanset^M \Psi
				\right\|_{L^2(\Sigma_t^u)}^2,
					\,
				\frac{1}{2}
				\left\|
					\sqrt{\upmu} \angdiff \Tanset^M \Psi
				\right\|_{L^2(\Sigma_t^u)}^2,
				\\
			& 
				C^{-1}
				\left\|
					\Tanset^M \Psi
				\right\|_{L^2(\Sigma_t^u)}^2,
					\,
				\left\|
					\Lunit \Tanset^M \Psi
				\right\|_{L^2(\mathcal{P}_u^t)}^2,
					\,
				\left\|
					\sqrt{\upmu} \angdiff \Tanset^M \Psi
				\right\|_{L^2(\mathcal{P}_u^t)}^2,
					\notag \\
		& 	C^{-1}
				\left\|
					\Tanset^M \Psi
				\right\|_{L^2(\ell_{t,u})}^2
				\Big\rbrace.
				\notag
	\end{align}

	Moreover,
	\begin{subequations}
	\begin{align}
		\left\|
			\Psi
		\right\|_{L^2(\Sigma_t^u)}
		& \leq
			C \mathring{\upepsilon}
			+ 
			C
			\totTanmax{1}^{1/2}(t,u),
			 \label{E:PSIL2ESTIMATELOSSOFONEDERIVATIVE} \\
		\left\|
			\Rad \Psi
		\right\|_{L^2(\Sigma_t^u)}
		& \leq 
		C
		\left\|
			\Rad \Psi
		\right\|_{L^2(\Sigma_0^u)}
		+ 
		C \mathring{\upepsilon}
		+ 
		C \totTanmax{1}^{1/2}(t,u).
		\label{E:RADPSIL2ESTIMATELOSSOFONEDERIVATIVE}
	\end{align}
	\end{subequations}
\end{lemma}

\begin{proof}
	We first prove \eqref{E:COERCIVENESSOFCONTROLLING}.
	We prove the estimates for
	$
	\displaystyle
	\left\|
		\Tanset^M \Psi
	\right\|_{L^2(\Sigma_t^u)}^2
	$
	and
	$
	\displaystyle
	\left\|
		\Tanset^M \Psi
	\right\|_{L^2(\ell_{t,u})}^2
	$
	in detail; the other estimates in \eqref{E:COERCIVENESSOFCONTROLLING}
	follow easily from Lemma~\ref{L:ORDERZEROCOERCIVENESS}
	and we omit those details.
	To derive the estimate 
	\eqref{E:COERCIVENESSOFCONTROLLING}
	for 
	$
	\displaystyle
	\left\|
		\Tanset^M \Psi
	\right\|_{L^2(\Sigma_{t,u})}^2
	$
	and
	$
	\displaystyle
	\left\|
		\Tanset^M \Psi
	\right\|_{L^2(\ell_{t,u})}^2
	$,
	we first note that the estimates for the former quantities follow 
	easily from integrating the estimates for the latter quantities with respect to $u$.
	Hence, it suffices to prove the estimates for
	$
	\displaystyle
	\left\|
		\Tanset^M \Psi
	\right\|_{L^2(\ell_{t,u})}^2
	$,
	and for this,
	we rely on the identity \eqref{E:UDERIVATIVEOFLINEINTEGRAL}.
	Using \eqref{E:ONERADIALFORGSPHERE}
	and the $L^{\infty}$ estimates of Prop.~\ref{P:IMPROVEMENTOFAUX},
	we bound the factor 
	$(1/2)\mytr \angdeform{\Rad}$ in \eqref{E:UDERIVATIVEOFLINEINTEGRAL}
	as follows:
	$
	\displaystyle
	(1/2)
	\left|
		\mytr \angdeform{\Rad}
	\right|
	= 
	\left|
		\angLie_{\Rad} \gsphere
	\right|
	\lesssim 1
	$.
	Using the previous estimate,
	\eqref{E:UDERIVATIVEOFLINEINTEGRAL}
	with $f = (\Tanset^M \Psi)^2$, 
	Young's inequality,
	and the fact that the solution is trivial when $u = 0$,
	we deduce that
	\begin{align} \label{E:UDERIVATIVEOFSTUINTEGRALGRONWALLREADY}
	\left\|
		\Tanset^M \Psi
	\right\|_{L^2(\ell_{t,u})}^2
	& \leq 
		\int_{u'=0}^u
			\left\|
				\Rad \Tanset^M \Psi
			\right\|_{L^2(\ell_{t,u'})}^2
		\, du'
		+
		c
		\int_{u'=0}^u
			\left\|
				\Tanset^M \Psi
			\right\|_{L^2(\ell_{t,u'})}^2
		\, du'.
\end{align}
From \eqref{E:UDERIVATIVEOFSTUINTEGRALGRONWALLREADY} and Gronwall's inequality, 
we find that
\begin{align} \label{E:UDERIVATIVEOFSTUINTEGRALGRONWALLED}
	\left\|
		\Tanset^M \Psi
	\right\|_{L^2(\ell_{t,u})}^2
	& \leq 
		C e^{cu}
		\int_{u'=0}^u
			\left\|
				\Rad \Tanset^M \Psi
			\right\|_{L^2(\ell_{t,u'})}^2
		\, du'
		= C e^{cu}
			\left\|
					\Rad \Tanset^M \Psi
			\right\|_{L^2(\Sigma_t^u)}^2
		\leq C 
			\left\|
				\Rad \Tanset^M \Psi
			\right\|_{L^2(\Sigma_t^u)}^2.
\end{align}
The desired bound 
for 
$
\left\|
	\Tanset^M \Psi
\right\|_{L^2(\ell_{t,u})}^2
$
now follows from
\eqref{E:UDERIVATIVEOFSTUINTEGRALGRONWALLED}
and the already proven estimate
\eqref{E:COERCIVENESSOFCONTROLLING}
for 
$
\left\|
	\Rad \Tanset^M \Psi
\right\|_{L^2(\Sigma_t^u)}^2
$.

To derive 
	\eqref{E:RADPSIL2ESTIMATELOSSOFONEDERIVATIVE},
	we use
	\eqref{E:L2NORMSOFTIMEINTEGRATEDFUNCTIONS}
	with $f = \Lunit \Psi$
	and $F(t,u,\vartheta) = \Psi(t,u,\vartheta) - \mathring{\Psi}(u,\vartheta)$,
	where $\mathring{\Psi}(u,\vartheta) = \Psi(0,u,\vartheta)$.
	Also using the data bound
	$
	\left\|
		\mathring{\Psi}
	\right\|_{L^{\infty}(\Sigma_0^u)}
	\leq C \mathring{\upepsilon}
	$
	(see \eqref{E:PSIDATAASSUMPTIONS}),
	we find that
	$
	\displaystyle
	\left\|
		\Psi
	\right\|_{L^2(\Sigma_t^u)}
	\leq
	C \mathring{\upepsilon}
	\left\|
		1
	\right\|_{L^2(\Sigma_t^u)}
	+
	C
	\int_{s=0}^t
		\left\|
			\Lunit \Psi
		\right\|_{L^2(\Sigma_s^u)}
	\, ds
	$.
	The desired estimate now follows easily
	from this inequality,
	\eqref{E:L2NORMOFCONSTANT},
	and the estimate \eqref{E:COERCIVENESSOFCONTROLLING}
	for $\left\| \Lunit \Psi \right\|_{L^2(\Sigma_s^u)}$.

	To prove \eqref{E:RADPSIL2ESTIMATELOSSOFONEDERIVATIVE},
	we first use the commutator estimate
	\eqref{E:PURETANGENTIALFUNCTIONCOMMUTATORESTIMATE}
	and the $L^{\infty}$ estimates of Prop.~\ref{P:IMPROVEMENTOFAUX} to deduce that
	$
	\displaystyle
	\Lunit \Rad \Psi 
	= \Rad \Lunit \Psi 
		+ \mathcal{O}(\Tanset^{\leq 1} \Psi)
	$.
	Taking the norm $\left\| \cdot \right\|_{L^2(\Sigma_t^u)}$
	of this inequality, we find that
	$
	\displaystyle
	\left\| 
		\Lunit \Rad \Psi 
	\right\|_{L^2(\Sigma_t^u)}
	\leq 
		\left\| 
			\Rad \Lunit \Psi
		\right\|_{L^2(\Sigma_t^u)} 
		+
		C
		\left\| 
			\Tanset^{\leq 1} \Psi
		\right\|_{L^2(\Sigma_t^u)}
	$.
	We have already bounded all terms on the RHS of this inequality
	by 
	$\lesssim 
		\mathring{\upepsilon}
		+ 
		C
		\totTanmax{1}^{1/2}(t,u)
	$.
	Hence, much like in the previous paragraph, 
	the desired estimate 
	\eqref{E:RADPSIL2ESTIMATELOSSOFONEDERIVATIVE}
	follows from 
	\eqref{E:L2NORMSOFTIMEINTEGRATEDFUNCTIONS}
	with 
	$
		F(t,u,\vartheta) = \Rad \Psi(t,u,\vartheta) - \Rad \Psi(0,u,\vartheta)
	$
	and $f(t,u,\vartheta) = \Lunit \Rad \Psi(t,u,\vartheta)$
	and the estimate \eqref{E:LINEVOLUMEFORMCOMPARISON},
	which ensures that
	the norms
	$
	\left\|
		\cdot
	\right\|_{L^2(\Sigma_t^u)}
	$
	and 
	$
	\left\|
		\cdot
	\right\|_{L^2(\Sigma_0^u)}
	$
	are uniformly comparable
	when applied to the $t-$independent function
	$\Rad \Psi(0,u,\vartheta)$.
\end{proof}

\begin{corollary}[$L^{\infty}$ \textbf{bounds for} $\Psi$ \textbf{in terms of the fundamental controlling quantities}]
\label{C:PSILINFTYINTERMSOFENERGIES}
Under the assumptions of Lemma~\ref{L:SOBOLEV}, 
the following estimates hold for
$(t,u) \in [0,\Tboot) \times [0,U_0]$:
\begin{align}  \label{E:PSILINFTYINTERMSOFENERGIES}
	\left\|
		\Tanset^{\leq 11} \Psi
	\right\|_{L^{\infty}(\Sigma_t^u)}
	& \lesssim 
		\totTanmax{[1,12]}^{1/2}(t,u)
		+ 
		\mathring{\upepsilon}.
\end{align}

\end{corollary}

\begin{proof}
	The bound
	$
	\left\|
		\Tanset^{[1,11]} \Psi
	\right\|_{L^{\infty}(\Sigma_t^u)}
		\lesssim 
		\totTanmax{[1,12]}^{1/2}(t,u)
	$
	follows from
	Lemma~\ref{L:COERCIVENESSOFCONTROLLING}
	and Lemma~\ref{L:SOBOLEV}.
	In particular, we have 
	$\left\|
		\Lunit \Psi
	\right\|_{L^{\infty}(\Sigma_t^u)}
		\lesssim 
		\totTanmax{[1,12]}^{1/2}(t,u)
	$.
	Integrating
	along the integral curves of
	$\Lunit$ as in \eqref{E:INTEGRATINGALONGINTEGRALCURVES}
	and using this bound and the small-data assumption
	$
	\left\|
		\Psi
	\right\|_{L^{\infty}(\Sigma_0^u)}
	\leq \mathring{\upepsilon}
	$
	(see \eqref{E:PSIDATAASSUMPTIONS}),
	we deduce that
	$
	\left\|
		\Psi
	\right\|_{L^{\infty}(\Sigma_t^u)}
	\leq 
		C \totTanmax{[1,12]}^{1/2}(t,u)
		+ 
		C \mathring{\upepsilon}
	$.
	We have thus proved the corollary.
\end{proof}

\section{Energy estimates}
\label{S:ENERGYESTIMATES}
This section contains the most important technical estimates in the article:
a priori estimates
for the controlling quantities
$\totTanmax{[1,N]}$
from
Def.~\ref{D:MAINCOERCIVEQUANT}
and the coercive spacetime integrals
$\coerciveTanspacetimemax{[1,N]}$
from Def.~\ref{D:COERCIVEINTEGRAL}.
The main result is Prop.~\ref{P:MAINAPRIORIENERGY}.
To obtain the proposition,
we use the pointwise estimates of Sect.~\ref{S:POINTWISEESTIMATESFORWAVEEQUATIONERRORTERMS} 
to establish suitable estimates for
the error integrals on RHS~\eqref{E:E0DIVID},
where $\Tanset^N \Psi $ is in the role of $\Psi$
and the factor $\mathfrak{F}$ in \eqref{E:E0DIVID}
is the inhomogeneous term in the commuted wave equation
$\upmu \square_{g(\Psi)} (\Tanset^N \Psi) = \mathfrak{F}$.
We have divided the error integrals into various classes 
that we separately treat in the ensuing sections.

\subsection{Statement of the main a priori energy estimates}
\label{SS:STATEMENTOFMAINAPRIORIESTIMATES}
We start by stating the proposition featuring our main a priori
energy estimates, the proof of which is located 
in Subsect.~\ref{SS:PROOFOFPROPMAINAPRIORIENERGY}.

\begin{proposition}[\textbf{The main a priori energy estimates}]
	\label{P:MAINAPRIORIENERGY}
	Consider the fundamental $L^2-$controlling quantities 
	$\lbrace \totTanmax{[1,N]}(t,u) \rbrace_{N=1,\cdots,18}$
	from Def.~\ref{D:MAINCOERCIVEQUANT}.
	There exists a constant $C > 0$ such that
	under the data-size and bootstrap assumptions 
	of Subsects.~\ref{SS:SIZEOFTBOOT}-\ref{SS:PSIBOOTSTRAP}
	and the smallness assumptions of Subsect.~\ref{SS:SMALLNESSASSUMPTIONS}, 
	the following estimates hold
	for $(t,u) \in [0,\Tboot) \times [0,U_0]$:
	\begin{subequations}
	\begin{align}
		\totTanmax{[1,13+M]}^{1/2}(t,u)
		+ \coerciveTanspacetimemax{[1,13+M]}^{1/2}(t,u)
		& \leq C \mathring{\upepsilon} \upmu_{\star}^{-(M+.9)}(t,u),
			&& (0 \leq M \leq 5),
				\label{E:MULOSSMAINAPRIORIENERGYESTIMATES} \\
		\totTanmax{[1,1+M]}^{1/2}(t,u)
		+ \coerciveTanspacetimemax{[1,1+M]}^{1/2}(t,u)
		& \leq C \mathring{\upepsilon},
		&& (0 \leq M \leq 11).
			\label{E:NOMULOSSMAINAPRIORIENERGYESTIMATES}
	\end{align}
	\end{subequations}
\end{proposition}

We prove Prop.~\ref{P:MAINAPRIORIENERGY} through a long Gronwall
argument that relies on the sharp estimates for $\upmu$ derived in Sect.~\ref{S:SHARPESTIMATESFORUPMU}
as well as the energy inequalities provided by the following result, 
Prop.~\ref{P:TANGENTIALENERGYINTEGRALINEQUALITIES}.
The proof of the proposition is located in
Subsect.~\ref{SS:PROOFOFPROPTANGENTIALENERGYINTEGRALINEQUALITIES}.
See Remark~\ref{R:BOXEDCONSTANTS}
regarding the boxed constants on RHS~\eqref{E:TOPORDERTANGENTIALENERGYINTEGRALINEQUALITIES}.

\begin{proposition}[\textbf{Inequalities derived from energy identities}]
	\label{P:TANGENTIALENERGYINTEGRALINEQUALITIES}
	Assume that $1 \leq N \leq 18$ and $\varsigma > 0$.
	There exists a constant $C > 0$,
	independent of $\varsigma$, such that
	under the data-size and bootstrap assumptions 
	of Subsects.~\ref{SS:SIZEOFTBOOT}-\ref{SS:PSIBOOTSTRAP}
	and the smallness assumptions of Subsect.~\ref{SS:SMALLNESSASSUMPTIONS}, 
	the following pointwise estimates hold
	for $(t,u) \in [0,\Tboot) \times [0,U_0]$
	(where $2 \leq N \leq 18$ in \eqref{E:BELOWTOPORDERTANGENTIALENERGYINTEGRALINEQUALITIES}):
	\begin{subequations}
	\begin{align} \label{E:TOPORDERTANGENTIALENERGYINTEGRALINEQUALITIES}
		&
		\max\left\lbrace
			\totTanmax{[1,N]}(t,u), \coerciveTanspacetimemax{[1,N]}(t,u)
		\right\rbrace
				\\
		& \leq C (1 + \varsigma^{-1}) \mathring{\upepsilon}^2 \upmu_{\star}^{-3/2}(t,u)
				\notag \\
		& \ \ 
			+ 
			\boxed{6}
			\int_{t'=0}^t
					\frac{\left\| [\Lunit \upmu]_- \right\|_{L^{\infty}(\Sigma_{t'}^u)}} 
							 {\upmu_{\star}(t',u)} 
				  \totTanmax{[1,N]}(t',u)
				\, dt'
				\notag \\
		& \ \
			+ 
			\boxed{8.1}
			\int_{t'=0}^t
				\frac{\left\| [\Lunit \upmu]_- \right\|_{L^{\infty}(\Sigma_{t'}^u)}} 
								 {\upmu_{\star}(t',u)} 
						\totTanmax{[1,N]}^{1/2}(t',u) 
						\int_{s=0}^{t'}
							\frac{\left\| [\Lunit \upmu]_- \right\|_{L^{\infty}(\Sigma_s^u)}} 
									{\upmu_{\star}(s,u)} 
							\totTanmax{[1,N]}^{1/2}(s,u) 
						\, ds
				\, dt'
			\notag	\\
		& \ \
			+ 
			\boxed{2}
			\frac{1}{\upmu_{\star}^{1/2}(t,u)}
			\totTanmax{[1,N]}^{1/2}(t,u)
			\left\| 
				\Lunit \upmu 
			\right\|_{L^{\infty}(\Sigmaminus{t}{t}{u})}
			\int_{t'=0}^t
				\frac{1}{\upmu_{\star}^{1/2}(t',u)} \totTanmax{[1,N]}^{1/2}(t',u)
			\, dt'
			\notag \\
		& \ \
			+ 
			C \varepsilon
			\int_{t'=0}^t
				\frac{1} {\upmu_{\star}(t',u)} 
						\totTanmax{[1,N]}^{1/2}(t',u) 
						\int_{s=0}^{t'}
							\frac{1} 
									{\upmu_{\star}(s,u)} 
							\totTanmax{[1,N]}^{1/2}(s,u) 
						\, ds
				\, dt'
			\notag	\\
		& \ \
			+ 
			C \varepsilon
			\int_{t'=0}^t
					\frac{1} 
						{\upmu_{\star}(t',u)} 
				  \totTanmax{[1,N]}(t',u)
				\, dt'
				\notag
				\\
		& \ \
			+ 
			C \varepsilon
			\frac{1}{\upmu_{\star}^{1/2}(t,u)}
			\totTanmax{[1,N]}^{1/2}(t,u)
			\int_{t'=0}^t
				\frac{1}{\upmu_{\star}^{1/2}(t',u)} \totTanmax{[1,N]}^{1/2}(t',u)
			\, dt'
				\notag \\
		& \ \
			+ 
			C
			\totTanmax{[1,N]}^{1/2}(t,u)
			\int_{t'=0}^t
				\frac{1}{\upmu_{\star}^{1/2}(t',u)} \totTanmax{[1,N]}^{1/2}(t',u)
			\, dt'
				\notag \\
		& \ \
			+
			C
			\int_{t'=0}^t
				\frac{1}{\sqrt{\Tboot - t'}} \totTanmax{[1,N]}(t',u)
			 \, dt'
			 \notag \\
		& \ \
			+ C (1 + \varsigma^{-1})
					\int_{t'=0}^t
					\frac{1} 
							 {\upmu_{\star}^{1/2}(t',u)} 
				  \totTanmax{[1,N]}(t',u)
				\, dt'
				\notag	\\
		& \ \
			+ C
					\int_{t'=0}^t
					\frac{1} 
							 {\upmu_{\star}(t',u)} 
				  \totTanmax{[1,N]}^{1/2}(t',u)
				  \int_{s = 0}^{t'}
				  	\frac{1} 
							 {\upmu_{\star}^{1/2}(s,u)} 
							 \totTanmax{[1,N]}^{1/2}(s,u)
				  \, ds
				\, dt'
				\notag	\\
		& \ \
			+ C
					\int_{t'=0}^t
					\frac{1} 
							 {\upmu_{\star}(t',u)} 
				  \totTanmax{[1,N]}^{1/2}(t',u)
				  \int_{s = 0}^{t'}
				  	\frac{1}{\upmu_{\star}(s,u)}
				  	\int_{s' = 0}^s
				  	\frac{1} 
							 {\upmu_{\star}^{1/2}(s',u)} 
							 \totTanmax{[1,N]}^{1/2}(s',u)
						\, ds'	 
				  \, ds
				\, dt'
				\notag	\\
		& \ \
			+ C (1 + \varsigma^{-1})
					\int_{u'=0}^u
						\totTanmax{[1,N]}(t,u')
					\, du'
				\notag	\\
		& \ \
			+ C \varepsilon \totTanmax{[1,N]}(t,u)
			+ C \varsigma \totTanmax{[1,N]}(t,u)
			+ C \varsigma \coerciveTanspacetimemax{[1,N]}(t,u)
				\notag \\
		& \ \
			+ C
				\int_{t'=0}^t
					\frac{1} 
							 {\upmu_{\star}^{5/2}(t',u)} 
				  \totTanmax{[1,N-1]}(t',u)
				\, dt',
			\notag 
	\end{align}

	\begin{align} \label{E:BELOWTOPORDERTANGENTIALENERGYINTEGRALINEQUALITIES}
		&
		\max\left\lbrace
			\totTanmax{[1,N-1]}(t,u), \coerciveTanspacetimemax{[1,N-1]}(t,u)
		\right\rbrace
			\\
		& \leq C \mathring{\upepsilon}^2 
			\notag \\
		& \ \
			+
			C
			\int_{t'=0}^t
				\frac{1}{\upmu_{\star}^{1/2}(t',u)} 
						\totTanmax{[1,N-1]}^{1/2}(t',u) 
						\int_{s=0}^{t'}
							\frac{1}{\upmu_{\star}^{1/2}(s,u)} 
							\totTanmax{[1,N]}^{1/2}(s,u) 
						\, ds
				\, dt'
			\notag	\\
		& \ \
			+
			C
			\int_{t'=0}^t
				\frac{1}{\sqrt{\Tboot - t'}} \totTanmax{[1,N-1]}(t',u)
			 \, dt'
			 \notag \\
		& \ \
			+ C (1 + \varsigma^{-1})
					\int_{t'=0}^t
					\frac{1} 
							 {\upmu_{\star}^{1/2}(t',u)} 
				  \totTanmax{[1,N-1]}(t',u)
				\, dt'
				\notag \\
		& \ \
			+ C \mathring{\upepsilon}
					\int_{t'=0}^t
					\frac{1} 
							 {\upmu_{\star}^{1/2}(t',u)} 
				  \totTanmax{[1,N-1]}^{1/2}(t',u)
				\, dt'
				\notag \\
		& \ \
			+ C (1 + \varsigma^{-1})
					\int_{u'=0}^u
						\totTanmax{[1,N-1]}(t,u')
					\, du'
				\notag	\\
		& \ \
			+ C \varsigma \coerciveTanspacetimemax{[1,N-1]}(t,u).
				\notag
	\end{align}

	\end{subequations}

\end{proposition}

\begin{remark}[\textbf{Less degeneracy at the cost of one derivative}]
Note that the estimate \eqref{E:BELOWTOPORDERTANGENTIALENERGYINTEGRALINEQUALITIES}
\emph{
does not involve any of the difficult ``boxed-constant'' error integrals appearing on 
}
RHS~\eqref{E:TOPORDERTANGENTIALENERGYINTEGRALINEQUALITIES}.
The price paid is that the term $\totTanmax{[1,N]}^{1/2}$ 
on RHS~\eqref{E:BELOWTOPORDERTANGENTIALENERGYINTEGRALINEQUALITIES}
corresponds to one derivative above the level of LHS~\eqref{E:BELOWTOPORDERTANGENTIALENERGYINTEGRALINEQUALITIES}
(that is, the estimate \eqref{E:BELOWTOPORDERTANGENTIALENERGYINTEGRALINEQUALITIES} loses one derivative).
\end{remark}

\subsection{Preliminary \texorpdfstring{$L^2$}{square integral} estimates for the eikonal function quantities that do not require modified quantities}
\label{SS:L2FOREIKONALNOMODIFIED}
In this subsection, we provide preliminary $L^2$ estimates for some error term factors.
The main result is Lemma~\ref{L:EASYL2BOUNDSFOREIKONALFUNCTIONQUANTITIES},
in which we bound the \emph{below-top-order} derivatives of 
$\upmu$, $\Lunit_{(Small)}^i$, and $\mytr \upchi$ 
in terms of the fundamental controlling quantities of Def.~\ref{D:MAINCOERCIVEQUANT}.
These estimates are not difficult to obtain because we allow them to lose one derivative relative to $\Psi$.
We also derive estimates for the top-order derivatives involving at least one $\Lunit$ differentiation.
These estimates are also not difficult because to obtain them, we do not need to 
rely on the modified quantities of Sect.~\ref{S:MODQUANTS}.

To derive the desired estimates, we will integrate
the transport equations of Lemma~\ref{L:UPMUANDLUNITIFIRSTTRANSPORT}
and their higher-order analogs
with respect to $t$ at fixed $(u,\vartheta)$
and apply Lemma~\ref{L:L2NORMSOFTIMEINTEGRATEDFUNCTIONS}.

\begin{lemma}[$L^2$ \textbf{bounds for the eikonal function quantities that do not require modified quantities}]
	\label{L:EASYL2BOUNDSFOREIKONALFUNCTIONQUANTITIES}
	Assume that $N \leq 18$.
	Under the data-size and bootstrap assumptions 
	of Subsects.~\ref{SS:SIZEOFTBOOT}-\ref{SS:PSIBOOTSTRAP}
	and the smallness assumptions of Subsect.~\ref{SS:SMALLNESSASSUMPTIONS}, 
	the following $L^2$ estimates hold for $(t,u) \in [0,\Tboot) \times [0,U_0]$
	(see Subsect.~\ref{SS:STRINGSOFCOMMUTATIONVECTORFIELDS} regarding the vectorfield operator notation):
	\begin{subequations}
	\begin{align}
		\left\|
			\Lunit \Tanset_*^{[1,N]} \upmu
		\right\|_{L^2(\Sigma_t^u)},
			\,
		\left\|
			\Lunit \Tanset^{\leq N} \Lunit_{(Small)}^i
		\right\|_{L^2(\Sigma_t^u)},
			\,
		\left\|
			\Lunit \Tanset^{\leq N-1} 
			\mytr \upchi
		\right\|_{L^2(\Sigma_t^u)}
		& \lesssim 
				\mathring{\upepsilon}
				+
				\frac{\totTanmax{[1,N]}^{1/2}(t,u)}{\upmu_{\star}^{1/2}(t,u)},
			\label{E:LUNITTANGENGITALEIKONALINTERMSOFCONTROLLING}
				 \\
		\left\|
			\Lunit \Fullset^{\leq N;1} \Lunit_{(Small)}^i
		\right\|_{L^2(\Sigma_t^u)},
			\,
		\left\|
			\Lunit \Fullset^{\leq N-1;1} \mytr \upchi
		\right\|_{L^2(\Sigma_t^u)}
		& \lesssim 
				\mathring{\upepsilon}
				+
				\frac{\totTanmax{[1,N]}^{1/2}(t,u)}{\upmu_{\star}^{1/2}(t,u)},
			\label{E:LUNITONERADIALEIKONALINTERMSOFCONTROLLING}
				 \\
		\left\|
			\Tanset_*^{[1,N]} \upmu
		\right\|_{L^2(\Sigma_t^u)},
			\,
		\left\|
			\Tanset^{\leq N} \Lunit_{(Small)}^i
		\right\|_{L^2(\Sigma_t^u)},
			\,
		\left\|
			\Tanset^{\leq N-1} \mytr \upchi
		\right\|_{L^2(\Sigma_t^u)}
		& \lesssim 
			\mathring{\upepsilon}
			+ 
			\int_{s=0}^t
				\frac{\totTanmax{[1,N]}^{1/2}(s,u)}{\upmu_{\star}^{1/2}(s,u)}
			\, ds,
				 \label{E:TANGENGITALEIKONALINTERMSOFCONTROLLING}
				 \\
		\left\|
			\Fullset_*^{\leq N;1} \Lunit_{(Small)}^i
		\right\|_{L^2(\Sigma_t^u)},
			\,
		\left\|
			\Fullset^{\leq N-1;1} \mytr \upchi
		\right\|_{L^2(\Sigma_t^u)}
		& \lesssim 
				\mathring{\upepsilon}
				+
				\int_{s=0}^t
				\frac{\totTanmax{[1,N]}^{1/2}(s,u)}{\upmu_{\star}^{1/2}(s,u)}
			\, ds.
			\label{E:ONERADIALEIKONALINTERMSOFCONTROLLING}
	\end{align}
	\end{subequations}

\end{lemma}

\begin{proof}
	See Subsect.~\ref{SS:OFTENUSEDESTIMATES} for some comments on the analysis.
	We set 
	\begin{align} \label{E:LITTLEQNDEF}
	q_N(t) 
		& 
		:= 
		\sum_{a=1}^2
		\left\|
			\Tanset^{\leq N} 
			\Lunit_{(Small)}^a
		\right\|_{L^2(\Sigma_t^u)}
		+ 
		\left\|
			\Tanset^{\leq N-1} \mytr \upchi
		\right\|_{L^2(\Sigma_t^u)}.
	\end{align}
	From 
	\eqref{E:LUNITTANGENTDIFFERENTIATEDLUNITSMALLIMPROVEDPOINTWISE},
	Lemma~\ref{L:COERCIVENESSOFCONTROLLING},
	\eqref{E:SIGMATVOLUMEFORMCOMPARISON},
	and Lemma~\ref{L:L2NORMSOFTIMEINTEGRATEDFUNCTIONS},
	we deduce that
	\begin{align} \label{E:GRONWALLREADYEIKONALFUNCTIONQUANTITIES}
		q_N(t)
		& \leq 
		C q_N(0)
		+ C \varepsilon
			\int_{t'=0}^t
				q_N(t')
			\, dt'
			+
			C
			\int_{s=0}^t
				\frac{\totTanmax{[1,N]}^{1/2}(s,u)}{\upmu_{\star}^{1/2}(s,u)}
			\, ds.
	\end{align}
	Next, we note that $q_N(0) \lesssim \mathring{\upepsilon}$,
	an estimate that follows from 
	the estimate 
	\eqref{E:POINTWISEESTIMATESFORGSPHEREANDITSTANGENTIALDERIVATIVES} for $\mytr \upchi$
	and Lemma~\ref{L:BEHAVIOROFEIKONALFUNCTIONQUANTITIESALONGSIGMA0}.
	We now apply Gronwall's inequality to \eqref{E:GRONWALLREADYEIKONALFUNCTIONQUANTITIES}
	to conclude that $q_N(t) \lesssim$ RHS~\eqref{E:TANGENGITALEIKONALINTERMSOFCONTROLLING} as desired.
	We have thus proved the desired estimates for 
	$\Tanset^{\leq N} \Lunit_{(Small)}^i$ and $\Tanset^{\leq N-1} \mytr \upchi$.

	Next, we consider 
	the first term on RHS~\eqref{E:PURETANGENTIALLUNITUPMUCOMMUTEDESTIMATE}.
	We use the commutation estimate 
	\eqref{E:ONERADIALTANGENTIALFUNCTIONCOMMUTATORESTIMATE}
	with $f=\Psi$,
	the $L^{\infty}$ estimates of Prop.~\ref{P:IMPROVEMENTOFAUX},
	and Cor.~\ref{C:SQRTEPSILONTOCEPSILON}
	to commute the (at most one) factor of $\Rad$
	in the operator $\Fullset_*^{N+1;1}$
	to the front, which allows us to write
	$|\Fullset_*^{\leq N+1;1} \Psi|
	\lesssim
	|\Rad \Tanset^{[1,N]} \Psi|
	+ 
	|\Tanset^{\leq N+1} \Psi|
	+
	\varepsilon |\Fullset_*^{\leq N;1} \GdVar|
	+ \varepsilon |\Tanset_*^{[1,N]} \BadVar|
	$. 
	Thanks to the previous estimate and inequality \eqref{E:PURETANGENTIALLUNITUPMUCOMMUTEDESTIMATE}, 
  we can use an argument similar to the one that we used to
	derive \eqref{E:GRONWALLREADYEIKONALFUNCTIONQUANTITIES}
	in order to deduce
	\begin{align} \label{E:GRONWALLREADYUPMUINL2}
		\left\|
			\Tanset_*^{[1,N]} \upmu
		\right\|_{L^2(\Sigma_t^u)}
		& \leq \left\|
			\Tanset_*^{[1,N]} \upmu
		\right\|_{L^2(\Sigma_0^u)}
		+
		C \mathring{\upepsilon}
			+ C 
			\int_{t'=0}^t
				q_N(t')
			\, dt'
				\\
	& \ \
		+ C \varepsilon
			\int_{t'=0}^t
				\left\|
					\Tanset_*^{[1,N]} \upmu
				\right\|_{L^2(\Sigma_{t'}^u)}
			\, dt'
			+
			C
			\int_{s=0}^t
				\frac{\totTanmax{[1,N]}^{1/2}(s,u)}{\upmu_{\star}^{1/2}(s,u)}
			\, ds.
			\notag
	\end{align}
	We clarify that the term
	$C \mathring{\upepsilon}$ on RHS~\eqref{E:GRONWALLREADYUPMUINL2}
	comes from the first term on RHS~\eqref{E:PSIL2ESTIMATELOSSOFONEDERIVATIVE}.
	Moreover, 
	Lemma~\ref{L:BEHAVIOROFEIKONALFUNCTIONQUANTITIESALONGSIGMA0}
	yields that  
	$\left\|
			\Tanset_*^{[1,N]} \upmu
		\right\|_{L^2(\Sigma_0^u)} 
	\lesssim \mathring{\upepsilon}
	$, 
	while the estimates we have already derived for $\Tanset^{\leq N} \Lunit_{(Small)}^i$
	imply that
	$
	\displaystyle
	C 
	\int_{t'=0}^t
		q_N(t')
	\, dt
	\lesssim 
	\mathring{\upepsilon}
	+ 
		\int_{s=0}^t
			\frac{\totTanmax{[1,N]}^{1/2}(s,u)}{\upmu_{\star}^{1/2}(s,u)}
		\, ds
	$.
	Also using Gronwall's inequality,
	we conclude the desired estimate for
	$
	\left\|
			\Tanset_*^{[1,N]} \upmu
	\right\|_{L^2(\Sigma_t^u)}
	$.

	To obtain the estimates \eqref{E:LUNITTANGENGITALEIKONALINTERMSOFCONTROLLING},
	we take the norm 
	$\left\|
	\cdot
	\right\|_{L^2(\Sigma_t^u)} 
	$
	of the inequalities 
	\eqref{E:PURETANGENTIALLUNITUPMUCOMMUTEDESTIMATE}
	and
	\eqref{E:LUNITTANGENTDIFFERENTIATEDLUNITSMALLIMPROVEDPOINTWISE} 
	and argue as above using the already proven estimates
	\eqref{E:TANGENGITALEIKONALINTERMSOFCONTROLLING}.
	In these estimates, we encounter the integrals
	$
	\displaystyle
	\int_{s=0}^t
				\frac{\totTanmax{[1,N]}^{1/2}(s,u)}{\upmu_{\star}^{1/2}(s,u)}
			\, ds
	$,
	which we bound by 
	$
	\lesssim \totTanmax{[1,N]}^{1/2}(t,u) 
	\leq \upmu_{\star}^{-1/2}(t,u) \totTanmax{[1,N]}^{1/2}(t,u)
	$
	with the help of inequality \eqref{E:LESSSINGULARTERMSMPOINTNINEINTEGRALBOUND}.

	The proofs of \eqref{E:LUNITONERADIALEIKONALINTERMSOFCONTROLLING} 
	and \eqref{E:ONERADIALEIKONALINTERMSOFCONTROLLING}
	are similar
	and are based on inequality
	\eqref{E:LUNITONERADIALTANGENTDIFFERENTIATEDLUNITSMALLIMPROVEDPOINTWISE}
	and the already proven estimates \eqref{E:TANGENGITALEIKONALINTERMSOFCONTROLLING}; 
	we omit the details.


\end{proof}

In the next corollary, we obtain $L^2$ estimates for $\Psi$
with a right-hand side that involves no explicit
degenerate factor involving a power of $\upmu_{\star}^{-1}$.
The price we pay is that the estimates lose one derivative, that
is, the left-hand side features $N$ derivatives of $\Psi$ 
but the right-hand side features a quantity
that is capable of controlling
up to $N+1$ derivatives of $\Psi$.

\begin{corollary}[\textbf{Non-degenerate $L^2$ estimates for $\Psi$ that lose one derivative}]
	\label{C:COMMUTATIONBOUNDSFORPSIINTERMSOFFUNDAMENTALCONTROLLINGQUANTITIES}
	Assume that $N \leq 18$.
	Under the data-size and bootstrap assumptions 
	of Subsects.~\ref{SS:SIZEOFTBOOT}-\ref{SS:PSIBOOTSTRAP}
	and the smallness assumptions of Subsect.~\ref{SS:SMALLNESSASSUMPTIONS}, 
	the following $L^2$ estimates hold for $(t,u) \in [0,\Tboot) \times [0,U_0]$
	(see Subsect.~\ref{SS:STRINGSOFCOMMUTATIONVECTORFIELDS} regarding the vectorfield operator notation):
	\begin{align}	\label{E:COMMUTATIONBOUNDSFORPSIINTERMSOFFUNDAMENTALCONTROLLINGQUANTITIES}
		\left\|
			\Fullset_*^{N;1} \Psi
		\right\|_{L^2(\Sigma_t^u)}
		& 
		\lesssim 
		\totTanmax{[1,N]}^{1/2}(t,u)
		+ \mathring{\upepsilon}.
	\end{align}
\end{corollary}

\begin{proof}
	In the second paragraph of 
	the proof of Lemma~\ref{L:EASYL2BOUNDSFOREIKONALFUNCTIONQUANTITIES},
	we obtained the following estimate
	(except that here we have $N$ in the role of $N+1$):
	$|\Fullset_*^{N;1} \Psi|
	\lesssim
	|\Rad \Tanset^{[1,N-1]} \Psi|
	+ 
	|\Tanset^{\leq N} \Psi|
	+
	\varepsilon |\Fullset_*^{\leq N-1;1} \GdVar|
	+ \varepsilon |\Tanset_*^{[1,N-1]} \BadVar|
	$,
	where only the second term is present when $N=0,1$.
	The desired bound 
	\eqref{E:COMMUTATIONBOUNDSFORPSIINTERMSOFFUNDAMENTALCONTROLLINGQUANTITIES}
	follows from this estimate, 
	Lemma~\ref{L:COERCIVENESSOFCONTROLLING},
	the estimates
	\eqref{E:TANGENGITALEIKONALINTERMSOFCONTROLLING}
	and
	\eqref{E:ONERADIALEIKONALINTERMSOFCONTROLLING}
	(with $N-1$ in the role of $N$ there),
	the fact that
	$\totTanmax{[1,N]}$
	is increasing in its arguments,
	and inequality \eqref{E:LESSSINGULARTERMSMPOINTNINEINTEGRALBOUND}
	(which we use to annihilate the factors of $\upmu_{\star}^{1/2}(s,u)$
	in the denominators of the integrands 
	on RHS~\eqref{E:TANGENGITALEIKONALINTERMSOFCONTROLLING} 
	and RHS~\eqref{E:TANGENGITALEIKONALINTERMSOFCONTROLLING}).
\end{proof}

\subsection{Estimates for the easiest error integrals}
\label{SS:ENERGYESTIMATESEASIESTERRORINTEGRALS}
In this subsection, 
we derive estimates for the simplest error integrals that appear in
our energy estimates, that is, for the simplest integrals on RHS~\eqref{E:E0DIVID}. 

We start with the following simple lemma, which shows that the fundamental
controlling quantities from Def.~\ref{D:MAINCOERCIVEQUANT}
are size $\mathcal{O}(\mathring{\upepsilon}^2)$
at time $0$.

\begin{lemma}[\textbf{The fundamental controlling quantities are initially small}]
\label{L:INITIALSIZEOFL2CONTROLLING}
	Assume that $1 \leq N \leq 18$.
	Under the data-size assumptions
	of Subsect.~\ref{SS:SIZEOFTBOOT},
	the following estimates hold for $u \in [0,U_0]$: 
	\begin{align} \label{E:INITIALSIZEOFL2CONTROLLING}
		\totTanmax{N}(0,u) \lesssim \mathring{\upepsilon}^2.
	\end{align}
\end{lemma}
\begin{proof}
	From Def.~\ref{D:MAINCOERCIVEQUANT},
	Lemma~\ref{L:ORDERZEROCOERCIVENESS},
	and
	Lemma~\ref{L:BEHAVIOROFEIKONALFUNCTIONQUANTITIESALONGSIGMA0}
	(which in particular implies that $\upmu \approx 1$ along $\Sigma_0^1$),
	we see that
	$
		\totTanmax{N}(0,u)
		\lesssim 
		\left\|
			\Fullset_*^{\leq N+1;1} \Psi
		\right\|_{L^2(\Sigma_0^u)}^2
	$.
	The estimate \eqref{E:INITIALSIZEOFL2CONTROLLING} now follows 
	from the initial data assumptions \eqref{E:PSIDATAASSUMPTIONS}.
\end{proof}

The next lemma provides
control over the error integrals corresponding to the deformation tensor
of the multiplier vectorfield \eqref{E:DEFINITIONMULT}, that is,
for the last integral on RHS~\eqref{E:E0DIVID}.
We stress that one of these error integrals is coercive in the geometric torus derivatives
and was treated separately in Lemma~\ref{L:KEYSPACETIMECOERCIVITY}.

\begin{lemma}[\textbf{Error integrals involving the deformation tensor of the multiplier vectorfield}]
	\label{L:MULTIPLIERVECTORFIELDERRORINTEGRALS}
	Assume that $1 \leq N  \leq 18$ and $\varsigma > 0$.
	Let $\basicenergyerrorarg{\Mult}{i}[\Tanset^N \Psi]$
	be the quantities defined by \eqref{E:MULTERRORINTEG1}-\eqref{E:MULTERRORINTEG5}
	(with $\Tanset^N \Psi$ in the role of $\Psi$).
	Under the data-size and bootstrap assumptions 
	of Subsects.~\ref{SS:SIZEOFTBOOT}-\ref{SS:PSIBOOTSTRAP}
	and the smallness assumptions of Subsect.~\ref{SS:SMALLNESSASSUMPTIONS}, 
	the following integral estimates hold for $(t,u) \in [0,\Tboot) \times [0,U_0]$,
	where the implicit constants are independent of $\varsigma$
	(and without any absolute value taken on the left):
	\begin{align}
		\int_{\mathcal{M}_{t,u}}
			\sum_{i=1}^5 \basicenergyerrorarg{\Mult}{i}[\Tanset^N \Psi]
		\, d \vol
		& \lesssim
		 	\int_{t'=0}^t
				\frac{1}{\sqrt{\Tboot - t'}} \totTanmax{[1,N]}(t',u)
			 \, dt'
			 \\
		& \ \
			+
		 	(1 + \varsigma^{-1})
		 		\int_{t'=0}^t
					\totTanmax{[1,N]}(t',u)
				\, dt'
			 \notag	\\
		 & \ \ 
		 	+
			(1 + \varsigma^{-1})
			\int_{u'=0}^u
				\totTanmax{[1,N]}(t,u')
			\, du'
			+ \varsigma 
				\coerciveTanspacetimemax{[1,N]}(t,u).
				\notag
	\end{align}
\end{lemma}

\begin{proof}
	We integrate \eqref{E:MULTIPLIERVECTORFIEDERRORTERMPOINTWISEBOUND} 
	(with $\Tanset^N \Psi$ in the role of $\Psi$)
	over $\mathcal{M}_{t,u}$ and use 
	Lemmas~\ref{L:KEYSPACETIMECOERCIVITY} and \ref{L:COERCIVENESSOFCONTROLLING}.
\end{proof}

The next lemma yields control over 
the simplest energy estimate error integrals
generated by the commutator terms.
These terms appear in the first error integral on RHS~\eqref{E:E0DIVID},
where $\Tanset^N \Psi $ is in the role of $\Psi$
and $\mathfrak{F}$ is the inhomogeneous term in the wave equation
$\upmu \square_{g(\Psi)} (\Tanset^N \Psi) = \mathfrak{F}$.

\begin{lemma}[$L^2$ \textbf{bounds for error integrals involving} $Harmless^{\leq N}$ \textbf{terms}]
	\label{L:STANDARDPSISPACETIMEINTEGRALS}
	Assume that $1 \leq N \leq 18$ and $\varsigma > 0$.
	Recall that the terms $Harmless^{\leq N}$ are defined in
	Def.~\ref{D:HARMLESSTERMS}.
	Under the data-size and bootstrap assumptions 
	of Subsects.~\ref{SS:SIZEOFTBOOT}-\ref{SS:PSIBOOTSTRAP}
	and the smallness assumptions of Subsect.~\ref{SS:SMALLNESSASSUMPTIONS}, 
	the following integral estimates hold for $(t,u) \in [0,\Tboot) \times [0,U_0]$,
	where the implicit constants are independent of $\varsigma$:
	\begin{subequations}
	\begin{align}  \label{E:STANDARDPSISPACETIMEINTEGRALS}
		\int_{\mathcal{M}_{t,u}}
		 	& \left|
					\myarray[(1 + \upmu) \Lunit \Tanset^N \Psi]
						{\Rad \Tanset^N \Psi}
				\right|
				\left|
					Harmless^{\leq N}
				\right|
		 \, d \vol
		   \\
		& \lesssim
		 	(1 + \varsigma^{-1})
		 	\int_{t'=0}^t
				\totTanmax{[1,N]}(t',u)
			\, dt'
			+
			(1 + \varsigma^{-1})
			\int_{u'=0}^u
				\totTanmax{[1,N]}(t,u')
			\, du'
				\notag \\
		& \ \
			+ \varsigma
			  \coerciveTanspacetimemax{[1,N]}(t,u)
			+ \mathring{\upepsilon}^2,
			\notag 
				\\
		\int_{\mathcal{M}_{t,u}}
		 	& \left|
					\angdiff \Tanset^N \Psi
				\right|
				\left|
					Harmless^{\leq N}
				\right|
		 \, d \vol
		  \label{E:ANGULARDERIVATIVESSTANDARDPSISPACETIMEINTEGRALS} \\
		& \lesssim
		 	\int_{t'=0}^t
				\totTanmax{[1,N]}(t',u)
			\, dt'
			+
		 	\int_{u'=0}^u
				\totTanmax{[1,N]}(t,u')
			\, du'
		 	+
		 	\coerciveTanspacetimemax{[1,N]}(t,u)
			+ 
			\mathring{\upepsilon}^2.
			\notag 
	\end{align}
	\end{subequations}
\end{lemma}

\begin{proof}
	See Subsect.~\ref{SS:OFTENUSEDESTIMATES} for some comments on the analysis.
	To prove 
	\eqref{E:STANDARDPSISPACETIMEINTEGRALS}
	and 
	\eqref{E:ANGULARDERIVATIVESSTANDARDPSISPACETIMEINTEGRALS},
	we must estimate the spacetime integrals of various quadratic terms.
	We derive the desired estimates for three representative
	quadratic terms. The remaining terms can be similarly bounded
	and we omit those details.
	We first bound the integral of
	$\left|
			\Lunit \Tanset^N \Psi
		\right|
	\left|
		\GeoAng \Tanset^{\leq N} \Psi
	\right|$.
	Using
	spacetime Cauchy-Schwarz,
	Lemmas~\ref{L:KEYSPACETIMECOERCIVITY} and \ref{L:COERCIVENESSOFCONTROLLING},
	and simple estimates of the form
	$ab \lesssim a^2 + b^2$,
	and separately treating the regions
	$\lbrace \upmu \geq 1/4 \rbrace$ 
	and $\lbrace \upmu < 1/4 \rbrace$
	when bounding the integral of
	$\left|
		\GeoAng \Tanset^{\leq N} \Psi
	\right|^2$,
	we derive the desired estimate as follows:
	\begin{align}
		& \int_{\mathcal{M}_{t,u}}
		 		\left|
					\Lunit \Tanset^N \Psi
				\right|
				\left|
					\GeoAng \Tanset^{\leq N} \Psi
				\right|
		 \, d \vol
		 	\label{E:FIRSTHARMLESSEXAMPLEINTEGRAL} \\
		 & \lesssim
		 \left\lbrace
		 \int_{\mathcal{M}_{t,u}}
		 		\left|
					\Lunit \Tanset^N \Psi
				\right|^2
			\, d \vol
			\right\rbrace^{1/2}
			\left\lbrace
			\int_{\mathcal{M}_{t,u}}
		 		\left|
					\GeoAng \Tanset^{\leq N} \Psi
				\right|^2
			\, d \vol
			\right\rbrace^{1/2}
				\notag \\
			& \lesssim
				(1 + \varsigma^{-1})
				\int_{u'=0}^u
					\int_{\mathcal{P}_{u'}^t}
						\left|
							\Lunit \Tanset^N \Psi
						\right|^2
					\, d \conevol
				\, du'
				\notag \\
		& \ \
				+
				\int_{u'=0}^u
					\int_{\mathcal{P}_{u'}^t}
						\upmu
						\left|
							\angdiff \Tanset^{\leq N} \Psi
						\right|^2
					\, d \conevol
				\, du'
				+
				\varsigma \TranminusdatasizeWithFactor
				\int_{\mathcal{M}_{t,u}}
				\textbf{1}_{\lbrace \upmu < 1/4 \rbrace}
		 		\left|
					\angdiff \Tanset^{\leq N} \Psi
				\right|^2
			\, d \vol
				\notag \\
			& \lesssim
				(1 + \varsigma^{-1})
				\int_{u'=0}^u
					\totTanmax{[1,N]}(t,u')
				\, du'
				+
				\varsigma
				\coerciveTanspacetimemax{[1,N]}(t,u), 
				\notag
	\end{align}
	which is $\lesssim$ RHS~\eqref{E:STANDARDPSISPACETIMEINTEGRALS} as desired.

	As our second example, we bound the integral of
	$\left|
			\Lunit \Tanset^N \Psi
	 \right|
	 \left|
	 		\Tanset_*^{[1,N]} \upmu
	 \right|$.
	Using spacetime Cauchy-Schwarz,
	Lemmas~\ref{L:KEYSPACETIMECOERCIVITY} and \ref{L:COERCIVENESSOFCONTROLLING},
	inequalities
	\eqref{E:LESSSINGULARTERMSMPOINTNINEINTEGRALBOUND}
	and 
	\eqref{E:TANGENGITALEIKONALINTERMSOFCONTROLLING},
	simple estimates of the form
	$ab \lesssim a^2 + b^2$,
	and the fact that
	$\totTanmax{[1,N]}$
	is increasing in its arguments,
	we derive the desired estimate as follows: 
	\begin{align}
		& \int_{\mathcal{M}_{t,u}}
		 		\left|
					\Lunit \Tanset^N \Psi
				\right|
				\left|
					\Tanset_*^{[1,N]} \upmu
				\right|
		 \, d \vol
		 	\label{E:SECONDHARMLESSEXAMPLEINTEGRAL} \\
		 & \lesssim
		 	\int_{u'=0}^u
					\int_{\mathcal{P}_{u'}^t}
						\left|
							\Lunit \Tanset^N \Psi
						\right|^2
					\, d \conevol
			\, du'
			+
			\int_{t'=0}^t
		 			\int_{\Sigma_{t'}^u}
		 				\left|
							\Tanset_*^{[1,N]} \upmu
						\right|^2
					\, d \tvol
			\, dt'
				\notag
					\\ 
		& \lesssim
				\int_{u'=0}^u
					\totTanmax{[1,N]}(t,u')
				\, du'
				+
				\int_{t'=0}^t
					\left\lbrace
						\int_{s=0}^{t'}
							\frac{\totTanmax{[1,N]}^{1/2}(s,u)}{\upmu_{\star}^{1/2}(s,u)}
						\, ds
					\right\rbrace^2
					+ 
					\mathring{\upepsilon}^2
				\, dt'
				\notag
					\\
		& \lesssim
				\int_{u'=0}^u
					\totTanmax{[1,N]}(t,u')
				\, du'
				+
				\int_{t'=0}^t
					\totTanmax{[1,N]}(t',u)
				\, dt'
				+ \mathring{\upepsilon}^2,
				\notag
	\end{align}
	which is $\lesssim$ RHS~\eqref{E:STANDARDPSISPACETIMEINTEGRALS} as desired.

	As our final example, we bound the integral of the product
	$\left|
			\Lunit \Tanset^N \Psi
	 \right|
	 \left|
	 		\Fullset^{N+1;1} \Psi
	 \right|
	$.
	We first recall the following estimate obtained in the second paragraph of 
	the proof of Lemma~\ref{L:EASYL2BOUNDSFOREIKONALFUNCTIONQUANTITIES}:
	$|\Fullset^{N+1;1} \Psi|
	\lesssim
	|\Rad \Tanset^{[1,N]} \Psi|
	+ 
	|\Tanset^{\leq N+1} \Psi|
	+
	\varepsilon |\Fullset_*^{\leq N;1} \GdVar|
	+ \varepsilon |\Tanset_*^{[1,N]} \BadVar|
	$. 
	Thus, we must bound the integral of the four corresponding products
	generated by the RHS of the previous inequality.
	To bound the integral of the first product,
	we argue as in the proof of 
	\eqref{E:FIRSTHARMLESSEXAMPLEINTEGRAL}
	to deduce that
	\begin{align}
		& \int_{\mathcal{M}_{t,u}}
		 		\left|
					\Lunit \Tanset^N \Psi
				\right|
				\left|
					\Rad \Tanset^{[1,N]} \Psi
				\right|
		 \, d \vol
		 	\label{E:FINALHARMLESSEXAMPLEINTEGRAL} \\
		& \lesssim
		 	\int_{u'=0}^u
					\int_{\mathcal{P}_{u'}^t}
						\left|
							\Lunit \Tanset^N \Psi
						\right|^2
					\, d \conevol
			\, du'
			+
			\int_{t'=0}^t
		 			\int_{\Sigma_{t'}^u}
		 				\left|
							\Rad \Tanset^{[1,N]} \Psi
						\right|^2
					\, d \tvol
			\, dt'
				\notag
					\\ 
		& \lesssim
				\int_{u'=0}^u
					\totTanmax{[1,N]}(t,u')
				\, du'
				+
				\int_{t'=0}^t
					\totTanmax{[1,N]}(t',u)
				\, dt',
				\notag
	\end{align}
	which is $\lesssim$ RHS~\eqref{E:STANDARDPSISPACETIMEINTEGRALS} as desired.
	Similar reasoning yields that the integral
	of the second product
	$|\Lunit \Tanset^N \Psi|
	 |\Tanset^{\leq N+1} \Psi|
	$
	is $\lesssim$ RHS~\eqref{E:FIRSTHARMLESSEXAMPLEINTEGRAL} 
	plus RHS~\eqref{E:SECONDHARMLESSEXAMPLEINTEGRAL}  as desired.
	We clarify that the factor $\mathring{\upepsilon}^2$ is generated by the square of RHS~\eqref{E:PSIL2ESTIMATELOSSOFONEDERIVATIVE}.
	Similar reasoning, together with inequalities 
	\eqref{E:TANGENGITALEIKONALINTERMSOFCONTROLLING}
	and
	\eqref{E:ONERADIALEIKONALINTERMSOFCONTROLLING},
	yields that
	the integral
	of the third product
	$	\varepsilon
		|\Lunit \Tanset^N \Psi|
		|\Fullset_*^{\leq N;1} \GdVar|
	$
	and
	the integral
	of the fourth product
	$	\varepsilon
		|\Lunit \Tanset^N \Psi|
	 	|\Tanset_*^{[1,N]} \BadVar|
	$
	are $\lesssim$ RHS~\eqref{E:FIRSTHARMLESSEXAMPLEINTEGRAL}
	plus RHS~\eqref{E:SECONDHARMLESSEXAMPLEINTEGRAL}  as desired.
	We clarify that we have used
	the fact that $\totTanmax{[1,N]}$ is increasing in its arguments
	and the estimate \eqref{E:LESSSINGULARTERMSMPOINTNINEINTEGRALBOUND}
	to bound the time integrals on RHSs
	\eqref{E:TANGENGITALEIKONALINTERMSOFCONTROLLING} 
	and
	\eqref{E:ONERADIALEIKONALINTERMSOFCONTROLLING}
	by $\lesssim \totTanmax{[1,N]}(t,u)$,
	as we did in passing to the last line of \eqref{E:SECONDHARMLESSEXAMPLEINTEGRAL}.

\end{proof}

\subsection{\texorpdfstring{$L^2$}{Square integral} bounds for the difficult top-order error integrals in terms of 
\texorpdfstring{$\totTanmax{[1,N]}$}{the fundamental controlling quantities}}
\label{SS:MOSTDIFFICULTENERGYESTIMATEINTEGRALS}
In the next lemma, we estimate, in the norm $\| \cdot \|_{L^2(\Sigma_t^u)}$, 
the most difficult product that appears in our energy estimates.

\begin{lemma}[$L^2$ \textbf{bound for the most difficult product}]
	\label{L:DIFFICULTTERML2BOUND}
		Assume that $1 \leq N \leq 18$.
		There exists a constant $C > 0$ such that
		under the data-size and bootstrap assumptions 
		of Subsects.~\ref{SS:SIZEOFTBOOT}-\ref{SS:PSIBOOTSTRAP}
		and the smallness assumptions of Subsect.~\ref{SS:SMALLNESSASSUMPTIONS},
		the following $L^2$ estimate holds for the difficult product $(\Rad \Psi) \GeoAng^N \mytr \upchi$ 
		from Prop.~\ref{P:KEYPOINTWISEESTIMATE}
		whenever $(t,u) \in [0,\Tboot) \times [0,U_0]$:
	\begin{align} \label{E:DIFFICULTTERML2BOUND}
		\left\|
			(\Rad \Psi) \GeoAng^N \mytr \upchi
		\right\|_{L^2(\Sigma_t^u)}
		& \leq
			\boxed{2} 
			\frac{\left\| [\Lunit \upmu]_- \right\|_{L^{\infty}(\Sigma_t^u)}} 
						{\upmu_{\star}(t,u)} 
				\totTanmax{[1,N]}^{1/2}(t,u)
					\\
	& \ \ +
			\boxed{4.05}
			\frac{\left\| [\Lunit \upmu]_- \right\|_{L^{\infty}(\Sigma_t^u)}} 
						{\upmu_{\star}(t,u)} 
			\int_{s=0}^t
				\frac{\left\| [\Lunit \upmu]_- \right\|_{L^{\infty}(\Sigma_s^u)}} 
				{\upmu_{\star}(s,u)} 
				\totTanmax{[1,N]}^{1/2}(s,u) 
			\, ds
			\notag \\
		& \ \ +
			C \varepsilon
			\frac{1} 
						{\upmu_{\star}(t,u)} 
			\int_{s=0}^t
				\frac{1} {\upmu_{\star}(s,u)} 
				\totTanmax{[1,N]}^{1/2}(s,u) 
			\, ds
			\notag \\
		& \ \
			+ C
				\frac{1}{\upmu_{\star}(t,u)}
				\int_{s'=0}^t
					\frac{1}{\upmu_{\star}(s',u)}
					\int_{s=0}^{s'}
						\frac{1}{\upmu_{\star}^{1/2}(s,u)}
						\totTanmax{[1,N]}^{1/2}(s,u)
					\, ds
				\, ds'
			\notag  \\
		& \ \
			+ C
				\frac{1}{\upmu_{\star}(t,u)}
				\int_{s=0}^t
					\frac{1}{\upmu_{\star}^{1/2}(s,u)}
					\totTanmax{[1,N]}^{1/2}(s,u)
				\, ds
		\notag  \\
		& \ \
			+ C \frac{1}{\upmu_{\star}^{1/2}(t,u)} \totTanmax{[1,N]}^{1/2}(t,u)
			+ C \frac{1}{\upmu_{\star}^{3/2}(t,u)} \totTanmax{[1,N-1]}^{1/2}(t,u)
				\notag  \\
		&  \ \
			+ C \frac{1}{\upmu_{\star}^{3/2}(t,u)} \mathring{\upepsilon}.
			\notag
	\end{align}

	Furthermore, we have the following less precise estimate:
	\begin{align} \label{E:LESSPRECISEDIFFICULTTERML2BOUND}
		\left\|
			\upmu \GeoAng^N \mytr \upchi
		\right\|_{L^2(\Sigma_t^u)}
		& \lesssim
			\totTanmax{[1,N]}^{1/2}(t,u)
			+
			\int_{s=0}^t
				\frac{1} 
				{\upmu_{\star}(s,u)} 
				\totTanmax{[1,N]}^{1/2}(s,u)
				\, ds
			\\
	& \ \ 
		+ \mathring{\upepsilon} 
				\left\lbrace 
					\ln \upmu_{\star}^{-1}(t,u) 
					+ 
					1 
				\right\rbrace.
				\notag
	\end{align}
\end{lemma}

\begin{proof}
	See Subsect.~\ref{SS:OFTENUSEDESTIMATES} for some comments on the analysis.
	We first prove \eqref{E:DIFFICULTTERML2BOUND}.
	We take the norm $\| \cdot \|_{L^2(\Sigma_t^u)}$
	of both sides of \eqref{E:KEYPOINTWISEESTIMATE}.
	Using \eqref{E:COERCIVENESSOFCONTROLLING},
	we see that the norm of the first term on RHS~\eqref{E:KEYPOINTWISEESTIMATE} is
	$\leq$ the first term on RHS~\eqref{E:DIFFICULTTERML2BOUND}
	as desired. Also using Lemma~\ref{L:L2NORMSOFTIMEINTEGRATEDFUNCTIONS},
	we see that the norm of the second term on RHS~\eqref{E:KEYPOINTWISEESTIMATE} is
	$\leq$ the second term on RHS~\eqref{E:DIFFICULTTERML2BOUND}.
	We now explain why
	the norm $\| \cdot \|_{L^2(\Sigma_t^u)}$
	of the term
	$\left|
			\mbox{\upshape Error}
	\right|
	$
	from \eqref{E:ERRORTERMKEYPOINTWISEESTIMATE}
	is $\leq$ the sum of the 
	terms on lines three to seven of RHS~\eqref{E:DIFFICULTTERML2BOUND}.
	With the exception of the bound for the first term
	$
		\displaystyle
		\frac{1}{\upmu_{\star}(t,u)}
			\left|\upchifullmodarg{\GeoAng^N} \right|(0,u,\vartheta)
	$
	on RHS~\eqref{E:ERRORTERMKEYPOINTWISEESTIMATE},
	the desired bounds follow from
	the same estimates used above
	together with those of Lemma~\ref{L:EASYL2BOUNDSFOREIKONALFUNCTIONQUANTITIES},
	Cor.~\ref{C:COMMUTATIONBOUNDSFORPSIINTERMSOFFUNDAMENTALCONTROLLINGQUANTITIES},
	inequalities
	\eqref{E:LOSSKEYMUINTEGRALBOUND},
	\eqref{E:LOGLOSSMUINVERSEINTEGRALBOUND},
	and
	\eqref{E:LESSSINGULARTERMSMPOINTNINEINTEGRALBOUND},
	the fact that $\totTanmax{[1,N]}$ is increasing in its arguments, 
	and simple inequalities of the form
	$ab \lesssim a^2 + b^2$.
Finally, we must bound 
$ \displaystyle
\left\| 
	\frac{1}{\upmu_{\star}(t,\cdot)}
	\left|\upchifullmodarg{\GeoAng^N} \right|(0,\cdot) 
\right\|_{L^2(\Sigma_t^u)}
$.
We first use 
\eqref{E:SIGMATVOLUMEFORMCOMPARISON}
with $s=0$
to deduce
$
\displaystyle
\left\| 
	\frac{1}{\upmu_{\star}(t,\cdot)}
	\left|\upchifullmodarg{\GeoAng^N} \right|(0,\cdot) 
\right\|_{L^2(\Sigma_t^u)}
\lesssim
\frac{1}{\upmu_{\star}(t,u)}
\left\|
	\upchipartialmodarg{\GeoAng^{N-1}}
\right\|_{L^2(\Sigma_0^u)}
$.
We now use definition \eqref{E:TRANSPORTPARTIALRENORMALIZEDTRCHIJUNK},
the simple inequality
$|G_{(Frame)}| = |\smoothfunction(\GdVar,\angdiff x^1,\angdiff x^2)| \lesssim 1$
(which follows from Lemmas~\ref{L:SCHEMATICDEPENDENCEOFMANYTENSORFIELDS}
and \ref{L:POINTWISEFORRECTANGULARCOMPONENTSOFVECTORFIELDS}
and the $L^{\infty}$ estimates of Prop.~\ref{P:IMPROVEMENTOFAUX}),
the estimates of Lemma~\ref{L:BEHAVIOROFEIKONALFUNCTIONQUANTITIESALONGSIGMA0},
the estimate \eqref{E:POINTWISEESTIMATESFORGSPHEREANDITSTANGENTIALDERIVATIVES},
and the assumptions on the data
to deduce the desired bound
$
\displaystyle
\frac{1}{\upmu_{\star}(t,u)}
\left\| 
	\upchifullmodarg{\GeoAng^N}
\right\|_{L^2(\Sigma_0^u)}
\lesssim \frac{1}{\upmu_{\star}(t,u)} \mathring{\upepsilon}
$.
We have thus proved \eqref{E:DIFFICULTTERML2BOUND}.

	The proof of \eqref{E:LESSPRECISEDIFFICULTTERML2BOUND}
	is based on inequality \eqref{E:LESSPRECISEKEYPOINTWISEESTIMATE}
	and is similar but much simpler; we omit the details, 
	noting only that 
	we use \eqref{E:COMMUTATIONBOUNDSFORPSIINTERMSOFFUNDAMENTALCONTROLLINGQUANTITIES}
	to bound the order $\leq N$ derivatives of $\Psi$ on RHS~\eqref{E:LESSPRECISEKEYPOINTWISEESTIMATE}
	and that
	inequality \eqref{E:LOGLOSSMUINVERSEINTEGRALBOUND}
	leads to the presence of the factor 
	$
		\ln \upmu_{\star}^{-1}(t,u) + 1 
	$.
\end{proof}

\subsection{\texorpdfstring{$L^2$}{Square integral} bounds for less degenerate top-order error integrals in terms of 
\texorpdfstring{$\totTanmax{[1,N]}$}{the fundamental controlling quantities}}
\label{SS:LESSDEGENERATEENERGYESTIMATEINTEGRALS}
In the next lemma, we bound some top-order error integrals
that appear in our energy estimates.
As in the proof of Lemma~\ref{L:DIFFICULTTERML2BOUND},  
we need to use the modified quantities to avoid
losing a derivative. However, 
the estimates of the lemma are much less degenerate than those of
Lemma~\ref{L:DIFFICULTTERML2BOUND} because of the availability
of a helpful factor of $\upmu$ in the integrands.

\begin{lemma}[\textbf{Bounds for less degenerate top-order error integrals}]
	\label{L:LESSDEGENERATEENERGYESTIMATEINTEGRALS}
	Assume that $1 \leq N \leq 18$.
	Under the data-size and bootstrap assumptions 
	of Subsects.~\ref{SS:SIZEOFTBOOT}-\ref{SS:PSIBOOTSTRAP}
	and the smallness assumptions of Subsect.~\ref{SS:SMALLNESSASSUMPTIONS}, 
	the following integral estimates hold for $(t,u) \in [0,\Tboot) \times [0,U_0]$:
	\begin{subequations}
	\begin{align} \label{E:FIRSTLESSDEGENERATEENERGYESTIMATEINTEGRALS}
		&
		\left|
			\int_{\mathcal{M}_{t,u}}
				\GeoAngFlatRadComponent
				(\Rad \GeoAng^N \Psi)	 
				(\Rad \Psi) 
				(\angdiffuparg{\#} \Psi)
				\cdot
				(\upmu \angdiff \GeoAng^{N-1} \mytr \upchi)
			\, d \vol
		\right|
			\\
		& \lesssim
			\int_{t'=0}^t
				\left\lbrace 
					\ln \upmu_{\star}^{-1}(t',u) 
					+ 
					1 
				\right\rbrace^2
				\totTanmax{[1,N]}(t',u)
			\, dt'
			+ \mathring{\upepsilon}^2,
			\notag
			\\
	&
		\left|
			\int_{\mathcal{M}_{t,u}}
				(1 + 2 \upmu)
				\GeoAngFlatRadComponent 
				(\Lunit \GeoAng^N \Psi)	 
				(\Rad \Psi)
				(\angdiffuparg{\#} \Psi)
				\cdot
				(\upmu \angdiff \GeoAng^{N-1} \mytr \upchi)
			\, d \vol
		\right|
			\label{E:SECONDLESSDEGENERATEENERGYESTIMATEINTEGRALS} \\
		& \lesssim
			\int_{t'=0}^t
				\left\lbrace 
					\ln \upmu_{\star}^{-1}(t',u) 
					+ 
					1 
				\right\rbrace^2
				\totTanmax{[1,N]}(t',u)
			\, dt'
			+
			\int_{u'=0}^u
				\totTanmax{[1,N]}(t,u')
			\, du'
				\notag \\
		& \ \
			+ \mathring{\upepsilon}^2.
			\notag
	\end{align}
	\end{subequations}
\end{lemma}

\begin{proof}
See Subsect.~\ref{SS:OFTENUSEDESTIMATES} for some comments on the analysis.
To prove \eqref{E:SECONDLESSDEGENERATEENERGYESTIMATEINTEGRALS},
we use the fact that $\GeoAngFlatRadComponent = \smoothfunction(\GdVar) \GdVar$
 (see \eqref{E:LINEARLYSMALLSCALARSDEPENDINGONGOODVARIABLES}),
	the $L^{\infty}$ estimates of Prop.~\ref{P:IMPROVEMENTOFAUX},
	Cauchy-Schwarz, 
	and \eqref{E:COERCIVENESSOFCONTROLLING}
	to deduce
\begin{align} \label{E:FIRSTBOUNDFORSECONDLESSDEGENERATEENERGYESTIMATEINTEGRALS}
	\mbox{LHS~\eqref{E:SECONDLESSDEGENERATEENERGYESTIMATEINTEGRALS}}
	& \lesssim
		\int_{\mathcal{M}_{t,u}}
			\left|
				\Lunit \GeoAng^N \Psi
			\right|^2 
		\, d \vol
		+
		\int_{\mathcal{M}_{t,u}}
			\left|
				\upmu \GeoAng^N \mytr \upchi
			\right|^2
		\, d \vol
		\\
		& 
		\lesssim
		\int_{u'=0}^u
			\left\|
				\Lunit \GeoAng^N \Psi
			\right\|_{L^2(\mathcal{P}_{u'}^t)}^2
		\, du'
		+
			\int_{t' = 0}^t
				\left\|
					\upmu \GeoAng^N \mytr \upchi 
				\right\|_{L^2(\Sigma_{t'}^u)}^2
			\, dt'
			\notag \\
		&
		\lesssim
			\int_{u'=0}^u
				\totTanmax{[1,N]}(t,u')
			\, du'
			+
			\int_{t' = 0}^t
				\left\|
					\upmu \GeoAng^N \mytr \upchi 
				\right\|_{L^2(\Sigma_{t'}^u)}^2
			\, dt'.
		\notag
\end{align}
To complete the proof of \eqref{E:SECONDLESSDEGENERATEENERGYESTIMATEINTEGRALS},
we must handle the final integral on RHS~\eqref{E:FIRSTBOUNDFORSECONDLESSDEGENERATEENERGYESTIMATEINTEGRALS}.
To bound the integral by $\leq$ RHS~\eqref{E:SECONDLESSDEGENERATEENERGYESTIMATEINTEGRALS}
we use inequality \eqref{E:LESSPRECISEDIFFICULTTERML2BOUND} 
(with $t'$ in place of $t$),
simple estimates of the form $ab \lesssim a^2 + b^2$,
and we in addition use \eqref{E:LOGLOSSMUINVERSEINTEGRALBOUND} 
and the fact that $\totTanmax{[1,N]}$ is increasing in its arguments 
to bound the time integral on RHS~\eqref{E:LESSPRECISEDIFFICULTTERML2BOUND} as follows:
\[
	\int_{s=0}^{t'}
				\frac{1} 
				{\upmu_{\star}(s,u)} 
				\totTanmax{[1,N]}^{1/2}(s,u)
	\, ds
	\lesssim
	\left\lbrace 
		\ln \upmu_{\star}^{-1}(t',u) 
		+ 
		1 
	\right\rbrace
	\totTanmax{[1,N]}^{1/2}(t',u).
\]
In carrying out this procedure, we encounter the following integral
generated by the next-to-last term on RHS~\eqref{E:LESSPRECISEDIFFICULTTERML2BOUND}:
\[
	\mathring{\upepsilon}^2
	\int_{t'=0}^t
		\left\lbrace 
					\ln \upmu_{\star}^{-1}(t',u) 
					+ 
					1 
				\right\rbrace^2
	\, dt'.
\]
Using \eqref{E:LESSSINGULARTERMSMPOINTNINEINTEGRALBOUND},
we deduce that the above term is $\lesssim \mathring{\upepsilon}^2$ as desired.
We have thus proved \eqref{E:SECONDLESSDEGENERATEENERGYESTIMATEINTEGRALS}.

The proof of \eqref{E:FIRSTLESSDEGENERATEENERGYESTIMATEINTEGRALS}
starts with the following analog of \eqref{E:FIRSTBOUNDFORSECONDLESSDEGENERATEENERGYESTIMATEINTEGRALS},
which can proved in the same way:
\[
\mbox{LHS~\eqref{E:FIRSTLESSDEGENERATEENERGYESTIMATEINTEGRALS}}
\lesssim
\int_{t' = 0}^t 
			\left\| 
				\Rad \GeoAng^N \Psi 
			\right\|_{L^2(\Sigma_{t'}^u)}^2 
		\, dt'
+
\int_{t' = 0}^t 
	\left\|
		\upmu \GeoAng^N \mytr \upchi 
	\right\|_{L^2(\Sigma_{t'}^u)}^2 
		\, dt'.
\]
The remaining details are similar to those given in the proof
of \eqref{E:SECONDLESSDEGENERATEENERGYESTIMATEINTEGRALS};
we therefore omit them.
\end{proof}

\subsection{Error integrals requiring integration by parts with respect to \texorpdfstring{$\Lunit$}{the rescaled null vectorfield}}
\label{SS:ERROINTEGRALSINVOLVINGIBPL}
In deriving top-order energy estimates, we encounter the error integral
\[
-           \int_{\mathcal{M}_{t,u}}
							(1 + 2 \upmu)
							(\Lunit \GeoAng^N \Psi)	 
							(\Rad \Psi) \GeoAng^N \mytr \upchi
							\, d \vol.
\]
It turns out that to suitably bound it, we must 
rely on the partially modified quantity 
$\upchipartialmodarg{\GeoAng^{N-1}}$ defined in \eqref{E:TRANSPORTPARTIALRENORMALIZEDTRCHIJUNK},
and we must also integrate by parts
via the identity \eqref{E:LUNITANDANGULARIBPIDENTITY}.
We derive the main estimate of interest for the above error integral in Lemma~\ref{L:BOUNDSFORDIFFICULTTOPORDERINTEGRALSINVOLVINGLUNITIBP}.
Before proving the lemma, we first establish some preliminary estimates for
various error integrals that arise from the integration by parts procedure. 
We bound the most difficult of these integrals, 
which is a $\Sigma_t^u$ boundary integral, 
in Lemma~\ref{L:ANNOYINGBOUNDRYSPATIALINTEGRALFACTORL2ESTIMATE}.

We start by deriving $\| \cdot \|_{L^2(\Sigma_t^u)}$ estimates
for the two second most difficult products that appear in our energy estimates.

\begin{lemma}[\textbf{A difficult hypersurface} $L^2$ \textbf{estimate}]
	\label{L:ANNOYINGBOUNDRYSPATIALINTEGRALFACTORL2ESTIMATE}
	Assume that $1 \leq N \leq 18$.
	Let $\upchipartialmodarg{\GeoAng^{N-1}}$
	be the partially modified quantity defined by \eqref{E:TRANSPORTPARTIALRENORMALIZEDTRCHIJUNK}.
	Under the data-size and bootstrap assumptions 
	of Subsects.~\ref{SS:SIZEOFTBOOT}-\ref{SS:PSIBOOTSTRAP}
	and the smallness assumptions of Subsect.~\ref{SS:SMALLNESSASSUMPTIONS},
	the following $L^2$ estimate holds for $(t,u) \in [0,\Tboot) \times [0,U_0]$:
	\begin{subequations}
	\begin{align} \label{E:ANNOYINGLDERIVATIVEBOUNDRYSPATIALINTEGRALFACTORL2ESTIMATE}
		\left\|
			\frac{1}{\sqrt{\upmu}} (\Rad \Psi) \Lunit \upchipartialmodarg{\GeoAng^{N-1}}
		\right\|_{L^2(\Sigma_t^u)}
		& \leq
			\boxed{\sqrt{2}}
			\frac{\left\| [\Lunit \upmu]_- \right\|_{C^0(\Sigma_t^u)}}
			{\upmu_{\star}(t,u)}
			\totTanmax{[1,N]}^{1/2}(t,u) 
			\\
		& \ \
			+ C \frac{1}{\upmu_{\star}^{1/2}(t,u)} \totTanmax{[1,N]}^{1/2}(t,u)
			+ C \varepsilon \frac{1}{\upmu_{\star}(t,u)} \totTanmax{[1,N]}^{1/2}(t,u)
				\notag \\
		& \ \
			+ C \mathring{\upepsilon} \frac{1}{\upmu_{\star}^{1/2}(t,u)},
			\notag
			\\
		\left\|
			\frac{1}{\sqrt{\upmu}} (\Rad \Psi) \upchipartialmodarg{\GeoAng^{N-1}}
		\right\|_{L^2(\Sigma_t^u)}
		& \leq
			\boxed{\sqrt{2}}
			\left\| \Lunit \upmu \right\|_{L^{\infty}(\Sigmaminus{t}{t}{u})}
			\frac{1}{\upmu_{\star}^{1/2}(t,u)}
			\int_{t'=0}^t
				\frac{1}{\upmu_{\star}^{1/2}(t',u)} \totTanmax{[1,N]}^{1/2}(t',u)
			\, dt'
				\label{E:ANNOYINGBOUNDRYSPATIALINTEGRALFACTORL2ESTIMATE} \\
		& \ \
			+ 
			C 
			\int_{t'=0}^t
				\frac{1}{\upmu_{\star}^{1/2}(t',u)} \totTanmax{[1,N]}^{1/2}(t',u)
			\, dt'
				\notag \\
		& \ \
			+ 
			C \varepsilon
			\frac{1}{\upmu_{\star}^{1/2}(t,u)}
			\int_{t'=0}^t
				\frac{1}{\upmu_{\star}^{1/2}(t',u)} \totTanmax{[1,N]}^{1/2}(t',u)
			\, dt'
			\notag \\
	& \ \
		+ C \mathring{\upepsilon}
			\frac{1}{\upmu_{\star}^{1/2}(t,u)}.
			\notag
	\end{align}
	\end{subequations}

	Moreover, we have the following less precise estimates:
	\begin{subequations}
	\begin{align}  \label{E:NOTTOOHARDLUNITAPPLIEDTOBOUNDRYSPATIALINTEGRALFACTORL2ESTIMATE}
		\left\|
			\Lunit \upchipartialmodarg{\GeoAng^{N-1}}
		\right\|_{L^2(\Sigma_t^u)}
		& \lesssim
			\frac{1}{\upmu_{\star}^{1/2}(t,u)} 
			\totTanmax{[1,N]}^{1/2}(t,u)
			+
			\mathring{\upepsilon},
				\\
		\left\|
			\upchipartialmodarg{\GeoAng^{N-1}}
		\right\|_{L^2(\Sigma_t^u)}
		& \lesssim
			\int_{t'=0}^t
				\frac{1}{\upmu_{\star}^{1/2}(t',u)} \totTanmax{[1,N]}^{1/2}(t',u)
			\, dt'
			+
			\mathring{\upepsilon}.
		\label{E:MUCHEASIERBOUNDRYSPATIALINTEGRALFACTORL2ESTIMATE}
	\end{align}
	\end{subequations}
\end{lemma}

\begin{proof}
	See Subsect.~\ref{SS:OFTENUSEDESTIMATES} for some comments on the analysis.
	We first prove \eqref{E:ANNOYINGBOUNDRYSPATIALINTEGRALFACTORL2ESTIMATE}.
	We take the norm $\| \cdot \|_{L^2(\Sigma_t^u)}$ of 
	$
	\displaystyle
	\frac{1}{\sqrt{\upmu}} (\Rad \Psi)
	$ 
	times \eqref{E:SHARPPOINTWISEPARTIALLYMODIFIED}.
	We bound the terms arising from the second line of RHS~\eqref{E:SHARPPOINTWISEPARTIALLYMODIFIED} 
	by $\leq$ the sum of the last two terms on RHS~\eqref{E:ANNOYINGBOUNDRYSPATIALINTEGRALFACTORL2ESTIMATE}
	with the help of
	Lemma~\ref{L:L2NORMSOFTIMEINTEGRATEDFUNCTIONS},
	Lemma~\ref{L:COERCIVENESSOFCONTROLLING},
	Lemma~\ref{L:EASYL2BOUNDSFOREIKONALFUNCTIONQUANTITIES},
	and the estimate $\| \Rad \Psi \|_{L^{\infty}(\Sigma_t^u)} \lesssim 1$
	(that is, \eqref{E:PSITRANSVERSALLINFINITYBOUNDBOOTSTRAPIMPROVED}).

	To bound the norm $\| \cdot \|_{L^2(\Sigma_t^u)}$ of the product
	$
	\displaystyle
	\frac{1}{2}
	\frac{1}{\sqrt{\upmu}} (\Rad \Psi)
	\left| G_{\Lunit \Lunit} \right|
	\int_{t'=0}^t
	\left|
		\angLap \GeoAng^{N-1} \Psi
	\right|
	\, dt'
	$,
	we first use equation \eqref{E:UPMUFIRSTTRANSPORT},
	the relations
	$G_{\Lunit \Lunit}, G_{\Lunit \Radunit} = \smoothfunction(\GdVar)$
	(see Lemma~\ref{L:SCHEMATICDEPENDENCEOFMANYTENSORFIELDS}),
	inequality \eqref{E:ANGDERIVATIVESINTERMSOFTANGENTIALCOMMUTATOR},
	the $L^{\infty}$ estimates of Prop.~\ref{P:IMPROVEMENTOFAUX},
	and Cor.~\ref{C:SQRTEPSILONTOCEPSILON}
	to pointwise bound the product by
	$
	\displaystyle
	\leq
	(1 + C \varepsilon)
	\left\lbrace
		\frac{|\Lunit \upmu(t,u,\vartheta)|}{\sqrt{\upmu(t,u,\vartheta)}}
		+
		C \varepsilon
	\right\rbrace
	\int_{t'=0}^t
		\left|
			\angdiff \GeoAng^{\leq N} \Psi
		\right|
		(t',u,\vartheta)
	\, dt'
	$.
	As above, we can bound the product involving the factor 
	$
	\displaystyle
	C \varepsilon
	$ 
	by $\leq$ the next-to-last term on RHS~\eqref{E:ANNOYINGBOUNDRYSPATIALINTEGRALFACTORL2ESTIMATE}
	by using
	Lemmas~\ref{L:L2NORMSOFTIMEINTEGRATEDFUNCTIONS}
	and
	\ref{L:COERCIVENESSOFCONTROLLING}.
	To bound the remaining (difficult) term
	$
	\displaystyle
	(1 + C \varepsilon)
	\left\|
	\left|
		\frac{\Lunit \upmu}{\sqrt{\upmu}}
	\right|
	\int_{t'=0}^t
		\left|
			\angdiff \GeoAng^{\leq N} \Psi
		\right|
		\, dt'
	\right\|_{L^2(\Sigma_t^u)}
	$,
	we first decompose $\Sigma_t^u = \Sigmaplus{t}{t}{u} \cup \Sigmaminus{t}{t}{u}$ 
	as in Def.~\ref{D:REGIONSOFDISTINCTUPMUBEHAVIOR}
	and use Lemmas~\ref{L:L2NORMSOFTIMEINTEGRATEDFUNCTIONS} 
	and \ref{L:COERCIVENESSOFCONTROLLING}
	to bound it by
	$\leq (1 + C \varepsilon)$ times
	\begin{align} \label{E:DIFFICULTTERMSANNOYINGBOUNDRYSPATIALINTEGRALFACTORL2ESTIMATE}
	&
	\boxed{\sqrt{2}}
	\left\|
		\frac{\Lunit \upmu}{\sqrt{\upmu}}
	\right\|_{L^{\infty}(\Sigmaminus{t}{t}{u})}
	\int_{t'=0}^t
		\frac{1}{\upmu_{\star}^{1/2}(t',u)}
		\totTanmax{[1,N]}^{1/2}(t',u)
	\, dt'
		\\
	& \ \
		+ C
		\left\|
			\frac{\Lunit \upmu}{\sqrt{\upmu}}
		\right\|_{L^{\infty}(\Sigmaplus{t}{t}{u})}
		\int_{t'=0}^t
			\frac{1}{\upmu_{\star}^{1/2}(t',u)}
			\totTanmax{[1,N]}^{1/2}(t',u)
		\, dt'.
		\notag
\end{align}
The fact that 
$
(1 + C \varepsilon) 
\times
\mbox{RHS~\eqref{E:DIFFICULTTERMSANNOYINGBOUNDRYSPATIALINTEGRALFACTORL2ESTIMATE}} 
\leq 
\mbox{RHS~\eqref{E:ANNOYINGBOUNDRYSPATIALINTEGRALFACTORL2ESTIMATE}}
$
follows from using 
\eqref{E:LUNITUPMULINFINITY},
\eqref{E:UPMULINFTY},
\eqref{E:POSITIVEPARTOFLMUOVERMUISBOUNDED},
and \eqref{E:KEYMUNOTDECAYINGMINUSPARTLMUOVERMUBOUND}
to deduce that
$
\displaystyle
(1 + C \varepsilon) 
\left\|
	\frac{\Lunit \upmu}{\sqrt{\upmu}}
\right\|_{L^{\infty}(\Sigmaminus{t}{t}{u})}
\leq
\left\|
	\Lunit \upmu
\right\|_{L^{\infty}(\Sigmaminus{t}{t}{u})}
\upmu_{\star}^{-1/2}(t,u)
+
C \varepsilon \upmu_{\star}^{-1/2}(t,u)
$
and
$
\displaystyle
\left\|
			\frac{\Lunit \upmu}{\sqrt{\upmu}}
		\right\|_{L^{\infty}(\Sigmaplus{t}{t}{u})}
\leq C 
$.

Finally, we must bound the norm $\| \cdot \|_{L^2(\Sigma_t^u)}$ of the product
arising from the first term $\left|
			\upchipartialmodarg{\GeoAng^{N-1}}
		\right|
		(0,\cdot)$ on RHS~\eqref{E:SHARPPOINTWISEPARTIALLYMODIFIED}.
We first use 
\eqref{E:PSITRANSVERSALLINFINITYBOUNDBOOTSTRAPIMPROVED}
and
\eqref{E:SIGMATVOLUMEFORMCOMPARISON}
with $s=0$
to deduce
$
	\displaystyle
	\left\|
	\frac{1}{\sqrt{\upmu}} (\Rad \Psi)
	\left|
		\upchipartialmodarg{\GeoAng^{N-1}}
	\right|
	(0,\cdot)
\right\|_{L^2(\Sigma_t^u)}
\lesssim
\frac{1}{\upmu_{\star}^{1/2}(t,u)}
\left\|
	\upchipartialmodarg{\GeoAng^{N-1}}
\right\|_{L^2(\Sigma_0^u)}
$.
Next, from definition \eqref{E:TRANSPORTPARTIALRENORMALIZEDTRCHIJUNK},
the simple inequality
$|G_{(Frame)}| = |\smoothfunction(\GdVar,\angdiff x^1,\angdiff x^2)| \lesssim 1$
(which follows from Lemmas~\ref{L:SCHEMATICDEPENDENCEOFMANYTENSORFIELDS}
and \ref{L:POINTWISEFORRECTANGULARCOMPONENTSOFVECTORFIELDS}
and the $L^{\infty}$ estimates of Prop.~\ref{P:IMPROVEMENTOFAUX}),
the estimates of Lemma~\ref{L:BEHAVIOROFEIKONALFUNCTIONQUANTITIESALONGSIGMA0},
the estimate \eqref{E:POINTWISEESTIMATESFORGSPHEREANDITSTANGENTIALDERIVATIVES},
and the assumptions on the data,
we find that 
$
\left\|
	\upchipartialmodarg{\GeoAng^{N-1}}
\right\|_{L^2(\Sigma_0^u)}
\lesssim \mathring{\upepsilon}
$.
In total, we conclude that the product under consideration 
is $\lesssim$ the last term on RHS~\eqref{E:ANNOYINGBOUNDRYSPATIALINTEGRALFACTORL2ESTIMATE} 
as desired. We have thus proved \eqref{E:ANNOYINGBOUNDRYSPATIALINTEGRALFACTORL2ESTIMATE}.

To prove \eqref{E:ANNOYINGLDERIVATIVEBOUNDRYSPATIALINTEGRALFACTORL2ESTIMATE},
we take the norm $\| \cdot \|_{L^2(\Sigma_t^u)}$ of
$
	\displaystyle
	\frac{1}{\sqrt{\upmu}} (\Rad \Psi)
$ 
times \eqref{E:SHARPPOINTWISELUNITPARTIALLYMODIFIED}.
We bound the terms arising from the last two terms on RHS~\eqref{E:SHARPPOINTWISELUNITPARTIALLYMODIFIED}
by $\leq$ the last two terms on RHS~\eqref{E:ANNOYINGBOUNDRYSPATIALINTEGRALFACTORL2ESTIMATE}
with the help of Lemma~\ref{L:COERCIVENESSOFCONTROLLING},
and the estimates
$\| \Rad \Psi \|_{L^{\infty}(\Sigma_t^u)} \lesssim 1$,
\eqref{E:LESSSINGULARTERMSMPOINTNINEINTEGRALBOUND},
\eqref{E:TANGENGITALEIKONALINTERMSOFCONTROLLING},
and \eqref{E:ONERADIALEIKONALINTERMSOFCONTROLLING}.
Note that we have used 
\eqref{E:LESSSINGULARTERMSMPOINTNINEINTEGRALBOUND}
and the fact that the $\totTanmax{[1,N]}$ are increasing in their arguments
to bound the time integrals 
on RHS~\eqref{E:TANGENGITALEIKONALINTERMSOFCONTROLLING}-\eqref{E:ONERADIALEIKONALINTERMSOFCONTROLLING}
by $\lesssim \totTanmax{[1,N]}(t,u)$.
To bound the norm $\| \cdot \|_{L^2(\Sigma_t^u)}$ of the product
$
\displaystyle
\frac{1}{2} 
\frac{1}{\sqrt{\upmu}} (\Rad \Psi)
\left| 
			G_{\Lunit \Lunit} 
		\right|
			\left|
					\angLap \GeoAng^{N-1} \Psi
			\right|
$,
we first use the reasoning from the second paragraph of this proof
to pointwise bound the product by
$
	\displaystyle
	\leq
	\frac{1}{\upmu_{\star}(t,u)}
	\left\|
		[\Lunit \upmu]_-
	\right\|_{L^{\infty}(\Sigma_t^u)}
	\left|
		\sqrt{\upmu} \angdiff \GeoAng^{\leq N} \Psi
	\right|
	+
	C
	\left\|
		\frac{[\Lunit \upmu]_+}{\upmu}
	\right\|_{L^{\infty}(\Sigma_t^u)}
	\left|
		\sqrt{\upmu} \angdiff \GeoAng^{\leq N} \Psi
	\right|
	+
	\frac{C \varepsilon}{\upmu_{\star}(t,u)}
	\left|
		\sqrt{\upmu} \angdiff \GeoAng^{\leq N} \Psi
	\right|
	$.
	Thus, using inequality \eqref{E:POSITIVEPARTOFLMUOVERMUISBOUNDED} and
	Lemma~\ref{L:COERCIVENESSOFCONTROLLING}, 
	we bound the norm $\| \cdot \|_{L^2(\Sigma_t^u)}$ of these products 
	by $\leq$ the sum of the first, second, and third 
	terms on RHS~\eqref{E:ANNOYINGLDERIVATIVEBOUNDRYSPATIALINTEGRALFACTORL2ESTIMATE}.

The proofs of 
\eqref{E:NOTTOOHARDLUNITAPPLIEDTOBOUNDRYSPATIALINTEGRALFACTORL2ESTIMATE}
and \eqref{E:MUCHEASIERBOUNDRYSPATIALINTEGRALFACTORL2ESTIMATE} are based on
a strict subset of the above arguments and are much simpler;
we omit the details.
\end{proof}

We now derive estimates for some error integrals that
are much easier to estimate than the ones
treated in Lemma~\ref{L:ANNOYINGBOUNDRYSPATIALINTEGRALFACTORL2ESTIMATE}.

\begin{lemma}[\textbf{Bounds connected to easy top-order error integrals requiring integration by parts with respect to} $\Lunit$]
	\label{L:HARMLESSIBPERRORINTEGRALS}
	Assume that $1 \leq N \leq 18$ and $\varsigma > 0$.
	Let 
	$\mbox{\upshape Error}_i[\GeoAng^N \Psi; \upchipartialmodarg{\GeoAng^{N-1}}]$
	be the error integrands defined in
	\eqref{E:LUNITIBPSPACETIMEERROR}
	and
	\eqref{E:LUNITIBPHYPERSURFACEERROR},
	where
	$\GeoAng^N \Psi$ is in the role of $\Tanset^N \Psi$
	and the partially modified quantity
	$\upchipartialmodarg{\GeoAng^{N-1}}$ 
	defined in \eqref{E:TRANSPORTPARTIALRENORMALIZEDTRCHIJUNK} is in role of $\ThirdSmoothFunction$.
	Under the data-size and bootstrap assumptions 
	of Subsects.~\ref{SS:SIZEOFTBOOT}-\ref{SS:PSIBOOTSTRAP}
	and the smallness assumptions of Subsect.~\ref{SS:SMALLNESSASSUMPTIONS},
	the following estimates hold for $(t,u) \in [0,\Tboot) \times [0,U_0]$,
	where the implicit constants are independent of $\varsigma$:
	\begin{subequations}
	\begin{align} \label{E:HARMLESSIBPSPACETIMEERRORINTEGRALS}
		\int_{\mathcal{M}_{t,u}}
			& 
			\left|
				\mbox{\upshape Error}_1[\GeoAng^N \Psi; \upchipartialmodarg{\GeoAng^{N-1}}]
			\right|
		\, d \vol
			\\
		& \lesssim 
			(1 + \varsigma^{-1})
			\int_{s=0}^t
				\frac{1}{\upmu_{\star}^{1/2}(s,u)} 
				\totTanmax{[1,N]}(s,u)
			\, ds
			+
			\int_{s=0}^t
				\frac{1}{\upmu_{\star}^{3/2}(s,u)} 
				\totTanmax{[1,N-1]}(s,u)
			\, ds
			\notag \\
		& \ \
			+
			\varsigma \coerciveTanspacetimemax{[1,N]}(t,u)
			+ (1 + \varsigma^{-1}) \mathring{\upepsilon}^2,
				\notag 
		\end{align}
		\begin{align}
		\int_{\Sigma_t^u}
		\left|
			\mbox{\upshape Error}_2[\GeoAng^N \Psi; \upchipartialmodarg{\GeoAng^{N-1}}]
		\right|
		\, d \vol
		& \lesssim 
			\mathring{\upepsilon}^2
			+ \varepsilon \totTanmax{[1,N]}(t,u)
			\label{E:HARMLESSIBPHYPERSURFACEERRORINTEGRALS},
			\\
		\int_{\Sigma_0^u}
		\left|
			\mbox{\upshape Error}_2[\GeoAng^N \Psi; \upchipartialmodarg{\GeoAng^{N-1}}]
		\right|
		\, d \vol
		& \lesssim \mathring{\upepsilon}^2,
			\label{E:DATAHARMLESSIBPHYPERSURFACEERRORINTEGRALS}
				\\
		\int_{\Sigma_0^u}
			\left|
				(1 + 2 \upmu) (\Rad \Psi) (\GeoAng \Tanset^N \Psi) \upchipartialmodarg{\GeoAng^{N-1}}
			\right|
		\, d \vol
		& \lesssim
			\mathring{\upepsilon}^2.
			\label{E:QUADRATICDATAHARMLESSIBPHYPERSURFACEERRORINTEGRALS}
	\end{align}
	\end{subequations}
\end{lemma}

\begin{proof}
	See Subsect.~\ref{SS:OFTENUSEDESTIMATES} for some comments on the analysis.
	We first prove \eqref{E:HARMLESSIBPSPACETIMEERRORINTEGRALS}.
	All products on RHS~\eqref{E:LUNITIBPSPACETIMEERROR}
	contain a quadratic factor of
	$(\angdiff \GeoAng^N \Psi) \upchipartialmodarg{\GeoAng^{N-1}}$,
	$(\GeoAng^{N+1} \Psi) \upchipartialmodarg{\GeoAng^{N-1}}$,
	$(\GeoAng^N \Psi) \upchipartialmodarg{\GeoAng^{N-1}}$,
	or 
	$(\GeoAng^N \Psi) \Lunit \upchipartialmodarg{\GeoAng^{N-1}}$.
	With the help of 
	the estimates \eqref{E:TANGENTDIFFERNTIATEDGEOANGDEFORMSPHERELSHARPPOINTWISE}
	and \eqref{E:TANGENTDIFFERNTIATEDANGDEFORMTANGENTPOINTWISE} 
	and the $L^{\infty}$ estimates of Prop.~\ref{P:IMPROVEMENTOFAUX},
	it is easy to see that the remaining factors are bounded in $L^{\infty}$ by $\lesssim 1$.
	Hence it suffices to bound the spacetime integrals of the 
	magnitude of the four quadratic terms
	by $\lesssim$ RHS~\eqref{E:HARMLESSIBPSPACETIMEERRORINTEGRALS}.
	To bound the spacetime integral of 
	$
	\left|
		(\GeoAng^{N+1} \Psi) \upchipartialmodarg{\GeoAng^{N-1}} 
	\right|
	$, 
	we use spacetime Cauchy-Schwarz,
	Lemmas~\ref{L:KEYSPACETIMECOERCIVITY} and \ref{L:COERCIVENESSOFCONTROLLING},
	inequalities
	\eqref{E:LESSSINGULARTERMSMPOINTNINEINTEGRALBOUND}
	and 
	\eqref{E:MUCHEASIERBOUNDRYSPATIALINTEGRALFACTORL2ESTIMATE},
	simple estimates of the form
	$ab \lesssim a^2 + b^2$,
	and the fact that
	$\totTanmax{[1,N]}$
	is increasing in its arguments
	to deduce
	\begin{align} \label{E:SPACETIMEINTEGRALBOUNDFORINTEGRATIONBYPARTISINL}
		\int_{\mathcal{M}_{t,u}}
			&
			\left|
				(\GeoAng^{N+1} \Psi) \upchipartialmodarg{\GeoAng^{N-1}}
			\right|
		\, d \vol
			\\  
		& 
		\lesssim
		\varsigma 
		\TranminusdatasizeWithFactor
		\int_{\mathcal{M}_{t,u}}
			\left|
				\angdiff \GeoAng^N \Psi
			\right|^2
		\, d \vol
		+
		\varsigma^{-1} \TranminusdatasizeWithFactor^{-1}
		\int_{s=0}^t
			\left\|
				\upchipartialmodarg{\GeoAng^{N-1}}
			\right\|_{L^2(\Sigma_s^u)}^2
		\, ds
			\notag \\
		& \lesssim
			\varsigma
			\coerciveTanspacetimemax{[1,N]}(t,u)
			+
			\varsigma^{-1}
			\int_{s=0}^t
				\left\lbrace
				\int_{t'=0}^s
					\frac{1}{\upmu_{\star}^{1/2}(t',u)} \totTanmax{[1,N]}^{1/2}(t',u)
				\, dt'
			  \right\rbrace^2
			  +
			  \varsigma^{-1}
			  \mathring{\upepsilon}^2
			\, ds
			\notag
			\\
		& \lesssim
			\varsigma
			\coerciveTanspacetimemax{[1,N]}(t,u)
			+
			\varsigma^{-1}
			\int_{s=0}^t
				\totTanmax{[1,N]}(s,u)
			\, ds
			+
			\varsigma^{-1} \mathring{\upepsilon}^2,
			\notag
	\end{align}
	which is $\leq$ RHS~\eqref{E:HARMLESSIBPSPACETIMEERRORINTEGRALS} as desired.
	We clarify that in passing to the last inequality in \eqref{E:SPACETIMEINTEGRALBOUNDFORINTEGRATIONBYPARTISINL},
	we have used the fact that $\totTanmax{[1,N]}$ is increasing in its arguments
	and the estimate \eqref{E:LESSSINGULARTERMSMPOINTNINEINTEGRALBOUND}
	to deduce that
	$
	\displaystyle
	\int_{t'=0}^s
		\frac{1}{\upmu_{\star}^{1/2}(t',u)} \totTanmax{[1,N]}^{1/2}(t',u)
	\, dt'
	\lesssim 
	\totTanmax{[1,N]}^{1/2}(s,u)
	$,
	as we did in passing to the last line of \eqref{E:SECONDHARMLESSEXAMPLEINTEGRAL}.

	The spacetime integral of
	$
	\left|
		(\angdiff \GeoAng^N \Psi) \upchipartialmodarg{\GeoAng^{N-1}} 
	\right|
	$
	can be bounded in the same way.

	The spacetime integral of
	$\left| (\GeoAng^N \Psi) \upchipartialmodarg{\GeoAng^{N-1}} \right|$
	can be bounded by $\leq$ RHS~\eqref{E:HARMLESSIBPSPACETIMEERRORINTEGRALS}
	by using essentially the same arguments; we omit the details.

	To bound the spacetime integral of
	$
	\left|
		(\GeoAng^N \Psi) \Lunit \upchipartialmodarg{\GeoAng^{N-1}}
	\right|
	$, 
	by $\leq$ RHS~\eqref{E:HARMLESSIBPSPACETIMEERRORINTEGRALS},
	we first use Cauchy-Schwarz,
	Lemmas~\ref{L:KEYSPACETIMECOERCIVITY} and \ref{L:COERCIVENESSOFCONTROLLING},
	and inequality \eqref{E:NOTTOOHARDLUNITAPPLIEDTOBOUNDRYSPATIALINTEGRALFACTORL2ESTIMATE}
	to deduce
	\begin{align} \label{E:GEOANGNLUNITPARTIALLMODIFIEDSPACETIMEINTEGRAL}
		\int_{\mathcal{M}_{t,u}}
			\left|
				(\GeoAng^N \Psi) \Lunit \upchipartialmodarg{\GeoAng^{N-1}}
			\right|
		\, d \vol
		& \lesssim
			\int_{s=0}^t
				\left\|
					\GeoAng^N \Psi
				\right\|_{L^2(\Sigma_s^u)}
				\left\|
					\Lunit \upchipartialmodarg{\GeoAng^{N-1}}
				\right\|_{L^2(\Sigma_s^u)}
		\, ds
			\\
		& \lesssim
			\int_{s=0}^t
				\frac{1}{\upmu_{\star}(s,u)} 
				\totTanmax{[1,N-1]}^{1/2}(s,u)
				\totTanmax{[1,N]}^{1/2}(s,u)
			\, ds
			\notag \\
		& \ \
			+
			\mathring{\upepsilon}
			\int_{s=0}^t
				\frac{1}{\upmu_{\star}^{1/2}(s,u)} 
				\totTanmax{[1,N-1]}^{1/2}(s,u)
			\, ds.
			\notag
		\end{align}
	Finally, using simple estimates of the form $ab \lesssim a^2 + b^2$,
	the estimate \eqref{E:LESSSINGULARTERMSMPOINTNINEINTEGRALBOUND},
	and the fact that
	$\totTanmax{[1,N]}$
	is increasing in its arguments,
	we bound RHS~\eqref{E:GEOANGNLUNITPARTIALLMODIFIEDSPACETIMEINTEGRAL} by $\lesssim$ RHS
	\eqref{E:HARMLESSIBPSPACETIMEERRORINTEGRALS} as desired. This concludes the proof of 
	\eqref{E:HARMLESSIBPSPACETIMEERRORINTEGRALS}.

	We now prove \eqref{E:HARMLESSIBPHYPERSURFACEERRORINTEGRALS}
	and \eqref{E:DATAHARMLESSIBPHYPERSURFACEERRORINTEGRALS}.
	We first note that RHS~\eqref{E:LUNITIBPHYPERSURFACEERROR}
	is in magnitude 
	$\lesssim \varepsilon 
		\left|
			\GeoAng^N \Psi 
		\right|
		\left|
			\upchipartialmodarg{\GeoAng^{N-1}} 
		\right|
	$,
	an estimate that can easily be verified 
	with the help of the estimate \eqref{E:TANGENTDIFFERNTIATEDANGDEFORMTANGENTPOINTWISE}
	and the $L^{\infty}$ estimates of Prop.~\ref{P:IMPROVEMENTOFAUX}.
	Next, using Cauchy-Schwarz on $\Sigma_t^u$,
	Lemma~\ref{L:COERCIVENESSOFCONTROLLING}, 
	\eqref{E:MUCHEASIERBOUNDRYSPATIALINTEGRALFACTORL2ESTIMATE},
	and the estimate \eqref{E:LESSSINGULARTERMSMPOINTNINEINTEGRALBOUND}, 
	we deduce that
	\begin{align} \label{E:NOTSODIFFICULTSIGMATUERRORINTEGRAL}
	\varepsilon 
		\int_{\Sigma_t^u}
			\left|
				\GeoAng^N \Psi 
			\right|
			\left|
				\upchipartialmodarg{\GeoAng^{N-1}} 
			\right|
		\, d \tvol
		& \lesssim
		\varepsilon
		\left\|
			\GeoAng^N \Psi 
		\right\|_{L^2(\Sigma_t^u)}
		\left\|
			\upchipartialmodarg{\GeoAng^{N-1}}  
		\right\|_{L^2(\Sigma_t^u)}
		 \\
		&
		\lesssim
		\varepsilon
		\totTanmax{[1,N]}^{1/2}(t,u)
		\left\lbrace
		\totTanmax{[1,N]}^{1/2}(t,u)
		+
		\mathring{\upepsilon}
		\right\rbrace
			\notag \\
	& \lesssim \mbox{RHS~\eqref{E:HARMLESSIBPHYPERSURFACEERRORINTEGRALS}},
		\notag
	\end{align}
	as desired. We clarify that in passing to the second line of \eqref{E:NOTSODIFFICULTSIGMATUERRORINTEGRAL},
	we have used \eqref{E:LESSSINGULARTERMSMPOINTNINEINTEGRALBOUND}
	and the fact that $\totTanmax{[1,N]}$ is increasing in its arguments to bound
	the time integral on RHS~\eqref{E:MUCHEASIERBOUNDRYSPATIALINTEGRALFACTORL2ESTIMATE} by
	$\lesssim \totTanmax{[1,N]}^{1/2}(t,u)$.
	\eqref{E:DATAHARMLESSIBPHYPERSURFACEERRORINTEGRALS} then follows from
	\eqref{E:HARMLESSIBPHYPERSURFACEERRORINTEGRALS} with $t=0$
	and Lemma~\ref{L:INITIALSIZEOFL2CONTROLLING}.

	The proof of \eqref{E:QUADRATICDATAHARMLESSIBPHYPERSURFACEERRORINTEGRALS} is similar.
	The main difference is that the
	$L^{\infty}$ estimates of Prop.~\ref{P:IMPROVEMENTOFAUX}  
	imply only that LHS~\eqref{E:QUADRATICDATAHARMLESSIBPHYPERSURFACEERRORINTEGRALS}
	is 
	$
	\lesssim
	\int_{\Sigma_0^u}
			\left|
				\GeoAng^{N+1} \Psi 
			\right|
			\left|
				\upchipartialmodarg{\GeoAng^{N-1}} 
			\right|
		\, d \tvol
	$,
	without a gain of a factor $\varepsilon$.
	However, this integral is quadratically small in the data
	parameter $\mathring{\upepsilon}$,
	as is easy to verify using 
	Lemma~\ref{L:INITIALSIZEOFL2CONTROLLING}
	and the arguments given in the previous paragraph. 
	We have thus proved \eqref{E:QUADRATICDATAHARMLESSIBPHYPERSURFACEERRORINTEGRALS}
	and established the lemma.
\end{proof}

We now combine the previous results to prove the main lemma of Subsect.~\ref{SS:ERROINTEGRALSINVOLVINGIBPL}.

\begin{lemma}[\textbf{Bounds for difficult top-order error integrals connected to integration by parts involving $\Lunit$}]
	\label{L:BOUNDSFORDIFFICULTTOPORDERINTEGRALSINVOLVINGLUNITIBP}
	Assume that $1 \leq N \leq 18$ and $\varsigma > 0$.
	Let $\upchipartialmodarg{\GeoAng^{N-1}}$  be the
	partially modified quantity
	defined in \eqref{E:TRANSPORTPARTIALRENORMALIZEDTRCHIJUNK}.
	There exists a constant $C > 0$,
	independent of $\varsigma$,
	such that under the data-size and bootstrap assumptions 
	of Subsects.~\ref{SS:SIZEOFTBOOT}-\ref{SS:PSIBOOTSTRAP}
	and the smallness assumptions of Subsect.~\ref{SS:SMALLNESSASSUMPTIONS}, 
	the following estimates hold for $(t,u) \in [0,\Tboot) \times [0,U_0]$:
	\begin{align} \label{E:DIFFICULTLUNITSPACETIMEIBPINTEGRALBOUND}
		&
		\left|
			\int_{\mathcal{M}_{t,u}}
				(1 + 2 \upmu) (\Rad \Psi) (\GeoAng^{N+1} \Psi) \Lunit \upchipartialmodarg{\GeoAng^{N-1}}
			\, d \vol
		\right|
			\\
		& \leq
			\boxed{2}
			\int_{t'=0}^t
					\frac{\left\| [\Lunit \upmu]_- \right\|_{L^{\infty}(\Sigma_{t'}^u)}} 
							 {\upmu_{\star}(t',u)} 
				  \totTanmax{[1,N]}(t',u)
				\, dt'
			\notag	\\
		& \ \
			+
			C \varepsilon
			\int_{t'=0}^t
				\frac{1} 
					{\upmu_{\star}(t',u)} 
				 \totTanmax{[1,N]}(t',u)
			\, dt'
			+
			C
			\int_{t'=0}^t
				\frac{1} 
							 {\upmu_{\star}^{1/2}(t',u)} 
				 \totTanmax{[1,N]}(t',u)
			\, dt'
			\notag	\\
		& \ \
			+ C \mathring{\upepsilon}^2,
			\notag
				\\
		&
		\left|
		\int_{\Sigma_t^u}
				(1 + 2 \upmu) (\Rad \Psi) (\GeoAng^{N+1} \Psi) \upchipartialmodarg{\GeoAng^{N-1}}
		\, d \vol
		\right|
			\label{E:DIFFICULTLUNITHYPERSURFACEIBPINTEGRALBOUND} \\
		& \leq 
			\boxed{2}
			\frac{1}{\upmu_{\star}^{1/2}(t,u)}
			\totTanmax{[1,N]}^{1/2}(t,u)
			\left\| \Lunit \upmu \right\|_{L^{\infty}(\Sigmaminus{t}{t}{u})}
			\int_{t'=0}^t
				\frac{1}{\upmu_{\star}^{1/2}(t',u)} \totTanmax{[1,N]}^{1/2}(t',u)
			\, dt'
				\notag \\
		& \ \
			+ 
			C \varepsilon
			\frac{1}{\upmu_{\star}^{1/2}(t,u)}
			\totTanmax{[1,N]}^{1/2}(t,u)
			\int_{t'=0}^t
				\frac{1}{\upmu_{\star}^{1/2}(t',u)} \totTanmax{[1,N]}^{1/2}(t',u)
			\, dt'
			\notag
				\\
		& \ \
			+ 
			C 
			\totTanmax{[1,N]}^{1/2}(t,u)
			\int_{t'=0}^t
				\frac{1}{\upmu_{\star}^{1/2}(t',u)} \totTanmax{[1,N]}^{1/2}(t',u)
			\, dt'
				\notag \\
		&  \ \
			+ C \varsigma \totTanmax{[1,N]}(t,u)
			+ C \varsigma^{-1} \mathring{\upepsilon}^2 \frac{1}{\upmu_{\star}(t,u)}.
			\notag
	\end{align}

	\end{lemma}

\begin{proof}
	See Subsect.~\ref{SS:OFTENUSEDESTIMATES} for some comments on the analysis.
	To prove \eqref{E:DIFFICULTLUNITSPACETIMEIBPINTEGRALBOUND}, we 
	first use Cauchy-Schwarz and 
	the estimate
		$
		\left|
			\GeoAng
		\right|
		\leq 1 + C \varepsilon
		$
		(which follows from \eqref{E:GEOANGPOINTWISE} and
		the $L^{\infty}$ estimates of Prop.~\ref{P:IMPROVEMENTOFAUX})
		and in particular the estimates
		$
		\left\| 
			\Rad \Psi
		\right\|_{L^{\infty}(\Sigma_t^u)} 
		\lesssim 1
		$
		and
		$
		\left\| 
			\upmu
		\right\|_{L^{\infty}(\Sigma_t^u)} 
		\lesssim 1
		$
		to bound the LHS by 
		\begin{align} \label{E:ANGDIFFTIMESPSIPARTIALMODDIFFICULTTOPORDERFIRSTESTIMATE}
		&	\leq
		(1 + C \varepsilon)
		\int_{t'=0}^t
		\left\|
			\sqrt{\upmu} \angdiff \GeoAng^N \Psi
		\right\|_{L^2(\Sigma_{t'}^u)}
		\left\|
			\frac{1}{\sqrt{\upmu}} (\Rad \Psi) \Lunit \upchipartialmodarg{\GeoAng^{N-1}}
		\right\|_{L^2(\Sigma_{t'}^u)}
		\, dt'
			\\
		& \	\
			+ 
		C 
		\int_{t'=0}^t
			\left\|
				\sqrt{\upmu} \angdiff \GeoAng^N \Psi
			\right\|_{L^2(\Sigma_{t'}^u)}
			\left\|
				\Lunit \upchipartialmodarg{\GeoAng^{N-1}}
			\right\|_{L^2(\Sigma_{t'}^u)}
		\, dt'.
	\notag
	\end{align}
	The desired estimate \eqref{E:DIFFICULTLUNITSPACETIMEIBPINTEGRALBOUND}
	now follows from \eqref{E:ANGDIFFTIMESPSIPARTIALMODDIFFICULTTOPORDERFIRSTESTIMATE},
	Lemma~\ref{L:COERCIVENESSOFCONTROLLING},
	and inequalities \eqref{E:ANNOYINGLDERIVATIVEBOUNDRYSPATIALINTEGRALFACTORL2ESTIMATE}
	and \eqref{E:NOTTOOHARDLUNITAPPLIEDTOBOUNDRYSPATIALINTEGRALFACTORL2ESTIMATE}.
	Note that to bound the integral 
	$
	\displaystyle
	C 
	\int_{t'=0}^t
		\mathring{\upepsilon}
		\frac{1}{\upmu_{\star}^{1/2}(t',u)}
		\totTanmax{[1,N]}^{1/2}(t',u)
	\, dt'
	$,
	which is generated by the last term 
	on RHS~\eqref{E:ANNOYINGLDERIVATIVEBOUNDRYSPATIALINTEGRALFACTORL2ESTIMATE},
	we first use Young's inequality to bound the integrand by
	$
	\displaystyle
	\lesssim
	\frac{\mathring{\upepsilon}^2}{\upmu_{\star}^{1/2}(t',u)}
	+ 
	\frac{\totTanmax{[1,N]}(t',u)}{\upmu_{\star}^{1/2}(t',u)}
	$.
	We then bound the time integral of the first term in the previous expression by 
	$\lesssim \mathring{\upepsilon}^2
	$ 
	with the help of the estimate \eqref{E:LESSSINGULARTERMSMPOINTNINEINTEGRALBOUND}
	and the time integral of the second by 
	$\leq$ the third term on RHS~\eqref{E:DIFFICULTLUNITSPACETIMEIBPINTEGRALBOUND}.

	The proof of \eqref{E:DIFFICULTLUNITHYPERSURFACEIBPINTEGRALBOUND}
	is similar but relies on 
	\eqref{E:ANNOYINGBOUNDRYSPATIALINTEGRALFACTORL2ESTIMATE}
	and \eqref{E:MUCHEASIERBOUNDRYSPATIALINTEGRALFACTORL2ESTIMATE}
	in place of 
	\eqref{E:ANNOYINGLDERIVATIVEBOUNDRYSPATIALINTEGRALFACTORL2ESTIMATE}
	and
	\eqref{E:NOTTOOHARDLUNITAPPLIEDTOBOUNDRYSPATIALINTEGRALFACTORL2ESTIMATE};
	we omit the details, noting only that we encounter the term
	$
	\displaystyle
	C \mathring{\upepsilon}
				\frac{1}{\upmu_{\star}^{1/2}(t,u)}
				\totTanmax{[1,N]}^{1/2}(t,u)
	$
	generated by the last term on RHS~\eqref{E:ANNOYINGBOUNDRYSPATIALINTEGRALFACTORL2ESTIMATE}.
		We bound this term by using Young's inequality as follows:
	$
	\displaystyle
	C \mathring{\upepsilon}
				\frac{1}{\upmu_{\star}^{1/2}(t,u)}
				\totTanmax{[1,N]}^{1/2}(t,u)
	\leq 
	C \varsigma^{-1} \mathring{\upepsilon}^2 \frac{1}{\upmu_{\star}(t,u)}
	+ C \varsigma \totTanmax{[1,N]}(t,u)
	$.

\end{proof}

\subsection{Estimates for error integrals involving a loss of one derivative}
The following lemma plays a central role 
in our proof that the energy estimates become successively less degenerate with respect to 
powers of $\upmu_{\star}^{-1}$ as we descend below top order.
In the lemma, we consider the two most difficult error integrals
that we encounter in our proof of Prop.~\ref{P:TANGENTIALENERGYINTEGRALINEQUALITIES}.
Here, we bound them in a much simpler way that incurs a loss of one derivative
(which is permissible below top order).
The main advantage of these estimates
compared to the ones that do not lose derivatives is:
\emph{the derivative-losing estimates are much less degenerate with respect to $\upmu_{\star}^{-1}$.}

\begin{lemma}[\textbf{Estimates for error integrals involving a loss of one derivative}]
\label{L:ERRORINTEGRALSLOSEONEDERIVATIVE}
	Assume that $2 \leq N \leq 18$.
	Under the data-size and bootstrap assumptions 
	of Subsects.~\ref{SS:SIZEOFTBOOT}-\ref{SS:PSIBOOTSTRAP}
	and the smallness assumptions of Subsect.~\ref{SS:SMALLNESSASSUMPTIONS},
	the following estimates hold for $(t,u) \in [0,\Tboot) \times [0,U_0]$:
\begin{subequations}
\begin{align}
	& 
	\left|
		\int_{\mathcal{M}_{t,u}}
			(\Rad \Tanset^{N-1} \Psi)	 
			(\Rad \Psi) \GeoAng^{N-1} \mytr \upchi
		\, d \vol
	\right|
		\label{E:RADPSIDERIVATIVELOSINGERRORINTEGRAL} \\
	& \lesssim
		\int_{t'=0}^t
			\totTanmax{[1,N-1]}^{1/2}(t',u)
			\left\lbrace
				\int_{s=0}^{t'}
					\frac{\totTanmax{[1,N]}^{1/2}(s,u)}{\upmu_{\star}^{1/2}(s,u)}
				\, ds
			\right\rbrace
		\, dt'
		+
		\mathring{\upepsilon}
		\int_{t'=0}^t
			\frac{\totTanmax{[1,N-1]}^{1/2}(t',u)}{\upmu_{\star}^{1/2}(t',u)}
		\, dt',
			\notag \\
		& 
		\left|
		\int_{\mathcal{M}_{t,u}}
				(1 + 2 \upmu)
				(\Lunit \Tanset^{N-1} \Psi)	 
				(\Rad \Psi) \GeoAng^{N-1} \mytr \upchi
			\, d \vol
		\right|
	\label{E:LUNITPSIDERIVATIVELOSINGERRORINTEGRAL}
		\\
	& \lesssim
		\int_{t'=0}^t
			\frac{\totTanmax{[1,N-1]}^{1/2}(t',u)}{\upmu_{\star}^{1/2}(t',u)}
			\left\lbrace
				\int_{s=0}^{t'}
					\frac{\totTanmax{[1,N]}^{1/2}(s,u)}{\upmu_{\star}^{1/2}(s,u)}
				\, ds
			\right\rbrace
		\, dt'
		+
		\mathring{\upepsilon}
		\int_{t'=0}^t
			\frac{\totTanmax{[1,N-1]}^{1/2}(t',u)}{\upmu_{\star}^{1/2}(t',u)}
		\, dt'.
		\notag
\end{align}
\end{subequations}

\end{lemma}

\begin{proof}
See Subsect.~\ref{SS:OFTENUSEDESTIMATES} for some comments on the analysis.
We first prove \eqref{E:LUNITPSIDERIVATIVELOSINGERRORINTEGRAL}.
We begin by using the $L^{\infty}$ estimates of Prop.~\ref{P:IMPROVEMENTOFAUX}
to bound two of the factors in the integrand on the LHS as follows:
$
\left\|
	(1 + 2 \upmu) (\Rad \Psi)
\right\|_{L^{\infty}(\Sigma_t^u)}
\lesssim 1
$.
Using the previous estimate, 
Cauchy-Schwarz,
Lemma~\ref{L:COERCIVENESSOFCONTROLLING},
and the estimate \eqref{E:TANGENGITALEIKONALINTERMSOFCONTROLLING},
we bound LHS~\eqref{E:LUNITPSIDERIVATIVELOSINGERRORINTEGRAL} by 
\begin{align} \label{E:FIRSTLUNITPSIDERIVATIVELOSINGERRORINTEGRAL}
	& \lesssim
	\int_{t'=0}^t
		\frac{1}{\upmu_{\star}^{1/2}(t,u)}
		\left\|
			\sqrt{\upmu} \Lunit \Tanset^{N-1} \Psi
		\right\|_{L^2(\Sigma_{t'}^u)}
		\left\|
			\GeoAng^{N-1} \mytr \upchi
		\right\|_{L^2(\Sigma_{t'}^u)}
	\, dt'
		\\
	& \lesssim
	\int_{t'=0}^t
		\frac{\totTanmax{[1,N-1]}^{1/2}(t',u)}{\upmu_{\star}^{1/2}(t',u)}
		\left\lbrace
			\mathring{\upepsilon}
			+
			\int_{s=0}^{t'}
				\frac{\totTanmax{[1,N]}^{1/2}(s,u)}{\upmu_{\star}^{1/2}(s,u)}
			\, ds
		\right\rbrace
	\, dt'
		\notag
		\\
& \lesssim
	\int_{t'=0}^t
		\frac{\totTanmax{[1,N-1]}^{1/2}(t',u)}{\upmu_{\star}^{1/2}(t',u)}
		\left\lbrace
		\int_{s=0}^{t'}
			\frac{\totTanmax{[1,N]}^{1/2}(s,u)}{\upmu_{\star}^{1/2}(s,u)}
		\, ds
		\right\rbrace
	\, dt'
	+
	\mathring{\upepsilon}
	\int_{t'=0}^t
		\frac{\totTanmax{[1,N-1]}^{1/2}(t',u)}{\upmu_{\star}^{1/2}(t',u)}
	\, dt'
	\notag
\end{align}
as desired.

	The proof of \eqref{E:RADPSIDERIVATIVELOSINGERRORINTEGRAL} is similar,
	the only difference being that we start by
	bounding LHS~\eqref{E:RADPSIDERIVATIVELOSINGERRORINTEGRAL} by
	$
	\displaystyle
	\lesssim
	\int_{t'=0}^t
		\left\|
			\Rad \Tanset^{N-1} \Psi
		\right\|_{L^2(\Sigma_{t'}^u)}
		\left\|
			\GeoAng^{N-1} \mytr \upchi
		\right\|_{L^2(\Sigma_{t'}^u)}
	\, dt'
	$.
\end{proof}

\subsection{Proof of Prop.~\ref{P:TANGENTIALENERGYINTEGRALINEQUALITIES}}
\label{SS:PROOFOFPROPTANGENTIALENERGYINTEGRALINEQUALITIES}
\ \\

\noindent \textbf{Proof of \eqref{E:TOPORDERTANGENTIALENERGYINTEGRALINEQUALITIES}:}
Assume that $1 \leq N \leq 18$
and let $\Tanset^N$ be an $N^{th}-$order 
$\mathcal{P}_u$-tangent vectorfield operator.
From \eqref{E:E0DIVID} with $\Tanset^N \Psi$ in the role of $\Psi$, 
the decomposition \eqref{E:MULTERRORINT} with $\Tanset^N \Psi$ in the role of $\Psi$,
and definition \eqref{E:COERCIVESPACETIMEDEF},
we have
\begin{align} \label{E:E0DIVIDMAINESTIMATES}
	\enzero[\Tanset^N \Psi](t,u)
	& + 
	\flzero[\Tanset^N \Psi](t,u)
	+
	\coercivespacetime[\Tanset^N \Psi](t,u)
	\\
	& 
		=
		\enzero[\Tanset^N \Psi](0,u)
		- 
		\int_{\mathcal{M}_{t,u}}
			\left\lbrace
				(1 + 2 \upmu) (\Lunit \Tanset^N \Psi)
				+ 
				2 \Rad \Tanset^N \Psi 
			\right\rbrace
			\upmu \square_g(\Tanset^N \Psi) 
		\, d \vol
							\notag \\
				& \ \ 
						+
						\sum_{i=1}^5 
						\int_{\mathcal{M}_{t,u}}
							\basicenergyerrorarg{\Mult}{i}[\Tanset^N \Psi]
						\, d \vol.
							\notag
\end{align}
We will show that RHS~\eqref{E:E0DIVIDMAINESTIMATES} $\leq$ RHS~\eqref{E:TOPORDERTANGENTIALENERGYINTEGRALINEQUALITIES}.
Then, taking the max over that estimate for all such operators of order in between $1$ and $N$ 
and appealing to Defs.~\ref{D:MAINCOERCIVEQUANT} and \ref{D:COERCIVEINTEGRAL},
we conclude \eqref{E:TOPORDERTANGENTIALENERGYINTEGRALINEQUALITIES}.

To show that RHS~\eqref{E:E0DIVIDMAINESTIMATES} $\leq$ RHS~\eqref{E:TOPORDERTANGENTIALENERGYINTEGRALINEQUALITIES},
we first use Lemma~\ref{L:INITIALSIZEOFL2CONTROLLING} to deduce that
$\enzero[\Tanset^N \Psi](0,u) \lesssim \mathring{\upepsilon}^2$,
which is $\leq$ the first term on RHS~\eqref{E:TOPORDERTANGENTIALENERGYINTEGRALINEQUALITIES} as desired.

To bound the last integral
$
\sum_{i=1}^5
\int_{\mathcal{M}_{t,u}} 
	\cdots
$
on RHS~\eqref{E:E0DIVIDMAINESTIMATES} by $\leq$ RHS~\eqref{E:TOPORDERTANGENTIALENERGYINTEGRALINEQUALITIES},
we use Lemma~\ref{L:MULTIPLIERVECTORFIELDERRORINTEGRALS}.

We now address the first integral 
$
- \int_{\mathcal{M}_{t,u}} \cdots
$
on RHS~\eqref{E:E0DIVIDMAINESTIMATES}.
If $N \geq 2$ and $\Tanset^N$ is not of the form 
$\GeoAng^{N-1} \Lunit$ or $\GeoAng^N$, 
then the desired bound follows from 
\eqref{E:HARMLESSORDERNCOMMUTATORS}
and \eqref{E:STANDARDPSISPACETIMEINTEGRALS},
which together allow us to bound error integrals involving 
$Harmless^{\leq N}$ factors.
Note that these bounds do not produce any of the difficult
``boxed-constant-involving'' terms on RHS~\eqref{E:TOPORDERTANGENTIALENERGYINTEGRALINEQUALITIES}.

We now consider the case $\Tanset^N = \GeoAng^N$. The case 
$\Tanset^N = \GeoAng^{N-1} \Lunit$ can be treated in an identical fashion
and we omit those details. We 
start by substituting RHS~\eqref{E:GEOANGANGISTHEFIRSTCOMMUTATORIMPORTANTTERMS}
for the term $\upmu \square_g(\Tanset^N \Psi)$ on RHS~\eqref{E:E0DIVIDMAINESTIMATES}.
It suffices for us to bound the
integrals corresponding to the terms
$(\Rad \Psi) \GeoAng^N \mytr \upchi$
and
$\GeoAngFlatRadComponent (\angdiffuparg{\#} \Psi) \cdot (\upmu \angdiff \GeoAng^{N-1} \mytr \upchi)$
from RHS~\eqref{E:GEOANGANGISTHEFIRSTCOMMUTATORIMPORTANTTERMS},
for the above argument has already addressed how to
bound the integrals generated by $Harmless^{\leq N}$
terms.
To bound the difficult integral
\begin{align} \label{E:RADDIFFICULTERRORINTEGRAL}
	-           2
							\int_{\mathcal{M}_{t,u}}
							(\Rad \GeoAng^N \Psi)	 
							(\Rad \Psi) \GeoAng^N \mytr \upchi
							\, d \vol
\end{align}
by $\leq$ RHS~\eqref{E:TOPORDERTANGENTIALENERGYINTEGRALINEQUALITIES},
we first use Cauchy-Schwarz and \eqref{E:COERCIVENESSOFCONTROLLING} to bound it by
\begin{align} \label{E:FIRSTSTEPDIFFICULTINTEGRALBOUND}
	\leq 2 
		\int_{t' = 0}^t 
			\totTanmax{[1,N]}^{1/2}(t',u) 
			\left\| 
				(\Rad \Psi) \GeoAng^N \mytr \upchi 
			\right\|_{L^2(\Sigma_{t'}^u)} 
		\, dt'.
\end{align}
We now substitute the estimate \eqref{E:DIFFICULTTERML2BOUND}
(with $t$ in \eqref{E:DIFFICULTTERML2BOUND} replaced by $t'$)
for the second factor in the integrand of \eqref{E:FIRSTSTEPDIFFICULTINTEGRALBOUND}.
Following this substitution,
the desired bound of \eqref{E:FIRSTSTEPDIFFICULTINTEGRALBOUND}
by $\leq$ RHS~\eqref{E:TOPORDERTANGENTIALENERGYINTEGRALINEQUALITIES}
follows easily with the help of simple estimates of the form
$ab \lesssim a^2 + b^2$.
Note that these estimates account for
the portion $\boxed{4} \cdots$ of the first boxed constant integral
$\boxed{6} \cdots$ on RHS~\eqref{E:TOPORDERTANGENTIALENERGYINTEGRALINEQUALITIES}
and the full portion of 
the boxed constant integral $\boxed{8.1} \cdots$ 
on RHS~\eqref{E:TOPORDERTANGENTIALENERGYINTEGRALINEQUALITIES}.

We now bound the error integral
\begin{align} \label{E:LUNITDIFFICULTERRORINTEGRAL}
	-           \int_{\mathcal{M}_{t,u}}
							(1 + 2 \upmu)
							(\Lunit \GeoAng^N \Psi)	 
							(\Rad \Psi) \GeoAng^N \mytr \upchi
							\, d \vol.
\end{align}
To proceed, we use \eqref{E:TRANSPORTPARTIALRENORMALIZEDTRCHIJUNK}-\eqref{E:TRANSPORTPARTIALRENORMALIZEDTRCHIJUNKDISCREPANCY} 
to decompose
$ \GeoAng^N \mytr \upchi 
= \GeoAng \upchipartialmodarg{\GeoAng^{N-1}}
	- \GeoAng \upchipartialmodinhomarg{\GeoAng^{N-1}}$.
Since RHS~\eqref{E:HARMLESSNATUREOFPARTIALLYMODIFIEDDISCREPANCY} $= Harmless^{\leq N}$,
we have already suitably bounded the error integrals generated by
$\GeoAng \upchipartialmodinhomarg{\GeoAng^{N-1}}$.
We therefore must bound
\begin{align} \label{E:TOPORDERLUNITIBPDIFFICULTTERMS}
	-
	\int_{\mathcal{M}_{t,u}}
		(1 + 2 \upmu)
		(\Lunit \GeoAng^N \Psi)	 
		(\Rad \Psi) 
		\GeoAng \upchipartialmodarg{\GeoAng^{N-1}}
	\, d \vol
\end{align}
by $\leq$ RHS~\eqref{E:TOPORDERTANGENTIALENERGYINTEGRALINEQUALITIES}.
To this end, we integrate by parts using \eqref{E:LUNITANDANGULARIBPIDENTITY}
with 
$\ThirdSmoothFunction := \upchipartialmodarg{\GeoAng^{N-1}}$.
We bound the error integrals on the last line of RHS~\eqref{E:LUNITANDANGULARIBPIDENTITY}
and the $\int_{\Sigma_0^u} \cdots$ integral on the second line
using Lemma~\ref{L:HARMLESSIBPERRORINTEGRALS}.
It remains for us to bound the first two (difficult) integrals on
RHS~\eqref{E:LUNITANDANGULARIBPIDENTITY} 
by $\leq$ RHS~\eqref{E:TOPORDERTANGENTIALENERGYINTEGRALINEQUALITIES}.
The desired bounds have been derived in Lemma~\ref{L:BOUNDSFORDIFFICULTTOPORDERINTEGRALSINVOLVINGLUNITIBP}.
Note that these estimates account for
the remaining portion $\boxed{2} \cdots$ of the first boxed constant integral
$\boxed{6} \cdots$ on RHS~\eqref{E:TOPORDERTANGENTIALENERGYINTEGRALINEQUALITIES}
and the full portion of 
the boxed constant integral $\boxed{2} \cdots$ 
on RHS~\eqref{E:TOPORDERTANGENTIALENERGYINTEGRALINEQUALITIES}.

To complete the proof of \eqref{E:TOPORDERTANGENTIALENERGYINTEGRALINEQUALITIES},
it remains for us to bound the two error integrals generated by the term
$\GeoAngFlatRadComponent (\angdiffuparg{\#} \Psi) \cdot (\upmu \angdiff \GeoAng^{N-1} \mytr \upchi)$
from RHS~\eqref{E:GEOANGANGISTHEFIRSTCOMMUTATORIMPORTANTTERMS}.
These two integrals were suitably bounded by $\leq$ RHS~\eqref{E:TOPORDERTANGENTIALENERGYINTEGRALINEQUALITIES}
in Lemma~\ref{L:LESSDEGENERATEENERGYESTIMATEINTEGRALS}
(note that we are using the simple bound
$\left\lbrace 
					\ln \upmu_{\star}^{-1}(t',u) 
					+ 
					1 
				\right\rbrace^2
\lesssim 
\upmu_{\star}^{-1/2}(t',u)
$
in order to bound the 
integrand factors in the first integrals
on RHS~\eqref{E:FIRSTLESSDEGENERATEENERGYESTIMATEINTEGRALS}
and RHS~\eqref{E:SECONDLESSDEGENERATEENERGYESTIMATEINTEGRALS}).
Note that these estimates do not contribute to the difficult boxed constant terms
on RHS~\eqref{E:TOPORDERTANGENTIALENERGYINTEGRALINEQUALITIES}.
We have thus proved \eqref{E:TOPORDERTANGENTIALENERGYINTEGRALINEQUALITIES}.

\medskip
\noindent \textbf{Proof of \eqref{E:BELOWTOPORDERTANGENTIALENERGYINTEGRALINEQUALITIES}:}
We repeat the proof of \eqref{E:TOPORDERTANGENTIALENERGYINTEGRALINEQUALITIES}
with $N-1$ in the role of $N$
and with one critically important change:
we bound the difficult error
integrals
\[
	-           2
							\int_{\mathcal{M}_{t,u}}
							(\Rad \Tanset^{N-1} \Psi)	 
							(\Rad \Psi) \GeoAng^{N-1} \mytr \upchi
							\, d \vol
\]
and
\[
	-           \int_{\mathcal{M}_{t,u}}
							(1 + 2 \upmu)
							(\Lunit \Tanset^{N-1} \Psi)	 
							(\Rad \Psi) \GeoAng^{N-1} \mytr \upchi
							\, d \vol
\]
using the derivative-losing Lemma~\ref{L:ERRORINTEGRALSLOSEONEDERIVATIVE} in place
of the arguments used in bounding 
\eqref{E:RADDIFFICULTERRORINTEGRAL} and
\eqref{E:LUNITDIFFICULTERRORINTEGRAL}.
Note that in the proof of \eqref{E:TOPORDERTANGENTIALENERGYINTEGRALINEQUALITIES},
the error integrals
\eqref{E:RADDIFFICULTERRORINTEGRAL} 
and
\eqref{E:LUNITDIFFICULTERRORINTEGRAL}
were the only ones that
resulted in the presence of 
very degenerate terms
on RHS~\eqref{E:TOPORDERTANGENTIALENERGYINTEGRALINEQUALITIES},
such as the ``boxed-constant-involving''
terms. This explains why
RHS~\eqref{E:BELOWTOPORDERTANGENTIALENERGYINTEGRALINEQUALITIES} 
features fewer (and less degenerate/simpler) terms
and why it features the one new factor
that represents the loss of a derivative
(namely, the factor $\totTanmax{[1,N]}^{1/2}(s,u)$
on the second line of RHS~\eqref{E:BELOWTOPORDERTANGENTIALENERGYINTEGRALINEQUALITIES}).

\hfill $\qed$

\subsection{Proof of Prop.~\ref{P:MAINAPRIORIENERGY}}
\label{SS:PROOFOFPROPMAINAPRIORIENERGY}
\ \\

\noindent \textbf{Estimates for} 
$\totTanmax{[1,18]}$, 
$\coerciveTanspacetimemax{[1,18]}$,
$\totTanmax{[1,17]}$,
and
$\coerciveTanspacetimemax{[1,17]}$: 
	We first derive the estimates \eqref{E:MULOSSMAINAPRIORIENERGYESTIMATES} for 
	$\totTanmax{[1,18]}$, 
	$\coerciveTanspacetimemax{[1,18]}$,
	$\totTanmax{[1,17]}$,
	and
	$\coerciveTanspacetimemax{[1,17]}$,
	which are highly coupled and must be treated as a system.
	To this end, we set
	\begin{align}
		F(t,u) 
		& := \sup_{(\hat{t},\hat{u}) \in [0,t] \times [0,u]} 
			\iota_F^{-1}(\hat{t},\hat{u})
			\max
			\left\lbrace
				\totTanmax{[1,18]}(\hat{t},\hat{u}),
				\coerciveTanspacetimemax{[1,18]}(\hat{t},\hat{u})
			\right\rbrace, 
			\label{E:TOPORDERENERGYRESCALED} \\
		G(t,u) 
		& := \sup_{(\hat{t},\hat{u}) \in [0,t] \times [0,u]} 
		\iota_G^{-1}(\hat{t},\hat{u})
		\max
		\left\lbrace
			\totTanmax{[1,17]}(\hat{t},\hat{u}),
			\coerciveTanspacetimemax{[1,17]}(\hat{t},\hat{u})
		\right\rbrace,
		\label{E:JUSTBELOWTOPORDERENERGYRESCALED}
	\end{align}
	where for $t' \leq \hat{t} \leq t$ 
	and $u' \leq \hat{u} \leq U_0$,
	we define
	\begin{align} \label{E:INTEGRATINGFACTOR}
		\iota_1(t')
		&:=
		\exp
		\left(
			\int_{s=0}^{t'}
				\frac{1}{\sqrt{\Tboot - s}}
			\, ds
		\right)
		= \exp
			\left(
				2 \sqrt{\Tboot}
				- 2 \sqrt{\Tboot - t'}
			\right),
				\\
	\iota_2(t',u')
		&:=
		\exp
		\left(
			\int_{s=0}^{t'}
				\frac{1}{\upmu_{\star}^{9/10}(s,u')}
			\, ds
		\right),
			\label{E:SECONDINTEGRATINGFACTOR} \\
	\iota_F(t',u')
	& :=
		\upmu_{\star}^{-11.8}(t',u')
		\iota_1^c(t')
		\iota_2^c(t',u')
		e^{c t'}
		e^{c u'},
			\\
	\iota_G(t',u')
	& :=
		\upmu_{\star}^{-9.8}(t',u')
		\iota_1^c(t')
		\iota_2^c(t',u')
		e^{c t'}
		e^{c u'},
		\label{E:JUSTBELOWTOPORDERTOTALINTEGRATINGFACTOR}
	\end{align}
	and $c$ is a sufficiently large positive constant that we choose below.
	The functions \eqref{E:INTEGRATINGFACTOR}-\eqref{E:JUSTBELOWTOPORDERTOTALINTEGRATINGFACTOR}
	are approximate integrating factors that
	will allow us to absorb all of the error integrals
	on the RHSs of the inequalities of Prop.~\ref{P:TANGENTIALENERGYINTEGRALINEQUALITIES}.
	We claim that to obtain the desired estimates for 
$\totTanmax{[1,18]}$, 
$\coerciveTanspacetimemax{[1,18]}$,
$\totTanmax{[1,17]}$,
and
$\coerciveTanspacetimemax{[1,17]}$,
	it suffices to show that 
	\begin{align} \label{E:DESIREDRESCALEDAPRIORIBOUND}
		F(t,u) 
		& \leq C \mathring{\upepsilon}^2,
		\qquad
		G(t,u) 
		\leq C \mathring{\upepsilon}^2,
	\end{align}
	where $C$ in \eqref{E:DESIREDRESCALEDAPRIORIBOUND} is allowed to depend on $c$.
	To justify the claim, we use the fact that for a fixed $c$,
	the functions
	$\iota_1^c(t)$,
	$\iota_2^c(t,u)$,
	$e^{c t}$,
	and
	$e^{c u}$
	are uniformly bounded from above by a positive constant
	for $(t,u) \in [0,\Tboot) \times [0,U_0]$;
	all of these estimates are simple to derive,
	except for \eqref{E:SECONDINTEGRATINGFACTOR},
	which relies on \eqref{E:LESSSINGULARTERMSMPOINTNINEINTEGRALBOUND}.

	To prove \eqref{E:DESIREDRESCALEDAPRIORIBOUND},
	it suffices to
	show that
	there exist positive constants 
	$\upalpha_1$,
	$\upalpha_2$, 
	$\upbeta_1$,
	and 
	$\upbeta_2$ 
	with
	\begin{align} \label{E:CONSTSBOUND}
		\displaystyle \upalpha_1 + \frac{\upalpha_2 \upbeta_1}{1 - \upbeta_2} 
		& < 1,
		\qquad	 
		\upbeta_2 < 1
	\end{align}
	such that if $c$ is sufficiently large,
	then
	\begin{align} \label{E:TOPORDERALMOSTESTIMATED}
		F(t,u) 
		& \leq C \mathring{\upepsilon}^2 
			+ \upalpha_1 F(t,u) 
			+ \upalpha_2 G(t,u),
			\\
		G(t,u) \label{E:BELOWTOPORDERALMOSTESTIMATED}
		& \leq C \mathring{\upepsilon}^2 
			+ \upbeta_1 F(t,u)
			+ \upbeta_2 G(t,u).
	\end{align}
 Once we have obtained \eqref{E:TOPORDERALMOSTESTIMATED}-\eqref{E:BELOWTOPORDERALMOSTESTIMATED},
 we easily deduce from those estimates that
	\begin{align} \label{E:MOREBELOWTOPORDERALMOSTESTIMATED}
		G(t,u)
		& \leq 
			C \mathring{\upepsilon}^2 
			+ \frac{\upbeta_1}{1 - \upbeta_2} F(t,u),
			\\
		F(t,u)
		& \leq
			C \mathring{\upepsilon}^2
			+ 
			\left\lbrace
				\upalpha_1
				+ \frac{\upalpha_2 \upbeta_1}{1 - \upbeta_2} 
			\right\rbrace
			F(t,u).
		\label{E:MORETOPORDERALMOSTESTIMATED}
	\end{align}
	The desired bounds \eqref{E:DESIREDRESCALEDAPRIORIBOUND}
	now follow easily from
	\eqref{E:CONSTSBOUND}
	and
	\eqref{E:MOREBELOWTOPORDERALMOSTESTIMATED}-\eqref{E:MORETOPORDERALMOSTESTIMATED}.

	It remains for us to derive \eqref{E:TOPORDERALMOSTESTIMATED}-\eqref{E:BELOWTOPORDERALMOSTESTIMATED}.
	To this end,
	we will use the critically important estimates of Prop.~\ref{P:MUINVERSEINTEGRALESTIMATES}
	as well as the following simple estimates, which are easy to derive:
	\begin{align}
		\int_{t'=0}^{\hat{t}}
				\iota_1^c(t') \frac{1}{\sqrt{\Tboot - t'}} 
		\, dt'
		=
		\frac{1}{c}
		\int_{t'=0}^{\hat{t}}
				\frac{d}{dt'}
				\left\lbrace
					\iota_1^c(t') 
				\right\rbrace
		\, dt'
		& \leq \frac{1}{c} \iota_1^c(\hat{t}),
			\label{E:IF1SMALLNESS} \\
		\int_{t'=0}^{\hat{t}}
				\iota_2^c(t',\hat{u}) \frac{1}{\upmu_{\star}^{9/10}(t',\hat{u})} 
		\, dt'
		=
		\frac{1}{c}
		\int_{t'=0}^{\hat{t}}
				\frac{d}{dt'}
				\left\lbrace
					\iota_2^c(t',\hat{u}) 
				\right\rbrace
		\, dt'
		& \leq \frac{1}{c} \iota_2^c(\hat{t},\hat{u}),
			\label{E:IF2SMALLNESS} \\
		\int_{t'=0}^{\hat{t}}
				e^{c t'}
		\, dt'
		& \leq \frac{1}{c} e^{c \hat{t}},
			\label{E:IF3SMALLNESS} \\
		\int_{u'=0}^{\hat{u}}
				e^{c u'} 
		\, du'
		& \leq \frac{1}{c} e^{c \hat{u}}.
		\label{E:IF4SMALLNESS}
		\end{align}
	The smallness needed to close our estimates will come from
	taking $c$ to be large and $\varepsilon$ to be small.

	\emph{
	We stress that from now through
	inequality \eqref{E:NEARLYBELOWTOPORDERALMOSTESTIMATED},
	the constants $C$ can be chosen to be independent of $c$.
	}

	We also use the fact that
	$\iota_1^c(\cdot)$,
	$\iota_2^c(\cdot)$
	$e^{c \cdot}$,
	and $e^{c \cdot}$ are non-decreasing in their arguments,
	and the estimate \eqref{E:MUSTARINVERSEMUSTGROWUPTOACONSTANT}, 
	which implies that for
	$t' \leq \hat{t}$
	and $u' \leq \hat{u}$, we have the approximate monotonicity inequality
	\begin{align} \label{E:APPROXIMATETIMEMONOTONICITYGRONWALLPROOF}
		(1 + C \varepsilon) \upmu_{\star}(t',u') 
		& \geq \upmu_{\star}(\hat{t},\hat{u}).
	\end{align}
	In our arguments below, we do not explicitly mention these monotonicity properties every time we use them.

	We now set $N=18$, 
	multiply both sides of inequality
	\eqref{E:TOPORDERTANGENTIALENERGYINTEGRALINEQUALITIES}
	by $\iota_F^{-1}(t,u)$ and then set $(t,u) = (\hat{t},\hat{u})$.
	Similarly, we multiply both sides of inequality
	\eqref{E:BELOWTOPORDERTANGENTIALENERGYINTEGRALINEQUALITIES}
	by $\iota_G^{-1}(t,u)$ and then set $(t,u) = (\hat{t},\hat{u})$.
	To deduce \eqref{E:TOPORDERALMOSTESTIMATED}-\eqref{E:BELOWTOPORDERALMOSTESTIMATED},
	the difficult step is to obtain suitable bounds
	for the terms generated by the integrals on RHSs
	\eqref{E:TOPORDERTANGENTIALENERGYINTEGRALINEQUALITIES}-\eqref{E:BELOWTOPORDERTANGENTIALENERGYINTEGRALINEQUALITIES}.
	Once we have such bounds, we can then take
	$\sup_{(\hat{t},\hat{u}) \in [0,t] \times [0,u]}$
	of both sides of the resulting inequalities, and by virtue of
	definitions
	\eqref{E:TOPORDERENERGYRESCALED}-\eqref{E:JUSTBELOWTOPORDERENERGYRESCALED},
	we will easily
	conclude \eqref{E:TOPORDERALMOSTESTIMATED}-\eqref{E:BELOWTOPORDERALMOSTESTIMATED}.

	We now show how to obtain suitable bounds for 
	the terms generated by
	the ``borderline'' terms
	$\boxed{6} \int \cdots$,
	$\boxed{8.1} \int \cdots$, 
	and $
	\displaystyle
	\boxed{2}
			\frac{1}{\upmu_{\star}^{1/2}(t,u)}
			\totTanmax{[1,N]}^{1/2}(t,u)
			\left\| 
				\Lunit \upmu 
			\right\|_{L^{\infty}(\Sigmaminus{t}{t}{u})} 
			\int \cdots$
	on RHS~\eqref{E:TOPORDERTANGENTIALENERGYINTEGRALINEQUALITIES}.
	The terms generated by the
	remaining ``non-borderline'' terms on RHS~\eqref{E:TOPORDERTANGENTIALENERGYINTEGRALINEQUALITIES} are easier to treat.
	We start with the term 
	$\boxed{6} 
		\iota_F^{-1}(\hat{t},\hat{u}) \int_{t'=0}^{\hat{t}} \cdots
	$.
	Multiplying and 
	dividing by
	$\upmu_{\star}^{11.8}(t', \hat{u})$ in the integrand, 
	taking $\sup_{t' \in [0,\hat{t}]} \upmu_{\star}^{11.8}(t',\hat{u}) \totTanmax{[1,18]}(t',\hat{u})$,
	pulling the $\sup-$ed quantity out of the integral,
	and using the critically important integral estimate
	\eqref{E:KEYMUTOAPOWERINTEGRALBOUND} with $\Contwo = 12.8$, 
	we find that
	\begin{align} \label{E:MOSTDANGEROUSINTEGRALESTIMATE}
	& \boxed{6}
			\iota_F^{-1}(\hat{t},\hat{u})
			\int_{t'=0}^{\hat{t}}
					\frac{\left\| [\Lunit \upmu]_- \right\|_{L^{\infty}(\Sigma_{t'}^{\hat{u}})}} 
							 {\upmu_{\star}(t',\hat{u})} 
				  \totTanmax{[1,18]}(t',\hat{u})
			\, dt'
				\\
		& \leq 
				\boxed{6}
				\iota_F^{-1}(\hat{t},\hat{u})
				\sup_{t' \in [0,\hat{t}]}
				\left\lbrace 
					\upmu_{\star}^{11.8}(t',\hat{u}) \totTanmax{[1,18]}(t',\hat{u})
				\right\rbrace
				\int_{t'=0}^{\hat{t}}
					\left\| [\Lunit \upmu]_- \right\|_{L^{\infty}(\Sigma_{t'}^{\hat{u}})}
					\upmu_{\star}^{-12.8}(t',\hat{u})
				\, dt'
				\notag \\
		& \leq 
				\boxed{6}
				\upmu_{\star}^{11.8}(\hat{t},\hat{u})
				\sup_{t' \in [0,\hat{t}]}
				\left\lbrace 
					\iota_1^{-c}(t')
					\iota_2^{-c}(t',\hat{u})
					e^{-c t'} 
					e^{-c \hat{u}}
					\upmu_{\star}^{11.8}(t',\hat{u}) 
					\totTanmax{[1,18]}(t',\hat{u})
				\right\rbrace
					\notag \\
			& \ \ \ \
					\times
					\int_{t'=0}^{\hat{t}}
					\left\| [\Lunit \upmu]_- \right\|_{L^{\infty}(\Sigma_{t'}^{\hat{u}})}
					\upmu_{\star}^{-12.8}(t',\hat{u})
				\, dt'
				\notag \\
		&   \leq
				\frac{6 + C \sqrt{\varepsilon}}{11.8}
				F(\hat{t},\hat{u})
				\leq
				\frac{6 + C \sqrt{\varepsilon}}{11.8}
				F(t,u).
				\notag
	\end{align}

To handle the integral 
$\boxed{8.1} \iota_F^{-1}(\hat{t},\hat{u}) \int \cdots$,
we use a similar argument, but this time taking into account
that there are two time integrations. We find that
\begin{align} \label{E:SECONDMOSTDANGEROUSINTEGRALESTIMATE}
	&
	\boxed{8.1}
			\iota_F^{-1}(\hat{t},\hat{u})
			\int_{t'=0}^{\hat{t}}
				\frac{\left\| [\Lunit \upmu]_- \right\|_{L^{\infty}(\Sigma_{t'}^{\hat{u}})}} 
								 {\upmu_{\star}(t',\hat{u})} 
				\totTanmax{[1,18]}^{1/2}(t',\hat{u}) 
						\int_{s=0}^{t'}
							\frac{\left\| [\Lunit \upmu]_- \right\|_{L^{\infty}(\Sigma_s^{\hat{u}})}} 
									{\upmu_{\star}(s,\hat{u})} 
							\totTanmax{[1,18]}^{1/2}(s,\hat{u}) 
						\, ds
				\, dt'
				\\
		& \leq
				\frac{8.1 + C \sqrt{\varepsilon}}{5.9 \times 11.8}
				F(t,u).
				\notag
\end{align}

To handle the integral 
$
	\displaystyle
	\boxed{2}
			\iota_F^{-1}(\hat{t},\hat{u})
			\frac{1}{\upmu_{\star}^{1/2}(t,u)}
			\totTanmax{[1,N]}^{1/2}(t,u)
			\left\| 
				\Lunit \upmu 
			\right\|_{L^{\infty}(\Sigmaminus{t}{t}{u})} 
			\int \cdots$,
we use a similar argument based on the critically important estimate \eqref{E:KEYHYPERSURFACEMUTOAPOWERINTEGRALBOUND}.
We find that
\begin{align} \label{E:THIRDMOSTDANGEROUSINTEGRALESTIMATE}
	&
	\boxed{2}
	\iota_F^{-1}(\hat{t},\hat{u})
	\frac{1}{\upmu_{\star}^{1/2}(\hat{t},\hat{u})}
	\totTanmax{[1,18]}^{1/2}(\hat{t},\hat{u})
	\left\| 
			\Lunit \upmu 
	\right\|_{L^{\infty}(\Sigmaminus{\hat{t}}{\hat{t}}{\hat{u}})}
	\int_{t'=0}^t
			\frac{1}{\upmu_{\star}^{1/2}(t',\hat{u})} \totTanmax{[1,18]}^{1/2}(t',\hat{u})
	\, dt'
		\\
	& \leq \frac{2 + C \sqrt{\varepsilon}}{5.4} F(t,u).
	\notag
\end{align}

The important point is that for small $\varepsilon$,
the factors 
$\displaystyle \frac{6 + C \sqrt{\varepsilon}}{11.8}$ on RHS
\eqref{E:MOSTDANGEROUSINTEGRALESTIMATE},
$\displaystyle \frac{8.1 + C \sqrt{\varepsilon}}{5.9 \times 11.8}$
on RHS~\eqref{E:SECONDMOSTDANGEROUSINTEGRALESTIMATE},
and $\displaystyle\frac{2 + C \sqrt{\varepsilon}}{5.4}$
on RHS~\eqref{E:THIRDMOSTDANGEROUSINTEGRALESTIMATE}
sum to 
$\displaystyle 
\frac{6}{11.8} 
+ \frac{8.1}{5.9 \times 11.8} 
+ \frac{2}{5.4} 
+ C \sqrt{\varepsilon} 
< 1
$.
This sum is the main contributor to the constant
$\upalpha_1$ on RHS~\eqref{E:TOPORDERALMOSTESTIMATED}.

The remaining integrals are easier to treat.
We now show how to bound the term arising from the integral
on the $12^{th}$ line of RHS~\eqref{E:TOPORDERTANGENTIALENERGYINTEGRALINEQUALITIES}, 
which involves three time integrations.
The term arising from the integrals on the $10^{th}$ and $11^{th}$ lines of RHS~\eqref{E:TOPORDERTANGENTIALENERGYINTEGRALINEQUALITIES}
can be handled using similar arguments, so we do not provide those details.
We claim that the following sequence of inequalities holds for the term of interest,
which yields the desired bound:
\begin{align} \label{E:ESTIMATEFORTRIPLETIMEINTEGRATEDERRORINTEGRAL}
					& 
					C
					\iota_F^{-1}(\hat{t},\hat{u})
					\int_{t'=0}^{\hat{t}}
					\frac{1} 
							 {\upmu_{\star}(t',\hat{u})} 
				  \totTanmax{[1,18]}^{1/2}(t',\hat{u})
				  \int_{s = 0}^{t'}
				  	\frac{1}{\upmu_{\star}(s,\hat{u})}
				  	\int_{s' = 0}^s
				  	\frac{1} 
							 {\upmu_{\star}^{1/2}(s',\hat{u})} 
							 \totTanmax{[1,18]}^{1/2}(s',\hat{u})
						\, ds'	 
				  \, ds
				\, dt'
					\\
				& \leq
					\frac{C}{c} 
					\iota_F^{-1}(\hat{t},\hat{u})
					\iota_2^{c/2}(\hat{t},\hat{u})
					\int_{t'=0}^{\hat{t}}
					\frac{1} 
							 {\upmu_{\star}(t',\hat{u})} 
				  \totTanmax{[1,18]}^{1/2}(t',\hat{u})
				  	\notag \\
				 & \ \
				 		\times
				 		\int_{s = 0}^{t'}
				  	\frac{1}{\upmu_{\star}(s,\hat{u})}
				  	\sup_{(s',u') \in [0,s] \times [0,\hat{u}]}
						\left\lbrace 
							\iota_2^{-c/2}(s',u') \totTanmax{[1,18]}^{1/2}(s',u')
						\right\rbrace
				  	\, ds
				\, dt'
					\notag \\
			& \leq
					\frac{C}{c} 
					\iota_F^{-1}(\hat{t},\hat{u})
					\iota_2^{c/2}(\hat{t},\hat{u})
					\int_{t'=0}^{\hat{t}}
					\frac{1} 
							 {\upmu_{\star}(t',\hat{u})} 
				  \totTanmax{[1,18]}^{1/2}(t',\hat{u})
				  \notag \\
			& \ \ 
						\times
						\sup_{(s',u') \in [0,t'] \times [0,\hat{u}]}
						\left\lbrace 
							\upmu_{\star}(s',u') \iota_2^{-c/2}(s',u') \totTanmax{[1,18]}^{1/2}(s',u')
						\right\rbrace
				  \int_{s = 0}^{t'}
				  	\frac{1}{\upmu_{\star}^2(s,\hat{u})}
				  \, ds
				\, dt'
					\notag \\
		& \leq
					\frac{C}{c} 
					\iota_F^{-1}(\hat{t},\hat{u})
					\iota_2^{c/2}(\hat{t},\hat{u})
					\sup_{(s',u') \in [0,\hat{t}] \times [0,\hat{u}]}
						\left\lbrace 
							\upmu_{\star}(s',u') \iota_2^{-c/2}(s',u') \totTanmax{[1,18]}^{1/2}(s',u')
						\right\rbrace
						\notag \\
		& \times 
					\sup_{(s',u') \in [0,\hat{t}] \times [0,\hat{u}]}
						\left\lbrace 
							\totTanmax{[1,18]}^{1/2}(s',u')
						\right\rbrace
					\int_{t'=0}^{\hat{t}}
						\frac{1} 
							 {\upmu_{\star}(t',\hat{u})} 
				  	\int_{s = 0}^{t'}
				  		\frac{1}{\upmu_{\star}^2(s,\hat{u})}
				  	\, ds
				\, dt'
				\notag \\
		& \leq
					\frac{C}{c} 
					\upmu_{\star}(\hat{t},\hat{u})
					\sup_{(s',u') \in [0,\hat{t}] \times [0,\hat{u}]}
						\left\lbrace 
							\iota_F^{-1}(s',u') \totTanmax{[1,18]}(s',u')
						\right\rbrace
				\times 
					\int_{t'=0}^{\hat{t}}
						\frac{1} 
							 {\upmu_{\star}(t',\hat{u})} 
				  	\int_{s = 0}^{t'}
				  		\frac{1}{\upmu_{\star}^2(s,\hat{u})}
				  	\, ds
				\, dt'
				\notag \\
		& \leq \frac{C}{c} F(\hat{t},\hat{u})
			\leq \frac{C}{c} F(t,u)
			\notag,
\end{align}
which yields the desired smallness factor $\displaystyle \frac{1}{c}$.
We now explain how to derive 
\eqref{E:ESTIMATEFORTRIPLETIMEINTEGRATEDERRORINTEGRAL}.
To deduce the first inequality,
we multiplied and divided by
$\iota_2^{c/2}(t',\hat{u})$ in the integral 
$\int \cdots ds'$,
then pulled
$\displaystyle 
\sup_{(s',u') \in [0,s] \times [0,\hat{u}]}
						\left\lbrace 
							\iota_2^{-c/2}(s',u') \totTanmax{[1,18]}^{1/2}(s',u')
						\right\rbrace
$
out of the integral, and finally used \eqref{E:IF2SMALLNESS}
to gain the smallness factor $\displaystyle \frac{1}{c}$
from the remaining terms
$
\displaystyle
\int_{s'=0}^s
	\frac{1} 
		{\upmu_{\star}^{1/2}(s',\hat{u})} 
		\iota_2^{c/2}(s',\hat{u})
\, ds'
$.
To derive the second
inequality in \eqref{E:ESTIMATEFORTRIPLETIMEINTEGRATEDERRORINTEGRAL},
we multiplied and divided by
$\upmu_{\star}(s,\hat{u})$ in the integral 
$\int \cdots ds$,
and used the approximate monotonicity property
\eqref{E:APPROXIMATETIMEMONOTONICITYGRONWALLPROOF}
to pull the factor 
$
\displaystyle
\sup_{(s',u') \in [0,t'] \times [0,\hat{u}]}
						\left\lbrace 
							\upmu_{\star}(s',u') \iota_2^{-c/2}(s',u') \totTanmax{[1,18]}^{1/2}(s',u')
						\right\rbrace
$
out of the $ds$ integral, which costs us a harmless multiplicative factor of $1 + C \varepsilon$.
The third inequality in \eqref{E:ESTIMATEFORTRIPLETIMEINTEGRATEDERRORINTEGRAL}
follows easily. To derive the fourth inequality,
we use the monotonicity of
$\iota_1^c(\cdot)$,
	$\iota_2^c(\cdot)$
	$e^{c \cdot}$,
	and $e^{c \cdot}$,
and
\eqref{E:APPROXIMATETIMEMONOTONICITYGRONWALLPROOF}.
To derive the fifth inequality, 
we use inequality \eqref{E:LOSSKEYMUINTEGRALBOUND} twice.
The final inequality follows easily.

Similarly, we claim that we can bound the terms on the 
$5^{th}$ through $8^{th}$ lines of RHS~\eqref{E:TOPORDERTANGENTIALENERGYINTEGRALINEQUALITIES} as follows:
\begin{align} 
	\label{E:FORGOTTENEASYHYPERSUFRACEINTEGRALENERGYESTIMATETERM}
			\iota_F^{-1}(\hat{t},\hat{u})
			C \varepsilon
			\int_{t'=0}^{\hat{t}}
				\frac{1} {\upmu_{\star}(t',\hat{u})} 
						\totTanmax{[1,N]}^{1/2}(t',\hat{u}) 
						\int_{s=0}^{t'}
							\frac{1} 
									{\upmu_{\star}(s,\hat{u})} 
							\totTanmax{[1,N]}^{1/2}(s,\hat{u}) 
						\, ds
				\, dt'
			& \leq 
				C \varepsilon F(\hat{t},\hat{u})
				\leq C \varepsilon F(t,u), 
				\\
		\label{E:FIRSTEASYHYPERSUFRACEINTEGRALENERGYESTIMATETERM}
			\iota_F^{-1}(\hat{t},\hat{u})
			C
			\varepsilon
			\int_{t'=0}^{\hat{t}}
					\frac{1} 
						{\upmu_{\star}(t',\hat{u})} 
				  \totTanmax{[1,18]}(t',\hat{u})
				\, dt'
			& \leq 
				C \varepsilon F(\hat{t},\hat{u})
				\leq C \varepsilon F(t,u), 
				\\
			\iota_F^{-1}(\hat{t},\hat{u})
			C
			\varepsilon
			\frac{1}{\upmu_{\star}^{1/2}(\hat{t},\hat{u})}
			\totTanmax{[1,18]}^{1/2}(\hat{t},\hat{u})
			\int_{t'=0}^{\hat{t}}
				\frac{1}{\upmu_{\star}^{1/2}(t',\hat{u})} \totTanmax{[1,18]}^{1/2}(t',\hat{u})
			\, dt'
			& \leq 
				C \varepsilon F(\hat{t},\hat{u})
				\leq C \varepsilon F(t,u), 
				\label{E:SECONDEASYHYPERSUFRACEINTEGRALENERGYESTIMATETERM} \\
			\iota_F^{-1}(\hat{t},\hat{u})
			C
			\totTanmax{[1,18]}^{1/2}(\hat{t},\hat{u})
			\int_{t'=0}^{\hat{t}}
				\frac{1}{\upmu_{\star}^{1/2}(t',\hat{u})} \totTanmax{[1,18]}^{1/2}(t',\hat{u})
			\, dt'
			& 
			\leq \frac{C}{c} F(\hat{t},\hat{u})
			\leq \frac{C}{c} F(t,u).
			\label{E:THIRDEASYHYPERSUFRACEINTEGRALENERGYESTIMATETERM}
\end{align}
To derive 
\eqref{E:FORGOTTENEASYHYPERSUFRACEINTEGRALENERGYESTIMATETERM},
we use arguments similar to the ones we used 
in deriving \eqref{E:SECONDMOSTDANGEROUSINTEGRALESTIMATE}, 
but in place of the delicate estimate \eqref{E:KEYMUTOAPOWERINTEGRALBOUND},
we use the estimate \eqref{E:LOSSKEYMUINTEGRALBOUND},
whose large constant $C$ is compensated for by the availability of the smallness 
factor $\varepsilon$.
Similarly, to derive 
\eqref{E:FIRSTEASYHYPERSUFRACEINTEGRALENERGYESTIMATETERM},
we use arguments similar to the ones we used 
in deriving \eqref{E:MOSTDANGEROUSINTEGRALESTIMATE}, 
using \eqref{E:LOSSKEYMUINTEGRALBOUND}
in place of \eqref{E:KEYMUTOAPOWERINTEGRALBOUND}.
To derive \eqref{E:SECONDEASYHYPERSUFRACEINTEGRALENERGYESTIMATETERM},
we use arguments similar to the ones we used above,
but we now multiply and divide by 
$\upmu_{\star}^{5.9}(t',\hat{u})$
in the time integral on LHS~\eqref{E:SECONDEASYHYPERSUFRACEINTEGRALENERGYESTIMATETERM}
and use \eqref{E:LOSSKEYMUINTEGRALBOUND}.
To derive \eqref{E:THIRDEASYHYPERSUFRACEINTEGRALENERGYESTIMATETERM},
we use similar arguments based on multiplying and dividing by
$\iota_2^{c/2}(t',\hat{u})$ in the time integral
and using \eqref{E:IF2SMALLNESS}.

Similarly, we derive the bound
\begin{align}
	&
	C
	\iota_F^{-1}(\hat{t},\hat{u})
	\int_{t'=0}^{\hat{t}}
		\frac{1}{\sqrt{\Tboot - t'}} \totTanmax{[1,18]}(t',\hat{u})
	\, dt'
	\leq \frac{C}{c} F(\hat{t},\hat{u})
	\leq \frac{C}{c} F(t,u)
\end{align}for the term on the $9^{th}$ line of RHS~\eqref{E:TOPORDERTANGENTIALENERGYINTEGRALINEQUALITIES}
by multiplying and dividing by $\iota_1^{c}(t')$
in the integrand and using
\eqref{E:IF1SMALLNESS} to gain the smallness factor
$
\displaystyle
\frac{1}{c}
$.

Similarly, we derive the bound
\begin{align}
	&
	C
	(1 + \varsigma^{-1})
	\iota_F^{-1}(\hat{t},\hat{u})
	\int_{u'=0}^{\hat{u}}
		\totTanmax{[1,18]}(\hat{t},u')
	\, du'
	\leq \frac{C}{c} F(\hat{t},\hat{u})
	\leq \frac{C}{c} F(t,u)
\end{align}
for the term on the $13^{th}$ line of RHS~\eqref{E:TOPORDERTANGENTIALENERGYINTEGRALINEQUALITIES}
by multiplying and dividing by $e^{cu'}$
in the integrand and using \eqref{E:IF4SMALLNESS} to 
gain the smallness factor
$
\displaystyle
\frac{1}{c}
$.

It is easy to see that the terms 
arising from the terms 
on the 
first and the
next-to-last lines of RHS~\eqref{E:TOPORDERTANGENTIALENERGYINTEGRALINEQUALITIES},
namely
$
C (1 + \varsigma^{-1}) \mathring{\upepsilon}^2 \upmu_{\star}^{-3/2}(t,u)
$,
$
C \varepsilon \totTanmax{[1,N]}(t,u)
$,
$
C \varsigma \totTanmax{[1,N]}(t,u)
$,
and
$
C \varsigma \coerciveTanspacetimemax{[1,N]}(t,u)
$,
are respectively bounded 
(after multiplying by $\iota_F^{-1}$ and taking the relevant sup)
by 
$\leq C (1 + \varsigma^{-1}) \mathring{\upepsilon}^2
$,
$
\leq C \varepsilon F(t,u)
$,
$
\leq
C \varsigma F(t,u)
$,
and 
$
\leq
C \varsigma F(t,u)
$.

To bound the term arising from the last integral on RHS
\eqref{E:TOPORDERTANGENTIALENERGYINTEGRALINEQUALITIES}, 
we argue as follows 
with the help of \eqref{E:IF2SMALLNESS} and \eqref{E:APPROXIMATETIMEMONOTONICITYGRONWALLPROOF}
(recall that $N=18$):
\begin{align} \label{E:BELOWTOPORDERTERMCOUPLEDINTOTOPORDERESTIMATE}
			& C
			\iota_F^{-1}(\hat{t},\hat{u})
			\int_{t'=0}^{\hat{t}}
					\frac{1} 
							 {\upmu_{\star}^{5/2}(t',\hat{u})} 
				  \totTanmax{[1,17]}(t',\hat{u})
			\, dt'
				\\
			& 
			\leq
			C
			\upmu_{\star}^{9.8}(\hat{t},\hat{u})
			\iota_1^{-c}(\hat{t})
			\iota_2^{-c}(\hat{t},\hat{u})
			e^{-c \hat{t}} 
			e^{-c \hat{u}}
			\sup_{t' \in [0,\hat{t}]}
				\left(
				\frac{\upmu_{\star}(\hat{t},\hat{u})} 
				{\upmu_{\star}(t',\hat{u})} 
				\right)^2
					\notag \\
	& \ \
			\times
			\sup_{t' \in [0,\hat{t}]}
				\left\lbrace 
					\iota_2^{-c}(t') \totTanmax{[1,17]}(t',\hat{u})
				\right\rbrace
			\times
			\int_{t'=0}^{\hat{t}}
				\frac{\iota_2^c(t')} 
				{\upmu_{\star}^{1/2}(t',\hat{u})} 
			\, dt'
				\notag 
				\\
		& 
			\leq
			C
			\iota_2^{-c}(\hat{t})
			\sup_{t' \in [0,\hat{t}]}
				\left\lbrace 
					\iota_G^{-1}(t',\hat{u}) \totTanmax{[1,17]}(t',\hat{u})
				\right\rbrace
			\times
			\int_{t'=0}^{\hat{t}}
				\frac{\iota_2^c(t')} 
				{\upmu_{\star}^{1/2}(t',\hat{u})} 
			\, dt'
				\notag 
				\\
		& \leq \frac{C}{c}
						G(\hat{t},\hat{u})
			\leq \frac{C}{c}G(t,u).
			\notag
\end{align}

We now bound the terms 
$\displaystyle \iota_G^{-1}(\hat{t},\hat{u}) \times \cdots$ arising from the terms on
RHS~\eqref{E:BELOWTOPORDERTANGENTIALENERGYINTEGRALINEQUALITIES}.
All terms except the one arising from the
integral involving the top-order factor $\totTanmax{[1,18]}^{1/2}$
(featured in the $ds$ integral on RHS~\eqref{E:BELOWTOPORDERTANGENTIALENERGYINTEGRALINEQUALITIES})
can be bounded by
$
\displaystyle 
\leq C \mathring{\upepsilon}^2 
+ \frac{C}{c}(1 + \varsigma^{-1})G(t,u) 
+ C \varsigma G(t,u)
$
by using essentially the same arguments given above.
To handle the remaining term involving the top-order factor 
$\totTanmax{[1,18]}^{1/2}$, 
we use arguments similar to the ones
we used to prove \eqref{E:ESTIMATEFORTRIPLETIMEINTEGRATEDERRORINTEGRAL}
(in particular, we use inequality \eqref{E:LOSSKEYMUINTEGRALBOUND} twice)
to bound it as follows:
\begin{align} \label{E:TOPORDERTERMCOUPLEDINTOBELOWTOPORDERESTIMATE}
			& C
			\iota_G^{-1}(\hat{t},\hat{u})
			\int_{t'=0}^{\hat{t}}
				\frac{1}{\upmu_{\star}^{1/2}(t',\hat{u})} 
						\totTanmax{[1,17]}^{1/2}(t',\hat{u}) 
						\int_{s=0}^{t'}
							\frac{1}{\upmu_{\star}^{1/2}(s,\hat{u})} 
							\totTanmax{[1,18]}^{1/2}(s,\hat{u}) 
						\, ds
				\, dt'
				\\
			& 
			\leq
			C
			\iota_G^{-1}(\hat{t},\hat{u})
			\sup_{(t',u') \in [0,\hat{t}] \times [0,\hat{u}]}
			\left\lbrace 
				\upmu_{\star}^{4.9}(t',\hat{u}) \totTanmax{[1,17]}^{1/2}(t',u')
			\right\rbrace
			\times
			\sup_{(t',u') \in [0,\hat{t}] \times [0,\hat{u}]}
			\left\lbrace 
				\upmu_{\star}^{5.9}(t',\hat{u}) \totTanmax{[1,18]}^{1/2}(t',u')
			\right\rbrace
				\notag \\
			& \times
				\int_{t'=0}^{\hat{t}}
						\frac{1}{\upmu_{\star}^{5.4}(t',\hat{u})} 
						\int_{s=0}^{t'}
							\frac{1}{\upmu_{\star}^{6.4}(s,\hat{u})}
						\, ds
				\, dt'
					\notag
					\\
		& \leq  C
						F^{1/2}(\hat{t},\hat{u})
						G^{1/2}(\hat{t},\hat{u})
			\leq C F(t,u)
				+ \frac{1}{2}G(t,u).
			\notag
\end{align}


Inserting all of these estimates into the RHSs of
$\iota_F^{-1}(\hat{t},\hat{u}) \times$ \eqref{E:TOPORDERTANGENTIALENERGYINTEGRALINEQUALITIES}$(\hat{t},\hat{u})$
and $\iota_G^{-1}(\hat{t},\hat{u}) \times$ \eqref{E:BELOWTOPORDERTANGENTIALENERGYINTEGRALINEQUALITIES}$(\hat{t},\hat{u})$
and taking $\sup_{(\hat{t},\hat{u}) \in [0,t] \times [0,u]}$ of both sides,
we deduce that
\begin{align} \label{E:NEARLYTOPORDERALMOSTESTIMATED}
		F(t,u) 
		& \leq C (1 + \varsigma^{-1})\mathring{\upepsilon}^2 
			+ \left\lbrace
					\frac{6}{11.8} 
					+ \frac{8.1}{5.9 \times 11.8} 
					+ \frac{2}{5.4} 
					+ C \sqrt{\varepsilon} 
				\right\rbrace
				F(t,u) 
					\\
		& \ \
			+ \frac{C}{c} (1 + \varsigma^{-1}) F(t,u)
			+ \frac{C}{c} G(t,u)
			+ C \varsigma F(t,u),
			\notag \\
		G(t,u) \label{E:NEARLYBELOWTOPORDERALMOSTESTIMATED}
		& \leq C \mathring{\upepsilon}^2 
				+ \frac{C}{c} F(t,u)
				+ C F(t,u)
				+ \frac{C}{c} (1 + \varsigma^{-1}) G(t,u)
				+ C \varsigma G(t,u)
				+ \frac{1}{2} G(t,u).
\end{align}
We remind the reader that the constants $C$ 
in \eqref{E:NEARLYTOPORDERALMOSTESTIMATED}-\eqref{E:NEARLYBELOWTOPORDERALMOSTESTIMATED}
can be chosen to be independent of $c$.
The desired estimates \eqref{E:TOPORDERALMOSTESTIMATED}-\eqref{E:BELOWTOPORDERALMOSTESTIMATED}
now follow from first choosing $\varsigma$ to be sufficiently small,
then choosing $c$ to be sufficiently large, then choosing
$\varepsilon$ to be sufficiently small,
and using the aforementioned fact that
$\displaystyle 
\frac{6}{11.8} 
+ \frac{8.1}{5.9 \times 11.8} 
+ \frac{2}{5.4} 
+ C \sqrt{\varepsilon} 
< 1
$.

\noindent \textbf{Estimates for} $\totTanmax{[1, \leq 16]}$ 
\textbf{and}
$\coerciveTanspacetimemax{[1,\leq 16]}$ \textbf{via a descent scheme}:
We now explain how to use inequality
\eqref{E:BELOWTOPORDERTANGENTIALENERGYINTEGRALINEQUALITIES}
to derive the estimates
for 
$\totTanmax{[1,\leq 16]}$ 
and
$\coerciveTanspacetimemax{[1,\leq 16]}$
by downward induction.
Unlike our analysis of the strongly coupled pair
$\max
\left\lbrace
	\totTanmax{[1,18]}, 
	\coerciveTanspacetimemax{[1,18]}
\right\rbrace
$
and 
$\max
\left\lbrace
	\totTanmax{[1,17]}, 
	\coerciveTanspacetimemax{[1,17]}
\right\rbrace
$,
we can derive the desired estimates
for
$\max
\left\lbrace
	\totTanmax{[1,16]}, 
	\coerciveTanspacetimemax{[1,16]}
\right\rbrace
$
by using only inequality \eqref{E:BELOWTOPORDERTANGENTIALENERGYINTEGRALINEQUALITIES}
and the already derived estimates
for 
$\max
\left\lbrace
	\totTanmax{[1,17]}, 
	\coerciveTanspacetimemax{[1,17]}
\right\rbrace
$.
At the end of the proof, we will describe
the minor changes needed to derive
the desired estimates for
$
\max
\left\lbrace
	\totTanmax{[1,15]}, 
	\coerciveTanspacetimemax{[1,15]}
\right\rbrace
$, 
$
\max
\left\lbrace
	\totTanmax{[1,14]}, 
	\coerciveTanspacetimemax{[1,14]}
\right\rbrace
$, 
$\cdots$,
$
\max
\left\lbrace
	\totTanmax{1}, 
	\coerciveTanspacetimemax{1}
\right\rbrace$.

To begin, we define the following analogs of 
\eqref{E:JUSTBELOWTOPORDERTOTALINTEGRATINGFACTOR}
and
\eqref{E:JUSTBELOWTOPORDERENERGYRESCALED}:
\begin{align} 	\label{E:MORETHANONEBELOWTOPORDERTOTALINTEGRATINGFACTOR}
	\iota_H(t',u')
	& :=
		\upmu_{\star}^{-7.8}(t',u')
		\iota_1^c(t')
		\iota_2^c(t',u')
		e^{c t'}
		e^{c u'},
\end{align}
\begin{align}
		H(t,u) 
		& := \sup_{(\hat{t},\hat{u}) \in [0,t] \times [0,u]} 
		\iota_H^{-1}(\hat{t},\hat{u})
		\max
		\left\lbrace
			\totTanmax{[1,16]}(\hat{t},\hat{u}),
			\coerciveTanspacetimemax{[1,16]}(\hat{t},\hat{u})
		\right\rbrace.
		\label{E:MORETHANONEBELOWTOPORDERENERGYRESCALED}
	\end{align}
Note that the power of $\upmu_{\star}^{-1}$
in the factor $\upmu_{\star}^{-7.8}$ has been reduced
by two in \eqref{E:MORETHANONEBELOWTOPORDERTOTALINTEGRATINGFACTOR}
compared to \eqref{E:JUSTBELOWTOPORDERTOTALINTEGRATINGFACTOR},
which corresponds to less singular behavior of
$
\max
\left\lbrace
			\totTanmax{[1,16]},
			\coerciveTanspacetimemax{[1,16]}
\right\rbrace
$
near the shock. As before, to prove the desired estimate
\eqref{E:MULOSSMAINAPRIORIENERGYESTIMATES} (now with $M=3$),
it suffices to prove
\begin{align} \label{E:MORETHANONEBELOWTOPORDERDESIREDBOUND}
	H(t,u)
	\leq C \mathring{\upepsilon}^2.
\end{align}

We now set $N=17$, 
multiply both sides of inequality
\eqref{E:BELOWTOPORDERTANGENTIALENERGYINTEGRALINEQUALITIES}
by $\iota_H^{-1}(t,u)$, 
and then set $(t,u) = (\hat{t},\hat{u})$
(note that $N=17$ in \eqref{E:BELOWTOPORDERTANGENTIALENERGYINTEGRALINEQUALITIES} 
corresponds to estimating
$
\max
\left\lbrace
			\totTanmax{[1,16]},
			\coerciveTanspacetimemax{[1,16]}
\right\rbrace
$).
With one exception,
we can bound
all terms arising from the integrals on RHS~\eqref{E:BELOWTOPORDERTANGENTIALENERGYINTEGRALINEQUALITIES}
by $\displaystyle \leq C \mathring{\upepsilon}^2 + \frac{C}{c}(1 + \varsigma^{-1})H(t,u) + \varsigma H(t,u)$
(where $C$ is independent of $c$)
by using the same arguments
that we used in deriving the estimate for
$
\max
\left\lbrace
			\totTanmax{[1,17]},
			\coerciveTanspacetimemax{[1,17]}
\right\rbrace
$.
The exceptional term is the one arising from the 
integral involving the above-present-order factor $\totTanmax{[1,17]}^{1/2}$.
We bound the exceptional term as follows by using
inequality \eqref{E:LOSSKEYMUINTEGRALBOUND}, 
the approximate monotonicity of $\iota_H$,
and the estimate
$
\totTanmax{[1,17]}^{1/2} \leq C_c \mathring{\upepsilon} \upmu_{\star}^{-4.9}(t,u)
$
(which follows from the already proven estimate \eqref{E:DESIREDRESCALEDAPRIORIBOUND} for $G(t,u)$):
\begin{align} \label{E:BELOWTOPORDERTERMCOUPLEDINTOMOREBELOWTOPORDERESTIMATE}
			& C
			\iota_H^{-1}(\hat{t},\hat{u})
			\int_{t'=0}^{\hat{t}}
				\frac{1}{\upmu_{\star}^{1/2}(t',\hat{u})} 
						\totTanmax{[1,16]}^{1/2}(t',\hat{u}) 
						\int_{s=0}^{t'}
							\frac{1}{\upmu_{\star}^{1/2}(s,\hat{u})} 
							\totTanmax{[1,17]}^{1/2}(s,\hat{u}) 
						\, ds
				\, dt'
				\\
			& 
			\leq
			C_c
			\mathring{\upepsilon}
			\iota_H^{-1/2}(\hat{t},\hat{u})
			\sup_{(t',u') \in [0,\hat{t}] \times [0,\hat{u}]}
			\left\lbrace 
				\iota_H^{-1/2}(t',u')
				\totTanmax{[1,16]}^{1/2}(t',u')
			\right\rbrace
			\times
			\int_{t'=0}^{\hat{t}}
						\frac{1}{\upmu_{\star}^{1/2}(t',\hat{u})} 
						\int_{s=0}^{t'}
							\frac{1}{\upmu_{\star}^{5.4}(s,\hat{u})}
						\, ds
				\, dt'
					\notag
					\\
		& \leq 
			C_c
			\mathring{\upepsilon}
			\iota_H^{-1/2}(\hat{t},\hat{u})
			\upmu_{\star}^{-3.9}(\hat{t},\hat{u})
			\sup_{(t',u') \in [0,\hat{t}] \times [0,\hat{u}]}
			\left\lbrace 
				\iota_H^{-1/2}(t',u')
				\totTanmax{[1,16]}^{1/2}(t',u')
			\right\rbrace
			\notag
				\\
		& \leq 
			C_c
			\mathring{\upepsilon}
			H^{1/2}(\hat{t},\hat{u})
			\leq 
			C_c \mathring{\upepsilon}^2
			+
			\frac{1}{2}
			H(t,u).
			\notag
\end{align}

In total, we have obtained the following
analog of \eqref{E:NEARLYBELOWTOPORDERALMOSTESTIMATED}:
\begin{align}
	H(t,u) \label{E:NEARLYMORETHANONEBELOWTOPORDERALMOSTESTIMATED}
		& \leq C_c \mathring{\upepsilon}^2 
				+ \frac{C}{c} (1 + \varsigma^{-1}) H(t,u)
				+ \frac{1}{2} H(t,u)
				+ C \varsigma H(t,u),
\end{align}
where $C_c$ is the only constant that depends on $c$.
The desired bound \eqref{E:MORETHANONEBELOWTOPORDERDESIREDBOUND}
easily follows from \eqref{E:NEARLYMORETHANONEBELOWTOPORDERALMOSTESTIMATED}
by first choosing $\varsigma$ to be sufficiently small
and then $c$ to be sufficiently large
so that we can absorb all factors of $H$ on RHS~\eqref{E:NEARLYMORETHANONEBELOWTOPORDERALMOSTESTIMATED}
into the LHS.

The desired bounds \eqref{E:NOMULOSSMAINAPRIORIENERGYESTIMATES} for
$
\max
\left\lbrace
			\totTanmax{[1,15]},
			\coerciveTanspacetimemax{[1,15]}
\right\rbrace
$, 
$
\max
\left\lbrace
			\totTanmax{[1,14]},
			\coerciveTanspacetimemax{[1,14]}
\right\rbrace
$, 
$\cdots$
can be (downward) inductively derived
by using an argument similar to the one
we used to bound
$
\max
\left\lbrace
			\totTanmax{[1,16]},
			\coerciveTanspacetimemax{[1,16]}
\right\rbrace
$,
which relied on the already available bounds for 
$
\max
\left\lbrace
			\totTanmax{[1,17]},
			\coerciveTanspacetimemax{[1,17]}
\right\rbrace
$.
The only difference is that
we define the analog of the approximating integrating factor \eqref{E:MORETHANONEBELOWTOPORDERTOTALINTEGRATINGFACTOR}
to be
$ \upmu_{\star}^{-p}
	\iota_1^c(t')
	\iota_2^c(t',u')
	e^{c t'}
	e^{c u'}
$,
where $p = 5.8$ for the 
$
\max
\left\lbrace
			\totTanmax{[1,15]},
			\coerciveTanspacetimemax{[1,15]}
\right\rbrace
$ 
estimate,
$p = 3.8$ for the 
$
\max
\left\lbrace
			\totTanmax{[1,14]},
			\coerciveTanspacetimemax{[1,14]}
\right\rbrace
$ 
estimate,
$p = 1.8$ for the
$
\max
\left\lbrace
			\totTanmax{[1,13]},
			\coerciveTanspacetimemax{[1,13]}
\right\rbrace
$
estimate,
and $p=0$ for the 
$
\max
\left\lbrace
			\totTanmax{[1,\leq 12]},
			\coerciveTanspacetimemax{[1,\leq 12]}
\right\rbrace
$ 
estimates;
these latter estimates \emph{do not involve any singular factor} of $\upmu_{\star}^{-1}$.
There is one important new detail that is relevant
for these estimates: in deriving the analog of the inequalities
\eqref{E:BELOWTOPORDERTERMCOUPLEDINTOMOREBELOWTOPORDERESTIMATE}
for 
$
\max
\left\lbrace
			\totTanmax{[1,\leq 12]},
			\coerciveTanspacetimemax{[1,\leq 12]}
\right\rbrace
$,
we use the estimate \eqref{E:LESSSINGULARTERMSMPOINTNINEINTEGRALBOUND}
in place of the estimate \eqref{E:LOSSKEYMUINTEGRALBOUND};
the estimate \eqref{E:LESSSINGULARTERMSMPOINTNINEINTEGRALBOUND}
is what allows us to break the $\upmu_{\star}^{-1}$ degeneracy.

\section{The Stable Shock Formation Theorem}
\label{S:STABLESHOCKFORMATION}
In this section, 
we state and prove our main stable shock formation theorem.

\subsection{The diffeomorphic nature of \texorpdfstring{$\Upsilon$}{the change of variables map} and continuation criteria}
\label{SS:PRELIMINARYLEMMASMAINTHEOREM}
We first provide a technical lemma concerning the change of variables map
$\Upsilon$ and a lemma providing continuation criteria.

\begin{lemma}[\textbf{Sufficient conditions for $\Upsilon$ to be a global diffeomorphism}]
\label{L:CHOVREMAINSADIFFEOMORPHISM}
	Assume the data-size and bootstrap assumptions 
	of Subsects.~\ref{SS:SIZEOFTBOOT}-\ref{SS:PSIBOOTSTRAP}
	and the smallness assumptions of Subsect.~\ref{SS:SMALLNESSASSUMPTIONS}.
	Assume in addition that 
	\begin{align}
		\inf_{(t,u) \in [0,\Tboot) \times [0,U_0]} \upmu_{\star}(t,u) > 0.
	\end{align}
	Then the change of variables map $\Upsilon$ extends to a global $C^{1,1}$ diffeomorphism 
	from $[0,\Tboot] \times [0,U_0] \times \mathbb{T}$
onto its image.

\end{lemma}

\begin{proof}
	First, from Lemma~\ref{L:CHOVMAPREGULARITY},
	we see that $\Upsilon$ extends as a $C^{1,1}$ function
	defined on $[0,\Tboot] \times [0,U_0] \times \mathbb{T}$.
	Hence, to prove the lemma, it remains for us to 
	show that $\Upsilon$ is a diffeomorphism
	from $[0,\Tboot] \times [0,U_0] \times \mathbb{T}$
	onto its image. To this end, we first use 
	\eqref{E:JACOBIAN},
	\eqref{E:POINTWISEESTIMATEFORGTANCOMP},
	the bootstrap assumptions \eqref{E:PSIFUNDAMENTALC0BOUNDBOOTSTRAP},
	the fact that $\gt_{ij} = \delta_{ij} + \mathcal{O}(\Psi)$,
	\eqref{E:UPMULINFTY},
	and the assumption $\inf_{(t,u) \in [0,\Tboot) \times [0,U_0]} \upmu_{\star}(t,u) > 0$
	to deduce that the Jacobian determinant of $\Upsilon$
	is uniformly bounded from above and from below strictly away from $0$.
	Hence, from the inverse function theorem, we deduce that
	$\Upsilon$ extends as a $C^{1,1}$ local diffeomorphism from
	$[0,\Tboot] \times [0,U_0] \times \mathbb{T}$
	onto its image. 

	To show that $\Upsilon$ is a global diffeomorphism
	on the domain under consideration, 
	it suffices to show that 
	for $u_1, u_2 \in [0,U_0]$ with $u_1 < u_2$,
	the distinct curves 
	$\ell_{\Tboot,u_1},
	\ell_{\Tboot,u_2}
	\subset \Sigma_{\Tboot}^{U_0}
	$ 
	do not intersect each other 
	and that
	for each $u \in [0,U_0]$,
	$\Upsilon(\Tboot,u,\cdot)$ is an injection from $\mathbb{T}$ onto its image.
	To rule out the intersection of two distinct curves,
	we use 
	\eqref{E:ALTERNATEEIKONAL},
	\eqref{E:UPMULINFTY},
	the assumption $\inf_{(t,u) \in [0,\Tboot) \times [0,U_0]} \upmu_{\star}(t,u) > 0$,
	the bootstrap assumptions \eqref{E:PSIFUNDAMENTALC0BOUNDBOOTSTRAP},
	and the fact that $\gt_{ij} = \delta_{ij} + \mathcal{O}(\Psi)$
	to deduce that $\sum_{a=1}^2 |\partial_a u|$ is uniformly bounded from
	above and strictly from below away from $0$.
	It follows that the (closed) null plane portions $\mathcal{P}_u^{\Tboot}$
	corresponding to two distinct values of $u \in [0,U_0]$ cannot intersect, which
	yields the desired result.
	It remains for us to show that when $u \in [0,U_0]$,
	$\Upsilon(\Tboot,u,\cdot)$ is a diffeomorphism from $\mathbb{T}$ onto its image.
	To this end, we note that for each fixed $u \in [0,U_0]$,
	the rectangular component $\Upsilon^2(\Tboot,u,\cdot)$
	(which can be identified with the local rectangular coordinate $x^2$),  
	viewed as a $\mathbb{T}$-valued function of $\vartheta \in \mathbb{T}$,
	is homotopic to the degree-one map $\Upsilon^2(0,u,\cdot)$
	by the homotopy $\Upsilon^2(\cdot,u,\cdot): [0,\Tboot] \times \mathbb{T} \rightarrow \mathbb{T}$.
	Hence, it is a basic result of degree theory 
	(see, for example, the Hopf Degree Theorem in \cite{vGaP2010})
	that $\Upsilon^2(\Tboot,u,\cdot)$ 
	is also a degree-one map.
	Next, we note that 
	\eqref{E:PERTURBEDPART},
	\eqref{E:PURETANGENTIALGEOANGSMALLIPOINTWISE} with $N=0$,
	the $L^{\infty}$ estimates of Prop.~\ref{P:IMPROVEMENTOFAUX},
	and Lemma~\ref{L:COORDCOMPESTIMATES} together imply that
	$\CoordAng \Upsilon^2(\Tboot,u,\vartheta) 
	=
	[\GeoAngCoordComp^{-1} \GeoAng^2](\Tboot,u,\vartheta) 
	=
	1 + \mathcal{O}(\varepsilon)
	$
	for $\vartheta \in \mathbb{T}$.
	From this estimate and the degree-one property of $\Upsilon^2(\Tboot,u,\cdot)$, 
	we deduce\footnote{Recall that if $f: \mathbb{T} \rightarrow \mathbb{T}$ 
	is a $C^1$ surjective map without critical points, then $f$ is degree-one if for $p,q \in \mathbb{T}$,
	$1 = \sum_{p \in f^{-1}(q)} \mbox{\upshape sign } (d_p f)$,
	where $d_p f$ denotes the differential of $f$ at $p$
	and the $d_p f$ are computed relative to an atlas corresponding to the 
	smooth orientation on $\mathbb{T}$ chosen at the beginning of the article.
	It is a basic fact of degree theory that the sum is independent of $q$.
	Note that in the context of the present argument, 
	the role of $d f(\cdot)$ is effectively played by
	$\CoordAng \Upsilon^2(\Tboot,u,\cdot)$.}
	that (for sufficiently small $\varepsilon$),
	$\Upsilon^2(\Tboot,u,\cdot)$ is a bijection\footnote{The surjective property of this map is easy to deduce.} 
	from $\mathbb{T}$ to $\mathbb{T}$.
	Hence, $\Upsilon(\Tboot,u,\cdot)$ is injective, which is the desired result.
\end{proof}

We now provide some continuation criteria, which we will
use to ensure that the solution survives until the shock forms.

\begin{lemma}[\textbf{Continuation criteria}]
\label{L:CONTINUATIONCRITERIA}
Let $(\mathring{\Psi},\mathring{\Psi}_0) := (\Psi|_{\Sigma_0},\partial_t \Psi|_{\Sigma_0}) 
\in H_{\Euct}^{19}(\Sigma_0^1) \times H_{\Euct}^{18}(\Sigma_0^1)$
be initial data for the covariant wave equation $\square_{g(\Psi)} \Psi = 0$
that are compactly supported in $\Sigma_0^1$
(see Remark~\ref{R:SOBOLEVSPACESMULTIPLECOORDINATESYSTEMS} regarding the Sobolev spaces $H_{\Euct}^N(\Sigma_0^1)$).
Let $T_{(Local)} > 0$ and $U_0 \in (0,1]$, and assume that
the corresponding classical solution $\Psi$ exists on an
(``open at the top'')
spacetime region $\mathcal{M}_{T_{(Local)},U_0}$
(see Def.~\ref{D:HYPERSURFACESANDCONICALREGIONS})
that is completely determined by the non-trivial data lying in 
$\Sigma_0^{U_0}$ and the trivial data 
lying to the right of the line $\lbrace x^1 = 0 \rbrace$ in $\Sigma_0$
(see Figure~\ref{F:REGION} on pg.~\pageref{F:REGION}).
Let $u$ be the eikonal function that verifies the eikonal equation \eqref{E:INTROEIKONAL}
with the initial data \eqref{E:INTROEIKONALINITIALVALUE}.
Assume that $\upmu > 0$ on $\mathcal{M}_{T_{(Local)},U_0}$ and that
the change of variables map $\Upsilon$ from geometric to 
rectangular coordinates (see Def.~\ref{D:CHOVMAP}) is 
a $C^1$ diffeomorphism from $[0,T_{(Local)}) \times [0,U_0] \times \mathbb{T}$ onto $\mathcal{M}_{T_{(Local)},U_0}$.
Let $\mathcal{H}$ be the set of real numbers 
$b$ such that the following conditions hold:
\begin{itemize}
	\item The rectangular components $g_{\mu \nu}(\cdot)$, $(\mu, \nu = 0,1,2)$, are smooth on a neighborhood of $b$.
	\item $g_{00}(b) < 0$.
	\item The eigenvalues of the $2 \times 2$ matrix $\gt_{ij}(b)$ 
	(see Def.~\ref{D:FIRSTFUND}),
	$(i,j = 1,2)$, are positive.
\end{itemize}

Assume that none of the following $4$ breakdown scenarios occur:
\begin{enumerate}
	\item $\inf_{\mathcal{M}_{T_{(Local)},U_0}} \upmu = 0$. 
	\item $\sup_{\mathcal{M}_{T_{(Local)},U_0}} \upmu = \infty$.
	\item There exists a sequence $p_n \in \mathcal{M}_{T_{(Local)},U_0}$
		such that $\Psi(p_n)$ escapes every compact subset of $\mathcal{H}$ as $n \to \infty$.
	\item $\sup_{\mathcal{M}_{T,U_0}} 
				\max_{\kappa=0,1,2,3}
				\left|
					\partial_{\kappa} \Psi
				\right|
		= \infty$.
\end{enumerate}

In addition, assume that the following condition is verified:
\begin{enumerate}
	\setcounter{enumi}{4}
	\item The change of variables map $\Upsilon$
		extends to the compact set $[0,T_{(Local)}] \times [0,U_0] \times \mathbb{T}$
		as a (global) $C^1$ diffeomorphism onto its image.
\end{enumerate}

Then there exists a $\Delta > 0$ such that 
$\Psi$, 
$u$, 
and all of the
other geometric quantities defined
throughout the article can be uniquely extended 
(where $\Psi$ and $u$ are classical solutions)
to a strictly larger region of the form
$\mathcal{M}_{T_{(Local)} + \Delta,U_0}$
into which their Sobolev regularity 
relative to both geometric and rectangular coordinates
is propagated.
Moreover, if $\Delta$ is sufficiently small,
then none of the four breakdown scenarios
occur in the larger region, and
$\Upsilon$ extends to $[0,T_{(Local)} + \Delta] \times [0,U_0] \times \mathbb{T}$
as a (global) $C^1$ diffeomorphism onto its image.
\end{lemma}

\begin{proof}[Sketch of a proof]
	Lemma~\ref{L:CONTINUATIONCRITERIA} is mostly standard.
	A sketch of the proof was provided in
	\cite{jS2014b}*{Proposition 21.1.1}, to which we refer the reader
	for more details. Here, we only mention the main ideas.
	Criterion $(3)$ is connected to avoiding a breakdown
	in hyperbolicity of the equation. 
	Criterion $(4)$ is a standard criterion used to locally continue
	the solution relative to the rectangular
	coordinates.
	Criteria $(1)$ and $(2)$ and the assumption on $\Upsilon$
	are connected to ruling out the blowup of $u$, 
	degeneracy of the change of variables map,
	and degeneracy of the region $\mathcal{M}_{T_{(Local)},U_0}$.
	In particular, criteria $(1)$ and $(2)$
	play a role in a proving that
	$\sum_{a = 1}^2 |\partial_a u|$ is uniformly bounded from
	above and strictly from below away from $0$
	on $\mathcal{M}_{T_{(Local)},U_0}$
	(the proof was essentially given in the proof of Lemma~\ref{L:CHOVREMAINSADIFFEOMORPHISM}).
\end{proof}

\subsection{The main stable shock formation theorem}
\label{SS:MAINTHEOREM}
We now state and prove the main result of the article.

\begin{theorem}[\textbf{Stable shock formation}]
\label{T:MAINTHEOREM}
Let $(\mathring{\Psi},\mathring{\Psi}_0) := (\Psi|_{\Sigma_0},\partial_t \Psi|_{\Sigma_0}) 
\in H_{\Euct}^{19}(\Sigma_0^1) \times H_{\Euct}^{18}(\Sigma_0^1)$
(see Remark~\ref{R:SOBOLEVSPACESMULTIPLECOORDINATESYSTEMS})
be initial data for the covariant wave equation 
$\square_{g(\Psi)} \Psi = 0$
that are compactly supported in $\Sigma_0^1$
and that verify the data-size assumptions\footnote{Recall that in Remark~\ref{R:EXISTENCEOFDATA}, 
we outlined a proof that such data exist.} 
of Subsect.~\ref{SS:DATAASSUMPTIONS}.
In particular, let
$\mathring{\upepsilon}$,
$\mathring{\updelta}$,
and $\TranminusdatasizeWithFactor$
be the data-size parameters from
\eqref{E:PSIDATAASSUMPTIONS} and \eqref{E:CRITICALBLOWUPTIMEFACTOR}.
Assume that the rectangular metric component functions verify
the structural assumptions 
\eqref{E:NONVANISHINGNONLINEARCOEFFICIENT}
and
\eqref{E:GINVERSE00ISMINUSONE}.
For each $U_0 \in [0,1]$, let $T_{(Lifespan);U_0}$
be the classical lifespan of the solution
in the region that is completely determined by the non-trivial data lying in 
$\Sigma_0^{U_0}$ and the trivial data 
lying to the right of the line $\lbrace x^1 = 0 \rbrace$ in $\Sigma_0$
(see Figure~\ref{F:REGION} on pg.~\pageref{F:REGION}).
If $\mathring{\upepsilon}$ is sufficiently small
relative to 
$\mathring{\updelta}^{-1}$
and 
$\TranminusdatasizeWithFactor$
(in the sense explained in Subsect.~\ref{SS:SMALLNESSASSUMPTIONS}),
then the following conclusions hold,
where all constants can be chosen to be independent of $U_0$.

\medskip

\noindent \underline{\textbf{Dichotomy of possibilities.}}
One of the following mutually disjoint possibilities must occur,
where $\upmu_{\star}(t,u)$ is defined in \eqref{E:MUSTARDEF}.
\begin{enumerate}
	\renewcommand{\labelenumi}{\textbf{\Roman{enumi})}}
	\item $T_{(Lifespan);U_0} > 2 \TranminusdatasizeWithFactor^{-1}$. 
		In particular, the solution exists classically on the spacetime
		region $\mbox{\upshape cl} \mathcal{M}_{2 \TranminusdatasizeWithFactor^{-1},U_0}$,
		where $\mbox{\upshape cl}$ denotes closure.
		Furthermore, $\inf \lbrace \upmu_{\star}(s,U_0) \ | \ s \in [0,2 \TranminusdatasizeWithFactor^{-1}] \rbrace > 0$.
	\item $T_{(Lifespan);U_0} \leq 2 \TranminusdatasizeWithFactor^{-1}$,
		and 
		\begin{align} \label{E:MAINTHEOREMLIFESPANCRITERION}
		T_{(Lifespan);U_0} 
		= \sup 
			\left\lbrace 
			t \in [0, 2 \TranminusdatasizeWithFactor^{-1}) \ | \
				\inf \lbrace \upmu_{\star}(s,U_0) \ | \ s \in [0,t) \rbrace > 0
			\right\rbrace.
		\end{align}
\end{enumerate}
In addition, case $\textbf{II)}$ occurs when $U_0 = 1$. In this case, we have
\begin{align} \label{E:CLASSICALLIFESPANASYMPTOTICESTIMATE}
	T_{(Lifespan);1} 
	= 
	\left\lbrace 1 + \mathcal{O}(\mathring{\upepsilon}) \right\rbrace
	\TranminusdatasizeWithFactor^{-1}.
\end{align}

\medskip

\noindent \underline{\textbf{What happens in Case I).}}
In case $\textbf{I)}$, 
all bootstrap assumptions,
the estimates of Props.~\ref{P:IMPROVEMENTOFAUX} and \ref{P:IMPROVEMENTOFHIGHERTRANSVERSALBOOTSTRAP},
and the energy estimates of Prop.~\ref{P:MAINAPRIORIENERGY}
hold on $\mbox{\upshape cl} \mathcal{M}_{2 \TranminusdatasizeWithFactor^{-1},U_0}$
with the factors of $\varepsilon$ on the RHS replaced by $C \mathring{\upepsilon}$.
Moreover, for $0 \leq M \leq 5$, the following estimates hold
for $(t,u) \in [0,2 \TranminusdatasizeWithFactor^{-1}] \times [0,U_0]$:
\begin{subequations}
	\begin{align}
		\left\|
			\Tanset_*^{[1,12]} \upmu
		\right\|_{L^2(\Sigma_t^u)},
			\,
		\left\|
			\Tanset^{\leq 12} 
			\Lunit_{(Small)}^i
		\right\|_{L^2(\Sigma_t^u)},
			\,
		\left\|
			\Tanset^{\leq 11} 
			\mytr \upchi
		\right\|_{L^2(\Sigma_t^u)}
		& \leq 
			C \mathring{\upepsilon},
				 \label{E:LOWORDERTANGENTIALEIKONALL2MAINTHEOREM}
				 \\
		\left\|
			\Tanset_*^{13 + M} \upmu
		\right\|_{L^2(\Sigma_t^u)},
			\,
		\left\|
			\Tanset^{13 + M} 
			\Lunit_{(Small)}^i
		\right\|_{L^2(\Sigma_t^u)},
			\,
		\left\|
			\Tanset^{12 + M} 
			\mytr \upchi
		\right\|_{L^2(\Sigma_t^u)}
		& \leq 
			C \mathring{\upepsilon} \upmu_{\star}^{-(M+.4)}(t,u),
				 \label{E:MIDORDERTANGENTIALEIKONALL2MAINTHEOREM} 
				 	\\
		\left\|
			\Lunit \Tanset^{18} \upmu
		\right\|_{L^2(\Sigma_t^u)},
			\,
		\left\|
			\Lunit \Fullset^{18;1} \Lunit_{(Small)}^i
		\right\|_{L^2(\Sigma_t^u)},
			\,
		\left\|
			\Lunit \Fullset^{17;1} \mytr \upchi
		\right\|_{L^2(\Sigma_t^u)}
		& \leq 
			C \mathring{\upepsilon} \upmu_{\star}^{-6.4}(t,u),
				 \label{E:HIGHORDERLUNITTANGENTIALEIKONALL2MAINTHEOREM} 
		\\
		\left\|
			\upmu \GeoAng^{18} \mytr \upchi
		\right\|_{L^2(\Sigma_t^u)}
		& \leq 
			C \mathring{\upepsilon} \upmu_{\star}^{-5.9}(t,u).
				 \label{E:HIGHORDERANGULARTRCHIL2MAINTHEOREM} 
	\end{align}
\end{subequations}

\medskip

\noindent \underline{\textbf{What happens in Case II).}}
In case $\textbf{II)}$, 
all bootstrap assumptions,
the estimates of Props.~\ref{P:IMPROVEMENTOFAUX} and \ref{P:IMPROVEMENTOFHIGHERTRANSVERSALBOOTSTRAP},
and the energy estimates of Prop.~\ref{P:MAINAPRIORIENERGY}
hold on $\mathcal{M}_{T_{(Lifespan);U_0},U_0}$
with the factors of $\varepsilon$ on the RHS replaced by $C \mathring{\upepsilon}$.
Moreover, for $0 \leq M \leq 5$, the estimates 
\eqref{E:LOWORDERTANGENTIALEIKONALL2MAINTHEOREM}-\eqref{E:HIGHORDERANGULARTRCHIL2MAINTHEOREM}
hold for $(t,u) \in [0,T_{(Lifespan);U_0}) \times [0,U_0]$.
In addition, the scalar functions
$\Fullset^{\leq 9;1} \Psi$,
$\Fullset^{\leq 4;2} \Psi$,
$\Rad \Rad \Rad \Psi$,
$\Fullset^{\leq 9;1} \Lunit^i$,
$\Rad \Rad \Lunit^i$,
$\Tanset^{\leq 9} \upmu$,
$\Fullset^{\leq 2;1} \upmu$,
and
$\Rad \Rad \upmu$
extend to 
$\Sigma_{T_{(Lifespan);U_0}}^{U_0}$ 
as functions of 
the geometric coordinates $(t,u,\vartheta)$ 
that are uniformly bounded in $L^{\infty}$.
Furthermore, the rectangular component functions
$g_{\alpha \beta}(\Psi)$ 
verify the estimate
$g_{\alpha \beta} = m_{\alpha \beta} + \mathcal{O}(\mathring{\upepsilon})$ 
(where $m_{\alpha \beta} = \mbox{\upshape diag}(-1,1,1)$ is the standard Minkowski metric)
and have the same extension properties as $\Psi$ and its derivatives with
respect to the vectorfields mentioned above.

Moreover,	let $\Sigma_{T_{(Lifespan);U_0}}^{U_0;(Blowup)}$
be the (non-empty) subset of $\Sigma_{T_{(Lifespan);U_0}}^{U_0}$ 
defined by
\begin{align} \label{E:BLOWUPPOINTS}
	\Sigma_{T_{(Lifespan);U_0}}^{U_0;(Blowup)}
	:= \left\lbrace
			(T_{(Lifespan);U_0},u,\vartheta)
			\ | \
			\upmu(T_{(Lifespan);U_0},u,\vartheta)
			= 0
		\right\rbrace.
\end{align}
Then for each point $(T_{(Lifespan);U_0},u,\vartheta) \in \Sigma_{T_{(Lifespan);U_0}}^{U_0;(Blowup)}$,
there exists a past neighborhood containing it such that the following lower bound holds in
the neighborhood:
\begin{align} \label{E:BLOWUPPOINTINFINITE}
	\left| \Radunit \Psi (t,u,\vartheta) \right|
	\geq \frac{\TranminusdatasizeWithFactor}{4 |G_{\Lunit_{(Flat)} \Lunit_{(Flat)}}(\Psi = 0)|} \frac{1}{\upmu(t,u,\vartheta)}.
\end{align}
In \eqref{E:BLOWUPPOINTINFINITE}, 
$
\displaystyle
\frac{\TranminusdatasizeWithFactor}{4 \left|G_{\Lunit_{(Flat)} \Lunit_{(Flat)}}(\Psi = 0) \right|}
$
is a \textbf{positive} data-dependent constant
(see \eqref{E:NONVANISHINGNONLINEARCOEFFICIENT}),
and the $\ell_{t,u}-$transversal vectorfield $\Radunit$ is near-Euclidean-unit length:
$\delta_{ab} \Radunit^a \Radunit^b = 1 + \mathcal{O}(\mathring{\upepsilon})$.
In particular, $\Radunit \Psi$ blows up like $1/\upmu$ at all points in $\Sigma_{T_{(Lifespan);U_0}}^{U_0;(Blowup)}$.
Conversely, at all points in
$(T_{(Lifespan);U_0},u,\vartheta) \in \Sigma_{T_{(Lifespan);U_0}}^{U_0} \backslash \Sigma_{T_{(Lifespan);U_0}}^{U_0;(Blowup)}$,
we have
\begin{align} \label{E:NONBLOWUPPOINTBOUND}
	\left| \Radunit \Psi (T_{(Lifespan);U_0},u,\vartheta) \right|
	< \infty.
\end{align}

\end{theorem}

\begin{proof}
	Let $C_* > 1$ be a constant (we will adjust $C_*$ throughout the proof).
	We define
	\begin{align}
		T_{(Max);U_0} &:= \mbox{ The supremum of the set of times } \Tboot \in [0,2 \TranminusdatasizeWithFactor^{-1}] \mbox{ such that:} 
			\label{E:LIFESPANPROOF}	\\
			& \bullet \mbox{$\Psi$, $u$, $\upmu$, $\Lunit_{(Small)}^i$, $\Upsilon$,
				and all of the other quantities} 
					\notag \\
			& \ \ \mbox{defined throughout the article
				  exist classically on } \mathcal{M}_{\Tboot,U_0}.
					\notag \\
		& \bullet 
			\mbox{The change of variables map $\Upsilon$ is a (global) $C^{1,1}$ diffeomorphism from }
				\notag \\
		& \ \ 
			\mbox{$[0,\Tboot) \times [0,U_0] \times \mathbb{T}$
			onto its image $\mathcal{M}_{\Tboot,U_0}$.}
				\notag \\
		& \bullet \inf \left\lbrace \upmu_{\star}(t,U_0) \ | \ t \in [0,\Tboot) \right\rbrace > 0.
				\notag \\
		& \bullet \mbox{The fundamental } L^{\infty} \mbox{ bootstrap assumptions } 
						\eqref{E:PSIFUNDAMENTALC0BOUNDBOOTSTRAP}
						\notag \\
		& \ \ \mbox{ hold with } \varepsilon := C_* \mathring{\upepsilon}
								\mbox{ for  } (t,u) \in \times [0,\Tboot) \times [0,U_0].
					\notag \\
			& \bullet \mbox{The $L^2-$type energy bounds }
				\notag \\
		& \ \ 
		\totTanmax{[1,13+M]}^{1/2}(t,u)
		+ 
		\coerciveTanspacetimemax{[1,13+M]}^{1/2}(t,u)
		\leq C_* \mathring{\upepsilon} \upmu_{\star}^{-(M+.9)}(t,u),
			\qquad (0 \leq M \leq 5),
				\label{E:PROOFMULOSSMAINAPRIORIENERGYESTIMATES} \\
		& \ \ 
		\totTanmax{[1,1+M]}^{1/2}(t,u)
		+ \coerciveTanspacetimemax{[1,1+M]}^{1/2}(t,u)
		\leq C_* \mathring{\upepsilon},
		\qquad (0 \leq M \leq 11)
			\label{E:PROOFNOMULOSSMAINAPRIORIENERGYESTIMATES}
								\\
			& \ \ \mbox{ hold for  } (t,u) \in \times [0,\Tboot) \times [0,U_0].
						\notag
		\end{align}
	It is a standard result that
	if $\mathring{\upepsilon}$ is sufficiently small and $C_*$ 
	is sufficiently large, then $T_{(Max);U_0}> 0$
	(this is a standard local well-posedness result
	combined with the initial smallness of the
	$L^2-$controlling quantities shown in Lemma~\ref{L:INITIALSIZEOFL2CONTROLLING}).

	We now show that the energy bounds
	\eqref{E:PROOFMULOSSMAINAPRIORIENERGYESTIMATES}-\eqref{E:PROOFNOMULOSSMAINAPRIORIENERGYESTIMATES}
	and the fundamental $L^{\infty}$ bootstrap assumption \eqref{E:PSIFUNDAMENTALC0BOUNDBOOTSTRAP}
	are not saturated for $(t,u) \in [0,T_{(Max);U_0}) \times [0,U_0]$.
	The non-saturation of the energy bounds (for $C_*$ sufficiently large)
	is provided by Prop.~\ref{P:MAINAPRIORIENERGY}.
	The non-saturation of the fundamental $L^{\infty}$ bootstrap assumptions
	\eqref{E:PSIFUNDAMENTALC0BOUNDBOOTSTRAP} then follows from
	Cor.~\ref{C:PSILINFTYINTERMSOFENERGIES}.
	Consequently, we conclude that all of the estimates proved throughout the article
	hold on $\mathcal{M}_{\Tboot,U_0}$ with
	the smallness parameter $\varepsilon$ replaced by $C \mathring{\upepsilon}$.
	We use this fact throughout the remainder of the proof without further remark. 

	Next, we show that \eqref{E:LOWORDERTANGENTIALEIKONALL2MAINTHEOREM}-\eqref{E:HIGHORDERANGULARTRCHIL2MAINTHEOREM}
	hold for $(t,u) \in [0,T_{(Max);U_0}) \times [0,U_0]$.
	To obtain \eqref{E:LOWORDERTANGENTIALEIKONALL2MAINTHEOREM}-\eqref{E:HIGHORDERLUNITTANGENTIALEIKONALL2MAINTHEOREM},
	we insert the energy estimates of Prop.~\ref{P:MAINAPRIORIENERGY}
	into the RHS of the inequalities of Lemma~\ref{L:EASYL2BOUNDSFOREIKONALFUNCTIONQUANTITIES}
	and use inequalities \eqref{E:LOSSKEYMUINTEGRALBOUND} and \eqref{E:LESSSINGULARTERMSMPOINTNINEINTEGRALBOUND}
	as well as the fact that $\totTanmax{[1,M]}$ is increasing in its arguments.
	Similarly, to obtain inequality \eqref{E:HIGHORDERANGULARTRCHIL2MAINTHEOREM}, 
	we insert the energy estimates of Prop.~\ref{P:MAINAPRIORIENERGY}
	into RHS~\eqref{E:LESSPRECISEDIFFICULTTERML2BOUND}
	and use inequality
	\eqref{E:LOSSKEYMUINTEGRALBOUND}.

	We now establish the dichotomy of possibilities.
	We first show that
	if 
	\[
		\inf \left\lbrace \upmu_{\star}(t,U_0) \ | \ t \in [0,T_{(Max);U_0}) \right\rbrace > 0,
	\]
	then $T_{(Max);U_0} = 2 \TranminusdatasizeWithFactor^{-1}$.
	To proceed, we assume for the sake of contradiction that
	the previous bound for $\upmu_{\star}$ holds but that
	$T_{(Max);U_0} < 2 \TranminusdatasizeWithFactor^{-1}$.
	To reach a contradiction, 
	we will use Lemmas~\ref{L:CHOVREMAINSADIFFEOMORPHISM}
	and \ref{L:CONTINUATIONCRITERIA}
	to deduce that we can classically extend the solution
	to a region of the form $\mathcal{M}_{T_{(Max);U_0} + \Delta,U_0}$,
	with $\Delta > 0$
	and
	$T_{(Max);U_0} + \Delta 
		< 2 \TranminusdatasizeWithFactor^{-1}
	$,
	such that all of the properties defining $T_{(Max);U_0}$ hold
	for the larger time
	$T_{(Max);U_0} + \Delta$.
	Since we have already shown that the energy bounds
	\eqref{E:PROOFMULOSSMAINAPRIORIENERGYESTIMATES}-\eqref{E:PROOFNOMULOSSMAINAPRIORIENERGYESTIMATES}
	are not saturated and that
	the fundamental $L^{\infty}$ bootstrap assumption \eqref{E:PSIFUNDAMENTALC0BOUNDBOOTSTRAP}
	are not saturated for $(t,u) \in [0,T_{(Max);U_0}) \times [0,U_0]$,
	the contradiction will follow once we show that
	the change of variables map $\Upsilon$ extends as a
	global $C^{1,1}$ diffeomorphism from
	$[0,T_{(Max);U_0}] \times [0,U_0] \times \mathbb{T}$
	onto its image
	and that none of the four breakdown scenarios
	of Lemma~\ref{L:CONTINUATIONCRITERIA} occur
	on $\mathcal{M}_{T_{(Max);U_0},U_0}$.
	Breakdown scenario $(1)$ 
	from Lemma~\ref{L:CONTINUATIONCRITERIA}
	is ruled out by assumption.
	Scenario $(2)$ is ruled out by the estimate
	\eqref{E:UPMULINFTY}. Scenario $(3)$ is ruled out 
	by the bootstrap assumptions \eqref{E:PSIFUNDAMENTALC0BOUNDBOOTSTRAP}
	and the fact that 
	$\gt_{ij} = \delta_{ij} + \mathcal{O}(\Psi)$, with $\delta_{ij}$
	the standard Kronecker delta.
	From Lemma~\ref{L:CHOVREMAINSADIFFEOMORPHISM},
	we obtain that $\Upsilon$ 
	extends as a
	global $C^{1,1}$ diffeomorphism from
	$[0,T_{(Max);U_0}] \times [0,U_0] \times \mathbb{T}$
	onto its image.
	Hence, 
	we can rule out the scenario $(4)$ 
	once we show that 
	\[
	\displaystyle
	\sup_{(t,u,\vartheta) \in [0,T_{(Max);U_0}] \times [0,U_0] \times \mathbb{T}}
		\sum_{i_1 + i_2 + i_3 \leq 1}
		\left| 
			\left(\frac{\partial}{\partial t}\right)^{i_1} 
			\left(\frac{\partial}{\partial u}\right)^{i_2}
			\left(\frac{\partial}{\partial \vartheta}\right)^{i_3}
			\Psi(t,u,\vartheta) 
		\right| 
	< \infty.
	\]
	This desired bound is a simple consequence of 
	the estimates \eqref{E:PSITRANSVERSALLINFINITYBOUNDBOOTSTRAPIMPROVED}-\eqref{E:PSIMIXEDTRANSVERSALTANGENTBOOTSTRAPIMPROVED}
	(which hold for $(t,u) \in [0,T_{(Max);U_0}] \times [0,U_0]$)
	and the fact that, 
	as we showed in the proof of Lemma~\ref{L:CHOVMAPREGULARITY}, 
	we have
	$\Lunit = \frac{\partial}{\partial t}$,
	$\GeoAng = (1 + \mathcal{O}(\mathring{\upepsilon})) \CoordAng$ 
	(recall that $\CoordAng := \frac{\partial}{\partial \vartheta}$)
	and $\frac{\partial}{\partial u} = \Rad + \mathcal{O}(1) \GeoAng$.
	We have thus reached a contradiction
	and established that either
	\textbf{I)} $T_{(Max);U_0} = 2 \TranminusdatasizeWithFactor^{-1}$
	or \textbf{II)} $\inf \left\lbrace \upmu_{\star}(t,U_0) \ | \ t \in [0,T_{(Max);U_0}) \right\rbrace = 0$.

	We now show that case \textbf{II)} corresponds to a singularity
	and that the classical lifespan is characterized by \eqref{E:MAINTHEOREMLIFESPANCRITERION}.
	To this end, we first use
	\eqref{E:LUNITUPMUDOESNOTDEVIATEMUCHFROMTHEDATA},
	\eqref{E:SMALLMUIMPLIESLMUISNEGATIVE},
	and the identity $\Rad = \upmu \Radunit$
	to deduce that inequality \eqref{E:BLOWUPPOINTINFINITE} holds.
	Furthermore, from 
	\eqref{E:LITTLEGDECOMPOSED},
	\eqref{E:DOWNSTAIRSUPSTAIRSSRADUNITPLUSLUNITISAFUNCTIONOFPSI},
	and the $L^{\infty}$ estimates of Prop.~\ref{P:IMPROVEMENTOFAUX},
	we deduce that
	$|\Radunit| 
	: =
	\sqrt{g_{ab} \Radunit^a \Radunit^b} 
	= 1 + \smoothfunction(\GdVar) \GdVar 
	=  1 
	+ \mathcal{O}(\mathring{\upepsilon})
	$.
	From this estimate and \eqref{E:BLOWUPPOINTINFINITE}, 
	we deduce that at points in $\Sigma_{T_{(Max);U_0},U_0}$
	where $\upmu$ vanishes, $|\Radunit \Psi|$ must blow up like $1/\upmu$.
	Hence, $T_{(Max);U_0}$ is the classical lifespan. That is, 
	we have $T_{(Max);U_0} = T_{(Lifespan);U_0}$
	as well as the characterization \eqref{E:MAINTHEOREMLIFESPANCRITERION} of the 
	classical lifespan. The estimate \eqref{E:NONBLOWUPPOINTBOUND}
	is an immediate consequence of
	the estimate \eqref{E:PSITRANSVERSALLINFINITYBOUNDBOOTSTRAPIMPROVED}
	and the identity $\Rad = \upmu \Radunit$.

	To obtain \eqref{E:CLASSICALLIFESPANASYMPTOTICESTIMATE},
	we use
	\eqref{E:MUSTARBOUNDS}
	and
	\eqref{E:CRUCIALLATETIMEDERIVATIVECOMPAREDTODATAPARAMETER}
	to deduce that $\upmu_{\star}(t,1)$
	vanishes for the first time when
	$t = 
	\left\lbrace 1 + \mathcal{O}(\mathring{\upepsilon}) \right\rbrace
	\TranminusdatasizeWithFactor^{-1}
	$.
	
	We now derive the statements regarding the quantities that extend to 
	$\Sigma_{T_{(Lifespan);U_0}}^{U_0}$ as $L^{\infty}$ functions.
	Let $q$ denote any of the quantities
	$\Fullset^{\leq 9;1} \Psi$,
	$\cdots$,
	$\Rad \Rad \upmu$
	that are stated in the theorem 
	to extend to $\Sigma_{T_{(Lifespan);U_0}}^{U_0}$ as an $L^{\infty}$ function
	of the geometric coordinates.
	The $L^{\infty}$ estimates of 
	Props.~\ref{P:IMPROVEMENTOFAUX}
	and \ref{P:IMPROVEMENTOFHIGHERTRANSVERSALBOOTSTRAP}
	imply that $\| \Lunit q \|_{L^{\infty}(\Sigma_t^{U_0})}$
	is uniformly bounded
	for $0 \leq t < T_{(Lifespan);U_0}$.
	Recalling that $\Lunit = \frac{\partial}{\partial t}$, we
	conclude that $q$ extends to $\Sigma_{T_{(Lifespan);U_0}}^{U_0}$
	as an element of $L^{\infty}(\Sigma_{T_{(Lifespan);U_0}}^{U_0})$
	as desired.
	The estimate 
	$g_{\alpha \beta} = m_{\alpha \beta} + \mathcal{O}(\mathring{\upepsilon})$
	and the extension properties of the $\Fullset-$derivatives of
	the scalar functions $g_{\alpha \beta}$ then follow from
	\eqref{E:LITTLEGDECOMPOSED},
	the already proven bound $\| \Psi \|_{L^{\infty}(\Sigma_t^{U_0})} \lesssim \mathring{\upepsilon}$,
	and the above extension properties of the $\Fullset-$derivatives of $\Psi$.

\end{proof}

\section*{Acknowledgments}
The authors thank the American Institute of Mathematics for funding three SQuaREs workshops
on the formation of shocks, which greatly furthered the development of many of the ideas
in this paper. They thank Sergiu Klainerman and Shiwu Yang for participating in the workshops 
and for their helpful contributions, and they are also grateful for the helpful suggestions offered by 
Jacques Smulevici. JS and WW gratefully acknowledge support from the Simons Center for Geometry and Physics, 
Stony Brook University, at which some of the research for this paper was performed.
Finally, we thank the four anonymous referees, who offered many useful insights 
that helped us improve and clarify various aspects of this work.

\appendix
\setcounter{equation}{0}
\numberwithin{equation}{subsection}

\section{Extending the results to 
\texorpdfstring{the equations $(g^{-1})^{\alpha \beta}(\partial \Phi) \partial_{\alpha} \partial_{\beta} \Phi = 0$}{non-covariant wave equations}
}
\label{A:EXTENDINGTONONGEOMETRICWAVEEQUATIONS}
In this appendix, we sketch how to extend our shock formation results to the Cauchy problem
\begin{subequations}
\begin{align} \label{E:APPENDIXNONGEOMETRICWAVE}
	(g^{-1})^{\alpha \beta}(\partial \Phi)
	\partial_{\alpha} \partial_{\beta} \Phi
	& = 0,
		\\
	(\Phi|_{\Sigma_0},\partial_t \Phi|_{\Sigma_0}) 
	&= (\mathring{\Phi},\mathring{\Phi}_0),
		\label{E:APPENDIXDATAFORNONGEOMETRICWAVE}
\end{align}
where equation \eqref{E:APPENDIXNONGEOMETRICWAVE} is written 
relative to the rectangular coordinates $\lbrace x^{\alpha} \rbrace_{\alpha = 0,1,2}$,
and
\begin{align} \label{E:APPENDIXMETRICEXPANDAROUNDMINKOWSKI}
	g_{\alpha \beta}(\partial \Phi) 
	& = m_{\alpha \beta}
		+ g_{\alpha \beta}^{(Small)}(\partial \Phi),
	\qquad g_{\alpha \beta}^{(Small)}(0) = 0.
\end{align}
Dividing the wave equation by 
$-(g^{-1})^{00}$ if necessary, we may assume as before that
\begin{align} \label{E:APPENDIXGINVERSE00ISMINUSONE}
	(g^{-1})^{00} 
	\equiv 
	- 1.
\end{align}
\end{subequations}

\subsection{Basic setup}
\label{SS:BASICSTEUP}
We start by defining $\Psi_{\nu}$, $(\nu = 0,1,2)$, and $\vec{\Psi}$ as follows:
	\begin{align} \label{E:APPENDIXPSICOMPONENTDEF} 
		\Psi_{\nu} 
		& := \partial_{\nu} \Phi,
			\qquad
		\vec{\Psi} 
		:= (\Psi_0,\Psi_1,\Psi_2).
	\end{align}

	The main strategy behind extending our results is to take rectangular derivatives of the equation
	\eqref{E:APPENDIXNONGEOMETRICWAVE} to form a system of wave equations
	in the unknowns $\vec{\Psi}$; see Lemma~\ref{L:ECOMMUTEDINTERMSOFHCON}. 
	The system has a special null structure that plays an important role in the analysis;
	see Lemma~\ref{L:SPECIALNULLSTRUCTUREINHOMOGENEOUS}.
	The vast majority of the proof of shock formation for the system 
	is the same as it is in the case of the 
	scalar equation \eqref{E:GEOWAVE}, but now with $\vec{\Psi}$ in the role of $\Psi$.
	We can treat the system
	using essentially the same methods that we used to treat the scalar equation \eqref{E:GEOWAVE}
	because the coupling between the $\Psi_{\nu}$ is not very difficult to handle 
	and because the tensorial
	structure of the equations matters only in a few key places.
	We devote the remainder of this appendix to highlighting those key places and to
	describing the handful of new ingredients that are needed. 

We first note that the analogs of the scalar functions \eqref{E:BIGGDEF} for equation \eqref{E:APPENDIXNONGEOMETRICWAVE} are
\begin{subequations}
\begin{align} \label{E:APPENDIXBIGGDEF}
		G_{\mu \nu}^{\lambda}
		= G_{\mu \nu}^{\lambda}(\vec{\Psi})
		& := \frac{\partial g_{\mu \nu}^{(Small)}(\vec{\Psi})}{\partial \Psi_{\lambda}},
			\\
		G^{\mu \alpha \beta}(\vec{\Psi})
	& := (g^{-1})^{\alpha \alpha'} (g^{-1})^{\beta \beta'} G_{\alpha' \beta'}^{\mu}.
		\label{E:APPENDIXRAISEDBIGGDEF}
\end{align}
\end{subequations}
For our proof to work, we assume 
an analog of \eqref{E:NONVANISHINGNONLINEARCOEFFICIENT},
specifically that there exist coordinates such that
$m_{\alpha \beta} = \mbox{\upshape diag}(-1,1,1)$
(that is, Minkowski-rectangular coordinates)
and such that with $\Lunit_{(Flat)} := \partial_t + \partial_1$, we have
\begin{align}  \label{E:APPENDIXNULLCONDITIONFAILS}
	m_{\kappa \lambda}
	G_{\mu \nu}^{\kappa}(\vec{\Psi} = 0)
	\Lunit_{(Flat)}^{\alpha} \Lunit_{(Flat)}^{\beta} \Lunit_{(Flat)}^{\lambda}
	\neq 0.
\end{align}
The assumption ensures that in the regime under study,
the term $\upomega^{(Trans-\vec{\Psi})} $
on RHS~\eqref{E:NEWUPMUTRANSPORT} is sufficiently strong to drive
$\upmu$ to $0$ in finite time.

We now provide the system of covariant wave equations
implied by equation \eqref{E:APPENDIXNONGEOMETRICWAVE}.
The proof is a straightforward but tedious computation that relies on the
identity $\partial_{\alpha} \Psi_{\beta} = \partial_{\beta} \Psi_{\alpha}$;
we omit the details.

\begin{lemma} [\textbf{The system of covariant wave equations}] 
\label{L:ECOMMUTEDINTERMSOFHCON}
As a consequence of equation \eqref{E:APPENDIXNONGEOMETRICWAVE}, 
the quantities $\Psi_{\nu} := \partial_{\nu} \Phi$
verify the following system of covariant wave equations
(where $\Psi_{\nu}$ is viewed to be a scalar-valued function under covariant differentiation):
\begin{align} \label{E:PSIGEOMETRIC}
	\square_{g(\vec{\Psi})} \Psi_{\nu} 
	& = \mathscr{Q}(\partial \vec{\Psi}, \partial \Psi_{\nu}), & (\nu & = 0,1,2),
\end{align}
where
$
\square_{g(\vec{\Psi})} \Psi := \frac{1}{\sqrt{|\mbox{\upshape{det}$g$}|}} 
\partial_{\alpha} \left(\sqrt{|\mbox{\upshape{det}$g$} |} (g^{-1})^{\alpha \beta} \partial_{\beta} \Psi \right)
$
is the covariant wave operator of $g$ applied to $\Psi$,
\begin{align} \label{E:APPENDIXCOVARIANTSYSTEMQUADRATICTERM}
	\mathscr{Q}(\partial \vec{\Psi}, \partial \Psi) 
	& :=	G^{\mu \alpha \beta}
				\left\lbrace
					\partial_{\beta} \Psi_{\alpha} \partial_{\mu} \Psi
					- \partial_{\mu} \Psi_{\alpha} \partial_{\beta} \Psi
				\right\rbrace
			+ (g^{-1})^{\alpha \beta} \Omega^{\lambda} \partial_{\alpha} \Psi_{\lambda} 
			\partial_{\beta} \Psi,
\end{align}
$G^{\mu \alpha \beta}$ is defined in \eqref{E:APPENDIXRAISEDBIGGDEF},
and
$
	\displaystyle
	\Omega^{\nu}(\vec{\Psi})
		:= \frac{1}{\sqrt{|\mbox{\upshape{det}$g$}|(\vec{\Psi})}} 
			\frac{\partial \sqrt{|\mbox{\upshape{det}$g$}|(\vec{\Psi})}}{\partial \Psi_{\nu}}
$.
\end{lemma}

$\hfill \qed$

The quadratic term $\mathscr{Q}$ on RHS~\eqref{E:PSIGEOMETRIC} has a special null
structure that is of critical importance for our proof.
We describe this structure in Lemma~\ref{L:SPECIALNULLSTRUCTUREINHOMOGENEOUS}
below. We first recall the definitions of the standard null (relative to $g$) forms 
$\mathscr{Q}_{(0)}$ and $\mathscr{Q}_{(\alpha \beta)}$:
	\begin{subequations}
	\begin{align}
		\mathscr{Q}_{(0)}(\partial \phi, \partial \widetilde{\phi})
		& := (g^{-1})^{\alpha \beta} \partial_{\alpha} \phi \partial_{\beta} \widetilde{\phi},
			\label{E:STANDARD0NULLFORM} \\
		\mathscr{Q}_{(\alpha \beta)}(\partial \phi, \partial \widetilde{\phi})
		& := \partial_{\alpha} \phi \partial_{\beta} \widetilde{\phi}
			- \partial_{\alpha} \widetilde{\phi} \partial_{\beta} \phi.
			\label{E:STANDARDALPHABETANULLFORM}
	\end{align}
	\end{subequations}

In the next lemma, 
we decompose the standard null forms
relative to the non-rescaled frame \eqref{E:UNITFRAME}
and exhibit their good geometric properties 
from the point of view of the shock formation problem.
The main point is that there is no term proportional to
$(\Radunit \phi) \Radunit \widetilde{\phi}$ on RHS~\eqref{E:STANDARDNULLFORMSTRUCTURE}.

	\begin{lemma}[\textbf{Good properties of the standard null forms}]
	\label{L:STANDARDNULLFORMSTRUCTURE}
	If $\mathscr{Q}$ is a standard null form, then we can decompose
	it as follows relative to the non-rescaled frame \eqref{E:UNITFRAME}:
	\begin{align} \label{E:STANDARDNULLFORMSTRUCTURE}
		\mathscr{Q}(\partial \phi, \partial \widetilde{\phi})
		& = f_1 (\Lunit \phi) \Lunit \widetilde{\phi}
			+ f_2 (\Lunit \phi) \Radunit \widetilde{\phi}
			+ f_3 (\Radunit \phi) \Lunit \widetilde{\phi}
				 \\
		& \ \ 
			+ (f_4 \cdot \angdiff \phi) \Lunit \widetilde{\phi}
			+ (f_5 \cdot \angdiff \phi) \Radunit \widetilde{\phi}
			+ (\Lunit \phi) f_6 \cdot \angdiff \widetilde{\phi} 
			+ (\Radunit \phi) f_7 \cdot \angdiff \widetilde{\phi} 
			+ f_8 \cdot \angdiff \phi \otimes \angdiff \widetilde{\phi},
			\notag
	\end{align}
	where $f_1$, $f_2$ and $f_3$ are scalar functions,
	$f_4$, $f_5$, $f_6$, $f_7$ are $\ell_{t,u}-$tangent vectorfields,
	and $f_8$ is a symmetric type $\binom{2}{0}$ $\ell_{t,u}-$tangent tensorfield
	with the following properties:
	$f_1$, $f_2$ and $f_3$
	and the rectangular components 
	$f_4^{\alpha}$, $f_5^{\alpha}$, $f_6^{\alpha}$, $f_7^{\alpha}$,
	and $f_8^{\alpha \beta}$
	are smooth scalar-valued functions of $\vec{\Psi}$
	and the rectangular components of the vectorfields $\Lunit$ and $\Radunit$.
\end{lemma}

\begin{proof}
	When $\mathscr{Q} = \mathscr{Q}_{(0)}$, \eqref{E:STANDARDNULLFORMSTRUCTURE} follows from
	Lemma~\ref{L:METRICDECOMPOSEDRELATIVETOTHEUNITFRAME}.
	When $\mathscr{Q} = \mathscr{Q}_{(\alpha \beta)}$,
	we view $\mathscr{Q}_{(\alpha \beta)}$ to be the rectangular components 
	of an anti-symmetric type $\binom{0}{2}$ spacetime tensor 
	which we decompose relative to the non-rescaled frame:
	$\mathscr{Q}_{(\alpha \beta)} 
	= 
		F_{\Lunit \Radunit}
		(
		\Lunit_{\alpha} \Radunit_{\beta}
		 -
		 \Radunit_{\alpha} \Lunit_{\beta}
		)
		+
		F_{\Lunit \CoordAng}
		(
		\Lunit_{\alpha} \CoordAng_{\beta}
		 -
		\CoordAng_{\alpha} \Lunit_{\beta}
		)
		+
		F_{\Radunit \CoordAng}
		(
		\Radunit_{\alpha} \CoordAng_{\beta}
		 -
		\CoordAng_{\alpha} \Radunit_{\beta}
		)
	$,
	where the $F_{\cdots}$ are scalar functions.
	To compute the $F_{\cdots}$, we contract both sides of the identity against pairs of elements
	of the non-rescaled frame $\lbrace \Lunit, \Radunit, \CoordAng \rbrace$.
	For example, contracting against 
	$\Lunit^{\alpha} \CoordAng^{\beta}$
	and using \eqref{E:METRICANGULARCOMPONENT}, 
	we find that 
	$(\Lunit \phi) \CoordAng \widetilde{\phi}
	 -
	 (\Lunit \widetilde{\phi}) \CoordAng \phi
	=
	- F_{\Radunit \CoordAng} \gtancomp^2
	$.
	This leads to a decomposition of the form
	$
	\mathscr{Q}_{(\alpha \beta)}
	= 
	\cdots
  -
  \gtancomp^{-2}
	\left\lbrace
		(\Lunit \phi) \CoordAng \widetilde{\phi}
		-
		(\Lunit \widetilde{\phi}) \CoordAng \phi
	\right\rbrace
  \left\lbrace
		\Radunit_{\alpha} \CoordAng_{\beta}
		 -
		\CoordAng_{\alpha} \Radunit_{\beta}
	\right\rbrace
	$.
	Using \eqref{E:METRICANGULARCOMPONENT},
	we can rewrite terms involving $\CoordAng$
	as in the following example:
	$\gtancomp^{-2} 
	(\Lunit \phi) (\CoordAng \widetilde{\phi})
	\Radunit_{\alpha} \CoordAng_{\beta}
	=
	(\Lunit \phi) (\Lineproject_{\beta}^{\ \#} \cdot \angdiff \widetilde{\phi})
	\Radunit_{\alpha} 
	$,
	where $\Lineproject_{\beta}^{\ \lambda}$ is defined in \eqref{E:LINEPROJECTION}.
	The desired decomposition
	\eqref{E:STANDARDNULLFORMSTRUCTURE}
	thus follows.
\end{proof}

In the next lemma, we characterize the good structure of the quadratic term $\mathscr{Q}$ on RHS~\eqref{E:PSIGEOMETRIC}.
The proof follows from observation.
\begin{lemma}[\textbf{Special null structure of the inhomogeneous terms}]
	\label{L:SPECIALNULLSTRUCTUREINHOMOGENEOUS}
	The quadratic term $\mathscr{Q}(\partial \vec{\Psi}, \partial \Psi)$ on the right-hand side of 
	\eqref{E:PSIGEOMETRIC} is a linear combination 
	of the standard null forms in $\partial \vec{\Psi}$ with coefficients depending on $\vec{\Psi}$.
\end{lemma}

$\hfill \qed$

\subsection{Additional smallness assumptions in the present context}
\label{SS:NEWSMALLNESS}
	To close the proof of shock formation for solutions to the system \eqref{E:PSIGEOMETRIC}, 
	we assume that each scalar function $\Psi_{\nu}$ has data verifying the same size assumptions as the 
	data for the scalar function $\Psi$, as described in Subsects.~\ref{SS:DATAASSUMPTIONS} and \ref{SS:SMALLNESSASSUMPTIONS}.
	Similarly, to derive estimates, we make the same $L^{\infty}$ bootstrap assumptions 
	for each $\Psi_{\nu}$ that we did for $\Psi$. 
	As we show below in 
	\eqref{E:SECONDEQUATIONNEWKINDOFSMALLNESSNEEDED}
	and the discussion surrounding
	\eqref{E:NEWSMALLNESSIDENTITYONE}-\eqref{E:NEWSMALLNESSIDENTITYTWO},
	these assumptions impose some subtle smallness restrictions on the data \eqref{E:APPENDIXDATAFORNONGEOMETRICWAVE}
	in the sense that they imply the $\mathcal{O}(\mathring{\upepsilon})$ smallness of 
	special combinations of the elements of 
	$\lbrace \Rad \Psi_0, \Rad \Psi_1, \Rad \Psi_2 \rbrace$.
	These smallness restrictions are consequences of our 
	size assumptions on the $\Psi_{\nu}$ and their derivatives
	and the symmetry property
	$
	\partial_{\alpha} \Psi_{\beta}
	=
	\partial_{\beta} \Psi_{\alpha}
	$. 
	As we will see in Subsect.~\ref{SS:MAINNEWESTIMATE}, 
	we especially rely on the following small-data estimates:
	\begin{align} \label{E:NEWPROBLEMADDITIONALSMALLNESSASSUMPTIONS}
		\left\|
			\Lunit \Rad \Phi
		\right\|_{L^{\infty}(\Sigma_0^1)},
			\,
		\left\|
			\GeoAng \Rad \Phi
		\right\|_{L^{\infty}(\Sigma_0^1)}
		& \lesssim \mathring{\upepsilon}.
	\end{align}
	We note that the $\mathcal{O}(\mathring{\upepsilon})$ smallness of
	$\left\|
		\Lunit \Rad \Phi
	\right\|_{L^{\infty}(\Sigma_0^1)}
	$
	is a simple consequence of the identity
	\begin{align} \label{E:LUNITRADPHISIMPLEID}
		\Lunit \Rad \Phi
		& = \Lunit (\upmu \Radunit^a \Psi_a) 
			= (\Lunit \upmu) \Radunit^a \Psi_a 
				+
				\upmu (\Lunit \Radunit_{(Small)}^a) \Psi_a
				+
				\upmu \Radunit^a \Lunit \Psi_a
	\end{align}
	and the $\mathcal{O}(\mathring{\upepsilon})$ smallness of 
	$
	\left\|
	 	\Psi_a
	\right\|_{L^{\infty}(\Sigma_0^1)}
	$
	and
	$
	\left\|
	 	\Lunit \Psi_a
	\right\|_{L^{\infty}(\Sigma_0^1)}
	$.
	Similar remarks apply to
	the term
	$
		\left\|
			\GeoAng \Rad \Phi
		\right\|_{L^{\infty}(\Sigma_0^1)}
	$
	on LHS~\eqref{E:NEWPROBLEMADDITIONALSMALLNESSASSUMPTIONS}.
	It is of course important that
	the smallness conditions \eqref{E:NEWPROBLEMADDITIONALSMALLNESSASSUMPTIONS} are propagated
	by the nonlinear flow.
	Specifically,
	in the analog of the proof of Prop.~\ref{P:IMPROVEMENTOFAUX},
	we could derive the estimates
	\begin{align} \label{E:NEWPROBLEMADDITIONALSMALNESSESTIMATE}
		\left\|
			\Lunit \Rad \Phi
		\right\|_{L^{\infty}(\Sigma_t^u)},
			\,
		\left\|
			\GeoAng \Rad \Phi
		\right\|_{L^{\infty}(\Sigma_t^u)}
		& \lesssim \varepsilon
	\end{align}
	at the end of the proof. For example, the estimate \eqref{E:NEWPROBLEMADDITIONALSMALNESSESTIMATE} for 
	$
	\left\|
		\Lunit \Rad \Phi
	\right\|_{L^{\infty}(\Sigma_t^u)}
	$
	would follow from the identity \eqref{E:LUNITRADPHISIMPLEID}
	and $L^{\infty}$ estimates for all of the terms on RHS~\eqref{E:LUNITRADPHISIMPLEID}, 
	which would already have been obtained in the proof of the proposition.

\subsection{The main new estimate needed at the top order}
\label{SS:MAINNEWESTIMATE}
We now explain how to extend Theorem~\ref{T:MAINTHEOREM}
to the system \eqref{E:PSIGEOMETRIC}. As we have suggested above, 
we can derive energy identities for each scalar function
$\Psi_{\nu}$
by using essentially the same arguments that we used to  
treat the scalar equation \eqref{E:GEOWAVE}. 
To derive inequalities that control $\vec{\Psi}$, 
we replace the controlling quantity $\totTanmax{N}$
from Def.~\ref{D:MAINCOERCIVEQUANT}
with
\begin{align} \label{E:NEWMAXEDENERGIESWITHOUTTENSORIALSTRUCTURE}
	\totTanmax{N}(t,u)
	& := \max_{|\vec{I}| = N} \max_{\nu=0,1,2} 
	\sup_{(t',u') \in [0,t] \times [0,u]} 
		\left\lbrace
			\enzero[\Tanset^{\vec{I}} \Psi_{\nu}](t',u')
			+ 
			\flzero[\Tanset^{\vec{I}} \Psi_{\nu}](t',u')
		\right\rbrace,
\end{align}
and similarly for the other controlling quantities of Sect.~\ref{S:FUNDAMENTALL2CONTROLLINGQUANTITIES}.

Thanks to Lemma~\ref{L:SPECIALNULLSTRUCTUREINHOMOGENEOUS},
the terms on RHS~\eqref{E:PSIGEOMETRIC} 
are easy to treat without invoking any new ideas.
In the remainder of this appendix,
we explain the one new ingredient that we
need to close the estimates. It is needed
for the top-order $L^2$ estimates for the $\Psi_{\nu}$.
To motivate the discussion, we first
recall a critically important aspect of our analysis of the scalar equation \eqref{E:GEOWAVE}.
At several points in our argument for deriving top-order $L^2$ estimates for solutions to \eqref{E:GEOWAVE},
we had to use equation \eqref{E:UPMUFIRSTTRANSPORT},
the fact that 
$G_{\Lunit \Lunit}, G_{\Lunit \Radunit} = \smoothfunction(\GdVar)$
(see Lemma~\ref{L:SCHEMATICDEPENDENCEOFMANYTENSORFIELDS}),
and the $L^{\infty}$ estimates of Prop.~\ref{P:IMPROVEMENTOFAUX} to obtain 
\begin{align} \label{E:KEYREPLACE}
	\left|G_{\Lunit \Lunit} \Rad \Psi \right| \leq 2 |\Lunit \upmu| 
		+ \upmu \mathcal{O}(\varepsilon).
\end{align}
For example, \eqref{E:KEYREPLACE} was used to derive\footnote{Actually, in deriving \eqref{E:INTROREPRESENTATIVEERRORINTEGRAL}, 
we used a version of \eqref{E:KEYREPLACE} in which the absolute value signs are missing and ``$\leq$'' is replaced with ``$=$.''
However, \eqref{E:KEYREPLACE} would have been sufficient for all of the arguments to go through.} 
equation \eqref{E:INTROREPRESENTATIVEERRORINTEGRAL}.
Since we are treating the coupled system $\vec{\Psi}$
by separately deriving energy identities for each scalar function $\Psi_{\nu}$,
our energy estimates rely on the following analog of \eqref{E:KEYREPLACE} 
\emph{for each of the three} $\Psi_{\nu}$:
\begin{align} \label{E:FIRSTSTATEMENTHARDSHARPRADPSIPOINTWISEESTIMATE}
	\left|G_{\Lunit \Lunit}^{\Lunit} \Rad \Psi_{\nu} \right| 
	& \leq 2 \left|\Lunit \upmu \right| 
	+ \upmu \mathcal{O}(\varepsilon),
\end{align}
where $G_{\Lunit \Lunit}^{\Lunit}:= G_{\alpha \beta}^{\kappa} \Lunit^{\alpha} \Lunit^{\beta} \Lunit_{\kappa}$.
The estimate \eqref{E:FIRSTSTATEMENTHARDSHARPRADPSIPOINTWISEESTIMATE} is the main new ingredient
that we need at the top order.
As we will see, it does not follow
directly from the evolution equation $\Lunit \upmu = \cdots$ 
and instead relies on a few new tensorial observations
and the estimate \eqref{E:NEWPROBLEMADDITIONALSMALNESSESTIMATE}.
Thus, we dedicate
the remainder of this appendix to sketching a proof of \eqref{E:FIRSTSTATEMENTHARDSHARPRADPSIPOINTWISEESTIMATE}.

We start by providing the evolution equation for $\upmu$ in the present context.

\begin{lemma} [\textbf{The transport equation verified by $\upmu$}] \label{L:NEWUPMUTRANSPORT}
In the case of equation \eqref{E:APPENDIXNONGEOMETRICWAVE},
the inverse foliation density $\upmu$ defined in \eqref{E:UPMUDEF} verifies the following transport equation:
\begin{align} \label{E:NEWUPMUTRANSPORT}
	\Lunit \upmu 
	& = \upomega^{(Trans-\vec{\Psi})} 
		+ \upmu \upomega^{(Tan-\vec{\Psi})}
		:= \upomega,
\end{align}
where
\begin{subequations}
\begin{align}
	\upomega^{(Trans-\vec{\Psi})} 
	& := - \frac{1}{2} G_{\Lunit \Lunit}^{\Lunit} \Radunit^a \Rad \Psi_a,
		\label{E:LITTLEOMEGATRANSVERSALPART} \\
	\upomega^{(Tan-\vec{\Psi})} 
	& := 
		- \frac{1}{2} G_{\Lunit \Lunit}^{\Lunit} \Radunit^a \Lunit \Psi_a
		- \frac{1}{2} G_{\Lunit \Lunit}^{\Radunit} \Radunit^a \Lunit \Psi_a
		+ \frac{1}{2} \angG_{\Lunit \Lunit}^{\#} \cdot \Radunit^a \angdiff \Psi_a
		- \frac{1}{2} G_{\Lunit \Lunit}^{\lambda} \Lunit \Psi_{\lambda}
		- G_{\Lunit \Radunit}^{\lambda} \Lunit \Psi_{\lambda}.
		\label{E:LITTLEOMEGATANGENTIALPART}
\end{align}
\end{subequations}
The scalar-valued function $G_{\mu \nu}^{\lambda}$ above is defined in \eqref{E:APPENDIXBIGGDEF},
$G_{\Lunit \Lunit}^{\Lunit}:= G_{\alpha \beta}^{\kappa} \Lunit^{\alpha} \Lunit^{\beta} \Lunit_{\kappa}$,
$\angG_{\Lunit \Lunit}^{\lambda} 
:= G_{\alpha \beta}^{\kappa} \Lunit^{\alpha} \Lunit^{\beta} \Lineproject_{\kappa}^{\ \lambda}$,
etc.
\end{lemma}
\begin{proof}
	The proof is very similar to the proof of \eqref{E:UPMUFIRSTTRANSPORT}.
	The main difference is that we use the identity
	$\partial_{\alpha} \Psi_{\beta} = \partial_{\beta} \Psi_{\alpha}$
	to rewrite 
	\begin{align}
		G_{\Lunit \Lunit}^{\lambda} \Rad \Psi_{\lambda}
		& = \upmu G_{\Lunit \Lunit}^{\lambda} \Radunit^a \partial_{\lambda} \Psi_a
			\\
		& = - G_{\Lunit \Lunit}^{\Lunit} \Radunit^a \Rad \Psi_a
			- \upmu G_{\Lunit \Lunit}^{\Lunit} \Radunit^a \Lunit \Psi_a
			- \upmu G_{\Lunit \Lunit}^{\Radunit} \Radunit^a \Lunit \Psi_a
			+ \upmu \angG_{\Lunit \Lunit}^{\#} \cdot \Radunit^a \angdiff \Psi_a.
			\notag
	\end{align}
\end{proof}

The $L^{\infty}$ estimates provided by the analogs of 
the estimate $|G_{(Frame)}^{\#}| = |\smoothfunction(\GdVar,\ginversesphere,\angdiff x^1,\angdiff x^2)| \lesssim 1$
(see Lemmas~\ref{L:SCHEMATICDEPENDENCEOFMANYTENSORFIELDS}
and \ref{L:POINTWISEFORRECTANGULARCOMPONENTSOFVECTORFIELDS}
and the $L^{\infty}$ estimates of Prop.~\ref{P:IMPROVEMENTOFAUX})
and Prop.~\ref{P:IMPROVEMENTOFAUX}
allow us to obtain the following bound for the term \eqref{E:LITTLEOMEGATRANSVERSALPART}:
$|\upomega^{(Tan-\vec{\Psi})}| = \mathcal{O}(\varepsilon)$.
Also using \eqref{E:NEWUPMUTRANSPORT}-\eqref{E:LITTLEOMEGATANGENTIALPART},
we see that the desired estimate \eqref{E:FIRSTSTATEMENTHARDSHARPRADPSIPOINTWISEESTIMATE} 
will follow once we show that
\begin{align} \label{E:SECONDEQUATIONNEWKINDOFSMALLNESSNEEDED}
		\left\|
			\Rad \Psi_0
			-
			\Radunit^a \Rad \Psi_a
		\right\|_{L^{\infty}(\Sigma_t^u)},
			\,
		\left\|
			\Rad \Psi_1
			+
			\Radunit^a \Rad \Psi_a
		\right\|_{L^{\infty}(\Sigma_t^u)},
			\,
		\left\|
			\Rad \Psi_2
		\right\|_{L^{\infty}(\Sigma_t^u)}
		& \lesssim \varepsilon.
	\end{align}
The proof of \eqref{E:SECONDEQUATIONNEWKINDOFSMALLNESSNEEDED} is not difficult.
Inequality \eqref{E:SECONDEQUATIONNEWKINDOFSMALLNESSNEEDED} for the second term
on the LHS follows from the identity (see \eqref{E:PERTURBEDPART})
\begin{align}
\Rad \Psi_1 
+
\Radunit^a \Rad \Psi_a 
= 
\Radunit_{(Small)}^1 \Rad \Psi_1 
+ \Radunit_{(Small)}^2 \Rad \Psi_2
\end{align}
and the bound $\| \Radunit_{(Small)}^i \|_{L^{\infty}(\Sigma_t^u)} \lesssim \varepsilon$
provided by the relation \eqref{E:LINEARLYSMALLSCALARSDEPENDINGONGOODVARIABLES},
the estimate \eqref{E:PURETANGENTIALLUNITAPPLIEDTOLISMALLANDLISMALLINFTYESTIMATE},
and the bootstrap assumptions \eqref{E:PSIFUNDAMENTALC0BOUNDBOOTSTRAP}
(the version for $\vec{\Psi}$).
The main idea of the proof of the other two estimates in \eqref{E:SECONDEQUATIONNEWKINDOFSMALLNESSNEEDED}
is to exploit the smallness \eqref{E:NEWPROBLEMADDITIONALSMALNESSESTIMATE}
and the following identities, which yield expressions for 
$\Rad \Psi_0 - \Radunit^a \Rad \Psi_a$
and
$\Rad \Psi_2$:
\begin{subequations}
\begin{align} \label{E:NEWSMALLNESSIDENTITYONE}
	\Lunit \Rad \Phi
	& = 
		\Rad \Psi_0 
		- \Radunit^a \Rad \Psi_a
		- (g^{-1})^{0a} \Rad \Psi_a
		+ (\Lunit \upmu) \Radunit^a \Psi_a
		+ \upmu (\Lunit \Radunit_{(Small)}^a) \Psi_a,
			\\
	\GeoAng \Rad \Phi
	& = 
		\Rad \Psi_2
		+ \GeoAng_{(Small)}^a \Rad \Psi_a
		+ (\GeoAng \upmu) \Radunit^a \Psi_a
 		+ \upmu (\GeoAng \Radunit_{(Small)}^a) \Psi_a.
		\label{E:NEWSMALLNESSIDENTITYTWO}
\end{align}
\end{subequations}
The identities \eqref{E:NEWSMALLNESSIDENTITYONE}-\eqref{E:NEWSMALLNESSIDENTITYTWO}
follow from the identity
$\partial_{\alpha} \Psi_{\beta} 
= 
\partial_{\beta} \Psi_{\alpha}
$,
\eqref{E:DOWNSTAIRSUPSTAIRSSRADUNITPLUSLUNITISAFUNCTIONOFPSI},
Def.~\ref{D:PERTURBEDPART},
and the fact that $\Lunit^0 = 1$.
From \eqref{E:LINEARLYSMALLSCALARSDEPENDINGONGOODVARIABLES},
the $L^{\infty}$ estimates of the analog of Prop.~\ref{P:IMPROVEMENTOFAUX} in the present context,
and the estimate \eqref{E:NEWPROBLEMADDITIONALSMALNESSESTIMATE},
we find that
$\| \Rad \Psi_0 - \Radunit^a \Rad \Psi_a \|_{L^{\infty}(\Sigma_t^u)}, 
	\, 
\| \Rad \Psi_2 \|_{L^{\infty}(\Sigma_t^u)}
\lesssim \varepsilon
$.
Also using the already proven estimate \eqref{E:SECONDEQUATIONNEWKINDOFSMALLNESSNEEDED} for the second term on
the LHS, we conclude the remaining two estimates stated in \eqref{E:SECONDEQUATIONNEWKINDOFSMALLNESSNEEDED}.

%

\section{Extending the Results to the Irrotational Euler Equations}
\label{A:EXTENSIONTOEULER}
\setcounter{equation}{0}
\numberwithin{equation}{section}
In this appendix, we sketch the minor changes needed to extend the shock formation results outlined in 
Appendix~\ref{A:EXTENDINGTONONGEOMETRICWAVEEQUATIONS} 
to the irrotational Euler equations of fluid mechanics in two spatial dimensions;
this is the content of Subsect.~\ref{SS:MASSAGINGEULER}.
Then, in Subsect.~\ref{SS:REMARKSONSMALLNESS},
we show that there exist initial data for the irrotational Euler equations
that verify the smallness-largeness hierarchy 
used in our proof of shock formation.

\subsection{Massaging the equations into the form of Appendix~\ref{A:EXTENDINGTONONGEOMETRICWAVEEQUATIONS}}
\label{SS:MASSAGINGEULER}
The necessary changes are all connected to normalization.
Under the assumption of irrotationality, 
the Euler equations reduce to a quasilinear wave equation 
for a potential function\footnote{In general, the potential function $\Phi$ can only be locally defined because
$\mathbb{R} \times \Sigma$ is not simply connected. 
However, the quasilinear wave equation for irrotational Euler flows is of the form
$(g^{-1})^{\alpha \beta}(\partial \Phi) \partial_{\alpha} \partial_{\beta} \Phi = 0$.
In particular, the equation depends only on the gradient of $\Phi$, 
which is ``globally'' defined throughout the maximal development of the data.
} 
$\Phi$ on the spacetime manifold 
$\mathbb{R} \times \Sigma$. The wave equation is the Euler-Lagrange equation 
(in particular it can be expressed in divergence form)
for a Lagrangian depending on $\partial \Phi$ 
that must satisfy various physical assumptions allowing for a fluid interpretation;
see \cite{dCsM2014} for the details in the case of the non-relativistic Euler equations and
\cite{dC2007} in the case of the (special) relativistic Euler equations.
A representative wave equation in the special relativistic case,
derivable from the Lagrangian\footnote{This Lagrangian corresponds to the fluid equation of state
$p = \rho/(2s+ 1)$, where $p$ is the pressure and $\rho$ is the proper energy density.}
$\mathscr{L} := [- (m^{-1})^{\kappa \lambda} \partial_{\kappa} \Phi \partial_{\lambda} \Phi]^{s+1}$,
is (see \cite{iRjS2013} for more details):
\begin{align} \label{E:MODELRELATIVISTICFLUIDWAVEQUATION}
\partial_{\alpha} 
\left\lbrace
	[- (m^{-1})^{\kappa \lambda} \partial_{\kappa} \Phi \partial_{\lambda} \Phi]^s
	(m^{-1})^{\alpha \beta}
	\partial_{\beta} \Phi
\right\rbrace
= 0,
\end{align}
where $s \in (0,\infty)$ is a constant
and $m_{\alpha \beta} = \mbox{\upshape diag}(-1,1,1)$ is the standard Minkowski metric.
The background solutions with perturbations that we are able to treat
correspond to constant solutions
with \emph{non-zero} energy density.\footnote{When the energy density vanishes, the wave equation becomes degenerate.}
In terms of the potential, these solutions are $\Phi = kt$, where $k > 0$ is a constant.
For some fluid wave equations, 
the values of $k$ that correspond to a physical fluid solution are restricted to a subset of $\mathbb{R}$;
this is not the case for equation \eqref{E:MODELRELATIVISTICFLUIDWAVEQUATION}.

Note that the spacetime metric corresponding to the background solution
is flat but typically not equal to $\mbox{\upshape diag}(-1,1,1)$.
We can remedy this by rescaling time. That is,
we can rescale the Minkowski time coordinate by $t \rightarrow \upalpha t$ 
(where the constant $\upalpha > 0$ generally depends on $k$ and the Lagrangian)
so that the metric corresponding to the quasilinear wave equation 
is equal to $\mbox{\upshape diag}(-1,1,1)$ for the background solution.
This is equivalent to choosing rescaled rectangular coordinates such that the
speed of sound (that is, the propagation speed) corresponding to the background solution is 
$1$. Note that after this rescaling,
the $00$ component of the tensorfield called ``$m$'' in \eqref{E:MODELRELATIVISTICFLUIDWAVEQUATION}
is no longer $-1$. The rescaling also changes $k$ to $\upalpha k$, 
but we will ignore that minor change here.
Moreover, in a slight abuse of notation, 
we also refer to the rescaled time variable as $x^0$ and/or $t$.
Having normalized the rectangular coordinates, 
we may now divide the wave equation by $-(g^{-1})^{00}$,
which allows us to assume that \eqref{E:APPENDIXGINVERSE00ISMINUSONE} 
holds. In total, we obtain a wave equation of the form \eqref{E:APPENDIXNONGEOMETRICWAVE}
verifying \eqref{E:APPENDIXGINVERSE00ISMINUSONE}.
For the rest of this appendix, we assume that this is the case.

We now define $\Psi_{\nu}$ and $\vec{\Psi}$ as in \eqref{E:APPENDIXPSICOMPONENTDEF},
except that we change the definition of $\Psi_0$ to $\Psi_0 := \partial_t \Phi - k$.
This is a good definition because for the kinds of perturbations of the background solutions 
that we consider, the (undifferentiated) $\Psi_{\nu}$ are small quantities.
The condition \eqref{E:APPENDIXMETRICEXPANDAROUNDMINKOWSKI} 
concerning the functional dependence of the metric on the wave variables
takes the following form in the present context:
\begin{align} \label{E:APPENDIXEULERMETRICEXPANDAROUNDMINKOWSKI}
	g_{\alpha \beta}(\vec{\Psi}) 
	& = m_{\alpha \beta}
		+ g_{\alpha \beta}^{(Small)}(\vec{\Psi}),
	\qquad g_{\alpha \beta}^{(Small)}(\vec{\Psi} = 0) = 0,
\end{align}
where
$m_{\alpha \beta} = \mbox{\upshape diag}(-1,1,1)$ in \eqref{E:APPENDIXEULERMETRICEXPANDAROUNDMINKOWSKI}.

To derive our main shock formation results, 
we again assume that \eqref{E:APPENDIXNULLCONDITIONFAILS} holds.
For all fluid Lagrangians in the regime of physically relevant $k$,
aside from one exceptional Lagrangian 
(mentioned in Footnote~\ref{FN:EXCEPTIONALLAGRANGIANS} on pg.~\pageref{FN:EXCEPTIONALLAGRANGIANS}),
it is possible to construct Minkowski-rectangular coordinates such
that \eqref{E:APPENDIXNULLCONDITIONFAILS} holds.
One can compute that the $\Psi_{\nu}$ verify the system
\eqref{E:PSIGEOMETRIC}.
We have thus massaged the wave equations of irrotational fluid mechanics 
into a form such that we can apply the shock formation proof outlined in 
Appendix~\ref{A:EXTENDINGTONONGEOMETRICWAVEEQUATIONS}.
We note in passing that for these wave equations,
the first product 
$G^{\mu \alpha \beta} \lbrace \cdots \rbrace$
on RHS~\eqref{E:APPENDIXCOVARIANTSYSTEMQUADRATICTERM}
vanishes. The vanishing is a consequence of the 
symmetry property $G^{\mu \alpha \beta} = G^{\beta \alpha \mu}$,
which holds for Euler-Lagrange equations
since $G^{\mu \alpha \beta}$ 
is proportional to the third partial derivative of
the Lagrangian with respect to its arguments
$\partial_{\mu} \Phi$,
$\partial_{\alpha} \Phi$,
$\partial_{\beta} \Phi$.

\subsection{The existence of data verifying the smallness assumptions}
\label{SS:REMARKSONSMALLNESS}
\setcounter{equation}{2}
In Subsect.~\ref{SS:NEWSMALLNESS},
we explained that in order to prove shock formation
for wave equations of the form
$(g^{-1})^{\alpha \beta}(\partial \Phi) \partial_{\alpha} \partial_{\beta} \Phi = 0$
using the framework of the present paper, we have to propagate the smallness
of various derivatives of $\Phi$ even though other derivatives are allowed to be large.
In this subsection, we study this 
smallness-largeness hierarchy at the level of the initial data 
in the case of the irrotational (special) relativistic Euler equations.
In particular, we show that there exist physically relevant initial data
exhibiting the desired size estimates. 
The point is that the desired smallness for the
appropriate higher derivatives of $\Phi$ is not immediate because 
some directional derivatives can be large.
By Cauchy stability,\footnote{Here we are referring to the continuous dependence of the solution on the initial data.} 
it suffices to exhibit plane symmetric data
(explained in the next paragraph)
verifying the desired smallness-largeness hierarchy.
It turns out that in plane symmetry,
the existence of such data is not difficult to see using Riemann invariants. 
In fact, the data-size assumptions
of Subsect.~\ref{SS:NEWSMALLNESS} can be realized at time $0$
by perturbations of plane symmetric data 
in which one Riemann invariant, denoted
by $\mathcal{R}_-$ below, completely vanishes, while the other one,
denoted by $\mathcal{R}_+$ below, 
is small with sufficiently large spatial derivatives.
The case of $\mathcal{R}_- \equiv 0$ corresponds to a simple outgoing
(that is, right-moving) plane wave solution.
In most of this subsection, we describe how to construct the Riemann invariants and how they
are related to other variables; at the end, we return to the issue of constructing
data verifying the desired smallness-largeness hierarchy.

In the analysis of this subsection, for simplicity, 
we restrict our attention to plane symmetric solutions to the relativistic Euler equations
on $\mathbb{R} \times \Sigma$, where (as in the rest of the article)
$\Sigma = \mathbb{R} \times \mathbb{T}$. 
Plane symmetric irrotational solutions are such that
the fluid potential function\footnote{In irrotational relativistic fluid mechanics,
all physical fluid variables are functions of $\partial \Phi$.} 
$\Phi$ is, relative to the rectangular coordinates, 
a function of only $t$ and $x^1$.
Plane symmetric solutions can of course be viewed as solutions on
$\mathbb{R} \times \mathbb{R}$.

We start by recalling some basic facts about irrotational special relativistic fluid mechanics.
The discussion in this paragraph is valid in all spatial dimensions.
Here we make many assertions without providing proofs; 
readers may consult \cite{dC2007} for more details.
The wave equations of irrotational fluid mechanics are 
Euler-Lagrange equations of the form
\[
	(g^{-1})^{\alpha \beta}(\partial \Phi) \partial_{\alpha} \partial_{\beta} \Phi = 0.
\]
The Lagrangian $\mathscr{L}$ may be identified with the fluid pressure $p$
and can be expressed as a function of 
$\upsigma$:
\begin{align} \label{E:IRROTATIONALLAG}
	p 
	=
	\mathscr{L}
	& := \mathscr{L}(\upsigma),
\end{align}
where
\begin{align} \label{E:RELENTHALPY}
	\upsigma
	& := - 
	(m^{-1})^{\alpha \beta}
	\partial_{\alpha} \Phi
	\partial_{\beta} \Phi
	> 0,
\end{align}
$m$ is the Minkowski metric,
and $\sqrt{\upsigma}$ is the enthalpy per particle.
The positivity assumption in \eqref{E:RELENTHALPY}
is a consequence of the timelike character of the fluid
velocity.\footnote{The fluid velocity $u^{\alpha}$
is equal to 
$
\displaystyle
- 
 \frac{(m^{-1})^{\alpha \beta} \partial_{\beta} \Phi}{\sqrt{\upsigma}}
$.}
Physically relevant Lagrangians satisfy 
various positivity assumptions ensuring, for example, 
that the pressure is positive,
the energy density is positive,
and that the speed of sound is real, positive, and less than the speed of light;
the following conditions ensure that the Lagrangian is physically relevant:
\begin{align}
	\mathscr{L}(\upsigma), 
		\,
	\frac{d \mathscr{L}}{d \upsigma},
		\,
	\frac{d}{d \upsigma}\left(\mathscr{L}/ \sqrt{\upsigma} \right),
		\,
	\frac{d^2\mathscr{L}}{d \upsigma^2} > 0. \label{E:FLUIDSINTERPRETATIONPOSITIVITY}
\end{align}
The acoustical metric $g$ and its inverse $g^{-1}$ can be expressed as
\begin{align}
	g_{\alpha \beta}(\partial \Phi)
	& = m_{\alpha \beta}
			+ 
			H 
			\partial_{\alpha} \Phi
			\partial_{\beta} \Phi,
			\label{E:ACOUSTICALMET} \\
	(g^{-1})^{\alpha \beta}(\partial \Phi)
	& = (m^{-1})^{\alpha \beta}
		- F
			(m^{-1})^{\alpha \kappa} (m^{-1})^{\beta \lambda} 
			\partial_{\kappa} \Phi \partial_{\lambda} \Phi,
			\label{E:INVERSEACOUSTICALMET} \\
	F= F(\upsigma) 
	& := \frac{2}{G} \frac{d G}{d \upsigma},
		\\
	G= G(\upsigma) 
	& := 2 \frac{d \mathcal{L}}{d \upsigma},
		\\
	H = H(\upsigma) 
	& := \frac{F}{1 + \upsigma F}.
\end{align}
The (generally non-constant) speed of sound is
\begin{align}
	c_s & = c_s (\upsigma) 
	= \sqrt{1 - \upsigma H}
	= \sqrt{\frac{1}{1 + \upsigma F}}.
		\label{E:CHRISTODOULOUSSPEED}
\end{align}
The assumptions in \eqref{E:FLUIDSINTERPRETATIONPOSITIVITY}
ensure in particular that $0 < c_s < 1$, where the speed of light is $1$.

In the remainder of this appendix, 
we will consider plane symmetric perturbations of the constant state
solution $\Phi = k t$, where $k > 0$ is a constant.
As in Subsect.~\ref{SS:MASSAGINGEULER},
we denote $\Psi_0 := \partial_t \Phi - k$ and $\Psi_1 = \partial_1 \Phi$.
We denote with a ``bar''
the value of a $\partial \Phi-$dependent variable evaluated at the constant
state solution. For example, $\bar{c}_s$ is the speed of sound
evaluated at the solution $\Phi = k t$. 
We now make the change of coordinates\footnote{In plane symmetry, the variable
$(x')^2$ does not play a role in the analysis.}
\begin{align} \label{E:RESCALEDCOORDINATES}
t' 
&:= \bar{c}_s t,
\qquad
(x')^i := x^i,
\qquad
i=1,2.
\end{align}
We refer to $(t',(x')^1,(x')^2)$ as the ``rescaled coordinates''
and we denote the corresponding partial
derivative vectorfields by
$\partial_{t'}$ and $\partial_{i'}$, $i=1,2$.
We set
\begin{align}
	\Psi_{0'} 
	& := \frac{\Psi_0}{\bar{c}_s},
	\qquad
	\Psi_{1'} = \Psi_1.
\end{align}
The rescaled coordinates are the exact analog of the rectangular coordinates
used in the bulk of the paper.
The quantities $\Psi_{0'}$ and $\Psi_{1'}$
are the exact analogs of the quantities 
$\lbrace \Psi_{\alpha} \rbrace$
appearing in Appendix~\ref{A:EXTENDINGTONONGEOMETRICWAVEEQUATIONS}.
Note that the metric 
denoted by $m$ in \eqref{E:ACOUSTICALMET} and \eqref{E:INVERSEACOUSTICALMET}
is equal to $\mbox{\upshape{diag}}(-\bar{c}_s^{-2},1,\cdots,1)$ in the rescaled coordinates.
Note also that in the rescaled coordinates, 
the constant state solution is 
$
\displaystyle
\Phi = k \frac{t'}{\bar{c}_s}
= k' t'
$,
where
\begin{align} \label{E:RESCALEDCONSTANT}
	k'
	& := \frac{k}{\bar{c}_s}.
\end{align}
We also set
\begin{align} \label{E:BACKGROUNDSOUNDSPEED1}
	c_s' 
	& = 
	\frac{c_s}{\bar{c}_s}.
\end{align} 
Note that $\bar{c}_s' = 1$. As is evident from the formulas
\eqref{E:CHARVECPLANESYMMETRY}, this 
condition implies that in the rescaled coordinates, 
the speed of sound of the constant state solution is precisely $1$;
in fact, this is the reason that we introduce the rescaled coordinates.
We may also rescale $g$ so that \eqref{E:APPENDIXGINVERSE00ISMINUSONE} holds
in the rescaled coordinates.

In plane symmetry, we can analyze solutions using the method of Riemann
invariants. The method is tied to the following vectorfields:
\begin{align} \label{E:CHARVECPLANESYMMETRY}
	\Lunit
	& = \partial_{t'}
			+
			\left(
			\frac{
				-\frac{
					\Psi_{1'}}{\bar{c}_s^2(\Psi_{0'} + k')} 
					+ 
					c_s'
					}
					{
					1 - \frac{\Psi_{1'}}{\Psi_{0'} + k'} c_s'
					}
			\right)
		\partial_{1'},
	\qquad 
	\uLunit
	=
			\partial_{t'}
			-
			\left(
			\frac{
				\frac{
					\Psi_{1'}}{\bar{c}_s^2(\Psi_{0'} + k')} 
					+ 
					c_s'
					}
					{
					1 + \frac{\Psi_{1'}}{\Psi_{0'} + k'} c_s'
					}
			\right)
		\partial_{1'}.
\end{align}
One may check that $\Lunit$ and $\uLunit$ are null: $g(\Lunit,\Lunit) = g(\uLunit,\uLunit) = 0$.
Having made the changes of variables and normalizations described above, 
we see that the restriction of $g$ to the $(t',(x')^1)$ plane, denoted by
$\widetilde{g}$, verifies
$
\widetilde{g}^{-1}
=
-
\frac{1}{2} 
\Lunit
\otimes
\uLunit
-
\frac{1}{2} 
\uLunit
\otimes
\Lunit
$.
The vectorfield $\Lunit$ is the analog of the vectorfield defined in 
\eqref{E:LUNITDEF}.
Moreover, the vectorfield
\begin{align} \label{E:RADUNITFORRELATIVSTICEULER}
	\Radunit 
	& := \frac{1}{2}
		\left\lbrace
			\uLunit
			-
			\Lunit
		\right\rbrace
\end{align}
is the analog of the one from Def.~\ref{D:RADANDXIDEFS}.
We also set
\begin{align} \label{E:UPMURESCALEDVECTORFIELDS}
	\uLgood 
	& := \upmu \uLunit,
	\qquad
	\Rad := \upmu \Radunit,
\end{align}
where $\upmu$ is defined by \eqref{E:FIRSTUPMU}, exactly as in the rest of the paper.
Note that $\Rad$ is the analog of the vectorfield defined in \eqref{E:RADDEF}.
Note also that
\begin{align} \label{E:ULGOODINTERMSOFLUNITANDRAD}
	\uLgood 
	& = \upmu \Lunit 
		+ 
		2 \Rad.
\end{align}

Making minor changes (including notational changes and normalization changes)
to the analysis presented in \cite{dCaL2016}, 
one finds that in plane symmetry,
the irrotational relativistic Euler wave equation 
is equivalent to the following system:
\begin{align} \label{E:RIEMANNINVARIANTEQUATIONS}
	\uLgood \mathcal{R}_- 
	& = 0,
	\qquad 
	\Lunit \mathcal{R}_+
	= 0,
\end{align}
where
$\mathcal{R}_-$ and $\mathcal{R}_+$ are Riemann invariants,
normalized so that 
\begin{align} \label{E:RIEMANNINVARIANTSINCONSTANTSTATE}
	\bar{\mathcal{R}}_-
	= 
	\bar{\mathcal{R}}_+
	= 0.
\end{align}
Specifically, we have
\begin{subequations}
\begin{align} \label{E:RMINUS}
	\mathcal{R}_- 
	& = 
		\frac{1}{\bar{c}_s} 
		\int \frac{1}{c_s' \sqrt{\upsigma}} \, d \sqrt{\upsigma} 
			+
			\frac{1}{2}
			\ln
			\left(
				\frac{\bar{c}_s + \frac{\Psi_{1'}}{\Psi_{0'} + k'}}{\bar{c}_s - \frac{\Psi_{1'}}{\Psi_{0'} + k'}}
			\right),
			 \\
	\mathcal{R}_+
	& = 
			\frac{1}{\bar{c}_s}
			\int \frac{1}{c_s' \sqrt{\upsigma}} \, d \sqrt{\upsigma}
			- 
			\frac{1}{2}
			\ln
			\left(
				\frac{\bar{c}_s + \frac{\Psi_{1'}}{\Psi_{0'} + k'}}{\bar{c}_s - \frac{\Psi_{1'}}{\Psi_{0'} + k'}}
			\right),
			\label{E:RPLUS}
\end{align}
\end{subequations}
where 
$c_s'$ is viewed as a function of $\sqrt{\upsigma}$ in the integrations
in \eqref{E:RMINUS}-\eqref{E:RPLUS}
and the constants of integration 
are chosen so that \eqref{E:RIEMANNINVARIANTSINCONSTANTSTATE} holds.
We can express $\Psi_{0'}$ and $\Psi_{1'}$ as follows:
\begin{subequations}
\begin{align} \label{E:PSI0INTERMSOFRIEMANNINVARIANTS}
	\Psi_{0'}
	& = 
			\frac{1}{\bar{c}_s}
			\sqrt{\upsigma} 
			\cosh
			\left(
			\frac{1}{2}
			\left(
				\mathcal{R}_+
				-
				\mathcal{R}_-
			\right)
			\right)
			-
			k',
			\\
	\Psi_{1'}
	& = \sqrt{\upsigma} 
			\sinh
			\left(
			\frac{1}{2}
			\left(
				\mathcal{R}_+
				-
				\mathcal{R}_-
			\right)
			\right),
			\label{E:PSI1INTERMSOFRIEMANNINVARIANTS}
\end{align}
\end{subequations}
where $\sqrt{\upsigma}$ can be viewed as a smooth function of
$\mathcal{R}_- + \mathcal{R}_+$.

We now achieve the main goal of this subsection:
explaining how to construct data so that the 
data-size assumptions stated in Subsect.~\ref{SS:NEWSMALLNESS} hold for $\Phi$.
As in the rest of the paper, we assume that the plane symmetric data 
$(\Phi|_{t'=0},\partial_{t'} \Phi|_{t'=0}) 
= (\mathring{\Phi},\mathring{\Phi}_0)$
are supported in the unit interval $[0,1]$ (which may be identified with $\Sigma_0^1$).
Note that the corresponding data 
$
(\mathcal{R}_-|_{t'=0},\mathcal{R}_+ |_{t'=0}) 
:= 
(\mathring{\mathcal{R}}_-,\mathring{\mathcal{R}}_+)
$
for the Riemann invariants are also supported in $[0,1]$.
We start by phrasing the data-size assumptions 
in terms of the Riemann invariants. Later,
we briefly overview how those conditions translate into the desired size assumptions 
for the data for $\Phi$. For brevity, we do not provide complete details here.

To proceed, we let $\mathring{\upepsilon}$ and $\mathring{\updelta}$
be data-size parameters satisfying the size
assumptions described in Subsect.~\ref{SS:SMALLNESSASSUMPTIONS}
(in particular, we assume that $\mathring{\upepsilon}$ is small relative to $\mathring{\updelta}^{-1}$).
We now simply take smooth data such that
$\mathring{\mathcal{R}}_- \equiv 0$,
such that
$\| \mathring{\mathcal{R}}_+ \|_{L^{\infty}(\Sigma_0^1)} = \mathcal{O}(\mathring{\upepsilon})$,
and such that
for\footnote{$M=1,2,3$ corresponds to our assumption that we control
up to three $\Rad$ derivatives of $\Psi$ in \eqref{E:PSIDATAASSUMPTIONS}.} 
$M=1,2,3$,
$\| \partial_{1'}^M \mathring{\mathcal{R}}_+ \|_{L^{\infty}(\Sigma_0^1)}$
is a relatively larger size $\mathcal{O}(\mathring{\updelta})$.
The smallness of $\mathring{\mathcal{R}}_-$ and 
$\mathring{\mathcal{R}}_+$ 
means that we are treating a perturbation of the 
data corresponding to the constant state $\Phi = k' t'$.
Moreover, with the help of equation \eqref{E:RIEMANNINVARIANTEQUATIONS}
and the commutation relation $[\Lunit, \Rad] = 0$
(which may be seen to be valid for the irrotational relativistic Euler equations
in plane symmetry),
we find that the higher derivatives 
of $\mathcal{R}_+$ with respect to $\Lunit$ and $\Rad$
\emph{completely vanish} as long as at least one $\Lunit$ differentiation is taken.
In addition, using the estimate $\upmu|_{t'=0} = 1 + \mathcal{O}(\mathring{\upepsilon})$,
which can be proved using arguments similar to the ones used 
in proving the first estimate stated in 
\eqref{E:UPMUDATATANGENTIALLINFINITYCONSEQUENCES},
and equations
\eqref{E:CHARVECPLANESYMMETRY}
and
\eqref{E:ULGOODINTERMSOFLUNITANDRAD},
we find that $\Rad|_{t=0} = -(1 + \mathcal{O}(\mathring{\upepsilon}))\partial_{1'}$.
Hence, we find that for $M=1,2,3$,
$\| \Rad^M \mathring{\mathcal{R}}_+ \|_{L^{\infty}(\Sigma_0^1)}$
is a relatively large size $\mathcal{O}(\mathring{\updelta})$.
Then using equations 
\eqref{E:PSI0INTERMSOFRIEMANNINVARIANTS} and \eqref{E:PSI1INTERMSOFRIEMANNINVARIANTS},
we may translate the $\mathring{\upepsilon}-\mathring{\updelta}$
hierarchy for the $\Lunit$ and $\Rad$ derivatives of $\mathcal{R}_-$ and $\mathcal{R}_+$
into a similar $\mathring{\upepsilon}-\mathring{\updelta}$ hierarchy
for the $\Lunit$ and $\Rad$ derivatives of $\Phi$, which yields
the desired data-size assumptions of Subsect.~\ref{SS:NEWSMALLNESS}.

We close this subsection by 
giving one concrete example illustrating the translation
mentioned at the end of the previous paragraph.
Specifically, we will show that
$
\| 
	\Rad (\Psi_{0'} + \Psi_{1'}) 
\|_{L^{\infty}(\Sigma_0^1)}
=
\mathcal{O}(\mathring{\upepsilon})
$.
As we will see, this estimate
is an easy consequence of our assumption\footnote{Actually, 
the $\mathcal{O}(\mathring{\upepsilon})-$smallness of $\Rad \mathring{\mathcal{R}}_-$ 
would suffice to obtain the desired bound.} 
that $\Rad \mathring{\mathcal{R}}_- = 0$.
The reason that this estimate is non-trivial is that
$
\| 
	\Rad \Psi_{0'}
\|_{L^{\infty}(\Sigma_0^1)}
$
and
$
\| 
	\Rad \Psi_{1'}
\|_{L^{\infty}(\Sigma_0^1)}
$
can be of a relatively large size
$\mathcal{O}(\mathring{\updelta})$.
The point of the estimate
$
\| 
	\Rad (\Psi_{0'} + \Psi_{1'}) 
\|_{L^{\infty}(\Sigma_0^1)}
=
\mathcal{O}(\mathring{\upepsilon})
$ 
is that, in view of the formula \eqref{E:NEWSMALLNESSIDENTITYONE}
and the fact that $\Radunit^1 = - 1 + \mathcal{O}(\mathring{\upepsilon})$,
it may be seen as a preliminary step 
(relevant for bounding the terms 
$
\Rad \Psi_0 
- \Radunit^a \Rad \Psi_a
$ 
on RHS~\eqref{E:NEWSMALLNESSIDENTITYONE})
in showing the desired data-size assumption
$
\| \Lunit \Rad \Phi \|_{L^{\infty}(\Sigma_0^1)}
=
\mathcal{O}(\mathring{\upepsilon})
$
stated in \eqref{E:NEWPROBLEMADDITIONALSMALLNESSASSUMPTIONS}.
To obtain the desired bound for
$\Rad (\Psi_{0'} + \Psi_{1'})$, 
we start by Taylor expanding RHS~\eqref{E:RMINUS} 
to first order around the constant state to obtain
\begin{align} \label{E:TAYLOREXPANDEDRMINUS}
	\mathcal{R}_-
	& = 
		\frac{1}{\bar{c}_s k}(\sqrt{\upsigma} - k)
		+ 
		\frac{\Psi_{1'}}{k}
		+ \mathcal{O}(|\mathcal{R}_-,\mathcal{R}_+|^2)
		= 
		\frac{\Psi_{0'}}{k}
		+ 
		\frac{\Psi_{1'}}{k}
		+ \mathcal{O}(|\mathcal{R}_-,\mathcal{R}_+|^2).
\end{align}
Applying $\Rad$ to both sides of 
\eqref{E:TAYLOREXPANDEDRMINUS}
and using our smallness assumptions
$\mathcal{R}_- \equiv 0$
and $\| \mathring{\mathcal{R}}_+ \|_{L^{\infty}(\Sigma_0^1)} = \mathcal{O}(\mathring{\upepsilon})$,
we conclude that
$
\| 
	\Rad (\Psi_{0'} + \Psi_{1'}) 
\|_{L^{\infty}(\Sigma_0^1)}
=
\mathcal{O}(\mathring{\upepsilon})
$
as desired.

\section{Notation}
\label{A:NOTATION}
In Appendix~\ref{A:NOTATION}, we collect some important notation 
and conventions that we use throughout the paper
so that the reader can refer to it as needed.

\subsection{Coordinates}
\label{SS:COORDINATES}
\begin{itemize}
	\item $(x^0,x^1,x^2)$ denote the rectangular spacetime coordinates.
	\item $(x^1,x^2$) denote the rectangular spatial coordinates.
	\item We often use the alternate notation $t = x^0$.
	\item $(t,u,\vartheta)$ are the geometric coordinates
		(where $t$ is the rectangular time coordinate, $u$ is the eikonal function,
		and $\vartheta$ is the geometric torus coordinate).
\end{itemize}

\subsection{Indices}
\label{SS:INDICES}
Lowercase Greek indices $\mu, \nu$, etc.\ correspond to components with respect
to the \emph{rectangular spacetime} coordinates $x^0,x^1,x^2$, and lowercase
Latin indices $i,j$, etc.\ correspond to components with respect
to the \emph{rectangular spatial} coordinates $x^1,x^2$. 
That is, lowercase Greek indices vary over $0,1,2$ and lowercase
Latin indices vary over $1,2$.
All lowercase Greek indices are lowered and raised with the spacetime metric
$g$ and its inverse $g^{-1}$, and \emph{not with the Minkowski metric}.
We use Einstein's summation convention in that repeated indices are summed 
over their respective ranges.

\subsection{Constants}
\label{S:CONSTANTS}
\begin{itemize}
\item $\mathring{\updelta}$ is the parameter corresponding to the initial size of the 
$\mathcal{P}_u-$transversal derivatives of the solution;
pg.~\pageref{E:PSIDATAASSUMPTIONS}.
\item $\mathring{\upepsilon}$ is a \emph{relatively small} parameter 
corresponding to the initial size of $\Psi$ and its derivatives involving at least one $\mathcal{P}_u-$tangential 
differentiation;
pg.~\pageref{E:PSIDATAASSUMPTIONS}.
We explain the kind of smallness that we impose on $\mathring{\upepsilon}$
in Subsect.~\ref{SS:SMALLNESSASSUMPTIONS}.
\item $\TranminusdatasizeWithFactor = \frac{1}{2} \sup_{\Sigma_0^1} \left[G_{\Lunit \Lunit} \Rad \Psi \right]_-$
is the key quantity that controls the blowup-time;
pg.~\pageref{E:CRITICALBLOWUPTIMEFACTOR}.
\item $C$ denotes a uniform constant that is free to vary from line to line.
\item The constants $C$ are allowed to depend on 
		the data-size parameters
		$\mathring{\updelta}$
		and 
		$\TranminusdatasizeWithFactor^{-1}$.
\item If we want to emphasize that the constant $C$ depends on an a quantity
	$Q$, then we use notation such as ``$C_Q$.''
\item We use the notation 
\[
	f_1 \lesssim f_2
\]
to indicate that there exists a uniform constant $C > 0$ such that
$f_1 \leq C f_2$. 
We sometimes use the alternate notation $\mathcal{O}(f_2)$ to denote a quantity $f_1$
that verifies $|f_1| \lesssim |f_2|$.
\item If we want to emphasize that the implicit constant $C$ depends on an a quantity
$Q$, then we use the alternate notation
\[
	f_1 \overset{Q}{\lesssim} f_2.
\]
\end{itemize}

\subsection{Spacetime subsets}
\label{SS:SPACETIMESUBSETS}
\begin{itemize}
	\item $\Sigma_{t'} := \lbrace (t,x^1,x^2) \in \mathbb{R} \times \mathbb{R} \times \mathbb{T} \ | \ t = t' \rbrace$; 
	pg.~\pageref{E:SIGMAT}.
	\item $\mathcal{P}_u =$ the outgoing null hypersurface equal to the corresponding level set of the eikonal function;
		pg.~\pageref{E:PUT}.
	\item $\Sigma_t^u =$ the portion of $\Sigma_t$ in between $\mathcal{P}_0$ and $\mathcal{P}_u$;
		pg.~\pageref{E:SIGMATU}.
	\item $\mathcal{P}_u^t =$ the portion of $\mathcal{P}_u$ in between $\Sigma_0$ and $\Sigma_t$;
		pg.~\pageref{E:PUT}.
	\item $\ell_{t',u'}=$ a topological one-dimensional torus equal to $\mathcal{P}_{u'}^{t'} \cap \Sigma_{t'}^{u'}$;
		pg.~\pageref{E:LTU}.
	\item $\mathcal{M}_{t,u} =$ the spacetime region trapped in between
		$\Sigma_0^u$, $\Sigma_t^u$, $\mathcal{P}_u^t$, and $\mathcal{P}_0^t$;
		pg.~\pageref{E:MTUDEF}.
\end{itemize}

\subsection{Metrics}
\label{SS:METRICS}
\begin{itemize}
	\item $g = g(\Psi)$ denotes the spacetime metric.
	\item Relative to rectangular coordinates,
		$g_{\mu \nu}(\Psi) = m_{\mu \nu} + g_{\mu \nu}^{(Small)}(\Psi)$,
		where $m_{\mu \nu} = \mbox{diag}(-1,1,1)$;
		pg.~\pageref{E:LITTLEGDECOMPOSED}.
	\item Relative to rectangular coordinates,
		$G_{\mu \nu}(\Psi) = \frac{d}{d \Psi} g_{\mu \nu}(\Psi)$
		and
		$G_{\mu \nu}'(\Psi) = \frac{d}{d \Psi} G_{\mu \nu}(\Psi)$;
 		pg.~\pageref{E:BIGGDEF}.
	\item $\gt$ denotes the first fundamental form of $\Sigma_t$, that is, $\gt_{ij} = g_{ij}$;
		pg.~\pageref{D:FIRSTFUND}.
	\item $\gt^{-1}$ denotes the inverse first fundamental form of $\Sigma_t$;
		pg.~\pageref{D:FIRSTFUND}. 
	\item $\gsphere$ denotes the first fundamental form of $\ell_{t,u}$;
		pg.~\pageref{D:FIRSTFUND}.
	\item $\gsphere^{-1}$ denotes the inverse first fundamental form of $\ell_{t,u}$;
		pg.~\pageref{D:FIRSTFUND}.
	\item $\gtancomp = \sqrt{\gsphere(\CoordAng,\CoordAng)}$ is a metric component;
		pg.~\pageref{E:METRICANGULARCOMPONENT}.
\end{itemize}

\subsection{Musical notation, contractions, and inner products}
\label{SS:MUSICALNOTATIONETC}
\begin{itemize}
	\item We denote the $\gsphere-$dual of an $\ell_{t,u}-$tangent one-form $\xi$
		by $\xi^{\#}$. 
		Similarly, if $Y$ is an $\ell_{t,u}-$tangent vector, then
		$Y_{\flat}$ denotes the $\gsphere-$dual of $Y$,
		which is an $\ell_{t,u}-$tangent covector. 
		Similarly, if $\xi$ is a symmetric type $\binom{0}{2}$ $\ell_{t,u}-$tangent tensor,
		then $\xi^{\#}$ denotes the type $\binom{1}{1}$ tensor that is
		$\gsphere-$dual to $\xi$, and $\xi^{\# \#}$
		denotes the type $\binom{2}{0}$ tensor that is
		$\gsphere-$dual to $\xi$.
		We use similar notation to denote the $\gsphere-$duals of general type 
		$\binom{m}{0}$ and type $\binom{0}{n}$ $\ell_{t,u}-$tangent tensors;
		pg.~\pageref{SS:NOTATION}.
	\item $g(X,Y) = g_{\alpha \beta} X^{\alpha} Y^{\beta}$ denotes the inner product
		of the vectors $X$ and $Y$ with respect to the metric $g$.
		Similarly, if $X$ and $Y$ are $\ell_{t,u}-$tangent,
		then $\gsphere(X,Y) = \gsphere_{ab} X^a Y^b$.
	\item $\cdot$ denotes the natural contraction between two tensors. 
		For example, if $\xi$ is a spacetime one-form and $V$ is a 
		spacetime vectorfield,
		then $\xi \cdot V := \xi_{\alpha} V^{\alpha}$.
		As a second example, if $T$ is a symmetric type $\binom{2}{0}$ 
		$\ell_{t,u}-$tangent tensorfield and $\xi$ is an $\ell_{t,u}-$tangent one-form, 
		then $(\angdiv T) \cdot \xi = (\angDarg{a} T^{ab}) \xi_b$;
		pg.~\pageref{SS:NOTATION}.
\item If $\xi$ is a one-form and $V$ is a vectorfield, then
	$\xi_V := \xi_{\alpha} V^{\alpha}$. 
	Similarly, if $W$ is a vectorfield, then
	$W_V := W_{\alpha} V^{\alpha} = g(W,V)$.
	We use similar notation when contracting higher-order tensorfields
	against vectorfields. Similarly, if $\Gamma_{\alpha \kappa \beta}$
	are the rectangular Christoffel symbols \eqref{E:CHRISTOFEELRECT}, then
	$\Gamma_{UVW} := U^{\alpha} V^{\kappa} W^{\beta} \Gamma_{\alpha \kappa \beta}$;
	pg.~\pageref{SS:NOTATION}.
	\item If $\xi$ is a symmetric type $\binom{0}{2}$ spacetime tensor and 
	$V$ is a vector, then $\xi_V$ is the one-form with rectangular components
	$(\xi_V)_{\nu} = \xi_{\alpha \nu}V^{\alpha}$;
	pg.~\pageref{E:TENSORVECTORANDSTUPROJECTED}.
\end{itemize}

\subsection{Tensor products and the trace of tensors}
\label{SS:TENSORPRODUCTS}
\begin{itemize}
	\item $(\xi \otimes \omega)_{\CoordAng \CoordAng} = \xi_{\CoordAng} \omega_{\CoordAng}$ 
		denotes the $\CoordAng \CoordAng$ component of the
		tensor product of the $\ell_{t,u}-$tangent one-forms $\xi$ and $\omega$.
	\item $\myspacetimetr \pi = (g^{-1})^{\alpha \beta} \pi_{\alpha \beta}$
		denotes the $g-$trace of the type $\binom{0}{2}$ spacetime tensor $\pi_{\mu \nu}$.
	\item $\mytr \xi = (\ginversesphere)^{\alpha \beta} \xi_{\alpha \beta}$ 
		denotes the $\gsphere-$trace 
		of the type $\binom{0}{2}$ $\ell_{t,u}-$tangent tensor $\xi$.
\end{itemize}

\subsection{Eikonal function quantities}
\label{SS:EIKONALFUNCTIONQUANTITIES}
	\begin{itemize}
		\item The eikonal function $u$ verifies the eikonal equation $(g^{-1})^{\alpha \beta}\partial_{\alpha} u \partial_{\beta} u = 0$
			and has the initial condition $u|_{t=0} = 1 - x^1$; 
			pg.~\pageref{E:INTROEIKONAL}.
		\item $\upmu = - \left\lbrace (g^{-1})^{\alpha \beta}\partial_{\alpha} u \partial_{\beta} t \right\rbrace^{-1}$ 
			denotes the inverse foliation density; pg.~\pageref{E:UPMUDEF}.
		\item $\Lgeo^{\nu} = - (g^{-1})^{\nu \alpha} \partial_{\alpha} u$ denotes the $\mathcal{P}_u-$tangent outgoing
			null geodesic vectorfield; pg.~\pageref{E:LGEOISGEODESIC}.
		\item $\Lunit^{\nu} = \upmu \Lgeo^{\nu}$ denotes a rescaled outgoing null vectorfield; pg.~\pageref{E:LUNITDEF}.
		\item $\Lunit_{(Small)}^i = \Lunit^i - \delta_1^i$ is the perturbed part of $\Lunit^i$;
			pg.~\pageref{E:PERTURBEDPART}.
		\item $\Radunit_{(Small)}^i = \Radunit^i + \delta_1^i$ is the perturbed part of $\Radunit^i$;
			pg.~\pageref{E:PERTURBEDPART}.
		\item $\upchi = \frac{1}{2} \angLie_{\Lunit} \gsphere$
			is the null second fundamental form of $\mathcal{P}_u$
			relative to $g$; pg.~\pageref{E:CHIDEF}.
\end{itemize}

\subsection{Additional \texorpdfstring{$\ell_{t,u}$}{toroidal} tensorfields related to the frame connection coefficients}
\label{SS:TOROIDALTENSORFIELDS}
\begin{itemize}
	\item $\upomega = \Lunit \upmu$; pg.~\pageref{E:UPMUFIRSTTRANSPORT}.
	\item $k = \frac{1}{2} \SigmatLie_{\Timenormal} \gt$ 
		is the second fundamental form of $\Sigma_t$ relative to $g$; pg.~\pageref{E:SECONDFUNDSIGMATDEF}.
	\item $\upzeta_{\CoordAng} = \angkdoublearg{\Radunit}{\CoordAng}
		= g(\D_{\CoordAng} \Lunit, \Radunit) $;
		pg.~\pageref{E:ZETADEF}.
	\item $\upzeta = \upmu^{-1} \upzeta^{(Trans-\Psi)} + \upzeta^{(Tan-\Psi)}$
		is a splitting of $\upzeta$;
		pg.~\pageref{E:ZETADECOMPOSED}.
	\item $\angk = \upmu^{-1} \angk^{(Trans-\Psi)} + \angk^{(Tan-\Psi)}$
		is a splitting of $\angk$;
		pg.~\pageref{E:ANGKDECOMPOSED}.
\end{itemize}

\subsection{Vectorfields}
\label{SS:VECTORFIELDS}
	\begin{itemize}
		\item $\Lunit$ is $\mathcal{P}_u-$tangent, outward pointing, 
			and verifies $g(\Lunit,\Lunit) = 0$ and $\Lunit t = 1$.
		\item $\Lunit = \frac{\partial}{\partial t}$ relative to the geometric coordinates;
			pg.~\pageref{E:LISDDT}.
		\item $\Rad$ is $\Sigma_t-$tangent, $\ell_{t,u}-$orthogonal,
			and verifies $g(\Rad,\Rad) = \upmu^2$, $\Rad u = 1$;
			pg.~\pageref{L:BASICPROPERTIESOFFRAME}.
		\item $\Rad = \frac{\partial}{\partial u} - \Xi$ relative to the geometric coordinates,
			where $\Xi$ is $\ell_{t,u}-$tangent;
			pg.~\pageref{E:RADSPLITINTOPARTTILAUANDXI}.
		\item $\Radunit = \upmu^{-1} \Rad$;
			pg.~\pageref{E:RADDEF}. 
		\item  
			$
			\displaystyle
			\CoordAng = \frac{\partial}{\partial \vartheta}
			$ 
			is the geometric torus coordinate partial derivative vectorfield;
			pg.~\pageref{D:ANGULARCOORDINATE}. 
		\item $\Timenormal = \Lunit + \Radunit$ is the 
			future-directed unit normal to $\Sigma_t$; 
			pg.~\pageref{E:TIMENORMAL}.
		\item $\lbrace \Lunit, \Rad, \CoordAng \rbrace$
			denotes the rescaled frame;
			pg.~\pageref{E:RESCALEDFRAME}.
		\item $\lbrace \Lunit, \Radunit, \CoordAng \rbrace$
			denotes the non-rescaled frame;
			pg.~\pageref{E:UNITFRAME}.
	\end{itemize}

\subsection{Projection operators and frame components}
\label{SS:PROJECTIONS}
	\begin{itemize}
		\item $\Sigmatproject$ denotes the type $\binom{1}{1}$ tensorfield that projects onto $\Sigma_t$;
			pg.~\pageref{E:SIGMATPROJECTION}.
		\item $\Lineproject$ denotes the type $\binom{1}{1}$ tensorfield that projects onto $\ell_{t,u}$;
			pg.~\pageref{E:LINEPROJECTION}.
		\item If $\xi$ is a spacetime tensor, then 
		$\angxi = \Lineproject \xi$
			is the projection of $\xi$ onto $\ell_{t,u}$;
			pg.~\pageref{E:TENSORSTUPROJECTED}.
		\item If $\xi$ is a type $\binom{0}{2}$ spacetime tensor, then 
			$\angxiarg{V} = \Lineproject (\xi_V)$;
			pg.~\pageref{E:TENSORSTUPROJECTED}.
		\item $G_{(Frame)} = 
			\left(
				G_{\Lunit \Lunit},
				G_{\Lunit \Radunit},
				G_{\Radunit \Radunit},
				\angGarg{\Lunit},
				\angGarg{\Radunit},
				\angG
			\right)$
			is the array of components of $G_{\mu \nu}$
			relative to the non-rescaled frame $\lbrace \Lunit, \Radunit, \CoordAng \rbrace$;
			pg.~\pageref{D:GFRAMEARRAYS}.
		\item $G_{(Frame)}' = 
			\left(
				G_{\Lunit \Lunit}',
				G_{\Lunit \Radunit}',
				G_{\Radunit \Radunit}',
				\angGprimearg{\Lunit},
				\angGprimearg{\Radunit},
				\angGprime
			\right)$
			is the array of components of $G_{\mu \nu}'$
			relative to the non-rescaled frame $\lbrace \Lunit, \Radunit, \CoordAng \rbrace$;
			pg.~\pageref{D:GFRAMEARRAYS}.
	\end{itemize}

\subsection{Arrays of solution variables and schematic functional dependence}
\label{SS:MOREARRAYS}
\begin{itemize}
	\item $\GdVar = \left(\Psi, \Lunit_{(Small)}^1, \Lunit_{(Small)}^2 \right)$;
		pg.~\pageref{E:GOODABBREIVATEDVARIABLES}.
	\item $\BadVar = \left(\Psi, \upmu - 1, \Lunit_{(Small)}^1, \Lunit_{(Small)}^2 \right)$;
		pg.~\pageref{E:BADABBREIVATEDVARIABLES}.
	\item $\smoothfunction(\xi_{(1)},\xi_{(2)},\cdots,\xi_{(m)})$ schematically denotes
		an expression depending smoothly on the $\ell_{t,u}-$tangent tensorfields $\xi_{(1)}, \xi_{(2)}, \cdots, \xi_{(m)}$;
		pg.~\pageref{R:SCHEMATICTENSORFIELDPRODUCTS}.
\end{itemize}

\subsection{Rescaled frame components of a vector}
\label{SS:FRAMECOMPONENTS}
	\begin{itemize}
		\item If $\mathscr{J}$ is a spacetime vector, then
			$\upmu \mathscr{J} 
			= 
			- \upmu \mathscr{J}_{\Lunit} \Lunit 
			- \mathscr{J}_{\Rad} \Lunit 
			- \mathscr{J}_{\Lunit} \Rad 
			+ \upmu \angJ$
		denotes its decomposition relative to the rescaled frame
		$\lbrace \Lunit, \Rad, \CoordAng \rbrace$,
		where
		$\mathscr{J}_{\Lunit} = \mathscr{J}^{\alpha} \Lunit_{\alpha}$,
		$\mathscr{J}_{\Rad} = \mathscr{J}^{\alpha} \Rad_{\alpha}$,
		and $\angJ = \Lineproject \mathscr{J}$;
		pg.~\pageref{L:DIVERGENCEFRAME}.
	\end{itemize}

\subsection{Energy-momentum tensorfield and multiplier vectorfields}
\begin{itemize}
	\item $\enmomtensor_{\mu \nu}[\Psi] = \D_{\mu} \Psi \D_{\nu} \Psi
	- \frac{1}{2} g_{\mu \nu} (g^{-1})^{\alpha \beta} \D_{\alpha} \Psi \D_{\beta} \Psi$
		denotes the energy-momentum tensorfield associated to $\Psi$;
		pg.~\pageref{E:ENERGYMOMENTUMTENSOR}.
	\item $\Mult = (1 + 2 \upmu) \Lunit + 2 \Rad$ denotes the timelike multiplier vectorfield;
		pg.~\pageref{E:DEFINITIONMULT}.
\end{itemize}

\subsection{Commutation vectorfields}
\begin{itemize}
	\item $\GeoAng_{(Flat)}$ denotes the $\Sigma_t-$tangent vectorfield
		with rectangular spatial components 
		$\GeoAng_{(Flat)}^i = \delta_2^i$;
		pg.~\pageref{E:GEOANGEUCLIDEAN}.
	\item $\GeoAng = \Lineproject \GeoAng_{(Flat)}$ 
		denotes the $\ell_{t,u}-$tangent commutation vectorfield;
		pg.~\pageref{E:GEOANGDEF}.
	\item $\GeoAng^i = \delta_2^i + \GeoAng_{(Small)}^i$
		is a splitting of $\GeoAng$ into
		$\GeoAng_{(Flat)}$ and a perturbation;
		pg.~\pageref{E:PERTURBEDPART}.
	\item $\Fullset = \lbrace \Lunit, \Rad, \GeoAng \rbrace$
		denotes the full set of commutation vectorfields;
		pg.~\pageref{E:COMMUTATIONVECTORFIELDS}.
	\item $\Tanset = \lbrace \Lunit, \GeoAng \rbrace$
		are the $\mathcal{P}_u-$tangent commutation vectorfields;
		pg.~\pageref{E:TANGENTIALCOMMUTATIONVECTORFIELDS}.
\end{itemize}

\subsection{Differential operators and commutator notation}
\label{SS:DIFFOPSANDCOMMUTATORNOTATION}
\begin{itemize}
	\item $\partial_{\mu}$ denotes the rectangular coordinate partial derivative 
	vectorfield $\displaystyle \frac{\partial}{\partial x^{\mu}}$.
	\item $\displaystyle \frac{\partial}{\partial t},\frac{\partial}{\partial u},\CoordAng = \frac{\partial}{\partial \vartheta}$ 
		denote the geometric coordinate partial derivative vectorfields.
	\item $V f = V^{\alpha} \partial_{\alpha} f$ 
		denotes the $V-$directional derivative of a function $f$.
	\item $d f$ denotes the standard differential of a function $f$ on spacetime.
	\item $\angdiff f = \Lineproject d f$, where $f$ is a function on spacetime and $\Lineproject$ denotes projection onto $\ell_{t,u}$;
		pg.~\pageref{D:ANGULARDIFFERENTIAL}. 
		Alternatively, $\angdiff f$ can be viewed as the inherent differential of a function $f$ defined on 
		$\ell_{t,u}$.
	\item $\D =$ Levi-Civita connection of $g$.
	\item $\angD =$ Levi-Civita connection of $\gsphere$.
	\item $\nabla = $ Levi-Civita connection of the Minkowski metric $m$.
	\item $\D_{XY}^2 = X^{\alpha} Y^{\beta} \D_{\alpha} \D_{\beta}$
		and similarly for other connections
		(contractions against $X$ and $Y$ are taken after the two covariant differentiations).
	\item $\angD^2$ denotes the second $\ell_{t,u}$ covariant derivative corresponding to $\gsphere$.
	\item $\angLap =  \mytr \angD^2 f$ denotes the covariant Laplacian on $\ell_{t,u}$ corresponding to $\gsphere$.
	\item If $\xi$ is an $\ell_{t,u}-$tangent one-form,
			then $\angdiv \xi$ is the scalar-valued function
			$\angdiv \xi := \ginversesphere \cdot \angD \xi$.
		Similarly, if $V$ is an $\ell_{t,u}-$tangent vectorfield,
			then $\angdiv V := \ginversesphere \cdot \angD V_{\flat}$,
			where $V_{\flat}$ is the one-form $\gsphere-$dual to $V$.
			If $\xi$ is a symmetric type $\binom{0}{2}$ 
		 $\ell_{t,u}-$tangent tensorfield, then
		 $\angdiv \xi$ is the $\ell_{t,u}-$tangent 
		 one-form $\angdiv \xi := \ginversesphere \cdot \angD \xi$,
		 where the two contraction indices in $\angD \xi$
		 correspond to the operator $\angD$ and the first index of $\xi$.
	\item $\Lie_V \xi$ denotes the Lie derivative of $\xi$ with respect to $V$;
		pg.~\pageref{E:LIEDERIVATIVE}.
	\item $[V,W] = \Lie_V W$ when $V$ and $W$ are vectorfields;
		pg.~\pageref{E:LIEDERIVATIVE}.
	\item More generally, if $P$ and $Q$ are two operators,
		then $[P,Q] = PQ - QP$ denotes their commutator.
	\item $\SigmatLie_V \xi = \Sigmatproject \Lie_V \xi$ is the $\Sigma_t-$projected 
		Lie derivative of $\xi$ with respect to $V$;
		pg.~\pageref{E:PROJECTIONS}.
	\item $\angLie_V \xi = \Lineproject \Lie_V \xi$ is the $\ell_{t,u}-$projected 
		Lie derivative of $\xi$ with respect to $V$;
		pg.~\pageref{E:PROJECTIONS}.
\end{itemize}

\subsection{Floor and ceiling functions and repeated differentiation}
\label{SS:FLOORANDCEILING}
\begin{itemize}
\item If $M$ is a non-negative integer,
	then $\lfloor M/2 \rfloor = M/2$ for $M$ even and 
	$\lfloor M/2 \rfloor = (M-1)/2$ for $M$ odd,
	while
	$\lceil M/2 \rceil = M/2$ for $M$ even and 
	$\lceil M/2 \rceil = (M+1)/2$ for $M$ odd.
\item We label the three vectorfields in 
	$\Fullset$ as follows: $Z_{(1)} = \Lunit, Z_{(2)} = \GeoAng, Z_{(3)} = \Rad$.
	Note that $\Tanset = \lbrace Z_{(1)}, Z_{(2)} \rbrace$.
\item If $\vec{I} = (\iota_1, \iota_2, \cdots, \iota_N)$ is a multi-index
		of order $|\vec{I}| := N$
		with $\iota_1, \iota_2, \cdots, \iota_N \in \lbrace 1,2,3 \rbrace$,
		then $\Fullset^{\vec{I}} := Z_{(\iota_1)} Z_{(\iota_2)} \cdots Z_{(\iota_N)}$ 
		denotes the corresponding $N^{th}$ order differential operator.
		We write $\Fullset^N$ rather than $\Fullset^{\vec{I}}$
		when we are not concerned with the structure of $\vec{I}$.
	\item Similarly, $\angLie_{\Fullset}^{\vec{I}} 
	:= \angLie_{Z_{(\iota_1)}} \angLie_{Z_{(\iota_2})} \cdots \angLie_{Z_{(\iota_N})}$
		denotes an $N^{th}$ order $\ell_{t,u}-$projected Lie derivative operator
		(see Def.~\ref{D:PROJECTEDLIE}),
		and we write $\angLie_{\Fullset}^N$
		when we are not concerned with the structure of $\vec{I}$.
	\item If $\vec{I} = (\iota_1, \iota_2, \cdots, \iota_N)$,
		then 
		$\vec{I}_1 + \vec{I}_2 = \vec{I}$ 
		means that
		$\vec{I}_1 = (\iota_{k_1}, \iota_{k_2}, \cdots, \iota_{k_m})$
		and
		$\vec{I}_2 = (\iota_{k_{m+1}}, \iota_{k_{m+2}}, \cdots, \iota_{k_N})$,
		where $1 \leq m \leq N$ and
		$k_1, k_2, \cdots, k_N$ is a permutation of 
		$1,2,\cdots,N$. 
	\item Sums such as $\vec{I}_1 + \vec{I}_2 + \cdots + \vec{I}_M = \vec{I}$
		have an analogous meaning.
	\item $\mathcal{P}_u-$tangent operators such as 
		$\Tanset^{\vec{I}}$ are defined analogously,  
		except in this case we clearly have
		$\iota_1, \iota_2, \cdots, \iota_N \in \lbrace 1,2 \rbrace$.
	\item $\Fullset^{N;M} f$ 
			denotes an arbitrary string of $N$ commutation
			vectorfields in $\Fullset$ (see \eqref{E:COMMUTATIONVECTORFIELDS})
			applied to $f$, where the string contains \emph{at most} $M$ factors of the $\mathcal{P}_u-$transversal
			vectorfield $\Rad$. 
		\item $\Tanset^N f$
			denotes an arbitrary string of $N$ commutation
			vectorfields in $\Tanset$ (see \eqref{E:TANGENTIALCOMMUTATIONVECTORFIELDS})
			applied to $f$.
		\item 
			For $N \geq 1$,
			$\Fullset_*^{N;M} f$
			denotes an arbitrary string of $N$ commutation
			vectorfields in $\Fullset$ 
			applied to $f$, where the string contains \emph{at least} one $\mathcal{P}_u-$tangent factor 
			and \emph{at most} $M$ factors of $\Rad$.
			We also set  $\Fullset_*^{0;0} f := f$.
		\item For $N \geq 1$,
					$\Tanset_*^N f$ 
					denotes an arbitrary string of $N$ commutation
					vectorfields in $\Tanset$ 
					applied to $f$, where the string contains
					\emph{at least one factor} of $\GeoAng$ or \emph{at least two factors} of $\Lunit$.
		\item For $\ell_{t,u}-$tangent tensorfields $\xi$, 
					we similarly define strings of $\ell_{t,u}-$projected Lie derivatives 
					such as $\angLie_{\Fullset}^{N;M} \xi$.
	\item $|\Fullset^{\leq N;M} f|$ is the sum over all terms of the form $|\Fullset^{N';M} f|$
			with $N' \leq N$ and $\Fullset^{N';M} f$ as defined above.
			When $N=M=1$, we sometimes write $|\Fullset^{\leq 1} f|$ instead of $|\Fullset^{\leq 1;1} f|$.
		\item $|\Fullset^{[1,N];M} f|$ is the sum over all terms of the form $|\Fullset^{N';M} f|$
			with $1 \leq N' \leq N$ and $\Fullset^{N';M} f$ as defined above.
		\item Sums such as 
			$|\Tanset_*^{[1,N]} f|$,
			$|\angLie_{\Fullset}^{\leq N;M} \xi|$,
			etc.,
			are defined analogously.
			We write $|\Tanset_* f|$ instead of $|\Tanset_*^{[1,1]} f|$.
			We also use the notation
			$|\Rad^{[1,N]} f| = |\Rad f| + |\Rad \Rad f| + \cdots + |\overbrace{\Rad \Rad \cdots \Rad}^{N \mbox{ \upshape copies}} f|$.
\end{itemize}
\

\subsection{Length, area, and volume forms}
\label{SS:LENGTHAREAANDVOLUMEFORM}
	\begin{itemize}
		\item $
		d \spherevol
		=
		d \argspherevol{(t,u,\vartheta)} = \gtancomp(t,u,\vartheta) d \vartheta$
		denotes the length form on $\ell_{t,u}$ induced by $g$;
		pg.~\pageref{E:RESCALEDVOLUMEFORMS}.
		\item $d \tvol = d \tvol(t,u',\vartheta) 
			= d \argspherevol{(t,u',\vartheta)} \, du'$
			denotes an area form on $\Sigma_t$;
			$\upmu \, d \tvol$ is the area form on $\Sigma_t$
				induced by $g$;
        pg.~\pageref{E:RESCALEDVOLUMEFORMS}.
		\item $d \conevol = d \conevol(t',u,\vartheta) = d \argspherevol{(t',u',\vartheta)} \, dt'$
			denotes an area form on $\mathcal{P}_u$;
			pg.~\pageref{E:RESCALEDVOLUMEFORMS}.
		\item $d \vol = d \vol(t',u',\vartheta) = d \argspherevol{(t',u',\vartheta)} \, du' \, dt'$
			denotes a volume form on $\mathcal{M}_{t,u}$;
			$\upmu \, d \vol$ is the volume form on $\mathcal{M}_{t,u}$
			induced by $g$;
			pg.~\pageref{E:RESCALEDVOLUMEFORMS}.
	\end{itemize}

\subsection{Norms}
\label{SS:MORENORMS}
	\begin{itemize}
		\item 
			$|\xi|^2 
				= 
				\gsphere_{\mu_1 \widetilde{\mu}_1} \cdots \gsphere_{\mu_m \widetilde{\mu}_m}
				(\ginversesphere)^{\nu_1 \widetilde{\nu}_1} \cdots (\ginversesphere)^{\nu_n \widetilde{\nu}_n}
				\xi_{\nu_1 \cdots \nu_n}^{\mu_1 \cdots \mu_m}
				\xi_{\widetilde{\nu}_1 \cdots \widetilde{\nu}_n}^{\widetilde{\mu}_1 \cdots \widetilde{\mu}_m}
			$ 
			denotes the square of the norm of the type $\binom{m}{n}$ $\ell_{t,u}$ tensor $\xi$;
			pg.~\pageref{E:POINTWISENORM}.
		\item $
		  \left\|
				\xi
			\right\|_{L^{\infty}(\ell_{t,u})}
				=
				\mbox{ess sup}_{\vartheta \in \mathbb{T}}
				|\xi|(t,u,\vartheta)
			$;
			pg.~\pageref{E:LINFTYNORMS}.
		\item 
		$ \left\|
				\xi
			\right\|_{L^{\infty}(\Sigma_t^u)}
			=
			\mbox{ess sup}_{(u',\vartheta) \in [0,u] \times \mathbb{T}}
				|\xi|(t,u',\vartheta)
		$;
		pg.~\pageref{E:LINFTYNORMS}.
		\item 
		$
		\left\|
				\xi
			\right\|_{L^{\infty}(\mathcal{P}_u^t)}
			=
			\mbox{ess sup}_{(t',\vartheta) \in [0,t] \times \mathbb{T}}
				|\xi|(t',u,\vartheta)
		$;
		pg.~\pageref{E:LINFTYNORMS}.
		\item $ \| f \|_{H_{\Euct}^N(\Sigma_0)}^2
					= \sum_{|\vec{I}| \leq N}
						\int_{\Sigma_0}
							(\partial_{\vec{I}} f)^2
						\, d^2 x
					$,
					where	$\partial_{\vec{I}}$ is a multi-indexed differential
					operator representing repeated differentiation with respect
					to the spatial coordinate partial derivatives
					and $d^2 x$
					is the area form corresponding to the standard Euclidean metric $\Euct$ on $\Sigma_0$;
					pg.~\pageref{R:SOBOLEVSPACESMULTIPLECOORDINATESYSTEMS}.
		\item $\| \xi \|_{L^2(\ell_{t,u})}^2 
			= \int_{\vartheta \in \mathbb{T}} |\xi|^2(t,u,\vartheta) d \spherevol
			= \int_{\ell_{t,u}} |\xi|^2 d \spherevol$;
			pg.~\pageref{E:L2NORMS}.
		\item $\| \xi \|_{L^2(\Sigma_t^u)}^2
			= \int_{u'=0}^u \int_{\vartheta \in \mathbb{T}} |\xi|^2(t,u',\vartheta) d \spherevol \, du'
			= \int_{\Sigma_t^u} |\xi|^2 \, d \tvol$;
			pg.~\pageref{E:L2NORMS}.
		\item $\| \xi \|_{L^2(\mathcal{P}_u^t)}^2
			= \int_{t'=0}^t \int_{\vartheta \in \mathbb{T}} |\xi|^2(t',u,\vartheta) d \spherevol \, dt'
			= \int_{\mathcal{P}_u^t} |\xi|^2 \, d \conevol$;
			pg.~\pageref{E:L2NORMS}.
		\item We use similar notation for the norms $\| \cdot \|_{L^2(\Omega)}$ 
		and
		$\| \cdot \|_{L^{\infty}(\Omega)}$
		of functions defined on subsets $\Omega$
		of $\ell_{t,u}$, $\Sigma_t^u$, or $\mathcal{P}_u^t$.
	\end{itemize}

\subsection{\texorpdfstring{$L^2-$}{Square-integral-}controlling quantities}
\begin{itemize}
	\item $\enzero[\Psi](t,u)$ denotes the energy of $\Psi$ along $\Sigma_t^u$ corresponding the multiplier $\Mult$;
		pg.~\pageref{E:ENERGYFLUX}.
	\item $\flzero[\Psi](t,u)$ denotes the null flux of $\Psi$ along $\mathcal{P}_u^t$ corresponding the multiplier $\Mult$;
		pg.~\pageref{E:ENERGYFLUX}.
	\item $\totTanmax{N}(t,u) = \max_{|\vec{I}| = N}
		\sup_{(t',u') \in [0,t] \times [0,u]} 
		\left\lbrace
			\enzero[\Tanset^{\vec{I}} \Psi](t',u')
			+ \flzero[\Tanset^{\vec{I}} \Psi](t',u')
		\right\rbrace$;
			pg.~\pageref{E:Q0TANNDEF}.
	\item $\totTanmax{[1,N]}(t,u) = \max_{1 \leq M \leq N} \totTanmax{M}(t,u)$;
		pg.~\pageref{E:MAXEDQ0TANLEQNDEF}.
	\item  $\coercivespacetime[\Psi](t,u) =
	 	\frac{1}{2}
	 	\int_{\mathcal{M}_{t,u}}
			[\Lunit \upmu]_-
			|\angdiff \Psi|^2
		\, d \vol$
		denotes the coercive spacetime integral associated to $\Psi$;
		pg.~\pageref{E:COERCIVESPACETIMEDEF}.
	\item $\coerciveTanspacetimemax{N}(t,u) = \max_{|\vec{I}| = N} \coercivespacetime[\Tanset^{\vec{I}} \Psi](t,u)$;
			pg.~\pageref{E:MAXEDCOERCIVESPACETIMEDEF}.
	\item $\coerciveTanspacetimemax{[1,N]}(t,u) = \max_{1 \leq M \leq N} \coerciveTanspacetimemax{M}(t,u)$;
			pg.~\pageref{E:MAXEDCOERCIVESPACETIMEDEF}.
\end{itemize}

\subsection{Modified quantities}
\begin{itemize}
	\item $\upchifullmodarg{\Tanset^N} = \upmu \Tanset^N \mytr \upchi + \Tanset^N \upchifullmodinhom$
		is the fully modified version of $\Tanset^N \mytr \upchi$; 
		pg.~\pageref{E:TRANSPORTRENORMALIZEDTRCHIJUNK}.
	\item $\upchipartialmodarg{\Tanset^N} = \Tanset^N \mytr \upchi + \upchipartialmodinhomarg{\Tanset^N}$
		is the partially modified version of $\Tanset^N \mytr \upchi$;
		pg.~\pageref{E:TRANSPORTPARTIALRENORMALIZEDTRCHIJUNK}.
\end{itemize}

\subsection{Curvature tensors}
\label{SS:CURVATURES}
\begin{itemize}
	\item $\Cur_{\mu \nu \alpha \beta} $ is the Riemann curvature tensor of $g$; 
		pg.~\pageref{E:SPACETIMERIEMANN}.
	\item $\Ric_{\alpha \beta} = (g^{-1})^{\kappa \lambda} \Cur_{\alpha \kappa \beta \lambda}$
		is the Ricci curvature tensor of $g$;
		pg.~\pageref{E:RICCIDEF}.
\end{itemize}

\subsection{Omission of the independent variables in some expressions}
\label{SS:OMISSION}
\begin{itemize}
	\item
	Many of our pointwise estimates are stated in the form
	\[
		|f_1| \lesssim h(t,u)|f_2|
	\]
	for some function $h$.
	Unless we otherwise indicate, it is understood that both $f_1$ and $f_2$
	are evaluated at the point with geometric coordinates
	$(t,u,\vartheta)$.
	\item Unless we otherwise indicate,
		in integrals $\int_{\ell_{t,u}} f \, d \spherevol$,
		the integrand $f$ and the length form $d \spherevol$ are viewed
		as functions of $(t,u,\vartheta)$ and $\vartheta$ is the integration variable.
	\item Unless we otherwise indicate,
		in integrals $\int_{\Sigma_t^u} f \, d \tvol$,
		the integrand $f$ and the area form $d \tvol$ are viewed
		as functions of $(t,u',\vartheta)$ and $(u',\vartheta)$ are the integration variables.
	\item Unless we otherwise indicate,
		in integrals $\int_{\mathcal{P}_u^t} f \, d \conevol$,
		the integrand $f$ and the area form $d \conevol$ are viewed
		as functions of $(t',u,\vartheta)$ and $(t',\vartheta)$ are the integration variables.
	\item Unless we otherwise indicate,
		in integrals $\int_{\mathcal{M}_{t,u}} f \, d \vol$,
		the integrand  $f$ and the volume form $d \vol$ are viewed
		as functions of $(t',u',\vartheta)$ and
		$(t',u',\vartheta)$ are the integration variables.
\end{itemize}

\bibliographystyle{amsalpha}
\bibliography{JBib}

\end{document}